\documentclass[11pt]{article}
\usepackage{blindtext}
\usepackage[a4paper, total={6in, 10in}]{geometry}

\usepackage[parfill]{parskip}    		
\usepackage{graphicx}				
\newcommand\smallO{
  \mathchoice
    {{\scriptstyle\mathcal{O}}}
    {{\scriptstyle\mathcal{O}}}
    {{\scriptscriptstyle\mathcal{O}}}
    {\scalebox{.5}{$\scriptscriptstyle\mathcal{O}$}}
  }

\usepackage{amsmath,amssymb,amsthm,amstext}
\numberwithin{equation}{section}

\DeclareMathOperator{\Tr}{Tr}

\usepackage{IEEEtrantools}

\usepackage[utf8]{inputenc}
\usepackage[english]{babel}
\usepackage{times}
\usepackage[Bjarne]{fncychap}

\usepackage{hyperref}
\usepackage{hypcap}
\theoremstyle{definition} \newtheorem{def2.1}{Definition}[section]
\theoremstyle{definition} \newtheorem{def2.2}[def2.1]{Definition}
\theoremstyle{definition} \newtheorem{def2.3}[def2.1]{Definition}
\theoremstyle{definition} \newtheorem{def2.4}[def2.1]{Definition}
\theoremstyle{definition} \newtheorem{def2.5}[def2.1]{Definition}
\theoremstyle{definition} \newtheorem{def2.6}[def2.1]{Definition}
\theoremstyle{definition} \newtheorem{def2.7}[def2.1]{Definition}
\theoremstyle{definition} \newtheorem{def2.8}[def2.1]{Definition}
\theoremstyle{definition} \newtheorem{def2.9}[def2.1]{Definition}
\theoremstyle{definition} \newtheorem{def2.10}[def2.1]{Definition}

\theoremstyle{definition} \newtheorem{theo2.1}{Theorem}[section]
\theoremstyle{definition} \newtheorem{theo2.2}[theo2.1]{Lemma}
\theoremstyle{definition} \newtheorem{theo2.3}[theo2.1]{Theorem}
\theoremstyle{definition} \newtheorem{theo2.4}[theo2.1]{Theorem}
\theoremstyle{definition} \newtheorem{theo2.5}[theo2.1]{Theorem}
\theoremstyle{definition} \newtheorem{theo2.6}[theo2.1]{Theorem}
\theoremstyle{definition} \newtheorem{theo2.7}[theo2.1]{Lemma}
\theoremstyle{definition} \newtheorem{theo2.8}[theo2.1]{Theorem}
\theoremstyle{definition} \newtheorem{theo2.9}[theo2.1]{Theorem}

\theoremstyle{definition} \newtheorem{def3.1}{Definition}[section]
\theoremstyle{definition} \newtheorem{def3.2}[def3.1]{Definition}

\theoremstyle{definition} \newtheorem{theo3.1}{Theorem}[section]
\theoremstyle{definition} \newtheorem{theo3.2}[theo3.1]{Theorem}
\theoremstyle{definition} \newtheorem{theo3.3}[theo3.1]{Theorem}
\theoremstyle{definition} \newtheorem{theo3.4}[theo3.1]{Theorem}
\theoremstyle{definition} \newtheorem{theo3.5}[theo3.1]{Theorem}
\theoremstyle{definition} \newtheorem{theo3.6}[theo3.1]{Theorem}
\theoremstyle{definition} \newtheorem{theo3.7}[theo3.1]{Theorem}
\theoremstyle{definition} \newtheorem{theo3.8}[theo3.1]{Theorem}
\theoremstyle{definition} \newtheorem{theo3.9}[theo3.1]{Theorem}
\theoremstyle{definition} \newtheorem{theo3.10}[theo3.1]{Theorem}
\theoremstyle{definition} \newtheorem{theo3.11}[theo3.1]{Theorem}
\theoremstyle{definition} \newtheorem{theo3.12}[theo3.1]{Theorem}
\theoremstyle{definition} \newtheorem{theo3.13}[theo3.1]{Theorem}
\theoremstyle{definition} \newtheorem{theo3.14}[theo3.1]{Theorem}

\theoremstyle{definition} \newtheorem{corol3.1}{Corollary}[section]
\theoremstyle{definition} \newtheorem{corol3.2}[corol3.1]{Corollary}
\theoremstyle{definition} \newtheorem{corol3.3}[corol3.1]{Corollary}
\theoremstyle{definition} \newtheorem{corol3.4}[corol3.1]{Corollary}

\theoremstyle{remark} \newtheorem{remark3.1}{Remark}[section]
\theoremstyle{remark} \newtheorem{remark3.2}[remark3.1]{Remark}
\theoremstyle{remark} \newtheorem{remark3.3}[remark3.1]{Remark}
\theoremstyle{remark} \newtheorem{remark3.4}[remark3.1]{Remark}
\theoremstyle{remark} \newtheorem{remark3.5}[remark3.1]{Remark}
\theoremstyle{remark} \newtheorem{remark3.6}[remark3.1]{Remark}
\theoremstyle{remark} \newtheorem{remark3.7}[remark3.1]{Remark}
\theoremstyle{remark} \newtheorem{remark3.8}[remark3.1]{Remark}
\theoremstyle{remark} \newtheorem{remark3.9}[remark3.1]{Remark}
\theoremstyle{remark} \newtheorem{remark3.10}[remark3.1]{Remark}

\theoremstyle{definition} \newtheorem{theo4.1}{Theorem}[section]
\theoremstyle{definition} \newtheorem{theo4.2}[theo4.1]{Theorem}
\theoremstyle{definition} \newtheorem{theo4.3}[theo4.1]{Theorem}

\theoremstyle{definition} \newtheorem{corol4.1}{Corollary}[section]
\theoremstyle{definition} \newtheorem{corol4.2}[corol4.1]{Corollary}

\theoremstyle{remark} \newtheorem{remark4.1}{Remark}[section]
\theoremstyle{remark} \newtheorem{remark4.2}[remark4.1]{Remark}
\theoremstyle{remark} \newtheorem{remark4.3}[remark4.1]{Remark}

\theoremstyle{definition} \newtheorem{def5.1}{Definition}[section]

\theoremstyle{definition} \newtheorem{theo5.1}{Theorem}[section]
\theoremstyle{definition} \newtheorem{theo5.2}[theo5.1]{Theorem}
\theoremstyle{definition} \newtheorem{theo5.3}[theo5.1]{Theorem}
\theoremstyle{definition} \newtheorem{theo5.4}[theo5.1]{Theorem}
\theoremstyle{definition} \newtheorem{theo5.5}[theo5.1]{Theorem}
\theoremstyle{definition} \newtheorem{theo5.6}[theo5.1]{Theorem}
\theoremstyle{definition} \newtheorem{theo5.7}[theo5.1]{Theorem}
\theoremstyle{definition} \newtheorem{theo5.8}[theo5.1]{Theorem}
\theoremstyle{definition} \newtheorem{theo5.9}[theo5.1]{Theorem}
\theoremstyle{definition} \newtheorem{theo5.10}[theo5.1]{Theorem}
\theoremstyle{definition} \newtheorem{theo5.11}[theo5.1]{Theorem}

\theoremstyle{definition} \newtheorem{corol5.1}{Corollary}[section]
\theoremstyle{definition} \newtheorem{corol5.2}[corol5.1]{Corollary}
\theoremstyle{definition} \newtheorem{corol5.3}[corol5.1]{Corollary}
\theoremstyle{definition} \newtheorem{corol5.4}[corol5.1]{Corollary}
\theoremstyle{definition} \newtheorem{corol5.5}[corol5.1]{Corollary}
\theoremstyle{definition} \newtheorem{corol5.6}[corol5.1]{Corollary}
\theoremstyle{definition} \newtheorem{corol5.7}[corol5.1]{Corollary}
\theoremstyle{definition} \newtheorem{corol5.8}[corol5.1]{Corollary}
\theoremstyle{definition} \newtheorem{corol5.9}[corol5.1]{Corollary}

\theoremstyle{remark} \newtheorem{remark5.1}{Remark}[section]
\theoremstyle{remark} \newtheorem{remark5.2}[remark5.1]{Remark}
\theoremstyle{remark} \newtheorem{remark5.3}[remark5.1]{Remark}
\theoremstyle{remark} \newtheorem{remark5.4}[remark5.1]{Remark}

\theoremstyle{plain}\newtheorem{gn_alg}{\underline{Gauss-Newton algorithms}}

\theoremstyle{plain}\newtheorem{lm_alg1}[gn_alg]{\underline{Levenberg-Marquardt algorithms}}
\theoremstyle{plain}\newtheorem{lm_alg2}[gn_alg]{\underline{Levenberg-Marquardt algorithms}}
\theoremstyle{plain}\newtheorem{lm_alg3}[gn_alg]{\underline{Levenberg-Marquardt algorithms}}

\theoremstyle{plain}\newtheorem{n_alg1}[gn_alg]{\underline{Newton algorithms}}
\theoremstyle{plain}\newtheorem{n_alg2}[gn_alg]{\underline{Newton algorithms}}
\theoremstyle{plain}\newtheorem{n_alg3}[gn_alg]{\underline{Newton algorithms}}

\theoremstyle{remark} \newtheorem{remark6.1}{Remark}[section]
\theoremstyle{remark} \newtheorem{remark6.2}[remark6.1]{Remark}
\theoremstyle{remark} \newtheorem{remark6.3}[remark6.1]{Remark}
\theoremstyle{remark} \newtheorem{remark6.4}[remark6.1]{Remark}

\title{Variable projection framework for the reduced-rank matrix approximation problem by weighted least-squares}
\author{Pascal Terray\footnote{Email: pascal.terray@locean.ipsl.fr} \footnote{Affiliation: Laboratoire d'Oc\'eanographie et du Climat: Exp\'erimentations et Approches Num\'eriques, Institut Pierre-Simon Laplace,
Sorbonne Universit\'e/CNRS/IRD/MNHN, Paris, France}}

\begin{document}

\maketitle

\pagestyle{plain}

\begin{abstract}
In this monograph, we review and develop  variable projection Gauss-Newton, Levenberg-Marquardt and Newton methods for the Weighted Low-Rank Approximation (WLRA) problem, which has now an increasing number of applications in many scientific fields. Particular attention is drawn at the robustness, efficiency and scalability of these variable projection second-order algorithms such that they can be used also on larger datasets now commonly found in many practical problems for which only first-order algorithms based on sequential repetitions of local optimization (e.g., majorization, Expectation-Maximization or alternating least-squares methods) or variations of gradient descent (e.g., conjugate, proximal or stochastic gradient descent methods), or hybrid algorithms from these two classes of methods, were only feasible due to their lower cost and memory requirement per iteration.

In parallel with this review of variable projection algorithms, we develop new formulae for the Jacobian and Hessian matrices involved in these variable projection methods and demonstrate their very specific properties such as the uniform rank deficiency of the Jacobian matrix or the rank deficiency of the Hessian matrix at the (local) minimizers of the cost function associated with the WLRA problem. These systematic deficiencies must be taken into account in any practical implementations of the algorithms. These different properties and the very particular geometry of the WLRA problem have not been well appreciated in the past and have been the main obstacles in the development of robust variable projection second-order algorithms for solving the WLRA problem.

In addition, we demonstrate that the variable projection framework gives original insights on the solvability, the landscape and the non-smoothness of the WLRA problem. It also helps to describe the tight links between previously unrelated methods, which have been proposed to solve it. Specifically, we illustrate the closed links between the variable projection framework and Riemannian optimization on the Grassmann manifold for the WLRA problem. We expect that software's developers and practitioners in different fields such as computer vision, signal processing, recommender systems, machine learning, multivariate statistics and geophysical sciences will benefit from the results in this monograph in order to devise more robust and accurate algorithms to solve the WLRA problem.
\end{abstract}

\newpage

\tableofcontents

\newpage

\section{Introduction} \label{intro:box}

Let $\mathbf{X}$ be a $p \times n$ real matrix and $\mathbf{W}$ be a $p \times n$ nonnegative real (weight) matrix (e.g., $\mathbf{W}_{ij} \ge 0$) associated with $\mathbf{X}$. This monograph is about the  Weighted Low-Rank Approximation (WLRA) problem:
\begin{equation} \label{eq:P0} \tag{P0}
\min_{\mathbf{Y} \in \mathbb{R}^{p \times n}_{\le k} } \, \quad\  \varphi( \mathbf{Y} ) =  \frac{1}{2}  \sum_{j=1}^n  \sum_{i=1}^p   { \mathbf{W}_{ij}.(  \mathbf{X}_{ij} -  \mathbf{Y}_{ij} )^{2}  }  = \frac{1}{2}  \Vert \sqrt{\mathbf{W}} \odot ( \mathbf{X} - \mathbf{Y} )  \Vert^{2}_{F} \ ,
\end{equation}
where $\mathbb{R}^{p \times n}_{\le k} = \big\{   \mathbf{Y} \in \mathbb{R}^{p \times n}  \text{ and } \emph{rank}( \mathbf{Y} ) \le{k}  \big\}$ and we assume that $k \le  \emph{rank}(   \mathbf{X} )  \le \emph{min}( p, n)$, $\odot$ denotes the Hadamard product (e.g., element-wise product) of two $p \times n$ matrices and $\Vert  \Vert_{F}$ is the Frobenius norm, i.e., the matrix norm induced by the standard inner product $  \langle \mathbf{Y}, \mathbf{X} \rangle = \emph{trace}( \mathbf{Y}^{T} \mathbf{X} ) $ on the Hilbert space of $p \times n$ real matrices. The factor  $\frac{1}{2}$ in the definition of $\varphi(.)$ has no effect on the minimizers of $\varphi(.)$, it is introduced only for notational convenience. Without this factor, we would have got an annoying factor of $2$ in many expressions of this monograph. Thus, a solution of the WLRA problem in its formulation~\eqref{eq:P0}, if it exists, is a $p \times n$ real matrix $\hat{\mathbf{X}}$ with $\emph{rank}(  \hat{\mathbf{X}} ) \le k \le  \emph{rank}(   \mathbf{X} )$. If $\ k =  \emph{rank}(  \hat{\mathbf{X}} ) = \emph{rank}(   \mathbf{X} )$ and $\mathbf{W}$  is a binary matrix, e.g., $\mathbf{W}_{ij}  \in  \big\{   0, 1  \big\} $, the WLRA problem is simply the so-called low-rank matrix completion problem (e.g., the problem of recovering matrices of low-rank when a large fraction of its elements are missing), which has been extensively studied in the past decades~\cite{NKS2019}\cite{DR2016}. In a slightly more general scenario, i.e., when $\ k =  \emph{rank}(  \hat{\mathbf{X}} ) <  \emph{rank}(  \mathbf{X} )$ and $\mathbf{W}$  is a binary matrix, a solution of the WLRA problem can be viewed as a robust generalisation of Principal Component Analysis (PCA) to incomplete, noisy or corrupted observations~\cite{J2002}\cite{IR2010}\cite{CLMW2011}\cite{NNSAJ2014}. In an even more general scenario when $\mathbf{W}$ is a general nonnegative matrix, a solution of the WLRA problem is very useful for denoising and revealing low-dimensional structures in incomplete and noisy datasets~\cite{PT1994}\cite{T2002}. Thus, in its general form, the WLRA problem can be considered as a robust generalization of (truncated) Singular Value Decomposition (SVD) analysis and extends significantly the usefulness and versatility of the classical low-rank approximation problem for many interesting applications arising from different fields including statistics~\cite{GZ1979}\cite{CLMW2011}\cite{TH2021}, computer vision~\cite{BF2005}\cite{C2008a}\cite{C2008b}\cite{HF2015}, machine learning for recommender systems~\cite{KBV2009}, signal processing and system identification\cite{MMH2003}\cite{MU2013}\cite{UM2014} and physical sciences~\cite{PT1994}\cite{T1995}\cite{T2002}\cite{BR2003}, to name a few.

Using general weights in the cost function $\varphi(.)$ allows us to take into account different confidence or sampling levels among the entries of the elements in $\mathbf{X}$ beyond the simple case of missing values, which corresponds to binary weights. As the error estimates of data are often widely varying, this is often better suited for many problems~\cite{PT1994}. Thus, weighted low-rank approximations of $\mathbf{X}$ can be used to deal with non-i.i.d. Gaussian noise in the data~\cite{T2002}\cite{MMH2003}\cite{C2008a} and to design robust versions of many multivariate statistical methods, which hinge on the classical low-rank matrix approximation in the Frobenius norm and are heavily used in data sciences. If the weight matrix $\mathbf{W}$ takes carefully into account the sampling properties of the dataset $\mathbf{X}$, the resulting weighted low-rank approximation $\hat{\mathbf{X}}$ is then defined to emphasize the better-observed aspects of the data~\cite{PT1994}\cite{T2002}. In other words, the nonnegative weights $\mathbf{W}_{ij}$ allow for a differential weighting of the accuracy of the measurements $\mathbf{X}_{ij}$ as well as for missing data if $\mathbf{W}_{ij} = 0$.  In particular, for the extreme case of zero sample size, an entry of the data matrix $\mathbf{X}$ should play no role in fitting the low-rank model; this can be done by assigning zero weight to such element of $\mathbf{X}$. 

Note, that we implicitly assume throughout the monograph that the weight matrix $\mathbf{W}$ is such that
\begin{equation*}
\sum_{i=1}^p { \mathbf{W}_{ij} } > 0  \text{ for } j = 1, \cdots, n  \text{ and  } \sum_{j=1}^n { \mathbf{W}_{ij} } > 0  \text{ for } i = 1, \cdots, p \ .
\end{equation*}
Stated more simply, these last two conditions imply that there is at least one nonzero weight in each column and row of $\mathbf{W}$ as otherwise the WLRA problem is not well-posed and tractable. Furthermore, we will demonstrate later that it is sometimes useful and necessary to impose stronger conditions on $\mathbf{W}$ such that each column and row of $\mathbf{W}$ have at least $k$ nonzero weights in order to avoid overfitting and obtain a meaningful approximate solution of the WLRA problem. In addition, as for the matrix completion problem, the WLRA problem may suffer from non-identifiability issues and is ill-posed without any incoherence type of conditions on the data matrix $\mathbf{X}$~\cite{CR2009}\cite{VMS2016}. As an illustration, with a sparse matrix $\mathbf{X}$, the matrix $\mathbf{W} \odot \mathbf{X}$ is likely to be a zero matrix if the number of non-zero weights $\mathbf{W}_{ij}$ is very small, and, in this case, the WLRA problem owns the zero matrix as a trivial solution, which obviously has no interest and is far from being optimal. To prevent this pathological case to occur, we need to impose some incoherent conditions on $\mathbf{X}$ with respect to the set of sparse matrices and assume that the number of samples is large enough, see~\cite{CR2009} or~\cite{VMS2016} for more formal definitions of these so-called low incoherence hypotheses, which provide reliable recoveries of the data matrix  $\mathbf{X}$ in the context of robust PCA, the matrix completion or WLRA problems.

If all the elements of $\mathbf{W}$ are all equal to $1$ (or more generally are all equal to a strictly positive real number), we have $\varphi( \mathbf{Y} ) =  \frac{1}{2}  \Vert  \mathbf{X} - \mathbf{Y}  \Vert^{2}_{F}$ up to a scaling constant, and this problem is well known and easily solved as the SVD theory provides the  best rank-$k$ approximation $\hat{\mathbf{X}}$ of a given  $p \times n$ real matrix $\mathbf{X}$ in terms of the Frobenius norm and also characterizes when this solution is unique or not  (see Theorem~\ref{theo2.1:box} below and~\cite{GVL1996} or~\cite{B2015} for details). Thus, in the simple case when all the elements of $\mathbf{W}$ are equal, but different from zero, it follows that once the SVD of  $\mathbf{X}$ is available, its best rank-$k$ approximation $\hat{\mathbf{X}}$ is readily computed. Moreover, if we are only interested in some $\hat{\mathbf{X}}$ with $k \ll \emph{min}( p, n)$, many less expensive alternatives than the computation of the complete SVD of $\mathbf{X}$ are available for computing $\hat{\mathbf{X}}$~\cite{GVL1996}\cite{STT2017}, including very fast and accurate randomized algorithms~\cite{HMT2011}\cite{LLSSKT2017}\cite{MDME2023}. Furthermore, under ideal conditions, i.e., $\mathbf{X}$ has no-missing values and the noise in all its elements can be modeled as zero-mean, independent and identically distributed (i.i.d.) Gaussian variables, the truncated SVD solution is the maximum likelihood solution and is, thus, the optimal one. However, this optimal property does not hold for non-i.i.d. Gaussian noise.

The more general case of uneven noisy observations (e.g., non-i.i.d. Gaussian noise) is in fact a particular instance of a WLRA problem in which we may assume that there is a ground truth low-rank matrix $\hat{\mathbf{X}}$, which we are trying to reconstruct and which is perturbed by non-i.i.d. Gaussian noise. Thus, implicit in the WLRA problem, is the statistical hypothesis that the input data consists of the observed (and also perturbed) data and weight matrices, $\mathbf{X}$ and $\mathbf{W}$, such that
\begin{equation}  \label{eq:stat_model}
\mathbf{X} = \mathbf{M} \odot ( \hat{\mathbf{X}} + \mathbf{E} ) \ ,
\end{equation}
where $\mathbf{M}$ is a boolean mask that indicates the observed elements of $\mathbf{X}$ (e.g., $\mathbf{M}_{ij} = 0$ if $\mathbf{W}_{ij} = 0$ and $\mathbf{M}_{ij} = 1$ otherwise), $\mathbf{E}$ is a noise matrix such that $\mathbf{E}_{ij} \sim \mathcal{N}( 0, \alpha_{ij}^{2})$ (e.g., $\mathbf{E}_{ij}$ is a Gaussian noise term) and $\mathbf{W}_{ij}$ is assumed to be modeled as a monotonically decreasing function of $\alpha_{ij}$, the noise level for each of the observed elements of $\mathbf{X}$. See~\cite{PT1994}\cite{T2002},~\cite{TH2021} and~\cite{C2008a} for examples, respectively, in the physical sciences, statistics and computer vision community on how such weight matrix $\mathbf{W}$ can be constructed in the case of non-i.i.d. Gaussian noise.

However, for a general choice of the weight matrix $\mathbf{W}$ and, even in the simple and very common case in which the weights are all $0$ or $1$ (e.g., the missing value or matrix completion problems~\cite{JHJ2009}\cite{IR2010}\cite{NKS2019}), the SVD of the masked observed matrix (e.g., set $\mathbf{X}_{ij} = 0$ if $\mathbf{W}_{ij} = 0$) may provide a useful and simple heuristic~\cite{MMW2021}, but does not give the desired closest fit to $\mathbf{X}$ in weighted 2-norm (or semi-norm if some weights are equal to zero) and the minimum of $\varphi(.)$. When general nonnegative weights are introduced, the problem of finding $\mathbf{Y} \in \mathbb{R}^{p \times n}_{\le k} $ so that $\varphi( \mathbf{Y} )$ is minimized is a NonLinear Least-Squares (NLLS) optimization problem in a finite-dimensional Hilbert space.

However, in its general setting, the WLRA problem is not convex (but only bi-convex) because of the nonconvexity and discontinuity of the rank function~\cite{HL2013}, has no closed-form solution because of the low-rank requirement, is known to be NP-hard~\cite{GG2011} and is, thus, not well understood~\cite{RSW2016}. Furthermore, for some matrices $\mathbf{X}$ and $\mathbf{W}$ and some integers $k$, the WLRA problem has no solution at all~\cite{GG2011} and, in other cases, the cost function $\varphi(.)$ may have several local minima~\cite{SJ2004}, a situation which can not occur in the classical low-rank approximation problem~\cite{SJ2004}\cite{H2010}. This hardness of the WLRA problem can be partly alleviated and some algorithms with provable guarantees have been proposed in the machine learning literature by making very strong assumptions such as incoherence of the ground truth low-rank matrix $\hat{\mathbf{X}}$, randomly sampled missing (or observed) entries in $\mathbf{X}$ or that the weight matrix $\mathbf{W}$ is spectrally closed to the all ones matrix~\cite{CR2009}\cite{KMO2010}\cite{JNS2013}\cite{BJ2014}\cite{LLR2016}. See the book by Vidal et al.~\cite{VMS2016} for a good introduction and discussion of these assumptions and the related algorithms in the case of binary weights (e.g., the matrix completion problem). However, for many applications these assumptions are unrealistic and violated, especially the assumption of randomly missing entries in the physical sciences, in which the statistical model~\eqref{eq:stat_model} is a more realistic framework, but does not provide any proven guarantees of success or provable time bounds for current WLRA algorithms. Taking into account this challenging background, the main objective of this paper is to discuss various efficient (pseudo) second-order iterative techniques for minimizing $\varphi(.)$, which exploit explicitly the separable properties of this cost function~\cite{R1974}\cite{W1976}\cite{RW1980}\cite{UM2014}\cite{HF2015}, and to show how to adapt standard NLLS algorithms to the special structure and geometry of the WLRA problem in its separable formulation.

The structure of the monograph is the following. Section~\ref{def:box} describes the notation used in the paper and gives an overview of some important definitions and preliminary results on linear algebra, multilinear algebra and  differentiation of vector and matrix functions relevant to the WLRA problem. In Section~\ref{seppb:box}, we study the geometry of the WLRA problem, the existence of solutions for it and we review several of its alternative formulations, which have been used in the literature, demonstrate their equivalence, which has not always been well appreciated in past studies and, finally, show that the WLRA problem can be reformulated as a separable NLLS problem~\cite{GP1973}\cite{RW1980}\cite{HPS2012}\cite{UM2014}\cite{HF2015}. This result was first used by Ruhe~\cite{R1974} for solving WLRA problems with binary weights and $k=1$ despite this separable formulation of the WLRA problem is often erroneously attributed to Wiberg~\cite{W1976} in the computer vision literature~\cite{SIR1995}\cite{OD2007}\cite{OYD2011}\cite{HF2015}. In fact, Wiberg~\cite{W1976} (who was a student of A. Ruhe) has extended the results of Ruhe~\cite{R1974} for an arbitrary integer $k$ and a slightly different component model specifically designed to the problem of estimating a principal components model when missing values are present in the data; see also~\cite{OD2007}\cite{VMS2016} for more details on this slightly different factor model used in~\cite{W1976}. Again in the computer vision literature, the separable NLLS algorithm originally proposed by Ruhe~\cite{R1974} and Wiberg~\cite{W1976} has been confused with the simplest Alternating Least-Squares (ALS) method~\cite{SIR1995}\cite{BF2005}\cite{VMS2016} as first noted by Okatani et al.~\cite{OD2007}. As a preamble to the variable projection algorithms, Section~\ref{nipals:box} gives a modern description of the block variant of this ALS method and its recent extensions.This ALS algorithm was perhaps the oldest and simplest method used to solve the WLRA problem in the statistical literature~\cite{W1966}\cite{WL1969}\cite{JHJ2009}\cite{GZ1979} and can be interpreted as a particular instance of the cyclic block-coordinate descent method for the WLRA problem. The ancestor of this block ALS algorithm is the Nonlinear Iterative PArtial Least Squares (NIPALS)  method devised originally by Wold and his collaborators for the missing value problem in PCA, i.e., in the case where the weight matrix is binary~\cite{W1966}\cite{WL1969}\cite{JHJ2009}\cite{IR2010}. Generalizations of the NIPALS algorithm to arbitrarily weighted least-squares have been first discussed in Gabriel and Zamir~\cite{GZ1979} and is now the  topic of many recent papers in different fields~\cite{SJ2004}\cite{BA2015}\cite{RSW2016}\cite{BWZ2019}\cite{TH2021}\cite{BRW2021}\cite{DLL2022}. However, nearly all the proposed algorithms dealing with general positive weights are first-order methods, excepted for the optimization approaches on the Grassmann manifold (e.g., the submanifold of fixed-rank matrices embedded in $\mathbb{R}^{p \times n}$) detailed in~\cite{MMH2003}\cite{C2008a}\cite{BA2015}. Section~\ref{varpro:box} is devoted to a detailed study of variable projection NLLS methods for solving the general WLRA problem, which use explicitly the separable property of this WLRA problem~\cite{GP1973}\cite{R1974}. Variable projection methods originate from numerical analysis and are efficient methods for solving separable NLLS problems in which some variables of the problem occur linearly and other nonlinearly; see Subsection~\ref{calculus:box} for a more formal definition. Explicit formulations of the gradient vector, Jacobian and Hessian matrices used in these variable projection second-order methods are given and their very specific mathematical properties are also derived in this Section~\ref{varpro:box}. Separable NLLS algorithms have a long history in applied mathematics and excellent reviews are offered in~\cite{RW1980}\cite{GP2003}\cite{HPS2012}. The closed relationships between the variable projection NLLS method and Riemannian optimization on the Grassmann manifold in the context of the WLRA problem are also explored in this Section~\ref{varpro:box}, extending and clarifying the results of Hong and Fitzgibbon~\cite{HF2015}\cite{HF2015b} who have focused on the binary weights case. Templates and implementation aspects of these variable projection second-order algorithms are detailed in Section~\ref{vpalg:box}. Finally, a summary of our contribution and perspectives for further advancing our understanding of the WLRA problem and methods for solving it are given in Section~\ref{conclu:box}.

\section{Definitions and preliminaries} \label{def:box}

We first collect in this section some basic notations, definitions and results concerning linear algebra, multilinear algebra, differentiation of vector and matrices and nonlinear optimization problems, which will be used frequently in the following sections.

Throughout this monograph, we have tried to adhere to the following conventions: bold capital letters will denote matrices and bold lower-case letters will indicate vectors. A lower-case letter in italic, but not in boldface, will indicate a scalar. The symbols $\mathbb{R}^{p}$ and $\mathbb{R}^{p \times n}$ denote, respectively, the linear spaces of the real $p$-vectors and of the real $p \times n$ matrices. In some occasions, the sizes of the vectors or the shapes of the matrices will be given as an upperscript. As an illustration, for  $a \in \mathbb{R}$, the symbols  $\mathbf{a}^p$ and $\mathbf{a}^{p \times n}$ represent, respectively, the $p$-vector and the $p \times n$ matrix composed of all $a$.  For $\mathbf{u} \in \mathbb{R}^{p}$, the symbol $\emph{diag}( \mathbf{u} )$ is used to represent a diagonal $p  \times p$ matrix with diagonal elements, $\lbrack \emph{diag}( \mathbf{u} ) \rbrack_{ii} = \mathbf{u}_{i}$ for  $i = 1, \cdots , p$. For any  $\mathbf{C}$ matrix, the symbol $\mathbf{C}_{.j}$ is used to represent the $j^{th}$ column vector of $\mathbf{C}$ and the symbol $\mathbf{C}_{i.}$  is used to represent the $i^{th}$  row vector of $\mathbf{C}$. The symbol $\mathbf{I}_p$  is used to denote the identity matrix of order $p$.

\subsection{Linear algebra} \label{lin_alg:box}

For a matrix $\mathbf{C} \in \mathbb{R}^{p \times n}$, we denote the transpose, the range and the null space of $\mathbf{C}$ by $\mathbf{C}^{T}$, $\emph{ran}(  \mathbf{C} )$  and $\emph{null}(  \mathbf{C} )$, respectively:
\begin{equation*}
\mathbf{C}^{T}_{ij} = \mathbf{C}_{ji}, \emph{ran}(  \mathbf{C} ) = \big\{ \mathbf{y} \in \mathbb{R}^{p}  \text{ }/\text{ }   \exists  \mathbf{x}   \in \mathbb{R}^{n}  \text{ with } \mathbf{y} = \mathbf{C}  \mathbf{x}  \big\}, \emph{null}(  \mathbf{C} ) = \big\{ \mathbf{x} \in \mathbb{R}^{n}  \text{ }/\text{ }  \mathbf{C}  \mathbf{x} = \mathbf{0}^p  \big\}.
\end{equation*}
$\emph{ran}(  \mathbf{C} )$ and $\emph{null}(  \mathbf{C} )$ are vector subspaces of $\mathbb{R}^{p}$ and $\mathbb{R}^{n}$, respectively. The rank of a matrix $\mathbf{C}$ is then defined by the dimension of the vector space $\emph{ran}(  \mathbf{C} )$, i.e.,  $\emph{rank}(  \mathbf{C} ) = \emph{dim}(  \emph{ran}( \mathbf{C} ) )$. Equivalently, the rank of a $p \times n$ matrix $\mathbf{C}$ can be defined as the smallest integer $k=\emph{rank}(  \mathbf{C} )$ such that it exists $\mathbf{A} \in \mathbb{R}^{p \times k}$ and $\mathbf{B} \in \mathbb{R}^{k \times n}$ such that $\mathbf{C}= \mathbf{A} \mathbf{B}$. From this definition, it is not difficult to show that $\emph{rank}(  \mathbf{C}^{T} ) =  \emph{rank}(  \mathbf{C})$. Then, it can been shown that
\begin{equation} \label{eq:rank}
\emph{dim} \big(  \emph{null}( \mathbf{C} )  \big) + \emph{dim} \big(  \emph{ran}( \mathbf{C} )  \big) = \emph{dim}  \big(  \emph{null}( \mathbf{C} )  \big) + \emph{rank}( \mathbf{C} ) = n \ ,
\end{equation}
which is known as the rank-nullity theorem or relationship, and also that  
\begin{equation} \label{eq:rank2}
\emph{rank}( \mathbf{A}\mathbf{B} ) \le min \big( \emph{rank}( \mathbf{A} ), \emph{rank}( \mathbf{B} )  \big) \ ,
\end{equation}
if the number of columns of $\mathbf{A}$  is equal to the number of rows of $\mathbf{B}$, and, finally, that
\begin{equation} \label{eq:rank3}
\emph{rank}( \mathbf{A} + \mathbf{B} )  \le \emph{rank}( \mathbf{A} ) + \emph{rank}( \mathbf{B} ) \ ,
\end{equation}
when $\mathbf{A}$ and $\mathbf{B}$ are matrices of the same dimensions. We further assume the following equalities
\begin{equation} \label{eq:null_ran}
\emph{null}( \mathbf{C}^{T} ) = \emph{ran}( \mathbf{C} )^\bot   \text{ and } \emph{ran}( \mathbf{C}^{T} ) = \emph{null}( \mathbf{C} )^\bot \ ,
\end{equation}
where  $\emph{ran}( \mathbf{C} )^\bot$ and $\emph{null}( \mathbf{C} )^\bot$ denote, respectively, the orthogonal complements of the range and null spaces of $\mathbf{C}$ with respect to the standard Euclidean inner products in $\mathbb{R}^{p}$ and $\mathbb{R}^{n}$, respectively.

We will use mostly the Euclidean norm for vectors and the Frobenius norm for matrices, i.e.,
\begin{equation} \label{eq:frob_norm}
 \Vert \mathbf{u}  \Vert_{2} =  \Big( \sum_{i=1}^p { \mathbf{u}^2_{i} } \Big)^{\frac{1}{2} } \text{ for }  \mathbf{u} \in \mathbb{R}^{p}  \text{  and  } \Vert \mathbf{C}  \Vert_{F} =  \Big( \sum_{i=1}^p { \sum_{j=1}^n { \mathbf{C}^2_{ij} } }  \Big)^{\frac{1}{2} } \text{ for }  \mathbf{C} \in \mathbb{R}^{p \times n} \ ,
\end{equation}
which are, respectively, associated to the Euclidean inner product in $\mathbb{R}^{p}$ 
\begin{equation}\label{eq:euclid_prod}
\langle \mathbf{u}  , \mathbf{v}  \rangle_{2} =   \sum_{i=1}^p { \mathbf{u}_{i} \mathbf{v}_{i} } \ ,
\end{equation}
and to the Frobenius inner product in $\mathbb{R}^{p \times n}$, defined for matrices $\mathbf{U}$  and $\mathbf{V}$ of identical sizes, by
\begin{equation}\label{eq:frob_prod}
\langle \mathbf{U} , \mathbf{V}  \rangle_{F} =  \Tr \big( \mathbf{U}^{T} \mathbf{V}  \big) = \sum_{i=1}^p { \sum_{j=1}^n { \mathbf{U}_{ij} \mathbf{V}_{ij} } } \ ,
\end{equation}
where for squared matrices $\Tr  \big( \mathbf{W} \big)  = \sum_{i=1}^p  { \mathbf{W}_{ii} }$. When we do not specify it, we implicitly mean these standard norms and inner products for vectors and matrices. Occasionally, especially in Section~\ref{seppb:box}, we will also use the spectral norm for matrices, which is the natural norm on the set of $p \times n$ real matrices induced by the Euclidean norm for vectors. For $\mathbf{C} \in \mathbb{R}^{p \times n}$, its spectral norm  $\Vert \mathbf{C}  \Vert_{S}$ can be computed as the squared root of the greatest eigenvalue of the matrix product $\mathbf{C}^{T}\mathbf{C}$~\cite{GVL1996}, i.e.,
\begin{equation}\label{eq:spect_norm}
\Vert \mathbf{C}  \Vert_{S} =  \max_{  \mathbf{x} \in  \mathbb{R}^{n} \text{ and }  \mathbf{x}  \ne \mathbf{0}^{n} }  \frac{\Vert \mathbf{C} \mathbf{x}  \Vert_{2} }{ \Vert  \mathbf{x}  \Vert_{2} } = ( \text{maximum eigenvalue of } \mathbf{C}^{T}\mathbf{C} )^{\frac{1}{2}} \ .
\end{equation}

A matrix $\mathbf{Q} \in \mathbb{R}^{p \times p}$ is said to be orthogonal if $\mathbf{Q}\mathbf{Q}^{T} = \mathbf{Q}^{T}\mathbf{Q} = \mathbf{I}_p$. It is easily verified that the product of two orthogonal matrices is also an orthogonal matrix. A matrix norm $\Vert   \Vert$ on $\mathbb{R}^{p \times n}$ is called unitarily invariant if $\Vert \mathbf{C} \Vert = \Vert \mathbf{Q}\mathbf{C}\mathbf{P} \Vert$ for all orthogonal matrices $\mathbf{Q}$ and $\mathbf{P}$ of order $p$ and $n$, respectively, and the Frobenius and spectral norms are unitarily invariant.

If $\mathbf{C} \in \mathbb{R}^{p \times n}$, then $\mathbf{C}^{+} \in \mathbb{R}^{n \times p}$ denotes the Moore-Penrose inverse (or pseudo-inverse) of $\mathbf{C}$ and is defined as the unique matrix which verifies the equalities
\begin{equation}\label{eq:ginv}
\mathbf{C}  \mathbf{C}^{+}  \mathbf{C}  = \mathbf{C}, \mathbf{C}^{+}  \mathbf{C}  \mathbf{C}^{+}  = \mathbf{C}^{+}, ( \mathbf{C}  \mathbf{C}^{+} )^{T} = \mathbf{C}  \mathbf{C}^{+} \text{ and }  ( \mathbf{C}^{+}  \mathbf{C} )^{T} = \mathbf{C}^{+}  \mathbf{C} \ .
\end{equation}
If $\mathbf{C}$ is of full column rank, it is easy to verify that
\begin{equation*}
\mathbf{C}^{+} = ( \mathbf{C}^{T} \mathbf{C} )^{-1} \mathbf{C}^{T} \ .
\end{equation*}
In addition, every matrix $\mathbf{C}^- \in \mathbb{R}^{n \times p}$ satisfying only the two equalities
\begin{equation}\label{eq:sginv}
\mathbf{C}  \mathbf{C}^-  \mathbf{C}  = \mathbf{C}  \text{ and }   ( \mathbf{C}  \mathbf{C}^- )^{T} = \mathbf{C}  \mathbf{C}^-
\end{equation}
is called a symmetric generalized inverse of $\mathbf{C}$.

An explicit formulation of the Moore-Penrose inverse $\mathbf{C}^{+}$ may be obtained with the help of the Singular Value Decomposition (SVD) of the matrix $\mathbf{C}$
\begin{equation}
\mathbf{C} = \mathbf{U}\Sigma\mathbf{V}^{T} \ ,
\end{equation}
where $\mathbf{U}$ and $\mathbf{V}$ are orthogonal matrices of order $p$ and $n$, respectively, and
\begin{equation*}
\Sigma = \left(
\begin{array}{ccccc}
\sigma_1 & 0              & \ldots    & 0     & 0      \\
0              & \sigma_2  & 0  &\ldots         & 0  \\
\vdots      & \ddots       & \ddots    & \ddots   & \vdots  \\
0            & \ldots         & 0           & \sigma_{n-1}   & 0  \\
0           & 0          &\ldots    & 0   & \sigma_n            \\
\vdots             &   \vdots & \ddots     &  \vdots  & \ddots    \\
0           &  0  & \ldots   &  0  &  0 
\end{array} \right) \ ,
\end{equation*}
where we have assumed for notational convenience that $p \ge n$. The existence of the SVD can be proved using the spectral theorem for symmetric matrices~\cite{GVL1996}\cite{B2015}. $\mathbf{U}$ and $\mathbf{V}$ consist of the orthonormal eigenvectors of $\mathbf{C} \mathbf{C}^{T}$ and of $\mathbf{C}^{T} \mathbf{C}$, respectively. $\mathbf{U}$ and $\mathbf{V}$ are called, respectively, the left and right singular vectors of $\mathbf{C}$. The diagonal elements of $\Sigma$ are called the singular values $\mathbf{C}$ and will always be taken to be nonnegative and ordered such that
\begin{equation*}
\sigma_1 \ge \sigma_2 \ge \cdots \sigma_{min(p,n)} \ge  0 \ .
\end{equation*}
These singular values are the non-negative square roots of the eigenvalues of $\mathbf{C}^{T} \mathbf{C}$ or $\mathbf{C} \mathbf{C}^{T}$. Then, in exact arithmetic, if $\emph{rank}( \mathbf{C} ) = k < n$, we have $\sigma_{k+1} = \sigma_{k+2} = \cdots = \sigma_{n} = 0$ and it is easy to verify that
\begin{equation}\label{eq:ginv_svd}
\mathbf{C}^{+} = \mathbf{V}\Lambda\mathbf{U}^{T} \ ,
\end{equation}
where $\Lambda$ is the $n \times p$ diagonal matrix with $\Lambda_{ii} = \sigma^{-1}_i$ for $i=1, \cdots, k$ and $\Lambda_{ii} = 0$ for $i=k+1, \cdots, n$. We also assume the following important property of the Moore-Penrose inverse $\mathbf{C}^{+}$ for all matrices $\mathbf{C}$: 
\begin{equation*}
\emph{null}( \mathbf{C}^{+} ) = \emph{null}( \mathbf{C}^{T} ) \ .
\end{equation*}

 A matrix $\mathbf{P}  \in \mathbb{R}^{p \times p}$ is an orthogonal projector if the following two conditions are satisfied:
 \begin{equation}\label{eq:projector}
\mathbf{P}\mathbf{P} = \mathbf{P} \text{ and } \mathbf{P}^{T} = \mathbf{P}
\end{equation}
and, given an orthogonal projector $\mathbf{P}  \in \mathbb{R}^{p \times p}$, its associated complementary projector is defined as $\mathbf{P}^{\bot} = \mathbf{I}_{p} - \mathbf{P}$ and is also an orthogonal projector. As for any matrix, an (orthogonal) projector $\mathbf{P}$ maps vectors into its range $\emph{ran}(  \mathbf{P} )$. However, an interesting and special property of any matrix $\mathbf{P}$ verifying $\mathbf{P}\mathbf{P} = \mathbf{P}$ is that it maps vectors of its range $\emph{ran}(  \mathbf{P} )$ to themselves. In addition, given an orthogonal projector $\mathbf{P}  \in \mathbb{R}^{p \times p}$ and a vector $\mathbf{x} \in \mathbb{R}^{p}$, the vector $\mathbf{P} \mathbf{x} \in \mathbb{R}^{p}$ uniquely solves the linear least-squares optimization problem
\begin{equation}\label{eq:projector2}
\mathbf{P} \mathbf{x}   = \text{Arg} \min_{  \mathbf{z} \in \emph{ran}(  \mathbf{P} )  }  \Vert \mathbf{x} - \mathbf{z} \Vert_{2} \ .
\end{equation}
In words, $\mathbf{P} \mathbf{x}$ is the unique closest point to $\mathbf{x}$ in $\emph{ran}(  \mathbf{P} )$. Note that $\emph{ran}(  \mathbf{P}^{\bot} ) = \emph{ran}(  \mathbf{P} )^{\bot}$, i.e., the range of $\mathbf{P}^{\bot}$ is the orthogonal complement of the range of  $\mathbf{P}$. Given a linear subspace $V$ of $\mathbb{R}^{p}$, we can decompose uniquely any vector $\mathbf{x} \in \mathbb{R}^{p}$ into the sum of one vector in $V$ and one vector  in $V^{\bot}$. This is easily verified as, given the (unique) orthogonal projector $\mathbf{P}$ onto $V$, we have immediately for any $\mathbf{x} \in \mathbb{R}^{p}$,
\begin{equation*}
\mathbf{x} = \mathbf{P} \mathbf{x} + (\mathbf{I}_{p} - \mathbf{P})\mathbf{x} = \mathbf{P} \mathbf{x} + \mathbf{P}^{\bot} \mathbf{x} \ ,
\end{equation*}
where $\mathbf{P} \mathbf{x}  \in  V$ and $ \mathbf{P}^{\bot} \mathbf{x}  \in V^{\bot}$. In such a case, we say that $\mathbb{R}^{p}$ is the direct sum of $V$ and $V^{\bot}$ and we write $\mathbb{R}^{p}=V \oplus V^{\bot}$. Finally, if the columns of $\mathbf{W}  \in \mathbb{R}^{p \times k}$ form an orthonormal basis of $V$, it is not difficult to verify that
\begin{equation*}
\mathbf{P} = \mathbf{W}\mathbf{W}^{T} \text{ and } \mathbf{P}^{\bot} = \mathbf{I}_{p} - \mathbf{W}\mathbf{W}^{T} \ .
\end{equation*}
Thus, provided that we have an orthonormal basis of $V$, we can also immediately project onto $V^{\bot}$ without constructing a basis for it. Furthermore, if we have such an orthonormal basis of $V$, we note that we have also a quick and efficient way of  applying orthogonal projectors to vectors as
\begin{equation*}
\mathbf{P}\mathbf{x} = \mathbf{W}( \mathbf{W}^{T} \mathbf{x} ) \text{ and } \mathbf{P}^{\bot} \mathbf{x}  = \mathbf{x}  - \mathbf{W}( \mathbf{W}^{T} \mathbf{x} ) \ .
\end{equation*}

The Moore-Penrose inverse and the SVD are also particularly useful to define and compute orthogonal projectors associated with the range of a matrix, especially if this matrix is rank deficient~\cite{GVL1996}\cite{B2015}. If the rank of the matrix $\mathbf{C} \in \mathbb{R}^{p \times n}$ is equal to $k$ (and looking at the distribution of the singular values of $\mathbf{C}$ is the best way to determine its numerical rank), the matrix
\begin{equation*}
\mathbf{P}_{\mathbf{C}} = \mathbf{C}\mathbf{C}^{+} = \mathbf{U}   \begin{bmatrix} \mathbf{I}_k   & \mathbf{0}^{k \times (p-k)}  \\  \mathbf{0}^{(p-k) \times k }  & \mathbf{0}^{(p-k) \times (p-k) }  \end{bmatrix}     \mathbf{U}^{T} \ ,
\end{equation*}
where $\mathbf{U}$ are the left singular vectors of $\mathbf{C}$, is the orthogonal projector onto $\emph{ran}(  \mathbf{C} )$. Furthermore, the matrix 
\begin{equation*}
\mathbf{P}_{\mathbf{C}}^{\bot} = \mathbf{I}_p - \mathbf{C}\mathbf{C}^{+} = \mathbf{U}   \begin{bmatrix}    \mathbf{0}^{k \times k}   & \mathbf{0}^{k \times (p-k)}  \\  \mathbf{0}^{(p-k) \times k }  & \mathbf{I}_{p-k}   \end{bmatrix}     \mathbf{U}^{T}
\end{equation*}
is the orthogonal projector onto the orthogonal complement of $\emph{ran}(  \mathbf{C} )$ (e.g., $\emph{ran}(  \mathbf{C} )^{\bot}$). It is easy to show that if $\mathbf{x} \in \mathbb{R}^{p}$ then $\mathbf{P}_{\mathbf{C}}\mathbf{x} \in \emph{ran}(  \mathbf{C} )$ and $\mathbf{P}_{\mathbf{C}}^{\bot}\mathbf{x} \in \emph{ran}(  \mathbf{C} )^{\bot}$. In the same conditions, the matrix
\begin{equation*}
\mathbf{P}_{\mathbf{C}^{T}} = \mathbf{C}^{+}\mathbf{C} = \mathbf{V}   \begin{bmatrix} \mathbf{I}_k   & \mathbf{0}^{k \times (n-k)}  \\  \mathbf{0}^{(n-k) \times k }  & \mathbf{0}^{(n-k) \times (n-k) }  \end{bmatrix}     \mathbf{V}^{T}
\end{equation*}
is the orthogonal projector onto the row space of $\mathbf{C}$, e.g., $\emph{ran}(  \mathbf{C}^{T} ) = \emph{null}(  \mathbf{C} )^\bot$ and the matrix
\begin{equation*}
\mathbf{P}_{\mathbf{C}^{T}}^{\bot} = \mathbf{I}_n - \mathbf{C}^{+}\mathbf{C} = \mathbf{V}   \begin{bmatrix}    \mathbf{0}^{k \times k}   & \mathbf{0}^{k \times (n-k)}  \\  \mathbf{0}^{(n-k) \times k }  & \mathbf{I}_{n-k}   \end{bmatrix}     \mathbf{V}^{T}
\end{equation*}
is the orthogonal projector onto $\emph{ran}(  \mathbf{C}^{T} )^\bot = \emph{null}(  \mathbf{C} )$. If $\emph{rank}(  \mathbf{C} ) = n \le p$, then $\mathbf{P}_{\mathbf{C}^{T}} = \mathbf{I}_{n}$ and $\mathbf{P}_{\mathbf{C}^{T}}^{\bot}$ is the $n \times n$ zero matrix $\mathbf{0}^{n \times n}$.

The Moore-Penrose inverse and SVD of a matrix are particularly useful for solving (rank-deficient) linear least-squares problems~\cite{LH1974}\cite{B2015}. For $\mathbf{C} \in \mathbb{R}^{p \times n}$ and $\mathbf{y} \in \mathbb{R}^{p}$, consider the linear least-squares problem
\begin{equation*}
\min_{\mathbf{x} \in \mathbb{R}^{n}}   \,  \Vert  \mathbf{y} - \mathbf{C}\mathbf{x} \Vert_{2} = \Vert  \mathbf{r} (\mathbf{x})  \Vert_{2} \text{ with }  \mathbf{r} (\mathbf{x}) = \mathbf{y} - \mathbf{C}\mathbf{x} \ .
\end{equation*}
The unique vector $\mathbf{\widehat{x}}$ of minimum Euclidean norm minimizing $\Vert  \mathbf{r}(\mathbf{x})  \Vert_{2}$ is given by $\mathbf{C}^{+} \mathbf{y}$ as $\mathbf{P}_{\mathbf{C}} \mathbf{y} = \mathbf{C}\mathbf{C}^{+}  \mathbf{y} = \mathbf{C} \mathbf{\widehat{x}}$ is the unique closest point to $\mathbf{y}$ in the range of $\mathbf{C}$. Even if $\mathbf{C}$ is of full column rank, in which case  $\mathbf{\widehat{x}} = ( \mathbf{C}^{T} \mathbf{C} )^{-1} \mathbf{C}^{T} \mathbf{y}$, we use the pseudo-inverse notation $\mathbf{C}^{+}$ for $( \mathbf{C}^{T} \mathbf{C} )^{-1} \mathbf{C}^{T}$ to indicate that the normal equations shall not be used to compute the solution of linear least-squares problems, especially if $\mathbf{C}$ is badly conditioned~\cite{LH1974}\cite{GVL1996}\cite{B2015}.

Note that linear least-squares problems can also be solved and orthogonal projectors be evaluated with the help of symmetric generalized inverses defined above~\cite{GP1976}\cite{HPS2012}. The advantage is that these symmetric generalized inverses can be computed much more cheaply than the SVD or the pseudo-inverse of $\mathbf{C}$  with the help of other matrix decompositions such as the standard QR decomposition with Column Pivoting (QRCP)~\cite{LH1974}\cite{GP1976}\cite{GVL1996}\cite{HPS2012}\cite{B2015}. According to this decomposition, there exist a $p \times p$ orthogonal matrix $\mathbf{Q}$ and a $n \times n$ permutation matrix $\mathbf{P}$ such that, for a given  $p \times n$ matrix $\mathbf{C}$ of rank $k$,
\begin{equation}\label{eq:qrcp}
\mathbf{Q}\mathbf{C}\mathbf{P} =  \begin{bmatrix} \mathbf{R}    &  \mathbf{S}   \\  \mathbf{0}^{(p-k) \times k }  & \mathbf{0}^{(p-k) \times (n-k) }  \end{bmatrix} \ ,
\end{equation}
where $\mathbf{R}$ is a $k \times k$ nonsingular upper triangular matrix (with diagonal elements of decreasing absolute magnitude) and $\mathbf{S}$ an $k \times (n-k)$ full matrix, which is vacuous if $k = n$. Several procedures are available to compute this QRCP, but the usual one is based on Householder transformations (e.g., elementary orthogonal reflectors), which are orthogonal matrices of the form
\begin{equation}\label{eq:hous_mat}
\mathbf{H}(i) = \mathbf{I}_{p} - 2 \mathbf{v}(i)  \mathbf{v}(i)^{T} \ ,
\end{equation}
where the $p$-vector $\mathbf{v}(i)$ has a 2-norm equal to one~\cite{LH1974}\cite{GVL1996}.  Premultiplication by $\mathbf{H}(i)$ is frequently used to zero out a sequence of entries in a given column $p$-vector. Thus, in order to compute the QR or QRCP decomposition, $\mathbf{C}$ is successively pre-multiplied by at most min$(n,p)$ Householder transformations $\mathbf{H}(i)$, permuting the columns of $\mathbf{C}$  if necessary (thus determining the permutation matrix $\mathbf{P}$). Also, the orthogonal matrix $\mathbf{Q}$ can be compactly stored (as only the vectors $\mathbf{v}(i)$ need to be stored) and explicitly computed as a product of $k$ elementary reflectors
\begin{equation*}
\mathbf{Q} = \mathbf{H}(k) \cdots  \mathbf{H}(2) \mathbf{H}(1) \ .
\end{equation*}
For more details concerning Householder transformations, see~\cite{LH1974}\cite{GVL1996}\cite{B2015}.

Furthermore, the rank $k$ of $\mathbf{C}$ can be efficiently estimated in an additional step from the upper triangular factor $\mathbf{R}$ computed during the QRCP, but we omit the details here~\cite{LH1974}\cite{GVL1996}\cite{HPS2012}\cite{B2015}. Note that the QRCP is not unique as the permutation matrix $\mathbf{P}$ is not unique. However, with the help of a QRCP of $\mathbf{C}$, the orthogonal projectors $\mathbf{P}_{\mathbf{C}} $ and $\mathbf{P}_{\mathbf{C}}^{\bot}$ can be efficiently computed as
\begin{equation} \label{eq:ginv_proj}
\mathbf{P}_{\mathbf{C}} =   \mathbf{C}\mathbf{C}^{-} = \mathbf{Q}^{T}   \begin{bmatrix} \mathbf{I}_k   & \mathbf{0}^{k \times (p-k)}  \\  \mathbf{0}^{(p-k) \times k }  & \mathbf{0}^{(p-k) \times (p-k) }  \end{bmatrix}     \mathbf{Q}
\end{equation}
and
\begin{equation}  \label{eq:ginv_proj_ortho}
\mathbf{P}_{\mathbf{C}}^{\bot} = \mathbf{I}_p -  \mathbf{C}\mathbf{C}^{-} = \mathbf{Q}^{T}   \begin{bmatrix}    \mathbf{0}^{k \times k}   & \mathbf{0}^{k \times (p-k)}  \\  \mathbf{0}^{(p-k) \times k }  & \mathbf{I}_{p-k}   \end{bmatrix}     \mathbf{Q} \ .
\end{equation}
Furthermore, a symmetric generalized inverse of $\mathbf{C}$ defined by the equations~\eqref{eq:sginv} can be represented as
\begin{equation}\label{eq:sginv_qrcp}
\mathbf{C}^- = \mathbf{P} \begin{bmatrix} \mathbf{R}^{-1}    &  \mathbf{0}^{ k \times (p-k) }     \\  \mathbf{0}^{(n-k) \times k }  & \mathbf{0}^{(n-k) \times (p-k) }  \end{bmatrix} \mathbf{Q} \ .
\end{equation}
Note that this particular symmetric generalized inverse also satisfies the additional equation
\begin{equation*}
\mathbf{C}^-  \mathbf{C}  \mathbf{C}^-  = \mathbf{C}^- \ .
\end{equation*}
Furthermore, if $k=n$ then $\mathbf{C}^{-} = \mathbf{C}^{+}$. Finally, the vector $\mathbf{C}^{-} \mathbf{y}$ is also a solution of the linear least-squares problem
\begin{equation*}
\min_{\mathbf{x} \in \mathbb{R}^{n}}   \,  \Vert  \mathbf{y} - \mathbf{C}\mathbf{x} \Vert_{2} = \Vert  \mathbf{r} (\mathbf{x})  \Vert_{2} \ ,
\end{equation*}
but not the solution of minimum Euclidean norm if $k<n$~\cite{LH1974}\cite{GP1976}\cite{HPS2012}\cite{B2015}. In other words, the pseudo-inverse $\mathbf{C}^{+}$ singles out the least-squares solution of minimum Euclidean length, which is not the case of $\mathbf{C}^{-}$.

If $k<n$, by applying  additional Householder  transformations (or, alternatively, Givens rotations) on the right of the QRCP to annihilate the submatrix $\mathbf{S}$, it is possible to obtain a Complete Orthogonal Decomposition (COD) of the matrix $\mathbf{C}$ of rank $k$ (see Chapter 5 of~\cite{GVL1996} or Theorems 29 and 30 of~\cite{HPS2012} and also~\cite{GP1976}). More precisely, by applying these additional Householder transformations, we obtain the following expression
\begin{equation*}
\mathbf{Q} \mathbf{C}\mathbf{P} =  \begin{bmatrix} \mathbf{R}    &  \mathbf{S}   \\  \mathbf{0}^{(p-k) \times k }  & \mathbf{0}^{(p-k) \times (n-k) }  \end{bmatrix} = \begin{bmatrix} \mathbf{T}    &  \mathbf{0}^{ k \times (n-k) }    \\  \mathbf{0}^{(p-k) \times k }  & \mathbf{0}^{(p-k) \times (n-k) }  \end{bmatrix} \mathbf{Z} \ ,
\end{equation*}
where $\mathbf{Z}$ is an $n \times n$ orthogonal matrix and is again implicitly represented by the product of elementary Householder matrices and $\mathbf{T}$ is an $k \times k$ upper triangular matrix of full rank (which is different from the triangular factor $\mathbf{R}$ in the QRCP). The COD of $\mathbf{C}$ is then defined as
\begin{equation}\label{eq:cod}
\mathbf{Q}\mathbf{C}\mathbf{O} =  \mathbf{Q}\mathbf{C} (\mathbf{P} \mathbf{Z}^{T}) =  \begin{bmatrix} \mathbf{T}    &  \mathbf{0}^{ k \times (n-k) }  \\  \mathbf{0}^{(p-k) \times k }  & \mathbf{0}^{(p-k) \times (n-k) }  \end{bmatrix} \ ,
\end{equation} 
where $\mathbf{Q}$  is the same $p \times p$ orthogonal matrix as in the QRCP, $\mathbf{T}$ is an $k \times k$ nonsingular upper triangular matrix and $\mathbf{O} = \mathbf{P} \mathbf{Z}^{T}$ is an $n \times n$ orthogonal matrix as the product of two orthogonal matrices. With the help of a COD of $\mathbf{C}$, its pseudo-inverse can be represented by
\begin{equation}\label{eq:ginv_cod}
\mathbf{C}^{+} = \mathbf{O} \begin{bmatrix} \mathbf{T}^{-1}    &  \mathbf{0}^{ k \times (p-k) }     \\  \mathbf{0}^{(n-k) \times k }  & \mathbf{0}^{(n-k) \times (p-k) }  \end{bmatrix} \mathbf{Q} \ .
\end{equation}
It is easily checked that this $n \times p$ matrix verifies the four equations~\eqref{eq:ginv} defining  the pseudo-inverse of $\mathbf{C}$ and since, for any matrix $\mathbf{C}$, there is only one matrix having these four properties, the above matrix is the pseudo-inverse of $\mathbf{C}$. This demonstrates that there is no need to compute a more costly SVD of $\mathbf{C}$ for this purpose. Importantly, with a COD, we also get the orthogonal projectors on the row space of $\mathbf{C}$ and its orthogonal complement as
\begin{equation*}
\mathbf{P}_{\mathbf{C}^{T}} =   \mathbf{C}^{+}\mathbf{C} = \mathbf{O}  \begin{bmatrix} \mathbf{I}_k   & \mathbf{0}^{k \times (n-k)}  \\  \mathbf{0}^{(n-k) \times k }  & \mathbf{0}^{(n-k) \times (n-k) }  \end{bmatrix}     \mathbf{O}^{T}
\end{equation*}
and
\begin{equation*}
\mathbf{P}_{\mathbf{C}^{T}}^{\bot} = \mathbf{I}_n -  \mathbf{C}^{+}\mathbf{C} = \mathbf{O}   \begin{bmatrix}    \mathbf{0}^{k \times k}   & \mathbf{0}^{k \times (n-k)}  \\  \mathbf{0}^{(n-k) \times k }  & \mathbf{I}_{n-k}   \end{bmatrix}     \mathbf{O}^{T} \ .
\end{equation*}
Finally, if we assume that $\mathbf{C}$ is of full column rank $k=n < p$, there is no need to compute a QRCP or COD of $\mathbf{C}$ to get the pseudo-inverse and the orthogonal projectors on the row or column spaces of $\mathbf{C}$ as a simple QR decomposition will do the job.

The SVD theory also provides the characterization of the best rank-$k$ approximation of a given $p \times n$ real matrix in terms of the Frobenius norm~\cite{GVL1996}. As the Frobenius norm is unitarily invariant, we first note that
\begin{equation} \label{eq:frob_norm2}
\Vert \mathbf{C}  \Vert_{F} = \Vert \Sigma  \Vert_{F} = \Big( \sum_{l=1} ^{min(p,n)} \sigma^2_l \Big)^{\frac{1}{2}} \ ,
\end{equation}
which shows that the Frobenius norm of a matrix is entirely defined by its singular values. From the SVD of a  $p \times n$ matrix $\mathbf{C} = \mathbf{U}\Sigma\mathbf{V}^{T}$, we can also obtain directly its spectral norm as
\begin{equation} \label{eq:spect_norm2}
\Vert \mathbf{C}  \Vert_{S} =  \sigma_1 = \mathbf{U}_{.1}^{T} \mathbf{C} \mathbf{V}_{.1} \ .
\end{equation}
Then, the following theorem is the reason for the importance of the SVD for applications involving low-rank approximation of matrices:
\begin{theo2.1} \label{theo2.1:box}
Let the SVD of $\mathbf{C} \in  \mathbb{R}^{p \times n}$ be $\mathbf{C} = \mathbf{U}\Sigma\mathbf{V}^{T}$ with $\sigma_1 \ge \sigma_2 \ge \cdots \ge \sigma_{min(p,n)}$. In addition, for $k$ such that $1 \le k \le min(p,n)$, defined the truncated SVD of $\mathbf{C}$ by
\begin{equation*}
\mathbf{C}_{k} =  \mathbf{U}_{k}\Sigma_{k}\mathbf{V}^{T}_{k} \ ,
\end{equation*}
where $\mathbf{U}_{k}$ and $\mathbf{V}_{k}$ are the submatrices formed by the $k$ first columns of $\mathbf{U}$ and $\mathbf{V}$, respectively, and $\Sigma_{k} = \emph{diag}( \lbrack \sigma_1, \cdots, \sigma_k \rbrack )$. Then, $\mathbf{C}_{k}$ provides a matrix of rank at most $k$ that is closest in Frobenius norm to $\mathbf{C}$ and this minimum distance is given by
\begin{equation*}
\Vert \mathbf{C} - \mathbf{C}_{k} \Vert_{F} = \min_{  \mathbf{B} \in  \mathbb{R}^{p  \times n } \text{ with }  \emph{rank}(  \mathbf{B} ) \le k }  \Vert \mathbf{C} - \mathbf{B} \Vert_{F} = \Big( \sum_{l=k+1} ^{min(p,n)} \sigma^2_l \Big)^{\frac{1}{2}} \ .
\end{equation*}
If $\sigma_k > \sigma_{k+1}$  or if $\sigma_k = 0$ then $\mathbf{C}_{k}$ is the unique best approximation of rank at most $k$ of $\mathbf{C}$.

$\Box$
\end{theo2.1}
Theorem~\ref{theo2.1:box} is often called the Eckart-Young Theorem and is in fact valid in any unitarily invariant norm, see~\cite{GVL1996}.

\subsection{Multilinear algebra} \label{multlin_alg:box}

In the next sections, we also need some operators and results from multilinear algebra~\cite{MN2019}. These tools will be particularly useful when we need to manipulate matrices as elements of a linear vector space and for computing derivatives of matrices (or matrix-matrix products) with respect to another matrix.

For any $\mathbf{C}$  and $\mathbf{D}$ matrices of the same dimensions, the expression $\mathbf{C}  \odot \mathbf{D}$ is used to mean the element-wise product of the $\mathbf{C}$  and $\mathbf{D}$ matrices (e.g., the Hadarmard product of two matrices):
\begin{equation}\label{eq:hadprod}
\big \lbrack \mathbf{C}  \odot \mathbf{D}  \big \rbrack_{ij} = \mathbf{C}_{ij} . \mathbf{D}_{ij} \ .
\end{equation}
The following property holds for matrices $\mathbf{B}$, $\mathbf{C}$ and $\mathbf{D}$ of the same shapes:
\begin{equation*}
\langle \mathbf{B} \odot  \mathbf{C} \; , \; \mathbf{D}  \rangle_{F} =  \langle  \mathbf{C} \; , \; \mathbf{B} \odot \mathbf{D}  \rangle_{F} \ .
\end{equation*}
Let $\mathbf{C} \in  \mathbb{R}^{q \times r}$ and $\mathbf{C}_{.j}$ denotes the $j^{th}$ column of $\mathbf{C}$, then the $\emph{vec}(.)$ function maps the $q \times r$ matrix $\mathbf{C}$ into a $q.r \times 1$ column vector by "stacking" the columns of $\mathbf{C}$ below one another
\begin{equation}\label{eq:vec}
\mathbf{C} \in  \mathbb{R}^{q \times r} \Longrightarrow \emph{vec}( \mathbf{C} ) = \begin{bmatrix} \mathbf{C}_{.1}  \\  \vdots   \\ \mathbf{C}_{.r} \end{bmatrix} \in  \mathbb{R}^{q.r} \ .
\end{equation}
The  $\emph{vec}(.)$  operator is an element of $\pounds ( \mathbb{R}^{q \times r}, \mathbb{R}^{q.r} )$, e.g., is a continuous linear mapping from $\mathbb{R}^{q \times r}$ into $\mathbb{R}^{q.r}$ and is also a bijection. The $\emph{mat}(.)$  operator is then the inverse mapping of  $\emph{vec}(.)$, which is a continuous linear bijection from $\mathbb{R}^{q.r}$ into $\mathbb{R}^{q \times r}$ such that
\begin{equation}\label{eq:mat}
\emph{mat} \big( \emph{vec}( \mathbf{C} ) \big) = \mathbf{C} \text{ , } \forall  \mathbf{C}  \in \mathbb{R}^{q \times r} \ .
\end{equation}
When it is not obvious from the context what is the shape of the image matrix for a given vector $\mathbf{c} \in \mathbb{R}^{q.r}$, we will use the notation $\emph{mat}_{q \times r}(.)$ instead.

A useful property involving the $\emph{vec}(.)$ and Hadamard operators is that the vectorized form of the Hadamard product of two matrices of the same dimensions can be written as a matrix-vector product
\begin{equation} \label{eq:vec_hadamard}
\emph{vec}( \mathbf{C} \odot \mathbf{D} ) = \emph{diag}\big( \emph{vec}( \mathbf{C} ) \big)  \emph{vec}( \mathbf{D} ) \ .
\end{equation}

Let further $\mathbf{D} \in  \mathbb{R}^{s \times t}$, then the Kronecker product $\mathbf{C} \otimes \mathbf{D}$ is the $q.s \times  r.t$ block matrix, whose $ij^{th}$ block is defined by
\begin{equation}\label{eq:kronprod}
\big \lbrack \mathbf{C} \otimes \mathbf{D} \big \rbrack^{ij} =   \mathbf{C}_{ij} \mathbf{D}   \text{ for } i = 1, \cdots , q \text{ and } j = 1, \cdots , r \ .
\end{equation}
The Kronecker product is a bilinear operator meaning that
\begin{align} \label{eq:bilin_kronprod}
( \mathbf{C}  + \mathbf{D} ) \otimes \mathbf{E}  & = ( \mathbf{C} \otimes \mathbf{E} )  + ( \mathbf{D} \otimes \mathbf{E} )  \ ,  \nonumber \\
\mathbf{E} \otimes  ( \mathbf{C}  + \mathbf{D} ) & = ( \mathbf{E} \otimes  \mathbf{C} ) + ( \mathbf{E} \otimes  \mathbf{D} )  \ ,   \\
\alpha ( \mathbf{C} \otimes \mathbf{D} )             & =  ( \alpha \mathbf{C} ) \otimes \mathbf{D} =  \mathbf{C} \otimes ( \alpha \mathbf{D} ) \ , \nonumber
\end{align}
where $\alpha \in \mathbb{R}$,  $\mathbf{E}$ is any matrix, and  $\mathbf{C}$ and $\mathbf{D}$ are two matrices of the same dimensions. We assume that the reader is familiar with the basic properties of Kronecker products (see Chapter 2 of~\cite{MN2019} for details). For easy reference, we only state the following relations for any matrices $\mathbf{C}$ and $\mathbf{D}$:
\begin{align} \label{eq:rank_kronprod}
( \mathbf{C} \otimes \mathbf{D} )^{T}                  & = \mathbf{C}^{T} \otimes \mathbf{D}^{T}  \nonumber \\
\emph{rank} ( \mathbf{C} \otimes \mathbf{D} ) & = \emph{rank} ( \mathbf{C} ) . \emph{rank} ( \mathbf{D} ) \ ;
\end{align}
for partitioned matrices:
\begin{equation} \label{eq:partmat_kronprod}
\begin{bmatrix} \mathbf{C}_1 &  \mathbf{C}_2  \end{bmatrix} \otimes \mathbf{D} = \begin{bmatrix} \mathbf{C}_1 \otimes \mathbf{D} & \mathbf{C}_2 \otimes \mathbf{D} \end{bmatrix} \ ;
\end{equation}
and for conforming matrices $\mathbf{C}$, $\mathbf{D}$, $\mathbf{E}$ and $\mathbf{F}$:
\begin{align}  \label{eq:mulmat_kronprod}
( \mathbf{C} \otimes \mathbf{D} )( \mathbf{E} \otimes \mathbf{F} ) & = \mathbf{C}\mathbf{E} \otimes \mathbf{D}\mathbf{F}  \ , \\
\emph{vec}( \mathbf{C} \mathbf{D} \mathbf{E} ) & = (  \mathbf{E}^{T}  \otimes \mathbf{C} ) \emph{vec}( \mathbf{D} ) \  .  \nonumber
\end{align}
This last equality is particularly useful to rearrange a matrix-matrix product as a simple matrix-vector product:
\begin{align} \label{eq:vec_kronprod}
\emph{vec}( \mathbf{C} \mathbf{D} ) & =  \emph{vec}( \mathbf{C} \mathbf{D} \mathbf{I}) = (  \mathbf{I}  \otimes \mathbf{C} ) \emph{vec}( \mathbf{D} )   \ ,  \\
\emph{vec}( \mathbf{C} \mathbf{D} ) & =  \emph{vec}( \mathbf{I} \mathbf{C} \mathbf{D} ) = (  \mathbf{D}^{T} \otimes  \mathbf{I} ) \emph{vec}( \mathbf{C} ) \ , \nonumber
\end{align}
where $\mathbf{I}$ is the identity matrix of appropriate order. These two relationships illustrate that we can evaluate the derivative of a matrix-matrix product with respect to one of the matrices by reshaping  the different matrices as vectors and computing the Jacobian matrix (see Subsection~\ref{calculus:box} for more details). Thus, we will use these two relations very frequently in the next sections, without explicit citation, when we need to compute derivatives of some matrices.

Let $\mathbf{X} \in  \mathbb{R}^{p \times n}$. We now introduce the $p.n \times p.n$ permutation matrix $\mathbf{K}_{(p,n)}$ uniquely defined by the relation
\begin{equation}\label{eq:commat}
\mathbf{K}_{(p,n)} \emph{vec}( \mathbf{X} ) = \emph{vec}( \mathbf{X}^{T} ) \ .
\end{equation}
This permutation matrix is well-known in statistics where it is called the commutation matrix, see Chapter 3 of~\cite{MN2019} for details. Its explicit form is given by
\begin{equation*}
\mathbf{K}_{(p,n)} =  \sum_{i=1}^p {  \sum_{j=1}^n { \mathbf{H}(i,j) \otimes \mathbf{H}(i,j)^{T} } } \ ,
\end{equation*}
where $\mathbf{H}(i,j)$ is an $p \times n$ matrix with a 1 in its $ij^{th}$ position and zeroes elsewhere. Important properties of the commutation matrix for our application are
\begin{equation}\label{eq:commat2}
\mathbf{K}_{(n,p)} = \mathbf{K}_{(p,n)}^{T} \text{  and  } \mathbf{K}_{(p,n)}^{T} \mathbf{K}_{(p,n)} = \mathbf{K}_{(p,n)} \mathbf{K}_{(p,n)}^{T} = \mathbf{I}_{p.n} \ .
\end{equation}
In other words, $ \mathbf{K}_{(p,n)}$ is an orthogonal matrix and its transpose and inverse is $\mathbf{K}_{(n,p)}$. Another useful property of the commutation matrix is that it can be used to reverse the order of a Kronecker product
\begin{equation}\label{eq:com_kron}
\mathbf{K}_{(s,p)} (  \mathbf{X}  \otimes \mathbf{Y} ) = (  \mathbf{Y}  \otimes \mathbf{X} ) \mathbf{K}_{(t,n)}  \text{ where } \mathbf{X}  \in \mathbb{R}^{p \times n}   \text{ and } \mathbf{Y}  \in \mathbb{R}^{s \times t} \ .
\end{equation}
This property will be also used frequently in the calculation of matrix derivatives. The following Lemma will also be useful later:
\\
\begin{theo2.2} \label{theo2.2:box}
Let $\mathbf{X} \in  \mathbb{R}^{p \times n}$, then
\begin{align}
& \mathbf{K}_{(p,n)} \emph{diag}( \emph{vec}( \mathbf{X} ) ) \mathbf{K}_{(n,p)}  = \emph{diag}( \emph{vec}( \mathbf{X}^{T} ) ) \ ,  \nonumber \\
& \emph{diag}( \emph{vec}( \mathbf{X} ) ) \mathbf{K}_{(n,p)}   = \mathbf{K}_{(n,p)} \emph{diag}( \emph{vec}( \mathbf{X}^{T} ) ) \ . \nonumber
\end{align}
\end{theo2.2}
\begin{proof}
Omitted.
\\
\end{proof}

\subsection{Topology of Euclidean vector or Frobenius matrix spaces} \label{topology:box}

For the sake of convenience, we define some notations first.
The following subsets of $\mathbb{R}^{p \times n}$, the set of $p \times n$ real matrices, and $\mathbb{R}$ will be used frequently in the following sections.
 \\
\begin{def2.1} \label{def2.1:box}
Let $p, n, k \in \mathbb{N}_*$ with $k \le \text{min}(p,n)$, then
\begin{IEEEeqnarray*}{rCl}
\mathbb{R}^{p \times n}_k         & = &  \big\{   \mathbf{Y} \in\mathbb{R}^{p \times n}  \text{ and } \emph{rank} ( \mathbf{Y} ) = {k}  \big\}  \ , \\
\mathbb{R}^{p \times n}_{\le k}  & = & \big\{   \mathbf{Y} \in\mathbb{R}^{p \times n}  \text{ and } \emph{rank}( \mathbf{Y} ) \le{k}  \big\} \ , \\
\mathbb{R}^{p \times n}_{> k}   & = & \big\{   \mathbf{Y} \in\mathbb{R}^{p \times n}  \text{ and } \emph{rank}( \mathbf{Y} ) >{k} \big\} \ , \\
\mathbb{R}^{p \times n}_{+}      & = & \big\{   \mathbf{Y} \in\mathbb{R}^{p \times n}  \text{ and }  \mathbf{Y}_{ij}  \ge {0}   \big\} \ , \\
\mathbb{R}^{p \times n}_{+*}     & = & \big\{   \mathbf{Y} \in\mathbb{R}^{p \times n}  \text{ and }  \mathbf{Y}_{ij}  > {0}   \big\} \ , \\
\mathbb{O}^{p \times k}            & = &  \big\{   \mathbf{U}\in\mathbb{R}^{p \times k}   \text{ }/\text{ }  \mathbf{U}^{T} \mathbf{U} = \mathbf{I}_k   \big\} \ , \\
\mathbb{O}^{k \times n}_{t}       & = &  \big\{   \mathbf{U}\in\mathbb{R}^{k \times n}   \text{ }/\text{ }  \mathbf{U} \mathbf{U}^{T} = \mathbf{I}_k   \big\} \ ,
\end{IEEEeqnarray*}
and
\begin{equation*}
\mathbb{R}_*  = \big\{   \mathbf{x}\in\mathbb{R}   \text{ }/\text{ }  \mathbf{x} \ne  {0}   \big\}  \ , 
\mathbb{R}_+  = \big\{   \mathbf{x}\in\mathbb{R}   \text{ }/\text{ }  \mathbf{x}  \ge {0}   \big\}  \text{  and }
\mathbb{R}_+*  = \big\{   \mathbf{x}\in\mathbb{R}   \text{ }/\text{ }  \mathbf{x}  > {0}   \big\} \ .
\end{equation*}

\end{def2.1}

The following definitions will also be useful:
\begin{def2.2} \label{def2.2:box}
Given $\mathbf{X}\in\mathbb{R}^{p \times n}$, $r\in\mathbb{R}_{+*}$, the open ball with center $\mathbf{X}$ and radius $r$ of $\mathbb{R}^{p \times n}$ is the set of $p \times n$ matrices defined by
\begin{equation*}
B_{p \times n}(\mathbf{X}, r)  = \big\{   \mathbf{Y} \in \mathbb{R}^{p \times n}  \text{ and }  \Vert  \mathbf{X} - \mathbf{Y} \Vert < r  \big\}
\end{equation*}
and the closed ball with center $\mathbf{X}$ and radius $r$ is the set
\begin{equation*}
\bar{B}_{p \times n}(\mathbf{X}, r)  = \big\{   \mathbf{Y} \in \mathbb{R}^{p \times n}  \text{ and }  \Vert  \mathbf{X} - \mathbf{Y} \Vert \le r  \big\} \  .
\end{equation*}
Note that in these definitions, $\Vert  \Vert$ can be the Frobenius norm or any norm defined on $\mathbb{R}^{p \times n}$ since  $\mathbb{R}^{p \times n}$ is a finite-dimensional vector space over $\mathbb{R}$ and, in this case, all norms on $\mathbb{R}^{p \times n}$ are equivalent and induce the same topology~\cite{B1993}\cite{C2017}. Similarly, given $\mathbf{x}\in\mathbb{R}^{n}$, $r\in\mathbb{R}_{+*}$, the open ball with center $\mathbf{x}$ and radius $r$ of $\mathbb{R}^{n}$ is the set of $n$-dimensional vectors defined by
\begin{equation*}
B_{n}(\mathbf{x}, r)  = \big\{   \mathbf{y} \in \mathbb{R}^{n}  \text{ and }  \Vert  \mathbf{x} - \mathbf{y} \Vert < r  \big\}
\end{equation*}
and the closed ball with center $\mathbf{x}$ and radius $r$ is the set
\begin{equation*}
\bar{B}_{n}(\mathbf{x}, r)  = \big\{   \mathbf{y} \in \mathbb{R}^{n}  \text{ and }  \Vert  \mathbf{x} - \mathbf{y} \Vert \le r  \big\} \ .
\end{equation*}
Again, here, $\Vert  \Vert$ can be the Euclidean norm or any norm defined on $\mathbb{R}^{n}$.
\\
\end{def2.2}
A set $U \subset \mathbb{R}^{p \times n}$ is open if every point of $U$ is contained in an open ball included in $U$. A set $U \subset \mathbb{R}^{p \times n}$ is closed if and only if its complement in $\mathbb{R}^{p \times n}$ is open. An arbitrary union of open sets is open and an arbitrary intersection of closed sets is closed. A finite union of closed sets is also closed. The closure of a set $U \subset \mathbb{R}^{p \times n}$ is the smallest closed set (in the sense of inclusion) of $\mathbb{R}^{p \times n}$ which contains $U$ and is denoted  $\bar{U}$. On the other hand, the interior of a set $U \subset \mathbb{R}^{p \times n}$ is the largest open set (in the sense of inclusion) of $\mathbb{R}^{p \times n}$ which is included in $U$ and is denoted $\mathring{U}$. A set $N$ is called a neighborhood of $\mathbf{X}$ in $\mathbb{R}^{p \times n}$  if there is an open set $U \subset N$ with $\mathbf{X} \in U$. A point $\mathbf{X}$ is a boundary point of a set $A$ if every neighborhood of $\mathbf{X}$ contains a point of $A$ and a point of its complement $B$ in $\mathbb{R}^{p \times n}$. The set of boundary points of $A$ is denoted $bd(A)$ and we have $bd(A) = \bar{A} \cap \bar{B}$. Thus, $bd(A)$ is closed as the intersection of two closed sets. The term frontier refers to the set of points of $bd(A)$ which are not in $A$ (e.g., $\bar{A}/A$). As an illustration, the boundary of both the open ball $B_{p \times n}(\mathbf{X}, r)$ and the closed ball $\bar{B}_{p \times n}(\mathbf{X}, r)$ is the sphere $S_{p \times n}(\mathbf{X}, r) = \big\{   \mathbf{Y} \in \mathbb{R}^{p \times n}  \text{ and }  \Vert  \mathbf{X} - \mathbf{Y} \Vert = r  \big\}$, but the frontier of $B_{p \times n}(\mathbf{X}, r)$ is  equal to $bd \big( B_{p \times n}(\mathbf{X}, r) \big)$ while the frontier of $\bar{B}_{p \times n}(\mathbf{X}, r)$ is empty. Let now $A$ and $B$ be two subsets of $\mathbb{R}^{p \times n}$ such that $A \subset B$. We say that $A$ is dense in $B$ if $B \subset \bar{A}$ and we say that $A$ is dense everywhere if  $\bar{A}=\mathbb{R}^{p \times n}$. Similar definitions hold for $\mathbf{x}$ in $\mathbb{R}^{n}$. In a finite-dimensional vector (or matrix) space over $\mathbb{R}$, a closed and bounded set is compact and all closed balls are compact. The preimage of a closed (open) set by a continuous function is a closed (open) set. The image of a compact set by a continuous function is compact.

Next, we collect some important topological results concerning certain subsets of $\mathbb{R}^{p \times n}$ in the following theorem that we will also use frequently in the following sections.
\begin{theo2.3} \label{theo2.3:box}
Let $p, n, k \in \mathbb{N}_*$ with $k \le \text{min}(p,n)$. The sets $\mathbb{O}^{p \times k}$ and $\mathbb{O}^{k \times n}_t$ are compact in $\mathbb{R}^{p \times k}$ and $\mathbb{R}^{k \times n}$, respectively. The set $\mathbb{R}^{p \times k}_k$  is open in $\mathbb{R}^{p \times k}$. If $k \ne n$, the interior of $\mathbb{R}^{p \times n}_k$ is empty and  $\mathbb{R}^{p \times n}_k$ is not closed or open in $\mathbb{R}^{p \times n}$. The sets $\mathbb{R}^{p \times n}_{> k}$ and $\mathbb{R}^{p \times n}_{\le k}$ are, respectively, open and closed in $\mathbb{R}^{p \times n}$. In all cases, the sets $\mathbb{R}^{p \times n}_{\le k}$ and $\mathbb{R}^{p \times n}_{<k}$  are, respectively, the closure and the frontier of $\mathbb{R}^{p \times n}_k$ in $\mathbb{R}^{p \times n}$, and $\mathbb{R}^{p \times n}_k$ is dense in  $\mathbb{R}^{p \times n}_{\le k}$. Furthermore, the set $\mathbb{R}^{p \times k}_k$ is an open subset dense everywhere in $\mathbb{R}^{p \times k}$.
\end{theo2.3}
\begin{proof}
Omitted. See  Section 3 and Theorem 2 of~\cite{HL2013},  Proposition 2.1 of~\cite{OA2022} and Section 3.1.5 of~\cite{AMS2008} for some details.
\end{proof}
 
 Importantly, since the $\emph{rank}(.)$ function defined on $\mathbb{R}^{p \times n}$ is an integer-valued and lower-semicontinuous function, an important result is that the rank function does not decrease in a sufficiently small neighborhood of any matrix $\mathbf{X}  \in \mathbb{R}^{p \times n}$~\cite{HL2013}. On the other hand, with the help of the SVD and Theorem~\eqref{theo2.1:box}, it is not difficult to see that if $\emph{rank}(\mathbf{X} ) = k < \text{min}(p,n)$ then any neighborhood of $\mathbf{X}$ contains matrices of rank $k+1, k+2, \dots, \text{min}(p,n)$.
 
To close this subsection, we finally recall the definition of a convex set for later reference. A subset $\mathcal{C}$ of a normed vector space $\mathcal{X}$ is called convex, if
\begin{equation} \label{eq:convex_set}
 \lambda.\mathbf{x} +  ( 1 - \lambda).\mathbf{y} \in \mathcal{C}  \text{ whenever } \mathbf{x} , \mathbf{y}  \in \mathcal{C}  \text{ and } 0 \le  \lambda \le 1 \  .
\end{equation}
Geometrically, a subset of a normed vector space is convex, if and only if, it contains the line segment $\lbrace  \lambda.\mathbf{x} +  ( 1 - \lambda).\mathbf{y}  \ /   \ 0 \le  \lambda \le 1 \rbrace$ joining each pair of its points $\mathbf{x} , \mathbf{y}$. As an illustration, the open and closed balls of $\mathcal{X}$ and the subspaces of $\mathcal{X}$  are all convex. Note, on the other hand, that the subsets $\mathbb{R}^{p \times n}_k$ and $\mathbb{R}^{p \times n}_{\le k}$ of $\mathbb{R}^{p \times n}$ are not convex if $k \ne n$ , which makes solving the WLRA problem~\eqref{eq:P0} challenging.

For more discussion about the topology of $\mathbb{R}^{p \times n}$ or arbitrary normed vector spaces, we refer the reader to~\cite{OR1970}\cite{B1993}\cite{C2017}.

\subsection{Differential calculus, variational geometry and optimization} \label{calculus:box}

We also assume that the reader has some familiarity with differentiation in a Euclidean space and derivatives of vectors and matrices, and their properties. Useful references on these topics are~\cite{OR1970}\cite{C2017}\cite{MN2019}.

Let $\mathcal{X}$ be a Euclidean space, e.g., a real vector space of finite dimension, say $k$, equipped with a scalar product $\langle . , .  \rangle_{\mathcal{X}}$ and the vector norm  $\Vert  . \Vert_{\mathcal{X}}$ induced by this scalar product. Let now $\phi(.)$ be a function from an open set $\Omega \subset \mathcal{X}$ to some other Euclidean space, say  $\mathcal{Y}$. In the following, we may have $\mathcal{Y} =  \mathbb{R}$, $\mathcal{Y} =  \mathbb{R}^{p}$ or  $\mathcal{Y} =  \mathbb{R}^{p \times n}$. We say that $\phi(.)$ is $\smallO( \Vert \mathbf{h} \Vert_{\mathcal{X}} )$ if
\begin{equation*}
\forall \varepsilon  \in  \mathbb{R}_{+*}, \exists \delta \in  \mathbb{R}_{+*} \text{ such that } \Vert \mathbf{h} \Vert_{\mathcal{X}} \le \delta \Longrightarrow  \Vert \phi( \mathbf{h} ) \Vert_{\mathcal{Y}}  \le \varepsilon  \Vert \mathbf{h} \Vert_{\mathcal{X}}  \  .
\end{equation*}
Similarly, we say that $\phi(.)$ is $\mathcal{O}( \Vert \mathbf{h} \Vert_{\mathcal{X}} )$ if
\begin{equation*}
\exists \lambda,  \eta \in  \mathbb{R}_{+*} \text{ such that } \Vert \mathbf{h} \Vert_{\mathcal{X}} \le \eta \Longrightarrow  \Vert \phi( \mathbf{h} ) \Vert_{\mathcal{Y}}  \le  \lambda  \Vert \mathbf{h} \Vert_{\mathcal{X}}  \  .
\end{equation*}
Notations like $\smallO( \Vert \mathbf{h} \Vert^{\alpha}_{\mathcal{X}} )$ or $\mathcal{O}( \Vert \mathbf{h} \Vert^{\alpha}_{\mathcal{X}} )$, for $\alpha \in \mathbb{N}_{*}$, will also be used to distinguish functions tending to zero faster than $\Vert \mathbf{h} \Vert^{\alpha}_{\mathcal{X}}$ instead of faster than $\Vert \mathbf{h} \Vert_{\mathcal{X}}$ .

With these notations, a function $\phi(.)$ from  the open set $\Omega \subset \mathcal{X}$ to the Euclidean space $\mathcal{Y}$ is said to be differentiable at  $\mathbf{a} \in \Omega$, if there exists a linear operator $\phi^{'} (\mathbf{a})$ from $\mathcal{X}$ to $\mathcal{Y}$ such that
\begin{equation*}
\phi(  \mathbf{a}  +  \mathbf{h} ) = \phi(  \mathbf{a}  ) + \phi^{'} (\mathbf{a})(  \mathbf{h} ) +  \smallO( \Vert \mathbf{h} \Vert_{\mathcal{X}} )  \  .
\end{equation*}
The set of (continuous) linear operators from $\mathcal{X}$ to $\mathcal{Y}$ is denoted by $\pounds (\mathcal{X}, \mathcal{Y}  )$. If $\mathcal{Y} =  \mathbb{R}$, then $\phi(.)$ is a real-valued function, $\phi^{'} (\mathbf{a})$ is a linear form and it can be represented by an unique element of $\mathcal{X}$, called the gradient of  $\phi(.)$ at $\mathbf{a}$ and denoted by $\nabla \phi( \mathbf{a}  )$, which verifies
\begin{equation*}
\phi^{'} (\mathbf{a})(  \mathbf{h} ) = \langle \nabla \phi( \mathbf{a}  ) , \mathbf{h}  \rangle_{\mathcal{X}} \text{ , } \forall  \mathbf{h} \in \mathcal{X}  \  .
\end{equation*}
In the same conditions, e.g., when $\mathcal{Y} =  \mathbb{R}$, we can consider the function from  $\Omega$ into $\pounds (\mathcal{X}, \mathbb{R}  )$, which at $\mathbf{a} \in \Omega$ associates the linear form $\phi^{'} (\mathbf{a})$. If this new function is itself differentiable, we get the second-order differential of $\phi(.)$ at $\mathbf{a}$, which is denoted by $\phi^{''} (\mathbf{a})$ and is an element of $\pounds (\mathcal{X},  \pounds (\mathcal{X}, \mathbb{R}  ) ) \simeq \pounds ( \mathcal{X} ,  \mathcal{X}  ; \mathbb{R} ) $. In other words, $\phi^{''} (\mathbf{a})$ can be identified with an unique bilinear form, also noted $\phi^{''} (\mathbf{a}) \in \pounds ( \mathcal{X} ,  \mathcal{X}  ; \mathbb{R} ) $ by an abuse of notation, and defined by
\begin{equation*}
\lbrack \phi^{''} (\mathbf{a}) (  \mathbf{h} ) \rbrack (  \mathbf{k} ) = \phi^{''} (\mathbf{a}) (  \mathbf{h} , \mathbf{k} ) \text{ , } \forall  (\mathbf{h} , \mathbf{k} ) \in \mathcal{X} \times \mathcal{X}  \  .
\end{equation*}
This bilinear form is also symmetric and yields the following second-order approximation of $\phi(.)$ at $\mathbf{a}$
\begin{equation*}
\phi(  \mathbf{a}  +  \mathbf{h} ) = \phi(  \mathbf{a}  ) + \phi^{'} (\mathbf{a})(  \mathbf{h} )  + \phi^{''} (\mathbf{a})(  \mathbf{h}, \mathbf{h} ) +  \smallO( \Vert \mathbf{h} \Vert^{2}_{\mathcal{X}} )  \  .
\end{equation*}
Again, using the Euclidean structure associated with $\mathcal{X}$, the symmetric bilinear form $\phi^{''} (\mathbf{a} )$ can be associated with an unique symmetric linear operator from $\mathcal{X}$  to $\mathcal{X}$, called the Hessian of $\phi(.)$ at $\mathbf{a}$, denoted by $\nabla^{2} \phi( \mathbf{a}  )$, and defined by
\begin{equation*}
\phi^{''} (\mathbf{a}) (  \mathbf{h} , \mathbf{k} ) =  \langle \nabla^{2} \phi( \mathbf{a}  ) ( \mathbf{h} ), \mathbf{k}  \rangle_{\mathcal{X}} =  \langle \mathbf{h} ,  \nabla^{2} \phi( \mathbf{a}  ) ( \mathbf{k} ) \rangle_{\mathcal{X}} \text{ , } \forall  (\mathbf{h} , \mathbf{k} ) \in \mathcal{X} \times \mathcal{X}  \  .
\end{equation*}
Note that  both $\nabla \phi( \mathbf{a}  )$ and $\nabla^{2} \phi( \mathbf{a}  )$ depend on the scalar product  $\langle . ,  .  \rangle_{\mathcal{X}}$, while $\phi^{'} (\mathbf{a})$ and $\phi^{''} (\mathbf{a})$ do not. When $\mathcal{X} =  \mathbb{R}^{k}$  and is equipped with the standard Euclidean inner product defined in Subsection~\ref{lin_alg:box} and the canonical basis of $\mathbb{R}^{k}$ is used to represent vectors in $\mathbb{R}^{k}$, the self-adjoint linear operator $\nabla^{2} \phi( \mathbf{a}  )$ is represented by a $k \times k$ symmetric real matrix, which is known as the Schwarz's theorem~\cite{C2017}. Then, by a slight abuse of notation, we will also use the symbol  $\nabla^{2} \phi( \mathbf{a}  )$ to represent this $k \times k$ symmetric matrix and we can write
\begin{equation*}
\phi^{''} (\mathbf{a}) (  \mathbf{h} , \mathbf{k} ) =  \langle \nabla^{2} \phi( \mathbf{a}  ) ( \mathbf{h} ), \mathbf{k}  \rangle_{2} =  \mathbf{h}^{T} \nabla^{2} \phi( \mathbf{a}  ) \mathbf{k}  \text{ , } \forall  (\mathbf{h} , \mathbf{k} ) \in \mathbb{R}^{k} \times \mathbb{R}^{k}  \  .
\end{equation*}

In the following sections, instead of the generic notations $\phi^{'} (\mathbf{a})$ and $\phi^{''} (\mathbf{a})$ for the first- and second-order derivatives of a (twice) differentiable function $\phi(.)$ from  $\Omega \subset \mathcal{X}$ to $\mathcal{Y}$ at a point $\mathbf{a} \in \Omega$, the symbols $\mathit{D}$, $\mathit{J}$, $\nabla$, $\nabla^2$ will be used for the (Euclidean) derivative of a real matrix, the Jacobian matrix (e.g., derivative) of a real vector function, the gradient (e.g., first derivative) and Hessian (e.g., second derivative) of a real functional, respectively.

As a first illustration, a $q \times r$ matrix function $\mathbf{C}(\mathbf{a})$ for $\mathbf{a} \in  \Omega  = \mathbb{R}^{k} $ can be interpreted as a (nonlinear) mapping from the linear space of parameters, $\mathbb{R}^{k}$, into the space of linear transformations $\pounds ( \mathbb{R}^{r}, \mathbb{R}^{q} ) = \mathcal{Y}$, which can be identified to the linear space $\mathbb{R}^{q \times r}$~\cite{C2017}. Consequently, the derivative of  the matrix function $\mathbf{C}(.)$ at a point $\mathbf{a}  \in  \mathbb{R}^{k}$ is an element of $\pounds ( \mathbb{R}^{k}, \pounds ( \mathbb{R}^{r}, \mathbb{R}^{q} ) )$, or equivalently $\pounds ( \mathbb{R}^{k}, \mathbb{R}^{q \times r} ) \simeq \mathbb{R}^{q \times r \times k}$, and can be interpreted as the tridimensional tensor $\mathit{D} (\mathbf{C}(\mathbf{a} ) ) \in \mathbb{R}^{q \times r \times k}$ defined by
\begin{equation} \label{eq:D_mat_tensor}
\big \lbrack \mathit{D}( \mathbf{C}(\mathbf{a}  ) ) \big \rbrack_{ijl} =  \frac{ \partial   \mathbf{C}_{ij}(\mathbf{a} )  }  {\partial\mathbf{a}_l}  \text{ for } i = 1, \cdots , q \text{ ; } j = 1, \cdots , r  \text{ ; } l = 1, \cdots , k \ ,
\end{equation}
following~\cite{GP1973}.

On the other hand, the first derivative of a real $q$-vector function  $\mathbf{r}(.)$ at a point $\mathbf{a}  \in  \mathbb{R}^{k}$ is an element of $\pounds ( \mathbb{R}^{k}, \mathbb{R}^{q} ) \simeq  \mathbb{R}^{q \times k}$. If $\mathbb{R}^{k}$ and $\mathbb{R}^{q}$ are equipped  with their usual Euclidean inner products and the canonical bases of $\mathbb{R}^{k}$ (e.g., the columns of the identity matrix $\mathbf{I}_{k}$) and $\mathbb{R}^{q}$ (e.g., the columns of the identity matrix $\mathbf{I}_{q}$) are used to represent vectors in these two linear spaces, the first derivative of $\mathbf{r}(.)$ at a point $\mathbf{a}  \in  \mathbb{R}^{k}$ can be  identified to the Jacobian matrix $\mathit{J} (\mathbf{r}(\mathbf{a} ) ) \in \mathbb{R}^{q \times k}$ defined by
\begin{equation}
\big \lbrack \mathit{J}( \mathbf{r}( \mathbf{a}  ) )\big \rbrack_{ij}  =  \frac{ \partial   \mathbf{r}_{i}(\mathbf{a} )  }  {\partial\mathbf{a}_j}  \text{ for } i = 1, \cdots , q \text{ ; } j = 1, \cdots , k  \  ,
\end{equation}
where $\mathbf{r}_{i}(\mathbf{a} )$ is the $i^{th}$ component of the real $q$-vector $\mathbf{r}( \mathbf{a}  )$. Note that each $i^{th}$ row of the Jacobian matrix $\mathit{J} (\mathbf{r}(\mathbf{a} ) ) \in \mathbb{R}^{q \times k}$ is equal to the transpose of the gradient of the real function $\mathbf{r}_{i}(.)$ at the point $\mathbf{a} \in  \mathbb{R}^{k}$,
\begin{equation*}
\big \lbrack \mathit{J}( \mathbf{r}( \mathbf{a}  ) )\big \rbrack_{i.}  =  \nabla \mathbf{r}_{i}( \mathbf{a}  )^{T}  \  .
\end{equation*}
In addition, if  $\mathbf{r}(.)$ is continuously differentiable at a point $\mathbf{a}  \in  \mathbb{R}^{k}$, we have the following first-order Taylor expansion
\begin{equation} \label{eq:taylor_expan1}
\mathbf{r}(  \mathbf{a}  +  d\mathbf{a} ) = \mathbf{r}(  \mathbf{a}  ) + \mathit{J}( \mathbf{r}(\mathbf{a} ) )d\mathbf{a} +  \mathcal{O}( \Vert d\mathbf{a} \Vert^{2}_{2} )  \  .
\end{equation}

As another illustration, if $\mathbb{R}^{k}$ is again equipped  with its usual Euclidean inner product and its canonical basis, the gradient of a real functional $\phi(.)$ at a point $\mathbf{a}  \in  \mathbb{R}^{k}$ forms a $k \times 1$ column vector, i.e.,
\begin{equation}  \label{eq:D_grad_vec}
\big \lbrack  \nabla \phi( \mathbf{a}  ) \big \rbrack_{i}  =   \frac{ \partial   \phi(\mathbf{a} )  }  {\partial\mathbf{a}_i} \text{ for } i = 1, \cdots , k \ ,
\end{equation}
and the Hessian of a real functional $\phi(.)$ at a point $\mathbf{a}  \in  \mathbb{R}^{k}$ can be identified with a $k \times k$ symmetric matrix, $\nabla^2 \phi( \mathbf{a} )$, defined by
\begin{equation} \label{eq:D_hess_mat}
\big \lbrack  \nabla^2 \phi( \mathbf{a}  ) \big \rbrack_{ij}  =  \frac{ \partial^2   \phi(\mathbf{a} )  }  {\partial\mathbf{a}_i \partial\mathbf{a}_j }  =  \frac{ \partial^2   \phi(\mathbf{a} )  }  {\partial\mathbf{a}_j \partial\mathbf{a}_i }  = \big \lbrack  \nabla^2 \phi( \mathbf{a}  ) \big \rbrack_{ji} \text{ for } i = 1, \cdots , k \text{ ; } j = 1, \cdots , k \ .
\end{equation}
Finally, if $\phi(.)$ is at least twice continuously differentiable at a point $\mathbf{a}  \in  \mathbb{R}^{k}$, we have the following second-order Taylor expansion
\begin{align} \label{eq:taylor_expan2}
\phi( \mathbf{a}  + d\mathbf{a} ) & = \phi( \mathbf{a}  ) + \langle d\mathbf{a} , \nabla \phi( \mathbf{a} ) \rangle_{2} + \frac{1}{2}  \langle \nabla^2 \phi( \mathbf{a}  )d\mathbf{a} ,  d\mathbf{a} \rangle_{2}  +  \mathcal{O}( \Vert d\mathbf{a} \Vert^{3}_{2} )   \nonumber \\
                                                   & = \phi( \mathbf{a}  ) + d\mathbf{a}^{T} \nabla \phi( \mathbf{a}  ) + \frac{1}{2}  d\mathbf{a}^{T} \nabla^2 \phi( \mathbf{a}  ) d\mathbf{a} +  \mathcal{O}( \Vert d\mathbf{a} \Vert^{3}_{2} ) \ .
\end{align}

Let $\mathcal{K}$ be a nonempty subset of $\mathbb{R}^{k}$ or more generally  a nonempty subset of an arbitrary normed vector space. We recall that a point $\hat{\mathbf{a}}  \in  \mathcal{K}$  is a global minimizer of a real function $\phi(.)$ defined over $\mathcal{K}$, if and only if, $\forall \mathbf{a} \in  \mathcal{K}$, we have $\phi( \mathbf{a} ) \ge \phi( \hat{\mathbf{a}} )$. On the other hand, $\hat{\mathbf{a}}$ is a local minimizer of $\phi(.)$ over $\mathcal{K}$, if and only if, $\exists r   \in  \mathbb{R}_{+*}$ such that $\forall \mathbf{a} \in  \mathcal{K}$ and $\Vert \mathbf{a} - \hat{\mathbf{a}} \Vert_{2} < r$ imply $\phi( \mathbf{a} ) \ge \phi( \hat{\mathbf{a}} )$. Similar definitions hold for strict global and local minimizers of  $\phi(.)$ over $\mathcal{K}$.

Let now  $\Omega$ be an open subset of $\mathbb{R}^{k}$ or more generally  an open subset of a  normed vector space of finite dimension. A necessary condition for a point $\hat{\mathbf{a}}   \in   \Omega$ to minimize a real function $\phi(.)$ defined and assumed to be twice continuously differentiable  on $\Omega$ is that the gradient of $\phi(.)$ at  $\hat{\mathbf{a}}$ is equal to the zero-vector of the ambient linear space, i.e.,
\begin{equation} \label{eq:D_first_kkt}
 \nabla \phi( \hat{\mathbf{a}}   ) = \mathbf{0}^{k}  \ ,
\end{equation}
and this condition defines the first-order Karush-Kuhn-Tucker (KKT) condition~\cite{OR1970}\cite{C2017}. If such KKT condition is satisfied then $\hat{\mathbf{a}}$ is said to be a first-order stationary or critical point of $\phi(.)$. However, first-order critical points of $\phi(.)$ can be minimizers, but also maximizers or saddle points (e.g., points for which the Hessian matrix has both positive and negative eigenvalues). A necessary condition for a first-order stationary point $\hat{\mathbf{a}}$ to be a local minimizer of $\phi(.)$ is that the Hessian  (bilinear form or matrix)  $\nabla^2 \phi( \hat{\mathbf{a}} )$ is positive semi-definite~\cite{OR1970}\cite{C2017}:
\begin{equation} \label{eq:D_second_kkt}
 \nabla^2 \phi( \hat{\mathbf{a}}  ) \big( d\mathbf{a} ,  d\mathbf{a}  \big)=   \langle \nabla^2 \phi(  \hat{\mathbf{a}}  )d\mathbf{a} ,  d\mathbf{a} \rangle_{2} \ge 0 \text{ , } \forall d\mathbf{a}  \in \mathbb{R}^{k}  \ . 
\end{equation}
Such first-order critical points for which the Hessian is positive semi-definite are called second-order stationary or critical points of $\phi(.)$. On the other hand, a sufficient condition for a first-order stationary point $\hat{\mathbf{a}}$ to be a strict local minimizer of $\phi(.)$ is that $\nabla^2 \phi( \hat{\mathbf{a}} )$ is positive definite  (second-order KKT condition). These assertions can be derived by noting that the second-order Taylor expansion of $\phi(.)$ at a first-order stationary point $\mathbf{\widehat{a}}$ reduces to
\begin{equation*}
\phi( \mathbf{\widehat{a}} + d\mathbf{a} ) = \phi( \mathbf{\widehat{a}} ) +  \frac{1}{2} \langle \nabla^2 \phi(  \hat{\mathbf{a}}  )d\mathbf{a} ,  d\mathbf{a} \rangle_{2}+  \mathcal{O}( \Vert d\mathbf{a} \Vert^{3}_{2} )  \ ,
\end{equation*}
see~\cite{OR1970}\cite{C2017} for details.

We now consider the case, where we seek to minimize a real function $\phi(.)$ on a linear subspace $\Upsilon \subset \Omega$ of dimension $s$, where $\Omega$ is an open subset of $\mathbb{R}^{k}$ on which $\phi(.)$ is defined and twice-differentiable. Obviously, we must have $\emph{dim}(  \Upsilon ) \le k$. In these conditions, we can consider $\phi(.)$ as a function from $\Omega$ to $\mathbb{R}$, but we can also consider its restriction to $\Upsilon$, $\phi_{\Upsilon}(.)$. If $\phi(.)$ is differentiable on $\Omega$ than $\phi_{\Upsilon}(.)$ will be differentiable on $\Upsilon$ as well and their differentials verify
\begin{equation*}
\phi^{'} (\mathbf{a}) (\mathbf{b}) = \phi^{'}_{\Upsilon} (\mathbf{a}) (\mathbf{b}) \text{ , } \forall \mathbf{a}, \mathbf{b} \in \Upsilon  \ .
\end{equation*}
In other words, the linear form $\phi^{'}_{\Upsilon} (\mathbf{a})$ is nothing else than the restriction of the linear form $\phi^{'} (\mathbf{a})$ to $\Upsilon$. Furthermore, if we equip both $\mathbb{R}^{k}$ and its linear subspace $\Upsilon$ with the same Euclidean structure induced by $\mathbb{R}^{k}$, we have, by definition, $\langle \mathbf{a} , \mathbf{b}  \rangle_{\Upsilon} = \langle \mathbf{a} ,\mathbf{b}  \rangle_{2}, \forall  \mathbf{a}, \mathbf{b} \in \Upsilon$.
In these conditions, the linear forms $\phi^{'} (\mathbf{a})$ and $\phi^{'}_{\Upsilon} (\mathbf{a})$ can be both represented by their own gradients, $\nabla \phi (\mathbf{a})$ and $\nabla \phi_{\Upsilon} (\mathbf{a})$, which are, respectively, elements of the linear spaces $\mathbb{R}^{k}$ and $\Upsilon$ such that
\begin{equation*}
 \phi^{'} (\mathbf{a}) (\mathbf{b}) = \langle  \nabla \phi (\mathbf{a}) , \mathbf{b}  \rangle_{2}  \text{ and } \phi^{'}_{\Upsilon} (\mathbf{a} ) (\mathbf{b}) = \langle \nabla \phi_{\Upsilon} (\mathbf{a} ),\mathbf{b}  \rangle_{\Upsilon} \text{ , }  \forall  \mathbf{a}, \mathbf{b} \in \Upsilon  \ .
\end{equation*}
Since the linear forms $\phi^{'} (\mathbf{a})$ and $\phi^{'}_{\Upsilon} (\mathbf{a} )$ coincide on $\Upsilon$, we deduce immediately that
\begin{equation*}
 \langle  \nabla \phi (\mathbf{a}) , \mathbf{b}  \rangle_{2}  = \langle \nabla \phi_{\Upsilon} (\mathbf{a} ),\mathbf{b}  \rangle_{\Upsilon} ,  \forall  \mathbf{a}, \mathbf{b} \in \Upsilon  \ .
\end{equation*}
Next, remember that $\nabla \phi (\mathbf{a}) \in \mathbb{R}^{k}$, while $\nabla \phi_{\Upsilon} (\mathbf{a} ) \in \Upsilon$, but we can easily expressed $\nabla \phi_{\Upsilon} (\mathbf{a} )$ as a function of $\nabla \phi (\mathbf{a})$. More precisely, using the two complementary orthogonal projectors on $\Upsilon$ and $\Upsilon^{\bot}$ (considered here as $k \times k$ symmetric and idempotent matrices rather than linear operators as discussed in Subsection~\ref{lin_alg:box}), denoted, respectively, by $\mathbf{P}_{\Upsilon}$ and $\mathbf{P}^{\bot}_{\Upsilon}$, and defined on $\mathbb{R}^{k}$, we have
\begin{equation*}
\nabla \phi (\mathbf{a}) = \mathbf{P}_{\Upsilon} \nabla \phi (\mathbf{a})  + \mathbf{P}^{\bot}_{\Upsilon} \nabla \phi (\mathbf{a}) , \text{ with } \mathbf{P}_{\Upsilon} \nabla \phi (\mathbf{a}) \in \Upsilon \text{ and } \mathbf{P}^{\bot}_{\Upsilon}  \nabla \phi (\mathbf{a})  \in \Upsilon^{\bot}  \  ,
\end{equation*}
and this implies immediately that
\begin{equation*}
 \langle  \nabla \phi (\mathbf{a}) , \mathbf{b}  \rangle_{2}  =  \langle  \mathbf{P}_{\Upsilon}  \nabla \phi (\mathbf{a}) , \mathbf{b}  \rangle_{2}  = \langle  \mathbf{P}_{\Upsilon}  \nabla \phi (\mathbf{a}) , \mathbf{b}  \rangle_{\Upsilon} = \langle \nabla \phi_{\Upsilon} (\mathbf{a} ),\mathbf{b}  \rangle_{\Upsilon} , \forall  \mathbf{a}, \mathbf{b} \in \Upsilon \  .
\end{equation*}
Then, by the unicity of the gradient of $\phi_{\Upsilon}(.)$ at $\mathbf{a} \in \Upsilon$, we get, by identification, the vector equality
\begin{equation} \label{eq:D_grad_vec2}
\nabla \phi_{\Upsilon} (\mathbf{a}) =  \mathbf{P}_{\Upsilon}  \nabla \phi (\mathbf{a}) , \forall  \mathbf{a} \in \Upsilon \  .
\end{equation}
In words, the gradient of $\phi_{\Upsilon} (.)$ at $\mathbf{a} \in \Upsilon$ is simply the orthogonal projection on $\Upsilon$ of the gradient of $\phi (.)$ at $\mathbf{a}$, considered as an element of $\mathbb{R}^{k}$ instead of $\Upsilon$. The interpretation of this result is simple and is that it is not necessary to check all the feasible directions in $\mathbb{R}^{k}$ to satisfy the first-order stationary condition for a point $\mathbf{a} \in \Upsilon$ if the search space is reduced to $\Upsilon$, only those belonging to $\Upsilon$ matter in that case. Similarly, it is not too difficult to verify that the linear operators $\nabla^{2} \phi_{\Upsilon} (\mathbf{a})$ and $\nabla^{2} \phi (\mathbf{a})$ are related by
\begin{equation}\label{eq:D_hess_mat2}
\nabla^{2} \phi_{\Upsilon} (\mathbf{a} )  \lbrack \mathbf{b} \rbrack = \mathbf{P}_{\Upsilon} \big ( \nabla^{2} \phi (\mathbf{a})   \lbrack \mathbf{b} \rbrack  \big) , \forall  \mathbf{b} \in \Upsilon \  ,
\end{equation}
where $\mathbf{P}_{\Upsilon}$ is now interpreted as an orthogonal projector operator rather than as a matrix. However, keep in mind that, in the above formulae, $\nabla^{2} \phi_{\Upsilon} (\mathbf{a} )$ is expressed as a non-symmetric linear operator from $\mathbb{R}^{k}$ to $\mathbb{R}^{k}$ rather than as a symmetric linear operator from $\Upsilon$ to $\Upsilon$ (or a $s \times s$ symmetric matrix), but both operators coincide on $\Upsilon$. Finally, $\Upsilon$ being a linear space of dimension $s$, the definitions of the first- and second-order stationary points of $\phi_{\Upsilon} (.)$ are exactly similar to those stated above. Namely, the first-order KKT condition is met if $\mathbf{P}_{\Upsilon}  \nabla \phi (\mathbf{a}) = \mathbf{0}^{k}$ and the second-order stationary (KKT) condition is equivalent to say that the bilinear form associated with the self-adjoint linear operator $\mathbf{P}_{\Upsilon} o \nabla^{2} \phi (\mathbf{a})$ is positive semi-definite (positive definite) over $\Upsilon$, for $ \mathbf{a} \in \Upsilon$.

Next, we consider the problem of the minimization of a smooth function $\phi(.)$ defined on a smooth submanifold $\mathcal{M}$ of dimension $r$ embedded in $\mathbb{R}^{k}$. This problem enters in the domain of differential geometry and optimization on Riemannian  manifolds described comprehensively in~\cite{AMS2008}\cite{RS2022}\cite{B2023}. We first precise what we mean by a smooth function and a smooth embedded manifold in $\mathbb{R}^{k}$ in the following definitions, which will be sufficient for our purpose.

\begin{def2.3} \label{def2.3:box}
Given any set $Q \subset \mathbb{R}^{k}$ and a mapping $\phi(.)$ from $Q$ to $\mathbb{R}^{m}$, we say that $\phi(.)$ is $C^{p}$ smooth if,  $\forall  \mathbf{a}  \in Q$, there is a neighborhood $U$ of $\mathbf{a}$ in $\mathbb{R}^{k}$ and a $C^{p}$ differentiable mapping $\hat{\phi}(.)$ from $U$ to $\mathbb{R}^{m}$ that agrees with $\phi(.)$ on $U \cap Q$. Here we assume that $p$ lies in $\mathbb{N}_{*} \cup  \lbrace \infty \rbrace$.
\end{def2.3}

\begin{def2.4} \label{def2.4:box}
Let $\mathcal{M}$ be a nonempty subset of $\mathbb{R}^{k}$. We say that $\mathcal{M}$ is a $C^{p}$ embedded submanifold of dimension $r$ of $\mathbb{R}^{k}$, with $r \le k$, if for each point $\mathbf{a} \in \mathcal{M}$, there is an open neighborhood $U$ around $\mathbf{a}$ in $\mathbb{R}^{k}$ such that $U  \cap \mathcal{M} = f^{-1} ( \mathbf{0}^{k-r} )$ for some $C^{p}$ differentiable map $f(.)$ from $U$ to $\mathbb{R}^{k-r}$, with its differential at $\mathbf{a}$, $f^{'} (\mathbf{a})$, being a surjective linear operator, which is equivalent to say that   $f^{'} (\mathbf{a})$ has full rank equal to $ k - r $.
\end{def2.4}
The mapping $f(.)$ is called a local defining function for $\mathcal{M}$ at  $\mathbf{a}$. This definition implies that, locally around  $\mathbf{a}$, a smooth embedded submanifold (of dimension $r$)  looks like a subspace of dimension $r$ of $\mathbb{R}^{k}$, which is also the dimension of the kernel of  $f^{'} (\mathbf{a})$, see Theorem 3.12 of Boumal~\cite{B2023} or Theorem 2.1.10 of Robbin and Salomon~\cite{RS2022} for details. More precisely, if $\mathcal{M}$ is a smooth submanifold of $\mathbb{R}^{k}$, it admits a tangent space noted $\mathcal{T_{\mathbf{a}} M}$, which is nothing else than the kernel of $f^{'} (\mathbf{a})$, where  $f(.)$ is any local defining function for $\mathcal{M}$ at $\mathbf{a}$ (see Theorem 3.15 of Boumal~\cite{B2023} or Theorem 2.2.3 of Robbin and Salomon~\cite{RS2022}) and this tangent space can be interpreted as a vector subspace of $\mathbb{R}^{k}$ that approximates the smooth submanifold locally. Thus, a smooth (sub)manifold of dimension $r$ is defined as a set that locally looks like a $r$-dimensional space, but can be very different globally.

Next, we clarify what we call a tangent vector to an arbitrary subset $C$ of a general vector space $\mathcal{X}$  at a point  $\mathbf{a} \in C$ in the following definition, which will also be sufficient for our purpose:
\begin{def2.5} \label{def2.5:box}
Let  $\mathcal{X}$ be a normed vector space and $C$ a nonempty subset of $\mathcal{X}$. A vector $\mathbf{d}  \in \mathcal{X}$ is a tangent vector to $C$ at $\mathbf{a} \in C$, if and only if, it exists $\alpha > 0$ and a mapping $\varepsilon(.)$ from $\lbrack -\alpha,  \alpha \rbrack$ to $\mathcal{X}$ such that
\begin{equation*}
\mathbf{a} + t.\mathbf{d} + t.\varepsilon(t)  \in C , \forall t \in \lbrack -\alpha,  \alpha \rbrack , \text{ and }  \lim_{t \to 0} \varepsilon(t) = 0 \ ,
\end{equation*}
or, equivalently, if it exists an open interval $I$ of $\mathbb{R}$ containing $t = 0$ and a function $\mathcal{E} : I \to C$ such that  $\mathcal{E}(.)$ is derivable at $t = 0$ with $\mathbf{d} = \mathcal{E}^{'} (0)$ and  $\mathcal{E}(0) = \mathbf{a}$.

The set of all tangent vectors to $C$  at $\mathbf{a} \in C$ is noted  $\mathcal{T_{\mathbf{a}}} C$. If $\mathcal{T_{\mathbf{a}}} C$ is a linear subspace of $\mathcal{X}$, it is called the tangent space to $C$ at  $\mathbf{a}$.
\\
\end{def2.5}
Note that, in this definition, it is only required that the function $\mathcal{E}(.)$ is derivable at $t = 0$, not on all $I$ and this definition is sufficient for many results stated in~\cite{B2023} or~\cite{RS2022} for a $C^{p}$ or $C^{\infty}$ differentiable function $\mathcal{E}(.)$ on all $I$. Furthermore, keep in mind that if $\mathcal{M}$ is a $C^{p}$ embedded submanifold in the sense of Definition~\ref{def2.4:box}, all the elements of its tangent space at a given point $\mathbf{a} \in \mathcal{M}$, defined as the kernel of $f^{'} (\mathbf{a})$ for any given local defining function $f(.)$  for $\mathcal{M}$ at  $\mathbf{a}$, verify Definition~\ref{def2.5:box} and the terminology is thus consistent~\cite{B2023}\cite{RS2022}.

Let $\mathcal{M}$ be a $C^{p}$ embedded submanifold of dimension $r$ of $\mathbb{R}^{k}$ in the sense of Definition~\eqref{def2.4:box}. If we now endow $\mathbb{R}^{k}$ with its standard Euclidean inner product, $\mathcal{T_{\mathbf{a}} M}$, which is a linear subspace of $\mathbb{R}^{k}$, admits an orthogonal supplementary subspace in $\mathbb{R}^{k}$, which is called the normal space of $\mathcal{M}$ at $\mathbf{a}$ and is denoted by $\mathcal{N_{\mathbf{a}} M}$ in the following. Both $\mathcal{T_{\mathbf{a}} M}$ and $\mathcal{N_{\mathbf{a}} M}$ are linear subspaces of $\mathbb{R}^{k}$ and we have the identity: $\mathbb{R}^{k} = \mathcal{T_{\mathbf{a}} M} \oplus \mathcal{N_{\mathbf{a}} M}$, which is equivalent to say that any vector of $\mathbb{R}^{k}$ can be written uniquely as the sum of an element of $\mathcal{T_{\mathbf{a}} M}$ and an element of $\mathcal{N_{\mathbf{a}} M}$.

Suppose now that we want to minimize a $C^{p}$ smooth function $\phi(.)$ from a smooth submanifold $\mathcal{M} \subset \mathbb{R}^{k}$ to $\mathbb{R}$. To define and also analyze Riemannian optimization methods on  $\mathcal{M}$ for solving this kind of problems, we need to define the notions of the Riemannian gradient and Hessian, which will be obviously different from their Euclidean analogs as  $\mathcal{M}$ is only locally homeomorphic to an Euclidean vector space. First, similarly to the standard case of a differentiable function from an open set $U$ to $\mathbb{R}$, the smooth function $\phi(.)$ admits a differential at $\mathbf{a} \in \mathcal{M}$, which is a linear mapping from $\mathcal{T_{\mathbf{a}} M}$ to $\mathbb{R}$ denoted also by $\phi^{'} (\mathbf{a})$~\cite{AMS2008}\cite{RS2022}\cite{B2023}. If, $\forall \mathbf{a} \in \mathcal{M}$, we equip $\mathcal{T_{\mathbf{a}} M}$ with the standard Euclidean inner product induced by $\mathbb{R}^{k}$, e.g.,
\begin{equation*}
\langle . ,  . \rangle_{\mathcal{T_{\mathbf{a}} M}} = \langle . ,  . \rangle_{2} , \forall \mathbf{a} \in \mathcal{M} \ ,
\end{equation*}
 $\mathcal{M}$ is then, by definition,  equipped with a smoothly varied inner product on all its tangent spaces and $(\mathcal{M}, \langle . , . \rangle_{2} )$ is a Riemannian manifold~\cite{AMS2008}\cite{RS2022}\cite{B2023}. In this setting, the Riemannian gradient of $\phi(.)$ at  $\mathbf{a} \in \mathcal{M}$, denoted here by $\nabla_{R} \phi (\mathbf{a})$, is then defined as the unique vector in $\mathcal{T_{\mathbf{a}} M}$ satisfying
\begin{equation*}
\phi^{'} (\mathbf{a})(  \mathbf{b} ) = \langle \nabla_{R} \phi( \mathbf{a}  ) , \mathbf{b}  \rangle_{\mathcal{T_{\mathbf{a}} M}} = \langle \nabla_{R} \phi( \mathbf{a}  ) , \mathbf{b}  \rangle_{2} \text{ , } \forall  \mathbf{b} \in \mathcal{T_{\mathbf{a}} M} \ ,
\end{equation*}
where $\phi^{'} (\mathbf{a})$ is the differential of the smooth mapping $\phi(.)$ at $\mathbf{a}$ in the sense defined above. We can also define the Riemannian Hessian of the smooth mapping $\phi(.)$ at $\mathbf{a}$, denoted by $\nabla^{2}_{R} \phi( \mathbf{a}  )$, which is  a self-adjoint linear operator from $\mathcal{T_{\mathbf{a}} M}$ to $\mathcal{T_{\mathbf{a}} M}$ defined by
\begin{equation*}
\nabla^{2}_{R} \phi( \mathbf{a}  ) \lbrack  \mathbf{b} \rbrack = \tilde{\nabla}_{\mathbf{b}} \nabla_{R} \phi( \mathbf{a} ) , \forall  \mathbf{b} \in \mathcal{T_{\mathbf{a}} M} \ ,
\end{equation*}
where $\tilde{\nabla}_{(.)} (.)$ denotes the so-called Levi-Civita connection on $\mathcal{M}$. The  Levi-Civita connection $\tilde{\nabla}_{\eta_{\mathbf{a}}}  \xi_{\mathbf{a}}$ on the Riemannian manifold $\mathcal{M}$ acting on two vector fields, $\eta_{\mathbf{a}}$ and $\xi_{\mathbf{a}}$, in the tangent bundle of $\mathcal{M}$ (the tangent bundle is the disjoint union of all the tangent spaces of the manifold $\mathcal{M}$, see Definition 3.42 in Boumal~\cite{B2023}) is a generalization of the notion of directional derivative of a vector field on the manifold $\mathcal{M}$. In this way, the Levi-Civita connection $\tilde{\nabla}_{\eta_{\mathbf{a}}}  \xi_{\mathbf{a}}$ can be interpreted as the directional derivative of the vector field $\xi_{\mathbf{a}} \in \mathcal{T_{\mathbf{a}} M}$ in the direction of $\eta_{\mathbf{a}}  \in \mathcal{T_{\mathbf{a}} M}$. Note further that the Riemannian gradient $\nabla_{R} \phi( . )$ defined for all $\mathbf{a} \in \mathcal{M}$ is a vector field from $\mathcal{M}$ to its tangent bundle and, in this condition, the Riemannian Hessian $\nabla^{2}_{R} \phi( \mathbf{a}  ) \lbrack  \mathbf{b} \rbrack $ can thus be interpreted as the directional derivative of the Riemannian gradient of $\phi(.)$ at $\mathbf{a} \in \mathcal{M}$ in the direction of $\mathbf{b} \in \mathcal{T_{\mathbf{a}} M}$. See Section 3.5 of Boumal~\cite{B2023}, Chapter 3 of Robbin and Salomon~\cite{RS2022} or Section 5.3 of Absil et al.~\cite{AMS2008} for more information.

Furthermore, if $\phi(.)$ is a $C^{p}$ smooth mapping, it can be extended to a $C^{p}$ differentiable function $\hat{\phi}(.)$ on an open neighborhood $U$ of $\mathbb{R}^{k}$ such that $\mathcal{M} \subset U$ (see Proposition 3.31 of Boumal~\cite{B2023}) and if, in addition, we equip the submanifold $\mathcal{M}$ with the Euclidean metric of the ambient linear space on all its tangent spaces, we have the following relationships between the Riemannian gradient and Hessian of $\phi(.)$  with the Euclidean gradient of  $\hat{\phi}(.)$, respectively:
\begin{align}  \label{eq:D_rgrad_vec}
 \nabla_{R} \phi( \mathbf{a}  ) =  \mathbf{P}_{\mathcal{T_{\mathbf{a}} M}}  \big( \nabla \hat{\phi} (\mathbf{a})  \big) , \forall  \mathbf{a} \in \mathcal{M} \ .
\end{align}
and
\begin{align}  \label{eq:D_rhess_mat}
\nabla^{2}_{R} \phi( \mathbf{a}  ) \lbrack  \mathbf{b} \rbrack  & = \mathbf{P}_{\mathcal{T_{\mathbf{a}} M}}  \Big ( \mathit{J}   \big ( \nabla_{R} \phi( \mathbf{a}  )   \big) \lbrack \mathbf{b} \rbrack  \Big)  \nonumber \\
                                                                                                & =  \mathbf{P}_{\mathcal{T_{\mathbf{a}} M}}  \Big ( \mathit{J} \big (\mathbf{P}_{\mathcal{T_{\mathbf{a}} M}}  \nabla \hat{\phi} (\mathbf{a})  \big)  \lbrack \mathbf{b} \rbrack  \Big),  \forall  \mathbf{a} \in \mathcal{M}, \forall  \mathbf{b} \in \mathcal{T_{\mathbf{a}} M} \ ,
\end{align}
where $\mathbf{P}_{\mathcal{T_{\mathbf{a}} M}}$ denotes the orthogonal projector operator onto $\mathcal{T_{\mathbf{a}} M}$ in  $\mathbb{R}^{k}$ and $\mathit{J}   \big( \nabla_{R} \phi( \mathbf{a}  )  \big)$ is the usual Euclidean derivative (e.g., Jacobian matrix operator) of the Riemannian gradient of $\phi(.)$ at $\mathbf{a}$. In words, if the metric on $\mathcal{M}$ is inherited from the ambient Euclidean space, the Riemannian gradient is just the tangent space projection of the embedded gradient in the ambient space and the Levi-Civita connection on $\mathcal{M}$ is the tangent space projection of the Levi-Civita connection on the ambient space, which is equivalent to the Euclidean (directional) derivative.

Alternatively, again in the case of an embedded submanifold, the Riemannian Hessian of $\phi(.)$ can be defined by means of so-called second-order retractions, which are  second-order approximations of the exponential map, see Propositions 5.5.4 and 5.5.5 in~\cite{AMS2008} and Proposition 3 in~\cite{AM2012} for details. This also allows to derive the Riemannian Hessian of a cost function defined on an embedded submanifold in terms of standard Euclidean derivatives as in equation~\eqref{eq:D_rhess_mat}. See Appendix A of~\cite{V2012} for an illustration with the Riemannian Hessian of the the cost function $\varphi(.)$ used in the formulation~\eqref{eq:P0} of the WLRA problem in the case of binary weights and also Proposition 2 in~\cite{LSX2019} for a generalization to an arbitrary twice differentiable cost function $\varphi(.)$ defined on  the smooth matrix  submanifold $\mathbb{R}^{p  \times n}_{k}$ embedded in $\mathbb{R}^{p  \times n}$. These results are useful in our WLRA context and will be used later, see equation~\eqref{eq:D_rhess_varphi} in Subsection~\ref{landscape_wlra:box}.

Finally, the first- and second-order stationary conditions for a $C^{p}$ smooth real function $\phi(.)$ defined on a submanifold $\mathcal{M} \subset \mathbb{R}^{k}$ are exactly similar to their standard Euclidean counterparts when the search space is reduced to a linear subspace embedded in $\mathbb{R}^{k}$~\cite{HLWY2020} : a vector $\hat{\mathbf{a}} \in \mathcal{M}$ is a first-order critical point for $\phi(.)$ if the vector $\nabla_{R} \phi( \hat{\mathbf{a}}   )  \in \mathcal{T_{\hat{\mathbf{a}}} M}$ is equal to the zero-vector. Using equation~\eqref{eq:D_rgrad_vec}, this is equivalent to say that the usual Euclidean gradient of  the differentiable extension $\hat{\phi} (.)$ at $\hat{\mathbf{a}}$, $\nabla \hat{\phi} (\hat{\mathbf{a}})$, is orthogonal to $\mathcal{T_{\hat{\mathbf{a}}} M}$, e.g., that  $\nabla \hat{\phi} (\hat{\mathbf{a}}) \in \mathcal{N_{ \hat{\mathbf{a}}} M}$. Thus, $\hat{\mathbf{a}} \in \mathcal{M}$ is a first-order stationary point of the $C^{p}$ smooth real function $\phi(.)$ if one of the following equivalent conditions are satisfied
\begin{equation} \label{eq:D_first_kkt_m}
\nabla \hat{\phi} (\hat{\mathbf{a}}) \in \mathcal{N_{ \hat{\mathbf{a}}} M} \iff \mathbf{P}_{\mathcal{T_{\hat{\mathbf{a}} } M}}  \big( \nabla \hat{\phi} (  \hat{\mathbf{a}} ) \big) = \mathbf{0}^{k} \iff \Vert \mathbf{P}_{\mathcal{T_{\hat{\mathbf{a}} } M}}  \big( \nabla \hat{\phi} (  \hat{\mathbf{a}} ) \big) \Vert_{2} = 0 \ ,
\end{equation}
where $\hat{\phi} (.)$ is a differentiable extension in the ambient linear space of the smooth function $\phi(.)$ at $\hat{\mathbf{a}}$. On the other hand, a vector $\hat{\mathbf{a}} \in \mathcal{M}$ is a second-order critical point for $\phi(.)$ if it is a first-order critical point  for $\phi(.)$ and if, in addition, the self-adjoint operator $\nabla^{2}_{R} \phi( \hat{\mathbf{a}})$ defines a (symmetric) positive semi-definite bilinear form on $\mathcal{T_{ \hat{\mathbf{a}} } M}$. Finally, a vector $\hat{\mathbf{a}} \in \mathcal{M}$ is a strict (local) minimum of the smooth real function $\phi(.)$ if it is a first-order critical point  for $\phi(.)$ and if, in addition, the self-adjoint operator $\nabla^{2}_{R} \phi( \hat{\mathbf{a}}   )$ defines a (symmetric) positive definite bilinear form on $\mathcal{T_{ \hat{\mathbf{a}}  } M}$.

In the following, we will also be concerned with the minimization of a smooth real mapping $\phi(.)$ defined over a smooth submanifold $\mathcal{M} \subset \mathbb{R}^{k}$ (or $\mathcal{M} \subset \mathbb{R}^{p \times k}$), where $\phi(.)$  is invariant under the action of a certain group $\mathcal{G}$, which allows us to define an equivalence relation $\sim$ in the total computational space $\mathcal{M}$. In these conditions, all the elements of a given equivalence class of  $\sim$ have the same value for  $\phi(.)$. The quotient $\mathcal{M}/\sim$ generated by this equivalence relation consists of elements that are equivalence classes. If $\mathring{\mathbf{a}} \in \mathcal{M}/\sim$ then its vector representation in $\mathcal{M}$ is $\mathbf{a}$. Because of the invariance property, we want to minimize $\phi(.)$ over the set of equivalence classes $\mathcal{M}/\sim$ instead on $\mathcal{M}$. This leads to the notion of smooth and Riemannian quotient manifolds if some conditions on the group $\mathcal{G}$ are satisfied~\cite{AMS2008}\cite{B2023}. An important example of quotient manifolds is the Grassmann manifold which is the collection of all linear subspaces of a given dimension $k$ in a particular Euclidean space of dimension $n>k$ and is denoted by $\text{Gr}(n,k)$; see Chapter 9 of~\cite{B2023} for a comprehensive overview of general quotient manifolds and $\text{Gr}(n,k)$. More precisely, each of these linear subspaces can be represented by a $n \times k$ matrix of rank $k$ whose columns form a basis of this given subspace and all the $n \times k$ matrices of rank $k$ which are associated with the same subspace of rank $k$ form obviously an equivalence class, which can be identified with each subspace of rank $k$ embedded in the Euclidean space of dimension $n$. See Section~\ref{seppb:box} for more concrete examples of Grassmann manifolds in the context of the WLRA problem.

On such smooth quotient manifolds, the concept of tangent space to the quotient manifold $\mathcal{M}/\sim$ at $\mathring{\mathbf{a}} \in \mathcal{M}/\sim$  can be also defined and this abstract tangent space will be denoted by  $\mathcal{T}_{ \mathring{\mathbf{a}} } \mathcal{M}/\sim$ or simply by  $\mathcal{T}_{ \mathring{\mathbf{a}} } \mathcal{M}$ by an abuse of notation. Furthermore, the notions of Riemannian gradient and Hessian of the smooth mapping $\phi(.)$ defined (again with a slight abuse of notation) on $\mathcal{M}/\sim$ and such that $\phi(\mathring{\mathbf{a}} ) =  \phi( \mathbf{a} ) , \forall  \mathring{\mathbf{a}} \in \mathcal{M}/\sim$ with $\mathbf{a} \in  \mathring{\mathbf{a}} \subset \mathcal{M}$, can be extended. First- and second-order optimality conditions of  $\phi(.)$ for an element of  $\mathcal{M}/\sim$ can also be formulated. More detailed information on the related backgrounds can be found in Section 3.4 of Absil et al.~\cite{AMS2008} or in Section 9.8 of Boumal~\cite{B2023}. Comprehensive introduction to these abstract notions are also provided in~\cite{MMBS2012}\cite{MMBS2014}\cite{MS2016}\cite{BA2015}. Fortunately, when $\mathcal{M}$ is an embedded submanifold of $\mathbb{R}^{k}$ (or $\mathbb{R}^{p \times k}$) and inherits of the Euclidean (or Frobenius) metric of the ambient linear space, each abstract element of $\mathcal{T}_{ \mathring{\mathbf{a}} } \mathcal{M}/\sim$ (where $\mathring{\mathbf{a}} \in \mathcal{M}/\sim$ and $\mathbf{a} \in \mathcal{M}$) can be uniquely represented by an element of the tangent space $\mathcal{T_{\mathbf{a}}  M}$ whose direction in the total space $\mathcal{M}$  does not induce a displacement (from $\mathbf{a}$) along the equivalence class $\mathring{\mathbf{a}}$. This is achieved by decomposing the tangent space $\mathcal{T_{\mathbf{a}}  M}$ to the total space $\mathcal{M}$ at  $\mathbf{a}$ in the following  complementary and orthogonal direct sum
\begin{equation*}
\mathcal{T_{\mathbf{a}}  M} = \mathcal{H_{\mathbf{a}} M} \oplus \mathcal{V_{\mathbf{a}} M} \ ,
\end{equation*}
where $\mathcal{H_{\mathbf{a}} M}$ and $\mathcal{V_{\mathbf{a}} M}$ are orthogonal (with respect to the inner product of the ambient linear space) linear subspaces of  $\mathcal{T_{\mathbf{a}}  M}$.  $\mathcal{V_{\mathbf{a}} M}$ is called the vertical space of $\mathcal{M}$ at $\mathbf{a}$ and is the set of tangent vectors to $\mathcal{M}$ at  $\mathbf{a}$, which do induce a  displacement along the equivalence class $\mathring{\mathbf{a}}$. The horizontal space  $\mathcal{H_{\mathbf{a}} M}$ is the orthogonal complement of $\mathcal{V_{\mathbf{a}} M}$ and provides a valid and one-to-one representation of the abstract tangent vectors to the quotient space $\mathcal{M}/\sim$ at $\mathring{\mathbf{a}}$; see Section 9.4 of~\cite{B2023} for more information. Displacements in the vertical space leave the vector $\mathbf{a}$, representing the equivalence class  $\mathring{\mathbf{a}}$, unchanged. This justifies to restrict both tangent vectors and metric to the horizontal space $\mathcal{H_{\mathbf{a}} M}$~\cite{AMS2008}\cite{B2023}.

Provided that the inherited Euclidean metric defined in the total space $\mathcal{M}$ is invariant along the equivalence classes in $\mathcal{M}/\sim$, the quotient space $\mathcal{M}/\sim$ endowed  with this (Riemannian) metric is called a Riemannian quotient manifold of  $\mathcal{M}$~\cite{AMS2008}\cite{B2023}. For such Riemannian quotient manifold $\mathcal{M}/\sim$ whose total space  $\mathcal{M}$ is a submanifold embedded in a Euclidean space, we can then obtain convenient practical representations for the abstract Riemannian gradient and Hessian of $\mathring{\phi}(.)$ (defined on  $\mathcal{M}/\sim$) at $\mathring{\mathbf{a}}$ by simply replacing the tangent space $\mathcal{T_{\mathbf{a}}  M}$ by its horizontal space $\mathcal{H_{\mathbf{a}} M}$ in expressions~\eqref{eq:D_rgrad_vec} and~\eqref{eq:D_rhess_mat}:
\begin{align}  \label{eq:D_rgrad_vec2}
 \nabla_{R} \phi( \mathring{\mathbf{a}}  ) \simeq  \mathbf{P}_{\mathcal{H_{\mathbf{a}} M}}  \big(  \nabla \hat{\phi} (\mathbf{a}) \big) , \forall  \mathbf{a} \in \mathcal{M} \ .
\end{align}
and
\begin{align}  \label{eq:D_rhess_mat2}
\nabla^{2}_{R} \phi( \mathring{\mathbf{a}}  ) \lbrack  \bar{\mathbf{b}} \rbrack  \simeq \mathbf{P}_{\mathcal{H_{\mathbf{a}} M}}  \Big ( \mathit{J} \big (\mathbf{P}_{\mathcal{H_{\mathbf{a}} M}}  \nabla \hat{\phi} (\mathbf{a})  \big)  \lbrack \mathbf{b} \rbrack  \Big) ,  \forall  \mathbf{a} \in \mathcal{M}, \forall  \mathbf{b} \in \mathcal{H_{\mathbf{a}} M} \ ,
\end{align}
where $\mathbf{P}_{\mathcal{H_{\mathbf{a}} M}}$ denotes now the orthogonal projector operator onto $\mathcal{H_{\mathbf{a}} M}$ in the ambient linear space $\mathbb{R}^{k}$, $\bar{\mathbf{b}}$ is an abstract tangent vector of the quotient manifold $\mathcal{M}/\sim$ at $\mathring{\mathbf{a}} \in \mathcal{M}/\sim$, which is uniquely represented by the so-called horizontal lift $\mathbf{b} \in \mathcal{H_{\mathbf{a}} M}$, and $\hat{\phi} (.)$ is a $C^{p}$ differentiable extension of the smooth function $\phi(.)$ defined on $\mathcal{M}$ to an open neighborhood $U$ of $\mathbb{R}^{k}$ (or of $\mathbb{R}^{p \times k}$) such that $\mathcal{M} \subset U$. See Mishra et al.~\cite{MMBS2012}\cite{MMBS2014}\cite{MS2016} or Boumal and Absil~\cite{BA2015} for concrete illustrations of these abstract objects in the context of the WLRA problem. Importantly, the first- and second-order critical conditions for  $\phi(.)$  on the quotient manifold $\mathcal{M}/\sim$ can now be expressed and evaluated concretely in terms of $\mathbf{P}_{\mathcal{H_{\mathbf{a}} M}}  \big( \nabla \hat{\phi} (\mathbf{a}) \big)$ and $\mathbf{P}_{\mathcal{H_{\mathbf{a}} M}}  \Big ( \mathit{J} \big (\mathbf{P}_{\mathcal{H_{\mathbf{a}} M}}  \nabla \hat{\phi} (\mathbf{a})  \big)  \lbrack \mathbf{b} \rbrack  \Big)$ as for a "standard" submanifold embedded in a Euclidean linear space. As an illustration, the first-order stationary condition for $\mathring{\mathbf{a}} \in \mathcal{M}/\sim$ becomes
\begin{equation} \label{eq:D_first_kkt_qm}
 \mathbf{P}_{\mathcal{H_{\mathbf{a} } M}}  \big( \nabla \hat{\phi} (  \mathbf{a} ) \big) = \mathbf{0}^{k} \iff \Vert \mathbf{P}_{\mathcal{H_{\mathbf{a} } M}}  \big( \nabla \hat{\phi} (  \mathbf{a} ) \big) \Vert_{2} = 0 \ .
\end{equation}

The minimization of a real function $\phi(.)$ over a nonempty (arbitrary) set $\mathcal{K}  \subset \mathbb{R}^{k}$ (or more generally over a subset of a Euclidean or Frobenius linear space) is more involved than solving the same problem over the whole linear space  $\mathbb{R}^{k}$, or over one of its linear subspaces or one of its embedded smooth submanifolds described above. The first difficulty  comes up in characterizing the optimality of feasible solutions itself, e.g., the necessary first- and second-order conditions for a point $\hat{\mathbf{a}} \in \mathcal{K}$ to be a local minimizer of $\phi(.)$ over $\mathcal{K}$. Here, the feasible set $\mathcal{K}$ and its topological properties play a role as important as the properties of the function $\phi(.)$ itself as it is first necessary to characterize which search directions are admissible around $\hat{\mathbf{a}} \in \mathcal{K}$. It is well known now that these admissible directions are related to the notions of tangent and normal cones to the set $\mathcal{K}$  at $\hat{\mathbf{a}}$, see Chapter 6 of Rockafellar and Wets~\cite{RW1998} and also Ruszczynski~\cite{R2006}. Moreover, minimizing an even $C^{\infty}$ differentiable real function $\phi(.)$ over a nonempty subset $\mathcal{K}  \subset \mathbb{R}^{k}$ leads to different and confusing notions of stationarity~\cite{RW1998}\cite{HLU2019}\cite{LKB2023}\cite{OW2024}\cite{P2024}.

We now introduce the required elements of variational geometry to characterize the first- and second-order stationarity conditions of a possible local solution $\hat{\mathbf{a}}$ to the minimization of $\phi(.)$ over a nonempty (arbitrary) set $\mathcal{K}  \subset \mathbb{R}^{k}$.

\begin{def2.6} \label{def2.6:box}
A subset $\mathcal{C}  \subset \mathbb{R}^{k}$ is called a cone if it contains the zero-vector and contains with each of its vectors, all positive multiples of that vector, e.g., if $\mathbf{a}  \in \mathcal{C} \Rightarrow \lambda.\mathbf{a}  \in \mathcal{C} \text{ , } \forall  \lambda \in \mathbb{R}_{+*}$.
\end{def2.6}
As an illustration, the set consisting of a nonzero-vector $\mathbf{a}  \in \mathbb{R}^{k}$ and all of its positive multiples $\lambda.\mathbf{a}$ (with $\lambda \ge 0$) is a particular cone, which is called a ray. In other words, a cone, which is distinct from $\lbrace \mathbf{0}^{k} \rbrace$, is therefore composed of the union of the rays it contains.

Next, if $\mathcal{K}$ is a nonempty subset of $\mathbb{R}^{k}$, the (dual) polar of $\mathcal{K}$, noted $\mathcal{K}^{o}$, is the set
\begin{equation}  \label{eq:D_polar_cone}
\mathcal{K}^{o} := \Big \lbrace  \mathbf{b}  \in \mathbb{R}^{k} \text{ / } \langle \mathbf{a} , \mathbf{b}  \rangle_{2} \le 0 \text{ , } \forall \mathbf{a}  \in \mathcal{K} \Big \rbrace \ .
\end{equation}
First of all, we see that the polar of $\mathcal{K}$ depends on the scalar product used in $\mathbb{R}^{k}$, if we changed this scalar product then $\mathcal{K}^{o}$ is also changed. Geometrically, $\mathcal{K}^{o}$ is the set of all vectors in $\mathbb{R}^{k}$, which have an angle of at least 90° with every vector in $\mathcal{K}$. Next, note that $\mathcal{K}^{o}$ is a cone, as it obviously contains $\mathbf{0}^{k}$, but also $\lambda.\mathbf{b}$ for any $\lambda \ge 0$ if $\mathbf{b}  \in \mathcal{K}^{o}$. $\mathcal{K}^{o}$ is further convex and closed in $\mathbb{R}^{k}$ as the above definition of $\mathcal{K}^{o}$ expresses $\mathcal{K}^{o}$ as the intersection of a family of closed half-spaces, which are also all convex:
\begin{equation*}
\mathcal{K}^{o} = \bigcap_{\mathbf{a}  \in \mathcal{K} }  \Big \lbrace  \mathbf{b}  \in \mathbb{R}^{k} \text{ / } \langle \mathbf{a} , \mathbf{b}  \rangle_{2} \le 0  \Big \rbrace \ .
\end{equation*}
If  $\mathcal{K}$ is a nonempty subset of $\mathbb{R}^{k}$, its orthogonal complement is the set:
\begin{equation}  \label{eq:D_ortho_comp}
\mathcal{K}^{\bot} :=   \mathcal{K}^{o}  \cap (-\mathcal{K})^{o}  = \Big \lbrace  \mathbf{b}  \in \mathbb{R}^{k} \text{ / } \langle \mathbf{a} , \mathbf{b}  \rangle_{2} = 0 \text{ , } \forall \mathbf{a}  \in \mathcal{K} \Big \rbrace \ .
\end{equation}
We deduce immediately that $\mathcal{K}^{\bot}$ is  a closed convex cone of $\mathbb{R}^{k}$ as the intersection of two closed convex cones. Obviously, $\mathcal{K}^{\bot}$ is also a linear subspace of $\mathbb{R}^{k}$. Interestingly, if $\mathcal{K}$ is a linear subspace of $\mathbb{R}^{k}$, we have $\mathcal{K} = -\mathcal{K}$ and, consequently, $\mathcal{K}^{\bot} = \mathcal{K}^{o}$. Thus, polarity generalises  the notion of orthogonality between linear subspaces discussed in Subsection~\ref{lin_alg:box} to arbitrary nonempty subsets of $\mathbb{R}^{k}$. If $\mathcal{K}_{1}$ and $\mathcal{K}_{2}$ are two nonempty cones of $\mathbb{R}^{k}$, we have
\begin{equation*}
( \mathcal{K}_{1} \cup  \mathcal{K}_{2} )^{o} =(  \mathcal{K}_{1} +  \mathcal{K}_{2} )^{o} =  \mathcal{K}_{1}^{o}  \cap  \mathcal{K}_{2}^{o}  \ . 
\end{equation*}
In addition, if $\mathcal{K}_{1}$ and $\mathcal{K}_{2}$ are two closed convex cones then
\begin{equation*}
( \mathcal{K}_{1} \cap  \mathcal{K}_{2} )^{o} = \mathcal{K}_{1}^{o} +  \mathcal{K}_{2}^{o}  \  ,
\end{equation*}
and, finally, the property $\mathcal{K}^{oo} = \mathcal{K}$ is true if $\mathcal{K}$ is a closed convex cone. See Deutsch~\cite{D2012} for more details on (convex) cones and their polars.

We now introduce the general concepts of tangent and normal vectors at a nonempty set  $\mathcal{K}  \subset \mathbb{R}^{k}$, which generalize  the notions of tangent and normal vectors at a smooth submanifold of $\mathbb{R}^{k}$ introduced above, following~\cite{RW1998}; see also~\cite{HL2004} or~\cite{R2006} for a more gentle introduction to these concepts.
\begin{def2.7} \label{def2.7:box}
For a nonempty set $\mathcal{K}  \subset \mathbb{R}^{k}$ and a point $\bar{\mathbf{a}} \in \mathcal{K}$, a vector $\mathbf{d} \in \mathbb{R}^{k}$ is said to be tangent to $\mathcal{K}$ at $\bar{\mathbf{a}}$, when there exists a sequence $(\mathbf{a}_{k})_{k \in \mathbb{N}_{*}}$ in $\mathcal{K}$ tending to $\bar{\mathbf{a}}$ and a sequence $(\mathbf{t}_{k})_{k \in \mathbb{N}_{*}}$ in $\mathbb{R}_{+*}$ tending to zero (e.g., decreasing to zero) such that the vectors $\mathbf{b}_{k} = \frac{ (\mathbf{a}_{k} - \bar{\mathbf{a}}) } { \mathbf{t}_{k} }$ tend to $\mathbf{d}$, e.g., if
\begin{equation*}
\lim_{k \rightarrow \infty} \frac{ ( \mathbf{a}_{k} - \bar{\mathbf{a}} ) } {\mathbf{t}_{k}} = \mathbf{d} \ .
\end{equation*}
\end{def2.7}
Note that, if $\lim_{k \rightarrow \infty} \frac{ ( \mathbf{a}_{k} - \bar{\mathbf{a}} ) } {\mathbf{t}_{k}} = \mathbf{d}$, it is implicit that the sequence  $(\mathbf{a}_{k})_{k \in \mathbb{N}_{*}}$ tends to $\bar{\mathbf{a}}$, as otherwise the above limit does not exist as the sequence $(\mathbf{t}_{k})_{k \in \mathbb{N}_{*}}$ tends to zero. Consequently, some authors define a tangent vector without the condition that the sequence  $(\mathbf{a}_{k})_{k \in \mathbb{N}_{*}}$ tends to $\bar{\mathbf{a}}$. Furthermore, different, but equivalent, definitions of a tangent vector are also used in the literature, see Guignard~\cite{G1969}, Equation 2.2 of Schneider and Uschmajew~\cite{SU2015} and Section 5.1 of Hiriart-Urruty and Le Marechal~\cite{HL2004} for details.

This new definition of tangency generalizes the classical Definition~\ref{def2.5:box} in which a tangent vector $\mathbf{d}$ to $\mathcal{K}$ at  $\bar{\mathbf{a}}$ is the derivative at  $\bar{\mathbf{a}}$ of some curve drawn on $\mathcal{K}$. This classical definition is not relevant here as $\mathcal{K}$ can be a subset of $\mathbb{R}^{k}$ of discrete type and also because half-derivatives are key here instead of full-derivatives as in standard differential geometry.

We observe immediately that $\mathbf{0}^{k}$ is always  a tangent vector at $\mathcal{K}$ for any $\bar{\mathbf{a}} \in \mathcal{K}$: it suffices to take $\mathbf{a}_{k} = \bar{\mathbf{a}} , \forall k \in \mathbb{N}_{*}$. Furthermore, if $\mathbf{d}$ is a tangent vector to $\mathcal{K}$  at $\bar{\mathbf{a}}$, then $\alpha.\mathbf{d}$ for $\alpha > 0$ is also a tangent vector  to $\mathcal{K}$ at $\bar{\mathbf{a}}$ since it suffices to change  $\mathbf{t}_{k}$ to $\frac {\mathbf{t}_{k}} {\alpha}$ in the Definition~\ref{def2.7:box} of a tangent vector. In other words, the set of all tangent vectors to $\mathcal{K}$ at $\bar{\mathbf{a}}$ in the sense of Definition~\ref{def2.7:box} is a cone. The next theorem further shows that the set of all tangent vectors to $\mathcal{K}$ at $\bar{\mathbf{a}}$ is in fact a closed cone, which is called the tangent cone (or the contingent or Bouligand's cone) to $\mathcal{K}$ at $\bar{\mathbf{a}}$ and is denoted by $\mathcal{T^{B}_{\bar{\mathbf{a}}} K}$.

\begin{theo2.4} \label{theo2.4:box}
Let  $\mathcal{K}$ be a nonempty subset of $\mathbb{R}^{k}$ and let $\bar{\mathbf{a}} \in \mathcal{K}$. The set $\mathcal{T^{B}_{\bar{\mathbf{a}}} K}$ of all tangent directions for $\mathcal{K}$ at  $\bar{\mathbf{a}}$ in the sense of Definition~\ref{def2.7:box} is a closed cone.
\end{theo2.4}
\begin{proof}
Omitted. See  Lemma 3.12 of~\cite{R2006} or Proposition 5.1.3 of~\cite{HL2004}.
\end{proof}
Furthermore, it is not difficult to see that if $\bar{\mathbf{a}}$ is an interior point of $\mathcal{K}$ (e.g., $\bar{\mathbf{a}}  \in \mathring{\mathcal{K}}$), we have  $\mathcal{T^{B}_{\bar{\mathbf{a}}} K} = \mathbb{R}^{k}$. Thus, "the interesting" points are those on $bd(\mathcal{K})$, the boundary of $\mathcal{K}$. We next define the notion of normal vectors or directions to a set $\mathcal{K}$ in the regular sense following~\cite{RW1998}:
\begin{def2.8} \label{def2.8:box}
For a nonempty set $\mathcal{K}  \subset \mathbb{R}^{k}$ and a point  $\bar{\mathbf{a}} \in \mathcal{K}$, a vector $\mathbf{d} \in \mathbb{R}^{k}$ is said to be normal to $\mathcal{K}$ at $\bar{\mathbf{a}}$ in the regular sense, or a regular normal, if
\begin{equation*}
\langle \mathbf{d} , \mathbf{a}  - \bar{\mathbf{a}}  \rangle_{2} \le   \smallO( \Vert \mathbf{a}  - \bar{\mathbf{a}} \Vert_{2} ) \text{ , } \forall \mathbf{a}  \in \mathcal{K} \ ,
\end{equation*}
where we denote  by $\smallO( \Vert \mathbf{a}  - \bar{\mathbf{a}} \Vert_{2})$, for $\mathbf{a}  \in \mathcal{K}$, a term with the property that $\frac {\smallO( \Vert \mathbf{a}  - \bar{\mathbf{a}} \Vert_{2})} {\Vert \mathbf{a}  - \bar{\mathbf{a}} \Vert_{2}}$ tends to zero when  $\mathbf{a}$ tends to $\bar{\mathbf{a}}$ in $\mathcal{K}$, with $\mathbf{a} \ne \bar{\mathbf{a}}$.

The set of normal vectors to $\mathcal{K}$ at $\bar{\mathbf{a}}$ in the regular sense is called the Frechet normal cone to $\mathcal{K}$ at $\bar{\mathbf{a}}$ and is denoted by $\mathcal{N^{F}_{\bar{\mathbf{a}}} K}$.

\end{def2.8}
This name is justified by the following result, which provides a more comprehensive interpretation of the set of normal vectors in the regular sense to $\mathcal{K}$ at $\bar{\mathbf{a}}$.
\begin{theo2.5} \label{theo2.5:box}
Let  $\mathcal{K}$ be a nonempty subset of $\mathbb{R}^{k}$ and let $\bar{\mathbf{a}} \in \mathcal{K}$. The set  $\mathcal{N^{F}_{\bar{\mathbf{a}}} K}$ of all regular normal vectors is characterized by
\begin{equation*}
\mathbf{d} \in \mathcal{N^{F}_{\bar{\mathbf{a}}} K}  \Longleftrightarrow \langle \mathbf{d} , \mathbf{a}  \rangle_{2} \le 0  \text{ , } \forall \mathbf{a}  \in \mathcal{T^{B}_{\bar{\mathbf{a}}} K} \ .
\end{equation*}
In other words, we have $\mathcal{N^{F}_{\bar{\mathbf{a}}} K} = (\mathcal{T^{B}_{\bar{\mathbf{a}}} K})^{o}$ and the Frechet normal cone to $\mathcal{K}$ at $\bar{\mathbf{a}}$ is the polar of the Bouligand tangent cone to $\mathcal{K}$ at $\bar{\mathbf{a}}$ and is, thus, a closed convex cone.
\end{theo2.5}
\begin{proof}
Omitted. See  Propostion 6.5 in Rockafellar and Wets~\cite{RW1998}.
\end{proof}
Thus, the normal vectors to $\mathcal{K}$ at $\bar{\mathbf{a}}$ in the regular sense, apart from $\mathbf{0}^{k}$, are simply the vectors $\mathbf{d}$ of $\mathbb{R}^{k}$ that make a right or obtuse angle with every tangent vector $\mathbf{a}$ to $\mathcal{K}$ at $\bar{\mathbf{a}}$. Importantly,  if the subset $\mathcal{K}$ is an embedded submanifold of $\mathbb{R}^{k}$, the Bouligand tangent  and Frechet normal cones to $\mathcal{K}$ at $\bar{\mathbf{a}}$ reduce, respectively, to the tangent and normal spaces to $\mathcal{K}$ at $\bar{\mathbf{a}}$~\cite{RW1998}, e.g.,
\begin{equation*}
\mathcal{T^{B}_{\bar{\mathbf{a}}} K} = \mathcal{T_{\bar{\mathbf{a}}} K} \text{ and } \mathcal{N^{F}_{\bar{\mathbf{a}}} K} = \mathcal{N_{\bar{\mathbf{a}}} K} \ .
\end{equation*}
Thus, in a sense, the notions of  Bouligand tangent  and Frechet normal cones generalize the concepts of tangent and normal spaces to a smooth submanifold, described above, to an arbitrary nonempty set $\mathcal{K}$ embedded in a given Euclidean vector or Frobenius matrix space. Furthermore, we will see now that the first- and second-order optimality conditions for mimimizing a real function $\phi(.)$ over $\mathcal{K}$ can also be interpreted as an extension of the first- and second-order optimality conditions required over a  smooth submanifold discussed above.

The motivation and interest for the above paragraphs about cones, tangent and normal directions are related to this task and come from the following Theorem~\ref{theo2.6:box}, which provides a first basic first-order necessary condition for a vector $\hat{\mathbf{a}}$ to be a solution of  the minimization of a real function $\phi(.)$ over a nonempty (arbitrary) subset $\mathcal{K}  \subset \mathbb{R}^{k}$ (or more generally a subset of a normed vector space of finite dimension).

To be more precise, consider a nonempty set $\mathcal{K}  \subset \mathbb{R}^{k}$, a differentiable function $\phi(.) : \Omega \longrightarrow \mathbb{R}$, where $\Omega$ is open in $\mathbb{R}^{k}$ and such that $\mathcal{K}  \subset \Omega$, and the constrained optimization problem $\min_{\mathbf{a} \in \mathcal{K} } \, \phi( \mathbf{a} )$. Note that  we don't assume here that $\mathcal{K}$ is open, so if the constrained problem has a (local) solution $\hat{\mathbf{a}}$, this solution $\hat{\mathbf{a}}$ can be a boundary point of the feasible set $\mathcal{K}$, in which case the necessary conditions of optimality formulated above in equations~\eqref{eq:D_first_kkt} and~\eqref{eq:D_second_kkt} do not have to be satisfied because the perturbations $d\mathbf{a}$ to the vector $\hat{\mathbf{a}}$ such that $\hat{\mathbf{a}} + d\mathbf{a} \notin \mathcal{K}$ do not have to be taken into account and therefore they may correspond to a decrease of the cost function $\phi(.)$. In order to obtain a correct first-order necessary condition for optimality in a such case, the next theorem shows that we can restrict the set of possible perturbations $d\mathbf{a}$  to the tangent directions to $\mathcal{K}$ at $\hat{\mathbf{a}}$ in the sense of Definition~\ref{def2.7:box}, e.g., to the elements of the Bouligand's cone to  $\mathcal{K}$ at $\hat{\mathbf{a}}$.
\begin{theo2.6} \label{theo2.6:box}
Let  $\mathcal{K}$ be a nonempty subset of $\mathbb{R}^{k}$ and assume that $\phi(.)$ is a differentiable real function from an open subset $\Omega$ of $\mathbb{R}^{k}$ to $\mathbb{R}$ such that $\mathcal{K}  \subset \Omega$. If $\phi(.)$ has a local minimum over $\mathcal{K}$ at $\hat{\mathbf{a}}$, then $\phi(.)$ has not descent vector $\mathbf{d} \in \mathcal{T^{B}_{\hat{\mathbf{a}}} K}$, i.e.,
\begin{equation} \label{eq:D_first_kkt_bouligand}
\langle  \nabla \phi (\hat{\mathbf{a}}) , \mathbf{d}  \rangle_{2} \ge 0  \text{ , } \forall \mathbf{d}  \in \mathcal{T^{B}_{\hat{\mathbf{a}}} K}  \ ,
\end{equation}
which is equivalent to say that
\begin{equation} \label{eq:D_first_kkt_frechet}
- \nabla \phi (\hat{\mathbf{a}})  \in \mathcal{N^{F}_{\hat{\mathbf{a}}} K} =  (\mathcal{T^{B}_{\hat{\mathbf{a}}} K})^{o}  \ .
\end{equation}
In words, if a vector $\hat{\mathbf{a}}$ is a local minimizer of $\phi(.)$ over $\mathcal{K}$, the anti-gradient $-\nabla \phi (\hat{\mathbf{a}})$ is a normal vector in the regular sense to $\mathcal{K}$ at $\hat{\mathbf{a}}$, which is equivalent to say that $-\nabla \phi (\hat{\mathbf{a}})$ is an element of the Frechet normal cone to $\mathcal{K}$ at $\hat{\mathbf{a}}$.
\end{theo2.6}
\begin{proof}
See  Theorem 3.24 of Ruszczyinski~\cite{R2006}, Theorem 6.12 of Rockafellar and Wets~\cite{RW1998} or Theorem 1 of Guignard~\cite{G1969} for a proof.
\end{proof}
Thus, Theorem~\ref{theo2.6:box} and equation~\eqref{eq:D_first_kkt_frechet} provides a first-order optimality condition for the problem of minimizing $\phi(.)$ over $\mathcal{K}$ at a point $\hat{\mathbf{a}}$ and we will say that $\hat{\mathbf{a}}$ is a Frechet first-order stationarity point for this minimizaion problem if such condition is fulfilled. However, beware that other first-order optimality conditions have been proposed in the literature by replacing the Frechet normal cone $\mathcal{N^{F}_{\hat{\mathbf{a}}} K}$ in equation~\eqref{eq:D_first_kkt_frechet} by other cones like the so-called Mordukhovich or Clarke normal cones depending on the assumed properties for the function $\phi(.)$;  see~\cite{RW1998}\cite{HLU2019}\cite{LSX2019}\cite{P2024}\cite{OW2024} for more information. However, if we only assume that  $\phi(.)$ is a continuously differentiable or twice continuously differentiable function, the above Frechet stationarity provides the strongest necessary condition~\cite{LSX2019}\cite{OW2024} and this is the first-order optimality condition we shall use in this monograph.

We now derive a more convenient expression to check that a given point $\hat{\mathbf{a}} \in \mathcal{K}$  is a Frechet first-order stationary point based on the notion of metric projection onto an arbitrary nonempty subset $\mathcal{K} \subset \mathbb{R}^{k}$ (or more generally a subset of an arbitrary normed vector space), which generalizes the concept of an orthogonal projection operator onto a linear subspace introduced in Subsection~\ref{lin_alg:box}.

Let first $\mathcal{K}$ be a linear subspace of $\mathbb{R}^{k}$ and denote by $\text{Proj}_{\mathcal{K}} (. )$ the orthogonal projector mapping onto the subspace $\mathcal{K}$. $\text{Proj}_{\mathcal{K}} (. )$ is linear, idempotent ($\text{Proj}_{\mathcal{K}} o \text{Proj}_{\mathcal{K}} = \text{Proj}_{\mathcal{K}} $), non-expansive ($ \Vert \text{Proj}_{\mathcal{K}} (  \mathbf{a} )  \Vert_{2}  \le \Vert  \mathbf{a}  \Vert_{2} \ ,  \ \forall  \mathbf{a}  \in \mathbb{R}^{k}$ ) and it defines a direct sum of $\mathbb{R}^{k}$ as $\mathbf{a} = \text{Proj}_{\mathcal{K}} (\mathbf{a} ) + \text{Proj}_{\mathcal{K}^{\bot}} (\mathbf{a} )  \ ,  \  \forall \mathbf{a}  \in \mathbb{R}^{k}$.

We now generalize this operator to the case where $\mathcal{K}$ is only a nonempty closed and, eventually, convex set in $\mathbb{R}^{k}$. We will also see that, if $\mathcal{K}$ is in addition a cone in the sense of Definition~\ref{def2.6:box}, almost all the above properties of an orthogonal projector can be conserved or extended to the metric projection operator. Let us first define precisely the metric projection operator with the following definition.
\begin{def2.9} \label{def2.9:box}
Let $\mathcal{K}$  be a nonempty subset of $\mathbb{R}^{k}$ and $\mathbf{a} \in \mathbb{R}^{k}$. An element $\mathbf{b} \in \mathcal{K}$ is called a nearest point to $\mathbf{a}$ from $\mathcal{K}$ if
\begin{equation*}
 \Vert \mathbf{a} - \mathbf{b}  \Vert_{2} = d( \mathbf{a}, \mathcal{K} ) \ ,
\end{equation*}
where $d( \mathbf{a}, \mathcal{K} ) := \inf_{\mathbf{d} \in \mathcal{K} } \Vert \mathbf{a} - \mathbf{d}  \Vert_{2}$. The number $d( \mathbf{a}, \mathcal{K} )$ always exists and is called the distance from $\mathbf{a}$ to $\mathcal{K}$. Next, the possibly empty, discrete or infinite set of all nearest points from $\mathbf{a}$ to $\mathcal{K}$ is denoted by $P_{\mathcal{K}} ( \mathbf{a} )$. In other words,
\begin{equation*}
P_{\mathcal{K}} ( \mathbf{a} ) := \big \lbrace  \mathbf{b} \in \mathcal{K} \ /  \  \Vert \mathbf{a} - \mathbf{b}  \Vert_{2} = d( \mathbf{a}, \mathcal{K} )  \big \rbrace \ .
\end{equation*}
This defines a mapping $P_{\mathcal{K}} ( . )$ from $\mathbb{R}^{k}$ to the subsets of $\mathcal{K}$ called the metric projection onto $\mathcal{K}$.
\\
\end{def2.9}
If each $\mathbf{a} \in \mathbb{R}^{k}$ has at least (respectively, exactly) one nearest point in  $\mathcal{K}$, then $\mathcal{K}$ is called a proximinal (respectively, Chebyshev) set~\cite{D2012}. In other words,  $\mathcal{K}$ is  proximinal if $P_{\mathcal{K}} ( \mathbf{a} ) \ne \emptyset \ , \forall \mathbf{a} \in \mathbb{R}^{k}$ and is Chebyshev, if and only if, $P_{\mathcal{K}} ( \mathbf{a} ) =  \big \lbrace \mathbf{b}   \big \rbrace \ ,   \text{ with } \mathbf{b} \in \mathcal{K} \ , \forall \mathbf{a} \in \mathbb{R}^{k}$. In this last case, $P_{\mathcal{K}} ( . )$ can be viewed simply as a mapping from $\mathbb{R}^{k}$ to $\mathcal{K}$ in the usual sense. This will be for example the case if  $\mathcal{K}$ is a linear subspace of $\mathbb{R}^{k}$ (in which case $P_{\mathcal{K}} ( . )$ is simply the orthogonal projector $\text{Proj}_{\mathcal{K}} (. )$) or, more generally, if $\mathcal{K}$ is a closed convex set, as we will show shortly.

First,  if we assume that $\mathcal{K}$ is a nonempty closed subset of $\mathbb{R}^{k}$ then all points $\mathbf{a}  \in \mathbb{R}^{k}$ have at least one nearest point in $\mathcal{K}$. To see this, define a real function $f_{\mathbf{a}} (.)$ from $\mathbb{R}^{k}$ to $\mathbb{R}, \forall \mathbf{a}  \in \mathbb{R}^{k}$, by
\begin{equation*}
f_{\mathbf{a}} (\mathbf{b}) = \Vert \mathbf{b} - \mathbf{a} \Vert_{2} \ , \  \forall \mathbf{b}  \in \mathbb{R}^{k} \ ,
\end{equation*}
take a point $\mathbf{c} \in \mathcal{K}$ and define the sublevel set
\begin{equation*}
S_{\mathbf{c}} = \lbrace  \mathbf{b}  \in \mathbb{R}^{k}   \ /  \  f_{\mathbf{a}} (\mathbf{b}) \le  f_{\mathbf{a}} (\mathbf{c}) \rbrace \ .
\end{equation*}
$S_{\mathbf{c}}$ is a compact set of $\mathbb{R}^{k}$ as $f_{\mathbf{a}} (.)$ is continuous, $\lbrack  -\infty  \ ,  \  f_{\mathbf{a}} (\mathbf{c}) \rbrack$ is closed in $\mathbb{R}$ and $S_{\mathbf{c}}$ is bounded by definition. Then, we have obviously
\begin{equation*}
d( \mathbf{a}, \mathcal{K} ) =  \inf_{\mathbf{b} \in \mathcal{K}  \cap S_{\mathbf{c}} } f_{\mathbf{a}} (\mathbf{b}) \ ,
\end{equation*}
which has a solution in $\mathcal{K}$ as $f_{\mathbf{a}} (.)$ is continuous and $\mathcal{K}  \cap S_{\mathbf{c}}$ is compact (since  $\mathcal{K}  \cap S_{\mathbf{c}}$ is closed and bounded) in $\mathbb{R}^{k}$. This implies, the existence of, at least, one nearest point in $\mathcal{K}$ to $\mathbf{a}$ for all $\mathbf{a} \in \mathbb{R}^{k}$ if $\mathcal{K}$ is closed.

On the other hand, if $\mathcal{K}$ is a convex subset of $\mathbb{R}^{k}$, then $\forall \mathbf{a} \in \mathbb{R}^{k}$, $\mathbf{a}$ has at most one nearest point in $\mathcal{K}$. To demonstrate this claim suppose that $\mathcal{K}$ is convex and that $\mathbf{a} \in \mathbb{R}^{k}$ has two distinct nearest points in $\mathcal{K}$, say $\mathbf{b}_{1}$ and  $\mathbf{b}_{2}$. By using the parallelogram law with $\mathbf{d}_{1} = \mathbf{b}_{1} - \mathbf{a}$ and $\mathbf{d}_{2} = \mathbf{b}_{2} - \mathbf{a}$, we get
\begin{align*}
\Vert \mathbf{d}_{1} + \mathbf{d}_{2} \Vert^{2}_{2} + \Vert \mathbf{d}_{1} - \mathbf{d}_{2} \Vert^{2}_{2} & = 2.\Vert \mathbf{d}_{1} \Vert^{2}_{2} + 2.\Vert \mathbf{d}_{2} \Vert^{2}_{2}  \nonumber \\
& \Longrightarrow \Vert \frac{ \mathbf{b}_{1} + \mathbf{b}_{2} } {2} - \mathbf{a} \Vert^{2}_{2} = \Vert \mathbf{b}_{1}  - \mathbf{a}  \Vert^{2}_{2} - \frac{1} {4}  \Vert   \mathbf{b}_{1}  -  \mathbf{b}_{2} \Vert^{2}_{2} \  .
\end{align*}
Since $\mathcal{K}$ is convex, $\frac{ \mathbf{b}_{1} + \mathbf{b}_{2} } {2}$ belongs to $\mathcal{K}$ and we have $\Vert \frac{ \mathbf{b}_{1} + \mathbf{b}_{2} } {2} - \mathbf{a} \Vert^{2}_{2} <  \Vert \mathbf{b}_{1}  - \mathbf{a}  \Vert^{2}_{2}$, which contradicts the fact that $\mathbf{b}_{1}$ is a nearest point to $\mathbf{a}$  in $\mathcal{K}$.

In summary, if $\mathcal{K}$ is a nonempty closed and convex subset of $\mathbb{R}^{k}$, $P_{\mathcal{K}} ( \mathbf{a} ) = \lbrace \mathbf{c}  \rbrace$ with  $\mathbf{c} \in \mathcal{K}  \ ,  \  \forall \mathbf{a} \in \mathbb{R}^{k}$, and, by an abuse of notation, the metric projection defines effectively a simple metric projection mapping $P_{\mathcal{K}} ( \mathbf{a} ) =   \mathbf{c}$, which to each $\mathbf{a} \in \mathbb{R}^{k}$ associates its unique nearest point in $\mathcal{K}$. Interestingly, when $\mathcal{K}$ is a nonempty closed and convex set, the point $P_{\mathcal{K}} ( \mathbf{a} ) =   \mathbf{c}$ is equivalently characterized by the following property:
\begin{equation} \label{eq:metric_proj}
P_{\mathcal{K}} ( \mathbf{a} ) =   \mathbf{c} \Longleftrightarrow \langle \mathbf{a} -  \mathbf{c}, \mathbf{b} -  \mathbf{c}  \rangle_{2} \le 0  \text{ , } \forall \mathbf{b}  \in \mathcal{K} \  ,
\end{equation}
see Theorem 3.1.1 of Hiriart-Urruty and Le Marechal~\cite{HL2004} for a proof. This equivalence can be obviously restated with the help of the polar cone of the set $( \mathcal{K} - \mathbf{c})$ as
\begin{equation*} 
P_{\mathcal{K}} ( \mathbf{a} ) =   \mathbf{c} \Longleftrightarrow  \mathbf{a} -  \mathbf{c}  \in  ( \mathcal{K} - \mathbf{c})^{o} \  ,
\end{equation*}
which generalizes the property  $\mathbf{a} -  \text{Proj}_{\mathcal{K}} (\mathbf{a} )  \in \mathcal{K} ^{\bot}$ when  $\mathcal{K}$  is a  subspace of $\mathbb{R}^{k}$ and   $\text{Proj}_{\mathcal{K}} (.)$ is the orthogonal projector onto $\mathcal{K}$. Furthermore, when $\mathcal{K}$ is a nonempty closed and convex subset of $\mathbb{R}^{k}$, we have the following additional properties~\cite{HL2004}\cite{D2012}:

- the set $\lbrace \mathbf{a} \in \mathbb{R}^{k}   \  /  \  P_{\mathcal{K}} (\mathbf{a} ) =  \mathbf{a} \rbrace$  of fixed points of  $P_{\mathcal{K}} ( .)$ is $\mathcal{K}$ itself;

- the  metric projection mapping is idempotent, e.g., $P_{\mathcal{K}} o P_{\mathcal{K}} = P_{\mathcal{K}}$ and this justifies the term metric projection for $P_{\mathcal{K}} ( .)$;

- The metric projection mapping $P_{\mathcal{K}} ( .)$ is nonexpansive in the sense that $ \Vert P_{\mathcal{K}} (  \mathbf{a} ) -  P_{\mathcal{K}} (  \mathbf{b} ) \Vert_{2}  \le \Vert  \mathbf{a}  -  \mathbf{b} \Vert_{2} \ ,  \ \forall  \mathbf{a},  \mathbf{b }   \in \mathbb{R}^{k}$, implying that the metric projection mapping $P_{\mathcal{K}} ( .)$ is uniformly continuous on $\mathbb{R}^{k}$. Furthermore, if $\mathcal{K}$ is also a cone, $\mathbf{0}^{k} \in \mathcal{K}$ and we have $ \Vert P_{\mathcal{K}} (  \mathbf{a} )  \Vert_{2}  \le \Vert  \mathbf{a}  \Vert_{2} \ ,  \ \forall  \mathbf{a}  \in \mathbb{R}^{k}$, as for the orthogonal projector $\text{Proj}_{\mathcal{K}} (. )$ when $\mathcal{K}$ is a subspace of $\mathbb{R}^{k}$;

- and, finally, $P_{\mathcal{K}} ( .)$ is a linear operator if and only if $\mathcal{K}$ is a linear subspace of $\mathbb{R}^{k}$.

Suppose now that $\mathcal{K}$  is a nonempty subspace of $\mathbb{R}^{k}$. Then, $\mathcal{K}$ is a closed convex set and the metric projection operator $P_{\mathcal{K}} ( . )$ is well defined as a mapping from $\mathbb{R}^{k}$ to $\mathcal{K}$. However, we also know from the results of Subsection~\ref{lin_alg:box} that
\begin{equation*}
\inf_{\mathbf{d} \in \mathcal{K} } \Vert \mathbf{a} - \mathbf{d}  \Vert_{2} = \text{min}_{\mathbf{d} \in \mathcal{K} } \Vert \mathbf{a} - \mathbf{d}  \Vert_{2} = \Vert \mathbf{a} - \text{Proj}_{\mathcal{K}} (\mathbf{a} )  \Vert_{2} \ ,  \  \forall \mathbf{a} \in \mathbb{R}^{k} \ ,
\end{equation*}
where  $\mathcal{K}$ is a nonempty linear subspace of $\mathbb{R}^{k}$ and $\text{Proj}_{\mathcal{K}} (. )$ is the unique orthogonal projector operator onto $\mathcal{K}$. Consequently, as the metric projection operator $P_{\mathcal{K}} ( . )$ also solves uniquely this minimization problem in $\mathcal{K}$, we deduce immediately that $\text{Proj}_{\mathcal{K}} (\mathbf{a} ) = P_{\mathcal{K}} (\mathbf{a} ) , \forall \mathbf{a}  \in \mathbb{R}^{k}$. Thus, when $\mathcal{K}$  is a linear subspace, the metric projection mapping $P_{\mathcal{K}} ( .)$ is nothing else than  the orthogonal projector operator onto $\mathcal{K}$, $\text{Proj}_{\mathcal{K}} (. )$, suggesting again that we can interpret the metric projection mapping as an extension of the orthogonal projector mapping.

All these different properties confirm  that we can somehow interpret the metric projection mapping as an extension of an orthogonal projector mapping when the set of fixed points is a closed and convex subset rather than a linear subspace. Furthermore, we come even closer to an orthogonal projector, if we further assume that $\mathcal{K}$ is also cone, since in that case we have
\begin{equation*}
\mathbf{a} = P_{\mathcal{K}} (\mathbf{a} ) + P_{\mathcal{K}^{o}} (\mathbf{a} ) \text{ with } \langle P_{\mathcal{K}} (\mathbf{a} ) ,  P_{\mathcal{K}^{o}} (\mathbf{a} ) \rangle_{2} = 0 \ , \ \forall  \mathbf{a}  \in \mathbb{R}^{k} \ ,
\end{equation*}
which generalizes the canonical orthogonal decomposition $\mathbf{a} = \text{Proj}_{\mathcal{K}} (\mathbf{a} ) + \text{Proj}_{\mathcal{K}^{\bot}} (\mathbf{a} )$ when $\mathcal{K}$ is a subspace, see Section 3.2 of~\cite{HL2004} for details.

Since, we will mainly use the metric projection to project onto closed cones (e.g., the Bouligand's tangent cone to $\mathcal{K}$ at $\mathbf{a}$ when $\mathcal{K}$ is a nonempty, eventually closed, subset of  $\mathbb{R}^{k}$), we focus now specifically on the properties of the metric projection operator, which are still valid in this case.

First, note that if $\mathcal{C}$ is a closed subset of $\mathbb{R}^{k}$, the distance function defined as $d_{\mathcal{C}}  ( \mathbf{a} ) = d( \mathbf{a} , \mathcal{C} )$ from $\mathbb{R}^{k}$ to $\mathbb{R}$ is well-defined (since $\mathcal{C}$ is closed) and continuous on $\mathbb{R}^{k}$, see example 1.20 of Rockafellar and Wets~\cite{RW1998} for a proof. Next, $\forall  \mathbf{a}  \in \mathbb{R}^{k}$, the set $P_{\mathcal{C}} (\mathbf{a} ) = d_{\mathcal{C}}^{-1} ( \mathbf{a} )$ is nonempty (as shown above), bounded and closed, and thus compact in $\mathbb{R}^{k}$. It is closed as the reciprocal image of the singleton $\lbrace d( \mathbf{a} , \mathcal{C} ) \rbrace$ of $\mathbb{R}$ by the continuous distance function $d_{\mathcal{C}} ( . )$. It is bounded, because if we take a fixed point  $\mathbf{c}  \in \mathcal{C}$, we have, $\forall  \mathbf{b}  \in P_{\mathcal{C}} (\mathbf{a} )$, by definition, the inequality $  \Vert \mathbf{a} - \mathbf{b}  \Vert_{2} \le   \Vert \mathbf{a} - \mathbf{c}  \Vert_{2}$ and the distance of $\mathbf{b}$  to $\mathbf{a}, \forall  \mathbf{b}  \in P_{\mathcal{C}} (\mathbf{a} )$ is bounded by $\Vert \mathbf{a} - \mathbf{c}  \Vert_{2}$.

We next state the following Lemma, which will be useful to prove our next Theorem:
\begin{theo2.7} \label{theo2.7:box}
Let  $\mathcal{C}$ be a closed cone in $\mathbb{R}^{k}$. $\forall \mathbf{a} \in \mathbb{R}^{k}$ and $\forall \mathbf{b} \in P_{\mathcal{C}} (\mathbf{a} )$, we have
\begin{equation*}
\Vert \mathbf{b}  \Vert_{2} = \text{max} \big( 0,   \text{max}_{\mathbf{c} \in \mathcal{C} , \Vert \mathbf{c}  \Vert_{2} = 1}   \langle \mathbf{a}  , \mathbf{c}  \rangle_{2} \big) = \sqrt{  \langle \mathbf{a}  , \mathbf{b}  \rangle_{2}  } \ .
\end{equation*}
\end{theo2.7}
\begin{proof}
Omitted. See  Proposition A.6 of Levin et al.~\cite{LKB2023} for a proof.
\end{proof}

\begin{theo2.8} \label{theo2.8:box}
Let  $\mathcal{C}$ be a closed cone in $\mathbb{R}^{k}$. $\forall \mathbf{a} \in \mathbb{R}^{k}$ and $\forall \mathbf{b} \in P_{\mathcal{C}} (\mathbf{a} )$, we have
\begin{equation*}
\Vert \mathbf{b}  \Vert_{2}^{2}  = \Vert \mathbf{a}  \Vert_{2}^{2} - d( \mathbf{a} , \mathcal{C} )^{2} \ , 
\end{equation*}
implying that all the elements of $P_{\mathcal{C}} (\mathbf{a} )$ have the same length, and 
\begin{equation*}
\mathbf{a}  \in \mathcal{C}^{o}  \iff    P_{\mathcal{C}} (\mathbf{a} ) =  \lbrace \mathbf{0}^{k} \rbrace \ .
\end{equation*}
In words, if $\mathbf{a}$ belongs to the polar of the closed cone $\mathcal{C}$, its metric projection over $\mathcal{C}$, $P_{\mathcal{C}} (\mathbf{a} )$, is reduced to the zero-vector of the ambiant linear space and reciprocally.
\end{theo2.8}
\begin{proof}
First, we have the equalities
\begin{align*}
\Vert \mathbf{b}  \Vert_{2}^{2}  & = \Vert  \mathbf{a} - ( \mathbf{a} - \mathbf{b}  ) \Vert_{2}^{2}  \\
                                                 & =  \Vert  \mathbf{a} \Vert_{2}^{2} +  \Vert  \mathbf{a} - \mathbf{b}  \Vert_{2}^{2} - 2 \langle \mathbf{a}  , \mathbf{a} - \mathbf{b}  \rangle_{2}  \\
                                                 & =  \Vert  \mathbf{a} \Vert_{2}^{2} +  \Vert  \mathbf{a} - \mathbf{b}  \Vert_{2}^{2} - 2 \Vert  \mathbf{a} \Vert_{2}^{2} +  2 \langle \mathbf{a}  , \mathbf{b}  \rangle_{2}   \\
                                                 & =   \Vert  \mathbf{a} - \mathbf{b}  \Vert_{2}^{2} -  \Vert  \mathbf{a} \Vert_{2}^{2} +  2 \langle \mathbf{a}  , \mathbf{b}  \rangle_{2} \ .
\end{align*}
Now, since $\mathcal{C}$ is  a closed cone by hypothesis, using Lemma~\eqref{theo2.7:box}, we have $\Vert \mathbf{b}  \Vert_{2}^{2} =  \langle \mathbf{a}  , \mathbf{b}  \rangle_{2}$, from which we get
\begin{equation*}
\Vert \mathbf{b}  \Vert_{2}^{2}  =   \Vert  \mathbf{a} - \mathbf{b}  \Vert_{2}^{2} -  \Vert  \mathbf{a} \Vert_{2}^{2} +  2 \Vert  \mathbf{b} \Vert_{2}^{2} \ ,
\end{equation*}
which is equivalent after simplification to
\begin{equation*}
\Vert \mathbf{b}  \Vert_{2}^{2}  =  \Vert  \mathbf{a} \Vert_{2}^{2}  - \Vert  \mathbf{a} - \mathbf{b}  \Vert_{2}^{2}  =  \Vert  \mathbf{a} \Vert_{2}^{2}  - d( \mathbf{a} , \mathcal{C} )^{2} \ ,
\end{equation*}
as claimed in the theorem.

We now demonstrate the implication $\mathbf{a}  \in \mathcal{C}^{o}  \Rightarrow P_{\mathcal{C}} (\mathbf{a} ) =  \lbrace \mathbf{0}^{k} \rbrace$ , $\forall \mathbf{a} \in \mathbb{R}^{k}$. If $\mathbf{a}  \in \mathcal{C}^{o}$, for $\mathbf{c}  \in \mathcal{C}$, we have first
\begin{equation*}
 \Vert  \mathbf{a} - \mathbf{c}  \Vert_{2}^{2}  =  \Vert  \mathbf{a} \Vert_{2}^{2}  +  \Vert  \mathbf{c} \Vert_{2}^{2} - 2 \langle \mathbf{a}  , \mathbf{c}  \rangle_{2}
\end{equation*}
and, as  $\mathbf{a}  \in \mathcal{C}^{o}$, also the inequality $\langle \mathbf{a}  , \mathbf{c}  \rangle_{2} \le 0$.
This implies that the term $\Vert  \mathbf{c} \Vert_{2}^{2} - 2 \langle \mathbf{a}  , \mathbf{c}  \rangle_{2}$ is strictly positive, $\forall \mathbf{c}  \in \mathcal{C}  \backslash  \lbrace \mathbf{0}^{k} \rbrace$, and we get the inequality
\begin{equation*}
 \Vert  \mathbf{a} - \mathbf{c}  \Vert_{2}^{2}  >   \Vert  \mathbf{a} \Vert_{2}^{2} \ , \ \forall \mathbf{c}  \in \mathcal{C}  \backslash  \lbrace \mathbf{0}^{k} \rbrace \ ,
\end{equation*}
and also 
\begin{equation*}
 \Vert  \mathbf{a} - \mathbf{c}  \Vert_{2}  >   \Vert  \mathbf{a} \Vert_{2} \ , \ \forall \mathbf{c}  \in \mathcal{C}  \backslash  \lbrace \mathbf{0}^{k} \rbrace \ ,
\end{equation*}
after simplification. In other words, $\mathbf{0}^{k}  \in \mathcal{C}$ is the unique nearest point in $\mathcal{C}$ to $\mathbf{a}$, e.g., $P_{\mathcal{C}} (\mathbf{a} ) =  \lbrace \mathbf{0}^{k} \rbrace$ if $\mathbf{a}  \in \mathcal{C}^{o}$, as claimed above.

Reciprocally, we now demonstrate the implication $\mathbf{0}^{k} \in  P_{\mathcal{C}} (\mathbf{a} )  \Rightarrow \mathbf{a}  \in \mathcal{C}^{o}$, $\forall \mathbf{a} \in \mathbb{R}^{k}$. First, as $\mathcal{C}$ is a closed cone by hypothesis, using the first assertion of the Theorem  demonstrated above, we have $\Vert \mathbf{a} \Vert_{2} = d( \mathbf{a} , \mathcal{C} )$ and, $\forall \mathbf{c} \in P_{\mathcal{C}} (\mathbf{a} )$, we have $\Vert \mathbf{c} \Vert_{2} = 0$ and, thus, $P_{\mathcal{C}} (\mathbf{a} ) =  \lbrace \mathbf{0}^{k} \rbrace$. In other words, $\mathbf{0}^{k}$ is the unique nearest point to $\mathbf{a}$ in $\mathcal{C}$. Furthermore, from Lemma~\eqref{theo2.7:box}, we have also 
\begin{equation*}
 \text{max}_{\mathbf{c} \in \mathcal{C} , \Vert \mathbf{c}  \Vert_{2} = 1}   \langle \mathbf{a}  , \mathbf{c}  \rangle_{2}  \le 0 \ .
\end{equation*}
In order to demonstrate that $\mathbf{a}  \in \mathcal{C}^{o}$, e.g., that $\langle \mathbf{a}  , \mathbf{d}  \rangle_{2}  \le 0 \ , \  \forall \mathbf{d}  \in \mathcal{C}$, we now proceed by contradiction. Suppose that it exists $\mathbf{d}  \in \mathcal{C}$ such that  $\langle \mathbf{a}  , \mathbf{d}  \rangle_{2}  > 0$. Then, $\mathbf{d}  \ne \mathbf{0}^{k}$ and $\Vert  \mathbf{d} \Vert_{2} \ne 0$, and we have
\begin{equation*}
 \text{max}_{\mathbf{c} \in \mathcal{C} , \Vert \mathbf{c}  \Vert_{2} = 1}   \langle \mathbf{a}  , \mathbf{c}  \rangle_{2}  \le 0 <   \langle \mathbf{a}  , \frac{ \mathbf{d} } {\Vert  \mathbf{d} \Vert_{2}} \rangle_{2} \ ,
\end{equation*}
which is a contradiction, since $\frac{ \mathbf{d} } {\Vert  \mathbf{d} \Vert_{2}}  \in \mathcal{C}$, because $\mathcal{C}$ is a cone, and $\Vert  \frac{ \mathbf{d} } {\Vert  \mathbf{d} \Vert_{2}} \Vert_{2} = 1$.

Summarizing, we have, $\forall \mathbf{a} \in \mathbb{R}^{k}$,
\begin{equation*}
\mathbf{a}  \in \mathcal{C}^{o}  \iff    P_{\mathcal{C}} (\mathbf{a} ) =  \lbrace \mathbf{0}^{k} \rbrace \ ,
\end{equation*}
as claimed in the Theorem and we are done.
\end{proof}

Now, we can come back to our problem of reformulating the first-order stationarity condition~\eqref{eq:D_first_kkt_frechet} for a point $\hat{ \mathbf{a} } \in \mathcal{K}$, $\mathcal{K}$ being an arbitrary subset of $\mathbb{R}^{k}$, to be a (local) minimizer of a real function $\phi(.)$ differentiable over an open neighborhood of $\mathcal{K}$. An application of Theorem~\ref{theo2.8:box} to the anti-gradient $-\nabla \phi (\hat{\mathbf{a}})$ and the tangent Bouligand's cone $\mathcal{T^{B}_{\hat{\mathbf{a}}} K}$, which is a closed cone, leads to the following equivalent first-order critical conditions
\begin{equation} \label{eq:D_first_kkt_c}
- \nabla \phi (\hat{\mathbf{a}})  \in \mathcal{N^{F}_{\hat{\mathbf{a}}} K} =  (\mathcal{T^{B}_{\hat{\mathbf{a}}} K})^{o}   \iff  P_{\mathcal{T^{B}_{\hat{\mathbf{a}}} K}} ( - \nabla \phi (\hat{\mathbf{a}}) ) =  \lbrace \mathbf{0}^{k} \rbrace \ .
\end{equation}
Moreover, by a small abuse of notation, we can write $P_{\mathcal{T^{B}_{\hat{\mathbf{a}}} K}} ( - \nabla \phi (\hat{\mathbf{a}}) ) =  \mathbf{0}^{k}$ and the first-order stationarity condition for $\hat{\mathbf{a}} \in \mathcal{K}$ to be a (local) minimizer becomes
\begin{equation} \label{eq:D_first_kkt_c2}
 \Vert  P_{\mathcal{T^{B}_{\hat{\mathbf{a}}} K}} ( - \nabla \phi (\hat{\mathbf{a}}) )  \Vert_{2} = 0 \ ,
\end{equation}
where $P_{\mathcal{T^{B}_{\hat{\mathbf{a}}} K}} ( - \nabla \phi (\hat{\mathbf{a}}) )$ designs now any of its elements since they have all the same length according to Theorem~\ref{theo2.8:box}. Note the similarity of this first-order condition~\eqref{eq:D_first_kkt_c}  or~\eqref{eq:D_first_kkt_c2} with the one stated above in equation~\eqref{eq:D_first_kkt_m} in the case where $\mathcal{K}$ is an embedded smooth submanifold of $\mathbb{R}^{k}$.

To conclude these paragraphs on optimality conditions for a real function $\phi(.)$ at a point $\hat{\mathbf{a}}  \in \mathcal{K}$, where $\mathcal{K}$ is an nonempty arbitrary subset of $\mathbb{R}^{k}$, we now recall in the following theorem the necessary second-order condition for a point $\hat{\mathbf{a}}  \in \mathcal{K}$ to be a (local) minimizer over $\mathcal{K}$ of a cost function $\phi(.)$ twice continuously differentiable over an open neighborhood of $\mathcal{K}$ in  $\mathbb{R}^{k}$.

\begin{theo2.9} \label{theo2.9:box}
Let  $\mathcal{K}$ be a nonempty subset of $\mathbb{R}^{k}$ and assume that $\phi(.)$ is a real function twice continuously differentiable over an open neighborhood of $\mathcal{K}$ in  $\mathbb{R}^{k}$ and that $\hat{\mathbf{a}}  \in \mathcal{K}$ is a (local) minimizer of $\phi(.)$ over $\mathcal{K}$. Then, for every $\mathbf{d}  \in \mathcal{T^{B}_{\hat{\mathbf{a}}} K}$ satisfying $ \langle \nabla \phi (\hat{\mathbf{a}})  , \mathbf{d}  \rangle_{2} = 0$ we have
\begin{equation} \label{eq:D_second_kkt_c}
 \langle \nabla \phi (\hat{\mathbf{a}})  , \mathbf{c}  \rangle_{2} +  \langle \big \lbrack \nabla^{2} \phi (\hat{\mathbf{a}}) \big \rbrack ( \mathbf{d} ) , \mathbf{d}  \rangle_{2} \ge 0   \ ,  \   \forall \mathbf{c}   \in \mathcal{T^{B}_{ (\hat{\mathbf{a}}, \mathbf{d}) } K} \ ,
\end{equation}
where $\mathcal{T^{B}_{\hat{\mathbf{a}}} K}$ is the Bouligand tangent cone to $\mathcal{K}$  at $\hat{\mathbf{a}}$ and $\mathcal{T^{B}_{ (\hat{\mathbf{a}}, \mathbf{d}) } K}$ is the second-order (Bouligand) tangent set to $\mathcal{K}$  at $\hat{\mathbf{a}}$ in the direction of $\mathbf{d}  \in \mathcal{T^{B}_{\hat{\mathbf{a}}} K}$ (see Definition 3.41 in Ruszczynski~\cite{R2006} for a precise definition of this second-order tangent set). Note, however, that $\mathcal{T^{B}_{ (\hat{\mathbf{a}}, \mathbf{d}) } K}$ is not a cone in general, nor it is convex.
\end{theo2.9}
\begin{proof}
Omitted. See  Theorem 3.45 of Ruszczynski~\cite{R2006} for a proof.
\end{proof}
Using Theorem~\ref{theo2.9:box}, we will say that $\hat{\mathbf{a}}  \in \mathcal{K}$ is a (Frechet) second-order stationarity point of $\phi(.)$ over $\mathcal{K}$ if it is a (Frechet) first-order stationarity point for $\phi(.)$  and if, in addition, the condition~\ref{eq:D_second_kkt_c} is fulfilled.

In the following sections, we will also manipulate (differentiable) scalar, vector or matrix functions with a matrix argument $\mathbf{A} \in  \mathbb{R}^{p  \times n}$. As an illustration, let $\phi(.)$ be a scalar function defined on $\mathbb{R}^{p  \times n}$. If $\mathbb{R}^{p  \times n}$ is equipped  with its usual Frobenius inner product, the gradient of $\phi(.)$ at a matrix variable  $\mathbf{A} \in  \mathbb{R}^{p  \times n}$ is also a $p  \times n$ matrix, i.e.,
\begin{equation} \label{eq:D_grad_mat}
\big \lbrack  \nabla \phi( \mathbf{A} ) \big \rbrack_{ij}  =   \frac{ \partial   \phi(\mathbf{A})  }  {\partial\mathbf{A}_{ij} } \text{ for } i = 1, \cdots , p \text{ ; }   j = 1, \cdots , n \ .
\end{equation}
Alternatively, we can interpret this gradient as a linear form  $  \big( \nabla \phi( \mathbf{A} ) \big) \in \pounds ( \mathbb{R}^{p  \times n} , \mathbb{R} )$ defined by
\begin{equation*}
\big( \nabla \phi( \mathbf{A} )  \big) ( \mathbf{C} ) =  \langle \nabla \phi( \mathbf{A} )  \; , \; \mathbf{C}   \rangle_{F} =  \Tr \big( \nabla \phi( \mathbf{A} )^{T} \mathbf{C}  \big) = \sum_{i=1}^p { \sum_{j=1}^n { \frac{ \partial   \phi(\mathbf{A})  }  {\partial\mathbf{A}_{ij} }  \mathbf{C}_{ij} } } \text{ , } \forall \mathbf{C} \in  \mathbb{R}^{p  \times n} \ .
\end{equation*}
On the other hand, the Hessian of $\phi(.)$  at $\mathbf{A} \in  \mathbb{R}^{p  \times n}$ can be viewed as a $4^{th}$ order tensor of dimension $p  \times n  \times p  \times n$, instead of a  symmetric matrix (see equation~\eqref{eq:D_hess_mat}) as in the case of a vector argument, which is equal to
\begin{equation} \label{eq:D_hess_mat3}
\big \lbrack  \nabla^2  \phi( \mathbf{A} ) \big \rbrack_{ijkl}  =   \frac{ \partial^2   \phi(\mathbf{A})  }  {\partial\mathbf{A}_{ij} \partial\mathbf{A}_{kl} } \text{ for } i = 1, \cdots , p \text{ ; }   j = 1, \cdots , n \text{ ; }  k = 1, \cdots , p \text{ ; }   l = 1, \cdots , n \ .
\end{equation}
Equivalently, we can view $\nabla^2  \phi( \mathbf{A} )$ as a bilinear form $\big(  \nabla^2  \phi( \mathbf{A} ) \big)$, from $\mathbb{R}^{p  \times n} \times \mathbb{R}^{p  \times n}$ to $\mathbb{R}$, defined by
\begin{equation*}
\big(  \nabla^2  \phi( \mathbf{A} ) \big) ( \mathbf{C},  \mathbf{D} ) =   \sum_{i,j,k,l}   { \frac{ \partial^2   \phi(\mathbf{A})  }  {\partial\mathbf{A}_{ij} \partial\mathbf{A}_{kl} }  \mathbf{C}_{ij}  \mathbf{D}_{kl} } \text{ , } \forall \mathbf{C} , \mathbf{D} \in  \mathbb{R}^{p  \times n} \ .
\end{equation*}
Finally, another very useful representation of $\nabla^2  \phi( \mathbf{A} )$, implicit in the preceding one, is as a huge $p.n  \times p.n$ symmetric matrix
\begin{equation*}
\big \lbrack  \nabla^2  \phi( \mathbf{A} ) \big \rbrack_{ij}  =  \frac{ \partial^2   \phi(\mathbf{A})  }  {\partial\mathbf{a}_{i} \partial\mathbf{a}_{j} }  \text{ for } i = 1, \cdots , p.n \text{ ; }   j = 1, \cdots , p.n \ ,
\end{equation*}
where $\mathbf{a}_{i}$ is the $i^{th}$ element of a vectorized form of $\mathbf{A}$, e.g., $\mathbf{a} = \emph{vec}( \mathbf{A} )$ or $\mathbf{a} = \emph{vec}( \mathbf{A}^{T} )$. For example, in Subsection~\ref{hess:box}  we will derive the Hessian of a real (variable projection) functional $\psi (.)$ of the matrix variable  $\mathbf{A} \in  \mathbb{R}^{p  \times k}$ (defined in the next section) using this specific representation.

The first and second derivatives of a matrix function $f(.)$ from $\mathbb{R}^{p  \times n}$ to $\mathbb{R}^{q  \times r}$ can also be viewed as higher order tensors. However, it is generally more convenient to represent them as linear or multi-linear operators~\cite{C2017}. For example, the first derivative of  $f(.)$ at $\mathbf{A} \in \mathbb{R}^{p  \times n}$ is a linear operator from $\mathbb{R}^{p  \times n}$ to $ \mathbb{R}^{q  \times r}$, e.g., $\mathit{D}f(\mathbf{A}) \in \pounds ( \mathbb{R}^{p  \times n}, \mathbb{R}^{q \times r} )$, and the second derivative of $f(.)$ at $\mathbf{A}$,  $\mathit{D}^{2} f(\mathbf{A})$, is an element of $\pounds \big(  \mathbb{R}^{p  \times n}, \pounds ( \mathbb{R}^{p  \times n}, \mathbb{R}^{q \times r} ) \big)$, which is isomorphic to $\pounds (  \mathbb{R}^{p  \times n}, \mathbb{R}^{p  \times n} ;\mathbb{R}^{q \times r} )$, the set of bilinear maps from $\mathbb{R}^{p  \times n}$ into $\mathbb{R}^{q \times r}$~\cite{C2017}. Thus, $\mathit{D}^{2} f(\mathbf{A})$ can be interpreted as a bilinear operator from $\mathbb{R}^{p  \times n} \times  \mathbb{R}^{p  \times n}$ to  $\mathbb{R}^{q  \times r}$. In this way,  the Hessian of  a scalar function $\phi(.)$ with a matrix argument $\mathbf{A} \in  \mathbb{R}^{p  \times n}$ discussed above is the first derivative of its gradient, which is a mapping from $\mathbb{R}^{p  \times n}$ to $\mathbb{R}^{p  \times n}$, and, thus, this Hessian can be viewed as a mapping from $\mathbb{R}^{p  \times n}$ to  $\pounds ( \mathbb{R}^{p  \times n},  \mathbb{R}^{p  \times n} )$, e.g., for $\mathbf{A} \in \mathbb{R}^{p  \times n}$, $\big \lbrack  \nabla^2  \phi( \mathbf{A} ) \big \rbrack \in  \pounds ( \mathbb{R}^{p  \times n},  \mathbb{R}^{p  \times n} )$ and is a linear operator  from $\mathbb{R}^{p  \times n}$ to $\mathbb{R}^{p  \times n}$. Furthermore, we can identify  the bilinear form $\big( \nabla^2  \phi( \mathbf{A} ) \big)$  with $\big \lbrack  \nabla^2  \phi( \mathbf{A} ) \big \rbrack$~\cite{C2017} and they verify the equality
\begin{equation*}
\langle  \big \lbrack  \nabla^2  \phi( \mathbf{A} ) \big \rbrack  ( \mathbf{C} )  \; , \;  \mathbf{D}  \rangle_{F} = \big( \nabla^2  \phi( \mathbf{A} ) \big) ( \mathbf{C},  \mathbf{D} )  \text{ , } \forall \mathbf{C} , \mathbf{D} \in  \mathbb{R}^{p  \times n}.
\end{equation*}

This identification of $\big( \nabla^2  \phi( \mathbf{A} ) \big)$ with $\big \lbrack  \nabla^2  \phi( \mathbf{A} ) \big \rbrack$ can be very useful in practice as evaluating directly  $\big \lbrack  \nabla^2  \phi( \mathbf{A} ) \big \rbrack  ( \mathbf{C} )$ (e.g., the directional derivative of the gradient of $\phi(.)$ in the direction of $\mathbf{C}$) can be much cheaper and efficient than computing  analytically the full Hessian $\nabla^2  \phi( \mathbf{A} )$. This is for example the approach followed by Boumal and Absil~\cite{BA2011}\cite{BA2015} in their Newton Riemannian trust-region method for solving the WLRA problem in a Grassmann manifold framework (recall that a Grassmann manifold is the collection of all linear subspaces of a given dimension in a particular Euclidean space as already discussed above).

Keep also in mind that all the above notions of a smooth function, smooth manifold, tangent space to a smooth manifold, tangent and normal cones to an arbitrary subset and metric projection onto an arbitrary subset can be defined without any difficulties in the case when the ambient linear space is $\mathbb{R}^{p  \times k}$ instead of $\mathbb{R}^{p}$ if the linear space $\mathbb{R}^{p  \times k}$ is equipped with the standard Frobenius inner product~\cite{LSX2019}. Moreover, the linear spaces $\mathbb{R}^{p  \times k}$  and $\mathbb{R}^{p.k}$ are isomorphic and the Frobenius metric on $\mathbb{R}^{p  \times k}$ is equivalent to the standard Euclidean metric on $\mathbb{R}^{p.k}$ thanks to this isomorphism.

We conclude that preliminary section by a few more definitions about nonlinear optimization, which will be useful for our next sections.

A function $\phi(.)$ is said to be nonlinear in some scalar parameter $\alpha$, vector parameter $\mathbf{a}$ or matrix parameter $\mathbf{A}$ if the derivatives $\frac{ \partial \phi(.)} {\partial \alpha}$,  $\frac{ \partial \phi(.)} {\partial\mathbf{a}}$ and $\frac{ \partial \phi(.)} {\partial\mathbf{A}}$ are functions of $\alpha$,  $\mathbf{a}$ and $\mathbf{A}$, respectively~\cite{HPS2012}. As an illustration, let $\mathbf{r}(.)$ be a real-vector function from $\mathbb{R}^{k}$ into  $\mathbb{R}^{q}$ and further assume that $\mathbf{r}(.)$  is at least twice continuously differentiable. Then, the real function $\phi(.)$ from $\mathbb{R}^{k}$ into  $\mathbb{R}$ defined by
\begin{equation} \label{def:NLLS_func}
\phi( \mathbf{a} )  =  \frac{1}{2} \Vert \mathbf{r}( \mathbf{a} )  \Vert^{2}_{2} =  \frac{1}{2} \mathbf{r}( \mathbf{a} )^{T} \mathbf{r}( \mathbf{a} )   \text{   for   } \mathbf{a} \in  \mathbb{R}^{k}
\end{equation}
is called  a Non-Linear Least-Squares (NLLS) functional.
If we differentiate $\phi(.)$ with respect to $\mathbf{a}  \in  \mathbb{R}^{k}$ (e.g., we compute its gradient at $\mathbf{a} $) and equate the derivative to zero, this leads to the following equation
\begin{equation} \label{eq:gradvec_nlls}
\nabla \phi( \mathbf{a}  ) =  \mathit{J} \big( \mathbf{r}( \mathbf{a}  ) \big)^{T} \mathbf{r}( \mathbf{a}  )  = \mathbf{0}^{k}  \  ,
\end{equation}
which may be used in practice to test the convergence of NLLS iterative algorithms employed for minimizing $\phi(.)$ over $\mathbb{R}^{k}$~\cite{OR1970}\cite{DS1983}\cite{MN2010}. This last equation shows that the vector $\mathbf{r}( \mathbf{a}  )$ is orthogonal to $\emph{ran}(  \mathit{J}( \mathbf{r}( \mathbf{a}  ) ) )$, the linear subspace spanned by the columns of the Jacobian matrix of the real $q$-vector function $\mathbf{r}(.)$  at $\mathbf{a}$, if $\mathbf{a}$ is a stationary point of $\phi(.)$. Furthermore, if $\phi(.)$ is a NLLS functional then its Hessian matrix is
\begin{equation}  \label{eq:hessian_nlls}
\nabla^2 \phi( \mathbf{a}  ) = \mathit{J} \big( \mathbf{r}(\mathbf{a} ) \big)^{T} \mathit{J} \big( \mathbf{r}(\mathbf{a} ) \big) + \sum_{l=1}^{q} \mathbf{r}_l (\mathbf{a} ) \nabla^2 \mathbf{r}_l (\mathbf{a}  ) \  ,
\end{equation}
where $\nabla^2 \mathbf{r}_l (\mathbf{a}  )$ is the Hessian matrix of the $l^{th}$ component of the $q$-vector function $\mathbf{r}(.)$ at $\mathbf{a} $ (i.e., $\mathbf{r}_l (\mathbf{a} )$)  given by
\begin{equation*}
   \big \lbrack \nabla^2 \mathbf{r}_l (\mathbf{a}  )  \big \rbrack_{ij} =   \frac{ \partial^2  \mathbf{r}_l (\mathbf{a} ) }{ \partial \mathbf{a}_i  \partial \mathbf{a}_j } \text{ for } i=1, \cdots, k \text{ and } j = 1, \cdots , k \  ,
\end{equation*}
for $l = 1, \cdots , q$. Note that the factor $\frac{1}{2}$  in the  definition~\ref{def:NLLS_func} of the NLLS functional $\phi(.)$ has been introduced here only for notational convenience as without it a factor $2$ will appear in the two preceding equations defining  $\nabla \phi(.)$ and $\nabla^2 \phi(.)$ and in many equations of this paper. Furthermore, the second-order Taylor expansion of the NLLS functional $\phi(.)$ at a point $\mathbf{a}  \in \mathbb{R}^{k}$ has the following form
\begin{align*}
\phi( \mathbf{a}  + d\mathbf{a} ) & =  \phi( \mathbf{a}  ) + d\mathbf{a}^{T}  \mathit{J}( \mathbf{r}( \mathbf{a}  ) )^{T} \mathbf{r}( \mathbf{a}  ) \text{ } +  \\
&     \frac{1}{2} d\mathbf{a}^{T} \left(   \mathit{J}( \mathbf{r}(\mathbf{a} ) )^{T} \mathit{J}( \mathbf{r}(\mathbf{a} ) )  + \sum_{l=1}^{q} \mathbf{r}_l (\mathbf{a} ) \nabla^2 \mathbf{r}_l (\mathbf{a}  ) \right) d\mathbf{a} +  \mathcal{O}( \Vert d\mathbf{a} \Vert^{3}_{2} )  \  .
\end{align*}
These special forms of the gradient, Hessian and Taylor expansion of $\phi(.)$ are exploited by methods for solving NLLS problems, see Subsection~\ref{opt:box} and~\cite{DS1983}\cite{MN2010}\cite{HPS2012} for details.

Finally,  we give the following definition, which will be also useful in the next sections:
 \\
\begin{def2.10} \label{def2.10:box}
Let $m, n, p,  k \in \mathbb{N}_*$ (e.g., the set of strictly positive integers). A NLLS problem associated with a cost function $\phi(.)$ from $\mathbb{R}^{k}$ into  $\mathbb{R}$ and a residual real-vector function $\mathbf{r}(.)$ from $\mathbb{R}^{k}$ into  $\mathbb{R}^{m}$ is said to be separable if the parameter vector  $\mathbf{a} \in  \mathbb{R}^{k}$ can be partitioned as
\begin{equation*}
\mathbf{a} =  \begin{bmatrix} \mathbf{b} \\  \mathbf{c}  \end{bmatrix} \text{ with } \mathbf{b} \in  \mathbb{R}^{n}, \mathbf{c} \in  \mathbb{R}^{p} \text{ and } n + p = k  \  ,
\end{equation*}
in such a way that the subproblem
\begin{equation*}
\min_{\mathbf{c} \in  \mathbb{R}^{p}  }   \phi \big( \begin{bmatrix} \mathbf{b} \\  \mathbf{c}  \end{bmatrix} \big) = \frac{1}{2} \big\Vert \mathbf{r} (  \begin{bmatrix} \mathbf{b} \\  \mathbf{c}  \end{bmatrix} )   \big\Vert^{2}_{2}
\end{equation*}
is easy to solve numerically for every fixed $\mathbf{b} \in  \mathbb{R}^{n}$~\cite{GP1973}\cite{RW1980}\cite{HPS2012}.
 \\
\end{def2.10}

In the following, we will be particularly interested in the particular case when $\mathbf{r} (  \begin{bmatrix} \mathbf{b} \\  \mathbf{c}  \end{bmatrix} )$ is linear in $\mathbf{c} \in  \mathbb{R}^{p}$, i.e.,
\begin{equation*}
\mathbf{r} (  \begin{bmatrix} \mathbf{b} \\  \mathbf{c}  \end{bmatrix} ) = \mathbf{F}(  \mathbf{b} )\mathbf{c} - \mathbf{g} ( \mathbf{b} ) \text{ with } \mathbf{F}(  \mathbf{b} ) \in  \mathbb{R}^{q \times p} \text{ and } \mathbf{g} ( \mathbf{b} ) \in  \mathbb{R}^{q}  \  .
\end{equation*}
Let $\mathbf{c}(\mathbf{b})$ denotes one solution of the above subproblem for a given $\mathbf{b} \in  \mathbb{R}^{n}$ and formulate the problem
\begin{equation*}
\min_{\mathbf{b} \in  \mathbb{R}^{n}  }   \psi \big(  \mathbf{b} \big) = \phi \big( \begin{bmatrix} \mathbf{b} \\  \mathbf{c}(\mathbf{b})   \end{bmatrix} \big) = \frac{1}{2} \big\Vert \mathbf{r} (  \begin{bmatrix} \mathbf{b} \\  \mathbf{c}(\mathbf{b})  \end{bmatrix} )   \big\Vert^{2}_{2}  \ .
\end{equation*}
In doing that we have replaced our initial $k$-dimensional NLLS minimization problem by a $n$-dimensional one and we have separated the vector variables $\mathbf{b}$ and $\mathbf{c}$~\cite{RW1980}\cite{GP2003}. This definition is also valid for a cost function $\phi(.)$ from $\mathbb{R}^{p \times k}$ into  $\mathbb{R}$ and a residual real-matrix function $\mathbf{r}(.)$ from $\mathbb{R}^{p \times k}$ into  $\mathbb{R}^{n \times m}$.  Algorithms for minimizing  a separable real function $\psi(.)$ with a vector or matrix argument are called variable projection methods~\cite{GP1973}\cite{K1974}\cite{K1975}\cite{RW1980}\cite{B2009}\cite{OR2013}.

\section{Alternative and separable forms of the weighted low-rank approximation problem} \label{seppb:box}

In this section, we first provide some theoretical insights into the WLRA problem and the existence of solutions for it. Of course, some information on the subject is already available in the literature~\cite{MMH2003}\cite{CFP2003}\cite{SJ2004}\cite{GG2011}\cite{RSW2016}, but further investigations are clearly needed both theoretically and numerically, especially about the solvability of the WLRA problem. Moreover, the WLRA problem in its general form is much less well understood that the matrix completion or low-rank approximation problems~\cite{GG2011}\cite{RSW2016}. We also explain how the WLRA problem can be reformulated in several different, but related, ways such that variable projection algorithms for separable NLLS problems~\cite{GP1973}\cite{K1975}\cite{RW1980}\cite{B2009} can be used to solve it efficiently even when the number of missing entries in the input matrix is high. Finally, we highlight the closed links between variable projection methods and  Riemannian optimization on Grassmann manifolds~\cite{AMS2008}\cite{B2023}, which are two seemingly different approaches often used to solve the WLRA problem numerically. Despite the similarity of the two frameworks has already been highlighted in some studies (e.g.,~\cite{HF2015b}), the near equivalence of these two approaches (from a numerical point of view) in the context of the WLRA problem has not been well appreciated in the literature, probably because these two approaches have been developed in different communities~\cite{EAS1998}\cite{MMH2003}\cite{C2008b}\cite{BA2015}\cite{HF2015}\cite{HZF2017}.

\subsection{Nonconvex formulations of the WLRA problem} \label{noconv_wlra:box}

A reasonable and efficient way to tackle the low-rank constraint in the formulation~\eqref{eq:P0} of the WLRA problem is to introduce a bilinear factorization model of the low-rank matrix solution as $\mathbf{Y} = \mathbf{A}\mathbf{B} \text{ with }\mathbf{A}\in\mathbb{R}^{p \times k}\text{ and }\mathbf{B}\in\mathbb{R}^{k \times n}$~\cite{G1978}\cite{GZ1979}\cite{SJ2004}. This non-convex bilinear formulation has a very long history in statistics~\cite{W1966}\cite{WL1969}\cite{JHJ2009} and has been revitalized recently for solving similar semi-definite problems~\cite{BM2003}. This re-parametrization technique is justified by the fact that any matrix  $\mathbf{Y}$ of rank at most $k$ can be written as $\mathbf{Y} = \mathbf{A}\mathbf{B} \text{, with }\mathbf{A}\in\mathbb{R}^{p \times k}\text{ and }\mathbf{B}\in\mathbb{R}^{k \times n}$  and that, reciprocally, any such matrix product $\mathbf{A}\mathbf{B}$ is of rank at most $k$ (see Subsection~\ref{lin_alg:box} for details). Note that a similar multiplicative formulation holds for the (Eckart-Young) Theorem~\ref{theo2.1:box}, which solves the WLRA problem in the simple case where all elements of $\mathbf{W}$ are equal to one~\cite{G1978}. In recent decades, this bilinear factorization approach for low-rank matrix decomposition (often called the Burer-Monteiro factorization in the machine learning literature~\cite{BM2003}) has also been the subject of intense research (for efficiency reasons) in solving large-scale convex optimization problems as this (nonconvex) reformulation of the original convex problems allows to drastically reduce  the number of optimization variables from $p.n$ to $(p+n).k$, when $k$ is small (e.g., $k \ll \emph{min}(p,n)$), and, thus, allowing it to scale to problems with  thousands or even millions of variables~\cite{HMLZ2015}\cite{PKCS2017}\cite{LZT2019}. However, as we will illustrate below, this increased efficiency comes with a price as the intrinsic bilinearity of the multiplicative (Burer-Monteiro) formulation makes the landscape and geometry of the factored objective functions much more complicated than the original (convex) ones with additional first-order critical and solution points that are not global optima of the factored optimization problems, which can be also badly-conditioned matrices~\cite{LZT2019}.
 
We begin with the following well-known and simple result:
\\
\begin{theo3.1} \label{theo3.1:box}
For $\mathbf{X}\in\mathbb{R}^{p \times n}$, $\mathbf{W}\in\mathbb{R}^{p \times n}_+$  (i.e., $\mathbf{W}_{ij} \ge 0$),  $\sqrt{\mathbf{W}}\in\mathbb{R}^{p \times n}_+$  with $\sqrt{\mathbf{W}}_{ij} =  \sqrt{\mathbf{W}_{ij}}$ and any fixed integer $ k \le \emph{rank}( \mathbf{X} ) \le \text{min}( {p},{n} )$, the problem~\eqref{eq:P0} is equivalent to the problem (P1):
\begin{equation} \label{eq:P1} \tag{P1}
\min_{\mathbf{A}\in\mathbb{R}^{p \times k}\text{, }\mathbf{B}\in\mathbb{R}^{k \times n} }   \, \quad\  \varphi^{*}( \mathbf{A},\mathbf{B} ) = \frac{1}{2}   \Vert  \sqrt{\mathbf{W}}  \odot ( \mathbf{X} - \mathbf{A}\mathbf{B} )  \Vert^{2}_{F}  \ .
\end{equation}
In other words, if we consider the range of $\varphi(.)$
\begin{equation*}
\text{C}_{\varphi} = \big\{  y\in\mathbb{R}_+ \text{ }/\text{ }  \exists \mathbf{Y} \in\mathbb{R}^{p \times n}_{\le k}  \text{ with } y= \varphi( \mathbf{Y} ) \big\} \ ,
\end{equation*}
and the range of $\varphi^{*}(.)$
\begin{equation*}
\text{C}_{\varphi^*} =  \big\{  y\in\mathbb{R}_+ \text{ }/\text{ }  \exists (\mathbf{A},\mathbf{B}) \in\mathbb{R}^{p \times k} \times \mathbb{R}^{k \times n} \text{ with } y= \varphi^{*}( \mathbf{A},\mathbf{B}  ) \big\} \ ,
\end{equation*}
these two subsets of $\mathbb{R}$ have the same infimum (e.g., greatest lower bound) and if this infimum is a minimum for one subset, the other subset also admits a minimum and these two minima are equal.
\end{theo3.1}
\begin{proof}
Since elements of the ranges $\text{C}_{\varphi}$ and $\text{C}_{\varphi^*}$ are sums of squares, they are bounded below by zero and both $\text{C}_{\varphi}$ and $\text{C}_{\varphi^*}$ admit an infimum greater or equal to zero, say $\bar{\mathbf{c}}_{\varphi}$ and $\bar{\mathbf{c}}_{\varphi^*}$, respectively. Now, we will demonstrate the stronger result $\text{C}_{\varphi} = \text{C}_{\varphi^*}$ in which case the assertions in the theorem are obvious.

Suppose first that $ y \in \text{C}_{\varphi}$. Then, $\exists \mathbf{Y}\in\mathbb{R}^{p \times n}_{\le k}$ such that $y = \varphi( \mathbf{Y} )$. Now let
\begin{equation*}
\mathbf{Y} = \mathbf{U}\mathbf{\Sigma}\mathbf{V}^{T}
\end{equation*}
be the SVD of $\mathbf{Y}$, where it is assumed that $\mathbf{\Sigma}$ is a diagonal matrix with the singular values of $\mathbf{Y}$ arranged in decreasing order of magnitude in the diagonal. Since $\mathbf{Y}$ is of rank less than or equal to $k$, this SVD will have no more than $k$ singular triplets with a singular value distinct from zero. Thus,
\begin{equation*}
\mathbf{Y} = \mathbf{U}_{k} \mathbf{\Sigma}_{k} \mathbf{V}^{T}_{k}  \ ,
\end{equation*}
where $\mathbf{U}_{k}$ and $\mathbf{V}_{k}$ stand for submatrices formed by the first $k$ columns of $\mathbf{U}$ and $\mathbf{V}$, respectively, and $\mathbf{\Sigma}_{k}$  is the submatrix  defined by the first $k$ columns and rows of $\mathbf{\Sigma}$. Defining $\mathbf{A} = \mathbf{U}_{k}$ and $\mathbf{B} = \mathbf{\Sigma}_{k} \mathbf{V}^{T}_{k}$, $\mathbf{Y}$ can be factorized as
\begin{equation*}
\mathbf{Y} = \mathbf{A} \mathbf{B}  \quad\  \text{ with } \mathbf{A}\in\mathbb{R}^{p \times k}\text{ and } \mathbf{B}\in\mathbb{R}^{k \times n} \ .
\end{equation*}
However, the equation $\mathbf{Y} = \mathbf{A} \mathbf{B}$ implies that $y = \varphi( \mathbf{Y} ) = \varphi^{*}( \mathbf{A},\mathbf{B} )$ and, thus, $y  \in \text{C}_{\varphi^*}$.

Reciprocally, assume that $y \in \text{C}_{\varphi^*}$. Then,  it exists $(\mathbf{A},\mathbf{B}) \in\mathbb{R}^{p \times k} \times \mathbb{R}^{k \times n}$ such that $y = \varphi^{*}( \mathbf{A},\mathbf{B}  )$. If we define $\mathbf{Y} = \mathbf{A}\mathbf{B}$, we have $\emph{rank}( \mathbf{Y} )  \le  \text{min}\big( \emph{rank}( \mathbf{A} ), \emph{rank}( \mathbf{B} ) \big) \le k$ according to equation~\eqref{eq:rank2} and we conclude that $\mathbf{Y} \in \mathbb{R}^{p \times n}_{\le k}$. In these conditions, $y  = \varphi^{*}( \mathbf{A},\mathbf{B} ) = \varphi( \mathbf{Y} )$ and $y \in \text{C}_{\varphi}$ and we are done.
\\
\end{proof}

\begin{remark3.1} \label{remark3.1:box}
Since any  $p \times n$ matrix  $\mathbf{Y}$ of rank at most $k$ can also be written as $\mathbf{Y} = \mathbf{A}\mathbf{B}$ with
\begin{align*}
1) &  \quad\ \mathbf{A}\in\mathbb{R}^{p \times k}_k \text{, } \mathbf{B}\in\mathbb{R}^{k \times n}   \ , \\
2) &  \quad\ \mathbf{A}\in\mathbb{R}^{p \times k} \text{, } \mathbf{B}\in\mathbb{R}^{k \times n}_k  \ , \\
3) &  \quad\ \mathbf{A}\in\mathbb{O}^{p \times k} \text{, } \mathbf{B}\in\mathbb{R}^{k \times n}   \ , \\
4) &  \quad\ \mathbf{A}\in\mathbb{R}^{p \times k} \text{, } \mathbf{B}\in \mathbb{O}^{k \times n}_t \ ,
\end{align*}
and, reciprocally, any of these $\mathbf{A}\mathbf{B}$ matrix products is also of rank at most $k$ and the range of $\varphi^{*}(.)$ is also equal to
\begin{align*}
\text{C}_{\varphi^*}  & =  \big\{  \mathbf{y}\in\mathbb{R}_+ \text{ }/\text{ }  \exists (\mathbf{A},\mathbf{B})  \in\mathbb{R}^{p \times k}_k \times \mathbb{R}^{k \times n}  \text{ and } \mathbf{y}= \varphi^{*}( \mathbf{A},\mathbf{B}  ) \big\} \\
                                & =  \big\{  \mathbf{y}\in\mathbb{R}_+ \text{ }/\text{ }  \exists (\mathbf{A},\mathbf{B})  \in\mathbb{R}^{p \times k} \times \mathbb{R}^{k \times n}_k  \text{ and } \mathbf{y}= \varphi^{*}( \mathbf{A},\mathbf{B}  ) \big\} \\
                                & =  \big\{  \mathbf{y}\in\mathbb{R}_+ \text{ }/\text{ }  \exists (\mathbf{A},\mathbf{B})  \in\mathbb{O}^{p \times k}  \times \mathbb{R}^{k \times n}  \text{ and } \mathbf{y}= \varphi^{*}( \mathbf{A},\mathbf{B}  ) \big\} \\
                                & =  \big\{  \mathbf{y}\in\mathbb{R}_+ \text{ }/\text{ }  \exists (\mathbf{A},\mathbf{B})  \in\mathbb{R}^{p \times k} \times  \mathbb{O}^{k \times n}_t \text{ and } \mathbf{y}= \varphi^{*}( \mathbf{A},\mathbf{B}  ) \big\} \ .
\end{align*}
In these conditions, it is immediate that the problems~\eqref{eq:P0} and~\eqref{eq:P1} are also equivalent to the problems:
\begin{align*}
1) &  \quad\    \min_{\mathbf{A}\in\mathbb{R}^{p \times k}_k \text{, } \mathbf{B}\in\mathbb{R}^{k \times n}  }   \, \quad\  \varphi^{*}( \mathbf{A},\mathbf{B} ) \ , \\
2) &  \quad\  \min_{\mathbf{A}\in\mathbb{R}^{p \times k} \text{, } \mathbf{B}\in\mathbb{R}^{k \times n}_k  }   \, \quad\  \varphi^{*}( \mathbf{A},\mathbf{B} )   \ , \\
3) &  \quad\    \min_{\mathbf{A}\in\mathbb{O}^{p \times k}  \text{, } \mathbf{B}\in\mathbb{R}^{k \times n}  }   \, \quad\  \varphi^{*}( \mathbf{A},\mathbf{B} )    \ , \\
4) &  \quad\  \min_{\mathbf{A}\in\mathbb{R}^{p \times k} \text{, } \mathbf{B}\in\mathbb{O}^{k \times n}_t  }   \, \quad\  \varphi^{*}( \mathbf{A},\mathbf{B} )   \ ,
\end{align*}
where $\varphi^{*}( \mathbf{A},\mathbf{B} ) = \frac{1}{2}   \Vert  \sqrt{\mathbf{W}} \odot ( \mathbf{X} - \mathbf{A}\mathbf{B} )  \Vert^{2}_{F}$ and we will use these alternative forms to demonstrate some important properties of the WLRA problem in this section and the followings.  $\blacksquare$
\\
\end{remark3.1}

\begin{remark3.2} \label{remark3.2:box}
By using the rank-nullity relationship~\eqref{eq:rank} in Subsection~\ref{lin_alg:box}, another way to tackle the low-rank constraint in the WLRA problem is to impose this low-rank constraint on the dimensions of the null space of $\mathbf{Y}$ (or $\mathbf{Y}^{T}$)  instead on the range of  $\mathbf{Y}$ (or $\mathbf{Y}^{T}$) as in the formulation~\eqref{eq:P0}~\cite{EAS1998}\cite{MMH2003}\cite{MU2013}\cite{UM2014}. Since, from equation~\eqref{eq:null_ran}, we have
\begin{equation*}
\emph{null}( \mathbf{Y} ) = \emph{ran}( \mathbf{Y}^{T} )^\bot  \text{ and }  \emph{null}( \mathbf{Y}^{T} ) = \emph{ran}( \mathbf{Y} )^\bot \ .
\end{equation*}
This is equivalent to impose the low-rank constraint on the dimensions of the orthogonal complements of $\emph{ran}( \mathbf{Y} )$ or $\emph{ran}( \mathbf{Y}^{T} )$ and leads to what we will call the formulation~\eqref{eq:P2} of the WLRA problem, which has the following form if the low-rank constraint is imposed on the dimension of $\emph{ran}( \mathbf{Y} )^\bot$
\begin{equation} \label{eq:P2}  \tag{P2}
\min_{\mathbf{N}\in\mathbb{R}^{p \times (p-k)}_{p-k} \text{, }  \mathbf{Y}\in\mathbb{R}^{p \times n} \text{ with } \mathbf{N}^{T}\mathbf{Y} =  \mathbf{0}^{(p-k) \times n} } \, \quad\  \varphi^{**}( \mathbf{N}, \mathbf{Y} ) = \frac{1}{2}  \Vert \sqrt{\mathbf{W}} \odot ( \mathbf{X} - \mathbf{Y} )  \Vert^{2}_{F} \ ,
\end{equation}
or its transpose formulation~\eqref{eq:P2t}, if the low-rank constraint is imposed on the dimension of $\emph{null}( \mathbf{Y} )=\emph{ran}( \mathbf{Y}^{T} )^{\bot}$,
\begin{equation}  \label{eq:P2t}  \tag{P2t}
\min_{\mathbf{N}\in\mathbb{R}^{n \times (n-k)}_{n-k} \text{, }  \mathbf{Y}\in\mathbb{R}^{p \times n} \text{ with } \mathbf{Y}\mathbf{N} =  \mathbf{0}^{p \times (n-k)} } \, \quad\  \varphi^{**}( \mathbf{N}, \mathbf{Y} ) = \frac{1}{2}  \Vert \sqrt{\mathbf{W}} \odot ( \mathbf{X} - \mathbf{Y} )  \Vert^{2}_{F} \ .
\end{equation}
If $p<n$, the formulation~\eqref{eq:P2} should be preferred as the number of parameters to be estimated is reduced and vice versa if $p>n$. Here, the rank constraint is imposed by the equalities
\begin{equation*}
\mathbf{N}^{T}\mathbf{Y} =  \mathbf{0}^{(p-k) \times n} \text{ and } \mathbf{Y}\mathbf{N} =  \mathbf{0}^{p \times (n-k)} \ ,
\end{equation*}
which are, respectively, equivalent to
\begin{equation*}
\emph{dim} \big(  \emph{null}( \mathbf{Y}^{T} ) \big) \ge p-k \text{ and }  \emph{dim}\big(  \emph{null}( \mathbf{Y} ) \big) \ge n-k \ ,
\end{equation*}
since all the columns of $\mathbf{N}$ belong to the null space of $\mathbf{Y}^{T}$, or $\mathbf{Y}$ in the second case, and $\mathbf{N}$ is of full column rank in both cases. Obviously, since by the rank-nullity relationship~\eqref{eq:rank} we have
\begin{equation*}
\emph{dim}\big(  \emph{null}( \mathbf{Y}^{T} ) \big) + \emph{rank}\big( \mathbf{Y}^{T} \big) = p, \emph{dim}\big(  \emph{null}( \mathbf{Y} ) \big) + \emph{rank}\big( \mathbf{Y} \big) = n \text{ and } \emph{rank}( \mathbf{Y}^{T} ) = \emph{rank}( \mathbf{Y} ) \ ,
\end{equation*}
this is equivalent in both cases to the rank constraint $\emph{rank}( \mathbf{Y} ) \le k$, which is used in the formulation~\eqref{eq:P0} of the WLRA problem. Further inspection along the same lines of Theorem~\ref{theo3.1:box} will demonstrate that this formulation~\eqref{eq:P2} is also equivalent to the formulations~\eqref{eq:P0} and~\eqref{eq:P1}. When $\mathbf{W}\in\mathbb{R}^{p \times n}_{+*}$, Edelman et al.~\cite{EAS1998} and Manton et al.~\cite{MMH2003} have proposed a Grassmann manifold framework to solve problem~\eqref{eq:P2} as the solution of this problem depends only on the span of the columns of $\mathbf{N}$. A Grassmann manifold is the collection of all linear subspaces of a given dimension in a particular Euclidean or Frobenius space, see Subsection~\ref{calculus:box} and~\cite{B2023}  for a good introduction on manifolds and optimization on manifolds. Furthermore, they have described a large variety of first- and second-order algorithms for minimizing the cost function $\varphi^{**}(.)$ in this framework. As we will illustrate below, the solutions of the problem~\eqref{eq:P1} also do not depend on the individual elements of the matrices $\mathbf{A}$ and $\mathbf{B}$, but only on the range of $\mathbf{A}$ and, thus, can also be formulated as an optimization problem on the Grassmann manifold~\cite{DKM2012}\cite{BA2015}.

In these conditions, it is not difficult to recognize that each algorithm develops for minimizing $\varphi^{**}(.)$ (when $\mathbf{W}\in\mathbb{R}^{p \times n}_{+*}$) has a dual formulation for minimizing $\varphi^{*}(.)$ and vice versa, as determining the range of $\mathbf{A}$ leads implicitly to determine its orthogonal complement. In practical applications, the choice between an algorithm to minimize $\varphi^{*}(.)$ or its dual version to minimize $\varphi^{**}(.)$ will depend on the values of $k$, $p$ and $n$. For small values of $k$, the formulation~\eqref{eq:P1} is likely to be more efficient as the size of the matrix variables  will be smaller and, conversely, the formulation~\eqref{eq:P2} can be a better choice for large values of $k$ as we will deal with smaller matrix variables when minimizing $\varphi^{**}(.)$. We will come back to these alternatives in the next sections. Finally, we mention that it is probably possible to extend the algorithms proposed by Manton et al.~\cite{MMH2003} to minimize the cost function $\varphi^{**}(.)$ to the case where $\mathbf{W}\in\mathbb{R}^{p \times n}_{+}$ instead of $\mathbb{R}^{p \times n}_{+*}$, see~\cite{C2008b} for work in this direction. But, this is not pursued here, as in most applications, we use values of $k$ which are much more smaller than min$(p,n)$ for which the formulation~\eqref{eq:P1} is likely more economical. $\blacksquare$
\\
\end{remark3.2}

\begin{remark3.3} \label{remark3.3:box}
A popular way to tackle the WLRA problem is also to consider the simpler problems:
\begin{equation}
\min_{\mathbf{Y} \in \mathbb{R}^{p \times n}_{k}  } \, \quad\  \varphi( \mathbf{Y} ) = \frac{1}{2}  \Vert \sqrt{\mathbf{W}} \odot ( \mathbf{X} - \mathbf{Y} )  \Vert^{2}_{F} \ ,
\end{equation}
or
\begin{equation}
\min_{\mathbf{A} \in \mathbb{R}^{p \times k}_{k} \text{, }\mathbf{B}\in\mathbb{R}^{k \times n}_{k} }   \, \quad\  \varphi^{*}( \mathbf{A},\mathbf{B} ) = \frac{1}{2}   \Vert  \sqrt{\mathbf{W}} \odot ( \mathbf{X} - \mathbf{A}\mathbf{B} )  \Vert^{2}_{F} \ ,
\end{equation}
which are equivalent (as $\emph{rank}( \mathbf{A}\mathbf{B} ) = k$ if  $\emph{rank}( \mathbf{A} ) = \emph{rank}( \mathbf{B} ) = k$, see Subsection~\ref{lin_alg:box}) and are also frequently solved by Riemannian optimization methods applied to smooth fixed-rank matrix manifolds~\cite{V2013}\cite{MMBS2012}\cite{MMBS2014}  as the cost functions $\varphi(.)$ and  $\varphi^{*}(.)$ are infinitely differentiable (e.g., of class $C^{\infty}$ ) and the set $\mathbb{R}^{p \times n}_{k}$ is a smooth ($C^{\infty}$) embedded submanifold of $\mathbb{R}^{p \times n}$ of dimension $(p+n-k)k$ (see Proposition 1.14 in Chap.~5 of~\cite{HM1996}, Example 8.14 of~\cite{L2003} or Section 7.5 in Chap.~7 of~\cite{B2023}). This approach is justified by the fact that $\mathbb{R}^{p \times n}_{k}$ is dense and open in $\mathbb{R}^{p \times n}_{\le k}$ (see Theorem~\ref{theo2.3:box}) meaning that with an initial guess in $\mathbb{R}^{p \times n}_{k}$, an iterate belonging to $\mathbb{R}^{p \times n}_{< k}$ or a non-smooth point of $\varphi(.)$ are both unlikely to occur in practice.

However, these two simpler problems are not mathematically equivalent to~\eqref{eq:P0} and~\eqref{eq:P1} for any choice of the weight matrix $\mathbf{W}$ as the submanifold $\mathbb{R}^{p \times n}_k$ is not closed in $\mathbb{R}^{p \times n}$  and a solution of these simpler problems may be on the frontier of $\mathbb{R}^{p \times n}_k$, which is $\mathbb{R}^{p \times n}_{<k}$, as stated in Theorem~\ref{theo2.3:box}. This implies that these simpler problems may not admit a global minimizer, while such global minimizer will exist for problems~\eqref{eq:P0} and~\eqref{eq:P1}~\cite{CFP2003}. Furthermore, closedness of the domain is important in (non-convex) nonlinear optimization to garantee that the limit point of the iterative sequence is still in the domain of interest. As the set $\mathbb{R}^{p \times n}_{k}$ is not closed, some matrices in $\mathbb{R}^{p \times n}_{<k}$ can be the limit points of the iterative sequences in $\mathbb{R}^{p \times n}_k$ leading to so-called spurious critical points which do not belong to the smooth fixed-rank manifold $\mathbb{R}^{p \times n}_{k}$~\cite{LKB2023}. Similarly, a sequence might also cross the frontier of  $\mathbb{R}^{p \times n}_{k}$ at a certain iterate and the rank might fall below $k$ breaking the sequence. For all these reasons, it is better to solve the WLRA problem over $\mathbb{R}^{p \times n}_{ \le k}$ rather than over $\mathbb{R}^{p \times n}_{k}$. Note, on the other hand, that optimization algorithms on smooth fixed-rank manifolds are not strictly applicable on $\mathbb{R}^{p \times n}_{\le k}$ as this set is a (non-smooth) real algebraic variety, not an embedded smooth submanifold of $\mathbb{R}^{p \times n}$ (see Proposition 1.1 in~\cite{{BV1988}}, Lecture 9 of~\cite{H1992} or~\cite{B2023}\cite{SU2015} for details). More precisely, $\mathbf{Y} \in \mathbb{R}^{p \times n}_{\le k}$ is, by definition, a space of matrices with a given upper bound on their ranks and is naturally an algebraic variety as the rank condition on a matrix is equivalent to the vanishing of all of its $(k+1, k+1)$-minors, which are  polynomials of degree $k+1$. $\mathbb{R}^{p \times n}_{\le k}$ is then defined as the solution set of polynomial equations therefore a so-called real determinantal variety~\cite{H1992}. Extending (first-order) optimization algorithms developed for smooth fixed-rank manifolds to real  determinantal varieties like $\mathbb{R}^{p \times n}_{\le k}$ is a very active area of research recently~\cite{SU2015}\cite{LKB2023}\cite{OGA2024}\cite{OA2024}, but variable projection techniques, which are the focus of this monograph, can also be used for that purpose. $\blacksquare$
\\
\end{remark3.3}

Thus, it is equivalent to minimize $\varphi( \mathbf{Y} )$ or $\varphi^{*}( \mathbf{A}, \mathbf{B} )$ for solving the WLRA problem. However, the WLRA problem (e.g., in the formulations~\eqref{eq:P0} and~\eqref{eq:P1}) has no known closed form solution in the general case and is known to be NP-hard~\cite{GG2011} as already discussed in the Introduction~\ref{intro:box}. For certain classes of weighting matrices, a globally optimal solution can be found and one such class is obviously the unweighted case (e.g., $\mathbf{W}_{i,j}=1$), since in that case the solution of the WLRA problem is given by the Eckart-Young Theorem~\ref{theo2.1:box}. Another very important specialization of this is the case where all the elements of $\mathbf{W}$ are greater than 0 in which case it is possible to demonstrate that the WLRA problem has a well-defined solution as demonstrated in Theorem~\ref{theo3.3:box} below. Moreover, in the case where all the elements of $\mathbf{W}$ are greater than 0 and the rank of  $\mathbf{W}$ is equal to $1$, the solution of the WLRA problem can also be found via a generalization of the SVD in which we use diagonal metrics and scalar products different from the identity matrix in $\mathbb{R}^{p}$ and $\mathbb{R}^{n}$ (see Theorem 3 of~\cite{MMH2003} and also~\cite{GG2011}\cite{RSW2016}). Finally, if $k = \emph{rank}( \mathbf{X} )$, the WLRA problem is equivalent to the $consistent$ matrix completion problem, which is to find one matrix $\widehat{\mathbf{X}}$ of rank at most $k$ consistent with the observed entries (e.g., $\mathbf{W}_{ij} \ne 0$) of $\mathbf{X}\in\mathbb{R}^{p \times n}_k$  (e.g., the problem of recovering large matrices of low rank when most of the entries are unknown). In this case,  the problem is  also well-posed since $\mathbf{X}$ is obviously a solution to the consistent  completion problem and we have $\varphi( \widehat{\mathbf{X}}  ) = 0$ for all solution matrices $\widehat{\mathbf{X}}$~\cite{DMK2011}\cite{DKM2012}. In the general case, a very large variety of iterative methods have been previously suggested to solve the WLRA problem or convex and smooth proxies of it, especially in the framework of low-rank matrix completion, which is also  NP-hard~\cite{CR2009}, and is the focus of lot of recent research~\cite{RS2005}\cite{SRJ2005}\cite{CR2009}\cite{CCS2010}\cite{WWY2012}\cite{KM2010}\cite{HMLZ2015}\cite{NKS2019}. Both the WLRA and matrix completion problems are also frequently recast as an optimization problem on smooth matrix manifolds as already noted above~\cite{AMS2008}\cite{MMH2003}\cite{SE2010}\cite{DKM2012}\cite{BA2015}\cite{B2023}.

The formulation~\eqref{eq:P0} of the WLRA problem is well suited to derived theoretical properties of the WLRA problem such as the existence of solutions for this problem. On the other hand, the interest of the alternative formulation~\eqref{eq:P1} and its variants (see Remark~\ref{remark3.1:box}), is that smaller matrices are manipulated and  the introduction of the (non-unique) parameterization $\mathbf{Y}=\mathbf{A}\mathbf{B}$ allows us to recast the WLRA problem as a standard unconstrained NLLS minimization problem as we will show below. This is particularly useful to derive practical algorithms to solve the WLRA problem as we will illustrate in the next sections.
\\
\begin{remark3.4} \label{remark3.4:box}
 The problem~\eqref{eq:P1} or its variants is over-parameterized. More precisely, if  $\mathbf{C}$ is a $k \times k$ invertible matrix, we have
\begin{equation*}
\mathbf{A} \mathbf{B} = \mathbf{A}( \mathbf{C}\mathbf{C}^{-1}  ) \mathbf{B} = ( \mathbf{A}\mathbf{C} )( \mathbf{C}^{-1}\mathbf{B} )  \quad\ \text{and} \quad\ \varphi^{*}( \mathbf{A}, \mathbf{B} ) = \varphi^{*}( \mathbf{A}\mathbf{C}, \mathbf{C}^{-1}\mathbf{B} )  \ .
\end{equation*}
Consequently, the set of global minimizers of $\varphi^{*}(.)$ can be empty or infinite, but never finite or an isolated minimum implying that the Hessian of $\varphi^{*}(.)$ is at best positive semi-definite, but never positive definite, see Subsection\ref{calculus:box} for details. This can severely degrade the performance of standard optimization algorithms, which are mostly developed for isolated optima ~\cite{DS1983}\cite{NW2006}. Furthermore, this scaling ambiguity tends to make the cost function $\varphi^{*}(.)$ of  problem~\eqref{eq:P1} badly-conditioned, especially when the matrix $\mathbf{C}$ or its inverse is nearly singular. To overcome this difficulty, many authors have proposed to add different regularizers to $\varphi^{*}(.)$ as we will discussed later in this section.

Notice also that, if $\mathbf{A} \in  \mathbb{R}^{p \times k}_k$ and  $\mathbf{B}  \in  \mathbb{R}^{k \times n}_k$, these two full rank matrices have $p.k$ and $k.n$ degrees of freedom, respectively. However, specifying the matrix product $\mathbf{A} \mathbf{B}$ in $\varphi^{*}(.)$   is equivalent to use the  matrix product $( \mathbf{A}\mathbf{C} )( \mathbf{C}^{-1}\mathbf{B} )$ for any $k  \times k$ matrix $\mathbf{C}$ of rank $k$, which is equivalent to specify the column space of $\mathbf{A} \mathbf{B}$. Hence, the matrix product $\mathbf{A} \mathbf{B}$, or its column space, has only $p.k + k.n - k.k = (p+n-k).k$ degrees of freedom in general, which is consistent with the fact that the set $\mathbb{R}^{p \times n}_k$ is a smooth submanifold of $\mathbb{R}^{p \times n}$ of dimension $ (n+p-k).k$ as already noted in Remark~\ref{remark3.3:box}  above.

More generally, as all the matrix products $( \mathbf{A}\mathbf{C} )( \mathbf{C}^{-1}\mathbf{B} )$  share the same column space, possibly remedies for the implicit over-parameterization in the formulation~\eqref{eq:P1} can be to recast the WLRA problem as an optimization problem on a Grassmann manifold~\cite{DKM2012}\cite{C2008b}\cite{BA2015}\cite{MMH2003}\cite{MMBS2012}\cite{MMBS2014} as discussed in Remark~\ref {remark3.3:box} or to use variable projection methods~\cite{R1974}\cite{C2008a}\cite{OYD2011}\cite{OD2007}. Moreover, these two seemingly different approaches for solving the WLRA problem are in fact tightly related as we will illustrate below. $\blacksquare$
\\
\end{remark3.4}

The cost functions $\varphi(.)$  and $\varphi^{*}(.)$ are the composition of several infinitely differentiable functions on their respective domain of definition and, consequently, are also infinitely differentiable as smoothness is preserved by composition thanks to the standard chain rule~\cite{C2017}. Since  $\varphi(.)$  and $\varphi^{*}(.)$ are smooth, they are also continuous on their respective domains. However, in the next theorem, we give a direct demonstration of the continuity of $\varphi(.)$  and $\varphi^{*}(.)$ by making clear that the WLRA problem differs from the standard low-rank approximation problem only by the choice of a different metric than the standard Frobenius metric on $\mathbb{R}^{p \times n}$. This metric is derived from the norm or seminorm induced by the choice of the weight matrix $\mathbf{W} \in \mathbb{R}^{p \times n}_{+}$.
\begin{theo3.2} \label{theo3.2:box}
Using the same notations and definitions as in Theorem~\ref{theo3.1:box}, the objective function defined in problem~\eqref{eq:P0}
\begin{equation*}
\varphi : \mathbb{R}^{p \times n}  \longrightarrow \mathbb{R} : \mathbf{Y} \mapsto \varphi( \mathbf{Y} ) =  \frac{1}{2}   \Vert  \sqrt{\mathbf{W}}  \odot ( \mathbf{X} - \mathbf{Y} )  \Vert^{2}_{F}  \ ,
\end{equation*}
and the objective function defined in problem~\eqref{eq:P1}
\begin{equation*}
\varphi^{*} : \mathbb{R}^{p \times k} \times  \mathbb{R}^{k \times n} \longrightarrow \mathbb{R} : ( \mathbf{A}, \mathbf{B} ) \mapsto \varphi^{*}( \mathbf{A}, \mathbf{B} ) =  \frac{1}{2}   \Vert \sqrt{\mathbf{W}}  \odot ( \mathbf{X} - \mathbf{A}\mathbf{B} )  \Vert^{2}_{F}  \ ,
\end{equation*}
are continuous on their respective domains of definition.
\end{theo3.2}
\begin{proof}
We first define a weighted norm or seminorm (if some of elements of $\mathbf{W}$ are equal to zero) of an $p \times n$ real matrix $\mathbf{Y}$ as
\begin{equation*}
 \Vert \mathbf{Y} \Vert^{2}_{\mathbf{W}} = \emph{vec}( \mathbf{Y} )^{T} \emph{diag}\big( \emph{vec}(\mathbf{W}) \big) \emph{vec}( \mathbf{Y} )  \ ,
\end{equation*}
where $\emph{vec}( \mathbf{Y} )$ stands for the vectorized form of $\mathbf{Y}$, i.e., a vector formed by stacking the consecutive columns of $\mathbf{Y}$ in one $p.n$-dimensional vector (see equation~\eqref{eq:vec} in Subsection~\ref{multlin_alg:box}). If none of the elements of $\mathbf{W}$ is equal to zero, $\Vert  \Vert_{\mathbf{W}}$ is obviously a norm on $\mathbb{R}^{p \times n}$ and, as $\mathbb{R}^{p \times n}$ is a finite-dimensional vector space over $\mathbb{R}$, all norms on $\mathbb{R}^{p \times n}$ are equivalent, induce the same topology and are continuous functions on $\mathbb{R}^{p \times n}$ with respect to this topology~\cite{C2017}\cite{B1993}.
On the other hand, if some of the elements of $\mathbf{W}$ are equal to zero, $\Vert  \Vert_{\mathbf{W}}$ is only a seminorm on $\mathbb{R}^{p \times n}$, e.g.,
$\Vert  \Vert_{\mathbf{W}}$ is a real-valued function : $\mathbb{R}^{p \times n} \longrightarrow \mathbb{R}$, which verifies, for all $\mathbf{Y}, \mathbf{Z} \in \mathbb{R}^{p \times n}$ and $\alpha \in \mathbb{R}$,
\begin{align*}
& \Vert  \mathbf{Y} \Vert_{\mathbf{W}} \ge 0    \ , \\
& \Vert  \alpha \mathbf{Y} \Vert_{\mathbf{W}} = \vert \alpha \vert \Vert  \mathbf{Y} \Vert_{\mathbf{W}}  \ , \\
& \Vert  \mathbf{Y} + \mathbf{Z} \Vert_{\mathbf{W}} \le \Vert  \mathbf{Y} \Vert_{\mathbf{W}} + \Vert  \mathbf{Z} \Vert_{\mathbf{W}}  \ .
\end{align*}
However, even if $\Vert  \Vert_{\mathbf{W}}$ is only a seminorm, it is still continuous with the respect to the unique topology on $\mathbb{R}^{p \times n}$ as demonstrated by Goldberg~\cite{G2017}.

Now, $\varphi( \mathbf{Y} )$ may be expressed as
\begin{align*}
 \varphi( \mathbf{Y} ) & =  \frac{1}{2}   \Vert  \sqrt{\mathbf{W}}  \odot ( \mathbf{X} - \mathbf{Y} )  \Vert^{2}_{F} \\
  & =  \frac{1}{2}   \Vert  \emph{vec} \big( \sqrt{\mathbf{W}}  \odot ( \mathbf{X} - \mathbf{Y} ) \big)  \Vert^{2}_{2} \\
  & =  \frac{1}{2}   \Vert   \emph{vec}( \sqrt{\mathbf{W}} ) \odot \emph{vec}( \mathbf{X} - \mathbf{Y} )  \Vert^{2}_{2} \\
  & =  \frac{1}{2}   \emph{vec}( \mathbf{X} - \mathbf{Y} )^{T} \emph{diag}\big( \emph{vec}(\mathbf{W}) \big) \emph{vec}( \mathbf{X} - \mathbf{Y} )  \\
  & =  \frac{1}{2} \Vert \mathbf{X} - \mathbf{Y} \Vert^{2}_{\mathbf{W}} \ .
\end{align*}
In other words, $\varphi(.)$ is the composition of the residual matrix function: $\mathbf{Y} \mapsto  \mathbf{X} - \mathbf{Y}$, the norm or seminorm: $\mathbf{Z}  \mapsto \Vert \mathbf{Z}  \Vert_{\mathbf{W}}$ and the square function: $y \mapsto y^2 $. As all these functions are continuous on their respective domain of definition, we conclude that  $\varphi(.)$ is also continuous on $\mathbb{R}^{p \times n}$.

Similarly, $\varphi^{*}( \mathbf{A}, \mathbf{B} )$ may be expressed as
\begin{align*}
 \varphi^{*}( \mathbf{A}, \mathbf{B} ) & =  \frac{1}{2}   \Vert  \sqrt{\mathbf{W}}  \odot ( \mathbf{X} - \mathbf{A}\mathbf{B} )  \Vert^{2}_{F} \\
    & =  \frac{1}{2}   \emph{vec}( \mathbf{X} - \mathbf{A}\mathbf{B} )^{T} \emph{diag}\big( \emph{vec}(\mathbf{W}) \big) \emph{vec}( \mathbf{X} - \mathbf{A}\mathbf{B} )  \\
    & =  \frac{1}{2} \Vert \mathbf{X} - \mathbf{A}\mathbf{B} \Vert^{2}_{\mathbf{W}}
\end{align*}
and $\varphi^{*}(.)$ is also the composition of several continuous functions on their respective domain of definition and, consequently,  $\varphi^{*}(.)$ is also continuous on $\mathbb{R}^{p \times k} \times  \mathbb{R}^{k \times n}$.
\\
\end{proof}
 As  $\varphi(.)$ is continuous on its domain of definition, it is not difficult to show that the problem~\eqref{eq:P0} has a well-defined solution when all the elements of the weight matrix $\mathbf{W}$ are strictly positive as stated in the next theorem.
 \\
 \begin{theo3.3} \label{theo3.3:box}
 For $\mathbf{X}\in\mathbb{R}^{p \times n}$ different of the zero matrix of $\mathbb{R}^{p \times n}$ and $\mathbf{W}\in\mathbb{R}^{p \times n}_{+*}$  (i.e., $\mathbf{W}_{ij} > 0$), and any fixed integer $ k \le \emph{rank}( \mathbf{X} ) \le \text{min}( {p},{n} )$, the set of global minimizers of $\varphi( \mathbf{Y} )$ on $\mathbb{R}^{p \times n}_{\le k}$ is nonempty and compact.
\end{theo3.3}
\begin{proof}
This theorem is a direct consequence of Theorem 3.1 stated without proof in Chu et al.~\cite{CFP2003}, but we give a direct proof for completeness.

As  $\mathbf{W}\in\mathbb{R}^{p \times n}_{+*}$ by hypothesis, we first observe that $\Vert  \Vert_{\mathbf{W}}$ defines a norm on  $\mathbb{R}^{p \times n}$.
Let us now consider the closed ball with center $\mathbf{X}$ and radius $r=\Vert  \mathbf{X} \Vert_{\mathbf{W}}$ with respect to this norm in $\mathbb{R}^{p \times n}$:
\begin{equation*}
\bar{B}_{p \times n}(\mathbf{X}, r)  = \big\{   \mathbf{Y} \in\mathbb{R}^{p \times n}  \text{ and }  \Vert  \mathbf{X} - \mathbf{Y} \Vert_{\mathbf{W}} \le r  \big\} \ .
\end{equation*}
$\bar{B}_{p \times n}(\mathbf{X}, r)$ is not empty as the zero matrix of $\mathbb{R}^{p \times n}$, which is also an element of $\mathbb{R}^{p \times n}_{\le k}$, is in this closed ball. As $\mathbb{R}^{p \times n}$ is a finite-dimensional vector space, this closed ball is also a compact set (as it is by definition a bounded set). Furthermore, as $\mathbb{R}^{p \times n}_{\le k}$ is closed in $\mathbb{R}^{p \times n}$ (see Theorem~\ref{theo2.3:box}), the intersection of $\mathbb{R}^{p \times n}_{\le k}$ and $\bar{B}_{p \times n}(\mathbf{X}, r)$ is also closed and bounded and, thus, compact in $\mathbb{R}^{p \times n}$.
Now, as $\varphi(.)$ is continuous on $\mathbb{R}^{p \times n}$ and the image of a compact set by a continuous function is also compact, we conclude that $\varphi \big(  \mathbb{R}^{p \times n}_{\le k}  \cap \bar{B}_{p \times n}(\mathbf{X}, r)  \big) \subset  \text{C}_{\varphi}$ is a compact set in $\mathbb{R}$ and, thus, a closed and bounded interval of $\mathbb{R}$.
Thus, $\varphi(.)$ attains its infimum on $\mathbb{R}^{p \times n}_{\le k}  \cap \bar{B}_{p \times n}(\mathbf{X}, r)$. In other words, it exists $\mathbf{\widehat{Y}} \in \mathbb{R}^{p \times n}_{\le k}  \cap \bar{B}_{p \times n}(\mathbf{X}, r)$ such that
\begin{equation*}
\varphi(\mathbf{\widehat{Y}}) \le \varphi(\mathbf{Y}) \text{ , } \forall \mathbf{Y} \in \mathbb{R}^{p \times n}_{\le k}  \cap \bar{B}_{p \times n}(\mathbf{X}, r) \ .
\end{equation*}
It remains to show that $\varphi(\mathbf{\widehat{Y}})=\bar{\mathbf{c}}_{\varphi}$ where $\bar{\mathbf{c}}_{\varphi}$ is the infimum of $\varphi(.)$ on $\mathbb{R}^{p \times n}_{\le k}$, i.e., that $\mathbf{\widehat{Y}}$ is also a global minimizer of $\varphi(.)$ on $\mathbb{R}^{p \times n}_{\le k}$.
By definition of $\bar{\mathbf{c}}_{\varphi}$, we already have $\bar{\mathbf{c}}_{\varphi} \le \varphi(\mathbf{\widehat{Y}})$ and it is sufficient to show that $\varphi(\mathbf{\widehat{Y}})  \le \bar{\mathbf{c}}_{\varphi}$ to demonstrate the theorem.
 \\
Suppose on the contrary that $\varphi(\mathbf{\widehat{Y}})  > \bar{\mathbf{c}}_{\varphi}$, then it exists $\mathbf{Y} \in \mathbb{R}^{p \times n}_{\le k}$ such that $\varphi(\mathbf{\widehat{Y}})  > \varphi(\mathbf{Y}) \ge\bar{\mathbf{c}}_{\varphi}$ by definition of $\bar{\mathbf{c}}_{\varphi}$. However, this implies that
\begin{equation*}
\frac{1}{2} \Vert \sqrt{\mathbf{W}} \odot ( \mathbf{X} - \mathbf{Y} )  \Vert^{2}_{F} = \varphi(\mathbf{Y}) < \varphi(\mathbf{\widehat{Y}}) = \frac{1}{2} \Vert \sqrt{\mathbf{W}} \odot ( \mathbf{X} - \mathbf{\widehat{Y}} ) \Vert^{2}_{F} \ ,
\end{equation*}
and it follows that 
\begin{equation*}
\Vert  \mathbf{X} - \mathbf{Y} \Vert_{\mathbf{W}}  < \Vert  \mathbf{X} - \mathbf{\widehat{Y}} \Vert_{\mathbf{W}}  \le \Vert  \mathbf{X} \Vert_{\mathbf{W}} = r \ .
\end{equation*}
In other words,  $\mathbf{Y} \in \mathbb{R}^{p \times n}_{\le k} \cap \bar{B}_{p \times n}(\mathbf{X}, r)$ and $\varphi(\mathbf{Y}) < \varphi(\mathbf{\widehat{Y}})$, which contradicts the assertion that $\mathbf{\widehat{Y}}$ is a minimizer of $\varphi$ on $\mathbb{R}^{p \times n}_{\le k} \cap \bar{B}_{p \times n}(\mathbf{X}, r)$ and we are done.
\\
\end{proof}

\begin{remark3.5} \label{remark3.5:box}
Using the equivalence between problems~\eqref{eq:P0} and~\eqref{eq:P1} stated in Theorem~\ref{theo3.1:box} above, we conclude that the set of global minimizers of  $\varphi^{*}(.)$, when the weight matrix is strictly positive, is also nonempty.  However, in the formulation~\eqref{eq:P1} of the WLRA problem, an important point to keep in mind is that, if the solution set is not empty, problem~\eqref{eq:P1} has an infinity of solutions as $\widehat{\mathbf{A}}$ and $\widehat{\mathbf{B}}$ are not determined uniquely and we can normalize them in an arbitrary manner without changing the value of $\varphi^{*}( \widehat{\mathbf{A}},\widehat{\mathbf{B}}  )$ (see Remark~\ref{remark3.4:box} above). Moreover, if $\alpha \in \mathbb{R}_*$, $( \alpha.\widehat{\mathbf{A}},  \frac{1}{\alpha}.\widehat{\mathbf{B}} )$ is also a solution of~\eqref{eq:P1}, which shows that the set of solutions in $\mathbb{R}^{p \times k} \times \mathbb{R}^{k \times n}$ is unbounded and, thus, not compact despite the set of global minimizers of $\varphi(.)$ is compact in $\mathbb{R}^{p \times n}$. $\blacksquare$
\\
\end{remark3.5}

\begin{remark3.6} \label{remark3.6:box}
In the unweighted case (and with no missing values), the WLRA problem has an unique global minimum and all critical points of  $\varphi(.)$ or $\varphi^{*}(.)$ which are not global minimizers are saddle points (e.g., critical points whose every neighborhood contains both "higher" and "smaller" points for $\varphi(.)$ or $\varphi^{*}(.)$), see Section 2.1 of~\cite{SJ2004} and Theorem 1.14 of~\cite{H2010} for details. In other words, $\varphi(.)$ or $\varphi^{*}(.)$ do not admit local minima in the unweighted case despite they are not convex functions.
 \\
While Theorem~\ref{theo3.3:box} shows that the WLRA problem has still well defined solutions when $\mathbf{W}\in\mathbb{R}^{p \times n}_{+*}$ because $\Vert  \Vert_{\mathbf{W}}$ is a norm, several authors have illustrated by examples that  $\varphi(.)$ or $\varphi^{*}(.)$ can have multiple local minima in addition to saddle points when the weights are all different of zero, but not uniform (see Section 2.1 of~\cite{SJ2004} and Example 1 of~\cite{GG2011}). Such local minima emerge especially when the weights become significantly non-uniform (see Figure 1 of~\cite{SJ2004} for illustration). When $\mathbf{W}$ has zero entries, the situation is even worse as $\varphi(.)$ or $\varphi^{*}(.)$ may have multiple local minima~\cite{IR2010}, but the infimum of $\varphi(.)$ or $\varphi^{*}(.)$ can also be unattained, see Example 2 of~\cite{GG2011} for illustration. $\blacksquare$
\\
\end{remark3.6}

An alternative and insightful demonstration of the above theorem can also be given using the notion of the level sets of a continuous real function as defined in Chapter 4 of Ortega and Rheinboldt~\cite{OR1970}. More precisely, for $\gamma \in \mathbb{R}$, the level set of $\varphi(.)$ at level $\gamma$ is the set $L(\gamma) =  \big\{  \mathbf{Y} \in \mathbb{R}^{p \times n}_{\le k}  \text{ / } \varphi(\mathbf{Y}) \le  \gamma \big\}$. In other words, $L(\gamma)$ is the subset of  $\mathbb{R}^{p \times n}_{\le k}$  whose elements $\mathbf{Y}$ verify the inequality  $\varphi(\mathbf{Y}) \le  \gamma$. Obviously, $L(\gamma)$ is empty if $\gamma <\bar{\mathbf{c}}_{\varphi}$, where $\bar{\mathbf{c}}_{\varphi}$ is the infimum of $\varphi(.)$, and is the set of the global minimizers of $\varphi(.)$ if $\gamma = \bar{\mathbf{c}}_{\varphi}$ (which can be also empty in the general case where $\mathbf{W}\in\mathbb{R}^{p \times n}_+$ as discussed above).
 \\
As $\varphi(.)$ is continuous on the closed set $\mathbb{R}^{p \times n}_{\le k}$ and the range of $\varphi(.)$, $\text{C}_{\varphi}$, is included in the nonnegative half-space of $\mathbb{R}$, then every level set of $\varphi(.)$ at level $\gamma$ for $\gamma \ge \bar{\mathbf{c}}_{\varphi}$ is closed  in $\mathbb{R}^{p \times n}$ as the reciprocal image of the closed interval $\lbrack \bar{\mathbf{c}}_{\varphi},  \gamma  \rbrack$ by a continuous and real function is also closed. Under these conditions, a necessary and sufficient condition for the set of global minimizers of $\varphi(.)$ to be nonempty and compact is that $\varphi(.)$ has a nonempty and bounded level set $L(\gamma)$ as this implies that $L(\gamma)$ is compact in $\mathbb{R}^{p \times n}$ (see Propositions 4.2.2 and 4.3.1 in Chap.~4 of~\cite{OR1970}). However, since $L(\gamma)$ is simply the intersection of $\mathbb{R}^{p \times n}_{\le k}$ and the closed ball with center $\mathbf{X}$ and radius $\sqrt{ 2 . \gamma }$ (with respect to the norm $\Vert  \Vert_{\mathbf{W}}$) if all the elements of $\mathbf{W}$ are strictly positive, $L(\gamma)$ is nonempty and bounded by definition for all $\gamma > \bar{\mathbf{c}}_{\varphi}$. This also proves that the set of global minimizers of $\varphi(.)$ is nonempty and compact if all the elements of $\mathbf{W}$ are strictly positive as stated in Theorem~\ref{theo3.3:box}.
 \\
 \\
In the more difficult case, where some elements of $\mathbf{W}$ are equal to zero, $\varphi(.)$ is still continuous as $\Vert  \Vert_{\mathbf{W}}$ defines a seminorm on  $\mathbb{R}^{p \times n}$ and every level set of $\varphi(.)$ is also automatically closed and the question of the existence of a global minimizer of $\varphi(.)$ reduces again to the existence of a bounded level set $L(\gamma)$ according to the previous discussion. However, in the case where some of the elements of $\mathbf{W}$ are equal to zero, the seminorm $\Vert  \Vert_{\mathbf{W}}$ does not define the topology of $\mathbb{R}^{p \times n}$~\cite{G2017} and the level set $L(\gamma)$ is not automatically bounded, so that the question of the existence of a nonempty and compact set of global minimizers is still unanswered in that case.
\\
\\
In order to discuss in more details, the existence of a nonempty and compact set of global minimizers of $\varphi(.)$ when some elements of $\mathbf{W}$ are equal to zero, let $\Omega \subset \lbrack p \rbrack \times  \lbrack n \rbrack$ be the set of indices of the elements of $\mathbf{W}$ such that $\mathbf{W}_{ij} \ne 0$, where $\lbrack L \rbrack = \lbrack 1, 2, ..., L \rbrack$. With this definition, from a weight matrix $\mathbf{W}$ with some zero elements and any $\lambda \in \mathbb{R}_{+*}$ (e.g., $\lambda > 0$), we can define a new $p \times n$ weight matrix $\mathbf{W}_\lambda$ as follows
 \begin{equation*}
    \big \lbrack \mathbf{W}_\lambda \big \rbrack_{ij} =
    \begin{cases}
        \displaystyle{ \mathbf{W}_{ij}  } & \text{if } (i,j)  \in \Omega\\
         \lambda                                    & \text{if } (i,j)  \notin  \Omega
    \end{cases} 
    \ .
\end{equation*}
\\
This new weight matrix $\mathbf{W}_\lambda$ induces a norm $\Vert  \Vert_{\mathbf{W}_\lambda}$ on $\mathbb{R}^{p \times n}$, which is closely related to the seminorm $\Vert  \Vert_{\mathbf{W}}$. More precisely, for any $\lambda \in \mathbb{R}_{+*}$ and $\mathbf{Y} \in \mathbb{R}^{p \times n}$, we have $\Vert \mathbf{Y} \Vert_{\mathbf{W}} \le \Vert \mathbf{Y} \Vert_{\mathbf{W}_\lambda}$, which provides another simpler and different proof that  $\Vert  \Vert_{\mathbf{W}}$ is a continuous real-valued function on $\mathbb{R}^{p \times n}$  (see Theorem~\ref{theo3.2:box}), and also 
\begin{equation*}
\lim_{\lambda \to 0 } \, \Vert \mathbf{Y} \Vert_{\mathbf{W}_\lambda} = \Vert \mathbf{Y} \Vert_{\mathbf{W}} \ .
\end{equation*}
Furthermore, for any $\gamma \in \mathbb{R}_{+*}$ with $\gamma \ge \bar{\mathbf{c}}_{\varphi}$, we have the implications
\begin{equation*}
\Vert   \mathbf{X} -  \mathbf{Y} \Vert_{\mathbf{W}_\lambda} \le \sqrt{ 2 . \gamma } \Rightarrow \Vert   \mathbf{X} -  \mathbf{Y} \Vert_{\mathbf{W}} \le \sqrt{ 2 . \gamma }  \Rightarrow \varphi( \mathbf{Y} )  \le \gamma \ .
\end{equation*}
This shows that $\mathbb{R}^{p \times n}_{\le k} \cap \bar{B}_{p \times n}(\mathbf{X}, \sqrt{ 2 . \gamma }) \subset L(\gamma)$ where $\bar{B}_{p \times n}(\mathbf{X}, \sqrt{ 2 . \gamma })$ is the closed ball of center $\mathbf{X}$ and radius $\sqrt{ 2 . \gamma }$ with respect to the norm  $\Vert  \Vert_{\mathbf{W}_\lambda}$ on $\mathbb{R}^{p \times n}$. While the reciprocal inclusion $L(\gamma) \subset \bar{B}_{p \times n}(\mathbf{X}, \sqrt{ 2 . \gamma })$ is obviously false in general, the fact that $\lim_{\lambda \to 0 } \, \Vert \mathbf{Y} \Vert_{\mathbf{W}_\lambda} = \Vert \mathbf{Y} \Vert_{\mathbf{W}}$ suggests that for some weight matrices $\mathbf{W}$, it may still exist $\gamma \ge \bar{\mathbf{c}}_{\varphi}$ and $\lambda \in \mathbb{R}_{+*}$ sufficiently small such that $L(\gamma) \subset \bar{B}_{p \times n}(\mathbf{X}, \sqrt{ 2 . \gamma })$ so that $L(\gamma) =  \mathbb{R}^{p \times n}_{\le k} \cap \bar{B}_{p \times n}(\mathbf{X}, \sqrt{ 2 . \gamma })$ because of the imposed rank constraint on the $p \times n$ matrix $\mathbf{Y}$ in the formulation~\eqref{eq:P0}. In such cases, $\varphi(.)$ will have a bounded level set and, consequently, the set of global minimizers of $\varphi(.)$ will be nonempty and compact.

\subsection{Landscape connections of formulations~\ref{eq:P0} and~\ref{eq:P1} of the WLRA problem} \label{landscape_wlra:box}

As noted above, the cost functions $\varphi(.)$  and $\varphi^{*}(.)$  are obviously infinitely differentiable as they are polynomial functions of the entries of $\mathbf{Y}$ or $(\mathbf{A},\mathbf{B})$, respectively. In these conditions, a natural and more modest question to ask, in addition of the existence of an absolute minimum of these cost functions, is the following: is there a connection between the first- and second-order critical points of $\varphi(.)$  and $\varphi^{*}(.)$?

To begin with, we first derive the gradient of $\varphi(.)$ at $\mathbf{Y} \in \mathbb{R}^{p \times n}$. We have the following differentiation rule for a differentiable function $g(.)$ defined from $\mathbb{R}^{p \times n}$ to  $\mathbb{R}^{p \times n}$ and $\forall \mathbf{Y} , \mathbf{H} \in \mathbb{R}^{p \times n}$:
\begin{equation*}
\mathit{D} \big ( \mathbf{Y}  \rightarrowtail  \frac{1}{2}  \Vert  g( \mathbf{Y} )  \Vert^{2}_{F}  \big ) ( \mathbf{Y} ) \lbrack  \mathbf{H} \rbrack = \langle \mathit{D} g( \mathbf{Y} )  \lbrack  \mathbf{H} \rbrack , g( \mathbf{Y} ) \rangle_{F} \ .
\end{equation*}
Here, we have $\varphi( \mathbf{Y} ) =  \frac{1}{2}  \Vert  g( \mathbf{Y} )  \Vert^{2}_{F}$ with $g( \mathbf{Y} ) = \sqrt{\mathbf{W}} \odot ( \mathbf{X} - \mathbf{Y} )$, and we get
\begin{align*}
\mathit{D} \varphi( \mathbf{Y} )  \lbrack  \mathbf{H} \rbrack  & =   \big \langle   \mathit{D} g( \mathbf{Y} )  \lbrack  \mathbf{H} \rbrack  , g( \mathbf{Y} )   \big \rangle_{F} \\
                                                                                              & =   \big \langle   \sqrt{\mathbf{W}}  \odot   -\mathbf{H} ,  \sqrt{\mathbf{W}} \odot ( \mathbf{X} - \mathbf{Y} )   \big \rangle_{F} \\
                                                                                              & =   \big  \langle   \mathbf{W} \odot   ( \mathbf{Y} - \mathbf{X} ) , \mathbf{H}    \big \rangle_{F} \ .
\end{align*}
By the unicity of the Frobenius gradient of $\varphi(.)$, this implies that 
\begin{equation} \label{eq:D_grad_varphi}
\nabla \varphi( \mathbf{Y} ) =  \mathbf{W} \odot   ( \mathbf{Y} - \mathbf{X} ) \ , \   \forall  \mathbf{Y}  \in \mathbb{R}^{p \times n} \ .
\end{equation}
In particular, the gradient of $\varphi(.)$ at $\mathbf{X}$ is $\nabla \varphi( \mathbf{X} ) =  \mathbf{W} \odot   ( \mathbf{X} - \mathbf{X} ) = \mathbf{0}^{p \times n}$, which implies that $\mathbf{X}$  is a first-order critical point of $\varphi(.)$ if the feasible set is the whole linear space $\mathbb{R}^{p \times n}$. However, in most cases, especially when $\mathbf{W}  \in \mathbb{R}^{p \times n}_{+*}$, $\mathbf{X}$ is the unique first-order critical point of $\varphi(.)$ considered as a function defined on the whole linear space $\mathbb{R}^{p \times n}$. In other words, and as expected from Subsection~\ref{calculus:box}, for $\mathbf{Y}  \in \mathbb{R}^{p \times n}_{\le k}$, $\nabla \varphi( \mathbf{Y} )$ cannot be used alone as a test of the optimality of $\mathbf{Y}$ in solving the WLRA problem~\eqref{eq:P0} because perturbations of $\mathbf{Y}$ which take it out of the feasible set $\mathbb{R}^{p \times n}_{\le k}$ are not allowed and they may correspond to a decrease of the cost function $\varphi(.)$.

In general term, $\nabla^{2} \varphi( \mathbf{Y} )$ is a $4^{th}$ order tensor of dimension $p \times n \times p \times n$, but $\nabla^{2} \varphi( \mathbf{Y} )$ can also be viewed  as a bilinear form $\big ( \nabla^{2} \varphi( \mathbf{Y} ) \big)$ from $\mathbb{R}^{p \times n} \times \mathbb{R}^{p \times n}$ to $\mathbb{R}$ and also as a self-adjoint linear operator $\big \lbrack \nabla^{2} \varphi( \mathbf{Y} ) \big \rbrack$ from $\mathbb{R}^{p \times n}$ to $\mathbb{R}^{p \times n}$  (see Subsection~\ref{calculus:box} for details), and we have the equality
\begin{equation*}
\big ( \nabla^{2} \varphi( \mathbf{Y} ) \big ) \big (  \mathbf{C}, \mathbf{D}  \big ) =  \langle \big \lbrack \nabla^{2} \varphi( \mathbf{Y} ) \big \rbrack ( \mathbf{C} ) , \mathbf{D}  \rangle_{F}  \ , \  \forall \mathbf{C} , \mathbf{D}  \in \mathbb{R}^{p \times n} \ .
\end{equation*}
Taking into account the particular form of $\nabla \varphi( \mathbf{Y} )$ derived in equation~\eqref{eq:D_grad_varphi}, we have simply
\begin{equation*}
\big \lbrack \nabla^{2} \varphi( \mathbf{Y} ) \big \rbrack ( \mathbf{C} ) =  \mathbf{W} \odot  \mathbf{C}  \ , \ \forall \mathbf{C}   \in \mathbb{R}^{p \times n} \ ,
\end{equation*}
and, thus, the bilinear form of $\nabla^{2} \varphi( \mathbf{Y} )$ is defined by
\begin{equation*}
\big ( \nabla^{2} \varphi( \mathbf{Y} ) \big ) \big (  \mathbf{C}, \mathbf{D}  \big ) =  \big \langle  \mathbf{W} \odot  \mathbf{C}  , \mathbf{D} \big \rangle_{F}  \ , \   \forall  \mathbf{C}  , \mathbf{D}  \in \mathbb{R}^{p \times n} \ .
\end{equation*}
In particular, the Hessian quadratic form $\big ( \nabla^{2} \varphi( \mathbf{Y} ) \big )$ for any $p \times n$ matrices $\mathbf{Y}$ and $\mathbf{C}$ is simply given by
\begin{equation} \label{eq:D_hess_varphi}
\big ( \nabla^{2} \varphi( \mathbf{Y} ) \big ) \big (  \mathbf{C}, \mathbf{C}  \big ) =  \Vert \sqrt{\mathbf{W}} \odot  \mathbf{C} \Vert^{2}_{F} \ge 0  \  .
\end{equation}
Thus, $\big ( \nabla^{2} \varphi( \mathbf{Y} ) \big )$ is always positive semi-definite and is even always positive definite when $\mathbf{W}  \in  \mathbb{R}^{p \times n}_{+*}$.

In summary, for $\mathbf{Y}  \in \mathbb{R}^{p \times n}_{\le k}$, $\nabla \varphi( \mathbf{Y} )$ and  $\nabla^{2} \varphi( \mathbf{Y} )$ cannot be used alone as test conditions for the global or local optimality of $\mathbf{Y}$ in solving the WLRA problem~\eqref{eq:P0} because in most settings the unconstrained local or global minimizers of $\varphi(.)$ do not satisfy the rank constrained $\emph{rank}(  \mathbf{Y} ) \le k$ and, also, for a given matrix $\mathbf{Y}$ of rank less than $k$, not all the search directions or perturbations have to be taken into account for determining the criticality conditions only those for which the rank constraint will be satisfied.

Thus, to continue with, we now characterize precisely the critical points of the rank-constrained minimization problem~\eqref{eq:P0} over the real-algebraic variety $\mathbb{R}^{p \times n}_{\le k}$, which is a closed subset of the matrix space $\mathbb{R}^{p \times n}$ as stated in Theorem~\eqref{theo2.3:box}. To this end, we first identify $\mathbb{R}^{p \times n}$ and $\mathbb{R}^{p.n}$ with the two isomorphisms $\emph{vec}(.)$ and $\emph{mat}(.)$, defined in equations~\eqref{eq:vec} and~\eqref{eq:mat}. Next, we note that the Euclidean scalar product in $\mathbb{R}^{p.n}$ and the Frobenius inner product in $\mathbb{R}^{p \times n}$ are intimately related since
\begin{equation*}
\langle \mathbf{C} , \mathbf{D}  \rangle_{F} =  \Tr \big( \mathbf{C}^{T} \mathbf{D}  \big) = \big \langle \emph{vec}(\mathbf{C}) , \emph{vec}(\mathbf{D})  \big \rangle_{2}  \ , \   \forall  \mathbf{C}  , \mathbf{D}  \in \mathbb{R}^{p \times n}
\end{equation*}
and, reciprocally, 
\begin{equation*}
\langle \mathbf{c} , \mathbf{d}  \rangle_{2} =  \big \langle \emph{mat}(\mathbf{c}) , \emph{mat}(\mathbf{d})  \big \rangle_{F}  \ , \   \forall  \mathbf{c}  , \mathbf{d}  \in \mathbb{R}^{p.n} \ .
\end{equation*}
Based on these considerations, it is rather straightforward to extend the notions of tangent vectors, Bouligand tangent and Frechet normal cones, and metric projection in $\mathbb{R}^{p.n}$ summarized in Subsection~\ref{calculus:box}, especially, Theorem~\eqref{theo2.6:box} and equations~\eqref{eq:D_first_kkt_c}, to the case of the matrix space $\mathbb{R}^{p \times n}$.

Thus, a matrix $\mathbf{D} \in \mathbb{R}^{p \times n}$ is said to be tangent to $\mathbb{R}^{p \times n}_{\le k}$ at $\bar{\mathbf{Y}} \in \mathbb{R}^{p \times n}_{\le k}$ if there exist  a matrix sequence $(\mathbf{Y}_{i})_{i \in \mathbb{N}_{*}}$ in $\mathbb{R}^{p \times n}_{\le k}$ tending to $\bar{\mathbf{Y}}$ and  a real sequence $(\mathbf{t}_{i})_{i \in \mathbb{N}_{*}}$ in $\mathbb{R}_{+*}$ tending to zero such that
\begin{equation*}
\lim_{i \rightarrow \infty} \frac{ ( \mathbf{Y}_{i} - \bar{\mathbf{Y}} ) } {\mathbf{t}_{i}} = \mathbf{D}.
\end{equation*}
The set of all tangent  matrices to $\mathbb{R}^{p \times n}_{\le k}$, at $\bar{\mathbf{Y}}$ is a closed cone (see Theorem~\eqref{theo2.4:box} for details), also called the Bouligand tangent cone to $\mathbb{R}^{p \times n}_{\le k}$ at $\bar{\mathbf{Y}}$, and denoted by $\mathcal{T}^{\mathcal{B}}_{\bar{\mathbf{Y}}} \mathbb{R}^{p \times n}_{\le k}$, similarly to the case of the vector space $\mathbb{R}^{p.n}$ discussed in Subsection~\ref{calculus:box}. Its polar is defined by
\begin{equation*}
(\mathcal{T}^{\mathcal{B}}_{\bar{\mathbf{Y}}} \mathbb{R}^{p \times n}_{\le k})^{o} = \big \lbrace  \mathbf{D} \in \mathbb{R}^{p \times n} \  /  \  \langle \mathbf{D} , \mathbf{Y}  \rangle_{F}   \le 0  \  , \   \forall  \mathbf{Y} \in  \mathcal{T}^{\mathcal{B}}_{\bar{\mathbf{Y}}} \mathbb{R}^{p \times n}_{\le k} \big \rbrace 
\end{equation*}
and is also a closed convex cone called the Frechet normal cone to $\mathbb{R}^{p \times n}_{\le k}$ at $\bar{\mathbf{Y}}$, noted as $\mathcal{N}^{\mathcal{F}}_{\bar{\mathbf{Y}}} \mathbb{R}^{p \times n}_{\le k}$, again similarly to the case of the vector space $\mathbb{R}^{p.n}$ discussed in Subsection~\ref{calculus:box}.

Finally, a point $\bar{\mathbf{Y}} \in \mathbb{R}^{p \times n}_{\le k}$ is a Frechet first-order stationary point for the WLRA problem~\eqref{eq:P0} if one of the following equivalent conditions is satisfied
\begin{align} \label{eq:D_first_kkt_frechet2}
& \langle \nabla \varphi( \bar{\mathbf{Y}} ) , \mathbf{Y}  \rangle_{F}   \ge 0  \  , \   \forall  \mathbf{Y} \in  \mathcal{T}^{\mathcal{B}}_{\bar{\mathbf{Y}}} \mathbb{R}^{p \times n}_{\le k} \ , \nonumber \\
& -\nabla \varphi( \bar{\mathbf{Y}} ) \in \mathcal{N}^{\mathcal{F}}_{\bar{\mathbf{Y}}} \mathbb{R}^{p \times n}_{\le k} \ , \\
& P_{\mathcal{T}^{\mathcal{B}}_{\bar{\mathbf{Y}}} \mathbb{R}^{p \times n}_{\le k}} ( -\nabla \varphi( \bar{\mathbf{Y}} ) ) = \lbrace \mathbf{0}^{p \times n}  \rbrace \ , \nonumber
\end{align}
where  $\nabla \varphi( \bar{\mathbf{Y}} )$ is given by equation~\eqref{eq:D_grad_varphi} and $P_{\mathcal{T}^{\mathcal{B}}_{\bar{\mathbf{Y}}} \mathbb{R}^{p \times n}_{\le k}} ( -\nabla \varphi( \bar{\mathbf{Y}} ) )$ is the metric projection of the antigradient $-\nabla \varphi( \bar{\mathbf{Y}} )$ onto  $\mathcal{T}^{\mathcal{B}}_{\bar{\mathbf{Y}}} \mathbb{R}^{p \times n}_{\le k}$ defined by
\begin{equation*}
P_{\mathcal{T}^{\mathcal{B}}_{\bar{\mathbf{Y}}} \mathbb{R}^{p \times n}_{\le k}} ( -\nabla \varphi( \bar{\mathbf{Y}} ) ) = \text{Arg}\min_{  \mathbf{Y} \in  \mathcal{T}^{\mathcal{B}}_{\bar{\mathbf{Y}}} \mathbb{R}^{p \times n}_{\le k} }   \Vert  -\nabla \varphi( \bar{\mathbf{Y}} ) - \mathbf{Y} \Vert^{2}_{F} \  .
\end{equation*}
Note that the set $P_{\mathcal{T}^{\mathcal{B}}_{\bar{\mathbf{Y}}} \mathbb{R}^{p \times n}_{\le k}} ( -\nabla \varphi( \bar{\mathbf{Y}} ) )$ is always nonempty as $\mathcal{T}^{\mathcal{B}}_{\bar{\mathbf{Y}}} \mathbb{R}^{p \times n}_{\le k}$ is a closed cone, but it is not necessarily reduced to a singleton as $\mathcal{T}^{\mathcal{B}}_{\bar{\mathbf{Y}}} \mathbb{R}^{p \times n}_{\le k}$  is not convex in $\mathbb{R}^{p \times n}$, see Subsection~\ref{calculus:box} for details.

However, $\forall  \mathbf{Z}  \in  P_{\mathcal{T}^{\mathcal{B}}_{\bar{\mathbf{Y}}} \mathbb{R}^{p \times n}_{\le k}} ( -\nabla \varphi( \bar{\mathbf{Y}} ) )$, we have
\begin{equation*}
\Vert  \mathbf{Z} \Vert_{F} = \sqrt{ \Vert  -\nabla \varphi( \bar{\mathbf{Y}} )  \Vert^{2}_{F} -  d(  -\nabla \varphi( \bar{\mathbf{Y}} ) , \mathcal{T}^{\mathcal{B}}_{\bar{\mathbf{Y}}} \mathbb{R}^{p \times n}_{\le k} )^{2}} \  ,
\end{equation*}
where the distance from $-\nabla \varphi( \bar{\mathbf{Y}} )$ to $\mathcal{T}^{\mathcal{B}}_{\bar{\mathbf{Y}}} \mathbb{R}^{p \times n}_{\le k}$ is given by
\begin{align*}
d(  -\nabla \varphi( \bar{\mathbf{Y}} ) , \mathcal{T}^{\mathcal{B}}_{\bar{\mathbf{Y}}} \mathbb{R}^{p \times n}_{\le k} ) & = \inf_{ \mathbf{T}  \in \mathcal{T}^{\mathcal{B}}_{\bar{\mathbf{Y}}} \mathbb{R}^{p \times n}_{\le k}} \,   \Vert  -\nabla \varphi( \bar{\mathbf{Y}} )   - \mathbf{T} \Vert_{F}    \\
                                                              & = \min_{ \mathbf{T}  \in \mathcal{T}^{\mathcal{B}}_{\bar{\mathbf{Y}}} \mathbb{R}^{p \times n}_{\le k}} \,   \Vert  -\nabla \varphi( \bar{\mathbf{Y}} )   - \mathbf{T} \Vert_{F} \  ,
\end{align*}
again because $\mathcal{T}^{\mathcal{B}}_{\bar{\mathbf{Y}}} \mathbb{R}^{p \times n}_{\le k}$ is a closed set. In other words, all elements of $P_{\mathcal{T}^{\mathcal{B}}_{\bar{\mathbf{Y}}} \mathbb{R}^{p \times n}_{\le k}} ( -\nabla \varphi( \bar{\mathbf{Y}} ) )$ have the same Frobenius norm and by the same small abuse of notation as used in equation~\eqref{eq:D_first_kkt_c2} of Subsection~\ref{calculus:box}, the Frechet first-order stationary condition for $\varphi(.)$ can be expressed as 
\begin{equation*}
  \Vert  P_{\mathcal{T}^{\mathcal{B}}_{\bar{\mathbf{Y}}} \mathbb{R}^{p \times n}_{\le k}} ( -\nabla \varphi( \bar{\mathbf{Y}} ) ) \Vert_{F}  = 0 \ ,
\end{equation*}
where $P_{\mathcal{T}^{\mathcal{B}}_{\bar{\mathbf{Y}}} \mathbb{R}^{p \times n}_{\le k}} ( -\nabla \varphi( \bar{\mathbf{Y}} ) )$ designs now any of its elements. However, to use these results, we first need to find convenient practical expressions for $\mathcal{T}^{\mathcal{B}}_{\bar{\mathbf{Y}}} \mathbb{R}^{p \times n}_{\le k}$ and the metric projection operator onto this closed set.

In order to derive a more convenient way for checking if a matrix $\bar{\mathbf{Y}} \in \mathbb{R}^{p \times n}_{\le k}$ is a Frechet first-order stationary point for $\varphi(.)$, we first note that the set $\mathbb{R}^{p \times n}_{\le k}$ stratifies into the set $\mathbb{R}^{p \times n}_{s}$ for $s= 1, \cdots, k$, e.g.,
\begin{equation*}
\mathbb{R}^{p \times n}_{\le k} = \bigcup^{k}_{s=1} \mathbb{R}^{p \times n}_{s} \ .
\end{equation*}
Furthermore, it is well-known, that each set $\mathbb{R}^{p \times n}_{s}$ is a smooth submanifold of dimension $(p+n-s).s$ embedded in $\mathbb{R}^{p \times n}$ and that  its tangent space at $\mathbf{Y} \in \mathbb{R}^{p \times n}_{s}$ is given by 
\begin{align*}
\mathcal{T}_{\mathbf{Y}} \mathbb{R}^{p \times n}_{s} & = \big \lbrace \mathbf{U}_{\mathbf{Y}} \mathbf{M} \mathbf{V}^{T}_{\mathbf{Y}} + \mathbf{U} \mathbf{V}^{T}_{\mathbf{Y}} + \mathbf{U}_{\mathbf{Y}}  \mathbf{V}^{T} \  /  \ \mathbf{M}  \in  \mathbb{R}^{s \times s}, \mathbf{U}  \in  \mathbb{R}^{p \times s}, \mathbf{V}  \in  \mathbb{R}^{n \times s}  \\
 & \qquad \text{ with }  \mathbf{U}^{T}_{\mathbf{Y}} \mathbf{U} = \mathbf{V}^{T}_{\mathbf{Y}} \mathbf{V} = \mathbf{0}^{s \times s} \big \rbrace \\
 & =  \big \lbrace    \lbrack  \mathbf{U}_{\mathbf{Y}}    \mathbf{U}^{\bot}_{\mathbf{Y}}  \rbrack 
 \begin{bmatrix}
\mathbf{A} & \mathbf{B} \\
\mathbf{C} & \mathbf{0}^{(p-s) \times (n-s)} 
\end{bmatrix}   
\lbrack  \mathbf{V}_{\mathbf{Y}}    \mathbf{V}^{\bot}_{\mathbf{Y}}  \rbrack^{T}  \\
 & \qquad  \text{ with } \mathbf{A}  \in  \mathbb{R}^{s \times s}, \mathbf{B}  \in  \mathbb{R}^{s \times (n-s)}   \text{ and } \mathbf{C}  \in  \mathbb{R}^{(p-s) \times s}  \big \rbrace \ ,
\end{align*}
where $\mathbf{Y} =  \mathbf{U}_{\mathbf{Y}}  \Sigma_{\mathbf{Y}} \mathbf{V}^{T}_{\mathbf{Y}}$ is the thin SVD of $\mathbf{Y}$ with $\mathbf{U}^{T}_{\mathbf{Y}} \mathbf{U}_{\mathbf{Y}}  = \mathbf{V}^{T}_{\mathbf{Y}} \mathbf{V}_{\mathbf{Y}}  = \mathbf{I}_{s}$ and $\Sigma_{\mathbf{Y}}$ is a $s \times s$ diagonal matrix with strictly positive diagonal elements  (e.g., the singular values of $\mathbf{Y}$), and $\lbrack  \mathbf{U}_{\mathbf{Y}}    \mathbf{U}^{\bot}_{\mathbf{Y}}  \rbrack$ and  $\lbrack  \mathbf{V}_{\mathbf{Y}}    \mathbf{V}^{\bot}_{\mathbf{Y}}  \rbrack$ are, respectively, $p \times p$ and $n \times n$ orthogonal matrices. See Example 8.14 of Lee~\cite{L2003}, Section 7.5 of Boumal~\cite{B2023} or Proposition 4.1 of Helmke and Shayman~\cite{HS1995} for proofs and further details. Furthermore, the equivalence of the two definitions of the tangent space $\mathcal{T}_{\mathbf{Y}} \mathbb{R}^{p \times n}_{s}$ can be easily verified by direct computations.

Interestingly, if  each $\mathbf{Y} \in \mathbb{R}^{p \times n}_{s}$ is identified by its singular triplets $(\mathbf{U}_{\mathbf{Y}},  \Sigma_{\mathbf{Y}}, \mathbf{V}_{\mathbf{Y}} )$, then the first formulation of $\mathcal{T}_{\mathbf{Y}} \mathbb{R}^{p \times n}_{s}$ shows that, to represent an element of $\mathcal{T}_{\mathbf{Y}} \mathbb{R}^{p \times n}_{s}$, we only need to store the small matrices $\mathbf{M}$, $\mathbf{U}$ and $\mathbf{V}$. Furthermore, this formulation also shows that the elements of  $\mathcal{T}_{\mathbf{Y}} \mathbb{R}^{p \times n}_{s}$ have a rank of at most $2.s$. On the other hand, the second formulation is useful for deriving the normal space to $\mathbb{R}^{p \times n}_{s}$ at $\mathbf{Y}$, noted $\mathcal{N}_{\mathbf{Y}} \mathbb{R}^{p \times n}_{s}$, which is the orthogonal complement of $\mathcal{T}_{\mathbf{Y}} \mathbb{R}^{p \times n}_{s}$ in $\mathbb{R}^{p \times n}$ with respect to the Frobenius inner product:
\begin{equation*}
\mathcal{N}_{\mathbf{Y}} \mathbb{R}^{p \times n}_{s} = ( \mathcal{T}_{\mathbf{Y}} \mathbb{R}^{p \times n}_{s} )^{\bot} \ ,
\end{equation*}
and also the orthogonal projectors on both $\mathcal{T}_{\mathbf{Y}} \mathbb{R}^{p \times n}_{s}$ and $\mathcal{N}_{\mathbf{Y}} \mathbb{R}^{p \times n}_{s}$ as we will see now.

First, the second formulation reveals immediately the dimension of $\mathcal{T}_{\mathbf{Y}} \mathbb{R}^{p \times n}_{s}$ as
\begin{equation*}
\emph{dim}( \mathcal{T}_{\mathbf{Y}} \mathbb{R}^{p \times n}_{s} ) = s.s + s.(p-s) + s.(n-s) = s.(p+n-s) \ .
\end{equation*}
Next, it is obvious from this formulation of $\mathcal{T}_{\mathbf{Y}} \mathbb{R}^{p \times n}_{s}$ that $\mathcal{N}_{\mathbf{Y}} \mathbb{R}^{p \times n}_{s}$ is equal to
\begin{equation}  \label{eq:D_normal_space_varphi}
\mathcal{N}_{\mathbf{Y}} \mathbb{R}^{p \times n}_{s} =  \Big \lbrace    \mathbf{U}^{\bot}_{\mathbf{Y}} \mathbf{N}  ( \mathbf{V}^{\bot}_{\mathbf{Y}} )^{T} \ \text{ with } \mathbf{N}  \in  \mathbb{R}^{(p-s) \times (n-s)} \Big \rbrace \ .
\end{equation}
Obviously and as expected, we have
\begin{equation*}
\emph{dim}( \mathcal{N}_{\mathbf{Y}} \mathbb{R}^{p \times n}_{s} ) = (p-s).(n-s) = p.n - s.(p+n-s) =  \emph{dim}(  \mathbb{R}^{p \times n} ) - \emph{dim}( \mathcal{T}_{\mathbf{Y}} \mathbb{R}^{p \times n}_{s} )
\end{equation*}
and the maximum rank of the matrix elements of $\mathcal{N}_{\mathbf{Y}} \mathbb{R}^{p \times n}_{s}$ is $\emph{min}( p, n) - s$ according to equation~\eqref{eq:rank2}. Next, by definition, the orthogonal projection of an arbitrary $\mathbf{Z} \in \mathbb{R}^{p \times n}$ onto $\mathcal{T}_{\mathbf{Y}} \mathbb{R}^{p \times n}_{s}$ satisfies both
\begin{equation*}
\mathbf{Z} - \mathbf{P}_{ \mathcal{T}_{\mathbf{Y}} \mathbb{R}^{p \times n}_{s} } ( \mathbf{Z} ) =  \mathbf{U}^{\bot}_{\mathbf{Y}} \mathbf{N}  ( \mathbf{V}^{\bot}_{\mathbf{Y}} )^{T} \ ,
\end{equation*}
for some $\mathbf{N}  \in  \mathbb{R}^{(p-s) \times (n-s)}$,  and 
\begin{equation*}
\mathbf{P}_{ \mathcal{T}_{\mathbf{Y}} \mathbb{R}^{p \times n}_{s} } ( \mathbf{Z} ) =  \mathbf{U}_{\mathbf{Y}} \mathbf{M} \mathbf{V}^{T}_{\mathbf{Y}} + \mathbf{U} \mathbf{V}^{T}_{\mathbf{Y}} + \mathbf{U}_{\mathbf{Y}}  \mathbf{V}^{T} \ ,
\end{equation*}
for some $\mathbf{M}  \in  \mathbb{R}^{s \times s}, \mathbf{U}  \in  \mathbb{R}^{p \times s}$ and $\mathbf{V}  \in  \mathbb{R}^{n \times s}$ with $\mathbf{U}^{T}_{\mathbf{Y}} \mathbf{U} = \mathbf{V}^{T}_{\mathbf{Y}} \mathbf{V} = \mathbf{0}^{s \times s}$. Combined, these two statements imply that
\begin{equation*}
\mathbf{Z} = \mathbf{U}_{\mathbf{Y}} \mathbf{M} \mathbf{V}^{T}_{\mathbf{Y}} + \mathbf{U} \mathbf{V}^{T}_{\mathbf{Y}} + \mathbf{U}_{\mathbf{Y}}  \mathbf{V}^{T} +  \mathbf{U}^{\bot}_{\mathbf{Y}} \mathbf{N}  ( \mathbf{V}^{\bot}_{\mathbf{Y}} )^{T} \ .
\end{equation*}
If we define now the orthogonal projectors associated with the column and row spaces of $\mathbf{Y}$ and their orthogonal complements
\begin{equation*}
\mathbf{P}_{\mathbf{U}} =  \mathbf{U}_{\mathbf{Y}}  \mathbf{U}^{T}_{\mathbf{Y}} \ , \ \mathbf{P}_{\mathbf{V}} =  \mathbf{V}_{\mathbf{Y}}  \mathbf{V}^{T}_{\mathbf{Y}} \ , \  \mathbf{P}^{\bot}_{\mathbf{U}} =  \mathbf{I}_{p} - \mathbf{P}_{\mathbf{U}}  \text{ and } \mathbf{P}^{\bot}_{\mathbf{V}} =  \mathbf{I}_{n} - \mathbf{P}_{\mathbf{V}} \ ,
\end{equation*}
we have, using orthogonal relationships,
\begin{align*}
\mathbf{P}_{\mathbf{U}} \mathbf{Z} \mathbf{P}_{\mathbf{V}}           & = \big (  \mathbf{U}_{\mathbf{Y}} \mathbf{M}  \mathbf{V}^{T}_{\mathbf{Y}} + \mathbf{U}_{\mathbf{Y}} \mathbf{V}^{T} \big ) \mathbf{P}_{\mathbf{V}} = \mathbf{U}_{\mathbf{Y}} \mathbf{M}  \mathbf{V}^{T}_{\mathbf{Y}} \ , \\
\mathbf{P}^{\bot}_{\mathbf{U}} \mathbf{Z} \mathbf{P}_{\mathbf{V}} & = \big ( \mathbf{U}  \mathbf{V}^{T}_{\mathbf{Y}} +  \mathbf{U}^{\bot}_{\mathbf{Y}} \mathbf{N} ( \mathbf{V}^{\bot}_{\mathbf{Y}} )^{T} \big ) \mathbf{P}_{\mathbf{V}} = \mathbf{U}  \mathbf{V}^{T}_{\mathbf{Y}} \ , \\
\mathbf{P}_{\mathbf{U}} \mathbf{Z} \mathbf{P}^{\bot}_{\mathbf{V}} & = \big ( \mathbf{U}_{\mathbf{Y}}  \mathbf{M}  \mathbf{V}^{T}_{\mathbf{Y}} +  \mathbf{U}_{\mathbf{Y}} \mathbf{V}^{T} \big )  \mathbf{P}^{\bot}_{\mathbf{V}} = \mathbf{U}_{\mathbf{Y}} \mathbf{V}^{T} \ .
\end{align*}
Using these results, we deduce that the orthogonal projector onto  $\mathcal{T}_{\mathbf{Y}} \mathbb{R}^{p \times n}_{s}$ is given, equivalently, by
\begin{align*}
\mathbf{P}_{ \mathcal{T}_{\mathbf{Y}} \mathbb{R}^{p \times n}_{s} } ( \mathbf{Z} )  & = \mathbf{P}_{\mathbf{U}} \mathbf{Z} \mathbf{P}_{\mathbf{V}} + \mathbf{P}^{\bot}_{\mathbf{U}} \mathbf{Z} \mathbf{P}_{\mathbf{V}} +  \mathbf{P}_{\mathbf{U}} \mathbf{Z} \mathbf{P}^{\bot}_{\mathbf{V}}  \\
                                                                                                                                   & = \mathbf{Z} \mathbf{P}_{\mathbf{V}} +  \mathbf{P}_{\mathbf{U}} \mathbf{Z} \mathbf{P}^{\bot}_{\mathbf{V}} = \mathbf{Z} ( \mathbf{V}_{\mathbf{Y}}  \mathbf{V}^{T}_{\mathbf{Y}} ) +( \mathbf{U}_{\mathbf{Y}}  \mathbf{U}^{T}_{\mathbf{Y}}  ) \mathbf{Z} (\mathbf{I}_{n} - \mathbf{V}_{\mathbf{Y}}  \mathbf{V}^{T}_{\mathbf{Y}}  ) \\
                                                                                                                                   & =   \mathbf{P}_{\mathbf{U}} \mathbf{Z} +  \mathbf{P}^{\bot}_{\mathbf{U}} \mathbf{Z} \mathbf{P}_{\mathbf{V}} = ( \mathbf{U}_{\mathbf{Y}}  \mathbf{U}^{T}_{\mathbf{Y}}  ) \mathbf{Z} +  ( \mathbf{I}_{p} - \mathbf{U}_{\mathbf{Y}}  \mathbf{U}^{T}_{\mathbf{Y}}  )  \mathbf{Z} ( \mathbf{V}_{\mathbf{Y}}  \mathbf{V}^{T}_{\mathbf{Y}} ) \ ,
\end{align*}
from which, we can also derive the orthogonal projector onto $\mathcal{N}_{\mathbf{Y}} \mathbb{R}^{p \times n}_{s}$ as
\begin{align*}
\mathbf{P}_{ \mathcal{N}_{\mathbf{Y}} \mathbb{R}^{p \times n}_{s} } ( \mathbf{Z} )  & = \mathbf{P}^{\bot}_{ \mathcal{T}_{\mathbf{Y}} \mathbb{R}^{p \times n}_{s} } ( \mathbf{Z} )  = \mathbf{Z} - \mathbf{P}_{ \mathcal{T}_{\mathbf{Y}} \mathbb{R}^{p \times n}_{s} } ( \mathbf{Z} )  \\
                                                                                                                                   & =   \mathbf{Z} - ( \mathbf{U}_{\mathbf{Y}}  \mathbf{U}^{T}_{\mathbf{Y}} )  \mathbf{Z} - ( \mathbf{I}_{p}  - \mathbf{U}_{\mathbf{Y}}  \mathbf{U}^{T}_{\mathbf{Y}} )  \mathbf{Z} ( \mathbf{V}_{\mathbf{Y}}  \mathbf{V}^{T}_{\mathbf{Y}} ) \\
                                                                                                                                   & =   ( \mathbf{I}_{p}  - \mathbf{U}_{\mathbf{Y}}  \mathbf{U}^{T}_{\mathbf{Y}} )  \mathbf{Z}   ( \mathbf{I}_{n}  - \mathbf{V}_{\mathbf{Y}}  \mathbf{V}^{T}_{\mathbf{Y}} ) = \mathbf{P}^{\bot}_{\mathbf{U}}  \mathbf{Z}   \mathbf{P}^{\bot}_{\mathbf{V}}  \ .
\end{align*}
In these conditions, if $\mathbf{Z} \in \mathcal{N}_{\mathbf{Y}} \mathbb{R}^{p \times n}_{s}$, we have $\mathbf{P}_{ \mathcal{N}_{\mathbf{Y}} \mathbb{R}^{p \times n}_{s} } ( \mathbf{Z} ) =   \mathbf{Z}$, which implies that
\begin{align*}
\mathbf{U}^{T}_{\mathbf{Y}}  \mathbf{Z} & = \mathbf{U}^{T}_{\mathbf{Y}} ( \mathbf{I}_{p}  - \mathbf{U}_{\mathbf{Y}}  \mathbf{U}^{T}_{\mathbf{Y}} )  \mathbf{Z}   ( \mathbf{I}_{n}  - \mathbf{V}_{\mathbf{Y}}  \mathbf{V}^{T}_{\mathbf{Y}} ) =  \mathbf{0}^{s \times n}  \ ,  \\
 \mathbf{Z}   \mathbf{V}_{\mathbf{Y}} & = ( \mathbf{I}_{p}  - \mathbf{U}_{\mathbf{Y}}  \mathbf{U}^{T}_{\mathbf{Y}} )  \mathbf{Z}   ( \mathbf{I}_{n}  - \mathbf{V}_{\mathbf{Y}}  \mathbf{V}^{T}_{\mathbf{Y}} ) \mathbf{V}_{\mathbf{Y}} = \mathbf{0}^{p \times s} \ .
\end{align*}
Reciprocally, if $\mathbf{Z} \in \mathbb{R}^{p \times n}$ with $\mathbf{U}^{T}_{\mathbf{Y}}  \mathbf{Z} =  \mathbf{0}^{s \times n}$ and $\mathbf{Z}   \mathbf{V}_{\mathbf{Y}} =  \mathbf{0}^{p \times s}$, we have
\begin{align*}
\mathbf{P}_{ \mathcal{T}_{\mathbf{Y}} \mathbb{R}^{p \times n}_{s} } ( \mathbf{Z} )  & =   ( \mathbf{Z} \mathbf{V}_{\mathbf{Y}} ) \mathbf{V}^{T}_{\mathbf{Y}} +  \mathbf{U}_{\mathbf{Y}} ( \mathbf{U}^{T}_{\mathbf{Y}}  \mathbf{Z} ) (\mathbf{I}_{n} - \mathbf{V}_{\mathbf{Y}}  \mathbf{V}^{T}_{\mathbf{Y}}  ) \\
                                                                                                                                    & = \mathbf{0}^{p \times n}   +  \mathbf{0}^{p \times n} = \mathbf{0}^{p \times n} 
\end{align*}
and $\mathbf{Z} \in \mathcal{N}_{\mathbf{Y}} \mathbb{R}^{p \times n}_{s}$. In other words, we get an alternative formulation of $\mathcal{N}_{\mathbf{Y}} \mathbb{R}^{p \times n}_{s}$ as
\begin{equation*}
\mathcal{N}_{\mathbf{Y}} \mathbb{R}^{p \times n}_{s} = \big \lbrace   \mathbf{Z} \in \mathbb{R}^{p \times n}  \text{ with }  \mathbf{U}^{T}_{\mathbf{Y}}  \mathbf{Z} =  \mathbf{0}^{s \times n}  \text{ and }   \mathbf{Z}   \mathbf{V}_{\mathbf{Y}} =  \mathbf{0}^{p \times s}  \big \rbrace \ ,
\end{equation*}
where the columns of $\mathbf{U}_{\mathbf{Y}}$ and $\mathbf{V}_{\mathbf{Y}}$ are, respectively, the leading $s$ left and right singular vectors of $\mathbf{Y}$, which is of rank $s$.

Armed with these various results on the smooth manifold $\mathbb{R}^{p \times n}_{s}$ embedded in $\mathbb{R}^{p \times n}$, we can now reformulate the definitions of the Bouligand tangent cone to  $\mathbb{R}^{p \times n}_{\le k}$ at a matrix $\mathbf{Y}$ of rank $s \le k$ and of the metric projection onto that closed set as follow.
 \begin{theo3.4} \label{theo3.4:box}
 Let  $\mathbf{Y} \in \mathbb{R}^{p \times n}_{\le k}$ with $\emph{rank}( \mathbf{Y} ) = s \le k$, the Bouligand tangent cone to  $\mathbb{R}^{p \times n}_{\le k}$ at  $\mathbf{Y}$ is given by
 \begin{equation*}
 \mathcal{T}^{\mathcal{B}}_{\mathbf{Y}} \mathbb{R}^{p \times n}_{\le k} =  \mathcal{T}_{\mathbf{Y}} \mathbb{R}^{p \times n}_{s}  \oplus \big (   \mathcal{N}_{\mathbf{Y}} \mathbb{R}^{p \times n}_{s}  \cap   \mathbb{R}^{p \times n}_{\le k-s}  \big ) \ ,
\end{equation*}
where $\oplus$ stands for a direct orthogonal sum with respect to the Frobenius inner product in $\mathbb{R}^{p \times n}$.

In addition, the metric projection of an arbitrary $\mathbf{Z} \in \mathbb{R}^{p \times n}$ onto $\mathcal{T}^{\mathcal{B}}_{\mathbf{Y}} \mathbb{R}^{p \times n}_{\le k}$ is given by
 \begin{equation*}
 \mathbf{P}_{ \mathcal{T}^{\mathcal{B}}_{\mathbf{Y}} \mathbb{R}^{p \times n}_{\le k} } ( \mathbf{Z} ) = \mathbf{P}_{ \mathcal{T}_{\mathbf{Y}} \mathbb{R}^{p \times n}_{s} } ( \mathbf{Z} ) +  \mathbf{P}_{  \mathbb{R}^{p \times n}_{\le k-s} } \big(  \mathbf{P}_{ \mathcal{N}_{\mathbf{Y}} \mathbb{R}^{p \times n}_{s} } ( \mathbf{Z} ) \big ) \ ,
\end{equation*}
where $\mathbf{P}_{ \mathcal{T}_{\mathbf{Y}} \mathbb{R}^{p \times n}_{s} } (.)$ and $\mathbf{P}_{ \mathcal{N}_{\mathbf{Y}} \mathbb{R}^{p \times n}_{s} } (.)$ are the two unique complementary orthogonal projectors onto the linear subspaces $\mathcal{T}_{\mathbf{Y}} \mathbb{R}^{p \times n}_{s}$ and $\mathcal{N}_{\mathbf{Y}} \mathbb{R}^{p \times n}_{s}$, which are orthogonal to each other with respect to the Frobenius inner product in $\mathbb{R}^{p \times n}$, and where $\mathbf{P}_{  \mathbb{R}^{p \times n}_{\le k-s} } (.)$ is the metric projection onto the closed set $\mathbb{R}^{p \times n}_{\le k-s}$.
\end{theo3.4}
\begin{proof}
For a proof, see Theorem 3.2 and Corollary 3.3 of Schneider and Uschmajew~\cite{SU2015}, Theorem 6.1 of Cason et al.~\cite{CAD2013} or Example 20.5 of Harris~\cite{H1992}.
\end{proof}
First note that, in Theorem~\ref{theo3.4:box}, $\mathbf{P}_{ \mathcal{T}^{\mathcal{B}}_{\mathbf{Y}} \mathbb{R}^{p \times n}_{\le k} } ( \mathbf{Z} )$ is always a nonempty set as $\mathcal{T}^{\mathcal{B}}_{\mathbf{Y}} \mathbb{R}^{p \times n}_{\le k}$ is a closed cone, but it is not neccessarily reduced to a singleton as $\mathcal{T}^{\mathcal{B}}_{\mathbf{Y}} \mathbb{R}^{p \times n}_{\le k}$ is not convex. More precisely, the cardinality of $\mathbf{P}_{ \mathcal{T}^{\mathcal{B}}_{\mathbf{Y}} \mathbb{R}^{p \times n}_{\le k} } ( \mathbf{Z} )$ relies on the cardinality of the set $ \mathbf{P}_{  \mathbb{R}^{p \times n}_{\le k-s} } \big(  \mathbf{P}_{ \mathcal{N}_{\mathbf{Y}} \mathbb{R}^{p \times n}_{s} } ( \mathbf{Z} ) \big )$.

For an arbitrary $\mathbf{T}  \in \mathbb{R}^{p \times n}$, the metric projection of $\mathbf{T}$ onto $\mathbb{R}^{p \times n}_{\le k-s}$ is the set defined by
 \begin{equation*}
\mathbf{P}_{  \mathbb{R}^{p \times n}_{\le k-s} } \big( \mathbf{T} \big ) = \text{Arg} \min_{ \mathbf{Z}  \in  \mathbb{R}^{p \times n}_{\le k-s} }    \Vert   \mathbf{T} -  \mathbf{Z}  \Vert_{F} \ .
\end{equation*}
Thus, the elements of $\mathbf{P}_{  \mathbb{R}^{p \times n}_{\le k-s} } ( \mathbf{T} )$ are easily determined with the help of the Eckart-Young Theorem~\ref{theo2.1:box} and are the best approximation of rank at most $k-s$ of $\mathbf{T}$ with respect to the Frobenius  norm. In other words,  $ \mathbf{P}_{  \mathbb{R}^{p \times n}_{\le k-s} } \big(  \mathbf{P}_{ \mathcal{N}_{\mathbf{Y}} \mathbb{R}^{p \times n}_{s} } ( \mathbf{Z} ) \big )$ is single-valued when
 \begin{equation*}
  \sigma_{ k-s} \big ( \mathbf{P}_{ \mathcal{N}_{\mathbf{Y}} \mathbb{R}^{p \times n}_{s} } ( \mathbf{Z} ) \big) > \sigma_{ k-s+1} \big( \mathbf{P}_{ \mathcal{N}_{\mathbf{Y}} \mathbb{R}^{p \times n}_{s} } ( \mathbf{Z} ) \big) \ ,
\end{equation*}
in which case its unique element is given by the truncated SVD of rank $k-s$ of  $\mathbf{P}_{ \mathcal{N}_{\mathbf{Y}} \mathbb{R}^{p \times n}_{s} } ( \mathbf{Z} )$ according to Theorem~\ref{theo2.1:box}), or, when,
 \begin{equation*}
  \sigma_{ k-s} \big ( \mathbf{P}_{ \mathcal{N}_{\mathbf{Y}} \mathbb{R}^{p \times n}_{s} } ( \mathbf{Z} ) \big) = 0 \ ,
\end{equation*}
in which case
 \begin{equation*}
  \mathbf{P}_{  \mathbb{R}^{p \times n}_{\le k-s} } \big(  \mathbf{P}_{ \mathcal{N}_{\mathbf{Y}} \mathbb{R}^{p \times n}_{s} } ( \mathbf{Z} ) \big ) =  \mathbf{P}_{ \mathcal{N}_{\mathbf{Y}} \mathbb{R}^{p \times n}_{s} } ( \mathbf{Z} ) \ .
\end{equation*}
In the above equations, $\sigma_{i} ( \mathbf{T} )$ denotes the $i^{th}$  largest singular value of the matrix $\mathbf{T}$. Furthermore, when $\mathbf{P}_{  \mathbb{R}^{p \times n}_{\le k-s} } \big(  \mathbf{P}_{ \mathcal{N}_{\mathbf{Y}} \mathbb{R}^{p \times n}_{s} } ( \mathbf{Z} ) \big )$  is single-valued then $ \mathbf{P}_{ \mathcal{T}^{\mathcal{B}}_{\mathbf{Y}} \mathbb{R}^{p \times n}_{\le k} } ( \mathbf{Z} )$ is also single-valued according to Theorem~\ref{theo3.4:box}.

In order to clarify the practical meaning of Theorem~\ref{theo3.4:box}, it is now useful to distinguish the two cases  $\emph{rank}( \mathbf{Y} ) = s = k$ and  $\emph{rank}( \mathbf{Y} ) = s < k$.

Obviously, in the first case, when $s=k$, we get
 \begin{equation*}
 \mathcal{T}^{\mathcal{B}}_{\mathbf{Y}} \mathbb{R}^{p \times n}_{\le k} =  \mathcal{T}_{\mathbf{Y}} \mathbb{R}^{p \times n}_{k} 
\end{equation*}
and 
 \begin{equation*}
 \mathbf{P}_{ \mathcal{T}^{\mathcal{B}}_{\mathbf{Y}} \mathbb{R}^{p \times n}_{\le k} } (.) = \mathbf{P}_{ \mathcal{T}_{\mathbf{Y}} \mathbb{R}^{p \times n}_{k} } (.) \ .
\end{equation*}
In words, when $\emph{rank}( \mathbf{Y} ) = k$, the Bouligand tangent cone to $\mathbb{R}^{p \times n}_{\le k}$ at $\mathbf{Y}$ coincides with the tangent linear space to  $\mathbb{R}^{p \times n}_{k}$ at $\mathbf{Y}$. Furthermore, the metric projection onto this Bouligand tangent cone is nothing else then the orthogonal projector onto the tangent linear space to $\mathbb{R}^{p \times n}_{k}$ at $\mathbf{Y}$. Finally, from these results, we deduce immediately that the Frechet normal cone to $\mathbb{R}^{p \times n}_{\le k}$ at $\mathbf{Y}$, which is defined as the polar of $ \mathcal{T}^{\mathcal{B}}_{\mathbf{Y}} \mathbb{R}^{p \times n}_{\le k}$, also coincides with the normal space to $\mathbb{R}^{p \times n}_{k}$ at $\mathbf{Y}$ (e.g., the orthogonal complement of  $\mathcal{T}_{\mathbf{Y}} \mathbb{R}^{p \times n}_{k}$ in $\mathbb{R}^{p \times n}$ with respect to the Frobenius inner product) when $\emph{rank}( \mathbf{Y} ) = k$.  If $\mathbf{Z} \in  \mathcal{N}^{\mathcal{F}}_{\mathbf{Y}} \mathbb{R}^{p \times n}_{\le k} =  (\mathcal{T}^{\mathcal{B}}_{\mathbf{Y}} \mathbb{R}^{p \times n}_{\le k} )^{o} = (\mathcal{T}_{\mathbf{Y}} \mathbb{R}^{p \times n}_{k} )^{o}$ then, by definition,
 \begin{equation*}
\langle \mathbf{Z} , \mathbf{Q}  \rangle_{F}   \le 0  \  , \   \forall  \mathbf{Q} \in  \mathcal{T}_{\mathbf{Y}} \mathbb{R}^{p \times n}_{k} .
\end{equation*}
However, since $\mathcal{T}_{\mathbf{Y}} \mathbb{R}^{p \times n}_{k}$ is a linear space, if $ \mathbf{Q} \in  \mathcal{T}_{\mathbf{Y}} \mathbb{R}^{p \times n}_{k}$ then $-\mathbf{Q}$ also belongs to $\mathcal{T}_{\mathbf{Y}} \mathbb{R}^{p \times n}_{k}$, from which we deduce 
 \begin{equation*}
\langle \mathbf{Z} , \mathbf{Q}  \rangle_{F}   \ge 0  \  , \   \forall  \mathbf{Q} \in  \mathcal{T}_{\mathbf{Y}} \mathbb{R}^{p \times n}_{k}
\end{equation*}
and we get the equivalences
 \begin{equation*}
 \mathbf{Z} \in  \mathcal{N}^{\mathcal{F}}_{\mathbf{Y}} \mathbb{R}^{p \times n}_{\le k} \Longleftrightarrow \langle \mathbf{Z} , \mathbf{Q}  \rangle_{F}  = 0  \  , \   \forall  \mathbf{Q} \in  \mathcal{T}_{\mathbf{Y}} \mathbb{R}^{p \times n}_{k}  \Longleftrightarrow \mathbf{Z} \in  ( \mathcal{T}_{\mathbf{Y}} \mathbb{R}^{p \times n}_{k} )^{\bot} = \mathcal{N}_{\mathbf{Y}} \mathbb{R}^{p \times n}_{k} \ .
\end{equation*}

Summarizing the preceding results, when $\mathbf{Y} \in  \mathbb{R}^{p \times n}_{k}$ and $\mathbf{Z} \in  \mathbb{R}^{p \times n}$, we have
 \begin{align}  \label{eq:D_riemann_bouligand_equi}
 \mathcal{T}^{\mathcal{B}}_{\mathbf{Y}} \mathbb{R}^{p \times n}_{\le k}  & =  \mathcal{T}_{\mathbf{Y}} \mathbb{R}^{p \times n}_{k}  \nonumber  \\
      & = \big \lbrace \mathbf{U}_{\mathbf{Y}} \mathbf{M} \mathbf{V}^{T}_{\mathbf{Y}} + \mathbf{U} \mathbf{V}^{T}_{\mathbf{Y}} + \mathbf{U}_{\mathbf{Y}}  \mathbf{V}^{T} \  /  \ \mathbf{M}  \in  \mathbb{R}^{k \times k}, \mathbf{U}  \in  \mathbb{R}^{p \times k}, \mathbf{V}  \in  \mathbb{R}^{n \times k}  \nonumber \\
      & \qquad \text{ with }  \mathbf{U}^{T}_{\mathbf{Y}} \mathbf{U} = \mathbf{V}^{T}_{\mathbf{Y}} \mathbf{V} = \mathbf{0}^{k \times } \big \rbrace ,  \nonumber \\
\mathcal{N}^{\mathcal{F}}_{\mathbf{Y}} \mathbb{R}^{p \times n}_{\le k}  & = ( \mathcal{T}_{\mathbf{Y}} \mathbb{R}^{p \times n}_{k} )^{\bot} = \mathcal{N}_{\mathbf{Y}} \mathbb{R}^{p \times n}_{k} \nonumber \\
      & = \big \lbrace   \mathbf{T} \in \mathbb{R}^{p \times n}  \text{ with }  \mathbf{U}^{T}_{\mathbf{Y}}  \mathbf{T} =  \mathbf{0}^{k \times n}  \text{ and }   \mathbf{T}   \mathbf{V}_{\mathbf{Y}} =  \mathbf{0}^{p \times k}  \big \rbrace \ ,  \\
  \mathbf{P}_{ \mathcal{T}^{\mathcal{B}}_{\mathbf{Y}} \mathbb{R}^{p \times n}_{\le k} } ( \mathbf{Z} ) & = \mathbf{P}_{ \mathcal{T}_{\mathbf{Y}} \mathbb{R}^{p \times n}_{k} } ( \mathbf{Z} ) \nonumber \\
      & = \mathbf{Z} ( \mathbf{V}_{\mathbf{Y}}  \mathbf{V}^{T}_{\mathbf{Y}} ) +( \mathbf{U}_{\mathbf{Y}}  \mathbf{U}^{T}_{\mathbf{Y}}  ) \mathbf{Z} (\mathbf{I}_{n} - \mathbf{V}_{\mathbf{Y}}  \mathbf{V}^{T}_{\mathbf{Y}}  ) \nonumber \\
      & =  ( \mathbf{U}_{\mathbf{Y}}  \mathbf{U}^{T}_{\mathbf{Y}}  ) \mathbf{Z} +  ( \mathbf{I}_{p} - \mathbf{U}_{\mathbf{Y}}  \mathbf{U}^{T}_{\mathbf{Y}}  )  \mathbf{Z} ( \mathbf{V}_{\mathbf{Y}}  \mathbf{V}^{T}_{\mathbf{Y}} ) \ ,  \nonumber \\
  \mathbf{P}_{ \mathcal{N}^{\mathcal{F}}_{\mathbf{Y}} \mathbb{R}^{p \times n}_{\le k} } ( \mathbf{Z} ) & =  \mathbf{Z} - \mathbf{P}_{ \mathcal{T}_{\mathbf{Y}} \mathbb{R}^{p \times n}_{k} } ( \mathbf{Z} )  = \mathbf{P}_{ \mathcal{N}_{\mathbf{Y}} \mathbb{R}^{p \times n}_{k} } ( \mathbf{Z} )  \nonumber \\
     & =  ( \mathbf{I}_{p}  - \mathbf{U}_{\mathbf{Y}}  \mathbf{U}^{T}_{\mathbf{Y}} )  \mathbf{Z}   ( \mathbf{I}_{n}  - \mathbf{V}_{\mathbf{Y}}  \mathbf{V}^{T}_{\mathbf{Y}} ) \ , \nonumber
\end{align}
where the columns of $\mathbf{U}_{\mathbf{Y}}$ and $\mathbf{V}_{\mathbf{Y}}$ are, respectively, the leading $k$ left and right singular vectors of $\mathbf{Y}$, which is of rank $k$.

These results are further consistent with the more general result that, when $\mathcal{M}$ is an arbitrary submanifold embedded in $\mathbb{R}^{p \times n}$ or $\mathbb{R}^{p}$, its tangent and normal spaces at an arbitrary $\mathbf{Z}  \in \mathcal{M}$ coincide exactly with the Bouligand tangent and Frechet normal cones to $\mathcal{M}$ at $\mathbf{Z}$, see Example 6.8 of Rockafellar and Wets~\cite{RW1998} or Theorem 3.15 of Ruszczynski~\cite{R2006} for details.

Furthermore, from the above results, we see that $\bar{\mathbf{Y}}  \in \mathbb{R}^{p \times n}_{k}$ is a Frechet first-order stationary point for the WLRA problem in its formulation~\eqref{eq:P0} if it satisfies one of the following equivalent conditions:
\begin{align*}
(1)& \langle \nabla \varphi( \bar{\mathbf{Y}} ) , \mathbf{Z}  \rangle_{F}   \ge 0  \  ,   \forall  \mathbf{Z} \in   \big \lbrace \mathbf{U}_{\bar{\mathbf{Y}}} \mathbf{M} \mathbf{V}^{T}_{\bar{\mathbf{Y}}} + \mathbf{U} \mathbf{V}^{T}_{\bar{\mathbf{Y}}} + \mathbf{U}_{\bar{\mathbf{Y}}}  \mathbf{V}^{T} /  \mathbf{M}  \in  \mathbb{R}^{k \times k}, \mathbf{U}  \in  \mathbb{R}^{p \times k}, \mathbf{V}  \in  \mathbb{R}^{n \times k}   \big \rbrace  \ , \\
& \\
(2)&  \mathbf{U}^{T}_{\bar{\mathbf{Y}}}  \nabla \varphi( \bar{\mathbf{Y}} ) =  \mathbf{0}^{k \times n}  \text{ and }  \nabla \varphi( \bar{\mathbf{Y}} )  \mathbf{V}_{\bar{\mathbf{Y}}} =  \mathbf{0}^{p \times k}  \ , \\
& \\
(3)&  \nabla \varphi( \bar{\mathbf{Y}} )  ( \mathbf{V}_{\bar{\mathbf{Y}}}  \mathbf{V}^{T}_{\bar{\mathbf{Y}}} ) +( \mathbf{U}_{\bar{\mathbf{Y}}}  \mathbf{U}^{T}_{\bar{\mathbf{Y}}}  ) \nabla \varphi( \bar{\mathbf{Y}} )  (\mathbf{I}_{n} - \mathbf{V}_{\bar{\mathbf{Y}}}  \mathbf{V}^{T}_{\bar{\mathbf{Y}}}  ) = \mathbf{0}^{p \times n}  \ ,
\end{align*}
where  the columns of $\mathbf{U}_{\bar{\mathbf{Y}}}$ and $\mathbf{V}_{\bar{\mathbf{Y}}}$ are, respectively, the leading $k$ left and right singular vectors of $\bar{\mathbf{Y}}$, which is of rank $k$.

Obviously, the second condition is the more convenient for our purpose and, as $\nabla \varphi( \bar{\mathbf{Y}} ) = \mathbf{W} \odot   ( \bar{\mathbf{Y}} - \mathbf{X} )$ according to equation~\eqref{eq:D_grad_varphi}, it translates to the simple statement
\begin{equation} \label{eq:D_first_kkt_frechet3}
\mathbf{U}^{T}_{\bar{\mathbf{Y}}} \big ( \mathbf{W} \odot   ( \bar{\mathbf{Y}} - \mathbf{X} ) \big )=  \mathbf{0}^{k \times n}  \text{ and }  \big ( \mathbf{W} \odot   ( \bar{\mathbf{Y}} - \mathbf{X} )  \big ) \mathbf{V}_{\bar{\mathbf{Y}}} =  \mathbf{0}^{p \times k}  \ .
\end{equation}

We now consider the case where $\bar{\mathbf{Y}}  \in \mathbb{R}^{p \times n}_{<k}$, e.g., when $\emph{rank}( \bar{\mathbf{Y}} ) = s < k$. In that case, we deduce from Theorem~\ref{theo3.4:box}  that the structure of  $\mathcal{T}^{\mathcal{B}}_{\bar{\mathbf{Y}}} \mathbb{R}^{p \times n}_{\le k}$ is more complex as it contains $\mathcal{T}_{\bar{\mathbf{Y}}} \mathbb{R}^{p \times n}_{s}$, but also matrices of rank less or equal to $k-s$ which intersect orthogonally $\mathbb{R}^{p \times n}_{\le k}$ (with respect the Frobenius inner product) and also sum of elements belonging to each of these two sets.

A key-remark for deriving a simple condition of Frechet first-order stationarity for $\varphi(.)$ at a point $\bar{\mathbf{Y}}  \in \mathbb{R}^{p \times n}_{<k}$ is the following. Assume that  $\emph{rank}( \bar{\mathbf{Y}} ) = s < k$ and consider an arbitrary matrix $\mathbf{Z} \in \mathbb{R}^{p \times n}_{\le k-s}$. We have
\begin{equation*}
\mathbf{Z} =  \mathbf{P}_{ \mathcal{T}_{\bar{\mathbf{Y}}} \mathbb{R}^{p \times n}_{s} } ( \mathbf{Z} ) + \mathbf{P}_{ \mathcal{N}_{\bar{\mathbf{Y}}} \mathbb{R}^{p \times n}_{s} } ( \mathbf{Z} )  \ ,
\end{equation*}
since  $\mathbf{P}_{ \mathcal{T}_{\bar{\mathbf{Y}}} \mathbb{R}^{p \times n}_{s} } (.)$ and $\mathbf{P}_{ \mathcal{N}_{\bar{\mathbf{Y}}} \mathbb{R}^{p \times n}_{s} } (.)$ are two complementary orthogonal projectors with respect to the Frobenius inner product in $\mathbb{R}^{p \times n}$. Clearly, by definition, $\mathbf{P}_{ \mathcal{T}_{\bar{\mathbf{Y}}} \mathbb{R}^{p \times n}_{s} } ( \mathbf{Z} ) \in  \mathcal{T}_{\bar{\mathbf{Y}}} \mathbb{R}^{p \times n}_{s}$ and  $\mathbf{P}_{ \mathcal{N}_{\bar{\mathbf{Y}}} \mathbb{R}^{p \times n}_{s} } ( \mathbf{Z} ) \in  \mathcal{N}_{\bar{\mathbf{Y}}} \mathbb{R}^{p \times n}_{s}$. Furthermore, as the orthogonal projector $\mathbf{P}_{ \mathcal{N}_{\bar{\mathbf{Y}}} \mathbb{R}^{p \times n}_{s} } (.)$ never increases the rank of a matrix, we also have $\mathbf{P}_{ \mathcal{N}_{\bar{\mathbf{Y}}} \mathbb{R}^{p \times n}_{s} } ( \mathbf{Z} )   \in \mathbb{R}^{p \times n}_{\le k-s}$ as $\mathbf{Z} \in \mathbb{R}^{p \times n}_{\le k-s}$, and we conclude that 
\begin{equation*}
\mathbf{P}_{ \mathcal{N}_{\bar{\mathbf{Y}}} \mathbb{R}^{p \times n}_{s} } ( \mathbf{Z} )   \in \mathcal{N}_{\bar{\mathbf{Y}}} \mathbb{R}^{p \times n}_{s} \cap \mathbb{R}^{p \times n}_{\le k-s}  \ ,
\end{equation*}
which implies finally that $\mathbf{Z} \in  \mathcal{T}^{\mathcal{B}}_{\bar{\mathbf{Y}}} \mathbb{R}^{p \times n}_{\le k}$. In other words, we have the inclusion $\mathbb{R}^{p \times n}_{\le k-s} \subset  \mathcal{T}^{\mathcal{B}}_{\bar{\mathbf{Y}}} \mathbb{R}^{p \times n}_{\le k}$. From this relationship and Theorem~\ref{theo3.4:box}, it is not difficult to see that an equivalent formulation of $\mathcal{T}^{\mathcal{B}}_{\bar{\mathbf{Y}}} \mathbb{R}^{p \times n}_{\le k}$ is
\begin{equation*}
\mathcal{T}^{\mathcal{B}}_{\bar{\mathbf{Y}}} \mathbb{R}^{p \times n}_{\le k} = \mathcal{T}_{\bar{\mathbf{Y}}} \mathbb{R}^{p \times n}_{s} + \mathbb{R}^{p \times n}_{\le k-s}  \ ,
\end{equation*}
where the direct orthogonal sum $\oplus$ is now replaced by an ordinary sum, see Hosseini et al.~\cite{HLU2019} for more details.

Now, if $\mathbf{Z} \in \mathbb{R}^{p \times n}_{\le k-s}$, $-\mathbf{Z}$ also belongs to $\mathbb{R}^{p \times n}_{\le k-s}$ and, thus, any element of $\mathcal{N}^{\mathcal{F}}_{\bar{\mathbf{Y}}} \mathbb{R}^{p \times n}_{\le k}$ needs to be orthogonal to $\mathbf{Z},  \forall \mathbf{Z} \in \mathbb{R}^{p \times n}_{\le k-s}$. Next, if $\mathbf{T} \in \mathcal{N}^{\mathcal{F}}_{\bar{\mathbf{Y}}} \mathbb{R}^{p \times n}_{\le k}$ and $\mathbf{T} \ne \mathbf{0}^{p \times n}$, this implies that $\mathbf{T}$ must be orthogonal (with respect to the Frobenius inner product in $\mathbb{R}^{p \times n}$) to its best approximation of rank $k-s$ given by the Eckart-Young Theorem~\ref{theo2.1:box}, which is absurd, and we conclude that
\begin{equation*}
\mathcal{N}^{\mathcal{F}}_{\bar{\mathbf{Y}}} \mathbb{R}^{p \times n}_{\le k} = \lbrace  \mathbf{0}^{p \times n} \rbrace \ .
\end{equation*}
In this condition, if $\bar{\mathbf{Y}}  \in \mathbb{R}^{p \times n}_{<k}$, the Frechet first-order stationary condition for $\varphi(.)$ at $\bar{\mathbf{Y}}  \in \mathbb{R}^{p \times n}_{<k}$ reduces to
\begin{equation*}
\nabla \varphi( \bar{\mathbf{Y}} ) = \mathbf{0}^{p \times n} \ ,
\end{equation*}
which translates to the simple matrix equality $ \mathbf{W} \odot   ( \bar{\mathbf{Y}} - \mathbf{X} ) = \mathbf{0}^{p \times n}$, using equation~\eqref{eq:D_grad_varphi}.

Collecting all the above developments, we have demonstrated  the following theorem, which is used without proof in Ha et al.~\cite{HLB2020} in a slightly larger setting where  $\varphi(.)$ is a continuously differentiable function instead of the objective function associated with the formulation~\eqref{eq:P0} of the WLRA problem.
\begin{theo3.5} \label{theo3.5:box}
Let $\bar{\mathbf{Y}} \in \mathbb{R}^{p \times n}_{\le k}$, with $\emph{rank}( \bar{\mathbf{Y}} ) = s \le k$. Then $\bar{\mathbf{Y}}$ is a Frechet first-order stationary point for $\varphi(.)$ if
\begin{equation*}
\mathbf{U}^{T}_{\bar{\mathbf{Y}}} \big ( \mathbf{W} \odot   ( \bar{\mathbf{Y}} - \mathbf{X} ) \big )=  \mathbf{0}^{k \times n}  \text{ and }  \big ( \mathbf{W} \odot   ( \bar{\mathbf{Y}} - \mathbf{X} )  \big ) \mathbf{V}_{\bar{\mathbf{Y}}} =  \mathbf{0}^{p \times k} \ ,
\end{equation*}
when $\emph{rank}( \bar{\mathbf{Y}} ) = s = k$, or if
\begin{equation*}
\mathbf{W} \odot   ( \bar{\mathbf{Y}} - \mathbf{X} ) = \mathbf{0}^{p \times n} \ ,
\end{equation*}
when $\emph{rank}( \bar{\mathbf{Y}} ) = s < k$ and the columns of $\mathbf{U}_{\bar{\mathbf{Y}}}$ and $\mathbf{V}_{\bar{\mathbf{Y}}}$ are, respectively, the leading $s$ left and right singular vectors of $\bar{\mathbf{Y}}$, which is of rank $s$.
\\
$\Box$
\end{theo3.5}

Any local minimizer $\bar{\mathbf{Y}}$ of $\varphi(.)$ in the set $\mathbb{R}^{p \times n}_{\le k}$ must satisfy the first-order conditions stated in Theorem~\ref{theo3.5:box}, though these conditions are not sufficient in general, see Theorem 6.12 in Rockafellar and Wets~\cite{RW1998} and also Ha et al.~\cite{HLB2020} for more details. However, in the case where $\emph{rank}( \bar{\mathbf{Y}} ) < k$ and $\mathbf{W} \odot   ( \bar{\mathbf{Y}} - \mathbf{X} ) = \mathbf{0}^{p \times n}$, we deduce immediately that 
\begin{equation*}
\sqrt{ \mathbf{W}}  \odot   ( \bar{\mathbf{Y}} - \mathbf{X} ) = \mathbf{0}^{p \times n}
\end{equation*}
and $\bar{\mathbf{Y}}$ is obviously a global minimizer of  $\varphi(.)$ and a solution of the WLRA problem in this particular case.

We now derive a more convenient expression than the one given in Theorem~\ref{theo2.9:box} to verify that a matrix $\bar{\mathbf{Y}} \in \mathbb{R}^{p \times n}_{\le k}$ is a Frechet second-order stationarity point of $\varphi(.)$ over  $\mathbb{R}^{p \times n}_{\le k}$. First, we observe   that, in the case where $\bar{\mathbf{Y}} \in \mathbb{R}^{p \times n}_{<k}$ is a Frechet first-order stationarity point of $\varphi(.)$,   $\nabla \varphi( \bar{\mathbf{Y}} ) = \mathbf{0}^{p \times n}$ and $\bar{\mathbf{Y}} \in \mathbb{R}^{p \times n}_{<k}$ is a global minimum of $\varphi(.)$ over $\mathbb{R}^{p \times n}$ and, thus, $\big ( \nabla^{2} \varphi( \bar{\mathbf{Y}} ) \big )$ is a positive semi-definite quadratic form over $\mathbb{R}^{p \times n}$ (note, alternatively, that $\big ( \nabla^{2} \varphi( \bar{\mathbf{Y}} ) \big )$ is always positive semi-definite according to equation~\eqref{eq:D_hess_varphi} ). From these results, when  $\emph{rank}( \bar{\mathbf{Y}} ) < k$ is a Frechet first-order stationarity point of $\varphi(.)$, we deduce immediately that the condition~\eqref{eq:D_second_kkt_c} in Theorem~\ref{theo2.9:box} is verified and consistently $\bar{\mathbf{Y}}$ is also a Frechet second-order stationarity point of  $\varphi(.)$ in the sense of Theorem~\ref{theo2.9:box}.

Next, in the case where  $\bar{\mathbf{Y}} \in \mathbb{R}^{p \times n}_{k}$, we first recall from equations~\eqref{eq:D_riemann_bouligand_equi} that
\begin{equation*}
\mathcal{T}^{\mathcal{B}}_{\bar{\mathbf{Y}}} \mathbb{R}^{p \times n}_{\le k} = \mathcal{T}_{\bar{\mathbf{Y}}} \mathbb{R}^{p \times n}_{k}  \    \text{ and }   \  \mathcal{N}^{\mathcal{F}}_{\bar{\mathbf{Y}}} \mathbb{R}^{p \times n}_{\le k} = \mathcal{N}_{\bar{\mathbf{Y}}} \mathbb{R}^{p \times n}_{k} 
\end{equation*}
and the Frechet first-order condition is thus equivalent to
\begin{equation*}
\nabla \varphi( \bar{\mathbf{Y}} )  \in \mathcal{N}_{\bar{\mathbf{Y}}} \mathbb{R}^{p \times n}_{k} \ ,
\end{equation*}
which is exactly similar to the statement that the Riemannian gradient of $\varphi(.)$ at $\bar{\mathbf{Y}} \in \mathbb{R}^{p \times n}_{k}$ is equal to zero by equation~\eqref{eq:D_first_kkt_m}. In other words, in the case  $\emph{rank}( \bar{\mathbf{Y}} ) = k$, the Frechet first-order condition for $\varphi(.)$,  considered as a function defined on $\mathbb{R}^{p \times n}_{\le k}$, at $\bar{\mathbf{Y}}$ stated in equation~\eqref{eq:D_first_kkt_frechet3} is equivalent to the Riemannian first-order condition for the restriction of $\varphi(.)$ over the embeddded smooth submanifold $\mathbb{R}^{p \times n}_{k}$ at $\bar{\mathbf{Y}}$ stated in equation~\eqref{eq:D_first_kkt_m}.

Furthermore, when $\emph{rank}( \bar{\mathbf{Y}} ) = k$ and $\bar{\mathbf{Y}}$ is a Frechet first-order stationary point, equation~\eqref{eq:D_second_kkt_c} in Theorem~\eqref{theo2.9:box},  specialized to the case  of $\mathbb{R}^{p \times n}_{\le k}$, simplifies to
\begin{equation*} 
 \langle \nabla \varphi( \bar{\mathbf{Y}} )  , \mathbf{Z}  \rangle_{F} +  \langle  \big \lbrack \nabla^{2} \varphi( \bar{\mathbf{Y}} )  \big \rbrack ( \mathbf{D} ) , \mathbf{D}  \rangle_{F} \ge 0   \ ,  \    \forall \mathbf{D}  \in  \mathcal{T}_{\bar{\mathbf{Y}}} \mathbb{R}^{p \times n}_{k} ,   \forall \mathbf{Z}   \in \mathcal{T}_{ (\bar{\mathbf{Y}}, \mathbf{D}) } \mathbb{R}^{p \times n}_{k} \ .
\end{equation*}
This condition is strictly equivalent to the statement that the Riemannian Hessian of $\varphi(.)$ at $\bar{\mathbf{Y}}$, $\big ( \nabla^{2}_{R}   \varphi( \bar{\mathbf{Y}} ) \big )$, is positive semi-definite over $\mathcal{T}_{\bar{\mathbf{Y}}} \mathbb{R}^{p \times n}_{k}$, as noted in~\cite{YZS2014} and~\cite{L2020}. Next, using the explicit formulation (in terms of standard Euclidean derivatives) of this Riemannian Hessian of the smooth function $\varphi(.)$ defined on the smooth submanifold $\mathbb{R}^{p \times n}_{k}$, derived in Proposition 2.2 of~\cite{V2013} and Proposition 2 of~\cite{LLZ2024}, the statement that $\big ( \nabla^{2}_{R}   \varphi( \bar{\mathbf{Y}} ) \big )$ is positive semi-definite is equivalent to
\begin{equation} \label{eq:D_rhess_varphi}
\big ( \nabla^{2}_{R}   \varphi( \bar{\mathbf{Y}} ) \big ) \lbrack \mathbf{D}, \mathbf{D}  \rbrack = \big ( \nabla^{2}  \varphi( \bar{\mathbf{Y}} ) \big ) \lbrack \mathbf{D}, \mathbf{D}  \rbrack +  2.\langle  \nabla \varphi( \bar{\mathbf{Y}} ) ,    \mathbf{U}^{\bot}_{\bar{\mathbf{Y}}}  \mathbf{C}   \Sigma^{-1}_{\bar{\mathbf{Y}}}  \mathbf{B}  (\mathbf{V}^{\bot}_{\bar{\mathbf{Y}}}  )^{T}  \rangle_{F} \ge 0 \ ,
\end{equation}
$\forall \mathbf{D}  \in \mathcal{T}_{\bar{\mathbf{Y}}} \mathbb{R}^{p \times n}_{k}$, and where the thin SVD of $\bar{\mathbf{Y}} \in \mathbb{R}^{p \times n}_{k}$ is given by $\bar{\mathbf{Y}} =  \mathbf{U}_{\bar{\mathbf{Y}}}  \Sigma_{\bar{\mathbf{Y}}} \mathbf{V}^{T}_{\bar{\mathbf{Y}}}$ with $\mathbf{U}^{T}_{\bar{\mathbf{Y}}}  \mathbf{U}_{\bar{\mathbf{Y}}}   = \mathbf{V}^{T}_{\bar{\mathbf{Y}}}  \mathbf{V}_{\bar{\mathbf{Y}}}   = \mathbf{I}_{k}$ and $\Sigma_{\bar{\mathbf{Y}}}$ is a $k \times k$ diagonal matrix with strictly positive diagonal elements  (e.g., the singular values of $\bar{\mathbf{Y}}$) and 
\begin{equation*}
\mathbf{D} =  \lbrack  \mathbf{U}_{\bar{\mathbf{Y}}}    \mathbf{U}^{\bot}_{\bar{\mathbf{Y}}}  \rbrack 
 \begin{bmatrix}
\mathbf{A} & \mathbf{B} \\
\mathbf{C} & \mathbf{0}^{(p-k) \times (n-k)}  
\end{bmatrix}   
\lbrack  \mathbf{V}_{\bar{\mathbf{Y}}}    \mathbf{V}^{\bot}_{\bar{\mathbf{Y}}}  \rbrack^{T} \in \mathcal{T}_{\bar{\mathbf{Y}}} \mathbb{R}^{p \times n}_{k} \ ,
\end{equation*}
where  $\mathbf{A}  \in  \mathbb{R}^{k \times k}, \mathbf{B}  \in  \mathbb{R}^{k \times (n-k)} , \mathbf{C}  \in  \mathbb{R}^{(p-k) \times k}$ and $\lbrack  \mathbf{U}_{\bar{\mathbf{Y}}}     \mathbf{U}^{\bot}_{\bar{\mathbf{Y}}}   \rbrack$ and  $\lbrack  \mathbf{V}_{\bar{\mathbf{Y}}}     \mathbf{V}^{\bot}_{\bar{\mathbf{Y}}}   \rbrack$ are, respectively, $p \times p$ and $n \times n$ orthogonal matrices.

Using equations~\eqref{eq:D_grad_varphi} and~\eqref{eq:D_hess_varphi}, the previous discussion leads to the following theorem, which characterizes more explicitly the Frechet second-order stationarity points of $\varphi(.)$.
\begin{theo3.6} \label{theo3.6:box}
Let $\bar{\mathbf{Y}} \in \mathbb{R}^{p \times n}_{\le k}$, with $\emph{rank}( \bar{\mathbf{Y}} ) \le k$. Then $\bar{\mathbf{Y}}$ is a Frechet second-order stationary point for $\varphi(.)$ if it is a Frechet first-order stationary point for $\varphi(.)$ and if, in addition, in the case of  $\emph{rank}( \bar{\mathbf{Y}} ) = k$, if
\begin{equation} \label{eq:D_rhess_varphi2}
\Vert \sqrt{\mathbf{W}}  \odot  \mathbf{D} \Vert^{2}_{F}  \ge  -2.\langle   \mathbf{W} \odot   (\bar{\mathbf{Y}} - \mathbf{X} ) ,    \mathbf{U}^{\bot}_{\bar{\mathbf{Y}}}  \mathbf{C}   \Sigma^{-1}_{\bar{\mathbf{Y}}}  \mathbf{B}  (\mathbf{V}^{\bot}_{\bar{\mathbf{Y}}}  )^{T}  \rangle_{F}  \ ,  \ \forall \mathbf{D}  \in \mathcal{T}_{\bar{\mathbf{Y}}} \mathbb{R}^{p \times n}_{k} \ ,
\end{equation}
 where the thin SVD of $\bar{\mathbf{Y}} \in \mathbb{R}^{p \times n}_{k}$ is given by $\bar{\mathbf{Y}} =  \mathbf{U}_{\bar{\mathbf{Y}}}  \Sigma_{\bar{\mathbf{Y}}} \mathbf{V}^{T}_{\bar{\mathbf{Y}}}$,  the columns of $\mathbf{U}^{\bot}_{\bar{\mathbf{Y}}} \in \mathbb{R}^{p \times (p-k) }$ and  $\mathbf{V}^{\bot}_{\bar{\mathbf{Y}}}  \in \mathbb{R}^{n \times (n-k) }$ form, respectively, orthonormal bases of  $\emph{ran}( \bar{\mathbf{Y}} )^{\bot}$ and $\emph{ran}( \bar{\mathbf{Y}}^{T} )^{\bot}$ and
\begin{equation*}
\mathbf{D} =  \lbrack  \mathbf{U}_{\bar{\mathbf{Y}}}    \mathbf{U}^{\bot}_{\bar{\mathbf{Y}}}  \rbrack 
 \begin{bmatrix}
\mathbf{A} & \mathbf{B} \\
\mathbf{C} & \mathbf{0}^{(p-k) \times (n-k)}  
\end{bmatrix}   
\lbrack  \mathbf{V}_{\bar{\mathbf{Y}}}    \mathbf{V}^{\bot}_{\bar{\mathbf{Y}}}  \rbrack^{T} \in \mathcal{T}_{\bar{\mathbf{Y}}} \mathbb{R}^{p \times n}_{k} \ ,
\end{equation*}
where  $\mathbf{A}  \in  \mathbb{R}^{k \times k}, \mathbf{B}  \in  \mathbb{R}^{k \times (n-k)} , \mathbf{C}  \in  \mathbb{R}^{(p-k) \times k}$.
\\
$\Box$
\end{theo3.6}
Interestingly, observe that, using equation~\eqref{eq:D_normal_space_varphi}, when $\bar{\mathbf{Y}}  \in \mathbb{R}^{p \times n}_{k}$ is a Frechet first-order stationary point for $\varphi(.)$, both $\nabla \varphi( \bar{\mathbf{Y}} ) =  \mathbf{W} \odot   (\bar{\mathbf{Y}} - \mathbf{X} )$ and $  \mathbf{U}^{\bot}_{\bar{\mathbf{Y}}}  \mathbf{C}   \Sigma^{-1}_{\bar{\mathbf{Y}}}  \mathbf{B}  (\mathbf{V}^{\bot}_{\bar{\mathbf{Y}}}  )^{T}$ are elements of $\mathcal{N}_{\bar{\mathbf{Y}}} \mathbb{R}^{p \times n}_{k}$ as $\mathbf{C}   \Sigma^{-1}_{\bar{\mathbf{Y}}}  \mathbf{B}  \in  \mathbb{R}^{(n-k) \times (n-k)}$.

We now characterize the critical points of the factorized cost function $\varphi^{*}(.)$, which is used in the formulation~\eqref{eq:P1} of the WLRA problem. $\varphi^{*}(.)$ is defined on the product space $\mathbb{R}^{p \times k} \times \mathbb{R}^{k \times n}$, which is a "standard" Euclidean linear (product) space. In other words, the gradient, Hessian and critical points of $\varphi^{*}(.)$ are defined in the usual way (see Subsection~\ref{calculus:box}) as there are no additional constraints on the matrix variables $\mathbf{A} \in \mathbb{R}^{p \times k}$ and $\mathbf{B} \in \mathbb{R}^{k \times n}$. Thus, the pair $(\mathbf{A} , \mathbf{B} )$ is a first-order stationary point of $\varphi^{*}(.)$, if and only if,
\begin{equation*}
\nabla \varphi^{*}( \mathbf{A} , \mathbf{B} ) = ( \mathbf{0}^{p \times k} , \mathbf{0}^{k \times n}) \ ,
\end{equation*}
and a second-order stationary point of $\varphi^{*}(.)$ if, in addition,
\begin{equation*}
\big ( \nabla^{2} \varphi^{*} ( \mathbf{A} , \mathbf{B} )  \big ) \big ( ( \mathbf{C} , \mathbf{D} ) ,  ( \mathbf{C} , \mathbf{D} ) \big ) \ge 0 \ , \ \forall ( \mathbf{C} , \mathbf{D} ) \in \mathbb{R}^{p \times k} \times \mathbb{R}^{k \times n} \ ,
\end{equation*}
where the second derivative (Hessian) $\big ( \nabla^{2} \varphi^{*} ( \mathbf{A} , \mathbf{B} ) \big )$ is a (symmetric) quadratic form mapping from $( \mathbb{R}^{p \times k} \times \mathbb{R}^{k \times n} ) \times ( \mathbb{R}^{p \times k} \times \mathbb{R}^{k \times n} )$ to $\mathbb{R}$.

By definition, $\varphi^{*}(.)$ is the composition of $\varphi(.)$, from $\mathbb{R}^{p \times n}$ to $\mathbb{R}$, with the bilinear mapping, from $\mathbb{R}^{p \times k} \times \mathbb{R}^{k \times n}$ to $\mathbb{R}^{p \times n}$, defined by $( \mathbf{A} , \mathbf{B} ) \longrightarrow \mathbf{A} \mathbf{B}  \ , \ \forall ( \mathbf{A} , \mathbf{B} ) \in \mathbb{R}^{p \times k} \times \mathbb{R}^{k \times n}$. Furthermore, $\varphi(.)$ and this bilinear mapping are $C^{\infty}$ differentiable. Thus, using the standard chain rule on the differential of the composition of two differentiable functions, we can easily obtain the two partial derivatives of $\varphi^{*}(.)$ since, $\forall ( \mathbf{C} , \mathbf{D} ) \in \mathbb{R}^{p \times k} \times \mathbb{R}^{k \times n}$, we have, using properties of the $\Tr(.)$ operator stated in Subsection~\ref{lin_alg:box},
\begin{align*}
\mathit{D} \varphi^{*}_{\mathbf{A}} ( \mathbf{A} , \mathbf{B} ) ( \mathbf{C} ) & = \big \langle \nabla \varphi( \mathbf{A} \mathbf{B} ) , \mathbf{C} \mathbf{B} \big \rangle_{F}  \\
                                                               & = \Tr \big(  \nabla \varphi( \mathbf{A} \mathbf{B} )^{T}  \mathbf{C} \mathbf{B} \big) \\
                                                               & = \Tr \big(   \mathbf{C} \mathbf{B}  \nabla \varphi( \mathbf{A} \mathbf{B} )^{T} \big) \\
                                                               & = \Tr \Big(   \mathbf{C} \big( \nabla \varphi( \mathbf{A} \mathbf{B} )\mathbf{B}^{T} \big)^{T} \Big) \\
                                                               & = \Tr \Big(   \big( \nabla \varphi( \mathbf{A} \mathbf{B} )\mathbf{B}^{T} \big)^{T} \mathbf{C}  \Big) \\
                                                               & =  \big \langle \nabla \varphi( \mathbf{A} \mathbf{B} )\mathbf{B}^{T} ,  \mathbf{C} \big \rangle_{F}
\end{align*}
and, similarly,
\begin{align*}
\mathit{D} \varphi^{*}_{\mathbf{B}} ( \mathbf{A} , \mathbf{B} ) ( \mathbf{D} ) & = \big \langle \nabla \varphi( \mathbf{A} \mathbf{B} ) , \mathbf{A} \mathbf{D} \big \rangle_{F}  \\
                                                               & = \Tr \big(  \nabla \varphi( \mathbf{A} \mathbf{B} )^{T}  \mathbf{A} \mathbf{D} \big) \\
                                                               & = \Tr \Big(   \big( \mathbf{A}^{T} \nabla \varphi( \mathbf{A} \mathbf{B} )  \big)^{T}  \mathbf{D}   \Big) \\
                                                               & =  \big \langle \mathbf{A}^{T} \nabla \varphi( \mathbf{A} \mathbf{B} ) ,  \mathbf{D} \big \rangle_{F} \ .
\end{align*}
Thus, by the unicity of the Frobenius gradients of the partial functions $\varphi^{*}_{\mathbf{A}} (.)$ and $\varphi^{*}_{\mathbf{B}} (.)$, and equation~\eqref{eq:D_grad_varphi}, we get
\begin{align} \label{eq:D_partial_grad_varphi*}
\nabla \varphi^{*}_{\mathbf{A}} ( \mathbf{A} , \mathbf{B} ) & = \nabla \varphi( \mathbf{A} \mathbf{B} )\mathbf{B}^{T} =  \big( \mathbf{W} \odot   ( \mathbf{A} \mathbf{B} - \mathbf{X} )  \big) \mathbf{B}^{T}  \in \mathbb{R}^{p \times k} \ ,  \nonumber \\
\nabla \varphi^{*}_{\mathbf{B}} ( \mathbf{A} , \mathbf{B} ) & = \mathbf{A}^{T} \nabla \varphi( \mathbf{A} \mathbf{B} ) = \mathbf{A}^{T}  \big( \mathbf{W} \odot   ( \mathbf{A} \mathbf{B} - \mathbf{X} )  \big)  \in \mathbb{R}^{k \times n}  \ ,
\end{align}
and, finally, we obtain the gradient of $\varphi^{*} (.)$ at any pair $( \mathbf{A} , \mathbf{B} ) \in \mathbb{R}^{p \times k} \times \mathbb{R}^{k \times n}$ as
\begin{align} \label{eq:D_grad_varphi*}
\nabla \varphi^{*} ( \mathbf{A} , \mathbf{B} ) & = \big( \nabla \varphi^{*}_{\mathbf{A}} ( \mathbf{A} , \mathbf{B} ) , \nabla \varphi^{*}_{\mathbf{B}} ( \mathbf{A} , \mathbf{B} ) \big) \nonumber \\
                                          & = \Big(   \big( \mathbf{W} \odot   ( \mathbf{A} \mathbf{B} - \mathbf{X} )  \big) \mathbf{B}^{T}  ,  \mathbf{A}^{T}  \big( \mathbf{W} \odot   ( \mathbf{A} \mathbf{B} - \mathbf{X} )  \big)  \Big) \ .
\end{align}
Consequently, the pair $( \mathbf{A} , \mathbf{B} )$ is a first-order stationary point of $\varphi^{*}(.)$ if
\begin{equation*}
 \big( \mathbf{W} \odot   ( \mathbf{A} \mathbf{B} - \mathbf{X} )  \big) \mathbf{B}^{T} =  \mathbf{0}^{p \times k}  \ \text{and} \   \mathbf{A}^{T}  \big( \mathbf{W} \odot   ( \mathbf{A} \mathbf{B} - \mathbf{X} )  \big) = \mathbf{0}^{k \times n} \ .
\end{equation*}
We now derive a convenient expression for the quadratic form $\big ( \nabla^{2} \varphi^{*} ( \mathbf{A} , \mathbf{B} ) \big )$ in order to characterize the second-order stationary points of $\varphi^{*}(.)$, which are defined by the conditions
\begin{equation*}
\nabla \varphi^{*} ( \mathbf{A} , \mathbf{B} ) = ( \mathbf{0}^{p \times k} ,  \mathbf{0}^{k \times n} )
\end{equation*}
and 
\begin{equation*}
\big ( \nabla^{2} \varphi^{*} ( \mathbf{A} , \mathbf{B} ) \big ) \big ( ( \mathbf{C} , \mathbf{D} ) ,  ( \mathbf{C} , \mathbf{D} ) \big ) \ge 0 \ , \ \forall ( \mathbf{C} , \mathbf{D} ) \in \mathbb{R}^{p \times k} \times \mathbb{R}^{k \times n} \ .
\end{equation*}
This will be useful to determine the relationships between the critical points of $\varphi(.)$ and $\varphi^{*}(.)$ in Theorem~\ref{theo3.7:box} below.

$\forall ( \mathbf{C} , \mathbf{D} ) \in \mathbb{R}^{p \times k} \times \mathbb{R}^{k \times n}$, we have by the bilinearity of $\nabla^{2} \varphi^{*} ( \mathbf{A} , \mathbf{B} )$ 
\begin{align}  \label{eq:D_quad_varphi*}
\big ( \nabla^{2} \varphi^{*} ( \mathbf{A} , \mathbf{B} ) \big ) \big ( ( \mathbf{C} , \mathbf{D} ) ,  ( \mathbf{C} , \mathbf{D} ) \big ) & = \big (  \nabla^{2} \varphi^{*} ( \mathbf{A} , \mathbf{B} ) \big ) \big ( ( \mathbf{C} , \mathbf{0}^{k \times n}  ) +  ( \mathbf{0}^{p \times k} ,  \mathbf{D} ) ,  \nonumber  \\
& \quad    \quad  \quad   \quad  \quad \quad  \quad  \quad ( \mathbf{C} , \mathbf{0}^{k \times n}  ) +  ( \mathbf{0}^{p \times k} , \mathbf{D} )\big ) \nonumber \\
& = \big ( \nabla^{2} \varphi^{*} ( \mathbf{A} , \mathbf{B} ) \big ) \big ( ( \mathbf{C} , \mathbf{0}^{k \times n}  )  ,   ( \mathbf{C} , \mathbf{0}^{k \times n}  )  \big ) \nonumber \\
& \quad + \big ( \nabla^{2} \varphi^{*} ( \mathbf{A} , \mathbf{B} ) \big ) \big (  ( \mathbf{0}^{p \times k} , \mathbf{D} )  ,    ( \mathbf{0}^{p \times k} , \mathbf{D} )  \big )  \nonumber \\
& \quad + 2. \big ( \nabla^{2} \varphi^{*} ( \mathbf{A} , \mathbf{B} ) \big ) \big (    ( \mathbf{0}^{p \times k} , \mathbf{D} )  ,  ( \mathbf{C} , \mathbf{0}^{k \times n}  ) \big ) \ .
\end{align}
The last equality resulting from the fact that $\nabla^{2} \varphi^{*} ( \mathbf{A} , \mathbf{B} )$ can also be considered as a self-adjoint (e.g., symmetric) mapping from $\mathbb{R}^{p \times k} \times \mathbb{R}^{k \times n}$ to $\mathbb{R}^{p \times k} \times \mathbb{R}^{k \times n}$ with respect to the inner product in $\mathbb{R}^{p \times k} \times \mathbb{R}^{k \times n}$ (see Subsection~\ref{calculus:box} for details).

Next, for the same reason, using the expression for $\nabla \varphi^{*}_{\mathbf{A}} ( \mathbf{A} , \mathbf{B} )$ given in equation~\eqref{eq:D_partial_grad_varphi*} and properties of the $\Tr(.)$ operator, notice that
\begin{align*}
  \big ( \nabla^{2} \varphi^{*} ( \mathbf{A} , \mathbf{B} )  \big ) \big ( ( \mathbf{C} , \mathbf{0}^{k \times n}  )  ,   ( \mathbf{C} , \mathbf{0}^{k \times n}  )  \big ) & =   \big \langle   \lbrack \nabla^{2} \varphi^{*} ( \mathbf{A} , \mathbf{B} ) \rbrack \big( ( \mathbf{C} , \mathbf{0}^{k \times n}  ) \big) , ( \mathbf{C} , \mathbf{0}^{k \times n}  )   \big \rangle_{\mathbb{R}^{p \times k} \times \mathbb{R}^{k \times n}}   \\
                    & =  \big ( \nabla^{2} \varphi^{*}_{\mathbf{A}} ( \mathbf{A} , \mathbf{B} ) \big ) ( \mathbf{C} , \mathbf{C} ) \\
                    & =  \big \langle   \lbrack \nabla^{2} \varphi^{*}_{\mathbf{A}} ( \mathbf{A} , \mathbf{B} ) \rbrack ( \mathbf{C} ) , \mathbf{C} \big \rangle_{F} \\
                    & =  \big \langle    \big( \mathbf{W} \odot   \mathbf{C} \mathbf{B}   \big) \mathbf{B}^{T}  , \mathbf{C} \big \rangle_{F} \\
                    & = \Tr \Big (   \big(  ( \mathbf{W} \odot   \mathbf{C} \mathbf{B}  )  \mathbf{B}^{T}  \big)^{T}  \mathbf{C}  \Big ) \\
                    & = \Tr \Big (    \mathbf{C}  \big(  ( \mathbf{W} \odot   \mathbf{C} \mathbf{B}  )  \mathbf{B}^{T}  \big)^{T}  \Big ) \\
                    & = \Tr \Big (    \mathbf{C}  \mathbf{B}  \big(  \mathbf{W} \odot   \mathbf{C} \mathbf{B}   \big)^{T}  \Big ) \\
                    & = \Tr \Big (   \big(  \mathbf{W} \odot   \mathbf{C} \mathbf{B}   \big)^{T}   \mathbf{C}  \mathbf{B}  \Big ) \\
                    & = \big \langle  \mathbf{W} \odot   \mathbf{C} \mathbf{B}  , \mathbf{C}  \mathbf{B} \big \rangle_{F} \\
                    & =  \big ( \nabla^{2} \varphi ( \mathbf{A} \mathbf{B} )  \big ) ( \mathbf{C}  \mathbf{B} , \mathbf{C}  \mathbf{B}  ) \ ,
\end{align*}
where the last equality results from equation~\eqref{eq:D_hess_varphi}. Similarly, we have
\begin{align*}
  \big ( \nabla^{2} \varphi^{*} ( \mathbf{A} , \mathbf{B} )  \big ) \big ( (  \mathbf{0}^{p \times k} , \mathbf{D} )  ,  (  \mathbf{0}^{p \times k} , \mathbf{D} )  \big ) & =   \big \langle   \lbrack \nabla^{2} \varphi^{*} ( \mathbf{A} , \mathbf{B} ) \rbrack  \big( (  \mathbf{0}^{p \times k} , \mathbf{D} ) \big), (  \mathbf{0}^{p \times k} , \mathbf{D} )   \big \rangle_{\mathbb{R}^{p \times k} \times \mathbb{R}^{k \times n}}   \\
                    & =  \big ( \nabla^{2} \varphi^{*}_{\mathbf{B}} ( \mathbf{A} , \mathbf{B} ) \big ) ( \mathbf{D} , \mathbf{D} ) \\
                    & =  \big \langle   \lbrack \nabla^{2} \varphi^{*}_{\mathbf{B}} ( \mathbf{A} , \mathbf{B} ) \rbrack ( \mathbf{D} ) , \mathbf{D} \big \rangle_{F} \\
                    & =  \big \langle    \mathbf{A}^{T} \big( \mathbf{W} \odot   \mathbf{A} \mathbf{D}   \big)   , \mathbf{D} \big \rangle_{F} \\
                    & = \Tr \Big (   \big(  \mathbf{A}^{T} ( \mathbf{W} \odot   \mathbf{A} \mathbf{D}  )   \big)^{T}  \mathbf{D}  \Big ) \\
                    & = \Tr \Big (   \big(  \mathbf{W} \odot   \mathbf{A} \mathbf{D}     \big)^{T}   \mathbf{A} \mathbf{D}  \Big ) \\
                    & = \big \langle  \mathbf{W} \odot   \mathbf{A} \mathbf{D}  , \mathbf{A}  \mathbf{D} \big \rangle_{F} \\
                    & =  \big ( \nabla^{2} \varphi ( \mathbf{A} \mathbf{B} )  \big ) ( \mathbf{A}  \mathbf{D} , \mathbf{A}  \mathbf{D}  ) \ .                    
\end{align*}
We now reformulate similarly the last factor in the right-hand side of equation~\eqref{eq:D_quad_varphi*} in terms of $\nabla \varphi ( \mathbf{A} \mathbf{B} )$ and $\big ( \nabla^{2} \varphi ( \mathbf{A} \mathbf{B} )  \big )$:
\begin{align*}
\big ( \nabla^{2} \varphi^{*} ( \mathbf{A} , \mathbf{B} ) \big ) \big (    ( \mathbf{0}^{p \times k} , \mathbf{D} )  ,  ( \mathbf{C} , \mathbf{0}^{k \times n}  ) \big )  & =  \Big \langle   \mathit{D}_{\mathbf{B}} \big ( \nabla \varphi^{*}_{\mathbf{A}} ( \mathbf{A} , \mathbf{B} ) \big ) ( \mathbf{D} ),  \mathbf{C} \Big \rangle_{F} \\
       & =  \Big \langle     \mathit{D}_{\mathbf{B}} \big (  (  \mathbf{W} \odot   ( \mathbf{A} \mathbf{B} - \mathbf{X} ) ) \mathbf{B}^{T} \big ) ( \mathbf{D} )  ,  \mathbf{C}  \Big \rangle_{F} \\
       & =  \Big \langle   \big (  \mathbf{W} \odot   ( \mathbf{A} \mathbf{B} - \mathbf{X} ) \big )  \mathbf{D}^{T}  +   (  \mathbf{W} \odot  \mathbf{A} \mathbf{D}  ) \mathbf{B}^{T}  ,  \mathbf{C}  \Big \rangle_{F}  \\
       & =  \Big \langle   \big (  \mathbf{W} \odot   ( \mathbf{A} \mathbf{B} - \mathbf{X} ) \big )  \mathbf{D}^{T}   ,  \mathbf{C}  \Big \rangle_{F} \\
       & \quad +  \Big \langle (  \mathbf{W} \odot  \mathbf{A} \mathbf{D}  ) \mathbf{B}^{T}  ,  \mathbf{C}  \Big \rangle_{F}  \\
       & = \Tr \Big ( \mathbf{D}  \big (  \mathbf{W} \odot   ( \mathbf{A} \mathbf{B} - \mathbf{X} ) \big )^{T}  \mathbf{C} \Big ) +  \Tr \Big (  \mathbf{B}  (  \mathbf{W} \odot  \mathbf{A} \mathbf{D}  )^{T}   \mathbf{C} \Big ) \\
       & = \Tr \Big ( \mathbf{C}  \mathbf{D}  \big (  \mathbf{W} \odot   ( \mathbf{A} \mathbf{B} - \mathbf{X} ) \big )^{T}  \Big ) + \Tr \Big (  \mathbf{C}  \mathbf{B}  (  \mathbf{W} \odot  \mathbf{A} \mathbf{D}  )^{T}  \Big ) \\
       & = \big \langle  \mathbf{W} \odot   ( \mathbf{A} \mathbf{B} - \mathbf{X} ) , \mathbf{C}  \mathbf{D}  \big \rangle_{F} +  \big \langle     \mathbf{W} \odot  \mathbf{A} \mathbf{D} ,  \mathbf{C}  \mathbf{B} \big \rangle_{F} \\
       & = \big \langle \nabla \varphi ( \mathbf{A} \mathbf{B} ) , \mathbf{C}  \mathbf{D} \big \rangle_{F} + \big ( \nabla^{2} \varphi ( \mathbf{A} \mathbf{B} )  \big ) ( \mathbf{A} \mathbf{D}  , \mathbf{C}  \mathbf{B} ) \ .
\end{align*} 
Summarizing the preceding results, we have
\begin{align}  \label{eq:D_partial_hess_varphi*}
 \big ( \nabla^{2} \varphi^{*} ( \mathbf{A} , \mathbf{B} )  \big ) \big ( ( \mathbf{C} , \mathbf{0}^{k \times n}  )  ,   ( \mathbf{C} , \mathbf{0}^{k \times n}  )  \big ) & =  \big ( \nabla^{2} \varphi^{*}_{\mathbf{A}} ( \mathbf{A} , \mathbf{B} ) \big ) ( \mathbf{C} , \mathbf{C} ) = \big ( \nabla^{2} \varphi ( \mathbf{A} \mathbf{B} )  \big ) ( \mathbf{C}  \mathbf{B} , \mathbf{C}  \mathbf{B}  )  ,  \nonumber  \\
  \big ( \nabla^{2} \varphi^{*} ( \mathbf{A} , \mathbf{B} )  \big ) \big ( (  \mathbf{0}^{p \times k} , \mathbf{D} )  ,  (  \mathbf{0}^{p \times k} , \mathbf{D} )  \big ) & = \big ( \nabla^{2} \varphi^{*}_{\mathbf{B}} ( \mathbf{A} , \mathbf{B} ) \big ) ( \mathbf{D} , \mathbf{D} )  = \big ( \nabla^{2} \varphi ( \mathbf{A} \mathbf{B} )  \big ) ( \mathbf{A}  \mathbf{D} , \mathbf{A}  \mathbf{D}  )  ,  \nonumber  \\
\big ( \nabla^{2} \varphi^{*} ( \mathbf{A} , \mathbf{B} ) \big ) \big (    ( \mathbf{0}^{p \times k} , \mathbf{D} )  ,  ( \mathbf{C} , \mathbf{0}^{k \times n}  ) \big )  & =   \big \langle \nabla \varphi ( \mathbf{A} \mathbf{B} ) , \mathbf{C}  \mathbf{D} \big \rangle_{F} + \big ( \nabla^{2} \varphi ( \mathbf{A} \mathbf{B} )  \big ) ( \mathbf{A} \mathbf{D}  , \mathbf{C}  \mathbf{B} ) \ ,
\end{align}
and this implies, finally, using the symmetry and bilinearity of the bilinear form $\big ( \nabla^{2} \varphi^{*} ( \mathbf{A} , \mathbf{B} ) \big )$ that
\begin{align} \label{eq:D_hess_varphi*}
\big ( \nabla^{2} \varphi^{*} ( \mathbf{A} , \mathbf{B} ) \big ) \big ( ( \mathbf{C} , \mathbf{D} ) ,  ( \mathbf{C} , \mathbf{D} ) \big )  & =    \big ( \nabla^{2} \varphi ( \mathbf{A} \mathbf{B} )  \big ) ( \mathbf{C}  \mathbf{B} , \mathbf{C}  \mathbf{B}  )  +  \big ( \nabla^{2} \varphi ( \mathbf{A} \mathbf{B} )  \big ) ( \mathbf{A}  \mathbf{D} , \mathbf{A}  \mathbf{D}  )  \nonumber \\
&  \quad + 2. \big ( \nabla^{2} \varphi ( \mathbf{A} \mathbf{B} )  \big ) ( \mathbf{A} \mathbf{D}  , \mathbf{C}  \mathbf{B} ) + 2.\big \langle \nabla \varphi ( \mathbf{A} \mathbf{B} ) , \mathbf{C}  \mathbf{D} \big \rangle_{F}  \nonumber \\
& = \big ( \nabla^{2} \varphi ( \mathbf{A} \mathbf{B} )  \big ) ( \mathbf{C}  \mathbf{B} , \mathbf{C}  \mathbf{B} +   \mathbf{A}  \mathbf{D} ) \nonumber \\
& \quad +  \big ( \nabla^{2} \varphi ( \mathbf{A} \mathbf{B} )  \big ) ( \mathbf{A}  \mathbf{D} , \mathbf{A}  \mathbf{D} + \mathbf{C}  \mathbf{B}  )  \nonumber \\
& \quad + 2.\big \langle \nabla \varphi ( \mathbf{A} \mathbf{B} ) , \mathbf{C}  \mathbf{D} \big \rangle_{F} \nonumber \\
& =  \big ( \nabla^{2} \varphi ( \mathbf{A} \mathbf{B} )  \big ) ( \mathbf{C}  \mathbf{B} +   \mathbf{A}  \mathbf{D} , \mathbf{C}  \mathbf{B} +   \mathbf{A}  \mathbf{D} ) \nonumber \\
& \quad + 2.\big \langle \nabla \varphi ( \mathbf{A} \mathbf{B} ) , \mathbf{C}  \mathbf{D} \big \rangle_{F} \ .
\end{align}
Thus, the second-order stationary condition for $\varphi^{*}(.)$ at $( \mathbf{A} , \mathbf{B} )$, e.g., that the quadratic form $\big ( \nabla^{2} \varphi^{*} ( \mathbf{A} , \mathbf{B} ) \big )$  is positive semi-definite, is equivalent to the inequality
\begin{equation*}
\big ( \nabla^{2} \varphi ( \mathbf{A} \mathbf{B} )  \big ) ( \mathbf{C}  \mathbf{B} +   \mathbf{A}  \mathbf{D} , \mathbf{C}  \mathbf{B} +   \mathbf{A}  \mathbf{D} )  \ge -2.\big \langle \nabla \varphi ( \mathbf{A} \mathbf{B} ) , \mathbf{C}  \mathbf{D} \big \rangle_{F} \ ,
\end{equation*}
or, using equations~\eqref{eq:D_grad_varphi}  and~\eqref{eq:D_hess_varphi}, to the more convenient inequality
\begin{equation*}
\Vert  \sqrt{\mathbf{W}} \odot ( \mathbf{C}  \mathbf{B} +   \mathbf{A}  \mathbf{D} ) \Vert_{F}  \ge -2.\big \langle  \mathbf{W} \odot   ( \mathbf{A} \mathbf{B} - \mathbf{X} ) , \mathbf{C}  \mathbf{D} \big \rangle_{F} \ , \ \forall ( \mathbf{C} , \mathbf{D} ) \in \mathbb{R}^{p \times k} \times \mathbb{R}^{k \times n} \ .
\end{equation*}
Note that, while these inequalities based on the quadratic expression of $\nabla^{2} \varphi ( \mathbf{A} \mathbf{B} ) $ and  the gradient $\nabla \varphi ( \mathbf{A} \mathbf{B} )$ are sufficient for our purpose in this section, it is rather straightforward to obtain the general bilinear form of $\nabla^{2} \varphi^{*} ( \mathbf{A} , \mathbf{B} ) $ since
\begin{align*}
\big ( \nabla^{2} \varphi^{*} ( \mathbf{A} , \mathbf{B} ) \big ) \big ( ( \mathbf{C} , \mathbf{D} ) ,  ( \mathbf{E} , \mathbf{F} ) \big )  & = \big ( \nabla^{2} \varphi^{*} ( \mathbf{A} , \mathbf{B} ) \big ) \big ( ( \mathbf{C} + \mathbf{E} , \mathbf{D} + \mathbf{F} ) , ( \mathbf{C} + \mathbf{E} , \mathbf{D} + \mathbf{F} ) \big ) \\
& \quad -  \frac{1}{4} \big ( \nabla^{2} \varphi^{*} ( \mathbf{A} , \mathbf{B} ) \big ) \big ( ( \mathbf{C} - \mathbf{E} , \mathbf{D} - \mathbf{F} ) , ( \mathbf{C} - \mathbf{E} , \mathbf{D} - \mathbf{F} ) \big ) \ .
\end{align*}
Furthermore, a vectorized formulation of the symmetric bilinear mapping $\big ( \nabla^{2} \varphi^{*} ( \mathbf{A} , \mathbf{B} ) \big )$ will be also derived later in Section~\ref{nipals:box}.

We are now in the position to characterize more precisely the connections between the critical points of $\varphi(.)$ and $\varphi^{*}(.)$ in the following theorem, which is a reformulation and a slight extension in our WLRA context of results first given in Ha et al.~\cite{HLB2020} and later refined in Levin et al.~\cite{LKB2025} and Luo et al.~\cite{LLZ2024}.
\begin{theo3.7} \label{theo3.7:box}
Let $( \mathbf{A} , \mathbf{B} ) \in \mathbb{R}^{p \times k} \times \mathbb{R}^{k \times n}$. Then:

$(1)$ If $\mathbf{A} \mathbf{B}  \in \mathbb{R}^{p \times n}_{\le k}$ is a Frechet first-order stationary point of $\varphi(.)$ in the sense of Theorem~\ref{theo3.5:box} then $( \mathbf{A} , \mathbf{B} )$ is a first-order stationary point of $\varphi^{*}(.)$.

$(2)$ Reciprocally, if $( \mathbf{A} , \mathbf{B} )$ is a first-order stationary point of $\varphi^{*}(.)$ such that $\mathbf{A} \mathbf{B}  \in \mathbb{R}^{p \times n}_{k}$ then $\mathbf{A} \mathbf{B}$ is a Frechet first-order stationary point of $\varphi(.)$ in the sense of Theorem~\ref{theo3.5:box}.

$(3)$ Moreover, if $( \mathbf{A} , \mathbf{B} )$ is a second-order stationary point of $\varphi^{*}(.)$ such that $\mathbf{A} \mathbf{B}  \in \mathbb{R}^{p \times n}_{<k}$, then $\mathbf{A} \mathbf{B}$ is a Frechet first-order stationary point of $\varphi(.)$ in the sense of Theorem~\ref{theo3.5:box} and, thus, also a Frechet second-order  stationary point of $\varphi(.)$  and even a solution of the WLRA problem in its formulation~\eqref{eq:P0}.

$(4)$ Reciprocally, if $\mathbf{A} \mathbf{B}  \in \mathbb{R}^{p \times n}_{<k}$ is a Frechet second-order  stationary point of $\varphi(.)$ then $( \mathbf{A} , \mathbf{B} )$ is a second-order stationary point of $\varphi^{*}(.)$ and also a solution of the WLRA problem in its formulation~\eqref{eq:P1}.

$(5)$ Finally, if  $( \mathbf{A} , \mathbf{B} )$ is a second-order stationary point of $\varphi^{*}(.)$ such that $\mathbf{A} \mathbf{B}  \in \mathbb{R}^{p \times n}_{k}$, then $\mathbf{A} \mathbf{B}$ is a Frechet second-order stationary point of $\varphi(.)$ in the sense of Theorem~\ref{theo3.5:box}.

$(6)$ Reciprocally, if $\mathbf{A} \mathbf{B}  \in \mathbb{R}^{p \times n}_{k}$ is a Frechet second-order  stationary point of $\varphi(.)$ then $( \mathbf{A} , \mathbf{B} )$ is a second-order stationary point of $\varphi^{*}(.)$.

\end{theo3.7}
\begin{proof}

$(1):$ In order to prove the first assertion, we assume that $\mathbf{A} \mathbf{B}  \in \mathbb{R}^{p \times n}_{\le k}$ is a Frechet first-order stationary point of $\varphi(.)$ and we consider separately the two cases $\emph{rank}( \mathbf{A} \mathbf{B} ) < k$ and $\emph{rank}( \mathbf{A} \mathbf{B} ) = k$.

If  $\emph{rank}( \mathbf{A} \mathbf{B} ) < k$, according to Theorem~\ref{theo3.5:box}, we have $\nabla \varphi ( \mathbf{A} \mathbf{B} ) =  \mathbf{0}^{p \times n}$ and we deduce immediately that
\begin{align*}
\nabla \varphi^{*}_{\mathbf{A}} ( \mathbf{A} , \mathbf{B} ) & = \nabla \varphi( \mathbf{A} \mathbf{B} )\mathbf{B}^{T} =    \mathbf{0}^{p \times k} \ ,  \\
\nabla \varphi^{*}_{\mathbf{B}} ( \mathbf{A} , \mathbf{B} ) & = \mathbf{A}^{T} \nabla \varphi( \mathbf{A} \mathbf{B} ) =   \mathbf{0}^{k \times n} \ .
\end{align*}
In other words, the pair $( \mathbf{A} , \mathbf{B} )$ is a first-order critical point of $\varphi^{*}(.)$.

On the other hand, if $\emph{rank}( \mathbf{A} \mathbf{B} ) = k$, again according to Theorem~\ref{theo3.5:box}, we have
\begin{equation*}
\nabla \varphi ( \mathbf{A} \mathbf{B} )^{T}  \mathbf{U}_{AB} =  \mathbf{0}^{n \times k} \ \text{ and } \   \nabla \varphi ( \mathbf{A} \mathbf{B} ) \mathbf{V}_{AB} =  \mathbf{0}^{p \times k} \ ,
\end{equation*}
where the columns of $\mathbf{U}_{AB}$ and  $\mathbf{V}_{AB}$  are, respectively, the first $k$ left and right singular vectors of the matrix product $\mathbf{A} \mathbf{B}$ in its thin SVD, e.g., $\mathbf{A} \mathbf{B} =  \mathbf{U}_{AB}  \Sigma_{AB} \mathbf{V}_{AB}$.

As  $\emph{rank}( \mathbf{A} \mathbf{B} ) =  \emph{rank}( \mathbf{A}  ) =  \emph{rank}( \mathbf{B}  ) = k$, we have  $\emph{ran}( \mathbf{A} \mathbf{B} ) = \emph{ran}( \mathbf{A}  )$ and $\emph{ran}( \mathbf{B}^{T}  \mathbf{A}^{T} ) = \emph{ran}( \mathbf{B}^{T}  )$, and also
\begin{align*}
\emph{ran}( \mathbf{U}_{AB} ) & = \emph{ran}( \mathbf{A} \mathbf{B} ) = \emph{ran}( \mathbf{A}  ) ,  \\
\emph{ran}( \mathbf{V}_{AB} ) & = \emph{ran}( \mathbf{B}^{T}  \mathbf{A}^{T} ) = \emph{ran}( \mathbf{B}^{T}  ) \ .
\end{align*}
This implies that it exists $\mathbf{C} \in \mathbb{R}^{k \times k}$ and $\mathbf{D} \in \mathbb{R}^{k \times k}$ such that
\begin{equation*}
 \mathbf{A} = \mathbf{U}_{AB}\mathbf{C}  \ \text{ and } \ \mathbf{B}^{T} = \mathbf{V}_{AB}\mathbf{D}  \ .
\end{equation*}
In these conditions, we have
\begin{align*}
\nabla \varphi^{*}_{\mathbf{A}} ( \mathbf{A} , \mathbf{B} ) & = \nabla \varphi( \mathbf{A} \mathbf{B} )\mathbf{B}^{T} = \big ( \nabla \varphi( \mathbf{A} \mathbf{B} ) \mathbf{V}_{AB} \big ) \mathbf{D}  =  \mathbf{0}^{p \times k}  \ ,  \\
\nabla \varphi^{*}_{\mathbf{B}} ( \mathbf{A} , \mathbf{B} ) & = \mathbf{A}^{T} \nabla \varphi( \mathbf{A} \mathbf{B} ) =    \mathbf{C}^{T}  \big (  \mathbf{U}_{AB}^{T} \nabla \varphi( \mathbf{A} \mathbf{B} ) \big )  =  \mathbf{0}^{k \times n}  \ ,
\end{align*}
as $\mathbf{A} \mathbf{B}$ is is a first-order critical point of $\varphi(.)$. In other words, we have $\nabla \varphi^{*} ( \mathbf{A} , \mathbf{B} ) = ( \mathbf{0}^{p \times k} , \mathbf{0}^{k \times n} )$ and the pair $( \mathbf{A} , \mathbf{B} )$ is a first-order stationary point of $\varphi^{*}(.)$.

$(2):$ Reciprocally, if the pair $( \mathbf{A} , \mathbf{B} )$ is a first-order stationary point of $\varphi^{*}(.)$ such that  $\mathbf{A} \mathbf{B}  \in \mathbb{R}^{p \times n}_{k}$, we have also $\emph{rank}( \mathbf{A} \mathbf{B} ) =  \emph{rank}( \mathbf{A}  ) =  \emph{rank}( \mathbf{B}^{T}  ) = k$, which implies again that $\mathbf{U}_{AB}$ and $\mathbf{A}$ span the same column space and that their columns form two bases of $\emph{ran}( \mathbf{U}_{AB} ) =  \emph{ran}( \mathbf{A}  )$. Similarly, $\mathbf{V}_{AB}$ and $\mathbf{B}^{T}$ span the same column space and  their columns form two bases of $\emph{ran}( \mathbf{V}_{AB} ) =  \emph{ran}( \mathbf{B}^{T}  )$. In these conditions, it exist  $\mathbf{C} \in \mathbb{R}^{k \times k}$ and $\mathbf{D} \in \mathbb{R}^{k \times k}$ such that
\begin{equation*}
 \mathbf{U}_{AB}  = \mathbf{A} \mathbf{C}  \ \text{ and } \ \mathbf{V}_{AB}  = \mathbf{B}^{T} \mathbf{D}  \ .
\end{equation*}
Using the first-order optimality conditions of $( \mathbf{A} , \mathbf{B} )$ for $\varphi^{*}(.)$, we have
\begin{equation*}
\nabla \varphi( \mathbf{A} \mathbf{B} )\mathbf{B}^{T} =  \mathbf{0}^{p \times k}  \ \text{ and } \  \mathbf{A}^{T} \nabla \varphi( \mathbf{A} \mathbf{B} ) =   \mathbf{0}^{k \times n} \ ,
\end{equation*}
which implies that
\begin{align*}
\nabla \varphi( \mathbf{A} \mathbf{B} ) \mathbf{V}_{AB}        & =  \big (  \nabla \varphi( \mathbf{A} \mathbf{B} )  \mathbf{B}^{T} \big ) \mathbf{D} = \mathbf{0}^{p \times k}  \ ,  \\
\nabla \varphi( \mathbf{A} \mathbf{B} )^{T} \mathbf{U}_{AB} & = \big (  \nabla \varphi( \mathbf{A} \mathbf{B} )^{T}  \mathbf{A} \big ) \mathbf{C}   =  \big ( \mathbf{A}^{T} \nabla \varphi( \mathbf{A} \mathbf{B} ) \big )^{T} \mathbf{C}  = \mathbf{0}^{k \times n}  \ ,
\end{align*}
and the matrix product $\mathbf{A} \mathbf{B}$ is a first-order critical point of  $\varphi(.)$ in the sense of Theorem~\ref{theo3.5:box}.

$(3):$ To demonstrate the next claim of the theorem, let $\mathbf{u}_{1} \in  \mathbb{R}^{p}$, $\mathbf{v}_{1} \in  \mathbb{R}^{n}$ and $\sigma_1  \in  \mathbb{R}_{+}$ be, respectively, the first left and right singular vectors and the first singular value of $\nabla \varphi( \mathbf{A} \mathbf{B} )  \in \mathbb{R}^{p \times n}$. We first recall from equation~\eqref{eq:spect_norm2} in Subsection~\ref{lin_alg:box} that the spectral norm of  $\nabla \varphi( \mathbf{A} \mathbf{B} )$ is given by
\begin{equation*}
\Vert \nabla \varphi( \mathbf{A} \mathbf{B} ) \Vert_{S} =  \sigma_1 = \mathbf{u}_{1}^{T} \nabla \varphi( \mathbf{A} \mathbf{B} ) \mathbf{v}_{1}  \ .
\end{equation*}
Moreover, as demonstrated just before Theorem~\ref{theo3.7:box}, the hypothesis that the pair $( \mathbf{A} , \mathbf{B} )$ is a second-order stationary point of $\varphi^{*}(.)$ is equivalent to the inequality
\begin{equation*}
\Vert   \sqrt{\mathbf{W}} \odot (  \mathbf{A}  \mathbf{D} + \mathbf{C}  \mathbf{B} ) \Vert_{F}  \ge -2.\big \langle \nabla \varphi ( \mathbf{A} \mathbf{B} ) , \mathbf{C}  \mathbf{D} \big \rangle_{F} \ , \ \forall ( \mathbf{C} , \mathbf{D} ) \in \mathbb{R}^{p \times k} \times \mathbb{R}^{k \times n}  \ .
\end{equation*}
Now, suppose that $\mathbf{A} \mathbf{B}  \in \mathbb{R}^{p \times n}_{<k}$ then $\mathbf{A}$ or $\mathbf{B}$ are not of full rank since $ k \le min( p, n)$. Without loss of generality suppose that $\emph{rank}( \mathbf{A} ) < k$. By the rank-nullity theorem~\eqref{eq:rank}, this implies that it exists a unit vector $\mathbf{w} \in  \mathbb{R}^{k}$ such that $\mathbf{A} \mathbf{w} = \mathbf{0}^{p}$. Let
\begin{equation*}
( \mathbf{C}_{c} , \mathbf{D}_{c} ) = ( - \mathbf{u}_{1} \mathbf{w}^{T} , c.\mathbf{w} \mathbf{v}_{1}^{T} ) \in \mathbb{R}^{p \times k} \times \mathbb{R}^{k \times n}  \ ,  \  \forall c  \in \mathbb{R}_{+*} \  .
\end{equation*}
By hypothesis, the pair $( \mathbf{A} , \mathbf{B} )$ is a second-order stationary point of $\varphi^{*}(.)$, which implies that
\begin{equation*}
\Vert   \sqrt{\mathbf{W}} \odot ( \mathbf{A}  \mathbf{D}_{c} +  \mathbf{C}_{c}  \mathbf{B} ) \Vert_{F}  \ge -2.\big \langle   \nabla \varphi ( \mathbf{A} \mathbf{B} ) , \mathbf{C}_{c} , \mathbf{D}_{c} \big \rangle_{F} \  .
\end{equation*}
Now, we have
\begin{equation*}
 \mathbf{A}  \mathbf{D}_{c} + \mathbf{C}_{c}  \mathbf{B}  = c. \mathbf{A} \mathbf{w} \mathbf{v}_{1}^{T} -  \mathbf{u}_{1} \mathbf{w}^{T} \mathbf{B} = -  \mathbf{u}_{1} \mathbf{w}^{T} \mathbf{B} \  ,
\end{equation*}
since $\mathbf{A} \mathbf{w} = \mathbf{0}^{p}$. Furthermore, as $ \Vert \mathbf{w} \Vert^{2}_{2} =  \mathbf{w}^{T}  \mathbf{w} = 1$, $\Vert \nabla \varphi( \mathbf{A} \mathbf{B} ) \Vert_{S} = \mathbf{u}_{1}^{T} \nabla \varphi( \mathbf{A} \mathbf{B} ) \mathbf{v}_{1} =  \sigma_1$ and
\begin{equation*}
\Tr \big(  \mathbf{E} \mathbf{F} \mathbf{G} \big) = \Tr \big( \mathbf{G}  \mathbf{E} \mathbf{F}  \big) \ , \ \forall \mathbf{E} \in \mathbb{R}^{p \times n} , \mathbf{F} \in \mathbb{R}^{n \times m} ,  \mathbf{G} \in \mathbb{R}^{m \times p}  \  ,
\end{equation*}
we deduce that
\begin{align*}
\big \langle   \nabla \varphi ( \mathbf{A} \mathbf{B} ) , \mathbf{C}_{c} , \mathbf{D}_{c} \big \rangle_{F}  & = \big \langle   \nabla \varphi ( \mathbf{A} \mathbf{B} ) , -c. \mathbf{u}_{1} \mathbf{w}^{T} \mathbf{w} \mathbf{v}_{1}^{T} \big \rangle_{F}  \\
 & = -c. \big \langle   \nabla \varphi ( \mathbf{A} \mathbf{B} ) ,  \mathbf{u}_{1} \mathbf{v}_{1}^{T} \big \rangle_{F}  \\
 & = -c. \Tr \big(    \nabla \varphi ( \mathbf{A} \mathbf{B} )^{T} \mathbf{u}_{1} \mathbf{v}_{1}^{T} \big)  \\
 & = -c. \Tr \big(   \mathbf{v}_{1}^{T}  \nabla \varphi ( \mathbf{A} \mathbf{B} )^{T} \mathbf{u}_{1}  \big)  \\
 & = -c.   \mathbf{v}_{1}^{T}  \nabla \varphi ( \mathbf{A} \mathbf{B} )^{T} \mathbf{u}_{1}   \\
 & = -c.   \mathbf{u}_{1}^{T}  \nabla \varphi ( \mathbf{A} \mathbf{B} ) \mathbf{v}_{1}^{T}   \\
  & = -c.\sigma_1 \  .
\end{align*}
Using these different results, the preceding inequality simplifies to 
\begin{equation*}
\Vert  \sqrt{\mathbf{W}} \odot ( \mathbf{u}_{1} \mathbf{w}^{T} \mathbf{B} ) \Vert_{F}  \ge  2.c.\sigma_1 = 2.c.\Vert \nabla \varphi( \mathbf{A} \mathbf{B} ) \Vert_{S} \  ,
\end{equation*}
which holds for any $ c > 0 $. On the other hand, since the left-hand side of the last inequality is the Frobenius norm of a fixed element of $\mathbb{R}^{p \times n}$, which is not a function of $c$, it must be finite and this implies that $\Vert \nabla \varphi( \mathbf{A} \mathbf{B} ) \Vert_{S} = \sigma_1 = 0$, i.e.,  $\nabla \varphi( \mathbf{A} \mathbf{B} ) =  \mathbf{0}^{p \times n}$. Consequently, since $\emph{rank}( \mathbf{A}\mathbf{B} ) < k$ by hypothesis, $\mathbf{A} \mathbf{B}$  is a Frechet first-order stationary point of $\varphi(.)$ in the sense of Theorem~\ref{theo3.5:box} and  even a solution of the WLRA problem in its formulation~\eqref{eq:P0}.

$(4):$  if $\mathbf{A} \mathbf{B}  \in \mathbb{R}^{p \times n}_{<k}$ is a Frechet second-order  stationary point of $\varphi(.)$ then this pair is a fortiori a Frechet first-order  stationary point of $\varphi(.)$ and, according to Theorem~\ref{theo3.5:box}, also a solution of the WLRA problem in its formulation~\eqref{eq:P1}. By an application of Theorem~\ref{theo3.1:box}, we deduce immediately that the pair $( \mathbf{A} , \mathbf{B} )$  is a solution of the WLRA problem in its formulation~\eqref{eq:P1} and, thus, also a second-order stationary point of $\varphi^{*}(.)$.

$(5)$ and $(6):$ the proofs of these two assertions can be found in Luo et al.~\cite{LLZ2024}, especially their Corollary 2, and we omit them here.
\\
\end{proof}
On the other hand, we highlight that, if the pair $( \mathbf{A} , \mathbf{B} )$ is a first-order stationary point of $\varphi^{*}(.)$ such that $\mathbf{A} \mathbf{B}  \in \mathbb{R}^{p \times n}_{<k}$, then $\mathbf{A} \mathbf{B}$ is not necessarily a Frechet first-order critical point of $\varphi(.)$, as noted by Ha et al.~\cite{HLB2020}. As an illustration, consider the pair $( \mathbf{0}^{p \times k} , \mathbf{0}^{k \times n} )$. Obviously, this pair is a  first-order critical point of $\varphi^{*}(.)$, but $ \mathbf{0}^{p \times k} \mathbf{0}^{k \times n} = \mathbf{0}^{p \times n} $ is not a Frechet first-order critical point of  $\varphi(.)$ in the sense of Theorem~\ref{theo3.5:box} as $\nabla \varphi ( \mathbf{0}^{p \times n} ) = - \mathbf{W} \odot \mathbf{X}$, which is not equal to $\mathbf{0}^{p \times n}$ as soon as we have for some pair of integers $(i,j)$, $\mathbf{X}_{ij} \ne 0$ and  $\mathbf{W}_{ij}  > 0$. Thus, in general, $\mathbf{0}^{p \times n}$ is not a Frechet first-order stationary point of $\varphi(.)$ and is obviously not a solution of of the WLRA problem in its formulation~\eqref{eq:P0}.

In addition, it is also possible to demonstrate that if the pair $( \mathbf{A} , \mathbf{B} )$ is not a second-order stationary point of $\varphi^{*}(.)$, then $\mathbf{A} \mathbf{B}  \in \mathbb{R}^{p \times n}_{ \le k}$ is not a (local) minimizer of $\varphi(.)$ over $\mathbb{R}^{p \times n}_{ \le k}$, see Ha et al.~\cite{HLB2020} and Levin et al.~\cite{LKB2025} for details.

\subsection{Approximate and regularized forms of the WLRA problem} \label{approx_wlra:box}

In practice, instead of an exact solution of the WLRA problem, which can even not exist if missing values are present as noted above, one often seeks an approximation of $\mathbf{X}$ such that
\begin{equation*}
\Vert \sqrt{\mathbf{W}} \odot ( \mathbf{X} - \widehat{\mathbf{X}} )  \Vert^{2}_{F}  \le ( 1 + \varepsilon ) \bar{\mathbf{c}}_{\varphi}  \ ,
\end{equation*}
where $\widehat{\mathbf{X}} \in \mathbb{R}^{p \times n}_{\le k}$ denotes the approximation, $\bar{\mathbf{c}}_{\varphi}$ is the infimum of $\varphi(.)$ and $\varepsilon \in ( 0, 1)$ is a tolerance parameter called the approximation error. In such framework, Razenshteyn et al.~\cite{RSW2016} recently show that in the case that $\mathbf{W}$ has at most $r$ distinct rows and $r$ distinct columns, there is an algorithm solving the above approximate version of the WLRA problem in $2^{O(k^2 . r/\varepsilon) } \text{poly}(n)$ time with probability of success at least $9/10$. In the case that $\mathbf{W}$ has at most $r$ distinct columns, but any number of distinct rows, there is also an algorithm solving the approximate version of the WLRA problem in $2^{O(k^2 . r^2 /\varepsilon) } \text{poly}(n)$ time with probability $9/10$. These bounds imply that for constant $r$ and $\varepsilon$, even if $r$ is as large as $\Theta\big(\text{log}(n)\big)$ in the first case, and $\Theta\big(\sqrt{ \text{log}(n)} \big)$ in the second case, the corresponding algorithms are polynomial time. Razenshteyn et al.~\cite{RSW2016} also consider the case when the rank of the weight matrix $\mathbf{W}$  is at most $r$, which includes as special cases the two above cases, and devise an $n^{O(k^2 . r/\varepsilon) }$ time algorithm for this more general case again with probability $9/10$. In other words, assuming that $\mathbf{W}$ has low rank, the algorithms of~\cite{RSW2016} achieve a $1 + \varepsilon$ multiplicative approximation to the infimum of $\varphi(.)$.

Alternatively, some authors have recently developed simple and greedy algorithms with additive error bounds for the WLRA problem which do not require any structural assumption on $\mathbf{W}$, see Bhaskara et al.~\cite{BRW2021} for general weights and also Musco et al.~\cite{MMW2021} in the case of binary weights. In such approach, one seeks an approximation $\widehat{\mathbf{X}}$ of $\mathbf{X}$ such that
\begin{equation*}
\Vert \sqrt{\mathbf{W}} \odot ( \mathbf{X} - \widehat{\mathbf{X}} )  \Vert^{2}_{F}  \le  \bar{\mathbf{c}}_{\varphi} + \varepsilon \Vert  \mathbf{X} \Vert^{2}_{F} \text{ or } \Vert \sqrt{\mathbf{W}} \odot ( \mathbf{X} - \widehat{\mathbf{X}} )  \Vert^{2}_{F}  \le   \gamma + \varepsilon \Vert  \mathbf{X} \Vert^{2}_{F} \ ,
\end{equation*}
where $\gamma$ is a (small) real constant of the order of $\bar{\mathbf{c}}_{\varphi}$ and the rank of $\widehat{\mathbf{X}}$ is of the order of $k$. Such methods with additive guarantees are interesting in applications (e.g., give sufficient matrix compression) when $\bar{\mathbf{c}}_{\varphi}$ is only a small fraction of the squared Frobenius norm of $\mathbf{X}$.

However, as these different algorithms with provable guarantees are inherently slow due the hardness of the WLRA problem and it is an open problem to determine when the WLRA problem has a closed form solution in general when some of the weights are zero, several authors have also proposed to minimize other related cost functions, which are convex, more smooth, and with a well-defined, nonempty and compact set of global minimizers, instead of problems~\eqref{eq:P0} or~\eqref{eq:P1} to address these issues~\cite{DKM2012}\cite{BA2015}\cite{MHT2010}\cite{MMBS2013}\cite{RS2005}\cite{SRJ2005}\cite{KM2010}\cite{KMO2010}\cite{BWZ2019}.
\\
\\
As a first illustration,~\cite{MHT2010}\cite{MMBS2013} have proposed the following convex relaxation to the rank constraint imposed in the formulation~\eqref{eq:P0}:
 \begin{equation*}
\min_{\mathbf{Y}\in\mathbb{R}^{p \times n}} \, \quad\  \varphi_{\lambda}( \mathbf{Y} ) = \frac{1}{2}  \Vert \sqrt{\mathbf{W}} \odot ( \mathbf{X} - \mathbf{Y} )  \Vert^{2}_{F} + \lambda \Vert \mathbf{Y} \Vert_{*} \ .
\end{equation*}
Here $\Vert \mathbf{Y} \Vert_{*}$ is the nuclear norm (also called the trace norm), which is equal to the sum of the singular values of the $p \times n$ matrix $\mathbf{Y}$ and $\lambda \in \mathbb{R}_{+*}$ is a regularization parameter controlling the nuclear norm of the minimizer $\widehat{\mathbf{Y}}(\lambda)$ of this Lagrange form of~\eqref{eq:P0}. $\varphi_{\lambda}(.)$ defines a convex function of its argument so that the above problem as an unique solution. Furthermore, it can be demonstrated that the rank of $\widehat{\mathbf{Y}}(\lambda)$ tends to zero when $\lambda$ grows unbounded so that this proxy can provide suboptimal low-rank minimizers of problem~\eqref{eq:P0} when this Lagrange form of~\eqref{eq:P0} is solved for a range of values of $\lambda$~\cite{MHT2010}\cite{MMBS2013}. Moreover, as the rank of $\widehat{\mathbf{Y}}(\lambda)$ increases when $\lambda$ decreases, if this problem is solved for a range of decreasing values of $\lambda$, the iterative algorithm can use efficiently the solution for the previous value of $\lambda$ as warm starts~\cite{MHT2010}\cite{HMLZ2015}.
\\
\\
Another class of related methods are maximum margin matrix factorization (MMMF) methods~\cite{RS2005}  ~\cite{SRJ2005}\cite{BWZ2019}\cite{LZT2019}, which use a factorization model of the matrix $\mathbf{Y}$, as in the formulation~\eqref{eq:P1} of the WLRA problem, but are also equipped with a regularization term $\lambda \in \mathbb{R}_{+*}$ as in the above Lagrange form of problem~\eqref{eq:P0}:
\begin{equation}  \label{eq:MMMF} \tag{MMMF}
\min_{\mathbf{A}\in\mathbb{R}^{p \times k}\text{, }\mathbf{B}\in\mathbb{R}^{k \times n} }   \, \quad\  \varphi^{*}_{\lambda}( \mathbf{A},\mathbf{B} ) = \frac{1}{2}   \Vert  \sqrt{\mathbf{W}}  \odot ( \mathbf{X} - \mathbf{A}\mathbf{B} )  \Vert^{2}_{F} + \frac{\lambda}{2} ( \Vert \mathbf{A} \Vert^{2}_{F} + \Vert \mathbf{B} \Vert^{2}_{F} ) \ .
\end{equation}
\\
Not surprisingly (e.g., taking into account the equivalence between the original problems~\eqref{eq:P0} and~\eqref{eq:P1} stated in Theorem~\ref{theo3.1:box}), there are closed relationships between the set of global minimizers of these Lagrange and regularized formulations of problems~\eqref{eq:P0} and~\eqref{eq:P1}, see Theorem 3 and Lemma 6 in Mazumder et al.~\cite{MHT2010} and also Hastie et al.~\cite{HMLZ2015} for details.
However, the above MMMF criterion is not convex in $(\mathbf{A},\mathbf{B})$, but only bi-convex as for the original problem~\eqref{eq:P1}, e.g., for a fixed $\mathbf{B}$ matrix, the modified function $ \varphi^{*}_{\lambda}(.)$ is convex in $\mathbf{A}$, and for a fixed $\mathbf{A}$ matrix, the function $\varphi^{*}_{\lambda}(.)$ is convex in $\mathbf{B}$. As the MMMF criterion is not convex, it can have possibly several local minima as the original problem~\eqref{eq:P1}~\cite{GZ1979}\cite{SJ2004}\cite{RS2005} and ALS algorithms (see Section~\ref{nipals:box}), which are very often used to solve these MMMF and~\eqref{eq:P1} problems, get frequently stuck in sub-optimal local minima for a small value of $k$ or a poorly chosen starting point, especially if some elements of the weight matrix $\mathbf{W}$ are equal to zero~\cite{GZ1979}\cite{SJ2004}. However, Ban et al.~\cite{BWZ2019} have demonstrated, extending the results of Razenshteyn et al.~\cite{RSW2016}, that it also exists polynomial time algorithms solving this weighted and regularized MMMF formulation of the WLRA problem, with provable guarantees, and also sharper time bounds than those proved in~\cite{RSW2016}.
\\
\\
Some other recent works have proposed to add to $\varphi^{*} ( \mathbf{A},\mathbf{B} )$, or similar regularized cost functions using the bilinear Burer-Monteiro approach, a balancing regularizer of the form
\begin{equation*}
R( \mathbf{A} ,\mathbf{B} ) = \frac{\lambda}{4}  \Vert  \mathbf{A}^{T} \mathbf{A} -   \mathbf{B}^{T} \mathbf{B} \Vert^{2}_{F} \ ,
\end{equation*}
where $\lambda$ controls the weight for the regularizer as before~\cite{PKCS2017}\cite{WZG2017}\cite{ZLTW2018}\cite{ZLTW2021}. $R( \mathbf{A} ,\mathbf{B} )$ implicitly forces the $\mathbf{A}$ and $\mathbf{B}$ matrices to have the same energy and, thus, helps to remove the scaling ambiguity which inherently affects the cost function $\varphi^{*} (.)$  and the minimization of  $\varphi^{*} ( \mathbf{A},\mathbf{B} )$ in the~\eqref{eq:P1} formulation of the WLRA problem as discussed in Remark~\ref{remark3.2:box} above. Moreover, for many cost functions which use the bilinear Burer-Monteiro approach, adding this balancing regularizer does not compromise the quality of the solutions~\cite{PKCS2017}\cite{LZT2019}\cite{ZLTW2021}\cite{OUV2023}.
\\
\\
Many of the proposed recent approaches also recast the WLRA problem as an optimization problem on the Grassmann manifold $\text{Gr}(p,k)$ or on the two Grassmann manifolds $\text{Gr}(p,k)$ and $\text{Gr}(n,k)$ (where $\text{Gr}(p,k)$ is the set of $k$-dimensional linear subspaces of $\mathbb{R}^p$) and introduce a regularization parameter $\lambda \in \mathbb{R}_{+*}$ as in the above Lagrange forms of problems~\eqref{eq:P0} and~\eqref{eq:P1} in order to ensure smoothness of the objective function and  hence obtain good convergence at the expense of slight increase of the objective~\cite{KM2010}\cite{KMO2010}\cite{DKM2012}\cite{MMBS2013}\cite{BA2015}. An interesting example in this class of methods, as it is closely related to the formulations~\eqref{eq:P0} or~\eqref{eq:P1} of the WLRA problem, is the unconstrained Riemannian optimization methods on a single Grassmann manifold $\text{Gr}(p,k)$ described in Boumal and Absil~\cite{BA2011}\cite{BA2015} for solving the matrix completion problem, which we now discussed in some details.
\\
\\
To this end, for any weight matrix $\mathbf{W}\in\mathbb{R}^{p \times n}_{+}$, let us define the set $\bar{\Omega} \subset \lbrack p \rbrack \times  \lbrack n \rbrack$, be the set of indices of the elements of $\mathbf{W}$ with $\mathbf{W}_{ij} = 0$ (e.g., $\bar{\Omega}$ is the complement of $\Omega$ in $\lbrack p \rbrack \times  \lbrack n \rbrack$) and the seminorms 
\begin{equation*}
\Vert  \mathbf{Y} \Vert^2_{\Omega} =  \sum_{(i,j) \in \Omega}  \mathbf{Y}^2_{ij}  \text{ and } \Vert  \mathbf{Y} \Vert^2_{\bar{\Omega}} =  \sum_{(i,j) \in \bar{\Omega}}  \mathbf{Y}^2_{ij} \ .
\end{equation*}
\\
With these definitions and in our notations, Boumal and Absil~\cite{BA2011}\cite{BA2015} proposed to solve the following optimization problem
\begin{equation*}
\min_{\mathbf{Y} \in \mathbb{R}^{p \times n}_{\le k} }  \, \quad\ g( \mathbf{Y} ) = \frac{1}{2}   \Vert  \sqrt{\mathbf{W}}  \odot ( \mathbf{X} - \mathbf{Y} )  \Vert^{2}_{\Omega} + \frac{\lambda}{2} \Vert \mathbf{Y} \Vert^{2}_{\bar{\Omega}} \ ,
\end{equation*}
where, as before, $\lambda \in \mathbb{R}_{+*}$ is a regularization parameter, which ensures that the solution to this problem exists and the cost function $g$ is smooth. They give the following interpretation for the minimization of the cost function $g$, which makes sense for the matrix completion problem: "we are looking for an optimal matrix $\widehat{\mathbf{X}}$ of rank at most $k$ and we have confidence $\sqrt{\mathbf{W}}_{ij}$ that $\widehat{\mathbf{X}}_{ij}$ should equal  $\mathbf{X}_{ij}$ for $(i,j)  \in \Omega$ and smaller confidence $\lambda$ that $\widehat{\mathbf{X}}_{ij}$ should equal  zero for $(i,j)  \in  \bar{\Omega}$". They have also illustrated that the solutions of this problem are largely insensitive to the value of $\lambda$  provided it is much smaller than the strictly positive values $\mathbf{W}_{ij}$. As an illustration, for matrix completion problems in their experiments, they used $\lambda = 10^{-6}$ and  $\mathbf{W}_{ij} = 1$ if $(i,j)  \in \Omega$. Finally, they describe and apply second-order Riemannian trust-region methods (RTRMC2) and Riemannian conjugate gradient methods (RCGMC)~\cite{B2023} to solve this problem efficiently and accurately, which are still state-of-the-art algorithms on a wide range of problem instances.
\\
\\
Interestingly, we now show that the minimization of the cost function $g(.)$ proposed by Boumal and Absil~\cite{BA2011}\cite{BA2015} is in fact a simple instance of  formulation~\eqref{eq:P0} of the WLRA problem so that the variable projection framework can also be used to solve this problem as we will illustrate in the following sections. More precisely, if, for any $p \times n$ weight matrix $\mathbf{W}$ with some zero elements and any $\lambda \in \mathbb{R}_{+*}$ (e.g., $\lambda > 0$), we define as above an $p \times n$ weight matrix $\mathbf{W}_\lambda \in \mathbb{R}^{p \times n}_{+*}$ as
 \begin{equation} \label{eq:weight_proj_op}
    \big \lbrack \mathbf{W}_\lambda \big \rbrack_{ij} =
    \begin{cases}
        \displaystyle{ \mathbf{W}_{ij}  } & \text{if } (i,j)  \in \Omega\\
         \lambda                                    & \text{if } (i,j)  \notin  \Omega
    \end{cases} \ ,
\end{equation}
and we introduce the projection operator associated with an $p \times n$ weight matrix $\mathbf{W}$ by $P_{\Omega} : \mathbb{R}^{p \times n}  \longrightarrow \mathbb{R}^{p \times n}$ with $P_{\Omega}(\mathbf{X}) = \mathbf{X}_{\Omega}$ where
 \begin{equation} \label{eq:proj_op}
    \big \lbrack \mathbf{X}_{\Omega} \big \rbrack_{ij} =
    \begin{cases}
        \displaystyle{ \mathbf{X}_{ij}  } & \text{if } (i,j)  \in \Omega\\
         0                                              & \text{if } (i,j)  \notin  \Omega
    \end{cases} \ ,
\end{equation}
we can rearrange the cost function $g(.)$ introduced by Boumal and Absil~\cite{BA2011}\cite{BA2015} as
\begin{equation} \label{eq:g_func}
g_{\lambda}( \mathbf{Y} ) = \frac{1}{2}   \Vert  \sqrt{\mathbf{W}_\lambda}  \odot ( \mathbf{X}_{\Omega} - \mathbf{Y} )  \Vert^{2}_{F} \ ,
\end{equation}
and it is readily observed that the minimization of this cost function $g_{\lambda}(.)$ w.r.t. $\mathbf{Y}\in\mathbb{R}^{p \times n}_{\le k}$ is equivalent to the form~\eqref{eq:P0}
\begin{equation*}
\min_{\mathbf{Y} \in \mathbb{R}^{p \times n}_{\le k} }   \quad\ g_\lambda( \mathbf{Y} ) = \frac{1}{2}   \Vert  \sqrt{\mathbf{W}_\lambda}  \odot ( \mathbf{X}_{\Omega} - \mathbf{Y} )  \Vert^{2}_{F} = \frac{1}{2}   \Vert   \mathbf{X}_{\Omega} - \mathbf{Y}  \Vert^{2}_{\mathbf{W}_\lambda}
\end{equation*}
of a standard WLRA problem in which we use the matrices $\mathbf{X}_{\Omega}$ and $\mathbf{W}_\lambda$ in place of $\mathbf{X}$ and $\mathbf{W}$, respectively. Furthermore, as all the elements of the weight matrix $\mathbf{W}_\lambda$ are greater than zero for any  $\lambda \in \mathbb{R}_{+*}$, $\Vert \Vert_{\mathbf{W}_\lambda}$ defines a norm on $\mathbb{R}^{p \times k}$ and Theorem~\ref{theo3.3:box} shows that the set of global minimizers of $g_{\lambda}(.)$ is nonempty and compact, so that the minimization of this cost function is a well-posed problem. In other words, for any $\lambda \in \mathbb{R}_{+*}$ there exists $\widehat{\mathbf{X}}_\lambda  \in \mathbb{R}^{p \times n}_{\le k}$ such that
\begin{equation*}
\widehat{\mathbf{X}}_\lambda  = \text{Arg}\min_{\mathbf{Y}\in\mathbb{R}^{p \times n}_{\le k} }  g_\lambda( \mathbf{Y} ) \ .
\end{equation*}
In addition, if we take a regularization parameter $\lambda$ (also called the Tikhonov parameter, see~\cite{GW2000}) sufficiently small, the following theorem shows that the minimization of  $g_{\lambda}(.)$ with a Tikhonov parameter tending to zero is an interesting alternative to the formulations~\eqref{eq:P0} and~\eqref{eq:P1} of the WLRA problem, which are not well-posed when some elements of the weight matrix $\mathbf{W}$ are equal to zero as discussed above.
 \begin{theo3.8} \label{theo3.8:box}
 Let $\mathbf{X} \in \mathbb{R}^{p \times n}$, $\mathbf{W} \in \mathbb{R}^{p \times n}_+$ (i.e., $\mathbf{W}_{ij} \ge 0$),  $k \in \mathbb{N}_*$ with $k \le \emph{rank}( \mathbf{X} ) \le \text{min}( {p},{n} )$ and $\lambda \in \mathbb{R}_{+*}$ (i.e., $\lambda > 0$). Furthermore, using definition~\eqref{eq:g_func} of the cost function, $g_{\lambda}(.)$, let
 \begin{equation*}
\widehat{\mathbf{X}}_\lambda = \text{Arg}\min_{\mathbf{Y}\in\mathbb{R}^{p \times n}_{\le k} }  g_\lambda( \mathbf{Y} )  \text{ and } f(\lambda) =  g_\lambda( \widehat{\mathbf{X}}_\lambda ) = \frac{1}{2} \Vert \mathbf{X}_{\Omega} - \widehat{\mathbf{X}}_\lambda \Vert^2_{\mathbf{W}_\lambda} \text{ for } \lambda \in \mathbb{R}_{+*}
\end{equation*}
then
\begin{equation*}
\lim_{\lambda \to 0 } \, f(\lambda) =  \bar{\mathbf{c}}_{\varphi} \ ,
\end{equation*}
where $\bar{\mathbf{c}}_{\varphi}$ is the infimum of the cost function $\varphi(.)$ used in the formulation~\eqref{eq:P0} of the WLRA problem and $\mathbf{X}_{\Omega} = P_{\Omega}(\mathbf{X})$ where $P_{\Omega}$ is the projection operator associated with the $p \times n$ weight matrix $\mathbf{W}$.
\end{theo3.8}
\begin{proof}
We first show that $f(.)$ has a well defined limit, $\bar{\mathbf{c}}_{f}$, when $f(.)$ tends to zero. To demonstrate this result, we first note that $f(.)$ is an increasing function. For $\alpha \in \mathbb{R}_{+*}$ and $\lambda \in \mathbb{R}_{+*}$ with $\alpha \ge \lambda$, let
\begin{equation*}
\widehat{\mathbf{X}}_\alpha = \text{Arg}\min_{\mathbf{Y}\in\mathbb{R}^{p \times n}_{\le k} }  g_\alpha( \mathbf{Y} )  \text{ and } \widehat{\mathbf{X}}_\lambda = \text{Arg}\min_{\mathbf{Y}\in\mathbb{R}^{p \times n}_{\le k} }  g_\lambda( \mathbf{Y} )
\end{equation*}
then we have
\begin{equation*}
 \Vert \mathbf{X}_{\Omega} - \widehat{\mathbf{X}}_\alpha \Vert^2_{\mathbf{W}_\alpha} \ge  \Vert \mathbf{X}_{\Omega} - \widehat{\mathbf{X}}_\alpha \Vert^2_{\mathbf{W}_\lambda} \ge  \Vert \mathbf{X}_{\Omega} - \widehat{\mathbf{X}}_\lambda \Vert^2_{\mathbf{W}_\lambda} \ ,
\end{equation*}
which implies that $f(\alpha) \ge f(\lambda)$. Furthermore, for all $\lambda \in \mathbb{R}_{+*}$, we have
\begin{equation*}
 f(\lambda) = \frac{1}{2}  \Vert \mathbf{X}_{\Omega} - \widehat{\mathbf{X}}_\lambda \Vert^2_{\mathbf{W}_\lambda} \ge \frac{1}{2} \Vert \mathbf{X}_{\Omega} - \widehat{\mathbf{X}}_\lambda \Vert^2_{\mathbf{W}} = \varphi( \widehat{\mathbf{X}}_\lambda ) \ge \bar{\mathbf{c}}_{\varphi} \ ,
\end{equation*}
which shows that $\lim_{\lambda \to 0 } \, f(\lambda)$ exists and that $\lim_{\lambda \to 0 } \, f(\lambda) = \bar{\mathbf{c}}_{f} \ge \bar{\mathbf{c}}_{\varphi}$.

It remains to show that $\bar{\mathbf{c}}_{\varphi} \ge  \bar{\mathbf{c}}_{f}$. To this end, suppose that $\bar{\mathbf{c}}_{\varphi} < \bar{\mathbf{c}}_{f}$, then it exists $\mathbf{Y} \in  \mathbb{R}^{p \times n}_{\le k}$ such that $\bar{\mathbf{c}}_{\varphi} \le \varphi(\mathbf{Y}) < \bar{\mathbf{c}}_{f}$, otherwise $\bar{\mathbf{c}}_{\varphi}$ is not the infimum of $\varphi(.)$. As
\begin{equation*}
\lim_{\lambda \to 0 } \, \Vert \mathbf{X}_{\Omega} - \mathbf{Y}  \Vert^2_{\mathbf{W}_\lambda} = \Vert \mathbf{X} - \mathbf{Y}  \Vert^2_{\mathbf{W}}  \text{ and }  \Vert \mathbf{X}_{\Omega} - \mathbf{Y}  \Vert^2_{\mathbf{W}_\lambda} \ge \Vert \mathbf{X} - \mathbf{Y}  \Vert^2_{\mathbf{W}} \text{ for all } \lambda \in \mathbb{R}_{+*} \ ,
\end{equation*}
it also exists $\alpha \in \mathbb{R}_{+*}$ such that
\begin{equation*}
\varphi(\mathbf{Y}) = \frac{1}{2} \Vert \mathbf{X} - \mathbf{Y}  \Vert^2_{\mathbf{W}} \le  \frac{1}{2} \Vert \mathbf{X}_{\Omega} - \mathbf{Y}  \Vert^2_{\mathbf{W}_\alpha} < \bar{\mathbf{c}}_{f}  \  .
\end{equation*}
However, we also have
\begin{equation*}
\bar{\mathbf{c}}_{f}  \le  f( \alpha )  = \frac{1}{2} \Vert \mathbf{X}_{\Omega} - \widehat{\mathbf{X}}_\alpha \Vert^2_{\mathbf{W}_\alpha}  \le \frac{1}{2} \Vert \mathbf{X}_{\Omega} - \mathbf{Y}  \Vert^2_{\mathbf{W}_\alpha} < \bar{\mathbf{c}}_{f}
\end{equation*}
and we obtain a contradiction.
\\
\end{proof}

Thus, one way of getting an useful approximate solution to the WLRA problem when missing values are present is to use a continuation Tikhonov method that approximately solves a sequence of  regularized WLRA problems for a sequence of decreasing Tikhonov parameter $\lambda$. The approximate solution of one regularized WLRA problem with Tikhonov parameter $\lambda_{t}$ (e.g., the minimization of $g_{\lambda_{t}}(.)$ ) is taken as the starting point for the next regularized WLRA problem with Tikhonov parameter $\lambda_{t+1}<\lambda_{t}$. This kind of Tikhonov methods has already been proposed in the context of ill-conditioned and uniformly rank-deficient NLLS problems~\cite{E1996}\cite{EW1996}\cite{EWGS2005}, see Section~\ref{vpalg:box} where such methods are further discussed.

\subsection{Variable projection formulation of the WLRA problem} \label{varpro_wlra:box}

We are now ready to show that the alternative formulation~\eqref{eq:P1} or its variants (see Remark~\ref{remark3.1:box}) of the WLRA problem
 \begin{equation*}
\min_{\mathbf{A}\in\mathbb{R}^{p \times k}\text{,}\mathbf{B}\in\mathbb{R}^{k \times n} }  \varphi^{*}( \mathbf{A},\mathbf{B}  ) = \frac{1}{2}   \Vert  \sqrt{\mathbf{W}} \odot ( \mathbf{X} - \mathbf{A}\mathbf{B} )  \Vert^{2}_{F}
\end{equation*}
is a separable NLLS problem as stated in the Definition~\ref{def2.10:box} of Subsection~\ref{calculus:box}~\cite{GP1973}\cite{RW1980}. This means that the minimization of $\varphi^{*}( \mathbf{A}, \mathbf{B} )$ is a mixed linear-nonlinear least-squares problem where the associated residual function $e(\mathbf{A}, \mathbf{B} )$ is linear in some variables and nonlinear in others.
\\
\\
In order to demonstrate this result, we first write $\varphi^{*}( \mathbf{A}, \mathbf{B} )$ as
\begin{equation*}
\varphi^{*}( \mathbf{A}, \mathbf{B} ) =  \frac{1}{2}   \Vert e(\mathbf{A},\mathbf{B}) \Vert^{2}_{2} = \frac{1}{2} {e(\mathbf{A},\mathbf{B})}^{T} e(\mathbf{A},\mathbf{B}) \ ,
\end{equation*}
where the residual vector function  $e(  \mathbf{A}, \mathbf{B} ) \in \mathbb{R}^{p.n}$ is defined by
\begin{equation}  \label{eq:res_func}
e( \mathbf{A},\mathbf{B} ) = \emph{vec} \big (  \sqrt{\mathbf{W}} \odot ( \mathbf{X} - \mathbf{A}\mathbf{B} )  \big) \ .
\end{equation}
Using equations~\eqref{eq:vec_hadamard} and~\eqref{eq:vec_kronprod}, the residual function $e(  \mathbf{A}, \mathbf{B} )$ can be further transformed as
\begin{align*}
    e(  \mathbf{A}, \mathbf{B} ) & =   \emph{diag}\big( \emph{vec}(\sqrt{\mathbf{W}}) \big) \emph{vec}( \mathbf{X} - \mathbf{A}\mathbf{B} )   \\
    & =    \emph{vec}( \sqrt{\mathbf{W}}  \odot \mathbf{X} ) - \emph{diag}\big( \emph{vec}( \sqrt{\mathbf{W}} ) \big) \emph{vec}( \mathbf{A}\mathbf{B} )   \\
    & =    \emph{vec}( \sqrt{\mathbf{W}}  \odot \mathbf{X} ) - \emph{diag}\big( \emph{vec}( \sqrt{\mathbf{W}} ) \big) \left(  \mathbf{I}_n  \otimes \mathbf{A}  \right) \emph{vec}( \mathbf{B} ) \ ,
\end{align*}
and $e(  \mathbf{A}, \mathbf{B} )$ is finally equal in explicit matrix form to
\begin{equation*}
      \left\lbrack
\begin{array}{ccccc}
\sqrt{\mathbf{W}}_{.1}\mathbf{X}_{.1}  \\
\vdots \\
\sqrt{\mathbf{W}}_{.j}\mathbf{X}_{.j} \\
\vdots \\
\sqrt{\mathbf{W}}_{.n}\mathbf{X}_{.n}
\end{array} \right\rbrack -
    \left\lbrack
\begin{array}{ccccc}
\emph{diag}(\sqrt{\mathbf{W}}_{.1})\mathbf{A} & 0          & \ldots             & 0          & 0          \\
0                      & \ddots &  0                    & \ldots   & 0          \\
\vdots              & 0         & \emph{diag}(\sqrt{\mathbf{W}}_{.j})\mathbf{A} & 0           & \vdots \\
0                      & \ldots  & 0                     & \ddots  & 0          \\
0                      & 0         & \ldots              & 0           & \emph{diag}(\sqrt{\mathbf{W}}_{.n})\mathbf{A}
\end{array} \right\rbrack
      \left\lbrack
\begin{array}{ccccc}
\mathbf{B}_{.1}  \\
\vdots \\
\mathbf{B}_{.j} \\
\vdots \\
\mathbf{B}_{.n}
\end{array} \right\rbrack \ .
\end{equation*}

In this residual function, we first note that all the lines corresponding to a zero weight (e.g., $\mathbf{W}_{ij}=0$) can be eliminated when evaluating this function in real computations. The same is true for all the equations of the following sections and in a practical computer implementation of the algorithms used to minimize  $\varphi^{*}(.)$. However, for notational simplicity and because we want to consider at the same time both the cases $\mathbf{W}\in\mathbb{R}^{p \times n}_+$ and $\mathbf{W}\in\mathbb{R}^{p \times n}_{+*}$, we do not introduce an incidence matrix in our equations to indicate which rows or columns must be eliminated as was done for example in~\cite{OD2007}\cite{C2008b}\cite{D2011}\cite{GM2011}\cite{BA2015}. Then, we may write
\begin{equation*}
\varphi^{*}(\mathbf{A}, \mathbf{B}) = \frac{1}{2} \Vert \mathbf{x}  - \mathbf{F}(  \mathbf{a} ) \mathbf{b} )  \Vert^{2}_{2}  \ ,
\end{equation*}
where $\mathbf{x} = \emph{vec}( \sqrt{\mathbf{W}} \odot \mathbf{X} )$, $\mathbf{a} = \emph{vec}( \mathbf{A}^{T} )$, $\mathbf{b} = \emph{vec}( \mathbf{B} )$ and $\mathbf{F}(  \mathbf{a} )$ is the block diagonal matrix
\begin{equation} \label{eq:F_mat}
    \mathbf{F}(  \mathbf{a} ) = \bigoplus_{j=1}^n \mathbf{F}_{j}(  \mathbf{a} ) =
    \left\lbrack
\begin{array}{ccccc}
\mathbf{F}_{1}(  \mathbf{a} )& 0          & \ldots             & 0          & 0          \\
0                      & \ddots &  0                    & \ldots   & 0          \\
\vdots              & 0         & \mathbf{F}_{j}(  \mathbf{a} ) & 0           & \vdots \\
0                      & \ldots  & 0                     & \ddots  & 0          \\
0                      & 0         & \ldots              & 0           & \mathbf{F}_{n}(  \mathbf{a} )
\end{array} \right\rbrack
= \emph{diag}\big( \emph{vec}( \sqrt{\mathbf{W}} ) \big)  \big(  \mathbf{I}_n  \otimes \mathbf{A}  \big) \ ,
\end{equation}
where
\begin{equation*}
    \mathbf{F}_{j}(  \mathbf{a} ) =  \emph{diag}(\sqrt{\mathbf{W}}_{.j})\mathbf{A} = \emph{diag}(\sqrt{\mathbf{W}}_{.j}) \big( \emph{mat}_{k \times p} (\mathbf{a}) \big )^{T} \ .
\end{equation*} 
The reason and interest of defining the vectorized form of $\mathbf{A}$ as
\begin{equation} \label{eq:d_veca}
\mathbf{a} = \emph{vec}( \mathbf{A}^{T} ) \ , 
\end{equation} 
instead of simply $\emph{vec}( \mathbf{A} )$ as usually done, will become clear in the next sections.
From this formulation, it is clear that minimizing $\varphi^{*}(.)$ is a separable NLLS problem, since for a fixed matrix $\mathbf{A}$, we have a linear least-squares problem to determine the optimal vector $\mathbf{\widehat{b}} = \emph{vec}( \mathbf{\widehat{B}} )$, i.e.,
\begin{equation*}
\mathbf{\widehat{b}} = \text{Arg}\min_{ \mathbf{b}\in\mathbb{R}^{n.k}} \, \varphi^{*}(\mathbf{A},\mathbf{B}  ) = \frac{1}{2}   \Vert  \mathbf{x}  - \mathbf{F}(  \mathbf{a} ) \mathbf{b} )  \Vert^{2}_{2}  \ .
\end{equation*}
Moreover, we observe that the residual function $e(  \mathbf{A}, \mathbf{B} )$ is linear in both $\mathbf{A}$ and $\mathbf{B}$, since
\begin{align*}
\mathbf{F}(  \mathbf{a} ) \mathbf{b}  & =  \Big(  \bigoplus_{j=1}^n \emph{diag}(\sqrt{\mathbf{W}}_{.j})\mathbf{A}  \Big)  \mathbf{b} \\
              & =    \emph{diag}\big( \emph{vec}( \sqrt{\mathbf{W}} ) \big)  \big(  \mathbf{I}_n  \otimes \mathbf{A}  \big)  \emph{vec}( \mathbf{B} ) \\
              & =   \emph{diag}\big( \emph{vec}( \sqrt{\mathbf{W}} ) \big)  \emph{vec}( \mathbf{A}\mathbf{B} ) \\
              & =   \emph{diag}\big( \emph{vec}( \sqrt{\mathbf{W}} ) \big)   \big(  \mathbf{B}^{T}  \otimes  \mathbf{I}_p  \big)  \emph{vec}( \mathbf{A} ) \\
              & =   \emph{diag}\big( \emph{vec}( \sqrt{\mathbf{W}} ) \big)   \big(  \mathbf{B}^{T}  \otimes  \mathbf{I}_p  \big)  \mathbf{K}_{(k,p)} \mathbf{K}_{(p,k)} \emph{vec}( \mathbf{A} ) \\
              & =   \emph{diag}\big( \emph{vec}( \sqrt{\mathbf{W}} ) \big)   \big(  \mathbf{B}^{T}  \otimes  \mathbf{I}_p  \big)  \mathbf{K}_{(k,p)}  \emph{vec}( \mathbf{A}^{T} ) \\
              & =   \emph{diag}\big( \emph{vec}( \sqrt{\mathbf{W}} ) \big)   \mathbf{K}_{(n,p)} \big( \mathbf{I}_p   \otimes  \mathbf{B}^{T} \big)  \mathbf{a} \\
              & =    \mathbf{K}_{(n,p)}  \emph{diag}\big( \emph{vec}( \sqrt{\mathbf{W}}^{T} ) \big)  \big( \mathbf{I}_p   \otimes  \mathbf{B}^{T} \big)  \mathbf{a} \\
              & =    \mathbf{K}_{(n,p)}   \Big(   \bigoplus_{i=1}^p   \emph{diag}(\sqrt{\mathbf{W}}_{i.})\mathbf{B}^{T}   \Big)  \mathbf{a} \ .
\end{align*}
Defining now
\begin{equation} \label{eq:G_mat}
    \mathbf{G}(  \mathbf{b} )  = \bigoplus_{i=1}^p   \mathbf{G}_{i}(  \mathbf{b} ) =
    \left\lbrack
\begin{array}{ccccc}
\mathbf{G}_{1}(  \mathbf{b} )  & 0          & \ldots             & 0          & 0          \\
0                      & \ddots &  0                    & \ldots   & 0          \\
\vdots              & 0         & \mathbf{G}_{i}(  \mathbf{b}  ) & 0           & \vdots \\
0                      & \ldots  & 0                     & \ddots  & 0          \\
0                      & 0         & \ldots              & 0           & \mathbf{G}_{p}(  \mathbf{b}  )
\end{array} \right\rbrack \ ,
\end{equation}
where
\begin{equation*}
    \mathbf{G}_{i}(  \mathbf{b} ) =  \emph{diag}(\sqrt{\mathbf{W}}_{i.})\mathbf{B}^{T} = \emph{diag}(\sqrt{\mathbf{W}}_{i.}) \big( \emph{mat}_{k \times n} ( \mathbf{b} ) \big)^{T} \ ,
\end{equation*} 
we note that the residual function $e( \mathbf{A},\mathbf{B}  )$ may then be written in the following alternative matrix form
\begin{align*}
e( \mathbf{A},\mathbf{B}  )  & =  \mathbf{x} - \mathbf{F}(  \mathbf{a} ) \mathbf{b}  \\
              & =    \mathbf{x} -  \mathbf{K}_{(n,p)}   \Big(   \bigoplus_{i=1}^p   \emph{diag}(\sqrt{\mathbf{W}}_{i.})\mathbf{B}^{T}   \Big)  \mathbf{a}\\
              & =    \mathbf{x} -  \mathbf{K}_{(n,p)}  \mathbf{G}(  \mathbf{b} ) \mathbf{a}  \\
              & =    \mathbf{K}_{(n,p)} \mathbf{K}_{(p,n)} \mathbf{x} -  \mathbf{K}_{(n,p)}  \mathbf{G}(  \mathbf{b} ) \mathbf{a}  \\
              & =    \mathbf{K}_{(n,p)} \mathbf{K}_{(p,n)} \emph{vec}( \sqrt{\mathbf{W}}  \odot \mathbf{X} ) -  \mathbf{K}_{(n,p)}  \mathbf{G}(  \mathbf{b} ) \mathbf{a}  \\
              & =    \mathbf{K}_{(n,p)}  \emph{vec} \big( (\sqrt{\mathbf{W}}  \odot \mathbf{X})^{T} \big) -  \mathbf{K}_{(n,p)}  \mathbf{G}(  \mathbf{b} ) \mathbf{a}  \\
              & =    \mathbf{K}_{(n,p)}  \big( \mathbf{z}  -  \mathbf{G}(  \mathbf{b} ) \mathbf{a} \big) \ ,
\end{align*}
where $\mathbf{z} = \emph{vec} \big( (\sqrt{\mathbf{W}}  \odot \mathbf{X})^{T} \big)$. This implies that $\varphi^{*}(\mathbf{A}, \mathbf{B})$ may also be expressed as 
\begin{align*}
\varphi^{*}(\mathbf{A}, \mathbf{B})  & =  \frac{1}{2}   \left(  \mathbf{z} - \mathbf{G}(  \mathbf{b} ) \mathbf{a} \right)^{T} \mathbf{K}_{(p,n)} \mathbf{K}_{(n,p)} \left(  \mathbf{z} - \mathbf{G}(  \mathbf{b} ) \mathbf{a} \right)  \\
              & =     \frac{1}{2}   \left(   \mathbf{z} - \mathbf{G}(  \mathbf{b} ) \mathbf{a} \right)^{T}  \left(  \mathbf{z} - \mathbf{G}(  \mathbf{b} ) \mathbf{a} \right)  \\
              & =     \frac{1}{2}   \Vert  \mathbf{z} - \mathbf{G}(  \mathbf{b} ) \mathbf{a} \Vert^2_2 \ ,
\end{align*}
which shows that the roles of $\mathbf{A}$ and $\mathbf{B}$ are interchangeable in $\varphi^{*}(\mathbf{A}, \mathbf{B})$ as already noted in the case of binary weights for example in~\cite{OD2007}. As for the choice between the formulations~\eqref{eq:P1} and~\eqref{eq:P2} of the WLRA problem (see Remark~\ref{remark3.2:box}), the choice between the formulations
\begin{equation*}
\varphi^{*}(\mathbf{A}, \mathbf{B}) = \frac{1}{2} \Vert \mathbf{x}  - \mathbf{F}(  \mathbf{a} ) \mathbf{b} )  \Vert^{2}_{2}  \quad \text{or}  \quad \varphi^{*}(\mathbf{A}, \mathbf{B}) = \frac{1}{2}   \Vert  \mathbf{z} - \mathbf{G}(  \mathbf{b} ) \mathbf{a} \Vert^2_2
\end{equation*} 
depends on the values of $p$ and $n$, and the first one should be preferred if $p<n$ as the number of parameters to estimate (e.g., $\mathbf{A}$) will be smaller once the other matrix variable  (e.g., $\mathbf{B}$) has been eliminated as we will show below, and vice-versa if $p>n$. Furthermore, in what follows, we note that the matrices $\mathbf{A}$ and $\mathbf{B}$ can be used in an interchangeable manner with their vectorized forms $\mathbf{a}$ and $\mathbf{b}$, respectively, as the mapping $\emph{vec}(.)$ is a bijective homeomorphism (see Subsection~\ref{multlin_alg:box} for details). 

Thus, the problem of minimizing  $\varphi^{*}(.)$ is separable and this property can be exploited in a least-squares estimation, and a number of special purpose algorithms have been proposed in this context~\cite{K1974}\cite{K1975}\cite{GP1973}\cite{RW1980}\cite{B2009}\cite{OR2013}\cite{BL2020}. Moreover, it has been demonstrated that these special algorithms provide greater stability than standard NLLS methods, besides reducing both the dimensionality of the optimization problem and the necessary number of iteration steps~\cite{N2000}\cite{GP2003}\cite{D2011}\cite{BL2020}. In most cases, the total computational work decreases with separable methods even though the code describing the separable problem is slightly more complicated than in standard NLLS algorithms. We now discuss how to reformulate the problems~\eqref{eq:P0} and~\eqref{eq:P1} so that we can exploit the separation property by eliminating one of the matrix variables (e.g.,  $\mathbf{B}$ if $p<n$) and devise more efficient algorithms to solve these problems.

Assuming that $p<n$, for a fixed $\mathbf{A}$ matrix, we have a linear least-squares problem to determine the optimal vector $\mathbf{b} = \emph{vec}( \mathbf{B} )$ which will minimize the cost function $\varphi^{*}(.)$ and the solution of this linear least-squares problem is $\mathbf{\widehat{b}} =  \mathbf{F}(  \mathbf{a} )^{+}\mathbf{x}$ where $\mathbf{x} = \emph{vec}( \sqrt{\mathbf{W}} \odot \mathbf{X} )$ and $\mathbf{F}(  \mathbf{a} )^{+}$ is the pseudo-inverse of the $p.n \times k.n$ matrix  $\mathbf{F}(  \mathbf{a} )$, see Subsection~\ref{lin_alg:box}. Inserting now $\mathbf{\widehat{b}}$ in $\varphi^{*}(.)$, we obtain a new nonlinear functional $\psi : \mathbb{R}^{p.k} \longrightarrow \mathbb{R}$ involving only the vectorized form of the $\mathbf{A}$ matrix
\begin{equation}  \label{eq:psi_func}
\psi(\mathbf{a} ) =  \frac{1}{2}  \Vert \left( \mathbf{I}_{p.n}  - \mathbf{F}(  \mathbf{a} )  \mathbf{F}(  \mathbf{a} )^{+} \right) \mathbf{x} \Vert^{2}_{2} = \frac{1}{2} \Vert \mathbf{P}^{\bot}_{\mathbf{F}(  \mathbf{a} ) }\mathbf{x} \Vert^{2}_{2} =  \frac{1}{2} \Vert \mathbf{r}(  \mathbf{a} )  \Vert^{2}_{2} \ ,
\end{equation}
where $\mathbf{P}^{\bot}_{\mathbf{F}(  \mathbf{a} ) }$ is the orthogonal projector onto the orthogonal complement of $\emph{ran}\big( \mathbf{F}(  \mathbf{a} ) \big)$ and $\mathbf{r}(.)$ is a nonlinear residual function of $\mathbf{a} = \emph{vec}( \mathbf{A}^{T})$ defined by
\begin{equation} \label{eq:r_vec}
\mathbf{r} : \mathbb{R}^{p.k} \longrightarrow \mathbb{R}^{p.n} : \mathbf{a} \mapsto \mathbf{P}^{\bot}_{\mathbf{F}(  \mathbf{a} ) }\mathbf{x} =  \mathbf{x} - \mathbf{F}(\mathbf{a}) \mathbf{\widehat{b}} = \mathbf{K}_{(n,p)} \big( \mathbf{z} - \mathbf{G}(\mathbf{\widehat{b}}) \mathbf{a} \big)  \  .
\end{equation}
$\mathbf{r}(  \mathbf{a} )$ is called the variable projection residual of $\mathbf{X}$ at  $\mathbf{A}$ (or equivalently of $\mathbf{x}$ at  $\mathbf{a}$) and the functional $\psi(.)$ can be termed a variable projection functional since $\mathbf{P}^{\bot}_{\mathbf{F}(  \mathbf{a} ) }$ is an orthogonal projector involving only the vectorized form of the $\mathbf{A}$ matrix~\cite{GP1973}. Again, if we take into account the block structure of $\mathbf{F}(  \mathbf{a} )$, we obtain an alternative formulation of $\psi(.)$, which is useful for computational purposes,
\begin{equation*}
\psi( \mathbf{a} ) = \frac{1}{2} \sum_{j=1}^n {  \psi_{j}( \mathbf{a} )}  \ ,
\end{equation*}
where $\psi_{j}(.)$ denotes the $j^{th}$ atomic function, which is defined for all $\mathbf{a} \in \mathbb{R}^{p.k}$,  by
\begin{align} \label{eq:psi_atomic_func}
\psi_{j}( \mathbf{a} ) & = \Vert  \mathbf{P}^{\bot}_{\mathbf{F}_{j}(  \mathbf{a} ) }  \mathbf{x}_{j}  \Vert^{2}_{2}  \nonumber \\
 & = \Vert  \big(  \mathbf{I}_p - \mathbf{F}_{j}(  \mathbf{a} ) \mathbf{F}_{j}(  \mathbf{a} )^{+}\big) \mathbf{x}_{j} \Vert^{2}_{2}   \nonumber \\
 & = \Vert  \left(  \mathbf{I}_p - \left( \emph{diag}(\sqrt{\mathbf{W}}_{.j}) \mathbf{A} \right)  \left( \emph{diag}(\sqrt{\mathbf{W}}_{.j})\mathbf{A}  \right)^{+} \right) ( \sqrt{\mathbf{W}}_{.j} \odot \mathbf{X}_{.j} ) \Vert^{2}_{2}  \ .
\end{align}
Here  $\mathbf{P}^{\bot}_{\mathbf{F}_{j}(  \mathbf{a} ) }$ is the orthogonal projector onto the orthogonal complement of $\emph{ran}\big( \mathbf{F}_{j}(  \mathbf{a} ) \big)$ and $\mathbf{x}_{j}= \sqrt{\mathbf{W}}_{.j} \odot \mathbf{X}_{.j} $.

This new formulation of our NLLS problem, based on the cost function $\psi(.)$, suggests that the minimization of $\varphi^{*}(.)$ can be separated in two steps. Once a $\mathbf{A}$ matrix has been obtained by minimizing $\psi( \mathbf{a} )$, the $\mathbf{B}$  matrix can be obtained by solving a large block diagonal least-squares problem, which is equivalent to solve $n$ independent smaller linear least-squares problems. The rational for employing this separation of variables to minimize $\varphi^{*}(.)$ is given by the following theorem, which is a slight variation of a theorem originally proved by Golub and Pereyra in a more general setting (see Theorem  2.1 in~\cite{GP1973}).
\\
\begin{theo3.9} \label{theo3.9:box}
With the same notations and definitions as in Theorem~\ref{theo3.1:box}, the problem~\eqref{eq:P1} is equivalent to the problem
\begin{equation} \label{eq:VP1} \tag{VP1}
\min_{\mathbf{A}\in\mathbb{R}^{p \times k}}   \,  \psi \big( \emph{vec}(\mathbf{A}^{T}) \big) = \psi(\mathbf{a} ) = \frac{1}{2} \Vert \mathbf{P}^{\bot}_{\mathbf{F}(  \mathbf{a} ) }\mathbf{x} \Vert^{2}_{2}  \ ,
\end{equation}
where $\mathbf{a} = \emph{vec}(\mathbf{A}^{T}) \in \mathbb{R}^{p.k}$ and $\mathbf{x} = \emph{vec}( \sqrt{\mathbf{W}} \odot \mathbf{X})  \in \mathbb{R}^{p.n}$.
In other words, if we consider the range of $\varphi^{*}(.)$, $\text{C}_{\varphi^*}$, and the range of $\psi(.)$,
\begin{equation*}
\text{C}_{\psi} =  \big\{  \mathbf{y}\in\mathbb{R} \text{ }/\text{ }  \exists \mathbf{A}\in\mathbb{R}^{p \times k}  \text{ with } \mathbf{y}= \psi \big( \emph{vec}(\mathbf{A}^{T}) \big)  \big\}  \ ,
\end{equation*}
these two subsets of $\mathbb{R}$ have the same infimum and if this infimum is a minimum for one set, the other set also admits a minimum and these two minima are equal.
\end{theo3.9}
\begin{proof}
As in Theorem~\ref{theo3.1:box}, $\text{C}_{\varphi^*}$ and $\text{C}_{\psi}$ are bounded below by zero and, thus, admit an infimum greater or equal to zero, say $\bar{\mathbf{c}}_{\varphi^*}$ and $\bar{\mathbf{c}}_{\psi}$, respectively.
\\
Suppose first that $\bar{\mathbf{c}}_{\psi}  < \bar{\mathbf{c}}_{\varphi^*}$. Then, it exists $\mathbf{A} \in \mathbb{R}^{p \times k}$ such that
\begin{equation*}
\bar{\mathbf{c}}_{\psi} \le  \psi \big( \emph{vec}(\mathbf{A}^{T}) \big)  =  \psi ( a ) < \bar{\mathbf{c}}_{\varphi^*}  \ ,
\end{equation*}
where $\mathbf{a} = \emph{vec}(\mathbf{A}^{T}) \in \mathbb{R}^{p.k}$. Now, let $\mathbf{b} = \mathbf{F}(  \mathbf{a} )^{+} \mathbf{x } \in \mathbb{R}^{k.n}$ and define $\mathbf{B} = \emph{mat} ( \mathbf{b} ) \in \mathbb{R}^{k \times n}$, we have
\begin{equation*}
\mathbf{P}^{\bot}_{\mathbf{F}(  \mathbf{a} ) }\mathbf{x} =\left( \mathbf{I}_{p.n} - \mathbf{F}(  \mathbf{a} ) \mathbf{F}(  \mathbf{a} )^{+} \right) \mathbf{x} =  \mathbf{x} - \mathbf{F}(  \mathbf{a} ) \mathbf{b}
\end{equation*}
and
\begin{equation*}
 \Vert \mathbf{P}^{\bot}_{\mathbf{F}(  \mathbf{a} ) }\mathbf{x}  \Vert^{2}_{2} =  \Vert \mathbf{x} - \mathbf{F}(  \mathbf{a} ) \mathbf{b}  \Vert^{2}_{2} = \Vert  \sqrt{\mathbf{W}} \odot  ( \mathbf{X} - \mathbf{A}\mathbf{B}) \Vert^{2}_{F} \ ,
\end{equation*}
which implies that $\psi( \mathbf{a} ) = \varphi^{*} ( \mathbf{A}, \mathbf{B} )$. In other words, we have $\varphi^{*} ( \mathbf{A}, \mathbf{B} ) < \bar{\mathbf{c}}_{\varphi^*}$, which contradicts the assertion that $\bar{\mathbf{c}}_{\varphi^*}$ is the infimum of $\varphi^{*}(.)$. This shows that $\bar{\mathbf{c}}_{\psi} \ge  \bar{\mathbf{c}}_{\varphi^{*}}$.
\\
Suppose now that $ \bar{\mathbf{c}}_{\varphi^{*}} < \bar{\mathbf{c}}_{\psi} $. Then,  it exists $( \mathbf{A}, \mathbf{B} )  \in \mathbb{R}^{p \times k} \times \mathbb{R}^{k \times n}$ such that $\bar{\mathbf{c}}_{\varphi^{*}}  \le \varphi^{*} ( \mathbf{A}, \mathbf{B} ) < \bar{\mathbf{c}}_{\psi} $, otherwise  $\bar{\mathbf{c}}_{\varphi^{*}}$ is not the infimum of $ \varphi^{*}(.)$. However, if we define $ \widehat{\mathbf{b}} =  \mathbf{F}(  \mathbf{a} )^{+} \mathbf{x } \in \mathbb{R}^{k.n}$ and $\widehat{\mathbf{B}} = \emph{mat} ( \widehat{\mathbf{b}} ) $, we have
\begin{equation*}
 \Vert  \sqrt{\mathbf{W}} \odot  ( \mathbf{X} - \mathbf{A}\widehat{\mathbf{B}} ) \Vert^{2}_{F} \le \Vert  \sqrt{\mathbf{W}} \odot  ( \mathbf{X} - \mathbf{A}\mathbf{B} ) \Vert^{2}_{F}  \ ,
\end{equation*}
as, for a fixed  $\mathbf{A}$ matrix, $\widehat{\mathbf{b}} = \emph{vec}( \widehat{\mathbf{B}} )$  is the solution of the least-squares problem
\begin{equation*}
\min_{\mathbf{b}\in\mathbb{R}^{k.n}  } \,   \Vert  \mathbf{x} - \mathbf{F}(  \mathbf{a} ) \mathbf{b}  \Vert^{2}_{2} = \min_{\mathbf{B}\in\mathbb{R}^{k \times n}  } \,  \Vert \sqrt{\mathbf{W}} \odot ( \mathbf{X} - \mathbf{A}\mathbf{B} )  \Vert^{2}_{F}
\end{equation*}
This implies that $\psi( \mathbf{a} ) = \varphi^{*} ( \mathbf{A}, \widehat{\mathbf{B}} ) \le \varphi^{*} ( \mathbf{A}, \mathbf{B} ) < \bar{\mathbf{c}}_{\psi} $, which contradicts the assertion that $ \bar{\mathbf{c}}_{\psi} $ is the infimum of $\psi(.)$. This demonstrates that  $ \bar{\mathbf{c}}_{\varphi^{*}} \ge \bar{\mathbf{c}}_{\psi} $.
\\
Finally, the inequalities $ \bar{\mathbf{c}}_{\varphi^{*}} \ge \bar{\mathbf{c}}_{\psi} $ and $ \bar{\mathbf{c}}_{\varphi^{*}} \le \bar{\mathbf{c}}_{\psi} $ imply that $ \bar{\mathbf{c}}_{\varphi^{*}} = \bar{\mathbf{c}}_{\psi} $, which proves the first part of the theorem.
\\
\\
Now assume that $ \widehat{\mathbf{A}}$ minimizes $\psi(.)$, e.g., $ \psi ( \emph{vec}(\widehat{\mathbf{A}}^{T}) ) = \psi ( \widehat{\mathbf{a}} ) = \bar{\mathbf{c}}_{\psi} $. If we let $\widehat{\mathbf{b}} = \mathbf{F}(  \widehat{\mathbf{a}} )^{+} \mathbf{x}$ and $\widehat{\mathbf{B}} = \emph{mat} ( \widehat{\mathbf{b}} )$, we have
\begin{equation*}
 \psi ( \widehat{\mathbf{a}} ) =   \frac{1}{2} \Vert \mathbf{P}^{\bot}_{\mathbf{F}(  \widehat{\mathbf{a}} ) }\mathbf{x}  \Vert^{2}_{2} =  \frac{1}{2} \Vert  \mathbf{x} - \mathbf{F}(  \widehat{\mathbf{a}} ) \widehat{\mathbf{b}}  \Vert^{2}_{2} = \frac{1}{2} \Vert  \sqrt{\mathbf{W}} \odot  ( \mathbf{X} -  \widehat{\mathbf{A}} \widehat{\mathbf{B}}) \Vert^{2}_{F} = \varphi^{*} ( \widehat{\mathbf{A}}, \widehat{\mathbf{B}} )
\end{equation*}
and the equalities $ \bar{\mathbf{c}}_{\varphi^{*}} = \bar{\mathbf{c}}_{\psi} $ and $\psi ( \widehat{\mathbf{a}} ) = \bar{\mathbf{c}}_{\psi} $ show that  $\varphi^{*} ( \widehat{\mathbf{A}}, \widehat{\mathbf{B}} ) = \bar{\mathbf{c}}_{\varphi^{*}}$ and we conclude that $( \widehat{\mathbf{A}}, \widehat{\mathbf{B}} )$ is a global minimizer of $\varphi^{*}(.)$.
\\
Reciprocally, assume that $( \widehat{\mathbf{A}}, \widehat{\mathbf{B}} )$ minimizes $\varphi^{*}(.)$, e.g., $ \varphi^{*} ( \widehat{\mathbf{A}}, \widehat{\mathbf{B}} ) = \bar{\mathbf{c}}_{\varphi^{*}}$. Let $\widehat{\mathbf{a}} = \emph{mat} ( \widehat{\mathbf{A}}^{T}) $, $\bar{\mathbf{b}} = \mathbf{F}(  \widehat{\mathbf{a}} )^{+} \mathbf{x}$ and  $\bar{\mathbf{B}} = \emph{mat} ( \bar{\mathbf{b}} )$, we have
\begin{equation*}
\bar{\mathbf{c}}_{\psi} \le \psi ( \widehat{\mathbf{a}} ) =  \varphi^{*} ( \widehat{\mathbf{A}}, \bar{\mathbf{B}} ) \le \varphi^{*} ( \widehat{\mathbf{A}}, \widehat{\mathbf{B}} )  = \bar{\mathbf{c}}_{\varphi^*} = \bar{\mathbf{c}}_{\psi}  \ ,
\end{equation*}
which implies that  $\psi ( \widehat{\mathbf{a}} ) = \bar{\mathbf{c}}_{\psi} $ and $\widehat{\mathbf{a}}$ is a global minimizer of $\psi(.)$ and we are done.
\\
\end{proof}

An alternative to the minimization of  $\psi(.)$ in the variable projection approach can be introduced with the help of a QRCP (see equation~\eqref{eq:qrcp} in Subsection~\ref{lin_alg:box}) of the matrix $\mathbf{F}( \mathbf{a} )$ of rank $r  \le k.n$
\begin{equation*}
\mathbf{Q}( \mathbf{a} )\mathbf{F}( \mathbf{a} )\mathbf{P} =  \begin{bmatrix} \mathbf{R}    &  \mathbf{S}   \\  \mathbf{0}^{(p.n-r) \times r }  & \mathbf{0}^{(p.n-r) \times (k.n-r) }  \end{bmatrix}  \ ,
\end{equation*}
where $\mathbf{Q}( \mathbf{a} )$ is an $p.n \times p.n$ orthogonal matrix, $\mathbf{P}$ is an $k.n \times k.n$ permutation matrix, $\mathbf{R}$ is an $r \times r$ nonsingular upper triangular matrix and $\mathbf{S}$ an $r \times (k.n - r)$ full matrix. Then, if $\mathbf{Q}( \mathbf{a} )$ is partitioned into
\begin{equation*}
\mathbf{Q}( \mathbf{a} ) =  \begin{bmatrix} \mathbf{Q}_{1}( \mathbf{a} )   \\  \mathbf{Q}_{2}( \mathbf{a} ) \end{bmatrix}  \ ,
\end{equation*}
where $\mathbf{Q}_{1}( \mathbf{a} )$ and $\mathbf{Q}_{2}( \mathbf{a} )$ are, respectively, $r \times p.n$ and $(p.n - r) \times p.n$ submatrices, using results in Subsection~\ref{lin_alg:box}, we have
\begin{equation*}
\mathbf{Q}( \mathbf{a} ) \mathbf{P}_{\mathbf{F}( \mathbf{a} )}^{\bot} = \begin{bmatrix}   \mathbf{0}^{ r \times p.n }  \\  \mathbf{Q}_{2}( \mathbf{a} ) \end{bmatrix}  \ .
\end{equation*}
This implies that, for all  $\mathbf{a} \in \mathbb{R}^{k.p}$ and $\mathbf{x} = \emph{vec}( \sqrt{\mathbf{W}} \odot \mathbf{X})  \in \mathbb{R}^{p.n}$,
\begin{equation*}
 \Vert \mathbf{P}_{\mathbf{F}( \mathbf{a} )}^{\bot} \mathbf{x}  \Vert^{2}_{2} =   \Vert  \mathbf{Q}( \mathbf{a} ) \mathbf{P}_{\mathbf{F}( \mathbf{a} )}^{\bot} \mathbf{x}  \Vert^{2}_{2} =   \Vert   \mathbf{Q}_{2}( \mathbf{a} ) \mathbf{x} \Vert^{2}_{2}  \ ,
\end{equation*}
as first noted by Krogh~\cite{K1974} and Kaufman~\cite{K1975} and later by Shen and Ypma~\cite{SY2019}. Thus, instead of minimizing the variable projection functional $\psi(.)$ to solve the~\eqref{eq:VP1} problem, we can minimize the variable orthogonal functional $\psi^{*}(.)$ defined by
\begin{equation}  \label{eq:psi*_func}
\psi^{*}( \mathbf{a} ) =  \Vert   \mathbf{Q}_{2}( \mathbf{a} ) \mathbf{x} \Vert^{2}_{2} \ ,
\end{equation}
assuming that the rank of $\mathbf{F}( \mathbf{a} )$ stays constant in a neighborhood of a solution $\mathbf{\widehat{a}}$ of the~\eqref{eq:VP1} problem. Note that this condition is also implicit when using the variable projection functional $\psi(.)$ as this condition is required both for the differentiation of  $\mathbf{P}_{\mathbf{F}(.)}^{\bot}$ and $\mathbf{Q}_{2}(.)$ in a neighborhood of $\mathbf{\widehat{a}}$ as we will illustrate below. Thus, in the first step, it is mathematically equivalent to minimize $\psi(.)$ or $\psi^{*}(.)$, even though minimizing  $\psi^{*}(.)$ may involved slightly different numerical algorithms~\cite{K1975}\cite{L2009}\cite{SY2019}. Once a minimum of $\psi^{*}(.)$ has been determined, one can again determine $\mathbf{\widehat{b}}$ by solving the linear least-squares problem
\begin{equation*}
\min_{\mathbf{b}\in\mathbb{R}^{k.n}} \,  \Vert \mathbf{x} - \mathbf{F}(\mathbf{a})\mathbf{b} \Vert^{2}_{2} \ ,
\end{equation*}
according to Theorem~\ref{theo3.9:box}. 

As there is no constraint on the rank of $\mathbf{A} \in \mathbb{R}^{p \times k}$ in the~\eqref{eq:VP1} problem stated in Theorem~\ref{theo3.9:box}, the search space for minimizing the cost functions $\psi(.)$ or $\psi^{*}(.)$ is at first sight the linear space $\mathbb{R}^{p \times k}$. However, the following corollaries demonstrate that we can restrict this search space to the submanifold $\mathbb{R}^{p \times k}_k$ or even to $\mathbb{O}^{p \times k}$, the set of $p \times k$ matrices with orthonormal columns, which is called the Stiefel manifold~\cite{B2023}. Moreover, in many practical instances, for example in the matrix completion problem in which we are looking for a matrix  $\widehat{\mathbf{X}}$ of specified and fixed rank $k$, which agrees with the observed entries of the input matrix  $\mathbf{X}$, or when we solve the WLRA problem to estimate a consistent factor or principal component model, restricting the search space to the submanifold $\mathbb{R}^{p \times k}_k$ or even to the Stiefel submanifold is fully justified.
\\
\begin{corol3.1} \label{corol3.1:box}
With the same notations and definitions as in Theorem~\ref{theo3.9:box}, the problem~\eqref{eq:P1} is also equivalent to the following alternative formulations of the problem~\eqref{eq:VP1} in which the search space is restricted to $\mathbb{R}^{p \times k}_k$ and $\mathbb{O}^{p \times k}$, respectively:
\begin{equation*}
\min_{\mathbf{A} \in \mathbb{R}^{p \times k}_k }   \,  \psi \big( \emph{vec}(\mathbf{A}^{T}) \big) = \psi(\mathbf{a} ) = \frac{1}{2} \Vert \mathbf{P}^{\bot}_{\mathbf{F}(  \mathbf{a} ) }\mathbf{x} \Vert^{2}_{2}
\end{equation*}
and
\begin{equation*}
\min_{\mathbf{A} \in \mathbb{O}^{p \times k} }  \,  \psi \big( \emph{vec}(\mathbf{A}^{T}) \big) = \psi(\mathbf{a} ) = \frac{1}{2} \Vert \mathbf{P}^{\bot}_{\mathbf{F}(  \mathbf{a} ) }\mathbf{x} \Vert^{2}_{2}  \ ,
\end{equation*}
where $\mathbf{a} = \emph{vec}(\mathbf{A}^{T}) \in \mathbb{R}^{p.k}$ and $\mathbf{x} = \emph{vec}( \sqrt{\mathbf{W}} \odot \mathbf{X})  \in \mathbb{R}^{p.n}$.
\end{corol3.1}
\begin{proof}
As in Remark~\ref{remark3.1:box}, we note that for all $(\mathbf{A}, \mathbf{B}) \in \mathbb{R}^{p \times k} \times \mathbb{R}^{k \times n} $, the matrix product $\mathbf{A}\mathbf{B}$ is of rank at most $k$ and can be factored as $\mathbf{A}\mathbf{B} = \mathbf{C}\mathbf{D}$ with $(\mathbf{C}, \mathbf{D}) \in \mathbb{R}^{p \times k}_k \times \mathbb{R}^{k \times n} $, which implies that the range of $\varphi^{*}(.)$, $\text{C}_{\varphi^{*}}$, is equal to the set  $\big\{  y \in\mathbb{R} \text{ }/\text{ }  \exists \mathbf{C}\in\mathbb{R}^{p \times k}_k \text{, }\exists \mathbf{D}\in\mathbb{R}^{k \times n} \text{ with } y = \varphi^{*}( \mathbf{C},\mathbf{D}  ) \big\}$. Taking into account this property, it is easy to verify that  a slight modification of the demonstration of Theorem~\ref{theo3.9:box} leads to the assertion that  the problem~\eqref{eq:P1} is also equivalent to the problem
\begin{equation*}
\min_{\mathbf{A} \in \mathbb{R}^{p \times k}_k }   \,  \psi \big( \emph{vec}(\mathbf{A}^{T}) \big) = \psi( \mathbf{a} ) = \frac{1}{2} \Vert \mathbf{P}^{\bot}_{\mathbf{F}(  \mathbf{a} ) }\mathbf{x} \Vert^{2}_{2}  \ .
\end{equation*}
We omit the details.
\\
Now, since $\mathbb{O}^{p \times k} \subset  \mathbb{R}^{p \times k}_k $ and any element $\mathbf{A}$ of $\mathbb{R}^{p \times k}_k $ can also be written as $\mathbf{A} = \mathbf{U}\mathbf{R}$ where $\mathbf{U} \in \mathbb{O}^{p \times k}  $ and $\mathbf{R} \in \mathbb{R}^{k \times k}_k $, for example by using the QR or SVD decompositions of $\mathbf{A}$, we have 
\begin{equation*}
\psi ( \mathbf{a} ) =  \varphi^{*} ( \mathbf{A}, \widehat{\mathbf{B}} ) = \varphi^{*} ( \mathbf{U}\mathbf{R}, \widehat{\mathbf{B}} ) = \varphi^{*} ( \mathbf{U}, \mathbf{R}\widehat{\mathbf{B}} )  \ge \psi ( \mathbf{u} )  \ ,
\end{equation*}
where $ \mathbf{a} = \emph{vec}( \mathbf{A}^{T} )$, $\widehat{\mathbf{B}} = \emph{mat}( \widehat{\mathbf{b}} )$ with $\widehat{\mathbf{b}} = \mathbf{F}(  \mathbf{a} )^{+} \mathbf{x} $ and $\mathbf{u} = \emph{vec}( \mathbf{U}^{T} )$. Reciprocally, if $\mathbf{A} = \mathbf{U}\mathbf{R}$  with $\mathbf{U} \in \mathbb{O}^{p \times k}  $ and $\mathbf{R} \in \mathbb{R}^{k \times k}_k $, we have
\begin{equation*}
\psi ( \mathbf{u} ) =  \varphi^{*} ( \mathbf{U}, \widehat{\mathbf{D}} )  =  \varphi^{*} ( \mathbf{A} \mathbf{R}^{-1}, \widehat{\mathbf{D}} ) = \varphi^{*} ( \mathbf{A}, \mathbf{R}^{-1} \widehat{\mathbf{D}} ) \ge \psi ( \mathbf{a} )  \ ,
\end{equation*}
where $\widehat{\mathbf{D}} = \emph{mat} ( \widehat{\mathbf{d}} )$ with $\widehat{\mathbf{d}} = \mathbf{F}(  \mathbf{u} )^{+} \mathbf{x} $. This implies that $\psi ( \mathbf{a} ) = \psi ( \mathbf{u} )$ if  $\mathbf{A} = \mathbf{U}\mathbf{R}$, which demonstrates that the range of $\psi(.)$, $\text{C}_{\psi}$, is equal to the set
\begin{equation*}
\big\{  y \in\mathbb{R} \text{ }/\text{ }  \exists \mathbf{U}\in \mathbb{O}^{p \times k} \text{ with } y = \psi \big( \emph{vec}(\mathbf{U}^{T}) \big)  \big\}  \ .
\end{equation*}
As a consequence, these two sets have the same infimum and the same minimum, if this minimum exists, and the problems~\eqref{eq:P1} or~\eqref{eq:VP1} are also equivalent to 
\begin{equation*}
\min_{\mathbf{A} \in \mathbb{O}^{p \times k} }  \,  \psi \big( \emph{vec}(\mathbf{A}^{T}) \big) = \psi(\mathbf{a} ) = \frac{1}{2} \Vert \mathbf{P}^{\bot}_{\mathbf{F}(  \mathbf{a} ) }\mathbf{x} \Vert^{2}_{2}  \ ,
\end{equation*}
as claimed in the Corollary.
\\
\end{proof}

\begin{corol3.2} \label{corol3.2:box}
With the same definitions and notations as in Theorem~\ref{theo3.9:box} and Corollary~\ref{corol3.1:box}, the following two assertions are true:
\begin{align*}
1) & \quad\  \psi ( \emph{vec}( \mathbf{A}^{T} )  ) =  \psi ( \emph{vec}( \mathbf{C}^{T} )  )  \text{ if }  \mathbf{A} = \mathbf{C} \mathbf{D}   \text{ with }  \mathbf{A}  \in \mathbb{R}^{p \times k}_k, \mathbf{C}  \in \mathbb{R}^{p \times k}_k   \text{ and }  \mathbf{D}  \in \mathbb{R}^{k \times k}_k  \  , \\
2) &  \quad\  \psi ( \emph{vec}( \mathbf{A}^{T} )  ) =  \psi ( \emph{vec}( \mathbf{C}^{T} )  )  \text{ if }  \mathbf{A} = \mathbf{C} \mathbf{D}   \text{ with }  \mathbf{A}  \in \mathbb{O}^{p \times k}, \mathbf{C}  \in \mathbb{O}^{p \times k}   \text{ and }  \mathbf{D}  \in \mathbb{O}^{k \times k} \ .
\end{align*}
\end{corol3.2}
\begin{proof}
The proof is very similar to the one used in Corollary~\ref{corol3.1:box}. Suppose first that $\mathbf{A} = \mathbf{C} \mathbf{D}$ with $\mathbf{A}  \in \mathbb{R}^{p \times k}_k, \mathbf{C}  \in \mathbb{R}^{p \times k}_k   \text{ and }  \mathbf{D}  \in \mathbb{R}^{k \times k}_k$ and let $\mathbf{a} = \emph{vec}( \mathbf{A}^{T} )$ and $\mathbf{c} = \emph{vec}( \mathbf{C}^{T} )$. We have
\begin{equation*}
\psi ( \mathbf{a} ) =  \varphi^{*} ( \mathbf{A}, \widehat{\mathbf{B}} ) = \varphi^{*} ( \mathbf{C}\mathbf{D}, \widehat{\mathbf{B}} ) = \varphi^{*} ( \mathbf{C}, \mathbf{D}\widehat{\mathbf{B}} )  \ge \psi ( \mathbf{c} )  \ ,
\end{equation*}
where $\widehat{\mathbf{B}} = \emph{mat}( \widehat{\mathbf{b}} )$ and $\widehat{\mathbf{b}} = \mathbf{F}(  \mathbf{a} )^{+} \mathbf{x} $. Reciprocally, we have
\begin{equation*}
\psi ( \mathbf{c} ) =  \varphi^{*} ( \mathbf{C}, \widehat{\mathbf{E}} )  =  \varphi^{*} ( \mathbf{A} \mathbf{D}^{-1}, \widehat{\mathbf{E}} ) = \varphi^{*} ( \mathbf{A}, \mathbf{D}^{-1} \widehat{\mathbf{E}} ) \ge \psi ( \mathbf{a} )  \ ,
\end{equation*}
where $\widehat{\mathbf{E}} = \emph{mat} ( \widehat{\mathbf{e}} )$ and $\widehat{\mathbf{e}} = \mathbf{F}(  \mathbf{c} )^{+} \mathbf{x} $. Finally, we obtain $\psi ( \emph{vec}( \mathbf{A}^{T} )  ) =  \psi ( \emph{vec}( \mathbf{C}^{T} )  ) $, and the first part of the corollary is demonstrated.
\\
\\
The proof of the second assertion is exactly similar and thus omitted.
\\
\end{proof}
\begin{remark3.7} \label{remark3.7:box}
Theorem~\ref{theo3.9:box} and Corollaries~\ref{corol3.1:box} and~\ref{corol3.2:box} also explain why the WLRA problem can be recast as an optimization problem on the Grassmann manifold $\text{Gr}(p,k)$~\cite{DKM2012}\cite{BA2015}, the collection of all linear subspaces of fixed dimension $k$ of the Euclidean space $\mathbb{R}^p$, which is a smooth (quotient) manifold of dimension $k.(p-k)$~\cite{AMS2008}\cite{B2023}. As stated in Theorems~\ref{theo3.1:box} and~\ref{theo3.9:box}, the formulations~\eqref{eq:P0},~\eqref{eq:P1} and~\eqref{eq:VP1} of the WLRA problem are equivalent and Corollary~\ref{corol3.1:box} shows that the search space for the~\eqref{eq:VP1} problem can be restricted to $\mathbb{R}^{p \times k}_k$ or $\mathbb{O}^{p \times k}$ without loss of generality. Next, Corollary~\ref{corol3.2:box} demonstrates that $\psi ( \emph{vec}( \mathbf{A}^{T} )  ) =  \psi ( \emph{vec}( \mathbf{C}^{T} )  )$ as soon as the linear subspaces $\emph{ran}( \mathbf{A} )$ and $\emph{ran}( \mathbf{C} )$ are equal. In other words, the value  of $\psi ( \emph{vec}( \mathbf{A}^{T} )  )$ for $\mathbf{A}  \in \mathbb{R}^{p \times k}_k$ depends only on the linear subspace  $\emph{ran}( \mathbf{A} )$ as for any matrix $\mathbf{C}  \in \mathbb{R}^{p \times k}_k$ such that the columns of $\mathbf{C}$ is a  (orthonormal or not) basis of $\emph{ran}( \mathbf{A} )$, we have $\psi ( \emph{vec}( \mathbf{A}^{T} )  ) =  \psi ( \emph{vec}( \mathbf{C}^{T} )  )$.
\\
As such, an element $\mathcal{U}$ of $\text{Gr}(p,k)$ can be represented by any $p  \times k$ matrix $\mathbf{U}$ of full column-rank such that $\emph{ran}(\mathbf{U}) = \mathcal{U}$, e.g., if the columns of $\mathbf{U}$ form a basis of $\mathcal{U}$. For numerical reasons, elements of $\text{Gr}(p,k)$ are very often represented by elements of $\mathbb{O}^{p \times k}$~\cite{EAS1998}\cite{MMH2003} \cite{C2008b}\cite{DMK2011}\cite{DKM2012}\cite{BA2015}, but any matrix with the same column space can be used, and we will see below that representing $\mathcal{U}$ by elements of $\mathbb{R}^{p \times k}_k$ instead of $\mathbb{O}^{p \times k}$ can also be useful to demonstrate important properties of the cost function $\psi(.)$ used in the~\eqref{eq:VP1} form of the WLRA problem, especially when the associated weight matrix $\mathbf{W}$ is not strictly positive.
\\
Stated differently, we can say that two $p \times k$ (orthonormal) matrices $\mathbf{U}$ and $\mathbf{V}$ are equivalent if and only if they have the same column space or, equivalently, if it exists some $\mathbf{Q} \in \mathbb{R}^{k \times k}_k$ (or $\mathbf{Q} \in \mathbb{O}^{k \times k}$) such that $\mathbf{U}=\mathbf{V}\mathbf{Q}$. Because of this equivalence relationship in  $\mathbb{R}^{p \times k}_k$, $\text{Gr}(p,k)$ can be described as the quotient of $\mathbb{R}^{p \times k}_k$ by the action of $\mathbb{R}^{k \times k}_k$. Alternatively, $\text{Gr}(p,k)$ can also be described as the quotient of $\mathbb{O}^{p \times k}$ by the action of $\mathbb{O}^{k \times k}$, see Subsection~\ref{calculus:box} for details. See~\cite{B2023}\cite{BZA2023} for a geometrical and comprehensive description of these different approaches of representing elements of $\text{Gr}(p,k)$ with matrices.
\\
Moreover, the fact that the formulation~\eqref{eq:P1} of the WLRA problem never has a unique or finite set of global minimizers is also related to the preceding discussion. If $\mathbf{Y}=\mathbf{A}\mathbf{B}$ with $\mathbf{A} \in \mathbb{R}^{p \times k}_k$, the elements of the columns of $\mathbf{B} \in \mathbb{R}^{k \times p}$ are the coordinates of the corresponding columns of $\mathbf{Y}$ in the particular basis of $\emph{ran}(\mathbf{A})$ provided by the columns of $\mathbf{A}$ and, thus, these coordinates depend on the choice of the basis.
\\
Finally, in a similar fashion that the formulation~\eqref{eq:VP1} and the associated variable projection functionals $\psi(.)$ and $\psi^{*}(.)$ are derived from the formulation~\eqref{eq:P1} of the WLRA problem, it is also possible to reformulate the problem~\eqref{eq:P2}  (see Remark~\ref{remark3.2:box}) as a double-minimization problem
 \begin{equation*}
\min_{\mathbf{N}\in\mathbb{R}^{p \times (p-k)}_{p-k}  } \, \quad\  \left(   \min_{\mathbf{Y}\in\mathbb{R}^{p \times n} \text{ with } \mathbf{N}^{T}\mathbf{Y} =  \mathbf{0}^{(p-k) \times n} }  \,  \frac{1}{2}  \Vert \sqrt{\mathbf{W}} \odot ( \mathbf{X} - \mathbf{Y} )  \Vert^{2}_{F} \right) \ ,
\end{equation*}
which will lead to a dual formulation~\eqref{eq:VP2} of the WLRA problem and its associated variable projection functional $\psi^{**}(.)$ on the Grassmann manifold $\text{Gr}(p,p-k)$ ~\cite{MMH2003}. More precisely, Manton et al.~\cite{MMH2003} have demonstrated that the above inner minimization problem has a closed form solution (see Theorem 1 in~\cite{MMH2003}), which can be calculated analytically and depends only on the range of $\mathbf{N}$ and not on the particular matrix $\mathbf{N}$ when $\mathbf{W}\in\mathbb{R}^{p \times n}_{+*}$ (and these results can be probably extended to the case $\mathbf{W}\in\mathbb{R}^{p \times n}_{+}$). Thus, problem~\eqref{eq:P2} is also a separable NLLS problem as stated in the Definition~\ref{def2.10:box} of Subsection~\ref{calculus:box} despite the variable matrix $\mathbf{Y}$ does not  occur linearly in the residual function associated with problem~\eqref{eq:P2}. This situation is exactly similar to the one described above for the problem~\eqref{eq:VP1} and its associated variable projection functional $\psi(.)$ where we proceed in two steps, namely, first find the matrix $\mathbf{\widehat{A}}$ such that $\psi(\mathbf{\widehat{a}})$ is minimum and, second, determine $\mathbf{\widehat{B}}$ by solving a large block diagonal least-squares problem. In other words, for $\mathbf{N}\in\mathbb{R}^{p \times (p-k)}_{p-k}$, the dual variable projection functional
\begin{equation}\label{eq:psi**_func}
\psi^{**}(\mathbf{N}) = \min_{\mathbf{Y}\in\mathbb{R}^{p \times n} \text{ with } \mathbf{N}^{T}\mathbf{Y} =  \mathbf{0}^{(p-k) \times n} }  \,  \frac{1}{2}  \Vert \sqrt{\mathbf{W}} \odot ( \mathbf{X} - \mathbf{Y} )  \Vert^{2}_{F}
\end{equation}
is well defined and we can attempt to find a solution $\mathbf{\widehat{N}}$ (more precisely a subspace $\mathcal{\widehat{N}}$, which is represented by  $\mathbf{\widehat{N}}$) of the problem
\begin{equation} \label{eq:VP2} \tag{VP2}
\min_{\mathbf{N}\in\mathbb{R}^{p \times (p-k)}_{p-k}  }  \, \quad\  \psi^{**}(\mathbf{N})
\end{equation}
in a first step by any iterative first- or second-order method working on the Grassmann manifold $\text{Gr}(p,p-k)$~\cite{EAS1998}\cite{MMH2003}. Once a minimum $\mathbf{\widehat{N}}$ has been found, the best rank-$k$ approximation matrix $\mathbf{\widehat{Y}}$ can be determined in a second step by solving the inner minimization problem analytically for the matrix $\mathbf{\widehat{N}}$. As for the~\eqref{eq:VP1} problem, the search space for the~\eqref{eq:VP2} problem can be restricted to $\mathbb{O}^{p \times (p-k)}$ without loss of generality for numerical reasons~\cite{EAS1998}\cite{MMH2003}. Furthermore, it will be shown in the following sections, that minimizing the cost function $\psi(.)$ associated with the~\eqref{eq:VP1} problem reduces to one of dimension $k.(p-k)$ instead of $k.p$ as for the minimization of the cost function $\psi^{**}(.)$ associated with the~\eqref{eq:VP2} problem~\cite{EAS1998}\cite{MMH2003}, highlighting again the duality between the~\eqref{eq:VP1} and~\eqref{eq:VP2} formulations of the WLRA problem. $\blacksquare$
 \\
\end{remark3.7}

We address now the question of the continuity of the cost function $\psi(.)$, which must be minimized in the~\eqref{eq:VP1} problem, as was done above for the cost functions $\varphi(.)$ and $\varphi^{*}(.)$ associated, respectively, with the formulations~\eqref{eq:P0} and~\eqref{eq:P1} of the WLRA problem. Taking into account that, for all $\mathbf{a} \in \mathbb{R}^{p.k}$, we have
\begin{equation}  \label{eq:psi_func2}
\psi(\mathbf{a} ) = \frac{1}{2} \Vert \mathbf{P}^{\bot}_{\mathbf{F}(  \mathbf{a} ) }\mathbf{x} \Vert^{2}_{2} \ ,
\end{equation}
we see that the continuity of $\psi(.)$ is closely associated to the continuity of the orthogonal projector $\mathbf{P}^{\bot}_{\mathbf{F}(.)}$ (or equivalently $\mathbf{P}_{\mathbf{F}(.)}$) as a function of $\mathbf{a} \in \mathbb{R}^{p.k}$. Furthermore, due to the block diagonal structure of $\mathbf{F}(.)$, we observe that the continuity of $\mathbf{P}^{\bot}_{\mathbf{F}(.)}$ as a function of $\mathbf{a}$ is equivalent to the continuity of the $n$ atomic orthogonal projectors $\mathbf{P}^{\bot}_{\mathbf{F}_{j}(.) }$, for $j=1  \text{ to } n$, since, for all $\mathbf{a} \in \mathbb{R}^{p.k}$, we have
\begin{equation*}
\mathbf{P}^{\bot}_{\mathbf{F}(  \mathbf{a} )} =  \bigoplus_{j=1}^n \mathbf{P}^{\bot}_{\mathbf{F}_{j}(  \mathbf{a} ) } \ .
\end{equation*}

The next theorem gives necessary and sufficient conditions for the continuity of a general orthogonal projector $\mathbf{P}_{\Phi( .)}$, which is associated with a $l \times t$ matrix function $\Phi(.)$ of a vector $\mathbf{a} \in  \mathbb{R}^{m}$, but let us first give the following definition:
 \\
\begin{def3.1} \label{def3.1:box}
Let $\Phi(.)$ be a matrix function defined as
\begin{equation*}
\Phi : \mathbb{R}^m  \longrightarrow \mathbb{R}^{l \times t}  : \mathbf{a}  \mapsto \Phi(  \mathbf{a}  )  \ .
\end{equation*}
We say that the matrix function $\Phi(.)$  has a local constant rank at a point $\mathbf{a}_0 \in  \mathbb{R}^m$ if there exists an open neighborhood $\Upsilon$ of $\mathbf{a}_0$ in $\mathbb{R}^m$ such that the matrix $\Phi(  \mathbf{a} )$ has a constant rank $q \le \text{min}(l,t)$ for all $\mathbf{a} \in \Upsilon$.
\end{def3.1}
We then restate the following fundamental result about  the continuity of the Moore-Penrose inverse and orthogonal projectors for real matrix functions, which can be traced back to the seminal works of Wedin~\cite{W1973} and Stewart~\cite{SS1990}.
\\
\begin{theo3.10} \label{theo3.10:box}
Let $\Phi(.)$ be a matrix function : $\mathbb{R}^m \longrightarrow \mathbb{R}^{l \times t}$, which is continuous at a point $\mathbf{a}_0 \in  \mathbb{R}^m$. The following conditions are equivalent.
\begin{align*}
1) &  \quad\  \Phi(.)  \text{ has a local constant rank at }  \mathbf{a}_0   \\
2) &  \quad\  \Phi(.)^{+}   \text{ is continuous at }  \mathbf{a}_0   \\
3) &  \quad\  \Phi(.)\Phi(.)^{+} =  \mathbf{P}_{\Phi(.)}  \text{ is continuous at }  \mathbf{a}_0   \\
4) &  \quad\   \Phi(.)^{+}\Phi(.)  =  \mathbf{P}_{\Phi(.)^{T}}   \text{ is continuous at }  \mathbf{a}_0
\end{align*}
In other words, the continuity of the pseudo-inverse of a continuous matrix function $\Phi(.)$  at a point $\mathbf{a}_0 \in \mathbb{R}^m$ is equivalent to the continuity of the orthogonal projectors onto the column or row spaces of this matrix function at $\mathbf{a}_0$ and all these conditions are equivalent to the assertion that $\Phi(.)$ has local constant rank at $\mathbf{a}_0$ if $\Phi(.)$ is itself continuous at $\mathbf{a}_0$.
\end{theo3.10}
\begin{proof}
See Propositions 8.1 and 8.2 in Chapter 8 of Magnus and Neudecker~\cite{MN2019} or Chapter 10 of Campbell and Meyer~\cite{CM2009}.

\end{proof}

To ease the notation burden in the rest of this section and the following sections, we define the following linear mapping $h(.)$ and its inverse mapping  $h^{-1}(.)$ :
\begin{align} \label{eq:h_func}
h  & : \mathbb{R}^{p.k}  \longrightarrow \mathbb{R}^{p \times k}, \mathbf{a} \mapsto \lbrack \emph{mat}_{k  \times p}(  \mathbf{a}  ) \rbrack^{T}  = \mathbf{A} \ ,   \\
h^{-1}  & : \mathbb{R}^{p \times k} \longrightarrow \mathbb{R}^{p.k} , \mathbf{A} \mapsto \emph{vec}(  \mathbf{A}^{T} ) = \mathbf{a} \ ,  \nonumber 
\end{align}
which are homeomorphisms from $\mathbb{R}^{p.k}$ to $\mathbb{R}^{p \times k}$  and from $\mathbb{R}^{p \times k}$ to $\mathbb{R}^{p.k}$, respectively, allowing to identify the elements and the topologies of these two finite dimensional vector spaces. These notations seem cumbersome, but are related to our definition of the vectorized form  of the matrix $\mathbf{A}$ as $\emph{vec}(  \mathbf{A}^{T} )$ instead of $\emph{vec}(  \mathbf{A} )$, which will be justified in the next sections. In other words, with these definitions, we have $h^{-1} (\mathbf{A} ) = \emph{vec}(  \mathbf{A}^{T} ) = \mathbf{a}$ for all $\mathbf{A} \in \mathbb{R}^{p \times k}$.

Armed with Theorem~\ref{theo3.10:box}, we now consider the continuity of the orthogonal projector $\mathbf{P}^{\bot}_{\mathbf{F}(.)}$ (or equivalently $\mathbf{P}_{\mathbf{F}(.)}$), which is used in the cost function $\psi(.)$ of the~\eqref{eq:VP1} problem. We first observe that $\mathbf{F}(.)$, as a function of $\mathbf{a}$, is a  linear mapping from $\mathbb{R}^{p.k}$ to $\mathbb{R}^{p.n \times n.k}$ since the $\emph{mat}$ and transpose operators are linear mappings, and, the Kronecker and matrix products are bilinear mappings. As $\mathbb{R}^{p.k}$ and $\mathbb{R}^{p.n \times n.k}$ are finite dimensional vector spaces, $\mathbf{F}(.)$ is thus continuous for all  $\mathbf{a} \in \mathbb{R}^{p.k}$. Similarly, the $n$ atomic matrix functions $\mathbf{F}_{j}(.)$ are also continuous linear mappings from $\mathbb{R}^{p.k}$ to $\mathbb{R}^{p \times k}$. In these conditions, Theorem~\ref{theo3.9:box} shows that the continuity of the orthogonal projectors $\mathbf{P}^{\bot}_{\mathbf{F}(.)}$ and $\mathbf{P}^{\bot}_{\mathbf{F}_{j}(.)}$ (or equivalently $\mathbf{P}_{\mathbf{F}(.)}$ and $\mathbf{P}_{\mathbf{F}_{j}(.)}$) at a point $\mathbf{a} \in \mathbb{R}^{p.k}$ is equivalent, respectively, to the propositions that $\mathbf{F}(.)$ and $\mathbf{F}_{j}(.)$  have a local constant rank at $\mathbf{a}$. Furthermore, the proposition that $\mathbf{F}(.)$ has a local constant rank at $\mathbf{a}$ is equivalent to the proposition that all the $n$ $\mathbf{F}_{j}(.)$ functions,  for $j=1  \text{ to } n$, have also a local constant rank at $\mathbf{a}$.

In the special case where the weight matrix $\mathbf{W} \in \mathbb{R}^{p \times n}_{+*}$  (e.g., $\mathbf{W}_{ij} > 0$), we have the following more precise result:
\begin{theo3.11} \label{theo3.11:box}
For $\mathbf{X} \in \mathbb{R}^{p \times n} \text{ and } \mathbf{W} \in \mathbb{R}^{p \times n}_{+*}$, and any fixed integer $k \le \emph{rank}( \mathbf{X} ) \le \text{min}( {p},{n} )$, define the matrix function $\mathbf{P}_{\mathbf{F}(.)}$, from $\mathbb{R}^{p.k}$ to $\mathbb{R}^{p.n \times p.n}$, by
\begin{equation*}
\mathbf{a} \mapsto \mathbf{P}_{\mathbf{F}(  \mathbf{a}  )} = \mathbf{F}(  \mathbf{a} ) \mathbf{F}(  \mathbf{a} )^{+} \ ,
\end{equation*}
where $\mathbf{F}(  \mathbf{a} )^{+}$ is the pseudo-inverse of $\mathbf{F}(  \mathbf{a} )$ and $\mathbf{F}(  \mathbf{a} )$ is the $ p.n \times n.k  $ block diagonal matrix 
\begin{equation*}
\mathbf{F}(  \mathbf{a} ) =  \bigoplus_{j=1}^n \emph{diag}(\sqrt{\mathbf{W}}_{.j}) h(  \mathbf{a}  ) = \bigoplus_{j=1}^n \emph{diag}(\sqrt{\mathbf{W}}_{.j})\mathbf{A}  \ .
\end{equation*}
$\mathbf{P}_{\mathbf{F}(.)}$ is continuous at all $\mathbf{a} \in h^{-1}( \mathbb{R}^{p \times k}_k )$ and discontinuous at all $\mathbf{a} \in h^{-1}( \mathbb{R}^{p \times k}_{<k} )$.

\end{theo3.11}
\begin{proof}
As noted above, the continuity of $\mathbf{P}_{\mathbf{F}(.)}$ at a point $\mathbf{a}_0 \in \mathbb{R}^{p.k}$ is equivalent to the existence of an open neighborhood of $\mathbf{a}_0$  in $\mathbb{R}^{p.k}$ in which $\mathbf{F}(.)$ has a local constant rank.
However, since
\begin{align*}
\mathbf{F}(  \mathbf{a} )  & =    \bigoplus_{j=1}^n \emph{diag}(\sqrt{\mathbf{W}}_{.j})\mathbf{A}  \\
  & =  \emph{diag}\big( \emph{vec}( \sqrt{\mathbf{W}} ) \big)  \big[  \bigoplus_{j=1}^n \mathbf{A}  \big]
\end{align*}
and $\emph{rank}\big( \emph{diag}( \emph{vec}( \sqrt{\mathbf{W}} ) ) \big) = p.n$ if $\mathbf{W}$  is strictly positive, then
\begin{equation*}
\emph{rank}\big( \mathbf{F}(  \mathbf{a} ) \big) = \emph{rank} ( \bigoplus_{j=1}^n \mathbf{A} ) = n.\emph{rank} ( \mathbf{A} )  \ .
\end{equation*}
Thus, the rank of $\mathbf{F}(  \mathbf{a}  )$ is entirely determined by the rank of $\mathbf{A} = h(  \mathbf{a}  )$ when the weight matrix is strictly positive. In other words, the continuity of $\mathbf{P}_{\mathbf{F}(.)}$ at $\mathbf{a}_0 \in \mathbb{R}^{p.k}$ is equivalent to the existence of an open neighborhood $\Upsilon$ of $\mathbf{A}_0 = h(  \mathbf{a}_0  ) \in \mathbb{R}^{p \times k} $ such that for all $\mathbf{A} \in \Upsilon$, the rank of $\mathbf{A} = h(  \mathbf{a}  )$ is constant.

Now, we have $\mathbb{R}^{p \times k} = \mathbb{R}^{p \times k}_k \bigcup \mathbb{R}^{p \times k}_{  <k}$ and let us consider separately the two cases $\mathbf{A}_0 \in \mathbb{R}^{p \times k}_k$ and $\mathbf{A}_0 \in \mathbb{R}^{p \times k}_{  <k}$. 

Suppose first that $\mathbf{A}_0  \in \mathbb{R}^{p \times k}_k$. Per definition, the rank of $\mathbf{A}_0 $ is constant and equal to $k$. Note further that $\mathbb{R}^{p \times k}_k$ is an open set of $\mathbb{R}^{p \times k}$ as the preimage of the open set $\mathbb{R}  \backslash  \lbrace 0 \rbrace$ under the continuous mapping $\mathbf{A} \mapsto det( \mathbf{A}^{T}\mathbf{A} )$ as stated in Theorem~\ref{theo2.3:box}. In other words, for all $\mathbf{A}_0  \in \mathbb{R}^{p \times k}_k$, there is an open neighborhood $\Upsilon$ of $\mathbf{A}_0$ included in $\mathbb{R}^{p \times k}_k$ and we deduce immediately that $\mathbf{P}_{\mathbf{F}(.)}$ is a continuous mapping for all $\mathbf{a}_0 \in h^{-1} ( \mathbb{R}^{p \times k}_k )$ using Theorem~\ref{theo3.10:box}.

Suppose now that $\mathbf{A}_0  \in \mathbb{R}^{p \times k}_{<k}$. Since $\mathbb{R}^{p \times k}_{<k}$ is the frontier of the open set $\mathbb{R}^{p \times k}_k$ in $\mathbb{R}^{p \times k}$ according to Theorem~\ref{theo2.3:box}, every open  neighborhood $\Upsilon$ of $\mathbf{A}_0$ in $\mathbb{R}^{p \times k}_{<k}$ also contains some points $\mathbf{A} \in \mathbb{R}^{p \times k}_k$ and as such if $\mathbf{A}_0  \in \mathbb{R}^{p \times k}_{<k}$ there is no open neighborhood $\Upsilon$ of $\mathbf{A}_0$ in $\mathbb{R}^{p \times k}$ in which the rank of $\mathbf{A}$ is constant for all $\mathbf{A} \in \Upsilon$. This implies that $\mathbf{P}_{\mathbf{F}(.)}$ (and, thus, also $\mathbf{P}^{\bot}_{\mathbf{F}(.)}$ and $\psi(.)$) is discontinuous at all points $\mathbf{a}_0  \in \mathbb{R}^{p.k}$ such that $h( \mathbf{a}_0 ) =  \mathbf{A}_0  \in \mathbb{R}^{p \times k}_{<k}$.

\end{proof}

\begin{corol3.3} \label{corol3.3:box}
With the same definitions and notations as in Theorem~\ref{theo3.11:box}, the cost function $\psi(.)$ of the~\eqref{eq:VP1} problem
\begin{equation*}
\psi(\mathbf{a} ) = \frac{1}{2} \Vert \mathbf{P}^{\bot}_{\mathbf{F}(  \mathbf{a} ) }\mathbf{x} \Vert^{2}_{2}
\end{equation*}
is continuous at all points $\mathbf{a}  \in h^{-1}(  \mathbb{R}^{p \times k}_k )$, if the weight matrix $\mathbf{W} \in \mathbb{R}^{p \times n}_{+*}$ (e.g.,  if $\mathbf{W}$ is strictly positive). Furthermore, the sets of global minimizers of $\psi(.)$ over the subsets $\mathbb{R}^{p \times k}_k$  and $\mathbb{O}^{p \times k}$ of $\mathbb{R}^{p \times k}$ are not empty if $\mathbf{W} \in \mathbb{R}^{p \times n}_{+*}$.
\end{corol3.3}
\begin{proof}
From Theorem~\ref{theo3.11:box} above, we know that $\mathbf{P}^{\bot}_{\mathbf{F}(.)}$ is a continuous mapping over the set $\mathbb{R}^{p \times k}_k$ if the weight matrix $\mathbf{W}$ is strictly positive. Thus, in the same conditions, the restriction of $\psi(.)$ to  $\mathbb{R}^{p \times k}_k$  is the composition of several continuous mappings on their respective domains of definition and, consequently, the restriction of $\psi(.)$ to $\mathbb{R}^{p \times k}_k$ is continuous at all points in $\mathbb{R}^{p \times k}_k$.

The second part of the corollary results immediately from Theorems~\ref{theo3.3:box},~\ref{theo3.9:box} and Corollary~\ref{corol3.1:box}, which show, respectively, that the set of solutions of the (WLRA) problem is not empty if the weight matrix is strictly positive and that the (WLRA) and~\eqref{eq:VP1} problems, or their variants, are equivalent. Alternatively, it results from the first part of the corollary and the fact that $\mathbb{O}^{p \times k}$ is included in $\mathbb{R}^{p \times k}_k$ and is a compact set for the topology of $\mathbb{R}^{p \times k}_k$ induced by the topology of $\mathbb{R}^{p \times k}$ (as $\mathbb{O}^{p \times k}$ is compact in  $\mathbb{R}^{p \times k}$ according to Theorem~\ref{theo2.3:box}). In these conditions, $\psi(.)$ is a continuous mapping over $\mathbb{R}^{p \times k}_k$ and attains its infimum over the compact set $\mathbb{O}^{p \times k}$, which implies directly that the set of global minimizers of the ~\eqref{eq:VP1} problem is not empty and we are done.

\end{proof}

Consider the preimage of $\mathbb{R}^{p \times k}_k$ by $h(.)$, e.g., $h^{-1} ( \mathbb{R}^{p \times k}_k )$. As $h(.)$ is continuous and $\mathbb{R}^{p \times k}_k$ is an open set for the topology of  $\mathbb{R}^{p \times k}$, $h^{-1} ( \mathbb{R}^{p \times k}_k )$ is also an open set of $\mathbb{R}^{p.k}$. Similarly, as $\mathbb{O}^{p \times k}$ is a compact set for the topology of  $\mathbb{R}^{p \times k}$, $h^{-1} ( \mathbb{O}^{p \times k} )$ is a compact set of $\mathbb{R}^{p.k}$ (because $h(.)$ is a homeomorphism or more simply because the reciprocal image of the closed set $\mathbb{O}^{p \times k}$ by $h(.)$ is a closed set of  $\mathbb{R}^{p.k}$ and $\Vert  h( \mathbf{a} ) \Vert_F = \Vert \mathbf{a} \Vert_{2}$ for all $\mathbf{a} \in \mathbb{R}^{p.k}$).
Then the preceding results suggest that it is much more convenient to restrict the domain of definition of the cost function $\psi(.)$ to the open set $h^{-1} ( \mathbb{R}^{p \times k}_k )$ or even to the compact set $h^{-1} ( \mathbb{O}^{p \times k} )$ when we try to solve the~\eqref{eq:VP1} problem or its convex variants with a strictly positive weight matrix. In these conditions, the~\eqref{eq:VP1} problem is a well-posed problem with a non empty set of solutions and the cost function $\psi(.)$ is continuous and even smooth over $h^{-1} ( \mathbb{R}^{p \times k}_k )$ or $h^{-1} ( \mathbb{O}^{p \times k} )$. More precisely, the fact that $\mathbf{F}(.)$  is a continuous linear mapping and has the local constant rank property for all neighborhoods included in $h^{-1} ( \mathbb{R}^{p \times k}_k )$  implies also that $\mathbf{P}^{\bot}_{\mathbf{F}(.)}$ (and thus also $\psi(.)$) is continuously and infinitely differentiable at each point of $h^{-1} ( \mathbb{R}^{p \times k}_k )$. See Theorem 8.4  in Chapter 8 of Magnus and Neudecker~\cite{MN2019} and Subsection~\ref{jacob:box} for details.

We now discuss the continuity of $\mathbf{P}_{\mathbf{F}(.)}$, $\mathbf{P}^{\bot}_{\mathbf{F}(.)}$ and $\psi(.)$ in the case where some elements of the weight matrix $\mathbf{W}$ are equal to zero. We already know that the~\eqref{eq:VP1} problem is not well-posed in these conditions as the set of global minimizers of $\psi(.)$ over $h^{-1} ( \mathbb{R}^{p \times k}_k )$ or $h^{-1} ( \mathbb{O}^{p \times k} )$ can be empty, as already noted for the equivalent forms~\eqref{eq:P0} or~\eqref{eq:P1} of the WLRA problem. In this more difficult case, $\mathbf{P}_{\mathbf{F}(.)}$ and $\mathbf{P}^{\bot}_{\mathbf{F}(.)}$ are also generally discontinuous at all $\mathbf{a} \in h^{-1} (  \mathbb{R}^{p \times k}_{<k} )$, as illustrated by the following theorem.
\\
\begin{theo3.12} \label{theo3.12:box}
Let $\mathbf{X} \in \mathbb{R}^{p \times n} \text{ and } \mathbf{W} \in \mathbb{R}^{p \times n}_{+}$, and $k \le \emph{rank}( \mathbf{X} )$. If $\mathbf{W}_{.j} \in \mathbb{R}^{p}_{+*}$ (e.g., if the $j^{th}$ column of the weight matrix has no zero element), then the matrix function $\mathbf{P}_{\mathbf{F}_{j}(.)}$ from $\mathbb{R}^{p.k}$ to $\mathbb{R}^{p \times p}$ defined by
\begin{equation*}
\mathbf{a} \mapsto \mathbf{P}_{\mathbf{F}_{j}(  \mathbf{a}  )} = \mathbf{F}_{j}(  \mathbf{a} ) \mathbf{F}_{j}(  \mathbf{a} )^{+}  \ ,
\end{equation*}
where $\mathbf{F}_{j}(  \mathbf{a} )^{+}$ is the pseudo-inverse of $\mathbf{F}_{j}(  \mathbf{a} )$ and $\mathbf{F}_{j}(  \mathbf{a} )$ is the $p \times k$ matrix 
\begin{equation*}
\mathbf{F}_{j}(  \mathbf{a} ) = \emph{diag}(\sqrt{\mathbf{W}}_{.j}) h(  \mathbf{a}  ) = \emph{diag}(\sqrt{\mathbf{W}}_{.j})\mathbf{A}  \ ,
\end{equation*}
 is continuous at all $\mathbf{a} \in h^{-1}(\mathbb{R}^{p \times k}_k)$ and discontinuous at all $\mathbf{a} \in h^{-1}(\mathbb{R}^{p \times k}_{<k})$. Furthermore, in the same conditions, the matrix function $\mathbf{P}_{\mathbf{F}(.)}$ from $\mathbb{R}^{p.k}$ to $\mathbb{R}^{p.n \times p.n}$ defined by
\begin{equation*}
\mathbf{a} \mapsto \mathbf{P}_{\mathbf{F}(  \mathbf{a}  )} =  \bigoplus_{j=1}^n \mathbf{P}_{\mathbf{F}_{j}(  \mathbf{a}  )}
\end{equation*}
is also discontinuous at all $\mathbf{a} \in  h^{-1}(\mathbb{R}^{p \times k}_{<k})$.
\end{theo3.12}
\begin{proof}
The proof is similar to the one of Theorem~\ref{theo3.11:box} and is thus omitted.

\end{proof}

However, when  some weights are equal to zero (e.g., when $\mathbf{W} \in \mathbb{R}^{p \times n}_{+}$ instead of $\mathbb{R}^{p \times n}_{+*}$), the condition that $h(\mathbf{a})  \in \mathbb{R}^{p \times k}_k$ is still necessary, but is not sufficient to ensure the continuity of the orthogonal projectors $\mathbf{P}_{\mathbf{F}(.)}$ and $\mathbf{P}^{\bot}_{\mathbf{F}(.)}$.
More precisely, we will demonstrate now that the set of points of $h^{-1} ( \mathbb{R}^{p \times k}_k )$ for which $\mathbf{P}_{\mathbf{F}(.)}$ and $\mathbf{P}^{\bot}_{\mathbf{F}(.)}$ are discontinuous, is not always empty, can be even infinite and grows in size with the number of zero-elements of the weight matrix (see Theorem~\ref{theo3.13:box} below).
Finally, we will show that the $j^{th}$ atomic function $\psi_{j} (.)$ is also discontinuous at all points of $h^{-1} ( \mathbb{R}^{p \times k}_k )$ for which $\mathbf{P}^{\bot}_{\mathbf{F}_{j}(.)}$ is discontinuous (see again Theorem~\ref{theo3.13:box}) and, in these conditions, $\psi(.) = \frac{1}{2} \sum_{j=1}^n \psi_{j} (.)$ can be hardly continuous or differentiable at those points. These results generalize the examples given in Dai et al.~\cite{DMK2011}\cite{DKM2012} about the discontinuity of some of the $j^{th}$ atomic functions $\psi_{j} (.)$ when some entries of the matrix  $\mathbf{X}$ are missing in the case of binary weights and provide a systematic characterization of the subset of points of $h^{-1} ( \mathbb{R}^{p \times k}_k )$ for which some of the atomic functions $\psi_{j} (.)$ are discontinuous when the weight matrix $\mathbf{W}$ is not strictly positive. This systematic assessment of the discontinuities of $\psi(.)$ is possible because we consider $h^{-1} ( \mathbb{R}^{p \times k}_k )$ as the domain of definition of $\psi(.)$  instead of $h^{-1} (  \mathbb{O}^{p \times k} )$ as in Dai et al.~\cite{DMK2011}\cite{DKM2012}, who used a Grassmann manifold's framework to minimize $\psi(.)$.

In order to identify precisely the points of $h^{-1} (\mathbb{R}^{p \times k}_k )$ for which the $j^{th}$ atomic function defined by
\begin{equation*}
\psi_{j}(\mathbf{a} ) = \big\Vert \mathbf{P}^{\bot}_{\mathbf{F}_{j}(  \mathbf{a} ) }\mathbf{x}_{j} \big\Vert^{2}_{2} = \big\Vert \left(  \mathbf{I}_p - \mathbf{F}_{j}(  \mathbf{a} ) \mathbf{F}_{j}(  \mathbf{a} )^{+} \right)\mathbf{x}_{j} \big\Vert^{2}_{2}  \ ,  \forall  \mathbf{a} \in \mathbb{R}^{p.k} \ ,
\end{equation*}
where  $\mathbf{P}^{\bot}_{\mathbf{F}_{j}(  \mathbf{a} ) }$ is the orthogonal projector onto the orthogonal complement of $\emph{ran}\big( \mathbf{F}_{j}(  \mathbf{a} ) \big)$ and $\mathbf{x}_{j}= \sqrt{\mathbf{W}}_{.j} \odot \mathbf{X}_{.j} $, is discontinuous, let us first define what we call the $j^{th}$ barrier set $\mathcal{B}_j$ associated with this  $j^{th}$ atomic function $\psi_{j}(.)$ and the corresponding $p \times k$ matrix function $\mathbf{F}_{j}(.)$.
\\
\begin{def3.2} \label{def3.2:box}
Let $\mathbf{W} \in \mathbb{R}^{p \times n}_{+}$ and define for all $\mathbf{a} \in \mathbb{R}^{p.k}$, the matrix function $\mathbf{F}(.)$
\begin{equation*}
\mathbf{F} : \mathbb{R}^{p.k} \longrightarrow \mathbb{R}^{p.n  \times k.n } : \mathbf{a} \mapsto \mathbf{F}(  \mathbf{a} ) = \bigoplus_{j=1}^n \mathbf{F}_{j}(  \mathbf{a} )  \ ,
\end{equation*}
where $\mathbf{F}_{j}(  \mathbf{a} ) = \emph{diag}(\sqrt{\mathbf{W}}_{.j})  h(\mathbf{a})  \in \mathbb{R}^{p \times k}$ is called the  $j^{th}$ matrix function. The $j^{th}$ barrier set $\mathcal{B}_j$ associated with the $j^{th}$ atomic and matrix functions, $\psi_{j}(.)$ and $\mathbf{F}_{j}(.)$, is the subset of $\mathbb{R}^{p.k}$ defined by
\begin{equation*}
\mathcal{B}_j = \big\{   \mathbf{a} \in \mathbb{R}^{p.k}  \text{ / }  \exists  \mathbf{A} \in  \mathbb{R}^{p \times k}_k  \text{ with } \mathbf{A} = h(  \mathbf{a}  ) \text{ and } \mathbf{F}_{j}(  \mathbf{a} ) = \mathbf{0}^{p \times k}   \big\}  \ ,
\end{equation*}
where $\mathbf{0}^{p \times k}$ is the zero $p \times k$ matrix.
\\
\end{def3.2}

\begin{remark3.8} \label{remark3.8:box}

This terminology is due to Dai et al.~\cite{DMK2011}\cite{DKM2012}, who illustrated by a few examples that some points belonging to the intersection of a barrier set $\mathcal{B}_j$ with $h^{-1} (\mathbb{O}^{p \times k} )$ act as "barriers", which may prevent gradient descent algorithms used to solve the~\eqref{eq:VP1} problem from converging to a global minimum or infimum. $\blacksquare$
\\
\end{remark3.8}

\begin{remark3.9} \label{remark3.9:box}

The subset $\mathcal{B}_j$ of $\mathbb{R}^{p.k}$ introduced in Definition~\ref{def3.2:box} can also be defined as follows. Consider again the preimage of $\mathbb{R}^{p \times k}_k$ by $h(.)$, e.g., $h^{-1} ( \mathbb{R}^{p \times k}_k )$, which is an open set  of $\mathbb{R}^{p.k}$ and also the continuous linear mapping $\mathbf{F}_{j} : \mathbb{R}^{p.k}  \longrightarrow \mathbb{R}^{p \times k}, \mathbf{a} \mapsto \mathbf{F}_{j}(  \mathbf{a} )$. The subset of $\mathbb{R}^{p.k}$ such that $\mathbf{F}_{j}(  \mathbf{a} ) = \mathbf{0}^{p \times k}$ is the null space of $\mathbf{F}_{j}(.)$, which is a closed linear subspace of $\mathbb{R}^{p.k}$. With these results, we have
\begin{equation*}
\mathcal{B}_j = h^{-1} ( \mathbb{R}^{p \times k}_k )  \cap  \mathbf{F}^{-1}_{j} (  \lbrace  \mathbf{0}^{p \times k }  \rbrace ) = h^{-1} ( \mathbb{R}^{p \times k}_k )  \cap  null( \mathbf{F}_{j} )  \ ,
\end{equation*}
e.g., $\mathcal{B}_j$ is the intersection of the preimages $h^{-1} ( \mathbb{R}^{p \times k}_k )$ and $\mathbf{F}^{-1}_{j} (  \lbrace \mathbf{0}^{p \times k } \rbrace )$, which is not open and nor closed in $\mathbb{R}^{p.k}$. $\blacksquare$
\\
\end{remark3.9}

We first observe that $\mathcal{B}_j$ is empty if all elements of the column vector $\mathbf{W}_{.j}$ are greater than zero as in that case we have $ \emph{rank}\big(  \emph{diag}( \sqrt{\mathbf{W}}_{.j} ) \big) = p$ and, thus, $ \emph{rank}(   \mathbf{F}_{j}(  \mathbf{a} ) ) = \emph{rank}( \mathbf{A} ) = k$ if  $h(  \mathbf{a}  ) = \mathbf{A} \in  \mathbb{R}^{p \times k}_k$, and in these conditions  $\mathbf{F}_{j}(  \mathbf{a} ) \ne \mathbf{0}^{p \times k}$. On the other hand, if the number of zero elements of the column vector $\mathbf{W}_{.j}$ is greater or equal to $k$, $\mathcal{B}_j$ is not empty and is an infinite subset of $h^{-1} ( \mathbb{R}^{p \times k}_k )$ as demonstrated in Theorem~\ref{theo3.13:box} below.

To demonstrate this proposition, let us introduce again some notations. Suppose that some elements of the column vector $\mathbf{W}_{.j}$ are equal to zero, which is equivalent to say that the corresponding elements of the column vector $\mathbf{X}_{.j}$ are missing. Then, let $p_u$ be the number of zero elements of $\mathbf{W}_{.j}$ and $p_o = p - p_u$ the number of elements of $\mathbf{W}_{.j}$ which are different from zero, e.g., $p_u$ and $p_o$ are, respectively, the number of "unobserved" and "observed" entries in the column vector $\mathbf{X}_{.j}$. In this case, we can partition the vectors $\mathbf{X}_{.j}$ and $\mathbf{W}_{.j}$ as
\begin{equation*}
\begin{bmatrix}
\mathbf{X}^{u}_{.j}  \\
\mathbf{X}^{o}_{.j}
\end{bmatrix}  = \mathbf{P}_j \mathbf{X}_{.j}
 \text{ and } 
 \begin{bmatrix}
\mathbf{W}^{u}_{.j}  \\
\mathbf{W}^{o}_{.j}
\end{bmatrix}  = 
 \begin{bmatrix}
\mathbf{0}^{p_u}  \\
\mathbf{W}^{o}_{.j}
\end{bmatrix}  = 
\mathbf{P}_j \mathbf{W}_{.j}  \ .
\end{equation*}
Here $\mathbf{X}^{u}_{.j} \in  \mathbb{R}^{p_u}$ is the "unobserved" part of the column vector $\mathbf{X}_{.j}$, $\mathbf{X}^{o}_{.j} \in  \mathbb{R}^{p_o}$ is the "observed" part of this vector and $\mathbf{P}_j$ is any $p \times  p$ permutation matrix, which reorders the elements of $\mathbf{W}_{.j}$ so that the zero elements of $\mathbf{W}_{.j}$ appear first. Obviously, $\mathbf{P}_j$ is not unique, but we can use any such permutation in what follows. Note also that this permutation matrix $\mathbf{P}_j$ could be different for each pair of column vectors $\mathbf{X}_{.j}$ and $\mathbf{W}_{.j}$, such that some of the elements of $\mathbf{W}_{.j}$ are equal to zero, if the patterns of "unobserved" entries differ among the columns of the matrix $\mathbf{X}$. Similarly, $\mathbf{A} \in  \mathbb{R}^{p \times k}$ can be partitioned as
\begin{equation*}
\mathbf{P}_j \mathbf{A} = 
\begin{bmatrix}
\mathbf{A}^{u} \\
\mathbf{A}^{o}
\end{bmatrix}  \ ,
\end{equation*}
where $\mathbf{A}^{u}  \in  \mathbb{R}^{p_u \times k}$ and $\mathbf{A}^{o}  \in  \mathbb{R}^{p_o \times k}$. With these notations, we have the following theorem:
\\
\begin{theo3.13} \label{theo3.13:box}
Let $p_u \in \mathbb{N}_{*}$ and $p_o \in \mathbb{N}_{*}$ with $p = p_u + p_o$ design, respectively, the numbers of "unobserved" and "observed" entries in the $j^{th}$ column of $\mathbf{X}$. If $p_u \ge k$, the $j^{th}$ barrier set $\mathcal{B}_j$ is equal to the (nonempty) subset of $\mathbb{R}^{p.k}$
\begin{equation*}
\mathcal{B}^{*}_j = \big\{   \mathbf{a} \in \mathbb{R}^{p.k}  \text{ / }   \exists  \mathbf{A}^{u} \in  \mathbb{R}^{p_u \times k}_k \text{ and }  \mathbf{P}_j h(  \mathbf{a}  )  =   \begin{bmatrix} \mathbf{A}^{u}  \\   \mathbf{0}^{p_o \times k} \end{bmatrix} \big\}  \ ,
\end{equation*}
where $\mathbf{0}^{p_o \times k}$ is the zero $p_o \times k$ matrix and $\mathbf{P}_j$ is any $p \times  p$ permutation matrix, which reorders the elements of $\mathbf{X}_{.j}$ so that the missing elements of $\mathbf{X}_{.j}$ appear first. On the other hand, if $p_u < k$ then $\mathcal{B}_j$ is empty.
\end{theo3.13}
\begin{proof}
Suppose first that $p_u \ge k$ and let $\mathbf{a} \in \mathcal{B}_j$. Then, it exists $h(  \mathbf{a}  ) = \mathbf{A} \in \mathbb{R}^{p \times k}_k$ and we have the implications
\begin{align*}
\mathbf{F}_{j}(  \mathbf{a} ) = \mathbf{0}^{p \times k}  &  \Rightarrow  \emph{diag}( \sqrt{\mathbf{W}}_{.j} ) \mathbf{A} = \mathbf{0}^{p \times k}  \\
                                 &  \Rightarrow   \emph{diag}( \mathbf{P}^{T}_j  \mathbf{P}_j \sqrt{\mathbf{W}}_{.j} ) \mathbf{P}^{T}_j  \begin{bmatrix} \mathbf{A}^{u}  \\   \mathbf{A}^{o }\end{bmatrix}  = \mathbf{0}^{p \times k} \\
                                 &  \Rightarrow   \mathbf{P}^{T}_j \emph{diag}( \mathbf{P}_j \sqrt{\mathbf{W}}_{.j} )  \begin{bmatrix} \mathbf{A}^{u}  \\   \mathbf{A}^{o }\end{bmatrix}  = \mathbf{0}^{p \times k}  \\
                                 &  \Rightarrow   \emph{diag}(  \begin{bmatrix} \mathbf{0}^{p_u}  \\   \sqrt{\mathbf{W}}^{o }_{.j}  \end{bmatrix} )   \begin{bmatrix} \mathbf{A}^{u}  \\   \mathbf{A}^{o} \end{bmatrix}  = \mathbf{0}^{p \times k}  \\
                                 &  \Rightarrow   \emph{diag}(   \sqrt{\mathbf{W}}^{o }_{.j}  )   \mathbf{A}^{o }  = \mathbf{0}^{p_o \times k} \\
                                 &  \Rightarrow    \mathbf{A}^{o }  = \mathbf{0}^{p_o \times k} \ .
\end{align*}
The last implication results from the fact that all elements of the vector $\mathbf{W}^{o }_{.j}$ are different from zero by definition. This implies that $\mathbf{A} =  \mathbf{P}^{T}_j \begin{bmatrix} \mathbf{A}^{u}  \\   \mathbf{0}^{p_o \times k}\end{bmatrix}$. Furthermore, since $\mathbf{P}^{T}_j$ is a permutation matrix and, thus, of full rank $p$, we have $k = \emph{rank}( \mathbf{A} ) = \emph{rank}( \mathbf{A}^{u} )$ and $\mathbf{A}^{u} \in  \mathbb{R}^{p_u \times k}_k$, which shows  that $\mathbf{a}  \in \mathcal{B}^{*}_j$.

Reciprocally, suppose that $\mathbf{a} \in \mathcal{B}^{*}_j$. Then, it exists $\mathbf{A} = h(  \mathbf{a}  )  \in  \mathbb{R}^{p \times k}$  and  $\mathbf{A}^{u} \in  \mathbb{R}^{p_u \times k}_k$ such that $\mathbf{A} =  \mathbf{P}^{T}_j   \begin{bmatrix} \mathbf{A}^{u}  \\   \mathbf{0}^{p_o \times k} \end{bmatrix}$ and $\mathbf{A}$ is of rank $k$ as $\mathbf{P}^{T}_j$ is of full rank $p$ and  $\mathbf{A}^{u}$ is of rank $k$. Then, we have
\begin{align*}
\mathbf{F}_{j}(  \mathbf{a} ) & = \emph{diag}( \sqrt{\mathbf{W}}_{.j} ) \mathbf{A}   \\
      & = \emph{diag}( \sqrt{\mathbf{W}}_{.j} ) \mathbf{P}^{T}_j  \begin{bmatrix} \mathbf{A}^{u}  \\   \mathbf{0}^{p_o \times k} \end{bmatrix}   \\
      & = \emph{diag}( \mathbf{P}^{T}_j  \mathbf{P}_j  \sqrt{\mathbf{W}}_{.j} ) \mathbf{P}^{T}_j  \begin{bmatrix} \mathbf{A}^{u}  \\   \mathbf{0}^{p_o \times k} \end{bmatrix}   \\
      & = \mathbf{P}^{T}_j  \emph{diag}(   \mathbf{P}_j  \sqrt{\mathbf{W}}_{.j} )  \begin{bmatrix} \mathbf{A}^{u}  \\   \mathbf{0}^{p_o \times k} \end{bmatrix}   \\
      & = \mathbf{P}^{T}_j  \emph{diag}(  \begin{bmatrix} \mathbf{0}^{p_u}  \\   \sqrt{\mathbf{W}}^{o }_{.j}  \end{bmatrix} )  \begin{bmatrix} \mathbf{A}^{u}  \\   \mathbf{0}^{p_o \times k} \end{bmatrix}   \\
      & = \mathbf{P}^{T}_j   \mathbf{0}^{p_o \times k}  = \mathbf{0}^{p_o \times k}  \ ,
\end{align*}
and $\mathbf{a} \in \mathcal{B}_j$, which demonstrates the first part of the theorem.

Suppose now that $p_u < k$, then if  $\mathbf{a} \in  \mathbb{R}^{p.k}$ and $h(  \mathbf{a}  ) = \mathbf{A}  \in  \mathbb{R}^{p \times k}$, we have,
\begin{align*}
\mathbf{F}_{j}(  \mathbf{a} ) = \mathbf{0}^{p \times k}  &  \Rightarrow   \emph{diag}( \sqrt{\mathbf{W}}_{.j} ) \mathbf{A} = \mathbf{0}^{p \times k}  \\
                                   \Rightarrow  &   \emph{diag}(  \begin{bmatrix} \mathbf{0}^{p_u}  \\   \sqrt{\mathbf{W}}^{o }_{.j}  \end{bmatrix} )   \begin{bmatrix} \mathbf{A}^{u}  \\   \mathbf{A}^{o} \end{bmatrix}  = \mathbf{0}^{p \times k}  \\
                                   \Rightarrow  &    \emph{diag}(   \sqrt{\mathbf{W}}^{o }_{.j}  )   \mathbf{A}^{o }  = \mathbf{0}^{p_o \times k} \\
                                   \Rightarrow  &     \mathbf{A}^{o }  = \mathbf{0}^{p_o \times k} \ .
\end{align*}
This implies that $\mathbf{A} =  \mathbf{P}^{T}_j \begin{bmatrix} \mathbf{A}^{u}  \\   \mathbf{0}^{p_o \times k}\end{bmatrix}$ with $\mathbf{A}^{u}  \in \mathbb{R}^{p_u \times k}$. But, as $\mathbf{P}^{T}_j$ is a permutation matrix and, thus, of full rank $p$, we have $\emph{rank}( \mathbf{A} ) = \emph{rank}( \mathbf{A}^{u} ) \le \text{min}( p_u , k) = p_u < k$ and we conclude that $\mathbf{a} \notin  \mathcal{B}_j$. In other words, if $p_u < k$, for $\mathbf{a} \in  \mathbb{R}^{p.k}$ we cannot have $\mathbf{F}_{j}(  \mathbf{a} ) = \mathbf{0}^{p \times k}$ and $h(  \mathbf{a}  ) =  \mathbf{A}  \in  \mathbb{R}^{p \times k}_k$.

\end{proof}

\begin{remark3.10} \label{remark3.10:box}
Using Remark~\ref{remark3.9:box}, it is easy to verify that the first part of Theorem~\ref{theo3.13:box} results from  $(i)$ the fact that  $null( \mathbf{F}_{j} )$ is a linear subspace of $\mathbb{R}^{p.k}$ of dimension $p_{u}.k$ and is equal to
\begin{equation*}
null( \mathbf{F}_{j} ) = \Big \{   \mathbf{a} \in \mathbb{R}^{p.k}  \text{ / }   \exists  \mathbf{A}^{u} \in  \mathbb{R}^{p_u \times k} \text{ such that }  \mathbf{P}_j h(  \mathbf{a}  )  =   \begin{bmatrix} \mathbf{A}^{u}  \\   \mathbf{0}^{p_o \times k}\end{bmatrix} \Big \}  \ ,
\end{equation*}
and  $(ii)$ the property: Let $\mathbf{A} \in  \mathbb{R}^{p \times k}$ and  $\mathbf{A}^{u} \in  \mathbb{R}^{p_u \times k}$ with $\mathbf{P}_j \mathbf{A} = \begin{bmatrix}
\mathbf{A}^{u} \\ \mathbf{0}^{p_o \times k} \end{bmatrix} $. If $p_u \ge k$ :
\begin{equation*}
 \emph{rank}( \mathbf{A} ) = {k} \Longleftrightarrow \emph{rank}( \mathbf{A}^{u} ) = {k} \ .
\end{equation*}
$\blacksquare$
\\
\end{remark3.10}

\begin{theo3.14} \label{theo3.14:box}
With the same notations as in Theorem~\ref{theo3.13:box} and above, if $p_u \ge k$ and $p_o \ge 1$ (e.g., if there is at least one observed entry in the column vector $\mathbf{X}_{.j}$), the following two assertions are true:
\begin{align*}
1) &  \quad\   \text{ The orthogonal projectors }  \mathbf{P}_{\mathbf{F}_{j}(.)}   \text{ and }  \mathbf{P}^{\bot}_{\mathbf{F}_{j}(.)}   \text{ are discontinuous at all points of } \mathcal{B}_j , \\
2) &  \quad\  \text{ If }   \Vert \mathbf{X}^{o}_{.j} \Vert_{2}  \ne 0,  \text{ then }  \psi_{j}(.) \text{  is also discontinuous at all points of } \mathcal{B}_j .
\end{align*}

\end{theo3.14}
\begin{proof}
According to Theorem~\ref{theo3.10:box}, in order to demonstrate the first assertion of the theorem it suffices to show that the matrix function $\mathbf{F}_{j}(.)$ has no local constant rank for all $\mathbf{a} \in  \mathcal{B}_j$. We first note that $\mathbf{F}_{j}(  \mathbf{a} ) = \mathbf{0}^{p \times k} $ and, thus, $\emph{rank}(  \mathbf{F}_{j}(  \mathbf{a} ) ) = 0$ if $\mathbf{a} \in  \mathcal{B}_j$. In these conditions, it suffices to show that, $\forall \alpha  \in  \mathbb{R}_{+*}$, if we consider the open ball $B_{p.k}(\mathbf{a}, \alpha)$, with center $\mathbf{a}$ and radius $\alpha$, in $\mathbb{R}^{p \times k}$, it exists $\mathbf{d} \in  B_{p.k}(\mathbf{a}, \alpha)  \cap  h^{-1}( \mathbb{R}^{p \times k}_k )$ such that $\emph{rank}( \mathbf{F}_{j}(  \mathbf{d} )  ) \ne 0$.

Since $\mathbf{a} \in  \mathcal{B}_j$ and  $p_u \ge k$ by hypothesis, according to Theorem~\ref{theo3.13:box}, it exists $\mathbf{A}^{u} \in  \mathbb{R}^{p_u \times k}_k$ such that   $\mathbf{P}_j \mathbf{A} = \begin{bmatrix} \mathbf{A}^{u} \\ \mathbf{0}^{p_o \times k} \end{bmatrix}$ where $\mathbf{A}  =  h( \mathbf{a} ) \in  \mathbb{R}^{p \times k}_k$.
Now, let $\beta  \in  \mathbb{R}_{+*}$ such that $\beta < \alpha$ and  $\mathbf{D}^{o}  \in  \mathbb{R}^{p_o \times k}$ such that $\mathbf{D}^{o}_{11} = \beta$ and $\mathbf{D}^{o}_{ij} = 0$ if $i \ne 1$ and  $j \ne 1$. If we define $\mathbf{D} = \mathbf{P}^{T}_j  \begin{bmatrix} \mathbf{A}^{u} \\ \mathbf{D}^{o}  \end{bmatrix}$ and $\mathbf{d} = \emph{vec}(  \mathbf{D}^{T} )$, we have $\mathbf{d}  \in h^{-1}( \mathbb{R}^{p \times k}_k )$ as $\mathbf{D} \in \mathbb{R}^{p \times k}_k$ (since $\mathbf{A}^{u} \in  \mathbb{R}^{p_u \times k}_k$) and also
\begin{equation*}
\Vert \mathbf{d} - \mathbf{a} \Vert^{2}_{2} = \Vert \mathbf{D} - \mathbf{A} \Vert^{2}_F = \beta^{2} \le \alpha^{2} \ .
\end{equation*}
Thus, $\mathbf{d} \in  B_{p.k}(\mathbf{a}, \alpha)  \cap  h^{-1}( \mathbb{R}^{p \times k}_k )$. Obviously, $\mathbf{F}_{j}(  \mathbf{d} ) \ne \mathbf{0}^{p \times k}$ as $\mathbf{D}^{o}  \ne \mathbf{0}^{p_o \times k}$ and  $\emph{rank}( \mathbf{F}_{j}(  \mathbf{d} )  ) = 1$ and we conclude that $\mathbf{F}_{j}(.)$ has no local constant rank at all points $\mathbf{a}$ of $\mathcal{B}_j$.

For demonstrating the second part of the theorem, we first remark that, for all points $\mathbf{a} \in \mathcal{B}_j$, we have $\mathbf{F}_{j}(  \mathbf{a} ) = \mathbf{0}^{p \times k}$  and $\mathbf{F}_{j}(  \mathbf{a} )^{+} = \mathbf{0}^{k \times p}$ and, thus,
\begin{equation*}
\psi_{j}(\mathbf{a} )  =   \Vert \mathbf{P}^{\bot}_{\mathbf{F}_{j}(  \mathbf{a} )}  \mathbf{x}_j \Vert^{2}_{2} = \Vert \mathbf{x}_j \Vert^{2}_{2} = \Vert  \mathbf{x}^{u}_j \Vert^{2}_{2} + \Vert  \mathbf{x}^{o}_j \Vert^{2}_{2} = \Vert  \mathbf{x}^{o}_j \Vert^{2}_{2} \ ,
\end{equation*}
since $\mathbf{W}^{u}_{.j} = \mathbf{0}^{p_u}$. By hypothesis, we have $ \Vert \mathbf{X}^{o}_{.j} \Vert_{2} \ne 0$ and, thus, $\Vert \mathbf{x}^{o}_j \Vert^{2}_{2} \ne 0$, so let us consider the open interval $\big \rbrack \frac{1}{2} \Vert \mathbf{x}^{o}_j \Vert^{2}_{2},   \frac{3}{2} \Vert \mathbf{x}^{o}_j \Vert^{2}_{2}  \big \lbrack$ of $\mathbb{R}$, e.g., the open ball $B_1\big( \psi_{j}(\mathbf{a} ), \frac{1}{2}  \psi_{j}(\mathbf{a} ) \big)$ of $\mathbb{R}$.
We have to show that, for all $\alpha  \in  \mathbb{R}_{+*}$, it exists $\mathbf{d} \in \mathbb{R}^{p.k}$ such that $\mathbf{d} \in  B_{p.k}(\mathbf{a}, \alpha)  \cap  h^{-1}( \mathbb{R}^{p \times k}_k )$, but $\psi_{j}(\mathbf{d} ) \notin B_1\big( \psi_{j}(\mathbf{a} ), \frac{1}{2}  \psi_{j}(\mathbf{a} ) \big)$.

Let $\beta  \in  \mathbb{R}_{+*}$ such that $\beta < \alpha$, and define $\mathbf{D}^{o} \in \mathbb{R}^{p_o \times k}$ by $\mathbf{D}^{o}_{.1} = \beta . \frac{\mathbf{X}^{o}_{.j}}{\Vert \mathbf{X}^{o}_{.j} \Vert_{2}}$ and $\mathbf{D}^{o}_{.i} = 0^{p_o}$ for $i=2$ to $k$, where $0^{p_o}$ is the zero vector of dimension $p_o$. Setting now
\begin{equation*}
\mathbf{D} = \mathbf{P}^{T}_j \begin{bmatrix}  \mathbf{A}^{u} \\ \mathbf{D}^{o}  \end{bmatrix} \in \mathbb{R}^{p \times k}  , \mathbf{d} = \emph{vec}(  \mathbf{D}^{T} ) \in \mathbb{R}^{p.k} \text{ and } \mathbf{d}^{o} = \emph{vec} \big (   ( \mathbf{D}^{o} )^{T} \big)  \in \mathbb{R}^{{p_o}.k}  \ ,
\end{equation*}
we have $\Vert \mathbf{d} - \mathbf{a} \Vert^{2}_2 = \Vert \mathbf{D} - \mathbf{A} \Vert^{2}_F = \beta^2 \le \alpha^2$ and rank($\mathbf{D}$) = $k$. Thus, $\mathbf{d} \in  B_{p.k}(\mathbf{a}, \alpha)  \cap  h^{-1}( \mathbb{R}^{p \times k}_k )$ and
\begin{align*}
\psi_{j}(\mathbf{d} )  & =  \Vert \mathbf{P}^{\bot}_{\mathbf{F}_{j}(  \mathbf{d} )}  \mathbf{x}_j \Vert^{2}_{2}   \\
                                & =   \Vert    \mathbf{x}_j  - \mathbf{F}_{j}(  \mathbf{d} ) \mathbf{F}_{j}(  \mathbf{d} )^{+} \mathbf{x}_j  \Vert^{2}_{2}  \\
                                & =   \Vert    \mathbf{x}^{o}_j  - \mathbf{F}_{j}(  \mathbf{d}^{o} ) \mathbf{F}_{j}(  \mathbf{d}^{o} )^{+} \mathbf{x}^{o}_j  \Vert^{2}_{2} \\
                                & =   \min_{\mathbf{b}\in\mathbb{R}^{k} } \,  \Vert    \mathbf{x}^{o}_j  - \mathbf{F}_{j}(  \mathbf{d}^{o} ) \mathbf{b}  \Vert^{2}_{2} \\
                                & =   \min_{\mathbf{b}_1 \in\mathbb{R} } \,  \Vert    \mathbf{x}^{o}_j  -  ( \sqrt{\mathbf{W}}^{o }_{.j}   \odot \mathbf{D}^{o}_{.1} ) \mathbf{b}_1  \Vert^{2}_{2} \\
                                & =   \min_{\mathbf{b}_1 \in\mathbb{R} } \,  \Vert    \mathbf{x}^{o}_j  -  \frac{\beta}{\Vert \mathbf{X}^{o}_{.j} \Vert_{2}} ( \sqrt{\mathbf{W}}^{o }_{.j}   \odot \mathbf{X}^{o}_{.j} ) \mathbf{b}_1  \Vert^{2}_{2} \\
                                & =   \min_{\mathbf{b}_1 \in\mathbb{R} } \,  \Vert    \mathbf{x}^{o}_j  -  \frac{\beta}{\Vert \mathbf{X}^{o}_{.j} \Vert_{2}} \mathbf{x}^{o}_j    \mathbf{b}_1  \Vert^{2}_{2} \\
                                & =   0 \, \text{ with } \mathbf{b}_1 = \frac{\Vert \mathbf{X}^{o}_{.j} \Vert_{2}}{\beta} \ .
\end{align*}
In other words, $\psi_{j}(\mathbf{d} ) \notin \big \rbrack \frac{1}{2} \Vert \mathbf{x}^{o}_j \Vert^{2}_{2},   \frac{3}{2} \Vert \mathbf{x}^{o}_j \Vert^{2}_{2}  \big \lbrack = B_1\big( \psi_{j}(\mathbf{a} ), \frac{1}{2}  \psi_{j}(\mathbf{a} ) \big)$, which demonstrates that $\psi_{j}(.)$ is discontinuous at all points of  $\mathcal{B}_j$ as claimed in the second part of the theorem.

\end{proof}

\begin{corol3.4} \label{corol3.4:box}
With the same definitions and notations as in Theorem~\ref{theo3.14:box}, the orthogonal projectors $\mathbf{P}_{\mathbf{F}(.)}$  and $\mathbf{P}^{\bot}_{\mathbf{F}(.)}$  are discontinuous at all points $\mathbf{a}  \in  \bigcup_{j=1}^{n} \mathcal{B}_j$.
\end{corol3.4}
\begin{proof}
This results immediately from Theorem~\ref{theo3.14:box} and the equality
\begin{equation*}
\mathbf{P}_{\mathbf{F}(  \mathbf{a} )} =  \bigoplus_{j=1}^n \mathbf{P}_{\mathbf{F}_{j}(  \mathbf{a} ) } ,
\end{equation*}
which shows that the continuity of $\mathbf{P}_{\mathbf{F}(.)}$ is equivalent to the continuity of the $n$ atomic orthogonal projectors $\mathbf{P}_{\mathbf{F}_{j}(.) }$, for $j=1  \text{ to } n$.

\end{proof}

Since $\mathbf{F}_{j}(.)$ is a continuous linear mapping different from the zero constant function, its null space, $null( \mathbf{F}_{j} ) = \mathbf{F}^{-1}_{j} (  \lbrace \mathbf{0}^{p \times k } \rbrace )$, is closed and not dense in $\mathbb{R}^{p.k}$ and its complement  $\mathbb{R}^{p.k}/null( \mathbf{F}_{j} )$ (in $\mathbb{R}^{p.k}$), is nonempty and open. Next, using Theorem~\ref{theo3.14:box} and Corollary~\ref{corol3.4:box}, we know that $\mathbf{P}^{\bot}_{\mathbf{F}(.)}$ is not continuous at all points $\mathbf{a} \in \bigcup_{j=1}^{n} \mathcal{B}_j$, implying that $\mathbf{P}^{\bot}_{\mathbf{F}(.)}$ will not be differentiable and $\psi(.)$ will also, in general, not be continuous and differentiable at these points. This implies that the "feasible" search set for a solution of the~\eqref{eq:VP1} problem will be severely restricted in the case of missing values, at least by standard first- and second-order optimization methods, which require that the objective function must be at least differentiable. Moreover, the "size" of this "feasible" search set will also decrease if the number of missing values increases as it is equal to  $\bigcap_{j=1}^{n}   \lbrack h^{-1}( \mathbb{R}^{p \times k}_k )/\mathcal{B}_j)  \rbrack$.

Hence, when the number of missing values in the data matrix  $\mathbf{X}$ is very large, one may prefer an alternative formulation of the WLRA problem that will allow for a continuous and differentiable objective function $\psi(.)$ for all points of $h^{-1}( \mathbb{R}^{p \times k}_k )$, no barrier sets $\mathcal{B}_j$ and no discontinuities for any of the atomic functions $\psi_{j}(.)$. In this context, when the weights are equal to one or zero (e.g., the missing value problem), Dai et al.~\cite{DKM2012} have proposed to replace the traditional Frobenius metric by what they called a "geometric performance metric" to avoid these discontinuities of the atomic functions when solving the matrix completion problem in a Grassmann manifold's setting. More generally, as discussed above and demonstrated in Theorem~\ref{theo3.8:box}, the minimization of the separable form of the cost function $g_{\lambda}(.)$ defined in equation~\eqref{eq:g_func} with a small regularization parameter $\lambda \in  \mathbb{R}_{+*}$ or continuation Tikhonov methods based on a family of such cost functions $g_{\lambda}(.)$ in which $\lambda$ tends to zero during the iterations are also promising alternatives in such difficult situation. Moreover, these alternatives work with nonuniform weight matrices including zero weights and not only for the missing value problem with binary weights.

We are now set to describe the different algorithms which may be used to minimize $\varphi^{*}(.)$ or $\psi(.)$. We start by a modern description of several ALS regression methods, which all originate from the NIPALS algorithm first introduced by Wold and his collaborators~\cite{W1966}\cite{WL1969}\cite{JHJ2009}, and alternate between minimization of the two sets of variables, $\mathbf{A}$ and $\mathbf{B}$, for solving the formulation~\eqref{eq:P1} of the WLRA problem in Section~\ref{nipals:box}. The more complicated separable NLLS algorithms (e.g., first- and second-order variable projection methods), which explicitly eliminate the linear parameters (for example $\mathbf{b} =  \emph{vec}( \mathbf{B} )$) obtaining a reduced, but somewhat more complicated, functional $\psi(.)$ that involves only the nonlinear parameters (e.g., $\mathbf{a} =  \emph{vec}( \mathbf{A}^{T} )$), are described in Section~\ref{varpro:box}.

\section{The block alternating least-squares method and its variants} \label{nipals:box}

As noted in Subsection~\ref{varpro_wlra:box}, if we fix $\mathbf{a} = h^{-1}( \mathbf{A} ) = \emph{vec}(  \mathbf{A}^{T}  )$, then the problem
\begin{equation*}
\min_{\mathbf{b}\in\mathbb{R}^{k.n}} \, \quad\  \frac{1}{2}   \Vert \mathbf{x} - \mathbf{F}(\mathbf{a})\mathbf{b} \Vert^{2}_2 = \varphi^{*}( \mathbf{A},\mathbf{B}  ) \ ,
\end{equation*}
where $\mathbf{F}(\mathbf{a})$  and $\mathbf{x}$ are also defined in Subsection~\ref{varpro_wlra:box}, is a linear least-squares problem with $k.n$ unknowns $\mathbf{B}_{ij}$. The unique minimum 2-norm solution of this linear least-squares problem for a fixed $\mathbf{A}$ matrix is $\mathbf{\widehat{b}} = \mathbf{F}(\mathbf{a})^{+} \mathbf{x}$ as stated in Subsection ~\ref{lin_alg:box}. More precisely, if we take into account the block structure of $\mathbf{F}(\mathbf{a})$, we observe that the best choice of $\mathbf{b} = \emph{vec}( \mathbf{B} )$ for a given $\mathbf{A}$  matrix is obtained by solving $n$ independent linear least-squares problems, each with $k$ unknowns, and $\mathbf{\widehat{B}}_{.j}$, for $j=1, \cdots, n$, can be calculated by
\begin{equation} \label{eq:B_nipals}
\mathbf{\widehat{B}}_{.j} = \left(  \emph{diag}(  \sqrt{\mathbf{W}}_{.j} ) \mathbf{A} \right)^{+}  (  \sqrt{\mathbf{W}}_{.j} \odot \mathbf{X}_{.j} ) = \mathbf{F}_{j}( \mathbf{a} )^{+} \mathbf{x}_{j} \  ,
\end{equation}
where  $\mathbf{F}_{j}( \mathbf{a} ) =  \emph{diag}(  \sqrt{\mathbf{W}}_{.j} ) \mathbf{A}$ and $ \mathbf{x}_{j} =  \sqrt{\mathbf{W}}_{.j} \odot \mathbf{X}_{.j}$. Likewise, if $\mathbf{b} = \emph{vec} (  \mathbf{B} )$ is fixed, the minimization problem
\begin{equation*}
\min_{\mathbf{a}\in\mathbb{R}^{p.k} } \, \quad\  \frac{1}{2} \Vert \mathbf{z} - \mathbf{G}(\mathbf{b})\mathbf{a} \Vert^{2}_2 = \varphi^{*}(\mathbf{A},\mathbf{B}  ) \ ,
\end{equation*}
where $\mathbf{G}(\mathbf{b})$ and $\mathbf{z}$ are again defined in Subsection~\ref{varpro_wlra:box}, is a linear least-squares problem with $k.p$ unknowns $\mathbf{A}_{ij}$. The unique minimum 2-norm solution of this linear least-squares problem for a fixed $\mathbf{B}$ matrix is then $\mathbf{\widehat{a}} = \mathbf{G}(\mathbf{b})^{+} \mathbf{z}$. Again, by taking into account the block structure of $\mathbf{G}(\mathbf{b})$, we observe that the best choice of $\mathbf{A}$ for a given $\mathbf{B}$  matrix is obtained by solving $p$ independent linear least-squares problems, each with $k$ unknowns, and $\mathbf{\widehat{A}}_{i.}$, for $i=1, \cdots, p$, can be calculated by
\begin{equation}\label{eq:A_nipals}
\mathbf{\widehat{A}}_{i.} = ( \sqrt{\mathbf{W}}_{i.}  \odot \mathbf{X}_{i.} ) \left( \mathbf{B} \emph{diag}( \sqrt{\mathbf{W}}_{i.} ) \right)^{+} =   \left( \mathbf{G}_{i}( \mathbf{b} )^{+} \mathbf{z}_{i} \right)^{T} \  ,
\end{equation}
where $\mathbf{G}_{i}( \mathbf{b} ) =  \emph{diag}( \sqrt{\mathbf{W}}_{i.} ) \mathbf{B}^{T}$ and $\mathbf{z}_{i} = \left( \sqrt{\mathbf{W}}_{i.}  \odot \mathbf{X}_{i.} \right)^{T}$.

The above results suggest that we can minimize the cost function  $\varphi^{*}(.)$ and solve the WLRA problem in its formulation~\eqref{eq:P1} by taking the block separable least-squares approach (e.g., NIPALS algorithm) of Wold and his collaborators~\cite{W1966}\cite{WL1969}\cite{JHJ2009}, a method rediscovered and studied many times after, particularly in the context of low-rank matrix completion and optimization problems~\cite{WWY2012}\cite{JNS2013}\cite{H2014}\cite{LZT2019b}\cite{OUV2023} or in the computer vision community~\cite{SIR1995}\cite{BF2005}\cite{D2011}. The idea is to minimize $\varphi^{*}(.)$ by alternatively improving the $\mathbf{A}$ and $\mathbf{B}$ matrices through a sequence of cyclic linear least-squares optimizations. One starts with some initial guess, say, $\mathbf{A}^{0}$ and iterates from $\mathbf{A}^{0}$ to $\mathbf{B}^{0}$ then from $\mathbf{B}^{0}$ to $\mathbf{A}^{1}$, etc \dots This method of iterations yields a decreasing sequence of functions values $\lbrace \varphi^{*}(\mathbf{A}^{i},\mathbf{B}^{i}  ) \rbrace_{i \in  \mathbb{N}}$ as the sandwich inequality
\begin{equation*}
\varphi^{*}(\mathbf{A}^{i},\mathbf{B}^{i}  ) \ge  \varphi^{*}(\mathbf{A}^{i},\mathbf{B}^{i+1}  ) \ge  \varphi^{*}(\mathbf{A}^{i+1},\mathbf{B}^{i+1}  )
\end{equation*}
holds for all $(\mathbf{A}^{i},\mathbf{B}^{i}  )$, $i \in  \mathbb{N}$. Since the continuous real-valued function $\varphi^{*}(.)$ is bounded below by zero, the sequence of function values should converge to an infimum. However, the convergence can be quite slow, especially in the presence of missing values~\cite{BF2005}, and the sequence of points $\lbrace (\mathbf{A}^{i},\mathbf{B}^{i}  ) \rbrace_{i \in  \mathbb{N}}$ may even cycle, stagnate and not converge to a stationary point of the WLRA problem~\cite{GS2000} as this block separable least-squares approach is a simple instance of a block coordinate descent method (also known as the block-nonlinear Gauss-Seidel method, see~\cite{OR1970}\cite{NW2006}) for minimizing $\varphi^{*}(.)$ and the cost function is not convex~\cite{P1973}. Importantly, subsequence convergence to a stationary point can still be obtained for this  block coordinate descent algorithm applied to a nonconvex function (as $\varphi^{*}(.)$) for special cases such as the existence of an unique minimizer per block of variables~\cite{L1973}\cite{B1999} or in the case of two block of variables with a nonempty set of minimizers per block, but without the unicity condition~\cite{GS1999}, as stated in the following theorem and corollary:
\\
\begin{theo4.1} \label{theo4.1:box}
Suppose that $f(.)$ is a continuously  differentiable function over a set $\Omega \subset \mathbb{R}^{l}$, which is a Cartesian product of closed convex sets $\Omega_1, \Omega_2, ..., \Omega_m$, where $\Omega_i \subset \mathbb{R}^{n_i}$ for $i = 1,. . . ,m$ and $\sum^m_{i = 1} n_i = l$. Suppose that the variable $\mathbf{x}$ is also partitioned accordingly as $\mathbf{x} = ( \mathbf{x}_1, ..., \mathbf{x}_m)$ where $\mathbf{x}_i \in \Omega_i$. Furthermore, suppose that for each $i$ and $\mathbf{x} \in \Omega$, the solution of
\begin{equation*}
\min_{\zeta\in\Omega_i } \, \quad\  f(\mathbf{x}_1, ...,\mathbf{x}_{i-1},\zeta, \mathbf{x}_{i+1},..., \mathbf{x}_m)
\end{equation*}
is uniquely attained. If $f(.)$ is minimized by a block coordinate descent algorithm, in which a single block of variables, $\mathbf{x}_{i}$, is optimized while the remaining variables are held fixed at each iteration, then any accumulation point of the sequence of points, $\lbrace \mathbf{x}^k \rbrace_{k \in  \mathbb{N}}$ generated by this block coordinate descent algorithm is also a first-order stationary point of $f(.)$.
\end{theo4.1}
\begin{proof}
Omitted. See  p.195 in~\cite{P1973} or Proposition 2.7.1 in~\cite{B1999} for details.
\\
\end{proof}
\begin{corol4.1} \label{corol4.1:box}
In the same conditions as in Theorem~\ref{theo4.1:box}, if $f(.)$ is defined only over a Cartesian product of two closed convex sets, $\Omega_1$ and $\Omega_2$, and the global minimization of $f(.)$ with respect to each component is well defined, but not necessarily unique, then any accumulation point of the sequence of points, $\lbrace (\mathbf{x}^k,\mathbf{y}^k)  \rbrace_{k \in  \mathbb{N}}$ generated by this two-block coordinate descent algorithm is also a first-order stationary point of $f(.)$.
\end{corol4.1}
\begin{proof}
Omitted. See Theorem 6.3 in~\cite{GS1999}.
\end{proof}

As $\mathbb{R}^{p \times k}$ and $\mathbb{R}^{k \times n}$ are closed convex sets and $\varphi^{*}(.)$ is continuous as stated in Theorem~\ref{theo3.2:box} and also continuously differentiable (as it is a polynomial in $(p \times k) + (k \times n)$ variables), Corollary~\ref{corol4.1:box} can be applied to the two-block separable least-squares approach described above. Note, on the other hand, that Theorem~\ref{theo4.1:box} cannot be used here because we cannot assume that all the regression problems for computing $\mathbf{A}$ and $\mathbf{B}$ cyclically in the  two-block separable least-squares method to minimize $\varphi^{*}(.)$ can be solved uniquely in general. For example, this will not be the case if some rows or columns of $\mathbf{X}$  have less than $k$ "observed" values. However, using Corollary~\ref{corol4.1:box}, we still obtain that any  accumulation point of the sequence  $\lbrace (\mathbf{A}^i,\mathbf{B}^i)  \rbrace_{i \in  \mathbb{N}}$, say $( \widehat{\mathbf{A}}, \widehat{\mathbf{B}}  )$,  is a first-order stationary point of $\varphi^{*}(.)$ and, thus, satisfies
\begin{equation*}
 \frac{ \partial\varphi^{*}( \widehat{\mathbf{A}}, \widehat{\mathbf{B}} )}  {\partial\mathbf{a}} = \nabla \varphi^{*}_{\mathbf{a}} ( \widehat{\mathbf{A}}, \widehat{\mathbf{B}} )  = \mathbf{0}^{k.p}  \quad \text{and}  \quad   \frac{ \partial\varphi^{*}( \widehat{\mathbf{A}}, \widehat{\mathbf{B}}  )} {\partial\mathbf{b} } = \nabla \varphi^{*}_{\mathbf{b}} ( \widehat{\mathbf{A}}, \widehat{\mathbf{B}} )  = \mathbf{0}^{k.n} \ ,
\end{equation*}
where the partial functions $\varphi^{*}_{\mathbf{a}}(.)$ and $\varphi^{*}_{\mathbf{b}}(.)$ are defined by
\begin{align*}
\varphi^{*}_{\mathbf{a}} : \mathbb{R}^{p.k} \longrightarrow \mathbb{R} : \mathbf{c} \mapsto \varphi^{*}( \emph{mat}_{k \times p}( \mathbf{c} )^{T}, \mathbf{B} ) =  \frac{1}{2}   \Vert \sqrt{\mathbf{W}}  \odot ( \mathbf{X} - \mathbf{C}\mathbf{B} )  \Vert^{2}_{F} \ , \\
\varphi^{*}_{\mathbf{b}} : \mathbb{R}^{k.n} \longrightarrow \mathbb{R} : \mathbf{d} \mapsto \varphi^{*}(  \mathbf{A} , \emph{mat}_{k \times n}( \mathbf{d} )) =  \frac{1}{2}   \Vert \sqrt{\mathbf{W}}  \odot ( \mathbf{X} - \mathbf{A}\mathbf{D} )  \Vert^{2}_{F} \ ,
\end{align*}
see Theorem~\ref{theo4.3:box} below for details. Note that this condition is however not sufficient to ensure that $\widehat{\mathbf{Y}} = \widehat{\mathbf{A}}\widehat{\mathbf{B}} \in  \mathbb{R}^{p \times n}_{\le k}$ is a Frechet first-order stationary point of $\varphi(.)$ according to Theorem~\ref{theo3.7:box} if $\widehat{\mathbf{Y}}  \in  \mathbb{R}^{p \times n}_{<k}$. However, from Theorem~\ref{theo4.3:box} below, we also get that the (partial) Hessian matrices of the vectorized form of $\varphi^{*}(.)$ are equal to
\begin{align*}
\frac{ \partial^2 \varphi^{*}( \mathbf{A}, \mathbf{B} )}  {\partial^2 \mathbf{a}}  & = \nabla^{2} \varphi^{*}_{\mathbf{a}} ( \mathbf{A}, \mathbf{B} )  = \mathbf{G}(\mathbf{b})^{T} \mathbf{G}(\mathbf{b}) \ , \\
\frac{ \partial^2 \varphi^{*}( \mathbf{A}, \mathbf{B} )}  {\partial^2 \mathbf{b}}  & = \nabla^{2} \varphi^{*}_{\mathbf{b}} ( \mathbf{A}, \mathbf{B} )  = \mathbf{F}(\mathbf{a})^{T} \mathbf{F}(\mathbf{a}) \ ,
\end{align*}
and are, thus, positive  semi-definite for all  $\mathbf{A} \in \mathbb{R}^{p \times k}$ and $\mathbf{B} \in \mathbb{R}^{k \times n}$, which implies that $\varphi^{*}(.)$ is bi-convex in its whole domain. In addition, if  the block matrices $\mathbf{F}(\widehat{\mathbf{a}} )$ and $\mathbf{G}( \widehat{\mathbf{b}} ) $ are of full column-rank, which will be the rule rather than the exception if there are at least $k$ "observed" values in each column and row of the incomplete data matrix $\mathbf{X}$, the (partial) Hessian matrices $\frac{ \partial^2 \varphi^{*}( \widehat{\mathbf{A}}, \widehat{\mathbf{B}} )}  {\partial^2 \mathbf{a}}$ and $\frac{ \partial^2 \varphi^{*}( \widehat{\mathbf{A}}, \widehat{\mathbf{B}}  )} {\partial^2 \mathbf{b} }$  will be further positive definite and this implies that 
\begin{equation*}
\widehat{\mathbf{A}}  =  \text{Arg}\min_{\mathbf{A} \in \mathbb{R}^{p \times k}} \, \varphi^{*}(\mathbf{A},\widehat{ \mathbf{B} } ) \quad\  \text{ , } \quad\  \widehat{\mathbf{B}}  =  \text{Arg}\ \min_{\mathbf{B} \in \mathbb{R}^{k \times n}} \, \varphi^{*}( \widehat{\mathbf{A} }, \mathbf{B} ) \ ,
\end{equation*}
are strict local minima for the partial functions $\varphi^{*}( . ,\widehat{ \mathbf{B} } )$ and $ \varphi^{*}( \widehat{\mathbf{A} }, . )$, respectively, and
\begin{equation*}
\varphi^{*}( \widehat{\mathbf{A}}, \widehat{\mathbf{B}}  ) = \min_{\mathbf{A} \in \mathbb{R}^{p \times k}} \, \varphi^{*}(\mathbf{A},\widehat{ \mathbf{B} } ) = \min_{\mathbf{B} \in \mathbb{R}^{k \times n}} \, \varphi^{*}( \widehat{\mathbf{A} }, \mathbf{B} ) \ ,
\end{equation*}
which is a necessary, though not sufficient condition, for the pair $( \widehat{\mathbf{A}}, \widehat{\mathbf{B}} )$'s being a minimum point of $\varphi^{*}(.)$. Finally, the resulting algorithm is globally convergent with a sublinear or linear convergence rate at best~\cite{RW1980}\cite{BF2005}\cite{BT2013}. However, little can be said in general about the convergence behaviour of the sequence $\lbrace (\mathbf{A}^{i},\mathbf{B}^{i}  )\rbrace_{i \in  \mathbb{N}}$ without additional assumptions or modifications (e.g., regularizations) of the cost function of the WLRA problem as we will discuss now in some details.

If we use vectorized matrix variables, e.g., $\mathbf{a} =  \emph{vec}(  \mathbf{A}^{T}  )$ and $\mathbf{b} =  \emph{vec}(  \mathbf{B} )$, then the iterations in the block ALS algorithm (e.g., NIPALS) take the following form
\begin{align*}
\mathbf{a}^{i+1} & =  \mathbf{G}(\mathbf{b}^{i})^{+} \mathbf{z} = \omega(\mathbf{b}^{i}) \ , \\
\mathbf{b}^{i+1} & =  \mathbf{F}(\mathbf{a}^{i+1})^{+} \mathbf{x} = \upsilon(\mathbf{a}^{i+1} ) \ ,
\end{align*}
for $i = 0, 1, 2, \ldots$, where $\upsilon(.)$ and $\omega(.)$ are two real-vector functions from $\mathbb{R}^{p.k}$ to $\mathbb{R}^{k.n}$ and from $\mathbb{R}^{k.n}$ to $\mathbb{R}^{p.k}$, respectively, defined by
\begin{align*}
\upsilon (\mathbf{a} ) & =
    \begin{cases}
        \text{Arg}\min_{\mathbf{b} \in \mathbb{R}^{k.n}} \,     \frac{1}{2}   \Vert \mathbf{x} - \mathbf{F}(\mathbf{a})\mathbf{b} \Vert^{2}_2 = \varphi^{*}( \mathbf{A},\mathbf{B}  ) \\
        \text{s.t. }  \text{Arg}\min_{\mathbf{b} \in \mathbb{R}^{k.n}} \,   \Vert \mathbf{b} \Vert_{2}
    \end{cases}   \\
\omega (\mathbf{b} ) & =
    \begin{cases}
        \text{Arg}\min_{\mathbf{a} \in \mathbb{R}^{p.k}} \,     \frac{1}{2}   \Vert \mathbf{z} - \mathbf{G}(\mathbf{b})\mathbf{a} \Vert^{2}_2 = \varphi^{*}( \mathbf{A},\mathbf{B}  ) \\
        \text{s.t. }  \text{Arg}\min_{\mathbf{a} \in \mathbb{R}^{p.k}} \,   \Vert \mathbf{a} \Vert_{2}
    \end{cases} \\
\end{align*}
That is, either subproblem has an unique minimum  2-norm minimizer (see Subsection~\ref{lin_alg:box}) and the functions $\upsilon(.)$ and $\omega(.)$ are thus well-defined. In these conditions, the composition map $\chi (.) =  \omega  (.) \circ \upsilon  (.) $ is also a well-defined function from $\mathbb{R}^{p.k}$ to $\mathbb{R}^{p.k}$ and the ALS algorithm takes the form of a standard fixed point iteration~\cite{OR1970}
\begin{equation*}
\mathbf{a}^{i+1} = \chi ( \mathbf{a}^{i} ) = \chi^{i} ( \mathbf{a}^{0} ) \text{ for } i = 0, 1, 2, \ldots
\end{equation*}
Clearly, if $\upsilon(.)$ and $\omega(.)$ are continuous, $\chi(.)$ is also continuous and if, in addition,
\begin{equation*}
\lim_{i \rightarrow \infty } \mathbf{a}^{i} = \widehat{\mathbf{a} } \ , 
\end{equation*}
then $\widehat{\mathbf{a} }$ solves the system $\mathbf{a} = \chi ( \mathbf{a} )$, (e.g., $\widehat{\mathbf{a} }$ is a fixed point of  $\chi(.)$) and
\begin{equation*}
\lim_{i \rightarrow \infty } \big ( \mathbf{a}^{i},  \upsilon (\mathbf{a}^{i} ) \big) =  \lim_{i \rightarrow \infty } \big ( \mathbf{a}^{i},  \mathbf{b}^{i} \big) = ( \widehat{\mathbf{a} }, \widehat{\mathbf{b} } )   \text{ with  }  \widehat{\mathbf{b} } =   \upsilon (  \widehat{\mathbf{a} } )  \text{ and  }  \widehat{\mathbf{a} } =   \omega (  \widehat{\mathbf{b} } ) \  .
\end{equation*}
Then, under these hypotheses, the ALS iterations converge to $( \widehat{\mathbf{a} }, \widehat{\mathbf{b} } )$ and we have
\begin{equation*}
 \widehat{\mathbf{a} } = \text{Arg}\min_{\mathbf{a} \in \mathbb{R}^{p.k}} \,  \frac{1}{2}   \Vert \mathbf{z} - \mathbf{G}(\widehat{\mathbf{b} })\mathbf{a} \Vert^{2}_2    \quad\  \text{ and } \quad\   \widehat{\mathbf{b} } = \text{Arg}\min_{\mathbf{b} \in \mathbb{R}^{k.n}} \,     \frac{1}{2}   \Vert \mathbf{x} - \mathbf{F}(\widehat{\mathbf{a} })\mathbf{b} \Vert^{2}_2 \ ,
\end{equation*}
which implies
\begin{equation*}
\nabla \varphi^{*}_{\mathbf{a}} ( \widehat{\mathbf{A}}, \widehat{\mathbf{B}} )  = \frac{ \partial\varphi^{*}( \widehat{\mathbf{A}}, \widehat{\mathbf{B}} )}  {\partial\mathbf{a}} = \mathbf{0}^{k.p}  \quad \text{,}  \quad  \nabla \varphi^{*}_{\mathbf{b}} ( \widehat{\mathbf{A}}, \widehat{\mathbf{B}} )  = \frac{ \partial\varphi^{*}( \widehat{\mathbf{A}}, \widehat{\mathbf{B}}  )} {\partial\mathbf{b} } = \mathbf{0}^{k.n} \ ,
\end{equation*}
and  $(\widehat{\mathbf{A}}, \widehat{\mathbf{B}} )$ is a first-order stationary point of $\varphi^{*}(.)$. Thus, in these conditions, we have  a one-to-one correspondence between the fixed points of  $\chi(.)$ and the first-order stationary points of $\varphi^{*}(.)$. However, the hypotheses that $\upsilon(.)$ and $\omega(.)$ are continuous cannot be proved here as the generalized inverse functions $\mathbf{F}(.)^{+}$ and $\mathbf{G}(.)^{+}$ are clearly not continuous on all points of  $\mathbb{R}^{p.k}$ and $\mathbb{R}^{k.n}$, respectively, according to Theorems~\ref{theo3.10:box},~\ref{theo3.11:box} and~\ref{theo3.12:box}. Consequently, this approach cannot be used to establish the general convergence of the whole sequence  $\lbrace (\mathbf{A}^{i},\mathbf{B}^{i}  )\rbrace_{i \in  \mathbb{N}}$. Similarly, $\chi(.)$ is not a contraction in any open ball $B_{p.k}(\mathbf{a}^0, r)$ of radius $r$ around the starting point $\mathbf{a}^0$ as otherwise the Contraction Mapping Theorem~\cite{OR1970} will imply that the equation  $\mathbf{a} = \chi ( \mathbf{a} )$ has an unique solution $\widehat{\mathbf{a}}$ in the closed ball $\bar{B}_{p.k}(\mathbf{a}^0, r)$, which is false according to Remark~\ref{remark3.4:box} and the over-parameterization of the  formulation~\eqref{eq:P1} of the WLRA problem. In other words, the convergence of the whole sequence $\lbrace (\mathbf{A}^{i},\mathbf{B}^{i}  )\rbrace_{i \in  \mathbb{N}}$ cannot be proved either with the help of the Contraction Mapping Theorem.

First, we observe that more precise and stronger results can be derived when all the weights are strictly positive, e.g., when $\mathbf{W} \in \mathbb{R}^{p \times n}_{+*}$, because the cost function $\varphi(.)$ is $\lambda$-smooth in that case, which means that the gradient mapping $\nabla \varphi(.)$ from $\mathbb{R}^{p \times n}$ to $\mathbb{R}^{p \times n}$ is  Lipschitz continuous with a Lipschitz constant $\lambda > 0$, e.g.,
\begin{equation*}
 \Vert \nabla \varphi( \mathbf{Y} ) - \nabla \varphi( \mathbf{Z} ) \Vert_{F} \le \lambda   \Vert \mathbf{Y} - \mathbf{Z} \Vert_{F}  \ ,  \  \forall  \mathbf{Y}, \mathbf{Z} \in \mathbb{R}^{p \times n} \ .
\end{equation*}
Using equation~\eqref{eq:D_grad_varphi} in Subsection~\ref{landscape_wlra:box}, we get immediately
\begin{equation*}
 \Vert \nabla \varphi( \mathbf{Y} ) - \nabla \varphi( \mathbf{Z} ) \Vert_{F} =  \Vert \mathbf{W} \odot   ( \mathbf{Y} - \mathbf{Z} )  \Vert_{F}  \le  \lambda  \Vert  \mathbf{Y} - \mathbf{Z}  \Vert_{F} \ ,  \  \forall  \mathbf{Y}, \mathbf{Z} \in \mathbb{R}^{p \times n} \ ,
\end{equation*}
with $\lambda =  \max_{ (i,j) \in  \lbrack p \rbrack \times  \lbrack n \rbrack} \mathbf{W}_{ij}$, implying that $\varphi(.)$ is effectively $\lambda$-smooth when  $\mathbf{W} \in \mathbb{R}^{p \times n}_{+*}$. As $\varphi(.)$ is also bounded from below, e.g., $\forall \mathbf{Y} \in \mathbb{R}^{p \times n}, \varphi( \mathbf{Y} ) \ge 0$, we have the following result, which is a direct application of Corollary 3.9 in Olikier et al.~\cite{OUV2023}.
\begin{theo4.2} \label{theo4.2:box}
Let $ \mathbf{Y}^{i}  =   \mathbf{A}^{i} \mathbf{B}^{i}  \in \mathbb{R}^{p \times n}_{\le k}  \ ,  \  \forall  i  \in \mathbb{N}$, where the sequence  $\lbrace (\mathbf{A}^{i},\mathbf{B}^{i}  ) \rbrace_{i \in  \mathbb{N}}$ is the iterates of the block ALS algorithm defined by equations~\eqref{eq:B_nipals} and~\eqref{eq:A_nipals}.

Then, the generated sequence $\lbrace \varphi( \mathbf{Y}^{i} ) \rbrace_{i \in  \mathbb{N}} = \lbrace \varphi^{*}( \mathbf{A}^{i},\mathbf{B}^{i}  ) \rbrace_{i \in  \mathbb{N}}$ of cost function values is monotonically decreasing and converges to some value $\varphi_{*} \ge \bar{\mathbf{c}}_{\varphi^{*}} = \bar{\mathbf{c}}_{\varphi}$, where $\bar{\mathbf{c}}_{\varphi^{*}} = \bar{\mathbf{c}}_{\varphi}$ is the infimum of $\varphi(.)$ on $\mathbb{R}^{p \times n}_{\le k}$, which is equal to the infimum of  $\varphi^{*}(.)$ on $\mathbb{R}^{p \times k} \times \mathbb{R}^{k \times n}$ (see Theorem~\ref{theo3.1:box} for details). Moreover, the Riemannian gradient of $\varphi(.)$  at $\mathbf{Y}^{i}$ tends to zero, e.g.,
\begin{equation*}
\lim_{i \rightarrow \infty } \nabla_{R} \varphi( \mathbf{Y}^{i} ) =  \mathbf{P}_{\mathcal{T}_{\mathbf{Y}^{i}}  \mathbb{R}^{p \times n}_{\emph{rank}( \mathbf{Y}^{i} )} }  \big( \nabla \varphi( \mathbf{Y}^{i} )  \big) = \mathbf{0}^{p \times n}
\end{equation*}
and every point of accumulation $\widehat{\mathbf{Y}}$ of the sequence $\lbrace  \mathbf{Y}^{i}  \rbrace_{i \in  \mathbb{N}}$ satisfies $\varphi( \widehat{\mathbf{Y}} ) = \varphi_{*}$ and $\nabla_{R} \varphi( \widehat{\mathbf{Y}}  ) =  \mathbf{0}^{p \times n}$, which means that $\widehat{\mathbf{Y}}$  is a Riemannian first-order stationarity point  of $\varphi(.)$  on the smooth manifold $\mathbb{R}^{p \times n}_{k^{'}}$ where $k^{'} = \emph{rank}( \widehat{\mathbf{Y}} ) \le k$. In particular, if $\emph{rank}( \widehat{\mathbf{Y}} ) = k$ then $\widehat{\mathbf{Y}}$ is also a Frechet first-order stationarity point  of $\varphi(.)$  on $\mathbb{R}^{p \times n}_{\le k}$ in the sense of Theorem~\ref{theo3.5:box} and the pair $( \widehat{\mathbf{A}}, \widehat{\mathbf{B}} )$ is a first-order critical point of $\varphi^{*}(.)$ on $\mathbb{R}^{p \times k} \times \mathbb{R}^{k \times n}$ according to Theorem~\ref{theo3.7:box}.

Furthermore, $\forall j  \in  \mathbb{N}$, it holds that
\begin{equation*}
\min_{ 0  \le i  \le j} \Vert \nabla_{R} \varphi( \mathbf{Y}^{i} )  \Vert_{F} \le \Big (   2. \lambda  . \frac{ \varphi( \mathbf{Y}^{0} ) - \varphi_{*}  }{2.j+1}  \Big )^{\frac{1}{2}  }  \ .
\end{equation*}
In particular, given $ \varepsilon > 0$, the algorithm returns a matrix satisfying $\Vert \nabla_{R} \varphi( \mathbf{Y}^{i} )  \Vert_{F} \le  \varepsilon$ after at most $\Big \lceil  \lambda.  \frac{ \varphi( \mathbf{Y}^{0} ) - \varphi_{*}  }{ \varepsilon^{2}}  - \frac{1}{2} \Big \rceil$ iterations.
\end{theo4.2}
\begin{proof}
Omitted. See  Corollary 3.9 of Olikier et al.~\cite{OUV2023} for details.
\end{proof}
Interestingly, this theorem also illustrated the impact of a "good" initialization of the block ALS algorithm on the required number iterations for the convergence of the sequence in terms of the norm of the Riemannian gradient of $\varphi(.)$. However, even in the case where  $\mathbf{W} \in \mathbb{R}^{p \times n}_{+*}$, we cannot ensure that $\mathbf{Y}^{i+1} - \mathbf{Y}^{i} \longrightarrow \mathbf{0}^{p \times n}$ and this result requires additional modifications of the algorithm or hypotheses.

In fact, many of the past works study when the ALS minimization algorithm converges to its infimum for the matrix completion problem in polynomial time under the additional assumptions that (i) there is a solution $\widehat{\mathbf{X}} = \widehat{\mathbf{A}} \widehat{\mathbf{B}}$, which is incoherent (e.g., the squared row norms of $\widehat{\mathbf{A}}$ and squared column norms of $\widehat{\mathbf{B}}$ are not small) and (ii) the non-missing entries of $\mathbf{X}$ are selected uniformly at random or have pseudorandom properties~\cite{JNS2013}\cite{H2014}. More precisely, these two studies have shown that with an appropriate SVD-based initialization, the ALS algorithm (with a few modifications) recovers the ground-truth in the case of random binary weights and under a resampling scheme. Convergence results with a relaxation of the random sampling hypothesis can be found in~\cite{BJ2014}\cite{LLR2016}\cite{SL2016}. However, all these past studies concern mainly the matrix completion problem with a binary weight matrix~\cite{H2014}\cite{JNS2013}\cite{BJ2014}\cite{SL2016} or assume that there are no zero weights and that $\mathbf{W}$ is spectrally closed to the all one matrix in the case of a nonuniform weight matrix~\cite{LLR2016}. Finally, there is some ongoing debate as to whether these different assumptions are valid for real-world datasets~\cite{SW2015}. Interestingly, the incoherency hypothesis of the solution pair $( \widehat{\mathbf{A}}, \widehat{\mathbf{B}} )$ stated above means that $\widehat{\mathbf{A}}$ is far away from any of the barrier sets $\mathcal{B}_j$, defined in the previous subsection (see Definition~\ref{def3.2:box}), illustrating how the variable projection framework shed also some lights on the solvability of the WLRA problem by other methods such as the block ALS algorithm described above.

The block ALS method can also be adapted to solve the MMMF formulation of the WLRA problem equipped with a regularization parameter $\lambda \in \mathbb{R}_{+*}$ already discussed in Subsection~\ref{approx_wlra:box} (see equation~\eqref{eq:MMMF}), since
\begin{align*}
\min_{\mathbf{A}\in\mathbb{R}^{p \times k}\text{, }\mathbf{B}\in\mathbb{R}^{k \times n} }   \, \quad\  \varphi^{*}_{\lambda}( \mathbf{A},\mathbf{B} ) & = \frac{1}{2}   \Vert  \sqrt{\mathbf{W}}  \odot ( \mathbf{X} - \mathbf{A}\mathbf{B} )  \Vert^{2}_{F} + \frac{\lambda}{2} ( \Vert \mathbf{A} \Vert^{2}_{F} + \Vert \mathbf{B} \Vert^{2}_{F} ) \\
        & =  \frac{1}{2}   \big\Vert  \begin{bmatrix}    \mathbf{z} - \mathbf{G}(\mathbf{b})\mathbf{a}   \end{bmatrix}  \big\Vert^{2}_{2} + \frac{\lambda}{2}  \Vert \mathbf{a} \Vert^{2}_{2}  + \frac{\lambda}{2} \Vert \mathbf{b} \Vert^{2}_{2}  \\
        & =  \frac{1}{2}   \big\Vert  \begin{bmatrix}    \mathbf{x} - \mathbf{F}(\mathbf{a})\mathbf{b}   \end{bmatrix}  \big\Vert^{2}_{2} + \frac{\lambda}{2}  \Vert \mathbf{a} \Vert^{2}_{2}  + \frac{\lambda}{2} \Vert \mathbf{b} \Vert^{2}_{2} \ ,
\end{align*}
where $\mathbf{a} = \emph{vec}(  \mathbf{A}^{T} )$, $\mathbf{b} = \emph{vec}(  \mathbf{B} )$, $\mathbf{x}, \mathbf{z}$ , $\mathbf{F}(\mathbf{a})$ and $\mathbf{G}(\mathbf{b})$ are defined as above.

In this case, the block ALS algorithm computes alternatively the solutions of the two regularized least-squares problems
\begin{equation*}
\text{Arg}\min_{\mathbf{a}\in\mathbb{R}^{p.k}} \,   \frac{1}{2}  \big\Vert \mathbf{z} - \mathbf{G}(\mathbf{b})\mathbf{a} \big\Vert^{2}_{2} +  \frac{\lambda}{2} \big\Vert \mathbf{a} \big\Vert^{2}_{2} =    \frac{1}{2}    \Big\Vert  \begin{bmatrix}    \mathbf{z} - \mathbf{G}(\mathbf{b})\mathbf{a}   \\    \sqrt{\lambda} . \mathbf{a}   \end{bmatrix}  \Big\Vert^{2}_{2} =    \frac{1}{2} \Big\Vert  \begin{bmatrix} \mathbf{z}   \\ \mathbf{0}^{k.p} \end{bmatrix} - \begin{bmatrix} \mathbf{G}(\mathbf{b}) \\  \sqrt{\lambda} . \mathbf{I}_{p.k}  \end{bmatrix} \mathbf{a}   \Big\Vert^{2}_{2} \ ,
\end{equation*}
and
\begin{equation*}
\text{Arg}\min_{\mathbf{b}\in\mathbb{R}^{k.n}} \,    \frac{1}{2} \big\Vert \mathbf{x} - \mathbf{F}(\mathbf{a})\mathbf{b} \big\Vert^{2}_2 + \frac{\lambda}{2}  \big\Vert \mathbf{b}\big \Vert^{2}_{2} =    \frac{1}{2}    \Big\Vert  \begin{bmatrix}    \mathbf{x} - \mathbf{F}(\mathbf{a})\mathbf{b}   \\     \sqrt{\lambda} . \mathbf{b}  \end{bmatrix}  \Big\Vert^{2}_{2} =    \frac{1}{2} \Big\Vert  \begin{bmatrix} \mathbf{x}   \\ \mathbf{0}^{k.n} \end{bmatrix} -  \begin{bmatrix} \mathbf{F}(\mathbf{a}) \\  \sqrt{\lambda} . \mathbf{I}_{k.n}  \end{bmatrix} \mathbf{b}  \Big\Vert^{2}_2 \ .
\end{equation*}
In other words, the MMMF ALS algorithm updates $\mathbf{a}$ and $\mathbf{b}$ at the $i+1$ iteration according to the rules
\begin{equation*}
\mathbf{a}^{i+1} = \left(  \mathbf{G}(\mathbf{b}^{i} )^{T} \mathbf{G}(\mathbf{b}^{i}) + \lambda . \mathbf{I}_{p.k} \right)^{-1} \mathbf{G}(\mathbf{b}^{i})^{T} \mathbf{z}
\end{equation*}
and
\begin{equation*}
\mathbf{b}^{i+1}  =  \left( \mathbf{F}(\mathbf{a}^{i+1})^{T} \mathbf{F}(\mathbf{a}^{i+1}) + \lambda . \mathbf{I}_{k.n} \right)^{-1} \mathbf{F}(\mathbf{a}^{i+1})^{T} \mathbf{x} \ .
\end{equation*}
Furthermore, Theorem~\ref{theo4.1:box} can now be applied directly to this regularized ALS algorithm  in order to show that any accumulation point of the sequence  $\lbrace (\mathbf{A}^{i},\mathbf{B}^{i})  \rbrace_{i \in  \mathbb{N}}$, say $( \widehat{\mathbf{A}}, \widehat{\mathbf{B}}  )$, is a stationary point of $\varphi^{*}_{\lambda}(.)$ as we are now sure that all the regression subproblems for computing $\mathbf{A}^{i}$ and $\mathbf{B}^{i}$ can be solved uniquely because of the presence of the regularization terms in $\varphi^{*}_{\lambda}(.)$.

Next, from Theorem~\ref{theo4.3:box} and its corollary (see below), we deduce that the partial Hessian matrices $ \frac{ \partial^2 \varphi^{*}_{\lambda}( \mathbf{A}, \mathbf{B} )}  {\partial^2 \mathbf{a}} =  \nabla^{2} ( \varphi^{*}_{\lambda} )_{\mathbf{a}} ( \mathbf{A}, \mathbf{B} ) $ and $ \frac{ \partial^2 \varphi^{*}_{\lambda}( \mathbf{A}, \mathbf{B} )}  {\partial^2 \mathbf{b}} = \nabla^{2} (\varphi^{*}_{\lambda} )_{\mathbf{b}} ( \mathbf{A}, \mathbf{B} ) $ are positive definite for all  $\mathbf{A} \in \mathbb{R}^{p \times k}$ and $\mathbf{B} \in \mathbb{R}^{k \times n}$ as soon as $\lambda>0$, which implies that $\varphi^{*}_{\lambda}(.)$ is now strongly bi-convex in its whole domain instead of only bi-convex as  $\varphi^{*}(.)$. Using the facts that  $\varphi^{*}_{\lambda}(.)$ is  also a coercive (thanks to the inclusion of the regularization term $\frac{\lambda}{2}  \Vert \mathbf{a} \Vert^{2}_{2}  + \frac{\lambda}{2} \Vert \mathbf{b} \Vert^{2}_{2}$) and real-analytic (as it is a polynomial in $(p \times k) + (k \times n)$ variables) function, it can be demonstrated that this strongly bi-convex cost function also verifies the so-called  Kurdyka-Lojasiewicz inequality, the sequence $( \mathbf{A}^{i}, \mathbf{B}^{i} )$ is bounded and that the whole sequence $( \mathbf{A}^{i} , \mathbf{B}^{i} )$ generated by the MMMF ALS algorithm converges to a first-order stationary point of $\varphi^{*}_{\lambda}(.)$, say $( \widehat{\mathbf{A}} , \widehat{\mathbf{B}}  )$~\cite{XY2013}, which is a much stronger result than the one delivered by Theorem~\ref{theo4.1:box} and its corollary.

Finally, Li et al.~\cite{LZT2019b}, using results from ~\cite{ABRS2010}\cite{XY2013}, were able to demonstrate recently that the sequence $( \mathbf{A}^{i} , \mathbf{B}^{i} )$ generated by the following proximal version of the MMMF ALS algorithm
\begin{align*}
\mathbf{a}^{i+1} & =  \text{Arg}\min_{\mathbf{a}\in\mathbb{R}^{p.k}} \,   \frac{1}{2}  \big\Vert \mathbf{z} - \mathbf{G}(\mathbf{b}^{i} )\mathbf{a} \big\Vert^{2}_2 +  \frac{\lambda}{2} \big\Vert \mathbf{a} \big\Vert^{2}_{2} + \frac{\beta}{2} \big\Vert \mathbf{a}^{i} - \mathbf{a} \big\Vert^{2}_{2}  \\
& =     \text{Arg}\min_{\mathbf{a}\in\mathbb{R}^{p.k}} \, \frac{1}{2}  \Big\Vert   \begin{bmatrix}    \mathbf{z} - \mathbf{G}(\mathbf{b}^{i} )\mathbf{a}   \\    \sqrt{\lambda} . \mathbf{a}   \\    \sqrt{\beta}  ( \mathbf{a}^{i}   - \mathbf{a} ) \end{bmatrix}  \Big\Vert^{2}_{2}    \\
& =     \text{Arg}\min_{\mathbf{a}\in\mathbb{R}^{p.k}} \, \frac{1}{2} \Big\Vert  \begin{bmatrix} \mathbf{z}   \\ \mathbf{0}^{p.k}  \\   \sqrt{\beta} . \mathbf{a}^{i}   \end{bmatrix} - \begin{bmatrix} \mathbf{G}(\mathbf{b}^{i} ) \\  \sqrt{\lambda} . \mathbf{I}_{p.k}  \\  \sqrt{\beta} .\mathbf{I}_{p.k} \end{bmatrix} \mathbf{a}   \Big\Vert^{2}_2 \\
& =     \left(  \mathbf{G}(\mathbf{b}^{i} )^{T} \mathbf{G}(\mathbf{b}^{i}) + ( \lambda + \beta ) \mathbf{I}_{p.k} \right)^{-1} \left( \mathbf{G}(\mathbf{b}^{i})^{T} \mathbf{z} + \beta . \mathbf{a}^{i} \right)
\end{align*}
and
\begin{align*}
\mathbf{b}^{i+1} & =  \text{Arg}\min_{\mathbf{b}\in\mathbb{R}^{k.n}} \,   \frac{1}{2}  \big\Vert \mathbf{x} - \mathbf{F}(\mathbf{a}^{i+1} )\mathbf{b} \big\Vert^{2}_2 +  \frac{\lambda}{2} \big\Vert \mathbf{b} \big\Vert^{2}_{2} + \frac{\beta}{2} \big\Vert \mathbf{b}^{i} - \mathbf{b} \big\Vert^{2}_{2}  \\
& =     \text{Arg}\min_{\mathbf{b}\in\mathbb{R}^{k.n}} \, \frac{1}{2}    \Big\Vert  \begin{bmatrix}    \mathbf{x} - \mathbf{F}(\mathbf{a}^{i+1} )\mathbf{b}   \\    \sqrt{\lambda} . \mathbf{b}   \\    \sqrt{\beta} ( \mathbf{b}^{i}   - \mathbf{b} ) \end{bmatrix}  \Big\Vert^{2}_{2}    \\
& =     \text{Arg}\min_{\mathbf{b}\in\mathbb{R}^{k.n}} \, \frac{1}{2} \Big\Vert  \begin{bmatrix} \mathbf{x}   \\ \mathbf{0}^{k.n}  \\   \sqrt{\beta} . \mathbf{b}^{i}   \end{bmatrix} - \begin{bmatrix} \mathbf{F}(\mathbf{a}^{i+1} ) \\  \sqrt{\lambda} . \mathbf{I}_{k.n}  \\  \sqrt{\beta} . \mathbf{I}_{k.n} \end{bmatrix} \mathbf{b}   \Big\Vert^{2}_2 \\
& =    \left(  \mathbf{F}(\mathbf{a}^{i+1} )^{T} \mathbf{F}(\mathbf{a}^{i+1}) + ( \lambda + \beta ) \mathbf{I}_{k.n} \right)^{-1} \left( \mathbf{F}(\mathbf{a}^{i+1})^{T} \mathbf{x} + \beta . \mathbf{b}^{i} \right) ,
\end{align*}
where
\begin{equation*}
\beta > 8. \Vert \mathbf{W} \Vert^{2}_{S}  \varphi^{*}_{\lambda}( \mathbf{A}^{0},\mathbf{B}^{0} )/\lambda + 4. \Vert \mathbf{W} \Vert_{S} \sqrt{ \varphi^{*}_{\lambda}( \mathbf{A}^{0},\mathbf{B}^{0} ) } + \lambda,
\end{equation*}
converges not only to a first-order stationary point, but in fact to a second-order stationary point of $\varphi^{*}_{\lambda}(.)$ (see Proposition 4 and example 3 in Section 4.3 of~\cite{LZT2019b}), e.g., to a point $( \widehat{\mathbf{A}} , \widehat{\mathbf{B}} )$ which verifies
\begin{equation*}
\nabla\varphi^{*}_{\lambda}( \widehat{\mathbf{A}} , \widehat{\mathbf{B}} ) = ( \mathbf{0}^{p \times k} , \mathbf{0}^{k \times n} )  \text{ and } \big ( \nabla^{2} \varphi^{*}_{\lambda}( \widehat{\mathbf{A}} , \widehat{\mathbf{B}} ) \big ) \big( ( \mathbf{C} , \mathbf{D} ) , ( \mathbf{C} , \mathbf{D} ) \big )\ge 0 ,
\end{equation*}
$\forall ( \mathbf{C} , \mathbf{D} ) \in \mathbb{R}^{p \times k} \times \mathbb{R}^{k \times n}$,
e.g.,  $\nabla^{2} \varphi^{*}_{\lambda}( \widehat{\mathbf{A}} , \widehat{\mathbf{B}} )$ is a positive semi-definite (symmetric) matrix. Importantly, if  $\varphi^{*}_{\lambda}(.)$ is well-conditioned (e.g., depending on the form of the weight matrix $\mathbf{W}$), these second-order stationary points may correspond to a local or even global optimal solution, see~\cite{ZLTW2018} and Theorem 3.10 of~\cite{OUV2023} for more information.

\begin{remark4.1} \label{remark4.1:box}
An interesting and open question is to determine if these strong first- and second-order convergence properties of the ALS method for solving the MMMF formulation of the WLRA problem may also extend to the cost function $g_{\lambda}(.)$ proposed by Boumal and Absil~\cite{BA2011}\cite{BA2015}  and discussed in Subsection~\ref{approx_wlra:box}  (see equation~\eqref{eq:g_func}). $\blacksquare$
\\
\end{remark4.1}

The block ALS algorithm and its MMMF variant have also been incorporated as a building block in various Expectation-Maximization or other first-order methods to increase their efficiency for large datasets by avoiding costly SVD computations in high dimensions~\cite{JHJ2009}\cite{HMLZ2015}\cite{TH2021}.

Interestingly, we note that Szlam et al.~\cite{STT2017} have recently demonstrated that only a few iterations of such ALS are sufficient to produce nearly optimal spectral- and Frobenius-norm accuracies of low-rank approximations to a matrix when all the weights $\mathbf{W}_{ij}$ are equal to one, provided that $\mathbf{A}_{0}$ is one of the random matrices used by~\cite{HMT2011} (for example, the entries of  $\mathbf{A}_{0}$ can be independent and identically distributed standard normal variates) and that iterating until convergence is unnecessary. Extending their demonstration to the case when the weights $\mathbf{W}_{ij}$ are unequal (and eventually with some equal to zero) is an interesting issue already discussed in~\cite{RSW2016}\cite{BWZ2019}, but is outside the scope of this paper. However, we highlight again that proper initialization of the ALS or variable projection methods described here is obviously an important topic, which also needs a careful attention~\cite{GZ1979}\cite{JNS2013}\cite{H2014}\cite{SL2016}\cite{SJ2004}. As an illustration,~\cite{JNS2013}\cite{H2014}\cite{SL2016}\cite{RSW2016} showed that given a good enough initialization, many simple local search algorithms, like ALS, succeed, a result which is consistent with Theorem~\ref{theo4.2:box} above.

Now, let us consider how to compute efficiently the first- and second-order derivatives of the vectorized form  of $\varphi^{*}(.)$  (and $\varphi^{*}_{\lambda}(.)$)  in order to obtain meaningful tests of convergence of these ALS methods to a (local) minimizer of this cost function. We already know from Subsection~\ref{noconv_wlra:box}  that the  objective function $\varphi^{*}(.)$ used in the~\eqref{eq:P1} formulation of the WLRA problem,
\begin{equation*}
\varphi^{*} : \mathbb{R}^{p \times k} \times  \mathbb{R}^{k \times n} \longrightarrow \mathbb{R} : ( \mathbf{A}, \mathbf{B} ) \mapsto \varphi^{*}( \mathbf{A}, \mathbf{B} ) =  \frac{1}{2}   \Vert \sqrt{\mathbf{W}}  \odot ( \mathbf{X} - \mathbf{A}\mathbf{B} )  \Vert^{2}_{F} \ ,
\end{equation*}
is  $C^{\infty}$  differentiable over its domain of definition. Furthermore, we have also already derived the first- and second-order derivatives of $\varphi^{*}(.)$ in equations~\eqref{eq:D_grad_varphi*} and~\eqref{eq:D_hess_varphi*}, respectively. As the vectorized form of  $\varphi^{*}(.)$ is defined by the composition of  $\varphi^{*}(.)$ with the linear mapping
\begin{equation*}
\mathbb{R}^{p.k} \times  \mathbb{R}^{k.n} \longrightarrow \mathbb{R}^{p \times k} \times  \mathbb{R}^{k \times n} : ( \mathbf{a}, \mathbf{b} ) \mapsto \big (  \emph{mat}_{k \times p}( \mathbf{a} )^{T},  \emph{mat}_{k \times n}( \mathbf{b} ) \big) = ( \mathbf{A}, \mathbf{B} ) \ ,
\end{equation*}
it is also  $C^{\infty}$  differentiable over its domain of definition, $\mathbb{R}^{p.k} \times  \mathbb{R}^{k.n} $, and we have the following results concerning the vectorized forms of the first- and second-order derivatives of  $\varphi^{*}(.)$, which offer more convenient expressions for checking the first- and second-KKT conditions of $\varphi^{*}(.)$  at a given pair of $\mathbb{R}^{p \times k} \times  \mathbb{R}^{k \times n}$ then equations~\eqref{eq:D_grad_varphi*} and~\eqref{eq:D_hess_varphi*}. 
\\
\begin{theo4.3} \label{theo4.3:box}
 For $\mathbf{X} \in \mathbb{R}^{p \times n}$,  $\sqrt{\mathbf{W}} \in \mathbb{R}^{p \times n}_+$ and any fixed integer $ k \le \emph{rank}( \mathbf{X} ) \le \text{min}( {p},{n} )$, the vectorized partial first-derivatives of $\varphi^{*}(.)$ with respect to $\mathbf{a} = \emph{vec}( \mathbf{A}^{T})$ and $\mathbf{b} = \emph{vec}( \mathbf{B} )$ are equal, respectively, to
\begin{equation} \label{eq:Da_varphi*}
\frac{ \partial   \varphi^{*} ( \mathbf{A}, \mathbf{B} ) }  {\partial\mathbf{a}} =   \nabla \varphi^{*}_{\mathbf{a}} ( \mathbf{A}, \mathbf{B} )  = \mathbf{G}(\mathbf{b})^{T} \mathbf{G}(\mathbf{b})\mathbf{a} - \mathbf{G}(\mathbf{b})^{T} \mathbf{z} 
\end{equation}
and
\begin{equation} \label{eq:Db_varphi*}
\frac{ \partial   \varphi^{*} ( \mathbf{A}, \mathbf{B} ) }  {\partial\mathbf{b}} =  \nabla \varphi^{*}_{\mathbf{b}} ( \mathbf{A}, \mathbf{B} )  = \mathbf{F}(\mathbf{a})^{T} \mathbf{F}(\mathbf{a})\mathbf{b} - {\mathbf{F}(\mathbf{a})} ^{T} \mathbf{x} \ ,
\end{equation}
where
\begin{align*}
\mathbf{F}(  \mathbf{a} ) & =  \bigoplus_{j=1}^n \mathbf{F}_{j}(  \mathbf{a} ) = \bigoplus_{j=1}^n  \emph{diag}(\sqrt{\mathbf{W}}_{.j}) \big(   \emph{mat}_{k \times p} (\mathbf{a})  \big)^{T} \ ,  \\
\mathbf{G}(  \mathbf{b} ) & =  \bigoplus_{i=1}^p   \mathbf{G}_{i}(  \mathbf{b} ) =  \bigoplus_{i=1}^p  \emph{diag}(\sqrt{\mathbf{W}}_{i.}) \big( \emph{mat}_{k \times n} ( \mathbf{b} ) \big)^{T}  \ , \\
\mathbf{x} & = \emph{vec} ( \sqrt{\mathbf{W}} \odot \mathbf{X}  )  \text{ and } \mathbf{z} = \emph{vec} \big( (\sqrt{\mathbf{W}}  \odot \mathbf{X})^{T}  \big) \  .
\end{align*}

Moreover, we have
\begin{align} \label{eq:DAB_varphi*}
 \nabla \varphi^{*}_{\mathbf{a}} ( \mathbf{A}, \mathbf{B} )  & =  \emph{vec}  \big(  \nabla \varphi^{*}_{\mathbf{A}} ( \mathbf{A}, \mathbf{B} )  \big)  \ , \nonumber \\
 \nabla \varphi^{*}_{\mathbf{b}} ( \mathbf{A}, \mathbf{B} )  & = \emph{vec}  \big(   \nabla \varphi^{*}_{\mathbf{B}} ( \mathbf{A}, \mathbf{B} ) \big)  \ ,
\end{align}
where $\nabla \varphi^{*}_{\mathbf{A}} ( \mathbf{A}, \mathbf{B} )$ and $\nabla \varphi^{*}_{\mathbf{B}} ( \mathbf{A}, \mathbf{B} )$ are defined in equation~\eqref{eq:D_partial_grad_varphi*}.

The vectorized second-derivative (symmetric) matrix form of $\varphi^{*}(.)$  is given by
\begin{align} \label{eq:D2_varphi*}
\big \lbrack \nabla^{2} \varphi^{*} ( \mathbf{A}, \mathbf{B} ) \big \rbrack & = 
 \begin{bmatrix}
\frac{ \partial^{2}   \varphi^{*} ( \mathbf{A}, \mathbf{B} ) }  {\partial^{2}\mathbf{a}} & \frac{ \partial^{2}   \varphi^{*} ( \mathbf{A}, \mathbf{B} ) }  {\partial\mathbf{a} \partial\mathbf{b} } \\
\frac{ \partial^{2}   \varphi^{*}( \mathbf{A}, \mathbf{B} ) }  {\partial\mathbf{b} \partial\mathbf{a} } & \frac{ \partial^{2}   \varphi^{*} ( \mathbf{A}, \mathbf{B} ) }  {\partial^{2} \mathbf{b}  }  
\end{bmatrix}  \nonumber \\
& =
 \begin{bmatrix}
\nabla^{2} \varphi^{*}_{\mathbf{a}} ( \mathbf{A}, \mathbf{B} ) & \nabla^{2} \varphi^{*}_{\mathbf{a} , \mathbf{b}} ( \mathbf{A}, \mathbf{B} ) \\
\nabla^{2} \varphi^{*}_{\mathbf{b} , \mathbf{a}} ( \mathbf{A}, \mathbf{B} ) & \nabla^{2} \varphi^{*}_{\mathbf{b}}( \mathbf{A}, \mathbf{B} )  
\end{bmatrix} \  ,
\end{align} 
where
\begin{align*}
\nabla^{2} \varphi^{*}_{\mathbf{a}} ( \mathbf{A}, \mathbf{B} )  & = \frac{ \partial^2 \varphi^{*}( \mathbf{A}, \mathbf{B} )}  {\partial^2 \mathbf{a}} = \mathbf{G}(\mathbf{b})^{T} \mathbf{G}(\mathbf{b}) \ , \\
\nabla^{2} \varphi^{*}_{\mathbf{b}} ( \mathbf{A}, \mathbf{B} )  & = \frac{ \partial^2 \varphi^{*}( \mathbf{A}, \mathbf{B} )}  {\partial^2 \mathbf{b}} = \mathbf{F}(\mathbf{a})^{T} \mathbf{F}(\mathbf{a}) \ , \\
\nabla^{2} \varphi^{*}_{\mathbf{b} , \mathbf{a}} ( \mathbf{A}, \mathbf{B} ) & =  \frac{ \partial^{2}   \varphi^{*} ( \mathbf{A}, \mathbf{B} ) }  {\partial\mathbf{a} \partial\mathbf{b} } =  \Big ( \big ( \mathbf{W} \odot ( \mathbf{A} \mathbf{B} - \mathbf{X} ) \big )^{T}  \otimes  \mathbf{I}_{k} \Big ) + \mathbf{F}(\mathbf{a})^{T}  \mathbf{K}_{(n,p)} \mathbf{G}( \mathbf{b}  ) \ , \\
\nabla^{2} \varphi^{*}_{\mathbf{a} , \mathbf{b}} ( \mathbf{A}, \mathbf{B} ) & = \big \lbrack \nabla^{2} \varphi^{*}_{\mathbf{b} , \mathbf{a}} ( \mathbf{A}, \mathbf{B} ) \big \rbrack^{T} \ .
\end{align*}

Finally, we have the following equalities, which precise the relationships between the quadratic forms $\big ( \nabla^{2} \varphi^{*} ( \mathbf{A}, \mathbf{B} ) \big )$ and $\big ( \nabla^{2} \varphi( \mathbf{A} \mathbf{B} ) \big )$ in complement of equations~\eqref{eq:D_partial_hess_varphi*} and~\eqref{eq:D_hess_varphi*}
\begin{align*}
\mathbf{c}^{T} \nabla^{2} \varphi^{*}_{\mathbf{a}} ( \mathbf{A}, \mathbf{B} ) \mathbf{c} & =  \big ( \nabla^{2} \varphi^{*}_{\mathbf{A}} ( \mathbf{A} , \mathbf{B} ) \big ) ( \mathbf{C} , \mathbf{C} ) = \big ( \nabla^{2} \varphi ( \mathbf{A} \mathbf{B} )  \big ) ( \mathbf{C}  \mathbf{B} , \mathbf{C}  \mathbf{B}  ) \ ,  \\
\mathbf{d}^{T} \nabla^{2} \varphi^{*}_{\mathbf{b}} ( \mathbf{A}, \mathbf{B} ) \mathbf{d} & = \big ( \nabla^{2} \varphi^{*}_{\mathbf{B}} ( \mathbf{A} , \mathbf{B} ) \big ) ( \mathbf{D} , \mathbf{D} )  = \big ( \nabla^{2} \varphi ( \mathbf{A} \mathbf{B} )  \big ) ( \mathbf{A}  \mathbf{D} , \mathbf{A}  \mathbf{D}  ) \ ,   \\
\mathbf{d}^{T} \nabla^{2} \varphi^{*}_{\mathbf{b} , \mathbf{a}} ( \mathbf{A}, \mathbf{B} )  \mathbf{c} & = \big \langle \nabla \varphi ( \mathbf{A} \mathbf{B} ) , \mathbf{C}  \mathbf{D} \big \rangle_{F} + \big ( \nabla^{2} \varphi ( \mathbf{A} \mathbf{B} )  \big ) ( \mathbf{A} \mathbf{D}  , \mathbf{C}  \mathbf{B} ) \ ,
\end{align*}
$\forall  \mathbf{C} \in \mathbb{R}^{p \times k}$  with $\mathbf{c} = \emph{vec} (  \mathbf{C}^{T} )$ and $\forall  \mathbf{D} \in \mathbb{R}^{k \times n}$  with $\mathbf{d} = \emph{vec} (  \mathbf{D} )$.
\end{theo4.3}
\begin{proof}
First, we observe that the matrix of first-derivatives of the vectorized residual function $\mathbf{e}( \mathbf{a}, \mathbf{b} ) = \mathbf{e}( \mathbf{A}, \mathbf{B} ) = \mathbf{x} - \mathbf{F}(\mathbf{a})\mathbf{b}$ with respect to $\mathbf{b}$ ($\mathbf{e}( \mathbf{A}, \mathbf{B} )$ is first defined in equation~\eqref{eq:res_func})  is simply
\begin{equation*}
\frac{ \partial  \mathbf{e}( \mathbf{A}, \mathbf{B} )}  {\partial\mathbf{b}} = - \mathbf{F}(\mathbf{a})
\end{equation*}
and is very sparse with only $k$ non-zero elements in each row as $\mathbf{F}(\mathbf{a})$ is a block diagonal matrix (see equation~\eqref{eq:F_mat}). Since
\begin{equation*}
\varphi^{*}( \mathbf{A}, \mathbf{B} ) = \frac{1}{2}  {\mathbf{e}( \mathbf{A}, \mathbf{B} )}^{T} {\mathbf{e}( \mathbf{A}, \mathbf{B} )} = \frac{1}{2} \Vert \mathbf{e}( \mathbf{a}, \mathbf{b} ) \Vert^{2} _{2}  \ ,
\end{equation*}
the derivative of  $\varphi^{*}( \mathbf{A}, \mathbf{B} )$ with respect to $\mathbf{b}$ is then easy to compute, using a standard differential rule for a mapping of the form  $\mathbb{R}^{k.n} \longrightarrow  \mathbb{R} : \mathbf{d} \mapsto  \frac{1}{2} \Vert g( \mathbf{d} ) \Vert^{2}_{2}$, where $g(.)$ is a differentiable mapping from $\mathbb{R}^{k.n}$ to $\mathbb{R}^{p.n}$~\cite{C2017},
\begin{equation*}
\frac{ \partial   \varphi^{*}( \mathbf{A}, \mathbf{B} ) }  {\partial\mathbf{b}} = - {\mathbf{F}(\mathbf{a})}^{T}  \big( \mathbf{x} - \mathbf{F}(\mathbf{a})\mathbf{b} \big) = {\mathbf{F}(\mathbf{a})} ^{T} \mathbf{F}(\mathbf{a})\mathbf{b} - {\mathbf{F}(\mathbf{a})} ^{T} \mathbf{x} \ .
\end{equation*}
For computing the derivative of $\varphi^{*}( \mathbf{A}, \mathbf{B} )$ with respect to $\mathbf{a}$, we first recall that the vectorized residual function $\mathbf{e}( \mathbf{a}, \mathbf{b} )$ may also be expressed in the alternative form
\begin{equation*}
\mathbf{e}( \mathbf{A}, \mathbf{B} ) = \mathbf{x} - \mathbf{K}_{(n,p)}\mathbf{G}(\mathbf{b})\mathbf{a} = \mathbf{K}_{(n,p)} \big(  \mathbf{z} - \mathbf{G}(\mathbf{b})\mathbf{a} \big) \ ,
\end{equation*}
see the paragraph after equation~\eqref{eq:G_mat} in Subsection~\ref{varpro_wlra:box} for details. Hence
\begin{equation*}
\frac{ \partial  \mathbf{e}( \mathbf{A}, \mathbf{B} )}  {\partial\mathbf{a}} = - \mathbf{K}_{(n,p)} \mathbf{G}(\mathbf{b})
\end{equation*}
and this matrix of derivatives with respect to $\mathbf{a}$ is also very sparse with only $k$ non-zero elements in each row. Now the derivative of $\varphi^{*}( \mathbf{A}, \mathbf{B} )$ with respect to $\mathbf{a}$, is also very simple to obtain, using the same differentiation rule as above and properties of the commutation matrix given in Subsection~\ref{multlin_alg:box},
\begin{align*}
\frac{ \partial   \varphi^{*}( \mathbf{A}, \mathbf{B} ) }  {\partial\mathbf{a}} & =  -\big(  \mathbf{K}_{(n,p)} \mathbf{G}(\mathbf{b})  \big)^{T}  \mathbf{K}_{(n,p)} \big(  \mathbf{z} - \mathbf{G}(\mathbf{b})\mathbf{a} \big)  \\
              & = - \mathbf{G}(\mathbf{b})^{T} \mathbf{K}_{(p,n)} \mathbf{K}_{(n,p)} \big( \mathbf{z} - \mathbf{G}(\mathbf{b})\mathbf{a} \big)  \\
              & =  \mathbf{G}(\mathbf{b})^{T} \mathbf{G}(\mathbf{b})\mathbf{a} - \mathbf{G}(\mathbf{b})^{T} \mathbf{z} \ .
\end{align*}
Next, to demonstrate that $\nabla \varphi^{*}_{\mathbf{b}} ( \mathbf{A}, \mathbf{B} )  =  \emph{vec}  \big(  \nabla \varphi^{*}_{\mathbf{B}} ( \mathbf{A}, \mathbf{B} )  \big)$, we observe that, by definition, we have
\begin{equation*}
\mathbf{F}(\mathbf{a}) \mathbf{b} =  \emph{diag}\big( \emph{vec}(\sqrt{\mathbf{W}}) \big) ( \mathbf{I}_{n} \otimes  \mathbf{A} ) \emph{vec}( \mathbf{B} ) =  \emph{diag}\big( \emph{vec}(\sqrt{\mathbf{W}}) \big)  \emph{vec}(  \mathbf{A} \mathbf{B} )
\end{equation*}
and, thus,
\begin{align*}
\mathbf{F}(\mathbf{a}) \mathbf{b} -  \mathbf{x} & =  \emph{diag} \big( \emph{vec}(\sqrt{\mathbf{W}}) \big)  \emph{vec}(  \mathbf{A} \mathbf{B} ) - \emph{vec} ( \sqrt{\mathbf{W}} \odot \mathbf{X}  ) \\
                                                                          & =  \emph{diag} \big( \emph{vec}(\sqrt{\mathbf{W}}) \big) \big(  \emph{vec}(  \mathbf{A} \mathbf{B} )  - \emph{vec} ( \mathbf{X} )  \big)  \\
                                                                          & =   \emph{diag} \big( \emph{vec}(\sqrt{\mathbf{W}}) \big)  \emph{vec}(  \mathbf{A} \mathbf{B}  - \mathbf{X} )  \ ,
\end{align*}
which implies that
\begin{align*}
\nabla \varphi^{*}_{\mathbf{b}} ( \mathbf{A}, \mathbf{B} )  & =  {\mathbf{F}(\mathbf{a})}^{T} \big ( \mathbf{F}(\mathbf{a})\mathbf{b} -  \mathbf{x} \big ) \\
                                                                                           & = ( \mathbf{I}_{n} \otimes  \mathbf{A} )^{T}   \emph{diag}\big( \emph{vec}( \mathbf{W}) \big)  \emph{vec}(  \mathbf{A} \mathbf{B}  - \mathbf{X} ) \\
                                                                                           & = ( \mathbf{I}_{n} \otimes  \mathbf{A}^{T} )   \emph{diag}\big( \emph{vec}( \mathbf{W}) \big)  \emph{vec}(  \mathbf{A} \mathbf{B}  - \mathbf{X} ) \\
                                                                                           & = ( \mathbf{I}_{n} \otimes  \mathbf{A}^{T} )   \emph{vec} \big(  \mathbf{W}  \odot ( \mathbf{A} \mathbf{B}  - \mathbf{X} ) \big) \\
                                                                                           & = \emph{vec} \Big(   \mathbf{A}^{T} \big(  \mathbf{W}  \odot ( \mathbf{A} \mathbf{B}  - \mathbf{X} ) \big)  \Big) ,
\end{align*}
and, using equation~\eqref{eq:D_partial_grad_varphi*}, we conclude that
\begin{equation*}
\nabla \varphi^{*}_{\mathbf{b}} ( \mathbf{A}, \mathbf{B} )  =  \emph{vec}  \big(  \nabla \varphi^{*}_{\mathbf{B}} ( \mathbf{A}, \mathbf{B} )  \big) \ .
\end{equation*}

Similarly, for demonstrating that $\nabla \varphi^{*}_{\mathbf{a}} ( \mathbf{A}, \mathbf{B} )  =  \emph{vec}  \big(  \nabla \varphi^{*}_{\mathbf{A}} ( \mathbf{A}, \mathbf{B} )  \big)$, we observe that
\begin{equation*}
\mathbf{G}(\mathbf{b}) \mathbf{a} =   \emph{diag} \big( \emph{vec}(\sqrt{\mathbf{W}}^{T} ) \big)  ( \mathbf{I}_{p}  \otimes \mathbf{B}^{T} ) \emph{vec}( \mathbf{A}^{T} ) = \emph{diag} \big( \emph{vec}(\sqrt{\mathbf{W}}^{T} ) \big)  \emph{vec} \big( (\mathbf{A}\mathbf{B})^{T}  \big)
\end{equation*}
and, thus,
\begin{align*}
\mathbf{G}(\mathbf{b}) \mathbf{a} -  \mathbf{z} & = \emph{diag} \big( \emph{vec}(\sqrt{\mathbf{W}}^{T} ) \big)  \emph{vec} \big( (\mathbf{A}\mathbf{B})^{T}  \big) -  \emph{vec} \big( (\sqrt{\mathbf{W}}  \odot \mathbf{X})^{T} \big) \\
                                                                          & =  \emph{diag} \big( \emph{vec}(\sqrt{\mathbf{W}}^{T} ) \big) \Big(   \emph{vec} \big( (\mathbf{A}\mathbf{B})^{T}  \big)  - \emph{vec} ( \mathbf{X}^{T} )  \Big)  \\
                                                                          & =   \emph{diag} \big( \emph{vec}(\sqrt{\mathbf{W}}^{T} ) \big)  \emph{vec} \big(  ( \mathbf{A} \mathbf{B}  - \mathbf{X} )^{T} \big) \ ,
\end{align*}
which implies that
\begin{align*}
\nabla \varphi^{*}_{\mathbf{a}} ( \mathbf{A}, \mathbf{B} )  & =  {\mathbf{G}(\mathbf{b})}^{T} \big ( \mathbf{G}(\mathbf{b})\mathbf{a} -  \mathbf{z} \big ) \\
                                                                                           & = ( \mathbf{I}_{p} \otimes  \mathbf{B} )   \emph{diag}\big( \emph{vec}( \mathbf{W}^{T}) \big)  \emph{vec} \big(  ( \mathbf{A} \mathbf{B}  - \mathbf{X} )^{T}  \big) \\
                                                                                           & = ( \mathbf{I}_{p} \otimes  \mathbf{B} )    \emph{vec} \big(  \mathbf{W}^{T}  \odot ( \mathbf{A} \mathbf{B}  - \mathbf{X} )^{T} \big) \\
                                                                                           & = \emph{vec} \Big(   \mathbf{B} \big(  \mathbf{W}  \odot ( \mathbf{A} \mathbf{B}  - \mathbf{X} ) \big)^{T}  \Big)  \\
                                                                                           & = \emph{vec} \Big(   \big(  \mathbf{W}  \odot ( \mathbf{A} \mathbf{B}  - \mathbf{X} ) )  \mathbf{B}^{T}   \big)^{T}  \Big) \ ,
\end{align*}
and, using again equation~\eqref{eq:D_partial_grad_varphi*}, we conclude that
\begin{equation*}
\nabla \varphi^{*}_{\mathbf{a}} ( \mathbf{A}, \mathbf{B} )  =  \emph{vec}  \big(  \nabla \varphi^{*}_{\mathbf{A}} ( \mathbf{A}, \mathbf{B} )^{T}  \big) \ .
\end{equation*}
Next, we immediately get that the (partial) Hessian matrices of the vectorized form of $\varphi^{*}(.)$ are equal to
\begin{align*}
\nabla^{2} \varphi^{*}_{\mathbf{a}} ( \mathbf{A}, \mathbf{B} )  & = \mathbf{G}(\mathbf{b})^{T} \mathbf{G}(\mathbf{b}) \ , \\
\nabla^{2} \varphi^{*}_{\mathbf{b}} ( \mathbf{A}, \mathbf{B} )  & = \mathbf{F}(\mathbf{a})^{T} \mathbf{F}(\mathbf{a}) \ ,
\end{align*}
since  the specific forms of  $\frac{ \partial \varphi^{*}( \mathbf{A}, \mathbf{B} )}  {\partial \mathbf{a}}$ and $\frac{ \partial \varphi^{*}( \mathbf{A}, \mathbf{B} )}  {\partial \mathbf{b}}$ derived above can be both interpreted as the sum of a linear mapping and a constant term, when they are considered as a function of  $\mathbf{a}$  and $\mathbf{b}$, respectively.

To derive an explicit formula for $\nabla^{2} \varphi^{*}_{\mathbf{b} , \mathbf{a}} ( \mathbf{A}, \mathbf{B} )$, we start from the equation
\begin{equation*}
\frac{ \partial   \varphi^{*}( \mathbf{A}, \mathbf{B} ) }  {\partial\mathbf{b}} = {\mathbf{F}(\mathbf{a})}^{T}  \big( \mathbf{F}(\mathbf{a})\mathbf{b} -  \mathbf{x}  \big) =  \big (   \mathbf{I}_{n}  \otimes  \mathbf{A}^{T}  \big )  \emph{vec} \big(   \mathbf{W}  \odot ( \mathbf{A} \mathbf{B}  - \mathbf{X} )  \big) \ ,
\end{equation*}
and apply the differential rule for a matrix product~\cite{C2017} to get
\begin{equation*}
\nabla^{2} \varphi^{*}_{\mathbf{b} , \mathbf{a}} ( \mathbf{A}, \mathbf{B} )  \mathbf{c} = (   \mathbf{I}_{n}  \otimes  \mathbf{C}^{T} )  \emph{vec} \big(   \mathbf{W}  \odot ( \mathbf{A} \mathbf{B}  - \mathbf{X} )  \big) + (   \mathbf{I}_{n}  \otimes  \mathbf{A}^{T}  )  \emph{vec} (   \mathbf{W}  \odot \mathbf{C} \mathbf{B}  ) \ ,
\end{equation*}
$\forall  \mathbf{C} \in \mathbb{R}^{p \times k}$  with $\mathbf{c} = \emph{vec} (  \mathbf{C}^{T} )$. On one hand, using equation~\eqref{eq:vec_kronprod}, we have
\begin{align*}
(   \mathbf{I}_{n}  \otimes  \mathbf{C}^{T} )  \emph{vec} \big(   \mathbf{W}  \odot ( \mathbf{A} \mathbf{B}  - \mathbf{X} )  \big) & = \emph{vec}  \Big(   \mathbf{C}^{T}   \big( \mathbf{W}  \odot ( \mathbf{A} \mathbf{B}  - \mathbf{X} )  \big)  \Big)  \\
     & = \Big(   \big(   \mathbf{W}  \odot ( \mathbf{A} \mathbf{B}  - \mathbf{X} )  \big)^{T}  \otimes  \mathbf{I}_{k} \Big) \emph{vec} ( \mathbf{C}^{T} )  \\
     & = \Big(   \big(   \mathbf{W}  \odot ( \mathbf{A} \mathbf{B}  - \mathbf{X} )  \big)^{T}  \otimes  \mathbf{I}_{k} \Big) \mathbf{c}  \ ,
\end{align*}
and, on the other hand, using equations~\eqref{eq:vec_kronprod}, \eqref{eq:commat}, \eqref{eq:commat2} and Lemma~\ref{theo2.2:box}, we get
\begin{align*}
(   \mathbf{I}_{n}  \otimes  \mathbf{A}^{T}  )  \emph{vec} (   \mathbf{W}  \odot \mathbf{C} \mathbf{B} ) & = (   \mathbf{I}_{n}  \otimes  \mathbf{A}^{T} ) \emph{diag} \big(  \emph{vec} ( \mathbf{W} ) \big) \emph{vec} (  \mathbf{C} \mathbf{B} )  \\
& =  \Big(  \emph{diag} \big(  \emph{vec} (  \sqrt{\mathbf{W}} ) \big) \mathbf{F}(\mathbf{a})  \Big)^{T} (    \mathbf{B}^{T} \otimes \mathbf{I}_{p} ) \emph{vec} (  \mathbf{C} )  \\
& =  \Big(  \emph{diag} \big(  \emph{vec} (  \sqrt{\mathbf{W}} ) \big) \mathbf{F}(\mathbf{a})  \Big)^{T} (    \mathbf{B}^{T} \otimes \mathbf{I}_{p} )  \mathbf{K}_{(k,p)} \mathbf{K}_{(p,k)} \emph{vec} (  \mathbf{C} )  \\
& =  \mathbf{F}(\mathbf{a})^{T}  \emph{diag} \big(  \emph{vec} (  \sqrt{\mathbf{W}} ) \big)  (    \mathbf{B}^{T} \otimes \mathbf{I}_{p} )  \mathbf{K}_{(k,p)}  \emph{vec} (  \mathbf{C}^{T} )  \\
& =   \mathbf{F}(\mathbf{a})^{T}  \emph{diag} \big(  \emph{vec} (  \sqrt{\mathbf{W}} ) \big)  \mathbf{K}_{(n,p)}  (  \mathbf{I}_{p} \otimes \mathbf{B}^{T}  ) \emph{vec} (  \mathbf{C}^{T} )  \\
& =   \mathbf{F}(\mathbf{a})^{T} \mathbf{K}_{(n,p)}  \emph{diag} \big(  \emph{vec} (  \sqrt{\mathbf{W}}^{T} ) \big)   (  \mathbf{I}_{p} \otimes \mathbf{B}^{T}  ) \emph{vec} (  \mathbf{C}^{T} )  \\
& =   \mathbf{F}(\mathbf{a})^{T} \mathbf{K}_{(n,p)}  \mathbf{G}(\mathbf{b}) \mathbf{c} \ .
\end{align*}
Together, these equalities imply, finally, that
\begin{equation*}
\nabla^{2} \varphi^{*}_{\mathbf{b} , \mathbf{a}} ( \mathbf{A}, \mathbf{B} )  = \Big(   \big(   \mathbf{W}  \odot ( \mathbf{A} \mathbf{B}  - \mathbf{X} )  \big)^{T}  \otimes  \mathbf{I}_{k} \Big) +  \mathbf{F}(\mathbf{a})^{T} \mathbf{K}_{(n,p)}  \mathbf{G}(\mathbf{b})  \ ,
\end{equation*}
as claimed in the theorem.

Finally, the equality $ \nabla^{2} \varphi^{*}_{\mathbf{a} , \mathbf{b}} ( \mathbf{A}, \mathbf{B} )  = \big \lbrack \nabla^{2} \varphi^{*}_{\mathbf{b} , \mathbf{a}} ( \mathbf{A}, \mathbf{B} ) \big \rbrack^{T}$ is a direct consequence of the fact that the Hessian $\nabla^{2} \varphi^{*} ( \mathbf{A}, \mathbf{B} )$ is a  $(p.k + k.n)  \times (p.k + k.n)$ symmetric matrix according to the Schwarz's theorem~\cite{C2017}, see Subsection~\ref{calculus:box} and Remark~\ref{remark4.3:box} below for details.

It remains to establish the equalities between the quadratic forms $\big ( \nabla^{2} \varphi^{*} ( \mathbf{A}, \mathbf{B} ) \big )$ and $\big ( \nabla^{2} \varphi( \mathbf{A} \mathbf{B} ) \big )$. First, note that
\begin{align*}
\nabla^{2} \varphi^{*}_{\mathbf{a}} ( \mathbf{A}, \mathbf{B} ) \mathbf{c} & =   \emph{vec}  \Big (  \big(  ( \mathbf{W}  \odot   \mathbf{C} \mathbf{B} )  \mathbf{B}^{T}  \big)^{T} \Big ) \ ,  \\
\nabla^{2} \varphi^{*}_{\mathbf{b}} ( \mathbf{A}, \mathbf{B} ) \mathbf{d} & =   \emph{vec} \big(  \mathbf{A}^{T} (  \mathbf{W}   \odot  \mathbf{A} \mathbf{D} ) \big) \ ,
\end{align*}
$\forall  \mathbf{C} \in \mathbb{R}^{p \times k}$  with $\mathbf{c} = \emph{vec} (  \mathbf{C}^{T} )$ and $\forall  \mathbf{D} \in \mathbb{R}^{k \times n}$  with $\mathbf{d} = \emph{vec} (  \mathbf{D} )$.

Using these equalities, we deduce
\begin{align*}
\mathbf{c}^{T} \nabla^{2} \varphi^{*}_{\mathbf{a}} ( \mathbf{A}, \mathbf{B} ) \mathbf{c} & =   \big \langle   \nabla^{2} \varphi^{*}_{\mathbf{a}} ( \mathbf{A}, \mathbf{B} ) \mathbf{c} ,  \mathbf{c}  \big \rangle_{2} \\
& = \big \langle  \emph{vec}  \Big (  \big(  ( \mathbf{W}  \odot   \mathbf{C} \mathbf{B} )  \mathbf{B}^{T}  \big)^{T} \Big )  ,  \emph{vec} (  \mathbf{C}^{T} )  \big \rangle_{2} \\
& = \big \langle  \emph{vec}  \big(  ( \mathbf{W}  \odot   \mathbf{C} \mathbf{B} )  \mathbf{B}^{T}  \big)  ,  \emph{vec} (  \mathbf{C} )  \big \rangle_{2} \\
& = \big \langle   ( \mathbf{W}  \odot   \mathbf{C} \mathbf{B} )  \mathbf{B}^{T}  ,   \mathbf{C} \big \rangle_{F} \\
& = \big \langle   \lbrack \nabla^{2} \varphi^{*}_{\mathbf{A}} ( \mathbf{A} , \mathbf{B} ) \rbrack ( \mathbf{C} ),   \mathbf{C} \big \rangle_{F} \\
& =  \big ( \nabla^{2} \varphi^{*}_{\mathbf{A}} ( \mathbf{A} , \mathbf{B} ) \big ) ( \mathbf{C} , \mathbf{C} ) \\
& =  \big ( \nabla^{2} \varphi ( \mathbf{A} \mathbf{B} )  \big ) ( \mathbf{C}  \mathbf{B} , \mathbf{C}  \mathbf{B}  ) \ ,
\end{align*}
and also
\begin{align*}
\mathbf{d}^{T} \nabla^{2} \varphi^{*}_{\mathbf{b}} ( \mathbf{A}, \mathbf{B} ) \mathbf{d} & =   \big \langle   \nabla^{2} \varphi^{*}_{\mathbf{b}} ( \mathbf{A}, \mathbf{B} ) \mathbf{d} ,  \mathbf{d}  \big \rangle_{2} \\
& = \big \langle   \emph{vec} \big(  \mathbf{A}^{T} (  \mathbf{W}   \odot  \mathbf{A} \mathbf{D} ) \big) ,  \emph{vec} (  \mathbf{D} )  \big \rangle_{2} \\
& = \big \langle   \mathbf{A}^{T} (  \mathbf{W}   \odot  \mathbf{A} \mathbf{D} )  ,   \mathbf{D}   \big \rangle_{F} \\
& = \big \langle   \lbrack \nabla^{2} \varphi^{*}_{\mathbf{B}} ( \mathbf{A} , \mathbf{B} ) \rbrack ( \mathbf{D} ),   \mathbf{D} \big \rangle_{F} \\
& =  \big ( \nabla^{2} \varphi^{*}_{\mathbf{B}} ( \mathbf{A} , \mathbf{B} ) \big ) ( \mathbf{D} , \mathbf{D} ) \\
& = \big ( \nabla^{2} \varphi ( \mathbf{A} \mathbf{B} )  \big ) ( \mathbf{A}  \mathbf{D} , \mathbf{A}  \mathbf{D}  )  \ ,
\end{align*}
where, in both cases, the last equality results from equation~\eqref{eq:D_partial_hess_varphi*}. The last equality in the theorem, 
\begin{equation*}
\mathbf{d}^{T} \nabla^{2} \varphi^{*}_{\mathbf{b} , \mathbf{a}} ( \mathbf{A}, \mathbf{B} )  \mathbf{c} = \big \langle \nabla \varphi ( \mathbf{A} \mathbf{B} ) , \mathbf{C}  \mathbf{D} \big \rangle_{F} + \big ( \nabla^{2} \varphi ( \mathbf{A} \mathbf{B} )  \big ) ( \mathbf{A} \mathbf{D}  , \mathbf{C}  \mathbf{B} )  \ ,
\end{equation*}
can be derived in a similar way, by a lengthy, but direct, computation.
\\
\end{proof}
\begin{corol4.2} \label{corol4.2:box}
 For $\mathbf{X}\in\mathbb{R}^{p \times n}$,  $\sqrt{\mathbf{W}}\in\mathbb{R}^{p \times n}_+$, $\lambda \in \mathbb{R}_{+*}$ and any fixed integer $ k \le \emph{rank}( \mathbf{X} ) \le \text{min}( {p},{n} )$, the  objective function $\varphi^{*}_{\lambda}(.)$ used in the~\eqref{eq:MMMF} formulation of the WLRA problem
\begin{equation*}
\varphi^{*}_{\lambda} : \mathbb{R}^{p \times k} \times  \mathbb{R}^{k \times n} \longrightarrow \mathbb{R} : ( \mathbf{A}, \mathbf{B} ) \mapsto \varphi^{*}_{\lambda}( \mathbf{A}, \mathbf{B} ) =  \frac{1}{2}   \Vert \sqrt{\mathbf{W}}  \odot ( \mathbf{X} - \mathbf{A}\mathbf{B} )  \Vert^{2}_{F} + \frac{\lambda}{2} ( \Vert \mathbf{A} \Vert^{2}_{F} + \Vert \mathbf{B} \Vert^{2}_{F} )
\end{equation*}
is $C^{\infty}$ differentiable over its domain of definition and the partial first-order derivatives of $\varphi^{*}_{\lambda}(.)$ with respect to $\mathbf{a} = \emph{vec}( \mathbf{A}^{T})$ and $\mathbf{b} = \emph{vec}( \mathbf{B} )$ are equal, respectively, to
\begin{equation*}
\frac{ \partial   \varphi^{*}_{\lambda}( \mathbf{A}, \mathbf{B} ) }  {\partial\mathbf{a}} =   \nabla ( \varphi^{*}_{\lambda} )_{\mathbf{a}} ( \mathbf{A}, \mathbf{B} )  = \mathbf{G}(\mathbf{b})^{T} \mathbf{G}(\mathbf{b})\mathbf{a} - \mathbf{G}(\mathbf{b})^{T} \mathbf{z} + \lambda \mathbf{a}
\end{equation*}
and
\begin{equation*}
\frac{ \partial   \varphi^{*}_{\lambda}( \mathbf{A}, \mathbf{B} ) }  {\partial\mathbf{b}} =  \nabla (\varphi^{*}_{\lambda})_{\mathbf{b}} ( \mathbf{A}, \mathbf{B} )  = \mathbf{F}(\mathbf{a})^{T} \mathbf{F}(\mathbf{a})\mathbf{b} - {\mathbf{F}(\mathbf{a})} ^{T} \mathbf{x} + \lambda \mathbf{b}  \ .
\end{equation*}
Furthermore,  the partial second-order derivatives of $\varphi^{*}_{\lambda}(.)$ with respect to $\mathbf{a} = \emph{vec}( \mathbf{A}^{T})$ and $\mathbf{b} = \emph{vec}( \mathbf{B} )$ are given by
\begin{equation*}
 \frac{ \partial^2 \varphi^{*}_{\lambda}( \mathbf{A}, \mathbf{B} )}  {\partial^2 \mathbf{a}} =  \nabla^{2} ( \varphi^{*}_{\lambda} )_{\mathbf{a}} ( \mathbf{A}, \mathbf{B} )  = \mathbf{G}(\mathbf{b} )^{T} \mathbf{G}(\mathbf{b} ) + \lambda . \mathbf{I}_{p.k}
\end{equation*}
and
\begin{equation*}
\frac{ \partial^2 \varphi^{*}_{\lambda}( \mathbf{A}, \mathbf{B} )}  {\partial^2 \mathbf{b}} = \nabla^{2} (\varphi^{*}_{\lambda} )_{\mathbf{b}} ( \mathbf{A}, \mathbf{B} )  = \mathbf{F}(\mathbf{a})^{T} \mathbf{F}(\mathbf{a}) + \lambda . \mathbf{I}_{k.n}  \ .
\end{equation*}
\end{corol4.2}
\begin{proof}
$\varphi^{*}_{\lambda}(.)$ is the sum of  three $C^{\infty}$ differentiable functions, e.g., $\varphi^{*}(.)$ and the mappings $\frac{\lambda}{2} \Vert \mathbf{A} \Vert^{2}_{F} = \frac{\lambda}{2} \Vert \mathbf{a} \Vert^{2}_{2}$ and  $\frac{\lambda}{2} \Vert \mathbf{B} \Vert^{2}_{F} = \frac{\lambda}{2} \Vert \mathbf{b} \Vert^{2}_{2}$ and is thus $C^{\infty}$ differentiable.
The formulas for $\frac{ \partial   \varphi^{*}_{\lambda}( \mathbf{A}, \mathbf{B} ) }  {\partial\mathbf{a}}$ and $\frac{ \partial   \varphi^{*}_{\lambda}( \mathbf{A}, \mathbf{B} ) }  {\partial\mathbf{b}}$ follow immediately from Theorem~\ref{theo4.3:box}, standard differentiation rules and the differential rule for a mapping of the form $ \mathbb{R}^{m} \longrightarrow  \mathbb{R} : \mathbf{x} \mapsto  \frac{1}{2} \Vert \mathbf{x} \Vert^{2}_{2} $.

The form of the partial second-order derivatives of $\varphi^{*}_{\lambda}(.)$ given in the theorem is a direct consequence of Theorem~\ref{theo4.3:box} and the fact that $\nabla ( \varphi^{*}_{\lambda} )_{\mathbf{a}} ( \mathbf{A}, \mathbf{B} )$ and $ \nabla (\varphi^{*}_{\lambda})_{\mathbf{b}} ( \mathbf{A}, \mathbf{B} )$ are both the sum of two linear mappings and of a constant term when they are considered as a function of $\mathbf{a}$ and $\mathbf{b}$, respectively.
\\
\end{proof}
\begin{remark4.2} \label{remark4.2:box}
The equations
\begin{equation*}
\frac{ \partial  \mathbf{e}( \mathbf{A}, \mathbf{B} )}  {\partial\mathbf{b}} = - \mathbf{F}(\mathbf{a})  \quad \text{and} \quad  \frac{ \partial  \mathbf{e}( \mathbf{A}, \mathbf{B} )}  {\partial\mathbf{a}} = - \mathbf{K}_{(n,p)} \mathbf{G}(\mathbf{b})
\end{equation*}
derived in the proof of Theorem~\ref{theo4.3:box} show that residual function $\mathbf{e}(.)$ is not a nonlinear function of its arguments as defined in Subsection~\ref{calculus:box}. However, despite of this, the cost function $\varphi^{*}(.)$ is still nonlinear as the partial derivatives of $\varphi^{*}(.)$ with respect to $\mathbf{a}$ and $\mathbf{b}$ are functions of  $\mathbf{b}$ and $\mathbf{a}$, respectively, as demonstrated in Theorem~\ref{theo4.3:box}. Furthermore, as the minimization of the cost function  $\varphi^{*}(.)$ has no closed form solution in general, we can still consider $\varphi^{*}(.)$ as a NLLS functional as defined in Subsection~\ref{calculus:box}. $\blacksquare$
\\
\end{remark4.2}

Due to the block diagonal structures of both $\mathbf{F}(\mathbf{a})$ and $\mathbf{G}(\mathbf{b})$, the evaluation of the partial derivatives of $\varphi^{*}(.)$ is fast, easy to implement and may be parallelized. Moreover, we already know that
\begin{equation*}
\frac{ \partial\varphi^{*}( \mathbf{A}, \mathbf{B} )}  {\partial\mathbf{a}} = \mathbf{0}^{k.p}  \quad \text{or}  \quad  \frac{ \partial\varphi^{*}( \mathbf{A}, \mathbf{B}  )} {\partial\mathbf{b} } = \mathbf{0}^{k.n}
\end{equation*}
if the ALS algorithm is used to minimize $\varphi^{*}( \mathbf{A}, \mathbf{B} )$ and the iterations are stopped after computing $\mathbf{A}$ or $\mathbf{B}$, respectively. Similar remarks apply to the partial derivatives of $\varphi^{*}_{\lambda}(.)$.

The main payoff of the two-block ALS method is its simplicity since it involves solving mainly two sequences of small (eventually regularized) linear least-squares problems. Moreover, compared to gradient-type algorithms, it has the advantage that there is no need to tune optimization parameters like step sizes~\cite{OUV2023}. However, practical experience with this algorithm shows that, in many cases, the "NIPALS" iterates do not converge to the closest fit (e.g., the infimum or minimum of  $\varphi^{*}( \mathbf{A}, \mathbf{B} )$ or  $\varphi^{*}_{\lambda}( \mathbf{A}, \mathbf{B} )$) and get frequently stuck in sub-optimal local minima for a small value of $k$ or a  poorly chosen starting point~\cite{GZ1979}\cite{RS2005}. This is especially true when some weights are equal to $0$ (i.e., when missing values are present in $\mathbf{X}$), even with the initialization procedure proposed by Gabriel and Zamir~\cite{GZ1979}. Moreover this initialization procedure is only applicable if $k=1$ and if there is one and only one missing cell ($\mathbf{W}_{ij}=0$)  for the matrix entries in the $i^{th}$ row and $j^{th}$ column of $\mathbf{X}$ for all $i$ and $j$ (see Gabriel and Zamir~\cite{GZ1979} for more details). Furthermore, it is known that the two-block ALS algorithm is vulnerable to flatlining~\cite{BF2005} and inherits in many cases of the very slow convergence of the block coordinate descent method~\cite{NW2006}. To overcome these difficulties, we describe in the next section, various first- and second-order separable NLLS algorithms for minimizing $\psi(.)$ instead of $\varphi^{*}(.)$.

\begin{remark4.3} \label{remark4.3:box}
If we concatenate the vectors $\mathbf{a}$ and $\mathbf{b}$ in $\mathbf{c} =  \begin{bmatrix} \mathbf{a}   \\ \mathbf{b} \end{bmatrix} \in \mathbb{R}^{k.(p+n)}$, we may define the following residual and objective functions:
\begin{equation*}
r( \mathbf{c} ) =   \mathbf{x} - \mathbf{F}(\mathbf{a})\mathbf{b} = \mathbf{K}_{(n,p)} \big (  \mathbf{z} - \mathbf{G}(\mathbf{b})\mathbf{a} \big )  = \mathbf{e}( \mathbf{A}, \mathbf{B} ) 
\end{equation*}
and
\begin{equation*}
\phi( \mathbf{c} ) = \frac{1}{2} r( \mathbf{c} )^{T} r( \mathbf{c} ) = \frac{1}{2} \Vert r( \mathbf{c} )  \Vert^{2} _{2} = \varphi^{*}( \mathbf{A}, \mathbf{B} ).
\end{equation*}
According to equation~\eqref{eq:gradvec_nlls}, the gradient of $\phi(.)$ is then equal to
\begin{equation*}
\nabla \phi( \mathbf{c} ) =  \mathit{J}( \mathbf{r}( \mathbf{c} ) )^{T} \mathbf{r}( \mathbf{c} ) \ , 
\end{equation*}
with the Jacobian matrix $\mathit{J}( \mathbf{r}( \mathbf{c} ) ) \in \mathbb{R}^{p.n \times k.(p+n)}$ having the form
\begin{equation*}
\mathit{J}( \mathbf{r}( \mathbf{c} ) ) = \begin{bmatrix}  \frac{ \partial \mathbf{r}( \mathbf{c} ) }  {\partial\mathbf{a}}  & \frac{ \partial \mathbf{r}( \mathbf{c} ) }  {\partial\mathbf{b}} \end{bmatrix} =  \begin{bmatrix}  \frac{ \partial \mathbf{e}( \mathbf{A}, \mathbf{B} ) }  {\partial\mathbf{a}}  & \frac{ \partial\mathbf{e}( \mathbf{A}, \mathbf{B} ) }  {\partial\mathbf{b}} \end{bmatrix} \  .
\end{equation*}
Now, using Remark~\ref{remark4.2:box} and Theorem~\ref{theo4.3:box}, we have:
\begin{equation*}
\mathit{J}( \mathbf{r}( \mathbf{c} ) ) = - \begin{bmatrix} \mathbf{K}_{(n,p)} \mathbf{G}(\mathbf{b}) &  \mathbf{F}(\mathbf{a}) \end{bmatrix}
\end{equation*}
and
\begin{equation*}
\nabla \phi( \mathbf{c} ) =  \begin{bmatrix}   \nabla \varphi^{*}_{\mathbf{a}}( \mathbf{A}, \mathbf{B} ) \\   \nabla \varphi^{*}_{\mathbf{b}}( \mathbf{A}, \mathbf{B} ) \end{bmatrix} = \begin{bmatrix}   \mathbf{G}(\mathbf{b})^{T} \mathbf{G}(\mathbf{b})\mathbf{a} - \mathbf{G}(\mathbf{b})^{T} \mathbf{z}  \\ \mathbf{F}(\mathbf{a})^{T} \mathbf{F}(\mathbf{a})\mathbf{b} - {\mathbf{F}(\mathbf{a})} ^{T} \mathbf{x}  \end{bmatrix}  \ .
\end{equation*}
Finally, if we differentiate again $\nabla \phi( \mathbf{c} )$ with respect to $\mathbf{c}$, an analytic formulae for the Hessian matrix $\nabla^2 \phi( \mathbf{c} )$ can be obtained, which is essentially equivalent to the results given in Theorem~\ref{theo4.3:box}, see~\cite{BF2005}\cite{HF2015b} for a derivation of this Hessian matrix.
Equipped with these exact formulas for $\nabla \phi( \mathbf{c} ), \mathit{J}( \mathbf{r}( \mathbf{c} ) )$ and $\nabla^2 \phi( \mathbf{c} )$, standard first- and second-order NLLS methods such as the steepest gradient, Gauss-Newton, Levenberg-Marquardt and Newton algorithms~\cite{DS1983}\cite{NW2006}\cite{MN2010} can also be used (and have been used) to minimize directly $\phi( \mathbf{c} ) = \varphi^{*}(\mathbf{A}, \mathbf{B} )$ or $\phi_{\lambda}( \mathbf{c} ) = \varphi^{*}_{\lambda}(\mathbf{A}, \mathbf{B} )$, and to solve the WLRA problem and its MMMF variant~\cite{BF2005}\cite{D2011}\cite{HF2015}. However, as it is arguably preferable to keep the dimension of the search space as much low as possible and because the joint optimization strategy of minimizing directly $\phi(.)$ has been found to be much less efficient and less robust than the variable projection framework (based on the minimization of $\psi(.)$) detailed in the next section~\cite{D2011}\cite{OYD2011}\cite{BA2015}\cite{HF2015}\cite{BL2020}, we don't focus here anymore on the direct minimization of $\phi(.)$ or  $\varphi^{*}(.)$ (or  alternatively $\phi_{\lambda}(.)$ or  $\varphi^{*}_{\lambda}(.)$) for solving the WLRA problem. $\blacksquare$
\end{remark4.3}

\section{The variable projection framework} \label{varpro:box}

We now explain how to minimize the cost function $\psi(.)$, which is used in the~\eqref{eq:VP1} formulation of the WLRA problem. In addition to the equivalence of the~\eqref{eq:P1} and~\eqref{eq:VP1} formulations of the WLRA problem stated in Theorem~\ref{theo3.9:box}, the variable projection approach is further justified by a theorem originally proved by Golub and Pereyra in~\cite{GP1973}, which shows, under some differentiability conditions, that if $\widehat{\mathbf{a}} = \emph{vec}(\widehat{\mathbf{A}}^{T})$ is a critical point of $\psi(.)$ and $\widehat{\mathbf{B}}$ is calculated by equation~\eqref{eq:B_nipals}, e.g., by solving $n$ independent linear least-squares problems as described in the beginning of Section~\ref{nipals:box}, then $(\widehat{\mathbf{A}}, \widehat{\mathbf{B}})$ is also a first-order critical point of $\varphi^{*}(.)$. We will give a demonstration of this result later in Theorem~\ref{theo5.7:box} (see Subsection~\ref{hess:box}) for completeness.

General optimization methods used to minimize a functional like $\psi(.)$ are termed variable projection algorithms and are described in~\cite{GP1973}\cite{RW1980}\cite{K1974}\cite{K1975}\cite{B2009}\cite{GP2003}\cite{OR2013}. Their advantages are that they usually solve mixed linear-nonlinear least-squares problems like $\varphi^{*}(.)$ in less time, fewer function evaluations and better global convergence than standard NLLS codes, and that no starting estimate of the linear variable (e.g., $\mathbf{B}$) is required~\cite{N2000}. In the context of the WLRA or matrix completion problems, they offer also other advantages as shown in~\cite{OD2007}\cite{OYD2011}\cite{D2011}\cite{BA2015}\cite{HF2015}\cite{HZF2017} and as we will illustrate in the next sections. However, many of them have also a major drawback as they expand considerably the dimensionality of the WLRA problem (see Subsection~\ref{varpro_wlra:box} for details). This limits severely their use for medium and large datasets, which are currently found now in many applications, beyond variations of the variable projection steepest (e.g., gradient) descent method or similar first-order methods~\cite{SJ2004}\cite{DMK2011}\cite{BA2015}\cite{BL2020}\cite{OUV2023}. In our WLRA context, the simplest variable projection steepest descent method can be written as
\begin{equation*}
\mathbf{a}_{i+1} = \mathbf{a}_{i} -  \alpha_{i} \nabla \psi( \mathbf{a}_{i} ) .
\end{equation*}
In words, with this basic method, we move by making a correction step that is proportional to the negative of the gradient of $\psi(.)$ and the positive scalar $\alpha_{i}$ can be used to control the size of the step without changing its direction~\cite{MN2010}\cite{OUV2023}. This basic method works fine for simple models, but is often too simplistic when there are many parameters to estimate like in our WLRA problem. Furthermore, its convergence can be very slow without cleaver strategies to control $\alpha_{i}$ or the use of second-order information, especially in the final stage~\cite{NW2006}\cite{BL2020}\cite{OUV2023}. Near a local minimizer, the steepest descent method converges at a linear rate  depending on the condition number in a neighborhood of this minimizer. However, this convergence rate deteriorates dramatically when the Hessian of the cost function is ill-conditioned and we will demonstrate later, in Subsection~\ref{hess:box}, that this always the case for the cost function $\psi(.)$. As another illustration, during the iterations, the curvature of $\psi(.)$ is usually not the same in all directions. If there is a long and narrow valley in the values of $\psi(.)$, which is not unusual when the weights are not uniform~\cite{SJ2004}\cite{ZLTW2018}, the component of the gradient in the direction that points along the bottom of the valley can be very small while the component perpendicular to the walls of the valley can be quite large even though we have to move a long distance along the base and a small distance perpendicular to the walls to move in the right direction. This is the so-called "error valley" problem, which can be alleviated only if we use some information about the curvature as well as the gradient of $\psi(.)$ in the design of the method~\cite{NW2006}. However, second-order derivatives of the cost function are very often prohibitively expensive to compute and we need to find a good compromise between accuracy and speed when the dimensions and the number of variables of the problem are large~\cite{MN2010}.

Thus, since the convergence of the steepest descent method or its variants, like conjugate gradient methods, can be very slow and second derivatives are expensive to evaluate, we concentrate our attention on (pseudo) second-order or quasi-Newton methods well adapted to NLLS problems~\cite{DS1983}\cite{NW2006}\cite{MN2010}\cite{HPS2012}. These methods aim to avoid the drawbacks of Newton’s methods while maintaining the benefits of using second-order
information and introduce also some suitable regularization to cup with the singularity of the Hessian. After a brief description of the Newton, Gauss-Newton, augmented Gauss-Newton and Levenberg-Marquardt algorithms in Subsection~\ref{opt:box}, we give in the next sections a detailed study of the Jacobian matrix $\mathit{J}( \mathbf{r}(\mathbf{a}) )$, gradient vector $\nabla \psi( \mathbf{a} )$ and Hessian matrix $\nabla^2 \psi( \mathbf{a} )$, which are pivotal in these variable projection quasi-Newton algorithms and whose specific properties in the context of the WLRA problem have not always been well appreciated in past studies, except in~\cite{R1974}\cite{OD2007}\cite{OYD2011}.

\subsection{Second-order NLLS optimization methods} \label{opt:box}

As discussed in Subsection~\ref{varpro_wlra:box}, the minimization of  $\psi(.)$ is equivalent to the standard NLLS problem
\begin{equation*}
\min_{\mathbf{a}\in\mathbb{R}^{p.k}}   \,   \psi(\mathbf{a} ) = \frac{1}{2} \Vert \mathbf{P}^{\bot}_{\mathbf{F}(  \mathbf{a} ) }\mathbf{x} \Vert^{2}_{2} = \frac{1}{2} \Vert  \mathbf{r}(  \mathbf{a} ) \Vert^{2}_{2} = \frac{1}{2} \mathbf{r}(  \mathbf{a} )^{T} \mathbf{r}(  \mathbf{a} ) \ ,
\end{equation*}
where $\mathbf{r}(  \mathbf{a} ) = \mathbf{P}^{\bot}_{\mathbf{F}(  \mathbf{a} ) }\mathbf{x}$. Numerous first- and second-order iterative methods are available for minimizing a sum of squares of nonlinear functions such as $\psi(.)$~\cite{DS1983}\cite{NW2006}\cite{MN2010}\cite{HPS2012}. However, for finding a solution of our~\eqref{eq:VP1} problem with these methods, we first note that a certain degree of smoothness of the objective function $\psi(.)$ is required, meaning that $\psi(.)$ must possess one or better two continuous derivatives and the results of Subsection~\ref{varpro_wlra:box} show that these smoothness conditions are not systematically verified if $\mathbf{W}$ has some zero elements  as the orthogonal projector $\mathbf{P}^{\bot}_{\mathbf{F}(  \mathbf{a} )}$ can be a discontinuous function of $\mathbf{a}$ even if $\mathbf{A} = \emph{mat}_{k \times p}(  \mathbf{a} )^{T} = h( \mathbf{a} )$ is of full column rank (see Theorem~\ref{theo3.14:box} and Corollary~\ref{corol3.4:box}). The degree of smoothness of  $\psi(.)$ and $\mathbf{P}^{\bot}_{\mathbf{F}(.)}$ will be further studied in Subsection~\ref{jacob:box}, but we note that, despite these caveats, some standard iterative NLLS algorithms have been used very successfully to solve the~\eqref{eq:VP1} problem even without proper regularization of $\psi(.)$ to insure its smoothness when missing values are present~\cite{C2008b}\cite{OYD2011}\cite{GM2011}\cite{HF2015}\cite{HZF2017}.

The recommended standard methods are the Gauss-Newton, Levenberg-Marquardt, trust-region Gauss-Newton and augmented Gauss-Newton algorithms if second-order derivatives are difficult or cumbersome  to evaluate~\cite{DS1983}\cite{NW2006}\cite{MN2010}\cite{HPS2012}. All these methods attempt to minimize $\psi(.)$ by finding a zero of $\nabla \psi(.)$, i.e., a point $\mathbf{a} = \emph{vec}(\mathbf{A}^{T})$ such that
\begin{equation*}
\nabla \psi( \mathbf{a} ) = \mathit{J} \big ( \mathbf{r}(\mathbf{a}) \big )^{T} \mathbf{r}(\mathbf{a}) = \mathbf{0}^{k.p} \ .
\end{equation*}
Moreover, all four methods may be interpreted as variations of Newton's method to find a zero of $\nabla \psi(.)$~\cite{DS1983}\cite{NW2006}\cite{MN2010}. In Newton's method, the correction vector $d \mathbf{a}_{n}$ for improving an approximate initial solution vector $\mathbf{a}$ of the equation $\nabla \psi( \mathbf{a} ) = \mathbf{0}^{k.p}$ is found as the solution to the linear system
\begin{equation} \label{eq:newton}
\nabla^2 \psi( \mathbf{a} ) d\mathbf{a}_{n} = -\mathit{J}  \big( \mathbf{r}(\mathbf{a})  \big)^{T} \mathbf{r}(\mathbf{a}) \ ,
\end{equation}
where $\nabla^2 \psi( \mathbf{a} )$ is the Hessian of $\psi(.)$ at $\mathbf{a}$ given by
\begin{equation} \label{eq:hessian}
\nabla^2 \psi(\mathbf{a} ) = \mathit{J}  \big ( \mathbf{r}(\mathbf{a})  \big )^{T} \mathit{J}  \big ( \mathbf{r}(\mathbf{a})  \big) + \sum_{l=1}^{n.p} \mathbf{r}_{l} (\mathbf{a}) \nabla^2 \mathbf{r}_{l} (\mathbf{a}) \ .
\end{equation}
In this last equation, $\nabla^2 \mathbf{r}_l (\mathbf{a})$ is the Hessian matrix of the $l^{th}$ component of the residual functional $\mathbf{r}(\mathbf{a})$ (i.e., $\mathbf{r}_l (\mathbf{a})$), which is a $p.k \times p.k$ symmetric matrix. The Newton method is based on the second-order Taylor expansion of $\psi(.)$ in a neighborhood of the current iterate $\mathbf{a}$ (see equation~\eqref{eq:taylor_expan2} in Subsection~\ref{calculus:box}), e.g.,
\begin{equation*}
\psi( \mathbf{a} + d\mathbf{a} ) \approx  N( d\mathbf{a} ) = \psi( \mathbf{a} ) +  d\mathbf{a}^{T} \nabla \psi( \mathbf{a}) +  \frac{1}{2} d\mathbf{a}^{T} \nabla^2 \psi( \mathbf{a} ) d\mathbf{a}  \ .
\end{equation*}
More precisely, the Newton method attempts to minimize $\psi(.)$ at each iteration by finding a first-order stationary point $d\mathbf{a}_{n}$ of this quadratic model $N(.)$. Setting the gradient of $N(.)$ to zero, e.g., $\nabla N(d\mathbf{a}_{n}) = \mathbf{0}^{k.p}$, we obtain the following equation
\begin{equation*}
\nabla \psi( \mathbf{a}) +  \nabla^2 \psi( \mathbf{a} ) d\mathbf{a}_{n} = \mathbf{0}^{k.p}  \ ,
\end{equation*}
from which we derived immediately equation~\eqref{eq:newton} defining the Newton iteration. Moreover, if the Hessian matrix $\nabla^2 \psi( \mathbf{a} )$, which is also equal to $\nabla^2 N( d\mathbf{a}_{n} )$, is positive definite then $d\mathbf{a}_{n}$ is a strict global minimizer of $N(.)$ and in a descent direction for $\psi(.)$. In other words, the Newton iteration is well defined as soon as $\nabla^2 \psi( \mathbf{a} )$ is positive definite, but runs into troubles when it is not, for example in regions of mixed curvature of $\psi(.)$. It may even happen during the iterations that $\nabla^2 \psi( \mathbf{a} )$ becomes definite negative in which case $d\mathbf{a}_{n}$ will be a strict global maximizer of $N(.)$ instead of a minimizer, which is a major drawback of the basic Newton method and explains why it lacks global convergence~\cite{MN2010}\cite{HPS2012}. Moreover, since Newton's method requires the computation of second-order derivatives, which can be cumbersome for large-scale problems (see equation~\eqref{eq:hessian}), it is rarely used in practice despite its quadratic convergence in a neighborhood of a first-order critical point of $\psi(.)$~\cite{DS1983}\cite{NW2006}\cite{MN2010}.

Importantly, the smallest eigenvalue of the positive (semi-definite) matrix $\mathit{J}( \mathbf{r}(\mathbf{a}) )^{T} \mathit{J}( \mathbf{r}(\mathbf{a}) )$ can be used to assess the relative importance of the two terms in $\nabla^2 \psi(\mathbf{a} )$~\cite{HPS2012}. More precisely, if for all $\mathbf{a}$ in a neighborhood of a minimizer of $\psi(.)$, the quantities $\vert \mathbf{r}_l (\mathbf{a})\vert \Vert \nabla^2 \mathbf{r}_l (\mathbf{a})\Vert_{2}$ for $l=1, \cdots, n.p$ are small relative to this eigenvalue, the term $\mathit{J}( \mathbf{r}(\mathbf{a}) )^{T} \mathit{J}( \mathbf{r}(\mathbf{a}) )$ will dominate the Hessian matrix~\cite{HPS2012}. Now, depending on the relative importance of these two terms in $\nabla^2 \psi(\mathbf{a} )$, the recommended methods are the Gauss-Newton, Levenberg-Marquardt, trust-region Gauss-Newton and augmented Gauss-Newton algorithms, which  involve different approximations of the second term in the Hessian of $\psi(.)$.

The Gauss-Newton method approximates $\nabla^2 \psi(\mathbf{a} ) $ with $\mathit{J}( \mathbf{r}(\mathbf{a}) )^{T} \mathit{J}( \mathbf{r}(\mathbf{a}) )$, i.e., drops the second term of the Hessian of $\psi(.)$, which contains products of the $\mathbf{r}_l (\mathbf{a})$ functions and their second-order derivatives. This approximation is exact only if the residual function $\mathbf{r}(.)$ is linear in $\mathbf{a}$, which is usually valid only in a neighborhood of a minimum of $\psi(.)$. The Gauss-Newton method is intended for problems in which the second term of the Hessian matrix is small relative to the first term. Thus, this Gauss-Newton approximation is based on the assumptions that the functions $\mathbf{r}_l (\mathbf{a})$ have small curvatures or that near the solution the magnitudes of the $\mathbf{r}_l (\mathbf{a})$ functions are small. If these conditions are satisfied the Gauss-Newton method will ultimately converge at the same rate as Newton's method despite full second-order derivatives are not used. In Gauss-Newton's method, the correction vector $d \mathbf{a}_{gn}$ for improving an approximate solution is then found as the solution to the linear system of equations
\begin{equation} \label{eq:gn_normal}
\mathit{J} \big( \mathbf{r}(\mathbf{a}) \big)^{T} \mathit{J} \big( \mathbf{r}(\mathbf{a}) \big) d\mathbf{a}_{gn} = -\mathit{J} \big( \mathbf{r}(\mathbf{a}) \big)^{T} \mathbf{r}(\mathbf{a})  \ .
\end{equation}
The Gauss-Newton method can be also introduced by a linearization argument. If, given $\mathbf{a}\in\mathbb{R}^{p.k}$, we could solve the problem
\begin{equation*}
\min_{d\mathbf{a}\in\mathbb{R}^{p.k}}   \,   \psi(\mathbf{a} +  d\mathbf{a} ) =  \frac{1}{2} \Vert  \mathbf{r}(  \mathbf{a} +  d\mathbf{a} ) \Vert^{2}_{2}  = \frac{1}{2} \mathbf{r}(  \mathbf{a} +  d\mathbf{a} )^{T} \mathbf{r}(  \mathbf{a} +  d\mathbf{a} )  \ ,
\end{equation*}
then $\mathbf{a} +  d\mathbf{a}$ is a minimizer of $\psi(.)$. Since $\mathbf{r}(.)$ is a nonlinear residual function, we must seek an approximate solution that can be improved iteratively. A natural way to find an approximate solution is to linearize the residual function $\mathbf{r}(.)$ around $\mathbf{a}$. If we assume that $\mathbf{r}(.)$  is twice continuously differentiable at $\mathbf{a}\in\mathbb{R}^{p.k}$, we have the first-order Taylor expansion 
\begin{equation*}
\mathbf{r}(  \mathbf{a} +  d\mathbf{a} ) = \mathbf{r}(  \mathbf{a} ) + \mathit{J}( \mathbf{r}(\mathbf{a}) )d\mathbf{a} +  \mathcal{O}( \Vert d\mathbf{a} \Vert^{2}_{2} )
\end{equation*}
and if we substitute this Taylor approximation for $\mathbf{r}(  \mathbf{a} +  d\mathbf{a} )$ in the definition of $\psi(\mathbf{a} +  d\mathbf{a} )$, this leads to the quadratic function approximation
\begin{equation*}
\psi( \mathbf{a} +  d\mathbf{a} ) \approx G( d\mathbf{a}  ) = \psi(\mathbf{a}) + d\mathbf{a}^{T} \mathit{J}  \big( \mathbf{r}(\mathbf{a})  \big)^{T} \mathbf{r}(\mathbf{a}) +  \frac{1}{2} d\mathbf{a}^{T}  \mathit{J}  \big( \mathbf{r}(\mathbf{a})  \big)^{T} \mathit{J}  \big( \mathbf{r}(\mathbf{a})  \big)  d\mathbf{a}  \  ,
\end{equation*}
which must be minimized at each iteration. As a model for the change of the cost function $\psi(.)$, the quadratic function $G(.)$ has two important advantages compared to the Newton quadratic model $N(.)$, first, it involves only first derivatives of the residual function $\mathbf{r}(.)$ and, second, the symmetric matrix $\mathit{J}( \mathbf{r}(\mathbf{a}) )^{T} \mathit{J}( \mathbf{r}(\mathbf{a}) )$ is always positive semi-definite and is positive definite  if $ \mathit{J}( \mathbf{r}(\mathbf{a}) )$ is of full column rank. Of course, the drawback is a lost of accuracy as full second-order information from the Hessian matrix is not taken into account. The gradient of this quadratic function is equal to
\begin{equation*}
\nabla G( d\mathbf{a} ) = \mathit{J} \big( \mathbf{r}(\mathbf{a}) \big)^{T} \mathbf{r}(\mathbf{a}) +   \mathit{J} \big( \mathbf{r}(\mathbf{a}) \big)^{T} \mathit{J} \big( \mathbf{r}(\mathbf{a}) \big) d\mathbf{a}  \  ,
\end{equation*}
and setting it to zero leads to the linear system~\eqref{eq:gn_normal}, which is also the normal equations of the linear least-squares problem
\begin{equation} \label{eq:gn_llsq}
d\mathbf{a}_{gn} = \text{Arg}\min_{d\mathbf{a} \in \mathbb{R}^{p.k}}   \,   \frac{1}{2} \Vert \mathbf{r}(  \mathbf{a} ) + \mathit{J} \big( \mathbf{r}(\mathbf{a}) \big) d\mathbf{a} \Vert^{2}_{2}  \  ,
\end{equation}
whose unique solution is
\begin{equation*}
d\mathbf{a}_{gn} = - \left( \mathit{J} \big( \mathbf{r}(\mathbf{a}) \big)^{T} \mathit{J} \big( \mathbf{r}(\mathbf{a}) \big) \right)^{-1} \mathit{J} \big( \mathbf{r}(\mathbf{a}) \big)^{T} \mathbf{r}(\mathbf{a})  \  ,
\end{equation*}
if $\mathit{J}( \mathbf{r}(\mathbf{a}) )$ has full column rank or, if this Jacobian matrix is rank-deficient or ill-conditioned, whose unique minimum 2-norm solution is
\begin{equation} \label{eq:gn_direction}
d\mathbf{a}_{gn} = - \mathit{J} \big( \mathbf{r}(\mathbf{a}) \big)^{+}  \mathbf{r}(\mathbf{a})  \  ,
\end{equation}
where $\mathit{J}( \mathbf{r}(\mathbf{a}) )^{+}$ is the pseudo-inverse of the Jacobian matrix of the residual function $\mathbf{r}(.)$ at $\mathbf{a}$. In the rest of this monograph, we will mostly use the pseudo-inverse notation $\mathit{J}( \mathbf{r}(\mathbf{a}) )^{+}$  to indicate that the normal equations shall not be used to compute $d\mathbf{a}_{gn}$ if $\mathit{J}( \mathbf{r}(\mathbf{a}) )$ is ill-conditioned or singular.

The linear least-squares problem~\eqref{eq:gn_llsq}  can be solved by stable orthogonalization methods or the SVD decomposition of the Jacobian matrix $\mathit{J}( \mathbf{r}(\mathbf{a}) )$, see Subsection~\ref{lin_alg:box} and~\cite{GVL1996}\cite{HPS2012} for details. Thus, the last equation becomes the iteration formula
\begin{equation} \label{eq:gn_iter}
\mathbf{a}_{i+1} = \mathbf{a}_{i} - \mathit{J} \big( \mathbf{r}(\mathbf{a}_i) \big)^{+}  \mathbf{r}(\mathbf{a}_i)  \ ,
\end{equation}
which is known as the Gauss-Newton algorithm. Given an initial estimate $\mathbf{a}_{0}$, the linear least-squares problem~\eqref{eq:gn_llsq} associated with the Taylorized equations are solved to yield a correction to this vector $\mathbf{a}_{0}$. This process is repeated and stops if and when the vectors $\mathbf{a}_{i}$ (or the values $\psi( \mathbf{a}_{i} )$) converge or the norm of $\nabla \psi( \mathbf{a}_{i}  )$  is sufficiently small to assume that we have reached a stationary point of $\psi(.)$.

Of course, the linearization argument used to derive the Gauss-Newton iteration is only valid in a neighborhood of $\mathbf{a}_{i}$ and it may happens that $\psi( \mathbf{a}_{i+1} ) > \psi( \mathbf{a}_{i} )$ meaning that the Gauss-Newton algorithm may compute bad corrections by taking steps that are too long, reaching points outside the region of validity of the affine model used to approximate $\mathbf{r}(.)$ around $\mathbf{a}_i$. Several cleaver variants have been proposed to overcome this problem in practice.

The first one is the damped Gauss-Newton algorithm which is defined as
\begin{equation*}
\mathbf{a}_{i+1} = \mathbf{a}_{i} - \alpha_i \mathit{J} \big( \mathbf{r}(\mathbf{a}_i) \big)^{+}  \mathbf{r}(\mathbf{a}_i)  \  .
\end{equation*}
In this equation, $\alpha_i$ is a damping parameter which is chosen at each iteration to make the algorithm a descent method (i.e, such that $\psi( \mathbf{a}_{i+1} ) < \psi( \mathbf{a}_{i} )$). The Gauss-Newton approximation of the Hessian is always  positive semi-definite and it is positive definite, if and only if, the Jacobian matrix has full column rank, and, in this case, $d\mathbf{a}_{gn}$ is the unique solution of the above linear least-squares problem and is also in a descent direction for $\psi(.)$ if $d\mathbf{a}_{gn} \neq \mathbf{0}^{k.p}$ or, equivalently, if  $ \nabla \psi( \mathbf{a} ) = \mathit{J}( \mathbf{r}(\mathbf{a}) )^{T} \mathbf{r}(\mathbf{a})  \neq  \mathbf{0}^{k.p} $ since in these conditions
\begin{equation*} 
0 < d\mathbf{a}^{T}_{gn} \mathit{J} \big( \mathbf{r}(\mathbf{a}) \big)^{T} \mathit{J} \big( \mathbf{r}(\mathbf{a}) \big) d\mathbf{a}_{gn} = - d\mathbf{a}^{T}_{gn} \mathit{J} \big( \mathbf{r}(\mathbf{a}) \big)^{T} \mathbf{r}(\mathbf{a}) = - d\mathbf{a}^{T}_{gn} \nabla \psi( \mathbf{a} )  \  .
\end{equation*}
This shows that, when the Jacobian matrix has full column rank, the Gauss-Newton method can always be complemented with a line search in order to enforce the descending condition $\psi( \mathbf{a}_{i+1} ) < \psi( \mathbf{a}_{i} )$ during the iterations~\cite{NW2006}\cite{MN2010}\cite{HPS2012}. Here, $\mathbf{a}_{i+1} = \mathbf{a}_{i} +   \widehat{ \alpha } d\mathbf{a}_{gn}$ and $\widehat{ \alpha }$ is found as a (approximate) solution to the problem
\begin{equation*}
\widehat{ \alpha } \approx \text{Arg}\min_{\alpha>0}  \,  \psi( \mathbf{a}_{i} +   \alpha d\mathbf{a}_{gn} )  \  .
\end{equation*}
Many strategies have been proposed to choose the damping parameter $\widehat{ \alpha }$~\cite{DS1983}\cite{NW2006}. The Gauss-Newton method with a line search can be shown to have guaranteed convergence, provided that the level set  $\big\{  \mathbf{a} \in \mathbb{R}^{p.k} \text{ }  \vert \text{ } \psi( \mathbf{a} ) \le  \psi( \mathbf{a}_{0} ) \big\}$  is bounded, and the Jacobian matrix  $\mathit{J}( \mathbf{r}(\mathbf{a}_i) )$ has full rank in all iterations~\cite{DS1983}\cite{NW2006}. Practical experience shows that the Gauss-Newton method may fail with or without a line search and that it usually has only linear convergence as opposed to the Newton's method, which exhibits quadratic convergence near a solution vector $\mathbf{\widehat{a}}$. However, if, at a solution $\mathbf{\widehat{a}}$, we have $\mathbf{r}( \mathbf{\widehat{a}} ) = \mathbf{0}^{k.p}$, then we have the equality
\begin{equation*}
\nabla^2 \psi(\mathbf{\widehat{a}} ) = \mathit{J} \big( \mathbf{r}(\mathbf{\widehat{a}}) \big)^{T} \mathit{J} \big( \mathbf{r}(\mathbf{\widehat{a}}) \big)
\end{equation*}
and we can also get quadratic convergence with the Gauss-Newton method. Similarly, if the component residual functions $\mathbf{r}_{l}(.)$ have small curvatures or if the $\vert \mathbf{r}_{l}(\mathbf{\widehat{a}})  \vert$ are small, we can also get superlinear convergence. For example, this will be the case for the matrix completion problem. This can also be observed if the values of the residual matrix $\mathbf{X} - \mathbf{\widehat{A}} \mathbf{\widehat{B}}$ behave like white noise, as in this case we can expect partial canceling in the sum
\begin{equation*}
\sum_{l=1}^{n.p} \mathbf{r}_{l} (\mathbf{\widehat{a}}) \nabla^{2} \mathbf{r}_{l} (\mathbf{\widehat{a}})  \ ,
\end{equation*}
in which case, we also get
\begin{equation*}
\nabla^2 \psi(\mathbf{\widehat{a}} ) \approx \mathit{J} \big( \mathbf{r}(\mathbf{\widehat{a}}) \big )^{T} \mathit{J} \big ( \mathbf{r}(\mathbf{\widehat{a}}) \big )  \  .
\end{equation*}
This situation also occurs in many applications, especially in climate science.

When the Hessian matrix $\nabla^2 \psi(\mathbf{a} )$ is positive definite, the full Newton direction $d\mathbf{a}_{n}$ is also a descent direction for $\psi(.)$ and, in this case, the full Newton method can also be complemented by a  line search to enforce the descending condition $\psi( \mathbf{a}_{i+1} ) < \psi( \mathbf{a}_{i} )$~\cite{DS1983}\cite{NW2006}\cite{MN2010}\cite{HPS2012}. However, contrary to the Gauss-Newton approximation of the Hessian, which is always positive semi-definite, the full Hessian matrix can be indefinite in some regions of mixed curvature of the search space or even negative definite, in which cases, further regularization of the Hessian matrix, such that inflating its diagonal elements, is required to transform it in a positive definite matrix before applying a line search (see~\cite{NW2006}\cite{MN2010} and Subsection~\ref{vp_n_alg:box}  for more details).

The second modification of the Gauss-Newton algorithm used in practice is the Levenberg-Marquardt method. This method approximates the second term in the Hessian of  $\psi(.)$ with $\lambda . \mathbf{D}^{T} \mathbf{D}$ where $\mathbf{D}$ is a full rank matrix and $\lambda$ a strictly positive real scalar (the Marquardt damping parameter). The standard choice for $\mathbf{D}$ is the identity matrix or a diagonal matrix $\mathbf{D} = \emph{diag}( \mathbf{d} )$ with appropriately chosen components $\mathbf{d}_j > 0$, which take into account the scaling of the problem and can be kept fixed or changed during the iterations~\cite{M1978}\cite{DS1983}\cite{NW2006}\cite{MN2010}.
In all cases, this implies that the approximate Hessian matrix
\begin{equation*}
\mathit{J} \big ( \mathbf{r}(\mathbf{a}) \big )^{T} \mathit{J} \big ( \mathbf{r}(\mathbf{a}) \big ) + \lambda . \mathbf{D}^{T} \mathbf{D}
\end{equation*}
is  positive definite if $\lambda>0$. Thus, the Levenberg-Marquardt's method is based on the following quadratic approximation model
\begin{equation*}
\psi( \mathbf{a} +  d\mathbf{a} ) \approx L_{\lambda}( d\mathbf{a}  ) = \psi(\mathbf{a}) + d\mathbf{a}^{T} \mathit{J} \big ( \mathbf{r}(\mathbf{a}) \big )^{T} \mathbf{r}(\mathbf{a}) +  \frac{1}{2} d\mathbf{a}^{T}  \left( \mathit{J} \big( \mathbf{r}(\mathbf{a}) \big)^{T} \mathit{J} \big( \mathbf{r}(\mathbf{a}) \big ) + \lambda . \mathbf{D}^{T} \mathbf{D} \right) d\mathbf{a}
\end{equation*}
and the correction vector $d\mathbf{a}_{lm}$ for improving an approximate solution $\mathbf{a}$ is found as the solution of the regularized normal system
\begin{equation} \label{eq:lm_normal}
\left( \mathit{J}  \big( \mathbf{r}(\mathbf{a})  \big)^{T} \mathit{J}  \big( \mathbf{r}(\mathbf{a})  \big) + \lambda . \mathbf{D}^{T} \mathbf{D} \right) d\mathbf{a}_{lm} = -\mathit{J}  \big( \mathbf{r}(\mathbf{a})  \big)^{T} \mathbf{r}(\mathbf{a})
\end{equation}
and is always in a descent direction for  $\psi(.)$, even when $\mathit{J}( \mathbf{r}(\mathbf{a}) )$ is not of full column rank, if $\lambda > 0$. Rather than dividing the steps when $\psi( \mathbf{a}_{i+1} ) > \psi( \mathbf{a}_{i} )$ as in the damped Gauss-Newton method, the Levenberg-Marquardt algorithm deflates the steps by inflating the diagonals of the cross-product Jacobian matrix (which is equivalent to shift positively its spectrum) before inverting it to solve for the correction vector. It may be demonstrated that a sufficiently large $\lambda$ always exists such that $\psi( \mathbf{a}_{i+1} ) < \psi( \mathbf{a}_{i} )$ will be satisfied unless $\mathbf{a}_{i}$ is already a stationary point of $\psi(.)$~\cite{DS1983}\cite{NW2006}\cite{MN2010}\cite{HPS2012}.

In other words, the Marquardt damping parameter $\lambda$ controls the nature of the iterations and limits the size of $d\mathbf{a}_{lm}$ at the same time. If we assume that $\mathbf{D}$ is the identity matrix and $\lambda$ is very large, then
\begin{equation*}
d\mathbf{a}_{lm} \approx - \frac{1}{\lambda} . \mathit{J} \big ( \mathbf{r}(\mathbf{a}) \big )^{T} \mathbf{r}(\mathbf{a}) =  - \frac{1}{\lambda} . \nabla \psi( \mathbf{a} )
\end{equation*}
is a short step in a direction very close to the steepest descent direction. If, on the other hand, $\lambda$ is very small, then $L_{\lambda}( \mathbf{a} ) \approx G( \mathbf{a} )$  and $d\mathbf{a}_{lm}$ is close to the Gauss-Newton step $d\mathbf{a}_{gn}$ described above. In other words, we can think of the Levenberg-Marquardt method as a hybrid method between the steepest descent and Gauss-Newton methods with the good performance of the steepest descent method in the initial stage and the faster convergence of the Gauss-Newton method at the final stage of the iterative process, assuming that the value of the Marquardt damping parameter decreases during the iterative process. Taking $\mathbf{D}$ as the identity matrix corresponds to the algorithm originally proposed by Levenberg~\cite{L1944}. Later, Marquardt~\cite{M1963} improved the method by choosing the diagonals of $\mathbf{D}$ to match the 2-norms of the columns of the Jacobian matrix $\mathit{J}( \mathbf{r}(\mathbf{a}) )$. This makes the algorithm invariant under diagonal scaling of the elements of the vector $\mathbf{a}$~\cite{M1978}\cite{NW2006}. This also allows to include local curvature information, even when $\lambda$ is large and we are essentially moving in the (negative) steepest gradient direction. This is, for example, useful to alleviate the "error valley" problem affecting the steepest gradient method discussed at the beginning of this section since, in that case, we are moving further in the directions in which the gradient is smaller. Later, many other choices for $\mathbf{D}$ have been proposed and tested~\cite{M1978}\cite{DGW1981}\cite{DS1983}.

The above equations defining the Levenberg-Marquardt's correction vector $d\mathbf{a}_{lm}$ are the normal equations for the regularized linear least-squares problem
\begin{equation} \label{eq:lm_llsq}
\min_{d\mathbf{a} \in \mathbb{R}^{p.k}}   \,    \frac{1}{2} \big\Vert  \begin{bmatrix} \mathbf{r}(  \mathbf{a} ) \\ \mathbf{0}^{p.k} \end{bmatrix}  +  \begin{bmatrix}  \mathit{J} \big( \mathbf{r}(\mathbf{a}) \big)   \\  \sqrt{\lambda} . \mathbf{D}  \end{bmatrix}  d\mathbf{a}  \big\Vert^{2}_{2} = \frac{1}{2} \big\Vert \mathbf{r}(  \mathbf{a} ) + \mathit{J} \big( \mathbf{r}(\mathbf{a}) \big) d\mathbf{a}  \big\Vert^{2}_{2} +  \frac{\lambda}{2}  \big\Vert \mathbf{D} d\mathbf{a} \big\Vert^{2}_{2}  \ ,
\end{equation}
which can also be solved accurately by stable methods as for the Gauss-Newton correction and there is no need to form nor to invert the symmetric matrix $\mathit{J}( \mathbf{r}(\mathbf{a}) )^{T} \mathit{J}( \mathbf{r}(\mathbf{a}) ) + \lambda . \mathbf{D}^{T} \mathbf{D}$~\cite{HPS2012}\cite{NW2006}. Moreover, this linear least-squares problem has always a unique solution if $\lambda > 0$.

The Levenberg-Marquardt algorithm is often considered superior to the (damped) Gauss-Newton algorithm since it is well defined even when the Jacobian matrix is rank deficient. Another advantage is that the Levenberg-Marquardt correction assures an optimal interpolation between a Gauss-Newton step and the steepest descent direction (e.g., negative gradient direction) when the Gauss-Newton step is much too long.

Similarly, we can define a Levenberg-Marquardt variant of the Newton method by computing the correction vector $d\mathbf{a}_{n}$ as
\begin{equation} \label{eq:dn_normal}
\left(   \nabla^2 \psi( \mathbf{a} ) + \lambda . \mathbf{I}_{k.p}  \right) d\mathbf{a}_{n} = - \mathit{J} \big( \mathbf{r}(\mathbf{a}) \big)^{T} \mathbf{r}(\mathbf{a})  \ ,
\end{equation}
where the term $\lambda . \mathbf{I}_{k.p}$ is included when $\nabla^2 \psi( \mathbf{a} )$ is not positive definite and hence the Newton direction may not be a descent direction. In such conditions, it is always possible to choose  $\lambda$ sufficiently large such that, first, the matrix $\nabla^2 \psi( \mathbf{a} ) + \lambda . \mathbf{I}_{k.p}$ becomes positive definite and, second, $\psi( \mathbf{a} + d\mathbf{a}_{n}) < \psi( \mathbf{a} )$~\cite{NW2006}\cite{MN2010}. This strategy is based on the quadratic approximation model
\begin{equation*}
\psi( \mathbf{a} + d\mathbf{a} ) \approx  N_{\lambda}( d\mathbf{a} ) = \psi( \mathbf{a} ) +  d\mathbf{a}^{T}  \mathit{J} \big( \mathbf{r}(\mathbf{a}) \big)^{T} \mathbf{r}(\mathbf{a}) +  \frac{1}{2} d\mathbf{a}^{T} \left(  \nabla^2 \psi( \mathbf{a} ) + \lambda . \mathbf{I}_{k.p} \right) d\mathbf{a}  \ .
\end{equation*}
As in the Levenberg-Marquardt  algorithm, the damping parameter $\lambda$ can be used to control both the size and direction of the  correction vector $d\mathbf{a}_{n}$ and, in this case, we can avoid the use of a line search to control the step size in order to get reasonable convergence in the Newton method. Thus, we can also think of this Levenberg-Marquardt variant of the Newton method as an hybrid between the steepest descent and Newton methods with the good performance of the steepest descent method in the initial stage, but the quadratic convergence of the Newton method at the final stage~\cite{NW2006}\cite{MN2010}. See the variable projection Newton algorithms~\eqref {n_alg1:box},~\eqref {n_alg2:box} and~\eqref {n_alg3:box} described in Subsection~\ref{vp_n_alg:box}, which all integrate a damping term   $\lambda . \mathbf{I}_{k.p}$ for some illustrations of this simple strategy in the context of the Newton method applied to the WLRA problem.

A variation of the Levenberg-Marquardt method is the trust-region Gauss-Newton algorithm where the correction vector $d\mathbf{a}_{t-gn}$  is defined as the solution of the constrained linear least-squares problem
\begin{equation*}
d\mathbf{a}_{t-gn} = \text{Arg}\min_{d\mathbf{a} \in \mathbb{R}^{p.k}}   \,   \frac{1}{2} \big\Vert \mathbf{r}(  \mathbf{a} ) + \mathit{J} \big( \mathbf{r}(\mathbf{a}) \big) d\mathbf{a}  \big\Vert^{2}_{2}  \text{ subject to }    \big\Vert \mathbf{D} d\mathbf{a} \big\Vert_{2}  \le \delta   \ .
\end{equation*}
Here, the set of feasible correction vectors $d\mathbf{a}$ is restricted to the ellipsoid $\lbrace  d\mathbf{a} \in \mathbb{R}^{p.k}  /  \Vert \mathbf{D} d\mathbf{a} \Vert_{2}  \le \delta \rbrace$ which is called the trust region. $\delta > 0$ is the trust region radius, which controls the size of the trust region and is updated recursively during the iterative process~\cite{DS1983}\cite{NW2006}. In this class of methods, the scaling matrix $\mathbf{D}$ generates the elliptic norm $\Vert d\mathbf{a} \Vert_{\mathbf{D}} = \Vert \mathbf{D} d\mathbf{a} \Vert_{2}$ in which the correction vector is measured~\cite{NW2006}. The trust region can then be thought of as a region of trust for the linear model
\begin{equation*}
\mathbf{r}( \mathbf{a} + d\mathbf{a} ) \approx \mathbf{r}( \mathbf{a} ) + \mathit{J} \big ( \mathbf{r}(\mathbf{a}) \big) d\mathbf{a}
\end{equation*}
and the idea in the trust-region Gauss-Newton method is to avoid using this linear model outside its range of validity. Note that the Gauss-Newton step $d\mathbf{a}_{gn}$ solves this constrained problem if $\Vert \mathbf{D} d\mathbf{a}_{gn} \Vert_{2}  \le \delta$. Otherwise, it can be shown that the trust-region Gauss-Newton correction vector is the unique solution $d\mathbf{a}(\lambda)$ of the unconstrained regularized linear least-squares problem
\begin{equation*}
d\mathbf{a}(\lambda) = \text{Arg}\min_{d\mathbf{a} \in \mathbb{R}^{p.k}}   \,   \frac{1}{2} \big\Vert \mathbf{r}(  \mathbf{a} ) + \mathit{J} \big( \mathbf{r}(\mathbf{a}) \big) d\mathbf{a}  \big\Vert^{2}_{2} +  \frac{\lambda}{2}  \big\Vert \mathbf{D} d\mathbf{a}  \ ,\big\Vert^{2}_{2}  \ ,
\end{equation*}
where $\lambda > 0$ is determined from the scalar equation $\Vert \mathbf{D} d\mathbf{a}(\lambda) \Vert_{2}  = \delta$ which is nonlinear in $\lambda$. In other words, when the correction vector $d\mathbf{a}$ is directly controlled by the Marquardt damping parameter $\lambda$ and not by $\delta$, we obtain the Levenberg-Marquardt algorithm, otherwise we have a trust region Gauss-Newton algorithm~\cite{NW2006}. We also observe that if $\mathbf{D}$ is nonsingular then a change of variables yields an equivalent linear least-squares problem with $\mathbf{D} = \mathbf{I}$ for computing both the Levenberg-Marquardt and trust-region Gauss-Newton corrections.

Finally, the augmented Gauss-Newton method partly takes second-order derivatives into account by approximating the second term of the Hessian of $\psi(.)$ by either finite differencing or a quasi-Newton update in order to improve the above NLLS methods in the large residuals case~\cite{DGW1981}\cite{DS1983}\cite{NW2006}\cite{MN2010}. The variable projection quasi-Newton algorithms discussed in Subsection~\ref{vp_n_alg:box} belong to this class of methods.

What has been described so far is well-known. In the following subsections, let us quantify the smoothness of $\psi(.)$ in more details and study the specific properties of the Jacobian matrix $\mathit{J}( \mathbf{r}(\mathbf{a}) )$, Hessian matrix $\nabla^2  \psi( \mathbf{a} )$ and vector gradient $\nabla \psi( \mathbf{a} )$, which need to be evaluated in the above second-order or pseudo second-order NLLS algorithms. Implementation details of these variable projection NLLS algorithms will be presented in Section~\ref{vpalg:box} after their main properties have been derived in the rest of this section. For small or medium sized NLLS or WLRA problems, the above methods will be much faster than variants of the steepest gradient method. However, for larger problems, the cost of solving a linear least-squares problem or a linear system with a huge coefficient matrix at each iteration is a major drawback compared to the (steepest) gradient methods as the time spent in each iteration scales as $\mathcal{O}((k.p)^3)$ and, thus, increases considerably for large and square data matrices and a large value of the $k$ parameter. Taking these difficulties in consideration, we propose also some parallel implementations of all our variable projection NLLS algorithms for the WLRA problem in Section~\ref{vpalg:box} so that they can also be used for larger sized problems also found now in many practical applications.

\subsection{Computation and properties of the Jacobian matrix} \label{jacob:box}

In order to use a (damped or trust-region) Gauss-Newton or Levenberg-Marquardt algorithm for minimizing $\psi(.)$ (and solve the WLRA problem), we must compute the Jacobian of the residual function
\begin{equation*}
\mathbf{r}(\mathbf{a})  = \mathbf{P}_{\mathbf{F}(\mathbf{a})}^{\bot} \mathbf{x} = ( \mathbf{I}_{n.p} - \mathbf{P}_{\mathbf{F}(\mathbf{a})} )\mathbf{x} \ ,
\end{equation*}
defined in equation~\eqref{eq:r_vec}. This requires computing the derivative of the orthogonal projector $\mathbf{P}_{\mathbf{F}(.)}$ with respect to $\mathbf{a} \in h^{-1} ( \mathbb{R}^{p \times k}_{k} )$ as shown in Subsection~\ref{varpro_wlra:box}. If $\mathbf{a} \in h^{-1} ( \mathbb{R}^{p \times k}_{<k} )$, keep in mind that $\mathbf{P}_{\mathbf{F}(.)}$ is not even continuous at $\mathbf{a}$ (see Theorems~\ref{theo3.11:box} and~\ref{theo3.12:box}) and cannot be differentiable either at this point.

A close formula for the derivative of orthogonal projectors has been derived first by Golub and Pereyra~\cite{GP1973} and Decell~\cite{D1974}  under the assumption that $\mathbf{F}(.)$ is of local constant rank at any point $\mathbf{a}$ (this means that $\mathbf{F}(.)$ is of constant rank in a neighborhood of $\mathbf{a}$, but not necessarily of full column-rank, see Definition~\ref{def3.1:box} for details) in which differentiation is to be performed as stated in the following theorem, which extends the results about the continuity of $\mathbf{P}_{\mathbf{F}(.)}$ given in Theorem~\ref{theo3.10:box}:
\begin{theo5.1} \label{theo5.1:box}
Let $\Phi(.)$ be a matrix function : $\mathbb{R}^{m} \longrightarrow \mathbb{R}^{l \times t}$, which is $q$ times continuously differentiable at a point $\mathbf{a} \in  \mathbb{R}^{m}$. The following conditions are equivalent:
\begin{align*}
1) &  \quad\  \Phi(.)   \text{ has a local constant rank at }  \mathbf{a}  \ , \\
2) &  \quad\  \Phi(.)^{+}   \text{ is q times continuously differentiable at }  \mathbf{a}  \ , \\
3) &  \quad\  \Phi(.) \Phi(.)^{+} =  \mathbf{P}_{\Phi(.)} \text{ is q times continuously differentiable at }  \mathbf{a}  \ , \\
4) &  \quad\  \Phi(.)^{+} \Phi(.) \text{ is q times continuously differentiable at }  \mathbf{a}  \ .
\end{align*}
In other words, the differentiability of the pseudo-inverse of a matrix function $\Phi(.)$  at a point $\mathbf{a} \in \mathbb{R}^{m}$ is equivalent to the differentiability of the orthogonal projectors onto the column or row spaces of this matrix function at $\mathbf{a}$ and all these conditions are equivalent to the assertion that this matrix function has local constant rank at $\mathbf{a}$ if $\Phi(.)$ is itself differentiable at  $\mathbf{a}$.
Furthermore, in these conditions, we have for any point $\mathbf{a} \in \mathbb{R}^{m}$ for which $\Phi(.)$ is differentiable
\begin{equation} \label{eq:d_projmat}
\mathit{D}( \mathbf{P}_{\Phi(\mathbf{a})} ) = \mathbf{P}_{\Phi(\mathbf{a})}^{\bot}  \mathit{D}( \Phi(\mathbf{a}) )  \Phi(\mathbf{a})^{+}   +  \big( \mathbf{P}_{\Phi(\mathbf{a})}^{\bot}  \mathit{D}( \Phi(\mathbf{a}) ) \Phi(\mathbf{a})^{+} \big)^{T} 
\end{equation}
and
\begin{align} \label{eq:d_ginvmat}
    \mathit{D}( \Phi(\mathbf{a})^{+} ) = & - \Phi(\mathbf{a})^{+} \mathit{D}( \Phi(\mathbf{a}) ) \Phi(\mathbf{a})^{+} + \Phi(\mathbf{a})^{+} (\Phi(\mathbf{a})^{+})^{T} \mathit{D}( \Phi(\mathbf{a})^{T} ) \mathbf{P}_{\Phi(\mathbf{a})}^{\bot}  \nonumber \\
    & + \big( \mathbf{I}_{m} - \Phi(\mathbf{a})^{+} \Phi(\mathbf{a}) \big) \mathit{D}( \Phi(\mathbf{a})^{T} ) (\Phi(\mathbf{a})^{+})^{T} \Phi(\mathbf{a})^{+}  \ .
\end{align}
Finally, note that, in the above equation defining the differential of the orthogonal projector $\mathbf{P}_{\Phi(.)}$, we can substitute in place of the pseudo-inverse $\Phi(\mathbf{a})^{+}$ any symmetric generalized inverse $\Phi(\mathbf{a})^{-}$ as defined in equations~\eqref{eq:sginv} or~\eqref{eq:sginv_qrcp} of Subsection~\ref{lin_alg:box}.
\end{theo5.1}
\begin{proof}
See Theorems 8.4 and 8.5 in Chapter 8 of~\cite{MN2019} and also~\cite{GP1973}\cite{GP1976}\cite{D1974}\cite{CM2009}.
\\
\end{proof}

As noted already in Subsection~\ref{varpro_wlra:box},  $\mathbf{F}(.)$ is a continuous linear mapping from $\mathbb{R}^{p.k}$ into $\mathbb{R}^{n.p \times n.k}$ (since the $\emph{mat}$ and transpose operators are linear mappings and the Kronecker product is a bilinear operator) and is, thus, continuously and infinitely differentiable at any point $\mathbf{a} \in  \mathbb{R}^{p.k}$~\cite{C2017}. Collecting the results from Theorems~\ref{theo3.10:box} and~\ref{theo5.1:box}, we then deduce that the proposition that  $\mathbf{P}_{\mathbf{F}(.)}$ is infinitely differentiable (e.g., of class $C^\infty$) at a point $\mathbf{a} \in h^{-1} ( \mathbb{R}^{p \times k}_k )$ is equivalent to its continuity at this point and to the proposition that $\mathbf{F}(.)$ is of constant rank in a neighborhood of $\mathbf{a}$. Next, using Theorem~\ref{theo3.11:box}, we obtain the following corollary in the case where $\mathbf{W} \in \mathbb{R}^{p \times n}_{+*}$:
\begin{corol5.1} \label{corol5.1:box}
For $\mathbf{X} \in \mathbb{R}^{p \times n} \text{ and } \mathbf{W} \in \mathbb{R}^{p \times n}_{+*}$, and any fixed integer $k \le \emph{rank}( \mathbf{X} ) $, the matrix function $\mathbf{P}_{\mathbf{F}(.)}$ from $\mathbb{R}^{p.k}$ to $\mathbb{R}^{p.n \times p.n}$ defined by
\begin{equation*}
\mathbf{a} \mapsto \mathbf{P}_{\mathbf{F}(  \mathbf{a}  )} = \mathbf{F}(  \mathbf{a} ) \mathbf{F}(  \mathbf{a} )^{+} \ ,
\end{equation*}
where $\mathbf{F}(  \mathbf{a} )^{+}$ is the pseudo-inverse of $\mathbf{F}(  \mathbf{a} )$ and $\mathbf{F}(  \mathbf{a} )$ is the $ p.n \times n.k  $ block diagonal matrix 
\begin{equation*}
\mathbf{F}(  \mathbf{a} ) =  \bigoplus_{j=1}^n \emph{diag}(\sqrt{\mathbf{W}}_{.j}) h(  \mathbf{a}  ) = \bigoplus_{j=1}^n \emph{diag}(\sqrt{\mathbf{W}}_{.j})\mathbf{A}  \ ,
\end{equation*}
 is of class $C^\infty$ (e.g., infinitely differentiable) at all points $\mathbf{a} \in h^{-1} ( \mathbb{R}^{p \times k}_k )$.
Furthermore, for all points $\mathbf{a} \in h^{-1} ( \mathbb{R}^{p \times k}_k )$, we have
\begin{equation} \label{eq:d_projmat2}
\mathit{D}( \mathbf{P}_{\mathbf{F}(\mathbf{a})} ) = \mathbf{P}_{\mathbf{F}(\mathbf{a})}^{\bot}  \mathit{D} \big( \mathbf{F}(\mathbf{a})  \big) \mathbf{F}(\mathbf{a})^{+}   +  \big( \mathbf{P}_{\mathbf{F}(\mathbf{a})}^{\bot}  \mathit{D} \big( \mathbf{F}(\mathbf{a})  \big) \mathbf{F}(\mathbf{a})^{+} \big)^{T}  \ .
\end{equation}
Here, as in equation~\eqref{eq:d_projmat} of Theorem~\ref{theo5.1:box}, we can substitute in place of $\mathbf{F}(\mathbf{a})^{+}$ any symmetric generalized inverse $\mathbf{F}(\mathbf{a})^{-}$ as defined in equations~\eqref{eq:sginv} or~\eqref{eq:sginv_qrcp} of Subsection~\ref{lin_alg:box}.
\\
\end{corol5.1}
As expected from Corollary~\ref{corol3.4:box}, the situation is much less favourable when $\mathbf{W}$ has some zero elements, as the condition that $\mathbf{a} \in h^{-1} ( \mathbb{R}^{p \times k}_k )$ is not sufficient to ensure that $\mathbf{F}(.)$ is of constant rank in a neighborhood of $\mathbf{a}$ and, thus, that $\mathbf{P}^{\bot}_{\mathbf{F}(.)}$ is differentiable at $\mathbf{a}$ in such situation:
\\
\begin{corol5.2} \label{corol5.2:box}
For $\mathbf{X} \in \mathbb{R}^{p \times n} \text{ and } \mathbf{W} \in \mathbb{R}^{p \times n}_{+}$, and any fixed integer $k \le \emph{rank}( \mathbf{X} ) $, the matrix function $\mathbf{P}_{\mathbf{F}(.)}$ from $\mathbb{R}^{p.k}$ into $\mathbb{R}^{p.n \times p.n}$ defined by
\begin{equation*}
\mathbf{a} \mapsto \mathbf{P}_{\mathbf{F}(  \mathbf{a}  )} = \mathbf{F}(  \mathbf{a} ) \mathbf{F}(  \mathbf{a} )^{+}
\end{equation*}
is not differentiable at all points $\mathbf{a}  \in  \bigcup_{j=1}^{n} \mathcal{B}_j$, where $\mathcal{B}_j$ is the  $j^{th}$ barrier set associated with the $j^{th}$ atomic and matrix functions, $\psi_{j}(.)$ and $\mathbf{F}_{j}(.)$, as defined, respectively, in equation~\eqref{eq:psi_atomic_func} and Definition~\ref{def3.2:box}.
\\
\end{corol5.2}

Despite the caveats stated in Corollary~\ref{corol5.2:box} when some elements of $\mathbf{W}$ are equal to zero, it is important to keep in mind that the general differential formula~\eqref{eq:d_projmat2} is still valid in that case as soon as $\mathbf{F}(.)$ has a local constant rank at $\mathbf{a} \in h^{-1} ( \mathbb{R}^{p \times k}_k )$. Furthermore, previous comparative studies have also demonstrated that first- and second-order variable projection methods used for minimizing $\psi(.)$ generally outperform other concurrent methods even for a large number of missing values in the case of binary weights and without any form of regularization to ensure the smoothness of $\psi(.)$ despite the non differentiability of $\mathbf{P}_{\mathbf{F}(.)}$ in some regions of the search space $h^{-1} ( \mathbb{R}^{p \times k}_k )$~\cite{C2008b}\cite{OYD2011}\cite{HF2015}\cite{HZF2017}.

Here, $\mathit{D}( \mathbf{F}(\mathbf{a}) )$ and $\mathit{D}( \mathbf{P}_{\mathbf{F}(\mathbf{a})} )$ are, for $\mathbf{a} \in h^{-1} ( \mathbb{R}^{p \times k}_k )$, elements of $\pounds ( \mathbb{R}^{p.k}, \pounds ( \mathbb{R}^{n.k}, \mathbb{R}^{n.p} ) )$ and $\pounds ( \mathbb{R}^{p.k}, \pounds ( \mathbb{R}^{n.p}, \mathbb{R}^{n.p} ) )$, respectively, and could be interpreted as tridimensional tensors (see equation~\eqref{eq:D_mat_tensor} in Subsection~\ref{calculus:box}). Now, since $\mathbf{P}_{\mathbf{F}(\mathbf{a})}^{\bot} = \mathbf{I}_{n.p} - \mathbf{P}_{\mathbf{F}(\mathbf{a})}$, we then have
\begin{equation*}
\mathit{D}( \mathbf{P}_{\mathbf{F}(\mathbf{a})}^{\bot} ) =  \mathit{D}( \mathbf{I}_{n.p} - \mathbf{P}_{\mathbf{F}(\mathbf{a})} ) = - \mathit{D}( \mathbf{P}_{\mathbf{F}(\mathbf{a})} )
\end{equation*}
and we deduce by the product differentiation rule~\cite{C2017} that
\begin{equation} \label{eq:r_jacob}
   \mathit{J}( \mathbf{r}(\mathbf{a}) ) = \mathit{J}( \mathbf{P}_{\mathbf{F}(\mathbf{a})}^{\bot}\mathbf{x} ) =  \mathit{D}( \mathbf{P}_{\mathbf{F}(\mathbf{a})}^{\bot} )\mathbf{x} + \mathbf{P}_{\mathbf{F}(\mathbf{a})}^{\bot} \mathit{J}( \mathbf{x} ) = - \mathit{D}( \mathbf{P}_{\mathbf{F}(\mathbf{a})} )\mathbf{x} \ .
\end{equation}
Substituting now for $\mathit{D}( \mathbf{P}_{\mathbf{F}(\mathbf{a})} )$ yields
\begin{IEEEeqnarray*}{rClClr}
      \mathit{J}( \mathbf{r}(\mathbf{a}) ) & =  & - \big( \thickspace \mathbf{P}_{\mathbf{F}(\mathbf{a})}^{\bot}  \mathit{D}( \mathbf{F}(\mathbf{a}) ) \mathbf{F}(\mathbf{a})^{+}\mathbf{x} & + & (\mathbf{F}(\mathbf{a})^{+})^{T} \mathit{D}( \mathbf{F}(\mathbf{a}) )^{T} (\mathbf{P}_{\mathbf{F}(\mathbf{a})}^{\bot})^{T}\mathbf{x}  & \thickspace \big) \\
      & = & - \big(\thickspace \mathbf{P}_{\mathbf{F}(\mathbf{a})}^{\bot}  \mathit{D}( \mathbf{F}(\mathbf{a}) ) \mathbf{\widehat{b}} & + & (\mathbf{F}(\mathbf{a})^{+})^{T} \mathit{D}( \mathbf{F}(\mathbf{a}) )^{T}  \mathbf{r}(\mathbf{a}) & \thickspace \big) \ ,
\end{IEEEeqnarray*}
where we have used the fact that $\mathbf{P}_{\mathbf{F}(\mathbf{a})}^{\bot}$ is a symmetric matrix (see Subsection~\ref{lin_alg:box}). In these equations, $\mathbf{\widehat{b}} = \mathbf{F}(\mathbf{a})^{+}\mathbf{x}$ and $\mathbf{r}(\mathbf{a})= \mathbf{P}_{\mathbf{F}(\mathbf{a})}^{\bot}\mathbf{x}$ are, respectively, the minimum Euclidean norm solution and residual vector of the following linear least-squares problem already encountered when describing the block ALS method in Section~\ref{nipals:box}
\begin{equation*}
\min_{\mathbf{b}\in\mathbb{R}^{n.k}} \enspace\ \frac{1}{2} \Vert \mathbf{x} - \mathbf{F}(\mathbf{a})\mathbf{b} \Vert^{2}_2 =  \varphi^{*}(\mathbf{A},\mathbf{B}  ) \ ,
\end{equation*}
where $\mathbf{B} = \emph{mat}( \mathbf{b} )$. Note that we can also use $\mathbf{\widehat{b}} = \mathbf{F}(\mathbf{a})^{-}\mathbf{x}$, which is cheaper to evaluate, in the above equations. Moreover, we recall that the linear mappings $\mathit{D}( \mathbf{F}(\mathbf{a}) ) \mathbf{\widehat{b}}$ and $\mathit{D}( \mathbf{F}(\mathbf{a}) )^{T}  \mathbf{r}(\mathbf{a})$ are elements of $\pounds ( \mathbb{R}^{p.k}, \mathbb{R}^{n.p} )$ and $\pounds ( \mathbb{R}^{p.k}, \mathbb{R}^{n.k} )$, respectively, since transposition in the tensor $\mathit{D}( \mathbf{F}(\mathbf{a}) )$ is performed on each slab $\partial\mathbf{F}(\mathbf{a})/ \partial\mathbf{a}_{i}$. See equation~\eqref{eq:D_mat_tensor} in Subsection~\ref{calculus:box} for details. Thus, these two factors correspond to $n.p \times p.k $ and $n.k \times p.k$ matrices, respectively.

We now derive an explicit formulation for the $n.p \times p.k$ matrix $\mathit{J}( \mathbf{r}(\mathbf{a}) )$, which is independent of the differentiability of the residual function $\mathbf{r}(.)$ and the existence of the "true" Jacobian matrix of this residual function. We first consider the first term in $\mathit{J}( \mathbf{r}(\mathbf{a}) )$, i.e.,
\begin{equation*}
\mathbf{M}(\mathbf{a}) = \mathbf{P}_{\mathbf{F}(\mathbf{a})}^{\bot} \mathit{D} \big( \mathbf{F}(\mathbf{a}) \big) \mathbf{\widehat{b}} \ ,
\end{equation*}
which is also a $n.p \times p.k$ matrix. As derived in equation~\eqref{eq:F_mat} of Subsection~\ref{varpro_wlra:box}, $\mathbf{F}(\mathbf{a})$ may be expressed in the form
\begin{equation*}
    \mathbf{F}(\mathbf{a}) = \emph{diag} \big( \emph{vec}(\sqrt{\mathbf{W}}) \big) (\mathbf{I}_{n}\otimes\mathbf{A}) = \emph{diag} \big(\emph{vec}(\sqrt{\mathbf{W}}) \big)  \big(\mathbf{I}_{n}\otimes \emph{mat}_{k \times p} (\mathbf{a})^{T} \big)
\end{equation*}
and it is clear that $\mathbf{F}(.)$ is a continuous linear mapping from $\mathbb{R}^{p.k}$ into $\mathbb{R}^{n.p \times n.k}$ since the $\emph{mat}$ and transpose operators are linear mappings and the Kronecker and matrix products are bilinear operators. Hence, $\forall \thickspace \mathbf{a}, \triangle\mathbf{a}  \in \mathbb{R}^{p.k}$ and $\triangle\mathbf{A} = h(\triangle\mathbf{a}) = \emph{mat}_{k \times p} (\triangle\mathbf{a})^{T}$, we have
\begin{equation*}
\mathit{D}\big ( \mathbf{F}(\mathbf{a}) \big) (\triangle\mathbf{a}) = \mathbf{F}(\triangle\mathbf{a}) = \emph{diag} \big( \emph{vec}(\sqrt{\mathbf{W}}) \big) (\mathbf{I}_{n}\otimes\triangle\mathbf{A})  \ .
\end{equation*}
Noting that (see equation~\eqref{eq:vec_kronprod} in Subsection~\ref{multlin_alg:box} )
\begin{equation*}
 \mathbf{F}(\triangle\mathbf{a})\mathbf{\widehat{b}} = \emph{diag} \big(\emph{vec}(\sqrt{\mathbf{W}}) \big) (\mathbf{I}_{n}\otimes\triangle\mathbf{A})\emph{vec}(\mathbf{\widehat{B}}) = \emph{diag} \big( \emph{vec}(\sqrt{\mathbf{W}}) \big) (\mathbf{\widehat{B}}^{T} \otimes \mathbf{I}_{p})\emph{vec}(\triangle\mathbf{A})  \ ,
\end{equation*}
where $\mathbf{\widehat{B}} = \emph{mat}_{k \times n} (\mathbf{\widehat{b}})$ and using the $\emph{p.k} \times \emph{p.k}$ commutation matrix $\mathbf{K}_{(k,p)}$ (see equation~\eqref{eq:commat} in Subsection~\ref{multlin_alg:box}), we deduce that
\begin{equation} \label{eq: deriv_F(a)_delta_b}
(\mathit{D} \big( \mathbf{F}(\mathbf{a})  \big)(\triangle\mathbf{a}))\mathbf{\widehat{b}} = \emph{diag}  \big(\emph{vec}(\sqrt{\mathbf{W}}) \big) (\mathbf{\widehat{B}}^{T} \otimes \mathbf{I}_{p}) \mathbf{K}_{(k,p)} \triangle\mathbf{a}  \ ,
\end{equation}
since
\begin{equation*}
\emph{vec}(\triangle\mathbf{A}) = \mathbf{K}_{(k,p)} \emph{vec}(\triangle\mathbf{A}^{T}) = \mathbf{K}_{(k,p)} \triangle\mathbf{a} \ ,
\end{equation*}
following our conventions for the vectorized form of the $\mathbf{A}$ matrix defined in equation~\eqref{eq:d_veca} of Subsection~\ref{varpro_wlra:box}. In view of this, we finally obtain the following explicit formulation for the $ n.p \times p.k $ matrix $\mathbf{M}(\mathbf{a})$
\begin{equation}\label{eq:M_mat}
\mathbf{M}(\mathbf{a}) = \mathbf{P}_{\mathbf{F}(\mathbf{a})}^{\bot} \emph{diag} \big(\emph{vec}(\sqrt{\mathbf{W}}) \big) ( \mathbf{\widehat{B}}^{T} \otimes \mathbf{I}_{p} ) \mathbf{K}_{(k,p)}  \ .
\end{equation}

An alternative useful formulation of the  $\mathbf{M}(\mathbf{a})$ matrix may be derived by noting that (see equation~\eqref{eq:com_kron} and Lemma~\ref{theo2.2:box} in Subsection~\ref{multlin_alg:box} )
\begin{align*}
\emph{diag} \big(\emph{vec}(\sqrt{\mathbf{W}}) \big) (\mathbf{\widehat{B}}^{T} \otimes \mathbf{I}_{p}) \mathbf{K}_{(k,p)} & = \emph{diag} \big(\emph{vec}(\sqrt{\mathbf{W}}) \big) \mathbf{K}_{(n,p)} ( \mathbf{I}_{p} \otimes  \mathbf{\widehat{B}}^{T} ) \\
     & =   \mathbf{K}_{(n,p)}  \emph{diag} \big(\emph{vec}(\sqrt{\mathbf{W}}^{T}) \big)  ( \mathbf{I}_{p} \otimes  \mathbf{\widehat{B}}^{T} )  \\
     & =   \mathbf{K}_{(n,p)}   \mathbf{G}(\mathbf{\widehat{b}}) \ ,
\end{align*}
where $\mathbf{G}(\mathbf{\widehat{b}})$ is defined in equation~\eqref{eq:G_mat} of Subsection~\ref{varpro_wlra:box}. Thus,
\begin{equation} \label{eq:M_mat2}
\mathbf{M}(\mathbf{a}) = \mathbf{P}_{\mathbf{F}(\mathbf{a})}^{\bot} \mathbf{K}_{(n,p)} \mathbf{G}(\mathbf{\widehat{b}}) \  ,
\end{equation}
which will be used later, in particular in Theorem~\ref{theo5.3:box} and for computing $\nabla \psi( \mathbf{a} )$ in Subsection~\ref{hess:box} (see Theorem~\ref{theo5.7:box}). As (see the demonstration of Theorem~\ref{theo4.3:box} for details)
\begin{equation*}
e(\mathbf{A},\mathbf{\widehat{B}}) =  \mathbf{x} - \mathbf{F}(  \mathbf{a} ) \mathbf{\widehat{b}} = \mathbf{K}_{(n,p)}  \big( \mathbf{z}  -  \mathbf{G}(  \mathbf{\widehat{b}} ) \mathbf{a} \big) \  ,
\end{equation*}
where $\mathbf{z} = \emph{vec} \big( (\sqrt{\mathbf{W}}  \odot \mathbf{X})^{T} \big)$, we can also write
\begin{equation} \label{eq:M_mat3}
\mathbf{M}(\mathbf{a}) = - \mathbf{P}_{\mathbf{F}(\mathbf{a})}^{\bot} \frac{\partial e(\mathbf{A},\mathbf{\widehat{B}})}{\partial\mathbf{a}}  \  .
\end{equation}
\\
In order to evaluate the second term in $\mathit{J}( \mathbf{r}(\mathbf{a}) )$, i.e.,
\begin{equation*}
\mathbf{L}(\mathbf{a}) =  \big( {\mathbf{F}(\mathbf{a})}^{+} \big)^{T}  \mathit{D} \big( \mathbf{F}(\mathbf{a})  \big)^{T} \mathbf{r}(\mathbf{a})  \ ,
\end{equation*}
which corresponds also to a $n.p \times p.k$ matrix, we first remark that, $\forall \thickspace \mathbf{a}, \triangle\mathbf{a}  \in \mathbb{R}^{p.k}$ and $\triangle\mathbf{A} = h(\triangle\mathbf{a}) = \emph{mat}_{k \times p} (\triangle\mathbf{a})^{T}$, we have
\begin{align*}
    \Big( \mathit{D} \big( \mathbf{F}(\mathbf{a}) \big) (\triangle\mathbf{a} ) \Big)^{T} & =  \mathbf{F}(\triangle\mathbf{a})^{T}  \\
    & = \Big( \emph{diag} \big( \emph{vec}(\sqrt{\mathbf{W}} ) \big) ( \mathbf{I}_{n} \otimes \triangle\mathbf{A} ) \Big)^{T}  \\
    & = ( \mathbf{I}_{n} \otimes \triangle\mathbf{A}^{T} ) \emph{diag} \big( \emph{vec}( \sqrt{\mathbf{W}} ) \big) \ ,
\end{align*}
since $\mathbf{F}(.)$  is a linear mapping and the transpose operator distributes over the Kronecker product. Now, $\forall \thickspace \mathbf{Z} \in \mathbb{R}^{p \times n}$ and $\forall \thickspace \mathbf{a}, \triangle\mathbf{a}  \in \mathbb{R}^{p.k}$, using equation~\eqref{eq:vec_kronprod}, we have
\begin{align*}
    \big( \mathit{D}( \mathbf{F}(\mathbf{a}) )( \triangle\mathbf{a} ) \big)^{T} \emph{vec}(\mathbf{Z}) & =  ( \mathbf{I}_{n} \otimes \triangle\mathbf{A}^{T} ) \emph{diag}  \big( \emph{vec}( \sqrt{\mathbf{W}} )  \big) \emph{vec}( \mathbf{Z} ) \\
          & =  ( \mathbf{I}_{n} \otimes \triangle\mathbf{A}^{T} ) \emph{vec}( \sqrt{\mathbf{W}} \odot \mathbf{Z} )  \\
          & =  \emph{vec} \big( \triangle\mathbf{A}^{T} ( \sqrt{\mathbf{W}} \odot \mathbf{Z} ) \big)  \\
          & =  \big( ( \sqrt{\mathbf{W}} \odot \mathbf{Z} )^{T} \otimes \mathbf{I}_{k} \big) \emph{vec}( \triangle\mathbf{A}^{T} ) \\
          & =  \big( ( \sqrt{\mathbf{W}} \odot \mathbf{Z} )^{T} \otimes \mathbf{I}_{k} \big) \triangle\mathbf{a}  \ ,
\end{align*}
and, thus, the $n.k \times p.k$ matrix representing the linear mapping $\mathit{D}( \mathbf{F}(\mathbf{a})(.) )^{T} \emph{vec}( \mathbf{Z} )$ is
\begin{equation*}
( \sqrt{\mathbf{W}} \odot \mathbf{Z} )^{T} \otimes \mathbf{I}_{k} \ .
\end{equation*}
Now, using the projection operator $P_{\Omega}(.)$ associated with the $p \times n$ weight matrix $\mathbf{W}$ defined in equation~\eqref{eq:proj_op}, we have
\begin{equation*}
     \big \lbrack P_{\Omega}(\mathbf{X} -\mathbf{A}\mathbf{\widehat{B}}) \big \rbrack_{ij} =
    \begin{cases}
        \displaystyle{ \mathbf{X}_{ij} - \sum_{l=1}^{k} \mathbf{A}_{il} \mathbf{\widehat{B}}_{lj} } & \text{if } \mathbf{W}_{ij} \ne 0 \\
         0                                                                                                                                      & \text{if } \mathbf{W}_{ij}  =   0 
    \end{cases} \ ,
\end{equation*}
and the variable projection residual vector of $\mathbf{x}$  at $\mathbf{A}$ can be written as
\begin{equation*}
    \mathbf{r}(\mathbf{a}) = \emph{vec} \big( \sqrt{\mathbf{W}} \odot P_{\Omega}(\mathbf{X} -\mathbf{A}\mathbf{\widehat{B}}) \big)
\end{equation*}
and it follows that, $\forall \thickspace \mathbf{a}, \triangle\mathbf{a}  \in \mathbb{R}^{p.k}$,
\begin{equation} \label{eq: deriv_F(a)_delta_T_r(a)}
    \Big(( \mathit{D} \big( \mathbf{F}(\mathbf{a}) \big)( \triangle\mathbf{a} ) \Big)^{T} \mathbf{r}(\mathbf{a}) =   \left( \big( \sqrt{\mathbf{W}}  \odot \sqrt{\mathbf{W}}  \odot P_{\Omega}(\mathbf{X} -\mathbf{A}\mathbf{\widehat{B}}) \big)^{T} \otimes \mathbf{I}_{k} \right) \triangle\mathbf{a}  \ ,
\end{equation}
hence 
\begin{equation}\label{eq:L_mat}
    \mathbf{L}(\mathbf{a}) = \big( \mathbf{F}(\mathbf{a})^{+} \big)^{T} \left( \big( \mathbf{W} \odot P_{\Omega}(\mathbf{X} -\mathbf{A}\mathbf{\widehat{B}}) \big)^{T} \otimes \mathbf{I}_{k} \right) \  .
\end{equation}

At this point, we will introduce two new intermediate quantities to simplify the notation going forward, especially in the computation of the Hessian matrix in the next section:
\begin{equation}\label{eq:U_mat}
    \mathbf{U} (\mathbf{a}) = \emph{diag} \big( \emph{vec}(\sqrt{\mathbf{W}} ) \big) (\mathbf{\widehat{B}}^{T} \otimes \mathbf{I}_{p}) \mathbf{K}_{(k,p)} = \mathbf{K}_{(n,p)} \mathbf{G}(\mathbf{\widehat{b}})
\end{equation}
and
\begin{equation}\label{eq:V_mat}
    \mathbf{V} (\mathbf{a}) = \big( \mathbf{W} \odot P_{\Omega}(\mathbf{X} -\mathbf{A}\mathbf{\widehat{B}}) \big)^{T} \otimes \mathbf{I}_{k} \ .
\end{equation} 

With these definitions, we have finally,
\begin{equation} \label{eq:J_mat}
    \mathit{J} \big( \mathbf{r}(\mathbf{a}) \big) = - \big( \mathbf{M}(\mathbf{a}) + \mathbf{L}(\mathbf{a}) \big) = - \Big( \mathbf{P}_{\mathbf{F}(\mathbf{a})}^{\bot} \mathbf{U}(\mathbf{a}) + \big( \mathbf{F}(\mathbf{a})^{+}  \big)^{T}  \mathbf{V}(\mathbf{a}) \Big) \ .
\end{equation}

We now demonstrate several important results concerning the ranges and null spaces associated with the $\mathbf{M}(\mathbf{a})$, $\mathbf{L}(\mathbf{a})$ and $\mathit{J}( \mathbf{r}(\mathbf{a}) )$ matrices, which result directly from the use of the variable projection method.

First, we have $ \mathbf{M}(\mathbf{a}) = \mathbf{P}_{\mathbf{F}(\mathbf{a})}^{\bot} \mathbf{U} (\mathbf{a}) $ and this leads to $ \emph{ran}(\mathbf{M}(\mathbf{a})) \subset  \emph{ran}( \mathbf{P}_{\mathbf{F}(\mathbf{a})}^{\bot} ) =  \emph{ran}( \mathbf{F}(\mathbf{a}) )^{\bot} $. Using the properties of the Moore-Penrose inverse (see equation~\eqref{eq:ginv} or more generally of any symmetric generalized inverse of the form~\eqref{eq:sginv_qrcp} defined in Subsection~\ref{lin_alg:box}), we also have
\begin{equation*}
    ( \mathbf{F}(\mathbf{a})^{+} )^{T} = \big( \mathbf{F}(\mathbf{a})^{+} \mathbf{F}(\mathbf{a}) \mathbf{F}(\mathbf{a})^{+} \big)^{T} = \big( \mathbf{F}(\mathbf{a}) \mathbf{F}(\mathbf{a})^{+} \big)^{T} \big( \mathbf{F}(\mathbf{a})^{+} \big)^{T} = \mathbf{P}_{\mathbf{F}(\mathbf{a})} \big( \mathbf{F}(\mathbf{a})^{+} \big)^{T}
\end{equation*}
and we deduce that
\begin{equation}\label{eq:L_mat2}
   \mathbf{L}(\mathbf{a}) = \mathbf{P}_{\mathbf{F}(\mathbf{a})} \big(  \mathbf{F}(\mathbf{a})^{+}  \big)^{T}  \mathbf{V} (\mathbf{a})
\end{equation}
and $\emph{ran}( \mathbf{L}(\mathbf{a}) ) \subset \emph{ran}( \mathbf{P}_{\mathbf{F}(\mathbf{a})} ) = \emph{ran}( \mathbf{F}(\mathbf{a}) )$. Hence the subspaces $\emph{ran}( \mathbf{M}(\mathbf{a}) )$ and $\emph{ran}( \mathbf{L}(\mathbf{a}) )$ of $\mathbb{R}^{p.n}$ are orthogonal and $\emph{ran}( \mathbf{M}(\mathbf{a}) ) \cap \emph{ran}( \mathbf{L}(\mathbf{a}) ) = \{ \mathbf{0}^{p.n} \}$. Now, since $\mathit{J}( \mathbf{r}(\mathbf{a}) ) = - ( \mathbf{M}(\mathbf{a}) + \mathbf{L}(\mathbf{a}) )$, any element of $\emph{ran}( \mathit{J}( \mathbf{r}(\mathbf{a}) ) ) $ may be written uniquely as a sum of an element of $\emph{ran}( \mathbf{M}(\mathbf{a}) )$ and an element of $\emph{ran}( \mathbf{L}(\mathbf{a}) )$ and it follows that
\begin{equation*}
   \emph{ran} \left( \mathit{J} \big( \mathbf{r}(\mathbf{a}) \big) \right) \subset \emph{ran}\big( \mathbf{M}(\mathbf{a}) \big)  \oplus \emph{ran}\big( \mathbf{L}(\mathbf{a}) \big)
\end{equation*}
where $ \oplus $ stands for the direct sum. From these results, it is then easy to show that
\begin{equation*}
   \emph{null}\left( \mathit{J}\big( \mathbf{r}(\mathbf{a}) \big) \right) = \emph{null}\big( \mathbf{M}(\mathbf{a}) \big)  \cap \emph{null}\big( \mathbf{L}(\mathbf{a}) \big) \ .
\end{equation*}
Since $ \mathit{J}( \mathbf{r}(\mathbf{a}) ) = - ( \mathbf{M}(\mathbf{a}) + \mathbf{L}(\mathbf{a}) ) $, we have, by definition,
\begin{equation*}
   \emph{null}( \mathbf{M} )  \cap \emph{null}( \mathbf{L} ) \subset \emph{null} \left( \mathit{J}\big( \mathbf{r}(\mathbf{a}) \big) \right) \ ,
\end{equation*}
and, reciprocally,
\begin{align*}
     \mathbf{c}  \in \emph{null} \left( \mathit{J}\big( \mathbf{r}(\mathbf{a}) \big) \right) & \Rightarrow   \mathbf{M}(\mathbf{a})\mathbf{c} + \mathbf{L}(\mathbf{a})\mathbf{c} = \mathbf{0}^{p.n}  \\
                                                                                                          & \Rightarrow   \mathbf{M}(\mathbf{a})\mathbf{c} = - \mathbf{L}(\mathbf{a})\mathbf{c}       \\
                                                                                                          & \Rightarrow   \mathbf{M}(\mathbf{a})\mathbf{c} \in \emph{ran}( \mathbf{M}(\mathbf{a}) )  \cap \emph{ran}( \mathbf{L}(\mathbf{a}) )  \\
                                                                                                          & \Rightarrow   \mathbf{M}(\mathbf{a})\mathbf{c} = \mathbf{L}(\mathbf{a})\mathbf{c} = \mathbf{0}^{p.n}  \\
                                                                                                          & \Rightarrow   \mathbf{c} \in \emph{null}\big( \mathbf{M}(\mathbf{a}) \big)  \cap \emph{null}\big( \mathbf{L}(\mathbf{a}) \big) \ .
\end{align*}
Now, we demonstrate that the matrices $\mathbf{M}(\mathbf{a})$, $\mathbf{L}(\mathbf{a})$ and $\mathit{J}( \mathbf{r}(\mathbf{a}) )$ are rank deficient $\forall \thickspace \mathbf{a}  \in \mathbb{R}^{p.k}$. This result for $\mathbf{M}(\mathbf{a})$ was first noted by Ruhe~\cite{R1974} for the case  $k=1$ and $\mathbf{W}_{ij} \in \{0,1\}$. It was proved later for general $k$, again only for $\mathbf{M}(\mathbf{a})$ and $\mathbf{W}_{ij} \in \{0,1\}$, by Okatani and Deguchi~\cite{OD2007}, but under the restrictive hypotheses that $\mathbf{A}$, $\mathbf{B}$, $\mathbf{F}(\mathbf{a})$ and $\mathbf{G}(\mathbf{b})$ are of full rank. See also Okatani et al.~\cite{OYD2011}, where these results are further developed. The next theorem and corollary extend this result for general $k$ and any nonnegative real matrix $\mathbf{W}$ and to $\mathbf{M}(\mathbf{a})$, $\mathbf{L}(\mathbf{a})$ and $\mathit{J}( \mathbf{r}(\mathbf{a}) )$ matrices without any restrictive assumptions.
\\
\begin{theo5.2} \label{theo5.2:box}
Let  $k_{\mathbf{A}} = \emph{rank}( \mathbf{A} )$ . If,
\begin{align*}
    \mathbf{M}(\mathbf{a})                    & =   \mathbf{P}_{\mathbf{F}(\mathbf{a})}^{\bot}   \mathbf{U} (\mathbf{a}) = \mathbf{P}_{\mathbf{F}(\mathbf{a})}^{\bot} \emph{diag} \left(\emph{vec} \big( \sqrt{\mathbf{W}} \big) \right) ( \mathbf{\widehat{B}}^{T} \otimes \mathbf{I}_{p} ) \mathbf{K}_{(k,p)}  \ , \\
     \mathbf{L}(\mathbf{a})                    & =   \big( \mathbf{F}(\mathbf{a})^{+} \big)^{T}  \mathbf{V} (\mathbf{a}) = \big( \mathbf{F}(\mathbf{a})^{+} \big)^{T} \big( ( \mathbf{W} \odot P_{\Omega}(\mathbf{X} -\mathbf{A}\mathbf{\widehat{B}}) )^{T} \otimes \mathbf{I}_{k} \big)  \ , \\
     \mathit{J}( \mathbf{r}(\mathbf{a}) ) & =   - \big( \mathbf{M}(\mathbf{a}) + \mathbf{L}(\mathbf{a}) \big)  \ ,
\end{align*}
where all the matrices and vectors have the same definitions as above, then the following relationships hold
\begin{align*}
    \emph{dim}\left( \emph{null}\big( \mathbf{M}(\mathbf{a}) \big) \right)                   & \geqslant  k_{\mathbf{A}}.k    \ ,\\
    \emph{dim}\left( \emph{null}\big( \mathbf{L}(\mathbf{a}) \big) \right)                    & \geqslant  k_{\mathbf{A}}.k    \ , \\
    \emph{dim}\left( \emph{null}\big( \mathit{J}( \mathbf{r}(\mathbf{a}) ) \big) \right) & \geqslant  k_{\mathbf{A}}.k  \ .
\end{align*}
\end{theo5.2}
 \begin{proof}
Consider first the matrix $\mathbf{N}$ defined by
\begin{equation*}
    \mathbf{N} = \mathbf{K}_{(p,k)} ( \mathbf{I}_{k} \otimes \mathbf{A} )  \ .
\end{equation*}
Since $ \mathbf{A} $ is of rank $ k_{\mathbf{A}} $, $ \mathbf{I}_{k} \otimes \mathbf{A} $ is of rank $ k.k_{\mathbf{A}} $ (see equation~\eqref{eq:rank_kronprod}) and $ \mathbf{N} $ is also of rank $ k.k_{\mathbf{A}} $ because $ \mathbf{K}_{(p,k)} $ is a permutation matrix and the rank of a matrix is unaltered by multiplication with a nonsingular square matrix.

Now, we  first demonstrate that the space spanned by the columns of $ \mathbf{N} $, which is of dimension $ k.k_{\mathbf{A}} $, is included in $ \emph{null}( \mathbf{M}(\mathbf{a}) ) $ and so $ \emph{dim}( \emph{null}( \mathbf{M}(\mathbf{a}) ) ) \geqslant k_{\mathbf{A}}.k $.

Let $\mathbf{t} \in \emph{ran}( \mathbf{N} )$, then $  \exists \mathbf{Z} \in \mathbb{R}^{k \times k} $ such that
\begin{equation*}
    \mathbf{t} = \mathbf{N} \emph{vec}( \mathbf{Z} ) = \mathbf{K}_{(p,k)} ( \mathbf{I}_{k} \otimes \mathbf{A} ) \emph{vec}( \mathbf{Z} ) = \mathbf{K}_{(p,k)} \emph{vec}( \mathbf{A} \mathbf{Z} ) \ .
\end{equation*}
From this equality, we deduce that
\begin{align*}
    \mathbf{M}(\mathbf{a}) \mathbf{t} & =  \mathbf{P}_{\mathbf{F}(\mathbf{a})}^{\bot} \emph{diag}( \emph{vec}( \sqrt{\mathbf{W}} ) ) ( \mathbf{\widehat{B}}^{T} \otimes \mathbf{I}_{p} ) \mathbf{K}_{(k,p)}\mathbf{K}_{(p,k)} \emph{vec}( \mathbf{A} \mathbf{Z} )  \\
                                        & =  \mathbf{P}_{\mathbf{F}(\mathbf{a})}^{\bot} \emph{diag}( \emph{vec}( \sqrt{\mathbf{W}} ) ) ( \mathbf{\widehat{B}}^{T} \otimes \mathbf{I}_{p} )  \emph{vec}( \mathbf{A} \mathbf{Z} )  \\
                                        & =  \mathbf{P}_{\mathbf{F}(\mathbf{a})}^{\bot} \emph{diag}( \emph{vec}( \sqrt{\mathbf{W}} ) ) \emph{vec}( \mathbf{A} \mathbf{Z} \mathbf{\widehat{B}} )  \\
                                        & =  \mathbf{P}_{\mathbf{F}(\mathbf{a})}^{\bot} \emph{diag}( \emph{vec}( \sqrt{\mathbf{W}} ) ) ( \mathbf{I}_{n} \otimes \mathbf{A} ) \emph{vec}(  \mathbf{Z} \mathbf{\widehat{B}} )  \\
                                        & =  \mathbf{P}_{\mathbf{F}(\mathbf{a})}^{\bot} \mathbf{F}(\mathbf{a})  \emph{vec}(  \mathbf{Z} \mathbf{\widehat{B}} ) \\
                                        & =  \mathbf{0}^{p.n} \ ,
\end{align*}
since $ \mathbf{P}_{\mathbf{F}(\mathbf{a})}^{\bot} $ is the orthogonal projector onto $ \emph{ran}( \mathbf{F}(\mathbf{a}) )^{\bot} $ and, finally, $ \mathbf{t} \in \emph{null}( \mathbf{M}(\mathbf{a}) )$. In other words, $ \emph{ran}( \mathbf{N} ) \subset \emph{null}( \mathbf{M}(\mathbf{a}) ) $ and, hence, $ \emph{dim}( \emph{null}( \mathbf{M}(\mathbf{a}) ) ) \geqslant \emph{dim}( \emph{ran}( \mathbf{N} ) ) = k_{\mathbf{A}}.k  $ .

We now demonstrate that the relation $ \emph{ran}( \mathbf{N} ) \subset \emph{null}( \mathbf{L}(\mathbf{a}) ) $ also holds. If $\mathbf{t} \in \emph{ran}( \mathbf{N} )$ and $\mathbf{Z} \in \mathbb{R}^{k \times k}$ with $\mathbf{t} = \mathbf{N} \emph{vec}( \mathbf{Z} )$, using equation~\eqref{eq:vec_kronprod} and Lemma~\ref{theo2.2:box}, we have
\begin{align*}
    \mathbf{L}(\mathbf{a}) \mathbf{t} & =  ( \mathbf{F}(\mathbf{a})^{+} )^{T} ( ( \mathbf{W} \odot P_{\Omega}(\mathbf{X} -\mathbf{A}\mathbf{\widehat{B}}) )^{T} \otimes \mathbf{I}_{k} ) \mathbf{K}_{(p,k)} \emph{vec}( \mathbf{A} \mathbf{Z} ) \\
                                                         & =  ( \mathbf{F}(\mathbf{a})^{+} )^{T} ( ( \mathbf{W} \odot P_{\Omega}(\mathbf{X} -\mathbf{A}\mathbf{\widehat{B}}) )^{T} \otimes \mathbf{I}_{k} )  \emph{vec}( \mathbf{Z}^{T} \mathbf{A}^{T}  ) \\
                                                         & =  ( \mathbf{F}(\mathbf{a})^{+} )^{T} \emph{vec}( \mathbf{Z}^{T} \mathbf{A}^{T} ( \mathbf{W} \odot P_{\Omega}(\mathbf{X} -\mathbf{A}\mathbf{\widehat{B}}) ) ) \\
                                                         & =  ( \mathbf{F}(\mathbf{a})^{+} )^{T} ( \mathbf{I}_{n}  \otimes \mathbf{Z}^{T}) \emph{vec}( \mathbf{A}^{T} ( \mathbf{W} \odot P_{\Omega}(\mathbf{X} -\mathbf{A}\mathbf{\widehat{B}}) ) ) \\
                                                         & =  ( \mathbf{F}(\mathbf{a})^{+} )^{T} ( \mathbf{I}_{n}  \otimes \mathbf{Z}^{T}) ( \mathbf{I}_{n}  \otimes \mathbf{A}^{T}) \emph{vec}( \mathbf{W} \odot P_{\Omega}(\mathbf{X} -\mathbf{A}\mathbf{\widehat{B}}) ) \\
                                                         & =  ( \mathbf{F}(\mathbf{a})^{+} )^{T} ( \mathbf{I}_{n}  \otimes \mathbf{Z}^{T}) ( \mathbf{I}_{n}  \otimes \mathbf{A}^{T}) \emph{diag}( \emph{vec}( \sqrt{\mathbf{W}} ) ) \emph{vec}( \sqrt{\mathbf{W}} \odot P_{\Omega}(\mathbf{X} -\mathbf{A}\mathbf{\widehat{B}}) ) \\
                                                         & =  ( \mathbf{F}(\mathbf{a})^{+} )^{T} ( \mathbf{I}_{n}  \otimes \mathbf{Z}^{T}) \mathbf{F}(\mathbf{a})^{T} \mathbf{r}(\mathbf{a}) \\
                                                         & =  \mathbf{0}^{p.n} \ ,
\end{align*}
since $ \mathbf{F}(\mathbf{a})^{T} \mathbf{r}(\mathbf{a}) = \mathbf{0}^{k.n} $. Thus,  $ \emph{ran}( \mathbf{N} ) \subset \emph{null}( \mathbf{L}(\mathbf{a}) ) $ and, hence, 
 \begin{equation*}
    \emph{dim}\left( \emph{null}\big( \mathbf{L}(\mathbf{a}) \big) \right) \geqslant \emph{dim}\big( \emph{ran}( \mathbf{N} ) \big) = k_{\mathbf{A}}.k \ .
\end{equation*}
Finally, we have $ \emph{ran}( \mathbf{N} ) \subset \emph{null}( \mathbf{M}(\mathbf{a}) ) \cap  \emph{null}( \mathbf{L}(\mathbf{a}) ) = \emph{null}( \mathit{J}( \mathbf{r}(\mathbf{a}) ) ) $ and, so, 
\begin{equation*}
    \emph{dim}\left( \emph{null}\big( \mathit{J}( \mathbf{r}(\mathbf{a}) ) \big) \right) \geqslant k_{\mathbf{A}}.k \  .
\end{equation*}
\\
\end{proof}

\begin{corol5.3} \label{corol5.3:box}
With the same notations as in Theorem~\ref{theo5.2:box}, we have
\begin{align*}
    \emph{rank}\big( \mathbf{M}(\mathbf{a}) \big)                      & \leqslant  ( p - k_{\mathbf{A}}).k   \\
    \emph{rank}\big( \mathbf{L}(\mathbf{a}) \big)                       & \leqslant  ( p - k_{\mathbf{A}}).k   \\
    \emph{rank}\big( \mathit{J}( \mathbf{r}(\mathbf{a}) ) \big)    & \leqslant  ( p - k_{\mathbf{A}}).k \ .
\end{align*}
\end{corol5.3}
\begin{proof}
These inequalities follow directly from Theorem~\ref{theo5.2:box} and the rank-nullity theorem (see equation~\eqref{eq:rank} in Subsection~\ref{lin_alg:box}):
\begin{align*}
    \emph{rank}\big( \mathbf{M}(\mathbf{a}) \big)                         & =  k.p - \emph{dim} \Big( \emph{null}\big( \mathbf{M}(\mathbf{a}) \big)  \Big) \\
    \emph{rank}\big( \mathbf{L}(\mathbf{a}) \big)                          & =  k.p - \emph{dim} \Big( \emph{null}\big( \mathbf{L}(\mathbf{a}) \big)  \Big)   \\
    \emph{rank}\big( \mathit{J}( \mathbf{r}(\mathbf{a}) ) \big)       & =  k.p - \emph{dim} \Big( \emph{null}\big( \mathit{J}( \mathbf{r}(\mathbf{a}) ) \big)  \Big) \ . 
\end{align*}
\\
\end{proof}

\begin{remark5.1} \label{remark5.1:box}
Theorem~\ref{theo5.2:box} and Corollary~\ref{corol5.3:box} are obviously connected to the fact that the matrix factorization $\mathbf{Y}=\mathbf{A}\mathbf{B}$ used in the~\eqref{eq:P1} formulation of the WLRA problem is overparameterized, that the minimization of  $\psi(.)$ is an optimization problem on $\text{Gr}(p,k)$, the set of  linear subspaces of fixed dimension $k$ of the Euclidean space $\mathbb{R}^{p}$, and that $\text{Gr}(p,k)$ is a smooth manifold of dimension $p.k - k.k$  (see Remarks~\ref{remark3.4:box} and~\ref{remark3.7:box} for details). Assuming that $\mathbf{A} \in  \mathbb{R}^{p \times k}_k$, the minimization of the cost function $\psi(.)$ is at first sight a $k.p$ dimensional problem. However, $\psi(\mathbf{a})$ depends only on the column space of $\mathbf{A}$ and not on its individual elements~\cite{EAS1998}\cite{MMH2003}\cite{C2008b}\cite{BA2015}. As an illustration, if  $\mathbf{A} \in  \mathbb{R}^{p \times k}_k$ and  $d\mathbf{A} \in  \mathbb{R}^{p \times k}$ is a perturbation matrix, for certain matrices $d\mathbf{A}$, we will have $\emph{ran}( \mathbf{A} + d\mathbf{A} ) = \emph{ran}( \mathbf{A} )$, which will imply that $\psi(\mathbf{A} + d\mathbf{A}) = \psi(\mathbf{A})$. This demonstrates that it is not useful to consider all $k.p$ search directions for minimizing $\psi(.)$ from a previous matrix estimate $\mathbf{A}$. As demonstrated in~\cite{EAS1998}\cite{MMH2003}\cite{C2008b}\cite{BA2015}, this symmetry can be exploited to reduce the dimension of the problem to $k.(p-k)$ parameters instead of $k.p$ in both the~\eqref{eq:VP1} and~\eqref{eq:VP2} formulations of the WLRA problem. Thus, in that sense, the column space of $\mathbf{A}$ has only $p.k - k.k$ degrees of freedom, which is consistent to the facts that the rank of $\mathit{J}( \mathbf{r}(\mathbf{a}) )$ is at most $p.k - k.k$ if $\emph{rank}( \mathbf{A} ) = k$ and that the the dimension of $\text{Gr}(p,k)$ is also $p.k - k.k$. $\blacksquare$
\end{remark5.1}

Theorem~\ref{theo5.2:box} demonstrates that the Jacobian matrix  $\mathit{J}( \mathbf{r}(\mathbf{a}) )$ is always rank-deficient. This implies that the linear least-squares problem
\begin{equation*}
\min_{d\mathbf{a} \in \mathbb{R}^{p.k}}   \,   \frac{1}{2} \Vert \mathbf{r}(  \mathbf{a} ) + \mathit{J}( \mathbf{r}(\mathbf{a}) ) d\mathbf{a} \Vert^{2}_{2} \ ,
\end{equation*}
which must be solved at each iteration of a Gauss-Newton type algorithm (see Subsection~\ref{opt:box} for details) has an infinite set of solutions~\cite{GVL1996}\cite{HPS2012}\cite{B2015} and we must remove this ambiguity in any practical implementation of the Gauss-Newton algorithm in a such way that the direction vector $d\mathbf{a}_{gn}$ can be determined uniquely at each iteration. The general solution $d\mathbf{\widehat{a}} \in \mathbb{R}^{k.p}$ of the above rank-deficient linear least-squares problem can be written as
\begin{equation*}
d\mathbf{\widehat{a}} = - \mathit{J}( \mathbf{r}(\mathbf{a}) )^{+} \mathbf{r}(  \mathbf{a} ) + \mathbf{c} =  d\mathbf{a}_{min} + \mathbf{c} \ ,
\end{equation*}
where, as before, $\mathit{J}( \mathbf{r}(\mathbf{a}) )^{+}$ is the pseudo-inverse of $\mathit{J}( \mathbf{r}(\mathbf{a}) )$, $d\mathbf{a}_{min}$ is the (unique) minimum 2-norm solution of the above linear least-squares problem (see Subsection~\ref{lin_alg:box} and~\cite{GVL1996}\cite{HPS2012}\cite{B2015}), and $\mathbf{c}$ is an arbitrary $k.p$ dimensional vector belonging to $\emph{null}( \mathit{J}( \mathbf{r}(\mathbf{a}) ) )$.

First, the pseudo-inverse solution $d\mathbf{a}_{min}$ is characterized uniquely by the two conditions
\begin{equation*}
\mathit{J}( \mathbf{r}(\mathbf{a}) )^{T} \mathit{J}( \mathbf{r}(\mathbf{a}) ) d\mathbf{a}_{min} = - \mathit{J}( \mathbf{r}(\mathbf{a}) )^{T} \mathbf{r}(  \mathbf{a} )  \text{ and } d\mathbf{a}_{min} \in \emph{null}( \mathit{J}( \mathbf{r}(\mathbf{a}) ) )^{\bot} \ .
\end{equation*}
The first condition states simply that $d\mathbf{a}_{min}$ is a solution of the normal equations of the linear-least-squares problem. Note, further, that
\begin{equation*}
d\mathbf{a}_{min} = - \mathit{J}( \mathbf{r}(\mathbf{a}) )^{+} \mathbf{r}(  \mathbf{a} ) = - \mathit{J}( \mathbf{r}(\mathbf{a}) )^{+} \mathit{J}( \mathbf{r}(\mathbf{a}) ) \mathit{J}( \mathbf{r}(\mathbf{a}) )^{+} \mathbf{r}(  \mathbf{a} ) = \mathbf{P}_{\mathit{J}( \mathbf{r}(\mathbf{a}) )^{T}} d\mathbf{a}_{min} \ ,
\end{equation*}
where $\mathbf{P}_{\mathit{J}( \mathbf{r}(\mathbf{a}) )^{T}}$ is the orthogonal projector onto the row space of $\mathit{J}( \mathbf{r}(\mathbf{a}) )$ (e.g., $\emph{ran}(\mathit{J}( \mathbf{r}(\mathbf{a}) )^{T}$), see Subsection~\ref{lin_alg:box} for details. Since $\emph{ran}(\mathit{J}( \mathbf{r}(\mathbf{a}) )^{T} ) = \emph{null}( \mathit{J}( \mathbf{r}(\mathbf{a}) ) )^\bot$, we deduce immediately that $d\mathbf{a}_{min} \in  \emph{null}( \mathit{J}( \mathbf{r}(\mathbf{a}) ) )^\bot$ as stated in the second condition.

Now, if $d\mathbf{\widehat{a}} = d\mathbf{a}_{min} + \mathbf{c}$, we have
\begin{equation*}
\mathit{J}( \mathbf{r}(\mathbf{a}) ) d\mathbf{\widehat{a}} =  \mathit{J}( \mathbf{r}(\mathbf{a}) ) ( d\mathbf{a}_{min} + \mathbf{c} ) = \mathit{J}( \mathbf{r}(\mathbf{a}) ) d\mathbf{a}_{min} \ ,
\end{equation*}
and, thus, $\Vert \mathbf{r}(  \mathbf{a} ) + \mathit{J}( \mathbf{r}(\mathbf{a}) ) d\mathbf{\widehat{a}} \Vert_{2} = \Vert   \mathbf{r}(  \mathbf{a} ) + \mathit{J}( \mathbf{r}(\mathbf{a}) ) d\mathbf{a}_{min} \Vert_{2}$, which implies that $d\mathbf{\widehat{a}}$ is also a solution of the above linear least-squares problem. Reciprocally, if $d\mathbf{\widehat{a}}$ is an arbitrary solution, we have $\Vert \mathbf{r}(  \mathbf{a} ) + \mathit{J}( \mathbf{r}(\mathbf{a}) ) d\mathbf{\widehat{a}} \Vert_{2} = \Vert   \mathbf{r}(  \mathbf{a} ) + \mathit{J}( \mathbf{r}(\mathbf{a}) ) d\mathbf{a}_{min} \Vert_{2}$, which implies that
\begin{equation*}
\mathit{J}( \mathbf{r}(\mathbf{a}) ) d\mathbf{\widehat{a}} =  - \mathbf{P}_{\mathit{J}( \mathbf{r}(\mathbf{a}) )} \mathbf{r}(\mathbf{a}) = \mathit{J}( \mathbf{r}(\mathbf{a}) ) d\mathbf{a}_{min} \ ,
\end{equation*}
as $-\mathbf{P}_{\mathit{J}( \mathbf{r}(\mathbf{a}) )} \mathbf{r}(  \mathbf{a} )$ is the unique closest point to $\mathbf{r}(  \mathbf{a} )$ in $\emph{ran}( \mathit{J}( \mathbf{r}(\mathbf{a}) )  )$, see equation~\eqref{eq:projector2} of Subsection~\ref{lin_alg:box} for details. Thus, $\mathit{J}( \mathbf{r}(\mathbf{a}) ) ( d\mathbf{\widehat{a}} - d\mathbf{a}_{min} ) = \mathbf{0}^{p.n}$ and we can write $d\mathbf{\widehat{a}}$ as
\begin{equation*}
d\mathbf{\widehat{a}} = d\mathbf{a}_{min} + (  d\mathbf{\widehat{a}} - d\mathbf{a}_{min} )  = d\mathbf{a}_{min} + \mathbf{c} \ ,
\end{equation*}
with $\mathbf{c} = d\mathbf{\widehat{a}} - d\mathbf{a}_{min}   \in \emph{null}( \mathit{J}( \mathbf{r}(\mathbf{a}) ) )$.

In other words, all solution vectors $d\mathbf{\widehat{a}}$ can be written uniquely as the sum of $d\mathbf{a}_{min}  \in  \emph{null}( \mathit{J}( \mathbf{r}(\mathbf{a}) ) )^{\bot}$ and a vector $\mathbf{c} \in  \emph{null}( \mathit{J}( \mathbf{r}(\mathbf{a}) ) )$ and finding all the solutions of the above rank-deficient linear least-squares problem requires computing both a generalized inverse and a basis of the null space of $\mathit{J}( \mathbf{r}(\mathbf{a}) )$. Obviously, this also implies to determine accurately the rank of $\mathit{J}( \mathbf{r}(\mathbf{a}) )$ or, equivalently, the rank of its null space.
More generally, proceeding in a similar manner, it is also easy to establish an one to one mapping between the elements of  $\emph{null}( \mathit{J}( \mathbf{r}(\mathbf{a}) ) )^\bot$ and those of $\emph{ran}( \mathit{J}( \mathbf{r}(\mathbf{a}) ) )$.

Now, the most natural choice is to select  $d\mathbf{a}_{gn} = d\mathbf{a}_{min}$ as  the solution of our linear least-squares problem since, with such minimum Euclidean norm solution, the first order Taylor's expansion
\begin{equation*}
\mathbf{r}(  \mathbf{a} +  d\mathbf{a}_{gn} ) = \mathbf{r}(  \mathbf{a} ) + \mathit{J}( \mathbf{r}(\mathbf{a}) )d\mathbf{a}_{gn} +  \mathcal{O}( \Vert d\mathbf{a}_{gn} \Vert^{2}_{2} ) \ ,
\end{equation*}
which is at the base  of the  Gauss-Newton algorithm is the most accurate. Selecting $d\mathbf{a}_{gn} = d\mathbf{a}_{min}$ has also a strong theoretical justification as, with this choice, a variable projection Gauss-Newton algorithm used to minimize  $\psi(.)$  is equivalent to a Riemannian optimization method operating directly on the Grassmann manifold  $\text{Gr}(p,k)$~\cite{AMS2008}\cite{B2023} as we will explain later in this subsection.

These considerations related to the uniform rank degeneracy of the Jacobian matrix $\mathit{J}( \mathbf{r}(\mathbf{a}) )$ apply also to the computation of the correction vector $d\mathbf{a}_{lm}$ in the Levenberg-Marquardt method as soon as the Marquardt damping parameter $\lambda$ approaches zero, as it is expected after some iterations of the Levenberg-Marquardt algorithm. Moreover, if the Marquardt parameter $\lambda$ is controlled so that it does not approach to zero in order to remove the uniform singularity of the Jacobian matrix  $\mathit{J}( \mathbf{r}(\mathbf{a}) )$, this may severely deteriorate the  global convergence as well as  the local convergence of the method in a neighborhood of a critical point. In other words, adding the additional constraint that $\Vert d\mathbf{a}_{lm}\Vert_{2}$ is minimum when $\lambda$ approaches zero, is also important for the robustness and efficiency of the Levenberg-Marquardt or similar regularized methods described in Subsection~\ref{opt:box} when they are used to solve  NNLS problems with an uniformly deficient Jacobian matrix, like the WLRA problem.

We now give sufficient conditions for the equalities:
\begin{align*} 
    \emph{dim}\Big( \emph{null}\big( \mathbf{M}(\mathbf{a}) \big) \Big)                   & =  k_{\mathbf{A}}.k    \\
    \emph{dim}\Big( \emph{null}\big( \mathit{J}( \mathbf{r}(\mathbf{a}) ) \big) \Big) & =  k_{\mathbf{A}}.k    \ ,
\end{align*}
which will be helpful to remove these ambiguities in determining uniquely and efficiently the correction vectors $d\mathbf{a}_{gn}$ and $d\mathbf{a}_{lm}$ in many practical applications.

Let us first introduce some definitions and notations. For any nonnegative real $ p \times n $ matrix $ \mathbf{W} $, we define the finite subset of $\mathbb{N}$
 \begin{equation*}
    \Theta ( \mathbf{W} ) = \big\{ j \in  \{1,2, \cdots ,n\}  \text{ } / \text{ } \mathbf{W}_{.j} \in  \mathbb{R}^{p}_{+*}  \big\} \ .
\end{equation*}
$  \Theta ( \mathbf{W} ) $ is the set of the column-vector indices of $ \mathbf{W} $ such that 0 is not an element of such column-vector of $\mathbf{W}$. Furthermore, let $  card( \Theta ( \mathbf{W} ) )$ be the number of elements of $  \Theta ( \mathbf{W} ) $ and, for any $ s \times n $ matrix $ \mathbf{C} $, define the $ s \times card( \Theta ( \mathbf{W} ) ) $ real submatrix $ \mathbf{C}^{'} $ obtained from $ \mathbf{C} $ by deleting the columns of $ \mathbf{C} $ whose indices do not belong to $  \Theta ( \mathbf{W} ) $. We then have the following result, which is new as far as we know.
\begin{theo5.3} \label{theo5.3:box}
With these definitions and the same notations as in Theorem~\ref{theo5.2:box}, if $\emph{card}( \Theta ( \mathbf{W} ) ) = n^{'}  \geqslant k$ and $\emph{rank}( \mathbf{\widehat{B}}^{'} )=k$ then the following equalities hold:
 \begin{align*}
    \emph{null}\big( \mathit{J}( \mathbf{r}(\mathbf{a}) ) \big)                            & =  \emph{null}\big( \mathbf{M}(\mathbf{a}) \big)   \ ,  \\
    \emph{dim}\big( \emph{null}( \mathit{J}( \mathbf{r}(\mathbf{a}) ) ) \big)     & =   k_{\mathbf{A}}.k    \ .
\end{align*}
\end{theo5.3}
\begin{proof}
First, consider the second formulation of the $\mathbf{M}(\mathbf{a})$ matrix (see equation~\eqref{eq:M_mat2}), e.g.,
\begin{equation*}
\mathbf{M}(\mathbf{a}) = \mathbf{P}_{\mathbf{F}(\mathbf{a})}^{\bot} \mathbf{K}_{(n,p)} \mathbf{G}(\mathbf{\widehat{b}})  \ ,
\end{equation*}
where $\mathbf{\widehat{b}}= \emph{vec}( \mathbf{\widehat{B}} )$ and $\mathbf{G}(\mathbf{\widehat{b}}) = \emph{diag}( \emph{vec}( \sqrt{\mathbf{W}}^{T} ) ) ( \mathbf{I}_{p} \otimes \mathbf{\widehat{B}}^{T} )$. Using the two hypotheses $\emph{card}( \Theta ( \mathbf{W} ) ) = n^{'}  \geqslant k$ and $\emph{rank}( \mathbf{\widehat{B}}^{'} )=k$, we first deduce that
\begin{align*}
\emph{rank} \Big( \emph{diag}\big( \emph{vec}( \sqrt{\mathbf{W}^{'}}^{T} ) \big) ( \mathbf{I}_{p} \otimes \mathbf{\widehat{B}}^{'T} ) \Big ) & =  \emph{rank} \big(  \mathbf{I}_{p} \otimes \mathbf{\widehat{B}}^{'T}  \big)   \ ,    \\
                                                                        & =  \emph{rank}(  \mathbf{I}_{p} ) . \emph{rank}( \mathbf{\widehat{B}}^{'T} )  \\
                                                                        & =  p.k   \ ,
\end{align*}
since $ \emph{diag}( \emph{vec}( \sqrt{\mathbf{W}^{'}}^{T} ) ) $ is a nonsingular diagonal matrix.
Now, using this equality, we have also
\begin{equation*}
\emph{rank}\big( \mathbf{K}_{(n,p)} \mathbf{G}(\mathbf{\widehat{b}}) \big) =  \emph{rank}( \mathbf{G}(\mathbf{\widehat{b}}) ) = k.p \ ,
\end{equation*}
as $\mathbf{K}_{(n,p)}$ is a (nonsingular) permutation matrix and $\emph{diag}( \emph{vec}( \sqrt{\mathbf{W}^{'}}^{T} ) ) ( \mathbf{I}_{p} \otimes \mathbf{\widehat{B}}^{'T} )$ is a submatrix of $\mathbf{G}(\mathbf{\widehat{b}})$ formed simply by eliminating some rows of $\mathbf{G}(\mathbf{\widehat{b}})$.

Now, for any matrix $ \mathbf{C} $ with $s$ columns, we have the basic rank-nullity relation (see equation~\eqref{eq:rank})
\begin{equation*}
    s = \emph{rank}( \mathbf{C} ) + \emph{dim} \big( \emph{null}( \mathbf{C} ) \big) \ .
\end{equation*}
Furthermore, for any matrix $ \mathbf{D} $ with $s$ rows, we also assume the equality
\begin{equation*}
    \emph{rank}( \mathbf{D} ) = \emph{rank}( \mathbf{C}\mathbf{D} ) + \emph{dim}\big( \emph{null}( \mathbf{C} ) \cap \emph{ran}( \mathbf{D} ) \big) \ ,
\end{equation*}
see Marsaglia and Styan~\cite{MS1974} for a proof.

Using these two relations, we deduce
\begin{equation*}
    \emph{rank}  \big( \mathbf{M}(\mathbf{a})  \big) + \emph{dim} \big( \emph{null}( \mathbf{M}(\mathbf{a}) )  \big) = k.p
\end{equation*}
and
\begin{equation*}
    k.p = \emph{rank} \big( \mathbf{M}(\mathbf{a})  \big) + \emph{dim} \big( \emph{null}( \mathbf{P}_{\mathbf{F}(\mathbf{a})}^{\bot} ) \cap \emph{ran}( \mathbf{K}_{(n,p)} \mathbf{G}(\mathbf{\widehat{b}}) )  \big) \ ,
\end{equation*}
and so
 \begin{align*}
    \emph{dim} \big( \emph{null}( \mathbf{M}(\mathbf{a}) )  \big) & = \emph{dim} \big( \emph{null}( \mathbf{P}_{\mathbf{F}(\mathbf{a})}^{\bot}  \big) \cap \emph{ran} \big( \mathbf{K}_{(n,p)} \mathbf{G}(\mathbf{\widehat{b}}) )  \big)  \\
                                                                                                   & =  \emph{dim} \big( \emph{ran}( \mathbf{F}(\mathbf{a})  \big) \cap \emph{ran} \big( \mathbf{K}_{(n,p)} \mathbf{G}(\mathbf{\widehat{b}}) )  \big) \ .
\end{align*}
Next, we consider the matrix $ \mathbf{H} $ defined by
\begin{equation*}
    \mathbf{H} =  \emph{diag}( \emph{vec}( \sqrt{\mathbf{W}} ) ) ( \mathbf{\widehat{B}}^{T}  \otimes  \mathbf{A} ) \ .
\end{equation*}
We have
\begin{equation*}
    \emph{rank}( \mathbf{H} ) \leqslant  \emph{rank}( \mathbf{\widehat{B}}^{T}  \otimes  \mathbf{A} ) = \emph{rank}( \mathbf{\widehat{B}} ).\emph{rank}( \mathbf{A} ) = k.k_{\mathbf{A}} \ ,
\end{equation*}
since the hypothesis $ \emph{rank}( \mathbf{\widehat{B}}^{'} )=k $ implies $  \emph{rank}( \mathbf{\widehat{B}} )=k $. We now demonstrate the inclusion
\begin{equation*}
    \emph{ran}( \mathbf{F}(\mathbf{a} ) ) \cap \emph{ran}( \mathbf{K}_{(n,p)} \mathbf{G}(\mathbf{\widehat{b}}) )  \subset  \emph{ran}( \mathbf{H} ) \ .
\end{equation*}
Let $ \mathbf{c} \in \emph{ran}( \mathbf{F}(\mathbf{a} ) ) \cap \emph{ran}( \mathbf{K}_{(n,p)} \mathbf{G}(\mathbf{\widehat{b}}) ) $, then $ \exists \; \mathbf{S} \in \mathbb{R}^{k \times n} $ and $ \mathbf{Z} \in \mathbb{R}^{p \times k}  $ such that
\begin{equation*}
   \mathbf{c} = \mathbf{F}(\mathbf{a} ) \emph{vec}( \mathbf{S} ) = \mathbf{K}_{(n,p)} \mathbf{G}(\mathbf{\widehat{b}}) \emph{vec}( \mathbf{Z}^{T} )
\end{equation*}
and we want to show that $ \exists  \mathbf{T} \in \mathbb{R}^{k \times k} $ such that $ \mathbf{c} =  \mathbf{H} \emph{vec}( \mathbf{T} ) $. But, $ \forall \mathbf{T} \in \mathbb{R}^{k \times k} $, we have
\begin{align*}
     \mathbf{H} \emph{vec}( \mathbf{T} ) & =   \emph{diag}( \emph{vec}( \sqrt{\mathbf{W}} ) ) ( \mathbf{\widehat{B}}^{T}  \otimes  \mathbf{A} ) \emph{vec}( \mathbf{T} ) \\
                                                                & =   \emph{diag}( \emph{vec}( \sqrt{\mathbf{W}} ) ) ( \mathbf{\widehat{B}}^{T}  \otimes  \mathbf{I}_{p} ) ( \mathbf{I}_{k}  \otimes  \mathbf{A} ) \emph{vec}( \mathbf{T} )       \\
                                                                & =   \emph{diag}( \emph{vec}( \sqrt{\mathbf{W}} ) ) ( \mathbf{\widehat{B}}^{T}  \otimes  \mathbf{I}_{p} ) \mathbf{K}_{(k,p)} \mathbf{K}_{(p,k)} \emph{vec}( \mathbf{A} \mathbf{T} )       \\
                                                                & =   \mathbf{K}_{(n,p)} \mathbf{G}(\mathbf{\widehat{b}}) \mathbf{K}_{(p,k)}  \emph{vec}( \mathbf{A} \mathbf{T} )  \ ,
\end{align*}
and so, using the facts that $ \mathbf{K}_{(n,p)} $ is a (nonsingular) permutation matrix and $ \mathbf{G}(\mathbf{\widehat{b}}) $ has full column rank demonstrated above, we have the equivalences
\begin{align*}
 \mathbf{c} =  \mathbf{H} \emph{vec}( \mathbf{T} ) & \Leftrightarrow  \mathbf{K}_{(n,p)} \mathbf{G}(\mathbf{\widehat{b}}) \emph{vec}( \mathbf{Z}^{T} ) = \mathbf{K}_{(n,p)} \mathbf{G}(\mathbf{\widehat{b}}) \mathbf{K}_{(p,k)}  \emph{vec}( \mathbf{A} \mathbf{T} )  \\
                                                                              & \Leftrightarrow   \mathbf{G}(\mathbf{\widehat{b}}) \emph{vec}( \mathbf{Z}^{T} )  =  \mathbf{G}(\mathbf{\widehat{b}}) \mathbf{K}_{(p,k)}  \emph{vec}( \mathbf{A} \mathbf{T} )   \\
                                                                              & \Leftrightarrow   \emph{vec}( \mathbf{Z}^{T} ) = \mathbf{K}_{(p,k)}  \emph{vec}( \mathbf{A} \mathbf{T} ) \\
                                                                              & \Leftrightarrow   \mathbf{K}_{(k,p)} \emph{vec}( \mathbf{Z}^{T} ) = \emph{vec}( \mathbf{A} \mathbf{T} ) \\
                                                                              & \Leftrightarrow   \emph{vec}( \mathbf{Z} )= \emph{vec}( \mathbf{A} \mathbf{T}  )  \\
                                                                              & \Leftrightarrow   \mathbf{Z} = \mathbf{A} \mathbf{T} \ .
\end{align*}
Thus, to demonstrate that $ \exists  \mathbf{T} \in \mathbb{R}^{k \times k} $ such that $ \mathbf{c} =  \mathbf{H} \emph{vec}( \mathbf{T} ) $, it suffices to show that $ \exists  \mathbf{T} \in \mathbb{R}^{k \times k} $ such that $ \mathbf{Z} = \mathbf{A} \mathbf{T} $. But,
\begin{IEEEeqnarray*}{*x+lCl}
   & \mathbf{c}  & = &   \mathbf{F}(\mathbf{a}) \emph{vec}( \mathbf{S}  )  \\
   &                   & = &   \emph{diag}( \emph{vec}( \sqrt{\mathbf{W}} ) ) ( \mathbf{I}_{n}  \otimes  \mathbf{A} )  \emph{vec}( \mathbf{S}  )  \\
   &                   & = &   \emph{diag}( \emph{vec}( \sqrt{\mathbf{W}} ) )  \emph{vec}( \mathbf{A} \mathbf{S}  ) \\
    \text{ and also }  \\
   & \mathbf{c}  & = &   \mathbf{K}_{(n,p)} \mathbf{G}(\mathbf{\widehat{b}}) \emph{vec}( \mathbf{Z}^{T} )   \\
   &                   & = &   \mathbf{K}_{(n,p)} \emph{diag}( \emph{vec}( \sqrt{\mathbf{W}}^{T} ) ) ( \mathbf{I}_{p}  \otimes  \mathbf{\widehat{B}}^{T} )  \emph{vec}( \mathbf{Z}^{T}  )  \\
   &                   & = &   \emph{diag}( \emph{vec}( \sqrt{\mathbf{W}} ) ) ( \mathbf{I}_{n}  \otimes  \mathbf{Z} )  \emph{vec}( \mathbf{\widehat{B}}  )  \\
   &                   & = &   \emph{diag}( \emph{vec}( \sqrt{\mathbf{W}} ) )  \emph{vec}( \mathbf{Z} \mathbf{\widehat{B}}  ) \ .
\end{IEEEeqnarray*}
From these equalities, we then deduce that
\begin{equation*}
    \emph{diag}( \emph{vec}( \sqrt{\mathbf{W}^{'}} ) )  \emph{vec}( \mathbf{A} \mathbf{S}^{'}  ) = \emph{diag}( \emph{vec}( \sqrt{\mathbf{W}^{'}} ) )  \emph{vec}( \mathbf{Z} \mathbf{\widehat{B}}^{'}  ) \ ,
\end{equation*}
and since $ \emph{diag}( \emph{vec}( \sqrt{\mathbf{W}^{'}} ) ) $ is a nonsingular diagonal matrix, we obtain
\begin{equation*}
    \emph{vec}( \mathbf{A} \mathbf{S}^{'}  ) =   \emph{vec}( \mathbf{Z} \mathbf{\widehat{B}}^{'}  ) \text { and, finally, } \mathbf{A} \mathbf{S}^{'} =   \mathbf{Z} \mathbf{\widehat{B}}^{'} \ .
\end{equation*}
Now, $ \mathbf{\widehat{B}}^{'} $ has full row-rank by hypothesis and, consequently, admits a right inverse $ \mathbf{R} $ such that $ \mathbf{\widehat{B}}^{'} \mathbf{R} = \mathbf{I}_{k} $ (see~\cite{MS1974} for details) and so $ \mathbf{A} \mathbf{S}^{'} \mathbf{R} = \mathbf{Z} $ and we can take $ \mathbf{T} = \mathbf{S}^{'} \mathbf{R} $ . Consequently, we have proved
\begin{equation*}
    \emph{ran}( \mathbf{F}(\mathbf{a} ) ) \cap \emph{ran}( \mathbf{K}_{(n,p)} \mathbf{G}(\mathbf{\widehat{b}}) )  \subset  \emph{ran}( \mathbf{H} ) \ ,
\end{equation*}
which implies
\begin{equation*}
    \emph{dim}( \emph{null}( \mathbf{M}(\mathbf{a}) ) ) = \emph{dim}( \emph{ran}( \mathbf{F}(\mathbf{a}) ) \cap \emph{ran}( \mathbf{K}_{(n,p)} \mathbf{G}(\mathbf{\widehat{b}}) ) )  \leqslant  \emph{rank}( \mathbf{H} ) \leqslant k_{\mathbf{A}}.k \  .
\end{equation*}
But, from Theorem~\ref{theo5.2:box}, we already know that
\begin{equation*}
    k_{\mathbf{A}}.k   \leqslant   \emph{dim}( \emph{null}( \mathbf{M}(\mathbf{a}) ) )
\end{equation*}
and we obtain $ \emph{dim}( \emph{null}( \mathbf{M}(\mathbf{a}) ) ) =  k_{\mathbf{A}}.k  $.

Again, from Theorem~\ref{theo5.2:box}, we also know that
\begin{equation*}
    k_{\mathbf{A}}.k   \leqslant    \emph{dim}( \emph{null}( \mathit{J}( \mathbf{r}(\mathbf{a}) ) ) ) =  \emph{dim}( \emph{null}( \mathbf{M}(\mathbf{a}) ) \cap  \emph{null}( \mathbf{L}(\mathbf{a}) ) )
\end{equation*}
and since
\begin{equation*}
    \emph{dim}( \emph{null}( \mathbf{M}(\mathbf{a}) ) \cap  \emph{null}( \mathbf{L}(\mathbf{a}) ) ) \leqslant  \emph{dim}( \emph{null}( \mathbf{M}(\mathbf{a}) ) ) =  k_{\mathbf{A}}.k \ ,
\end{equation*}
we also conclude that $ \emph{dim}( \emph{null}( \mathit{J}( \mathbf{r}(\mathbf{a}) ) ) ) =  k_{\mathbf{A}}.k  $.

Finally, from the propositions,
\begin{equation*}
\emph{dim}( \emph{null}( \mathit{J}( \mathbf{r}(\mathbf{a}) ) ) ) = \emph{dim}( \emph{null}( \mathbf{M}(\mathbf{a}) ) ) \text{ and }  \emph{null}(  \mathit{J}( \mathbf{r}(\mathbf{a}) ) ) \subset \emph{null}( \mathbf{M}(\mathbf{a}) ) \ ,
\end{equation*}
 we deduce that $ \emph{null}( \mathit{J}( \mathbf{r}(\mathbf{a}) ) ) = \emph{null}( \mathbf{M}(\mathbf{a}) ) $ as claimed in the theorem.
\\
\end{proof}
Before stating some consequences of Theorem~\ref{theo5.3:box}, it is important to highlight that the two hypotheses of this theorem are not very stringent and are easily checked in practice. Moreover, in most practical cases, these two conditions will be meet as long as $\mathbf{W}$ has $k$ column-vectors without zero elements and each other column-vector of $\mathbf{W}$ has at least $k$ nonzero elements in it. The following corollary is then obvious and is stated without proof:
\\
\begin{corol5.4} \label{corol5.4:box}
With the same notations and hypotheses as in Theorem~\ref{theo5.3:box}, we also have the relations
\begin{align*}
    \emph{rank} \big( \mathbf{M}(\mathbf{a})  \big)                                               & =  ( p - k_{\mathbf{A}}).k  \ , \\
    \emph{rank} \big( \mathit{J}( \mathbf{r}(\mathbf{a}) )  \big)                             & =  ( p - k_{\mathbf{A}}).k \ .
\end{align*}
$\Box$
\end{corol5.4}

Furthermore, the next corollary shows that the sufficient conditions stated in Theorem~\ref{theo5.3:box} and used in Corollary~\ref{corol5.4:box} can be simplified when $\mathbf{W}$ is a strictly positive $ p \times n $  real matrix. A similar result has been obtained by Chen~\cite{C2008a} for $\mathbf{M}(\mathbf{a})$ only, but with a different proof and the additional hypothesis that $\mathbf{A}$ is of full column-rank.
\\
\begin{corol5.5} \label{corol5.5:box}
With the same notations as in Theorem~\ref{theo5.2:box}, if $ \mathbf{W} \in \mathbb{R}^{p \times n}_{+*} $ and $ \emph{rank}( \mathbf{\widehat{B}} ) = k $ then the following equalities hold:
\begin{align*}
    \emph{dim} \big( \emph{null}( \mathbf{M}(\mathbf{a}) ) \big)                          & =   k_{\mathbf{A}}.k   \ ,  \\
   \emph{dim} \big( \emph{null}( \mathit{J}( \mathbf{r}(\mathbf{a}) ) ) \big)         & =   k_{\mathbf{A}}.k  \ ,
\end{align*}
and also
\begin{align*}
     \emph{rank} \big( \mathbf{M}(\mathbf{a}) \big)                                               & =  ( p - k_{\mathbf{A}}).k   \ ,  \\
     \emph{rank} \big( \mathit{J}( \mathbf{r}(\mathbf{a}) ) \big)                             & =  ( p - k_{\mathbf{A}}).k \ .
\end{align*}
\end{corol5.5}
\begin{proof}
It suffices to note that $ \mathbf{W} \in \mathbb{R}^{p \times n}_{+*} $ and $ \emph{rank}( \mathbf{\widehat{B}} ) = k $ lead to
\begin{equation*}
    \emph{card}( \Theta ( \mathbf{W} ) ) = n  \geqslant k \text{ and } \mathbf{\widehat{B}}^{'} = \mathbf{\widehat{B}}
\end{equation*}
and the results follow immediately from Theorem~\ref{theo5.3:box}  and Corollary~\ref{corol5.4:box}.

\end{proof}
\begin{remark5.2} \label{remark5.2:box}
More generally, if $ \mathbf{W} \in \mathbb{R}^{p \times n}_{+*} $ and $ \emph{rank}( \mathbf{\widehat{B}} ) = k_{\mathbf{\widehat{B}}} < k $ it is possible to demonstrate
\begin{align*}
       \emph{dim} \big( \emph{null}( \mathbf{M}(\mathbf{a}) )  \big)              & =  p.( k - k_{\mathbf{\widehat{B}}}) + k_{\mathbf{A}}.k_{\mathbf{\widehat{B}}}    \ , \\
       \emph{rank} \big( \mathbf{M}(\mathbf{a})  \big)                                   & =  k_{\mathbf{\widehat{B}}}.( p - k_{\mathbf{A}})  \ ,
\end{align*}
and also
\begin{IEEEeqnarray*}{rCcCl}
     k_{\mathbf{A}}.k                                            &  \leqslant &  \emph{dim} \big( \emph{null}( \mathit{J}( \mathbf{r}(\mathbf{a}) ) )  \big)   & \leqslant &  p.( k - k _{\mathbf{\widehat{B}}}) + k_{\mathbf{A}}.k_{\mathbf{\widehat{B}}}   \ , \\
     k_{\mathbf{\widehat{B}}}.(p-k_{\mathbf{A}}) &  \leqslant &  \emph{rank} \big( \mathit{J}( \mathbf{r}(\mathbf{a}) )  \big)                        & \leqslant &  k.(p-k_{\mathbf{A}}) \ ,
\end{IEEEeqnarray*}
but we omit the details since these results are not very useful in practical applications in which $ \emph{rank}( \mathbf{A} ) = \emph{rank}( \mathbf{\widehat{B}} ) = k$ is the rule. $\blacksquare$
\\
\end{remark5.2}

If the hypotheses of Theorem~\ref{theo5.3:box} and Corollaries~\ref{corol5.4:box} and~\ref{corol5.5:box} are satisfied, we know precisely the dimensions of $\emph{null}( \mathit{J}( \mathbf{r}(\mathbf{a}) ) ) = \emph{null}( \mathbf{M}(\mathbf{a}) )$ and of its orthogonal complement in $\mathbb{R}^{k.p}$, e.g.,
\begin{equation*}
\emph{dim}  \big( \emph{null}( \mathit{J}( \mathbf{r}(\mathbf{a}) ) )  \big) = k_{\mathbf{A}}.k \text{ and } \emph{dim}  \big( \emph{null}( \mathit{J}( \mathbf{r}(\mathbf{a}) ) )^{\bot}  \big) = ( p - k_{\mathbf{A}}).k \ .
\end{equation*}
Furthermore, while Theorem~\ref{theo5.3:box}, Corollaries~\ref{corol5.4:box} and~\ref{corol5.5:box} are valid for any $\mathbf{A} \in  \mathbb{R}^{p \times k}$, we recall that the condition $\emph{rank}( \mathbf{A} ) = k$ is required for $\psi(.)$ to be differentiable, which implies that for all practical WLRA applications, we will have
\begin{equation*}
\emph{dim}  \big( \emph{null}( \mathit{J}( \mathbf{r}(\mathbf{a}) ) )  \big) = k.k \text{ and } \emph{dim}  \big( \emph{null}( \mathit{J}( \mathbf{r}(\mathbf{a}) ) )^{\bot}  \big) = ( p - k).k \ .
\end{equation*}

Then, if  (orthonormal or not) bases of $\emph{null}( \mathit{J}( \mathbf{r}(\mathbf{a}) ) ) = \emph{null}(  \mathbf{M}(\mathbf{a}) ) $ and of its orthogonal complement are available, it is easy to obtain the minimal 2-norm solution and also all solutions of the linear least-squares problem
\begin{equation*}
\min_{d\mathbf{a} \in \mathbb{R}^{p.k}}   \,   \frac{1}{2} \Vert \mathbf{r}(  \mathbf{a} ) + \mathit{J}  \big( \mathbf{r}(\mathbf{a})  \big) d\mathbf{a} \Vert^{2}_{2} \ .
\end{equation*}
More, precisely, if  the columns of $\mathbf{N} \in \mathbb{R}^{k.p \times k_{\mathbf{A}}.k }$ and $\mathbf{N}^{\bot} \in  \mathbb{R}^{k.p \times ( p - k_{\mathbf{A}}).k }$ form, respectively, (orthonormal) bases of $\emph{null}( \mathit{J}( \mathbf{r}(\mathbf{a}) ) )$ and $\emph{null}( \mathit{J}( \mathbf{r}(\mathbf{a}) ) )^{\bot}$, then $d\mathbf{a}_{gn} = d\mathbf{a}_{min}$ can be computed in a two-step procedure. In the first step, we need to solve the following reduced (and nonsingular) linear least-squares problem
\begin{equation*}
\min_{d\mathbf{a} \in \mathbb{R}^{( p - k_{\mathbf{A}}).k } }  \,   \frac{1}{2} \Vert \mathbf{r}(  \mathbf{a} ) + \mathit{J} \big( \mathbf{r}(\mathbf{a}) \big) \mathbf{N}^{\bot} d\mathbf{a} \Vert^{2}_{2} \ ,
\end{equation*}
which has an unique solution as demonstrated above, say $d\mathbf{\bar{a}}_{gn} \in \mathbb{R}^{( p - k_{\mathbf{A}}).k }$. Next, in a second step, we obtain $d\mathbf{a}_{gn} = d\mathbf{a}_{min}$ as the matrix-vector product
\begin{equation*}
d\mathbf{a}_{gn}  = \mathbf{N}^{\bot} d\mathbf{\bar{a}}_{gn} \ .
\end{equation*}
Then, the general solution $d\mathbf{\widehat{a}}$ of the full linear least-squares problem can be also computed in a third step as
\begin{equation*}
d\mathbf{\widehat{a}} = d\mathbf{a}_{min} +  \mathbf{N}\mathbf{c} = \mathbf{N}^{\bot} d\mathbf{\bar{a}}_{gn} + \mathbf{N}\mathbf{c} \ ,
\end{equation*}
with $\mathbf{c} \in \mathbb{R}^{k_{\mathbf{A}}.k }$ is a vector of $k_{\mathbf{A}}.k$ arbitrary constants. Note that if we use the matrix $- \mathbf{M}(\mathbf{a})$ as an approximate Jacobian instead, we can also find all the solutions of the linear least-squares problem
\begin{equation*}
\min_{d\mathbf{a} \in \mathbb{R}^{p.k}}   \,   \frac{1}{2} \Vert \mathbf{r}(  \mathbf{a} ) - \mathbf{M} ( \mathbf{a}) d\mathbf{a} \Vert^{2}_{2}
\end{equation*}
in three steps using again the bases $\mathbf{N}$ and $\mathbf{N}^{\bot}$, but in the first step, we need to solve the reduced and nonsingular linear least-squares problem
\begin{equation*}
\min_{d\mathbf{a} \in \mathbb{R}^{( p - k_{\mathbf{A}}).k } }  \,   \frac{1}{2} \Vert \mathbf{r}(  \mathbf{a} ) - \mathbf{M}(\mathbf{a}) \mathbf{N}^{\bot} d\mathbf{a} \Vert^{2}_{2}  \ ,
\end{equation*}
instead of the previous one involving the rank-deficient matrix  $\mathit{J}( \mathbf{r}(\mathbf{a}) )$.

Alternatively, it is also possible to obtain $d\mathbf{a}_{gn}$ and all solutions of the above Gauss-Newton linear least-squares problems involving the Jacobian matrix  $\mathit{J}( \mathbf{r}(\mathbf{a}) )$ or its approximation $-\mathbf{M}(\mathbf{a})$ using only a basis $\mathbf{N}$ of $\emph{null}( \mathit{J}( \mathbf{r}(\mathbf{a}) ) ) = \emph{null}(  \mathbf{M}(\mathbf{a}) )$ as first noted by Okatani et al.~\cite{OYD2011}. Using the fact, demonstrated above, that $d\mathbf{a}_{gn} = d\mathbf{a}_{min}$ is the unique solution of these linear least-squares problems, which belongs to $\emph{null}( \mathit{J}( \mathbf{r}(\mathbf{a}) ) )^{\bot}$, it is not difficult to see that $d\mathbf{a}_{gn}$ is also the unique solution among the infinite set of solutions $d\mathbf{\widehat{a}}$ of these linear least-squares problems, which verifies the equality
\begin{equation*}
 \mathbf{N}^{T} d\mathbf{\widehat{a}} =  \mathbf{0}^{k_{\mathbf{A}}.k} \ ,
\end{equation*}
if the columns of $\mathbf{N}$ form a basis of $\emph{null}( \mathit{J}( \mathbf{r}(\mathbf{a}) ) ) = \emph{null}(  \mathbf{M}(\mathbf{a}) )$. 

In addition, using results in the demonstration of Theorem~\ref{theo5.2:box} and assuming for simplicity that $k_{\mathbf{A}} = k$, e.g., $\mathbf{A}$ is of full column-rank, we observe that the matrix defined by
\begin{equation*}
\mathbf{N} = \mathbf{K}_{(p,k)} ( \mathbf{I}_{k} \otimes \mathbf{A} )
\end{equation*}
is also of full column-rank (e.g., $\emph{rank} (\mathbf{N} ) = k.k$), and that
\begin{equation*}
\emph{null} \Big( \mathit{J} \big( \mathbf{r}(\mathbf{a}) \big) \Big) = \emph{null} \big( \mathbf{M}(\mathbf{a}) \big) = \emph{ran}( \mathbf{N} )  \ ,
\end{equation*}
if the hypotheses of Theorem~\ref{theo5.3:box} are satisfied. In other words, it is very easy to compute a basis $\mathbf{N}$ of  $\emph{null}\big( \mathit{J}( \mathbf{r}(\mathbf{a}) ) \big) = \emph{null}( \mathbf{M}(\mathbf{a}) )$ in practical applications or even an orthonormal basis of this linear subspace of $\mathbb{R}^{ p.k }$ with the help of Corollary~\ref{corol5.6:box}  demonstrated below.

Furthermore, it is also very easy to introduce the linear constraint $ \mathbf{N}^{T} d\mathbf{a} =  \mathbf{0}^{k_{\mathbf{A}}.k}$  in the linear least-squares problems, which must be solved for computing the correction vector $d\mathbf{a}_{gn}$ at each iteration of the Gauss-Newton algorithm, as follows
\begin{equation*}
d\mathbf{a}_{gn} = \text{Arg} \min_{d\mathbf{a} \in \mathbb{R}^{p.k}}   \,   \frac{1}{2} \Vert \mathbf{r}(  \mathbf{a} ) + \mathit{J}( \mathbf{r}(\mathbf{a}) ) d\mathbf{a} \Vert^{2}_{2} + \frac{1}{2} \Vert \mathbf{N}^{T} d\mathbf{a} \Vert^{2}_{2} \ ,
\end{equation*}
or                       
\begin{equation*}
d\mathbf{a}_{gn} = \text{Arg} \min_{d\mathbf{a} \in \mathbb{R}^{p.k}}   \,   \frac{1}{2} \Vert \mathbf{r}(  \mathbf{a} ) - \mathbf{M}(\mathbf{a}) d\mathbf{a} \Vert^{2}_{2} + \frac{1}{2} \Vert \mathbf{N}^{T} d\mathbf{a} \Vert^{2}_{2} \ ,
\end{equation*}
if we use the approximate Jacobian matrix $- \mathbf{M}(\mathbf{a})$. These two "constrained" linear least-squares problems are also, respectively, equivalent to the following standard linear least-squares problems
\begin{equation*}
d\mathbf{a}_{gn} = \text{Arg} \min_{d\mathbf{a} \in \mathbb{R}^{p.k}}   \,   \frac{1}{2} \big\Vert \begin{bmatrix} \mathbf{r}(  \mathbf{a} )  \\  \mathbf{0}^{k_{\mathbf{A}}.k}   \end{bmatrix} + \begin{bmatrix} \mathit{J}( \mathbf{r}(\mathbf{a}) ) \ , \\   \mathbf{N}^{T} \end{bmatrix} d\mathbf{a} \big\Vert^{2}_{2}
\end{equation*}
and
\begin{equation*}
d\mathbf{a}_{gn} = \text{Arg} \min_{d\mathbf{a} \in \mathbb{R}^{p.k}}   \,   \frac{1}{2} \big\Vert \begin{bmatrix} \mathbf{r}(  \mathbf{a} )  \\  \mathbf{0}^{k_{\mathbf{A}}.k}   \end{bmatrix} - \begin{bmatrix} \mathbf{M}(\mathbf{a})  \\   \mathbf{N}^{T} \end{bmatrix} d\mathbf{a} \big\Vert^{2}_{2} \ ,
\end{equation*}
which are both easily solved and have an unique solution as the associated coefficient matrices are nonsingular if the columns of $\mathbf{N}$ form a basis of $\emph{null}( \mathit{J}( \mathbf{r}(\mathbf{a}) ) ) = \emph{null}(  \mathbf{M}(\mathbf{a}) )$.

To demonstrate this result, it suffices to show that the null space of these matrices is reduced to the zero vector. To this end, we first observe that
\begin{equation*}
\emph{null} \big(  \begin{bmatrix} \mathit{J}( \mathbf{r}(\mathbf{a}) ) \\   \mathbf{N}^{T} \end{bmatrix} \big) = \emph{null} \big( \mathit{J}( \mathbf{r}(\mathbf{a}) ) \big) \cap  \emph{null}(  \mathbf{N}^{T} ) = \emph{null} \big(  \begin{bmatrix}  \mathbf{M}(\mathbf{a}) \\   \mathbf{N}^{T} \end{bmatrix} \big) \ ,
\end{equation*}
if the hypotheses of Theorem~\ref{theo5.3:box} are satisfied. Moreover, the following relationships hold (see equation~\eqref{eq:rank} in Subsection~\ref{lin_alg:box})
\begin{equation*}
\emph{null}(  \mathbf{N}^{T} )^\bot = \emph{ran}(  \mathbf{N} ) = \emph{null}  \big( \mathit{J}( \mathbf{r}(\mathbf{a}) )  \big) = \emph{null}(  \mathbf{M}(\mathbf{a}) )  \ ,
\end{equation*}
and this implies, finally, that
\begin{equation} \label{eq:N_mat}
\emph{null}(  \begin{bmatrix} \mathit{J}( \mathbf{r}(\mathbf{a}) ) \\   \mathbf{N}^{T} \end{bmatrix} ) =  \emph{null}(  \begin{bmatrix}  \mathbf{M}(\mathbf{a}) \\   \mathbf{N}^{T} \end{bmatrix} ) =  \emph{null}(  \mathbf{N}^{T} )^\bot \cap  \emph{null}(  \mathbf{N}^{T} ) = \lbrace \mathbf{0}^{p.k} \rbrace  \ ,
\end{equation}
which demonstrates that the corresponding matrices are effectively nonsingular if the columns of $\mathbf{N}$ form a basis of $\emph{null}( \mathit{J}( \mathbf{r}(\mathbf{a}) ) ) = \emph{null}(  \mathbf{M}(\mathbf{a}) )$.

The above results can also be used to improve the Levenberg-Marquardt algorithm. For example, using an exact Jacobian matrix, $\mathit{J}( \mathbf{r}(\mathbf{a}) )$, an accurate Levenberg-Marquardt's correction vector $d\mathbf{a}_{lm}$ can also be obtained in two steps. First, by solving the following reduced and regularized linear least-squares problem
\begin{equation*}
\min_{d\mathbf{a} \in \mathbb{R}^{( p - k_{\mathbf{A}}).k } }  \,    \frac{1}{2} \big\Vert  \begin{bmatrix} \mathbf{r}(  \mathbf{a} ) \\ \mathbf{0}^{( p - k_{\mathbf{A}}).k }  \end{bmatrix}  +  \begin{bmatrix}  \mathit{J}( \mathbf{r}(\mathbf{a}) )\mathbf{N}^{\bot}    \\  \sqrt{\lambda} \mathbf{D}  \end{bmatrix}  d\mathbf{a}  \big\Vert^{2}_{2}   =  \frac{1}{2} \big\Vert \mathbf{r}(  \mathbf{a} ) + \mathit{J}( \mathbf{r}(\mathbf{a}) )\mathbf{N}^{\bot}  d\mathbf{a}  \big\Vert^{2}_{2} +  \frac{\lambda}{2}  \big\Vert \mathbf{D} d\mathbf{a} \big\Vert^{2}_{2} \ ,
\end{equation*}
where $\lambda$ is the damping Marquardt parameter and  $\mathbf{D}$ is a diagonal scaling matrix of dimension $( p - k_{\mathbf{A}}).k$. This damped linear least-squares problem has always an unique solution, $d\mathbf{\bar{a}}_{lm}$, even when $\lambda$ tends to zero if the hypotheses of Theorem~\ref{theo5.3:box} are satisfied. Once $d\mathbf{\bar{a}}_{lm}$ has been found, $d\mathbf{a}_{lm}$ can be computed by the matrix-vector product
\begin{equation*}
d\mathbf{a}_{lm} = \mathbf{N}^{\bot} d\mathbf{\bar{a}}_{lm} \ ,
\end{equation*}
as for the correction vector $d\mathbf{a}_{gn}$ in the Gauss-Newton algorithm.

Alternatively, the linear constraint
\begin{equation*}
 \mathbf{N}^{T} d\mathbf{a} =  \mathbf{0}^{k_{\mathbf{A}}.k}
\end{equation*}
can also be introduced in the  linear least-squares problem, which must be solved for computing the correction vector at each iteration of the Levenberg-Marquardt algorithm. As an illustration, if we use an exact Jacobian matrix, a correction vector can be computed as
\begin{equation*}
d\mathbf{a}^{'}_{lm} = \text{Arg} \min_{d\mathbf{a} \in \mathbb{R}^{p.k}}   \,   \frac{1}{2} \Vert \mathbf{r}(  \mathbf{a} ) + \mathit{J}( \mathbf{r}(\mathbf{a}) ) d\mathbf{a} \Vert^{2}_{2} +  \frac{1}{2} \Vert \mathbf{N}^{T} d\mathbf{a} \Vert^{2}_{2} +  \frac{\lambda}{2}  \big\Vert \mathbf{D} d\mathbf{a} \big\Vert^{2}_{2} \ ,
\end{equation*}
where $\lambda$ is the damping Marquardt parameter and  $\mathbf{D}$ is now a diagonal scaling matrix of dimension $k.p$. This problem is also equivalent to the standard linear least-squares problem
\begin{equation*}
d\mathbf{a}^{'}_{lm} = \text{Arg} \min_{d\mathbf{a} \in \mathbb{R}^{p.k}}   \,   \frac{1}{2} \big\Vert \begin{bmatrix} \mathbf{r}(  \mathbf{a} )  \\  \mathbf{0}^{k_{\mathbf{A}}.k}  \\  \mathbf{0}^{k.p} \end{bmatrix} + \begin{bmatrix} \mathit{J}( \mathbf{r}(\mathbf{a}) )   \\   \mathbf{N}^{T}  \\    \sqrt{\lambda} \mathbf{D} \end{bmatrix} d\mathbf{a}  \big\Vert^{2}_{2} \ .
\end{equation*}
However, contrary to the case of the Gauss-Newton method, the correction vectors, $d\mathbf{a}_{lm}$ and $d\mathbf{a}^{'}_{lm}$, obtained by these two alternative formulations of the Levenberg-Marquardt algorithm will differ in general. More precisely, if $\lambda \neq 0$, we cannot assume that it always exists $\mathbf{c}  \in \mathbb{R}^{( p - k_{\mathbf{A}}).k }$ such that
\begin{equation*}
d\mathbf{a}^{'}_{lm} = \mathbf{N}^{\bot} \mathbf{c}  \ ,
\end{equation*}
as we only have $\mathbf{N}^{T} d\mathbf{a}^{'}_{lm} \approx \mathbf{0}^{k_{\mathbf{A}}.k}$,  but not exactly $\mathbf{N}^{T} d\mathbf{a}^{'}_{lm} = \mathbf{0}^{k_{\mathbf{A}}.k}$ as for $d\mathbf{a}_{lm}$. Thus, in these conditions, $d\mathbf{a}^{'}_{lm} \notin \emph{null}( \mathit{J}( \mathbf{r}(\mathbf{a}) ) )^{\bot}$, while $d\mathbf{a}_{lm} \in \emph{null}( \mathit{J}( \mathbf{r}(\mathbf{a}) ) )^{\bot}$, and these two correction vectors will not be equal in general. Moreover, the fact that $d\mathbf{a}_{lm} \in \emph{null}( \mathit{J}( \mathbf{r}(\mathbf{a}) ) )^{\bot}$ implies that the version of the Levenberg-Marquardt algorithm using this correction vector can also be considered as a Riemannian optimization algorithm operating directly on the Grassmann manifold $\text{Gr}(p,k)$~\cite{AMS2008}\cite{B2023} in the same way as the Gauss-Newton algorithm discussed above (see below for details), while the version using $d\mathbf{a}^{'}_{lm}$ as a correction step does not enjoy this theoretical property.

 At first sight, the approach using only a basis $\mathbf{N}$ of  $\emph{null}\big( \mathit{J}( \mathbf{r}(\mathbf{a}) ) \big) = \emph{null}( \mathbf{M}(\mathbf{a}) )$ and a linear constraint for computing the Gauss-Newton and Levenberg-Marquardt directions at each iteration, seems to be much cheaper than the first approach, which needs to compute a nonsingular matrix $\mathbf{N}^{\bot} \in  \mathbb{R}^{k.p \times ( p - k_{\mathbf{A}}).k }$ whose columns form a basis of $\emph{null}( \mathit{J}( \mathbf{r}(\mathbf{a}) ) )^{\bot}$, to multiply the huge Jacobian matrix $\mathit{J}( \mathbf{r}(\mathbf{a})$  (or its approximation) by this matrix $\mathbf{N}^{\bot}$ and, finally, to compute the matrix-vector products $d\mathbf{a}_{gn}  =  \mathbf{N}^{\bot} d\mathbf{\bar{a}}_{gn}$ or $d\mathbf{a}_{lm} = \mathbf{N}^{\bot} d\mathbf{\bar{a}}_{lm}$ in the case of the Levenberg-Marquardt algorithm~\cite{OYD2011}.

However, the next corollary shows that the overhead cost incurred by the first approach can be drastically reduced since it is easy to obtain orthonormal bases of  $\emph{null}\big( \mathit{J}( \mathbf{r}(\mathbf{a}) ) \big)$ and its orthogonal complement in $ \mathbb{R}^{k.p}$ if the conditions of Theorem~\ref{theo5.3:box} are fulfilled. Furthermore, thanks to the particular form of these orthonormal bases, the above matrix product between the  Jacobian matrix (or its approximation) and a basis of $\emph{null}\big( \mathit{J}( \mathbf{r}(\mathbf{a}) ) \big)^{\bot}$ can be computed very efficiently with almost the same cost as evaluating the Jacobian matrix or its approximation itself as also demonstrated in this corollary.
\\
\begin{corol5.6} \label{corol5.6:box}
With the same notations and hypotheses as in Theorem~\ref{theo5.3:box}, let $ \mathbf{O} \in \mathbb{O}^{p \times k_{\mathbf{A}} } $ be an orthonormal basis of $ \emph{ran}( \mathbf{A} ) $ and  $ \mathbf{O}^{\bot} \in \mathbb{O}^{p \times ( p - k_{\mathbf{A}} ) } $ be an orthonormal basis of $ \emph{ran}( \mathbf{A} )^{\bot} $, then
\begin{align*}
   \mathbf{\bar{O}} & =  \mathbf{K}_{(p,k)} ( \mathbf{I}_{k}  \otimes  \mathbf{O}) \text{ is an orthonormal basis of } \emph{null}\big( \mathit{J}( \mathbf{r}(\mathbf{a}) ) \big) = \emph{null}( \mathbf{M}(\mathbf{a}) ) \ ,  \\
   \mathbf{\bar{O}}^{\bot} & =  \mathbf{K}_{(p,k)} ( \mathbf{I}_{k}  \otimes  \mathbf{O}^{\bot}) \text{ is an orthonormal basis of } \emph{null}\big( \mathit{J}( \mathbf{r}(\mathbf{a}) ) \big)^{\bot} = \emph{null}( \mathbf{M}(\mathbf{a}) )^{\bot}  \ .
\end{align*}
Furthermore, we have
\begin{align*}
\mathbf{M}(\mathbf{a})\mathbf{\bar{O}}^{\bot} & =  \mathbf{P}_{\mathbf{F}(\mathbf{a})}^{\bot} \mathbf{U}(\mathbf{a}) \mathbf{\bar{O}}^{\bot} = \mathbf{P}_{\mathbf{F}(\mathbf{a})}^{\bot} \emph{diag}\big( \emph{vec}(\sqrt{\mathbf{W}})\big) \big( \mathbf{\widehat{B}}^{T} \otimes \mathbf{O}^{\bot} \big) \ ,  \\
\mathbf{L}(\mathbf{a})\mathbf{\bar{O}}^{\bot} & = \big( \mathbf{F}(\mathbf{a})^{+}  \big)^{T}  \mathbf{V}(\mathbf{a}) \mathbf{\bar{O}}^{\bot} = \big( \mathbf{F}(\mathbf{a})^{+} \big)^{T} \Big( \big( \mathbf{W} \odot P_{\Omega}(\mathbf{X} -\mathbf{A}\mathbf{\widehat{B}}) \big)^{T} \mathbf{O}^{\bot} \otimes \mathbf{I}_{k} \Big) \mathbf{K}_{(p-k,p)} \ .
\end{align*}
\end{corol5.6}
\begin{proof}
Using the results of Theorem~\ref{theo5.2:box}, we first observe that
\begin{equation*}
\emph{null}\big( \mathit{J}( \mathbf{r}(\mathbf{a}) ) \big) = \emph{null} \big( \mathbf{M}(\mathbf{a}) \big) = \emph{ran}( \mathbf{N} )  \ ,
\end{equation*}
where $ \mathbf{N} = \mathbf{K}_{(p,k)} ( \mathbf{I}_{k} \otimes \mathbf{A} ) $.

Now, since $ \mathbf{O} $ is an orthonormal basis of the column space of $ \mathbf{A} $, then $ \exists \mathbf{C} \in \mathbb{R}^{k_{\mathbf{A}} \times k} $ such that $ \mathbf{A} = \mathbf{O} \mathbf{C} $, $ \emph{rank}( \mathbf{C} ) = k_{\mathbf{A}} $ and $ \mathbf{C} $ is uniquely determined (see Theorem 1 of Marsaglia and Styan~\cite{MS1974}). Moreover, since $ \mathbf{C} $ has full row-rank, $ \mathbf{C} $ admits a right-inverse $\mathbf{R}$ such that $ \mathbf{C} \mathbf{R} = \mathbf{I}_{k_{\mathbf{A}}} $. From this equality, we deduce that
\begin{equation*}
    \mathbf{A} \mathbf{R} = \mathbf{O} \mathbf{C} \mathbf{R} = \mathbf{O} \ .
\end{equation*}
Using these properties and equation~\eqref{eq:vec_kronprod}, we have, $ \forall \mathbf{Z} \in \mathbb{R}^{k \times k} $,
\begin{align*}
    ( \mathbf{I}_{k}  \otimes  \mathbf{A} ) \emph{vec} ( \mathbf{Z} )   &  =    ( \mathbf{I}_{k}  \otimes  \mathbf{O} \mathbf{C} ) \emph{vec} ( \mathbf{Z} )  \\
                                                                                                           &  =      \emph{vec} ( \mathbf{O} \mathbf{C} \mathbf{Z} )  \\
                                                                                                           &  =      ( \mathbf{I}_{k}  \otimes  \mathbf{O} ) \emph{vec} ( \mathbf{C} \mathbf{Z} ) \ ,
\end{align*}
and, $ \forall \mathbf{T} \in \mathbb{R}^{k_{\mathbf{A} } \times k_{\mathbf{A}} } $,
\begin{align*}
    ( \mathbf{I}_{k}  \otimes  \mathbf{O} ) \emph{vec} ( \mathbf{T} )   &  =   ( \mathbf{I}_{k}  \otimes  \mathbf{A} \mathbf{R} ) \emph{vec} ( \mathbf{T} )  \\
                                                                                                           &  =    \emph{vec} ( \mathbf{A} \mathbf{R} \mathbf{T} )  \\
                                                                                                           &  =    ( \mathbf{I}_{k}  \otimes  \mathbf{A} ) \emph{vec} ( \mathbf{R} \mathbf{T} ) \ .
\end{align*}
From these equalities, it can be easily proved that
\begin{equation*}
    \emph{ran} ( \mathbf{N} ) = \emph{ran} \big( \mathbf{K}_{(p,k)} ( \mathbf{I}_{k} \otimes \mathbf{O} ) \big) = \emph{ran}( \mathbf{\bar{O}} ) \ .
\end{equation*}
Moreover, $  \mathbf{\bar{O}} = \mathbf{K}_{(p,k)} ( \mathbf{I}_{k} \otimes \mathbf{O} ) $ is a matrix with orthonormal columns since:
\begin{align*}
     \big( \mathbf{K}_{(p,k)} ( \mathbf{I}_{k} \otimes \mathbf{O} ) \big)^{T} \mathbf{K}_{(p,k)} ( \mathbf{I}_{k} \otimes \mathbf{O} )
                                                                                              & =  ( \mathbf{I}_{k} \otimes \mathbf{O}^{T} )( \mathbf{I}_{k} \otimes \mathbf{O} )   \\
                                                                                              & =  ( \mathbf{I}_{k} \otimes \mathbf{O}^{T} \mathbf{O} )  \\
                                                                                              & =  \mathbf{I}_{k} \otimes \mathbf{I}_{k_{\mathbf{A}}}            \\
                                                                                              & =  \mathbf{I}_{k.k_{\mathbf{A}}} \ .
\end{align*}
A similar argument applies to $ \mathbf{\bar{O}}^{\bot} = \mathbf{K}_{(p,k)} ( \mathbf{I}_{k} \otimes \mathbf{O}^{\bot} ) $. Finally, we have
\begin{align*}
         \big( \mathbf{K}_{(p,k)} ( \mathbf{I}_{k} \otimes \mathbf{O} ) \big)^{T} \mathbf{K}_{(p,k)} ( \mathbf{I}_{k} \otimes \mathbf{O}^{\bot} )
                                                                                              & =  ( \mathbf{I}_{k} \otimes \mathbf{O}^{T} )( \mathbf{I}_{k} \otimes \mathbf{O}^{\bot} )   \\
                                                                                              & =  ( \mathbf{I}_{k} \otimes \mathbf{O}^{T} \mathbf{O}^{\bot} )  \\
                                                                                              & =  \mathbf{I}_{k} \otimes \mathbf{0}^{k_{A} \times (p - k_{A}) }      \\
                                                                                              & =  \mathbf{0}^{k.k_{A} \times k.(p - k_{A}) } \  ,
\end{align*}
and it follows that $ \mathbf{\bar{O}}$  is an orthonormal basis of $  \emph{null}\big( \mathit{J}( \mathbf{r}(\mathbf{a}) ) \big)$ and  $ \mathbf{\bar{O}}^{\bot}$ is an orthonormal basis of $\emph{null}\big( \mathit{J}( \mathbf{r}(\mathbf{a}) ) \big)^{\bot} $ as stated in the corollary. Furthermore, the columns of the matrix $\begin{bmatrix} \mathbf{\bar{O}}   & \mathbf{\bar{O}}^{\bot} \end{bmatrix}$ form an orthonormal basis of $\mathbb{R}^{k.p}$.

Finally, to demonstrate the second part of the corollary, let us evaluate compactly the matrix products $\mathbf{U} \mathbf{\bar{O}}^{\bot}$ and $\mathbf{V} \mathbf{\bar{O}}^{\bot}$ using the properties of the Kronecker product, commutation matrix and $\emph{vec}$ operator stated in Subsection~\ref{multlin_alg:box}. We have
\begin{align*}
\mathbf{U} \mathbf{\bar{O}}^{\bot} & =  \emph{diag} \big( \emph{vec}(\sqrt{\mathbf{W}} )  \big) (\mathbf{\widehat{B}}^{T} \otimes \mathbf{I}_{p}) \mathbf{K}_{(k,p)} \mathbf{K}_{(p,k)} ( \mathbf{I}_{k}  \otimes  \mathbf{O}^{\bot} ) \\
                                                       & =  \emph{diag} \big( \emph{vec}(\sqrt{\mathbf{W}} )  \big) (\mathbf{\widehat{B}}^{T} \otimes \mathbf{I}_{p})  ( \mathbf{I}_{k}  \otimes  \mathbf{O}^{\bot})  \\
                                                       & =  \emph{diag} \big( \emph{vec}(\sqrt{\mathbf{W}} )  \big) (\mathbf{\widehat{B}}^{T} \otimes \mathbf{O}^{\bot} ) \  ,
\end{align*}
and also
\begin{align*}
\mathbf{V} \mathbf{\bar{O}}^{\bot} & =  \big( ( \mathbf{W} \odot P_{\Omega}(\mathbf{X} -\mathbf{A}\mathbf{\widehat{B}}) )^{T} \otimes \mathbf{I}_{k} \big) \mathbf{K}_{(p,k)} ( \mathbf{I}_{k}  \otimes  \mathbf{O}^{\bot} ) \\
                                                       & =  \big( ( \mathbf{W} \odot P_{\Omega}(\mathbf{X} -\mathbf{A}\mathbf{\widehat{B}}) )^{T} \otimes \mathbf{I}_{k} \big) ( \mathbf{O}^{\bot} \otimes \mathbf{I}_{k}  ) \mathbf{K}_{(p-k,k)} \\
                                                       & =  \big( ( \mathbf{W} \odot P_{\Omega}(\mathbf{X} -\mathbf{A}\mathbf{\widehat{B}}) \big)^{T} \mathbf{O}^{\bot} \otimes \mathbf{I}_{k} \big) \mathbf{K}_{(p-k,p)} \  ,
\end{align*}
which concludes the demonstration of the corollary.
\\
\end{proof}

Using Corollary~\ref{corol5.6:box} and the preceding results, we can write the correction vectors $d\mathbf{a}_{gn}$ and $d\mathbf{a}_{lm}$ of the Gauss-Newton and Levenberg-Marquardt algorithms as the matrix-vector products
\begin{equation*}
d\mathbf{a}_{gn} = \mathbf{\bar{O}}^{\bot} d\mathbf{\bar{a}}_{gn}  \text{ and } d\mathbf{a}_{lm} = \mathbf{\bar{O}}^{\bot} d\mathbf{\bar{a}}_{lm}  \ ,
\end{equation*}
where $d\mathbf{\bar{a}}_{gn}$ and $d\mathbf{\bar{a}}_{lm}$ are, respectively, the solutions of the problems
\begin{equation*}
d\mathbf{\bar{a}}_{gn}  =  \text{Arg}\min_{d\mathbf{a} \in \mathbb{R}^{( p - k_{\mathbf{A}}).k } }  \,   \frac{1}{2} \big\Vert \mathbf{r}(  \mathbf{a} ) + \mathit{J}\big( \mathbf{r}(\mathbf{a}) \big)\mathbf{\bar{O}}^{\bot} d\mathbf{a}  \big\Vert^{2}_{2}
\end{equation*}
and
\begin{equation*}
d\mathbf{\bar{a}}_{lm}  =  \text{Arg}\min_{d\mathbf{a} \in \mathbb{R}^{( p - k_{\mathbf{A}}).k } }  \,   \frac{1}{2} \big\Vert \mathbf{r}(  \mathbf{a} ) + \mathit{J}\big( \mathbf{r}(\mathbf{a}) \big)\mathbf{\bar{O}}^{\bot} d\mathbf{a}  \big\Vert^{2}_{2} +  \frac{\lambda}{2}  \big\Vert \mathbf{D} d\mathbf{a} \big\Vert^{2}_{2}  \ ,
\end{equation*}
or of similar linear least-squares problems involving the approximate Jacobian matrix $-\mathbf{M}(\mathbf{a})$ instead of $\mathit{J}( \mathbf{r}(\mathbf{a}) )$.
Defining  $d\mathbf{A}_{gn}  \in \mathbb{R}^{ p \times k  }$ and $d\mathbf{\bar{A}}_{gn}  \in \mathbb{R}^{ (p - k_\mathbf{A}) \times k  }$ such that  $d\mathbf{a}_{gn} = \emph{vec}(  d\mathbf{A}_{gn}^{T} )$ and $d\mathbf{\bar{a}}_{gn}  = \emph{vec}(  d\mathbf{\bar{A}}_{gn} )$, we have (using equation~\eqref{eq:vec_kronprod} and Lemma~\ref{theo2.2:box} in Subsection~\ref{multlin_alg:box})
\begin{equation*}
\emph{vec}(  d\mathbf{A}_{gn}^{T} ) = \mathbf{K}_{(p,k)} ( \mathbf{I}_{k}  \otimes  \mathbf{O}^{\bot}) \emph{vec}(  d\mathbf{\bar{A}}_{gn} ) = \mathbf{K}_{(p,k)}  \emph{vec}(  \mathbf{O}^{\bot} d\mathbf{\bar{A}}_{gn} ) = \emph{vec}\big(  ( \mathbf{O}^{\bot} d\mathbf{\bar{A}}_{gn} )^{T} \big)  \ ,
\end{equation*}
which implies that $d\mathbf{A}_{gn} = \mathbf{O}^{\bot} d\mathbf{\bar{A}}_{gn}$. Obviously, the equality $d\mathbf{A}_{lm} = \mathbf{O}^{\bot} d\mathbf{\bar{A}}_{lm}$ can be derived in a similar fashion.

Thus, the columns of the perturbation matrices $d\mathbf{A}_{gn}$ and $d\mathbf{A}_{lm}$ belong to $\emph{ran}( \mathbf{O}^{\bot} ) = \emph{ran}( \mathbf{A} )^{\bot}$. In other words, these variations of the variable projection Gauss-Newton and Levenberg-Marquardt methods described above to deal with the singularity of the Jacobian matrix $\mathit{J}( \mathbf{r}(\mathbf{a}) )$, or of its approximation $-\mathbf{M}(\mathbf{a})$, consider only search directions of the form $d\mathbf{A} = \mathbf{O}^{\bot} \mathbf{C}$ where $\mathbf{C}  \in \mathbb{R}^{(p - k_\mathbf{A}) \times k  }$. This is consistent with Remark~\ref{remark5.1:box} and the fact that we have only $k.( p - k )$ degrees of freedom to update $\mathbf{A}$ if $\emph{rank}( \mathbf{A} ) = k$ at each iteration of the Gauss-Newton or Levenberg-Marquardt algorithms.

In addition, if $\mathbf{W} \in  \mathbb{R}^{ p \times n  }_{+*}$, the unvectorized form of the cost function $\psi(.)$ (e.g., $\psi  \circ  h^{-1}(.)$ where $h^{-1} ( \mathbf{A} ) = \emph{vec}(  \mathbf{A}^{T} ) = \mathbf{a}, \forall \mathbf{A} \in \mathbb{R}^{p \times k}$, with $h(.)$ and $h^{-1}(.)$ defined in equation~\eqref{eq:h_func} of Subsection~\ref{varpro_wlra:box}), which is used in the~\eqref{eq:VP1} formulation of the WLRA problem, is smooth  (e.g., of class $C^{\infty}$) over the subdomains $\mathbb{R}^{ p \times k  }_{k}$ or $\mathbb{O}^{p \times k}$ according to Corollaries~\ref{corol3.3:box}  and~\ref{corol5.1:box}. Note that this (unvectorized) cost function $\psi(.)$ can be an instance of the~\eqref{eq:VP1} formulation of the cost function  $g_{\lambda}(.)$  introduced by Boumal and Absil~\cite{BA2011}\cite{BA2015}, defined in equation~\eqref{eq:g_func}, and already discussed in Subsection~\ref{approx_wlra:box} since its associated weight matrix $\mathbf{W}_\lambda  \in \mathbb{R}^{ p \times n  }_{+*}$ if $\lambda > 0$. In the same conditions, if $\mathbf{A} \in \mathbb{R}^{ p \times k  }_{k}$ or $\mathbf{A} \in \mathbb{O}^{p \times k}$, we have $  \emph{rank}( \mathbf{M}(\mathbf{a}) ) =  \emph{rank}( \mathit{J}( \mathbf{r}(\mathbf{a}) ) )  = k.( p - k) $ according to Corollary~\ref{corol5.5:box}. Using these different results,  we can recast the variable projection formulation of the WLRA problem as an optimization problem on the Grassmann  manifold $\text{Gr}(p,k)$~\cite{BA2015}\cite{B2023} and the above variable projection Gauss-Newton and Levenberg-Marquardt methods to solve the WLRA problem as Riemannian optimization algorithms operating on this Grassmann  manifold~\cite{AMS2008}\cite{B2023} as their numerical behavior only depends on $\mathring{\mathbf{A}}  = \emph{ran}( \mathbf{A} ) \in \text{Gr}(p,k)$, for $\mathbf{A} \in \mathbb{R}^{ p \times k  }_{k}$ or  $\mathbf{A} \in \mathbb{O}^{p \times k}$, and not on the arbitrarily chosen matrix $\mathbf{A}$ to represent  $\emph{ran}( \mathbf{A} )$ according to Corollaries~\ref{corol3.1:box},~\ref{corol3.2:box} and Remark~\ref{remark3.7:box}.

More precisely, the smooth function $\psi  \circ  h^{-1}(.)$ defined on the smooth submanifold $\mathbb{R}^{ p \times k  }_{k}$ embedded in $\mathbb{R}^{p \times k}$ is invariant on the equivalence classes of the equivalence relation $\sim$ defined on $\mathbb{R}^{ p \times k  }_{k}$ by, $\forall \mathbf{A}, \mathbf{C} \in  \mathbb{R}^{ p \times k  }_{k}$,
\begin{equation*}
 \mathbf{A} \sim  \mathbf{C} \text{ if and only if it exists } \mathbf{D}  \in \mathbb{R}^{k \times k}_{k}  \text{ such that }  \mathbf{A} = \mathbf{C} \mathbf{D} \ ,
\end{equation*}
according to Corollary~\ref{corol3.2:box}. In this setting, we can say that $\text{Gr}(p,k)$ is the quotient of $\mathbb{R}^{ p \times k  }_{k}$ by the action of the group  $\mathbb{R}^{k \times k}_{k}$ following the terminology introduced in Subsection~\ref{calculus:box}. Alternatively, if we prefer to work with orthogonal matrices (e.g., with the Stiefel manifold $\mathbb{O}^{p \times k}$), we can consider the restriction of  $\psi  \circ  h^{-1}(.)$ to $\mathbb{O}^{p \times k}$ and the equivalence relation  $\sim$ defined on $\mathbb{O}^{p \times k}$ by, $\forall \mathbf{A}, \mathbf{C} \in \mathbb{O}^{p \times k}$,
\begin{equation*}
 \mathbf{A} \sim  \mathbf{C} \text{ if and only if it exists } \mathbf{D}  \in \mathbb{O}^{k \times k} \text{ such that }  \mathbf{A} = \mathbf{C} \mathbf{D} \ ,
\end{equation*}
and, similarly,  $\psi  \circ  h^{-1}(.)$ is invariant on the equivalence classes of this equivalence relation according to Corollary~\ref{corol3.2:box} and  $\text{Gr}(p,k)$ is defined now as the quotient of the Stiefel manifold $\mathbb{O}^{  p \times k  }$ by the action of the orthogonal group  $\mathbb{O}^{k \times k}$.

Moreover, as a quotient manifold (see Subsection~\ref{calculus:box} and Chapter 9 of~\cite{B2023} ), the Grassmannian admits a tangent space at $\mathring{\mathbf{A}} = \emph{ran}( \mathbf{A} )  \in \text{Gr}(p,k)$, $\forall \mathbf{A} \in \mathbb{R}^{ p \times k  }_{k}$ or $\forall \mathbf{A}  \in \mathbb{O}^{p \times k}$, designed by $\mathcal{T}_{\mathring{\mathbf{A}}} \text{Gr}(p,k)$ (in the terminology of Subsection~\ref{calculus:box}), which can be identified uniquely with the linear subspace of $\mathbb{R}^{ p \times k  }$ of dimension $k.( p - k )$ defined by
\begin{equation*}
\mathcal{T}_{\mathbf{A}} \text{Gr}(p,k) =  \big\{ \mathbf{D} \in \mathbb{R}^{ p \times k  }  \text{ } / \text{ } \mathbf{A}^{T} \mathbf{D} =  \mathbf{0}^{ k \times k  } \big\}  \ .
\end{equation*}
With this identification, $\mathcal{T}_{\mathring{\mathbf{A}}} \text{Gr}(p,k)$ is nothing else than the horizontal space,  $\mathcal{H}_{\mathbf{A}} \mathbb{R}^{ p \times k  }_{k}$, of $\mathbb{R}^{ p \times k  }_{k}$ at $\mathbf{A} \in \mathbb{R}^{ p \times k  }_{k}$ or, alternatively, the horizontal space,  $\mathcal{H}_{\mathbf{A}} \mathbb{O}^{ p \times k  }$, of $\mathbb{O}^{ p \times k  }$ at $\mathbf{A} \in \mathbb{O}^{ p \times k  }$; see Subsection~\ref{calculus:box} for details. Furthermore, the orthogonal projector onto $\mathcal{T}_{\mathbf{A}} \text{Gr}(p,k)$ with respect to the Frobenius inner product in $\mathbb{R}^{ p \times k  }$ is given by
\begin{equation*}
\mathbf{P}_{\mathcal{T}_{\mathbf{A}} \text{Gr}(p,k)} : \mathbb{R}^{ p \times k  }  \longrightarrow \mathcal{T}_{\mathbf{A}} \text{Gr}(p,k) , \mathbf{D}  \mapsto \mathbf{P}_{\mathcal{H}_{\mathbf{A}}  \mathbb{R}^{ p \times k  }_{k}} ( \mathbf{D} ) = ( \mathbf{I}_{p} - \mathbf{A}\mathbf{A}^{+} )\mathbf{D}  \ ,
\end{equation*}
or, equivalently, if we prefer to work with the Stiefel submanifold, by
\begin{equation*}
\mathbf{P}_{\mathcal{T}_{\mathbf{A}} \text{Gr}(p,k)} ( \mathbf{D} ) =  \mathbf{P}_{\mathcal{H}_{\mathbf{A}}  \mathbb{O}^{p \times k}} ( \mathbf{D} ) = ( \mathbf{I}_{p} - \mathbf{O}\mathbf{O}^{T} ) \mathbf{D} =   \mathbf{O}^{\bot} ( \mathbf{O}^{\bot} )^{T} \mathbf{D}  \ ,
\end{equation*}
where the columns of $\mathbf{O}$ and $\mathbf{O}^{\bot}$ form, respectively, orthogonal bases of $\emph{ran}(\mathbf{A})$ and $\emph{ran}(\mathbf{A})^{\bot}$, and $\mathbf{P}_{\mathcal{H}_{\mathbf{A}}  \mathbb{R}^{ p \times k  }_{k}}$ and $\mathbf{P}_{\mathcal{H}_{\mathbf{A}}  \mathbb{O}^{p \times k}}$ design, respectively, the orthogonal projectors onto the horizontal subspaces of the tangent spaces $\mathcal{T}_{\mathbf{A}}  \mathbb{R}^{ p \times k  }_{k}$ and $\mathcal{T}_{\mathbf{A}}  \mathbb{O}^{p \times k}$. See Subsection~\ref{calculus:box} and~\cite{AMS2008}\cite{B2023} for more details on the geometry of smooth manifolds, including the (quotient) Grassmann manifold.

Thus, in this Grassmann manifold framework, we have $d\mathbf{A}_{gn}, d\mathbf{A}_{lm}  \in \mathcal{T}_{\mathbf{A}} \text{Gr}(p,k)$, since $d\mathbf{A}_{gn} = \mathbf{O}^{\bot} d\mathbf{\bar{A}}_{gn}$ and $d\mathbf{A}_{lm} = \mathbf{O}^{\bot} d\mathbf{\bar{A}}_{lm}$, and this implies that the above Gauss-Newton and Levenberg-Marquardt algorithms can be interpreted exactly as Riemannian optimization methods operating on the Grassmann manifold $\text{Gr}(p,k)$: at the $(i+1)^{th}$ iteration, these algorithms move on the Grassmann manifold from $\mathbf{A}^{i} \in \mathbb{R}^{ p \times k  }_{k}$ along some direction prescribed by the tangent vectors $d\mathbf{A}^{i}_{gn}$ or $d\mathbf{A}^{i}_{lm}$ to 
\begin{equation*}
\mathbf{A}^{i+1} = \mathbf{A}^{i} + d\mathbf{A}^{i}_{gn}   \text{ or } \mathbf{A}^{i+1} = \mathbf{A}^{i} + d\mathbf{A}^{i}_{lm} \  .
\end{equation*}
As at each iteration, $d\mathbf{A}^{i}_{gn}$ and $d\mathbf{A}^{i}_{lm}$ belong to  $\mathcal{T}_{\mathbf{A}^{i}} \text{Gr}(p,k)$, we have
\begin{equation*}
(\mathbf{A}^{i} )^{T} d\mathbf{A}^{i}_{gn} = \mathbf{0}^{ k \times k  } \text{ and } (\mathbf{A}^{i} )^{T} d\mathbf{A}^{i}_{lm} =  \mathbf{0}^{ k \times k  } \ ,
\end{equation*}
and, in these conditions, $\mathbf{A}^{i+1}$ is always of full column-rank and so $\mathbf{A}^{i+1}  \in  \mathbb{R}^{ p \times k  }_{k}$, e.g., $\mathring{\mathbf{A}}^{i+1} = \emph{ran}( \mathbf{A}^{i+1} )$ belongs to the Grassmann manifold $\text{Gr}(p,k)$, validating our claim about the nature of these Gauss-Newton and Levenberg-Marquardt algorithms.

Alternatively, if we require that each element of $\text{Gr}(p,k)$ must be represented by an element of the Stiefel manifold $\text{St}(p,k) = \mathbb{O}^{p \times k }$
as in~\cite{BA2011}\cite{BA2015}, it is necessary to perform an additional retraction step to the correct (Stiefel) manifold at each iteration of the above variable projection Gauss-Newton and Levenberg-Marquardt algorithms in order to consider these algorithms as Riemannian optimization methods operating on the Grassmann manifold $\text{Gr}(p,k)$~\cite{AMS2008}\cite{BA2015}\cite{B2023}. In general terms, a retraction on a manifold can be interpreted as a tool that transforms a tangent update vector at a point of this manifold into a new iterate on this manifold. In other words, at the $(i+1)^{th}$ iteration, in order to move from $\mathbf{O}^{i}  \in \text{St}(p,k) = \mathbb{O}^{p \times k }$ along the tangent vectors $d\mathbf{O}^{i}_{gn}$ or $d\mathbf{O}^{i}_{lm}  \in \mathcal{T}_{\mathbf{O}^{i}} \text{Gr}(p,k)$ while remaining on the Stiefel manifold, after computing
\begin{equation*}
\mathbf{A}^{i+1} = \mathbf{O}^{i} + d\mathbf{O}^{i}_{gn}   \text{ or } \mathbf{A}^{i+1} = \mathbf{O}^{i} + d\mathbf{O}^{i}_{lm} \ ,
\end{equation*}
we need to perform the retraction
\begin{equation} \label{eq:D_retraction} 
\text{Retraction}_{\mathbf{O}^{i} } (  d\mathbf{O}^{i}_{gn}  ) =  \text{Ortho} (\mathbf{O}^{i} + d\mathbf{O}^{i}_{gn} ) \text{ or } \text{Retraction}_{\mathbf{O}^{i} } (  d\mathbf{O}^{i}_{lm}  ) =  \text{Ortho} (\mathbf{O}^{i} + d\mathbf{O}^{i}_{lm} ) \ ,
\end{equation}
where $\text{Ortho}( \mathbf{H} ) \in \text{St}(p,k) = \mathbb{O}^{p \times k }$ designates the $p \times k$ orthonormal factor of a thin QR or polar decomposition of $\mathbf{H} \in \mathbb{R}^{ p \times k  }$; see~\cite{AMS2008}\cite{BA2011}\cite{AM2012}\cite{BA2015}\cite{B2023} for more details on the concept of retraction on manifolds and how these retractions can be computed and used as cheap ways of moving on a specific manifold in Riemannian optimization algorithms.

Interestingly, these results are also very similar to those concerning the algorithms derived in Edelman et al.~\cite{EAS1998} and Manton et al.~\cite{MMH2003} for minimizing the cost function $\psi^{**}(.)$ on the Grassmann manifold $\text{Gr}(p,p-k)$, defined in equation~\eqref{eq:psi**_func} and discussed in Remark~\ref{remark3.7:box}.

More generally, for $\mathbf{W} \in  \mathbb{R}^{ p \times n  }_{+}$,  $\mathbf{A} \in \mathbb{R}^{p \times k }_k$  and any iterative NLLS algorithms use to solve the~\eqref{eq:VP1} or~\eqref{eq:VP2} problems, it is possible to demonstrate that there is almost no loss of generality to restrict the search directions during the iterations to the subspace $\emph{ran}( \mathbf{A} )^{\bot} = \emph{ran}( \mathbf{O}^{\bot} )$ when minimizing the cost function $\psi(.)$, or to $\emph{ran}( \mathbf{A} ) = \emph{ran}( \mathbf{O} )$ when minimizing the cost function $\psi^{**}(.)$ as noted, respectively, in~\cite{C2008b} and~\cite{EAS1998}\cite{MMH2003}.
As an illustration, consider the minimization of $\psi(.)$ in the~\eqref{eq:VP1} problem and take an arbitrary matrix $\mathbf{A} \in \mathbb{R}^{p \times k }_k$ and an arbitrary perturbation matrix $d\mathbf{A} \in \mathbb{R}^{p \times k }$. With the orthonormal basis $\begin{bmatrix} \mathbf{O}   & \mathbf{O}^{\bot} \end{bmatrix}$ of $\mathbb{R}^{p}$ derived from $\mathbf{A}$ in Corollary~\ref{corol5.6:box}, where $\mathbf{O} \in  \mathbb{O}^{p \times k }$ and $\mathbf{O}^{\bot} \in  \mathbb{O}^{p \times (p - k) }$, the perturbation matrix $d\mathbf{A}$ can be decomposed uniquely as  \begin{equation*}
d\mathbf{A} = \mathbf{O}\mathbf{K} + \mathbf{O}^{\bot}\mathbf{K}^{\bot} \ ,
\end{equation*}
where $\mathbf{K} \in \mathbb{R}^{k \times k }$ and $\mathbf{K}^{\bot} \in \mathbb{R}^{(p - k) \times k }$, and we have the equalities
\begin{equation*}
\mathbf{A} + d\mathbf{A} = ( \mathbf{A} + \mathbf{O}\mathbf{K} ) + \mathbf{O}^{\bot}\mathbf{K}^{\bot} = \mathbf{A}\mathbf{Z} + \mathbf{O}^{\bot}\mathbf{K}^{\bot} \ ,
\end{equation*}
since $\mathbf{A} + \mathbf{O}\mathbf{K} \in \emph{ran}( \mathbf{A} ) = \emph{ran}( \mathbf{O} )$ implies that it exists $\mathbf{Z} \in \mathbb{R}^{k \times k }$ such that $\mathbf{A} + \mathbf{O}\mathbf{K} = \mathbf{A}\mathbf{Z}$. Assuming further that $\mathbf{Z}$ is non-singular, which will be the rule in most practical applications, we have the equalities
\begin{equation*}
\psi( \mathbf{A} + d\mathbf{A} ) = \psi( \mathbf{A}\mathbf{Z} + \mathbf{O}^{\bot}\mathbf{K}^{\bot} ) = \psi \big(  ( \mathbf{A}\mathbf{Z} + \mathbf{O}^{\bot}\mathbf{K}^{\bot} ) \mathbf{Z}^{-1} \big) = \psi \big( \mathbf{A} + \mathbf{O}^{\bot}\mathbf{K}^{\bot} \mathbf{Z}^{-1} \big)  \ ,
\end{equation*}
where the second equality results from Corollary~\ref{corol3.2:box}. In other words, as noted by Chen~\cite{C2008b}, for most perturbation matrices $d\mathbf{A}$ around $\mathbf{A}$, there exists a perturbation matrix $d\mathbf{A}^{'} = \mathbf{O}^{\bot}\mathbf{K}^{\bot} \mathbf{Z}^{-1}$ whose range is included in $\emph{ran}( \mathbf{A} )^{\bot}$ and which has exactly the same effect as $d\mathbf{A}$ since $\psi( \mathbf{A} + d\mathbf{A} ) = \psi( \mathbf{A} + d\mathbf{A}^{'} )$. Thus, for any iterative NLLS algorithm used to minimize $\psi(.)$, at each iteration, it suffices generally to search for a perturbation matrix $d\mathbf{A}^{'}$ such that $\emph{ran}( d\mathbf{A}^{'} ) \subset \emph{ran}( \mathbf{A} )^{\bot}$.

Since, in these conditions, $\exists \mathbf{T} \in \mathbb{R}^{(p - k) \times k }$ such that $d\mathbf{A}^{'} = \mathbf{O}^{\bot}\mathbf{T}$, this reduces the number of degrees of freedom, or equivalently the number of parameters, to estimate from $k.p$ to $k.(p - k )$. Moreover, as $\mathbf{A}^{T} d\mathbf{A}^{'} = \mathbf{0}^{k \times k }$, we are also sure that $\mathbf{A} + d\mathbf{A}^{'}$ is always of full column-rank (e.g., $\emph{rank}( \mathbf{A} + d\mathbf{A}^{'} ) = k$ ) across the iterations, which is required for the validity of the algorithms as otherwise the associated orthogonal projector $\mathbf{P}^{\bot}_{\mathbf{F}(.) }$ is not differentiable at $\mathbf{a}$. This is also a useful benefit of restricting the domain of the perturbation matrices $d\mathbf{A}^{'}$ to $\emph{ran}( \mathbf{A} )^{\bot}$ (e.g., such that $\emph{ran}( d\mathbf{A}^{'} ) \subset \emph{ran}( \mathbf{A} )^{\bot}$).

These different properties further justify all the variations of the Gauss-Newton and Levenberg-Marquardt algorithms already described in this subsection, which  restrict the search directions to the subspace $\emph{ran}( \mathbf{A} )^{\bot}$. Proceeding with similar arguments, Edelman et al.~\cite{EAS1998} and Manton et al.~\cite{MMH2003} have also demonstrated that it suffices to consider perturbation matrices whose range is included in $\emph{ran}( \mathbf{O} )$ at each iteration of any NLLS algorithm used to minimize the cost function $\psi^{**}(.)$ in the~\eqref{eq:VP2} problem.
More precisely, in our notations, these algorithms try to minimize the functional  $\psi^{**}(\mathbf{O}^{\bot})$ for $\mathbf{O}^{\bot} \in \mathbb{O}^{p \times ( p - k ) }$ and at each iteration of these algorithms, the search directions for updating $\mathbf{O}^{\bot}$ are restricted to perturbation matrices of the form $d\mathbf{O}^{\bot} = \mathbf{O}\mathbf{T}$ where $\mathbf{O} \in \mathbb{O}^{p \times k }$ with $\mathbf{O}^{T}\mathbf{O}^{\bot} = \mathbf{0}^{k \times ( p - k ) }$ and $\mathbf{T}  \in  \mathbb{R}^{k \times (p - k)}$.
Obviously, as for the minimization of $\psi(.)$ in the~\eqref{eq:VP1} problem, this reduces the dimension of the problem from $p.(p - k )$ to $k.(p - k )$ since  $\mathbf{T}  \in  \mathbb{R}^{k \times (p - k)}$. This confirms that the~\eqref{eq:VP1} and~\eqref{eq:VP2} formulations of the WLRA problem are dual of each other and should have a similar performance, but not necessarily the same cost depending on the values of $p, n$ and $k$.

In the previous linear algebra theorems and corollaries, the matrix $ \mathbf{A} $ is never assumed to have full column-rank; in other words, the rank of $ \mathbf{A} $, $ k_{ \mathbf{A} }$, is a free parameter and, consequently, it can be less than $k$. However, we also recall that the case $ k_{ \mathbf{A} } < k $ is an anomaly in the framework of the WLRA problem because the condition $ \emph{rank}( \mathbf{A} ) = k $ is a necessary condition for the continuity and differentiability of the orthogonal projector  $\mathbf{P}_{\mathbf{F}(.)}$ as demonstrated in Theorems~\ref{theo3.11:box},~\ref{theo3.12:box} and Corollaries~\ref{corol5.1:box},~\ref{corol5.2:box} at the beginning of this subsection. Thus, for the following theorem we will make the natural assumption that $ \mathbf{A} $ has full column-rank in order to derive stronger results. This theorem demonstrates that it is easy to localize the $ k.k $ linearly dependent columns in $\mathbf{M}(\mathbf{a})$, $\mathbf{L}(\mathbf{a})$ and $\mathit{J}( \mathbf{r}(\mathbf{a}) )$ in almost all practical cases if $ k_{ \mathbf{A} } = k $. This result is also new as far we know.
\\
 \begin{theo5.4} \label{theo5.4:box}
With the same notations as in Theorem~\ref{theo5.2:box}, if $ \mathbf{A} $,  $ \mathbf{M}(\mathbf{a}) $, $ \mathbf{L}(\mathbf{a}) $ and $ \mathit{J}( \mathbf{r}(\mathbf{a}) ) $ are partitioned, respectively, as
\begin{IEEEeqnarray*}{llll}
    \mathbf{A} =  \begin{bmatrix}
            \mathbf{A}_{1}  \\
            \mathbf{A}_{2}
        \end{bmatrix}  \thickspace
      &
   \text{ with } \mathbf{A}_{1} \in \mathbb{R}^{(p-k) \times k}    & \text{ and } & \mathbf{A}_{2} \in \mathbb{R}^{k \times k}  \ ,  \\
     \mathbf{M}(\mathbf{a}) = \begin{bmatrix}
            \mathbf{M}(\mathbf{a})_{1} & \mathbf{M}(\mathbf{a})_{2}
       \end{bmatrix}  \thickspace
       &
   \text{ with }  \mathbf{M}(\mathbf{a})_{1} \in \mathbb{R}^{n.p \times k.(p-k)}   & \text{ and }  & \mathbf{M}(\mathbf{a})_{2} \in \mathbb{R}^{n.p \times k.k} \ ,  \\
    \mathbf{L}(\mathbf{a}) = \begin{bmatrix}
            \mathbf{L}(\mathbf{a})_{1} &  \mathbf{L}(\mathbf{a})_{2}
        \end{bmatrix}  \thickspace
        &
   \text{ with } \mathbf{L}(\mathbf{a})_{1} \in \mathbb{R}^{n.p \times k.(p-k)}     & \text{ and }  & \mathbf{L}(\mathbf{a})_{2} \in \mathbb{R}^{n.p \times k.k}  \  , \\
    \mathit{J}( \mathbf{r}(\mathbf{a}) ) = \begin{bmatrix}
            \mathit{J}( \mathbf{r}(\mathbf{a}) )_{1} & \mathit{J}( \mathbf{r}(\mathbf{a}) )_{2}
        \end{bmatrix} \thickspace
         &
    \text{ with }  \mathit{J}( \mathbf{r}(\mathbf{a}) )_{1} \in \mathbb{R}^{n.p \times k.(p-k)}   & \text{ and }  & \mathit{J}( \mathbf{r}(\mathbf{a}) )_{2} \in \mathbb{R}^{n.p \times k.k} \ ,
\end{IEEEeqnarray*} 
and if $ \emph{rank}( \mathbf{A}_{2} ) = k $, then $ \exists  \mathbf{Z} \in \mathbb{R}^{k.(p-k) \times k.k } $ such that
\begin{equation*}
    \mathbf{M}(\mathbf{a})_{2} = \mathbf{M}(\mathbf{a})_{1} \mathbf{Z} \text{ , }  \mathbf{L}(\mathbf{a})_{2} = \mathbf{L}(\mathbf{a})_{1} \mathbf{Z} \text{ and } \mathit{J}( \mathbf{r}(\mathbf{a}) )_{2} = \mathit{J}( \mathbf{r}(\mathbf{a}) )_{1} \mathbf{Z} \ .
\end{equation*}
In other words, if $ \emph{rank}( \mathbf{A}_{2} ) = k $, the last $k.k$ columns of $\mathbf{M}(\mathbf{a})$, $\mathbf{L}(\mathbf{a})$ and $\mathit{J}( \mathbf{r}(\mathbf{a}) )$ are linearly dependent upon the first $ k.(p-k) $ columns of $\mathbf{M}(\mathbf{a})$, $\mathbf{L}(\mathbf{a})$ and $\mathit{J}( \mathbf{r}(\mathbf{a}) )$, respectively.
\end{theo5.4}
\begin{proof}
We will first demonstrate that $ \mathbf{M}(\mathbf{a})_{2} = \mathbf{M}(\mathbf{a})_{1} \mathbf{Z} $ for some $ \mathbf{Z} \in \mathbb{R}^{k.(p-k) \times k.k } $. To this end, let us derive an explicit expression for the $\mathbf{M}(\mathbf{a})_{1}$ and $\mathbf{M}(\mathbf{a})_{2}$ submatrices using the formulation of the $\mathbf{M}(\mathbf{a})$ matrix given by equation~\eqref{eq:M_mat}, e.g.,
\begin{equation*}
\mathbf{M}(\mathbf{a}) = \mathbf{P}_{\mathbf{F}(\mathbf{a})}^{\bot} \emph{diag}(\emph{vec}(\sqrt{\mathbf{W}}))(\mathbf{\widehat{B}}^{T} \otimes \mathbf{I}_{p}) \mathbf{K}_{(k,p)} \ .
\end{equation*}
Using equation~\eqref{eq:com_kron}, we have
\begin{equation*}
   (  \mathbf{\widehat{B}}^{T} \otimes \mathbf{I}_{p} ) \mathbf{K}_{(k,p)} = \mathbf{K}_{(n,p)}( \mathbf{I}_{p} \otimes \mathbf{\widehat{B}}^{T} ) = \mathbf{K}_{(n,p)}
   \begin{bmatrix}
         \mathbf{I}_{p-k} \otimes \mathbf{\widehat{B}}^{T}  &  \mathbf{0}^{n.(p-k) \times k.k }   \\
         \mathbf{0}^{k.n  \times n.(p-k)}  & \mathbf{I}_{k} \otimes \mathbf{\widehat{B}}^{T}
    \end{bmatrix} \ ,
\end{equation*}
hence
\begin{equation*}
   \mathbf{M}(\mathbf{a}) =  \mathbf{P}_{\mathbf{F}(\mathbf{a})}^{\bot} \emph{diag}(\emph{vec}(\sqrt{\mathbf{W}})) \mathbf{K}_{(n,p)}
   \begin{bmatrix}
         \mathbf{I}_{p-k} \otimes \mathbf{\widehat{B}}^{T}  &  \mathbf{0}^{n.(p-k) \times k.k }  \\
         \mathbf{0}^{k.n  \times n.(p-k)}  & \mathbf{I}_{k} \otimes \mathbf{\widehat{B}}^{T}
    \end{bmatrix}\  ,
\end{equation*}
and so
\begin{align*}
   \mathbf{M}(\mathbf{a})_{1} & =   \mathbf{P}_{\mathbf{F}(\mathbf{a})}^{\bot} \emph{diag}(\emph{vec}(\sqrt{\mathbf{W}})) \mathbf{K}_{(n,p)}
   \begin{bmatrix}
         \mathbf{I}_{p-k} \otimes \mathbf{\widehat{B}}^{T}    \\
          \mathbf{0}^{k.n  \times n.(p-k)}
    \end{bmatrix} \ ,   \\
   \mathbf{M}(\mathbf{a})_{2} & =   \mathbf{P}_{\mathbf{F}(\mathbf{a})}^{\bot} \emph{diag}(\emph{vec}(\sqrt{\mathbf{W}})) \mathbf{K}_{(n,p)}
   \begin{bmatrix}
         \mathbf{0}^{n.(p-k) \times k.k }  \\
         \mathbf{I}_{k} \otimes \mathbf{\widehat{B}}^{T}
    \end{bmatrix} \ .
\end{align*}
Now, we want to show that, $\forall j \in \{1,2, \cdots ,k.k\}$, there is $ \mathbf{z}_{j} \in \mathbb{R}^{k.(p-k)} $ such that $ \mathbf{M}(\mathbf{a})_{2} \mathbf{e}_{j} = \mathbf{M}(\mathbf{a})_{1} \mathbf{z}_{j} $, where $ \mathbf{e}_{j} $  is the $j^{th}$ column unit vector of order $k.k$ and $\mathbf{z}_{j}$ is the $j^{th}$ column of the matrix $\mathbf{Z}$ we are looking for. Using the preceding expressions of the $\mathbf{M}(\mathbf{a})_{1}$ and $\mathbf{M}(\mathbf{a})_{2}$ submatrices, $\forall j \in \{1,2, \cdots ,k.k\} $, we have the equivalences
\begin{align*}    
     &  \mathbf{M}(\mathbf{a})_{2} \mathbf{e}_{j} = \mathbf{M}(\mathbf{a})_{1} \mathbf{z}_{j}  \\   
     & \Leftrightarrow   \mathbf{P}_{\mathbf{F}(\mathbf{a})}^{\bot} \emph{diag}(\emph{vec}(\sqrt{\mathbf{W}})) \mathbf{K}_{(n,p)}
     \Bigg(
        \begin{bmatrix}
         \mathbf{0}^{n.(p-k) \times k.k }  \\
         \mathbf{I}_{k} \otimes \mathbf{\widehat{B}}^{T}
    \end{bmatrix} \mathbf{e}_{j}  - 
       \begin{bmatrix}
         \mathbf{I}_{p-k} \otimes \mathbf{\widehat{B}}^{T}    \\
         \mathbf{0}^{k.n  \times n.(p-k)}
    \end{bmatrix} \mathbf{z}_{j}  \Bigg) =  \mathbf{0}^{p.n}   \\
    & \Leftrightarrow   \mathbf{P}_{\mathbf{F}(\mathbf{a})}^{\bot} \emph{diag}(\emph{vec}(\sqrt{\mathbf{W}})) \mathbf{K}_{(n,p)}
        \begin{bmatrix}
         - ( \mathbf{I}_{p-k} \otimes \mathbf{\widehat{B}}^{T}  )  \mathbf{z}_{j} \\
            ( \mathbf{I}_{k} \otimes \mathbf{\widehat{B}}^{T}     )  \mathbf{e}_{j}
        \end{bmatrix}  = \mathbf{0}^{p.n}    \\
    & \Leftrightarrow   \emph{diag}(\emph{vec}(\sqrt{\mathbf{W}})) \mathbf{K}_{(n,p)}
        \begin{bmatrix}
         - ( \mathbf{I}_{p-k} \otimes \mathbf{\widehat{B}}^{T}  )  \mathbf{z}_{j} \\
            ( \mathbf{I}_{k} \otimes \mathbf{\widehat{B}}^{T}     )  \mathbf{e}_{j}
        \end{bmatrix}    \in \emph{ran}( \mathbf{F}(\mathbf{a}) )   \\
      & \Leftrightarrow   \exists  \mathbf{t}_{j}  \in \mathbb{R}^{n.k} \text{ / }  \emph{diag}(\emph{vec}(\sqrt{\mathbf{W}})) \mathbf{K}_{(n,p)}
        \begin{bmatrix}
         - ( \mathbf{I}_{p-k} \otimes \mathbf{\widehat{B}}^{T}  )  \mathbf{z}_{j} \\
            ( \mathbf{I}_{k} \otimes \mathbf{\widehat{B}}^{T}     )  \mathbf{e}_{j}
        \end{bmatrix}   = \mathbf{F}(\mathbf{a}) \mathbf{t}_{j} \ .
\end{align*}
For all $ j \in \{1,2, \cdots ,k.k\} $, we will look for a vector $ \mathbf{t}_{j} $ such that $ \mathbf{t}_{j} =  \emph{vec}(\mathbf{C}_{j} \mathbf{\widehat{B}})$ where $ \mathbf{C}_{j}  \in \mathbb{R}^{k \times k} $. But, using equations~\eqref{eq:vec_kronprod},~\eqref{eq:com_kron} and~\eqref{eq:commat}, $ \forall \mathbf{C}  \in \mathbb{R}^{k \times k} $, we have
\begin{align*}
    \mathbf{F}(\mathbf{a}) \emph{vec}(\mathbf{C\widehat{B}}) & =  \emph{diag}(\emph{vec}(\sqrt{\mathbf{W}})) ( \mathbf{I}_{n} \otimes \mathbf{A} ) \emph{vec}(\mathbf{C\widehat{B}})  \\
                                                                                          & =  \emph{diag}(\emph{vec}(\sqrt{\mathbf{W}})) \emph{vec}(\mathbf{AC\widehat{B}})  \\
                                                                                          & =  \emph{diag}(\emph{vec}(\sqrt{\mathbf{W}})) (  \mathbf{\widehat{B}}^{T}  \otimes \mathbf{I}_{p} )  \emph{vec}(\mathbf{AC})  \\
                                                                                          & =  \emph{diag}(\emph{vec}(\sqrt{\mathbf{W}})) \mathbf{K}_{(n,p)} (  \mathbf{I}_{p} \otimes  \mathbf{\widehat{B}}^{T}  ) \mathbf{K}_{(p,k)}  \emph{vec}(\mathbf{AC})  \\
                                                                                          & =  \emph{diag}(\emph{vec}(\sqrt{\mathbf{W}})) \mathbf{K}_{(n,p)} (  \mathbf{I}_{p} \otimes  \mathbf{\widehat{B}}^{T}  ) \emph{vec}( \mathbf{C}^{T} \mathbf{A}^{T} )  \\
                                                                                          & =  \emph{diag}(\emph{vec}(\sqrt{\mathbf{W}})) \mathbf{K}_{(n,p)}
                                                                                                        \begin{bmatrix}
                                                                                                            \mathbf{I}_{p-k} \otimes \mathbf{\widehat{B}}^{T}  &  \mathbf{0}^{n(p-k) \times k.k }  \\
                                                                                                             \mathbf{0}^{k.n  \times n(p-k)}  & \mathbf{I}_{k} \otimes \mathbf{\widehat{B}}^{T}
                                                                                                        \end{bmatrix}  \emph{vec}( \mathbf{C}^{T} \mathbf{A}^{T} )  \\
                                                                                          & =  \emph{diag}(\emph{vec}(\sqrt{\mathbf{W}})) \mathbf{K}_{(n,p)}
                                                                                                        \begin{bmatrix}
                                                                                                           (  \mathbf{I}_{p-k} \otimes \mathbf{\widehat{B}}^{T} )   \emph{vec}( \mathbf{C}^{T} \mathbf{A}^{T}_{1} )  \\
                                                                                                           (  \mathbf{I}_{k} \otimes \mathbf{\widehat{B}}^{T} )      \emph{vec}( \mathbf{C}^{T} \mathbf{A}^{T}_{2} )
                                                                                                        \end{bmatrix} \ ,
\end{align*}
and it is therefore sufficient to demonstrate that, $\forall j \in \{1,2, \cdots ,k.k\}$, there are $ \mathbf{z}_{j} \in \mathbb{R}^{k.(p-k)} $ and $ \mathbf{C}_{j}  \in \mathbb{R}^{k \times k} $ such that
\begin{equation*}
         \begin{bmatrix}
         - ( \mathbf{I}_{p-k} \otimes \mathbf{\widehat{B}}^{T}  )  \mathbf{z}_{j} \\
            ( \mathbf{I}_{k} \otimes \mathbf{\widehat{B}}^{T}     )  \mathbf{e}_{j}
        \end{bmatrix} =
        \begin{bmatrix}
            (  \mathbf{I}_{p-k} \otimes \mathbf{\widehat{B}}^{T} )   \emph{vec}( \mathbf{C}^{T}_{j} \mathbf{A}^{T}_{1} )  \\
            (  \mathbf{I}_{k} \otimes \mathbf{\widehat{B}}^{T} )      \emph{vec}( \mathbf{C}^{T} _{j}\mathbf{A}^{T}_{2} )
        \end{bmatrix}
\end{equation*}
in order to have $ \mathbf{M}(\mathbf{a})_{2} \mathbf{e}_{j} = \mathbf{M}(\mathbf{a})_{1} \mathbf{z}_{j}  $. Furthermore, if $ \mathbf{z}_{j} $ is set to $ -\emph{vec}( \mathbf{C}^{T}_{j} \mathbf{A}^{T}_{1} ) $ in the above equation, then it suffices to show that there is $ \mathbf{C}_{j}  \in \mathbb{R}^{k \times k} $ such that
\begin{equation*}
    \mathbf{e}_{j} = \emph{vec}( \mathbf{C}^{T} _{j}\mathbf{A}^{T}_{2} )
\end{equation*}
in order to obtain the desired result. But, $ \forall j \in \{1,2, \cdots ,k.k\} $, it is easily verified that there are two integers $t(j)$ and $u(j)$ such that
\begin{equation*}
    \mathbf{e}_{j} = \emph{vec}( \mathbf{i} _{t(j)} \mathbf{i}^{T}_{u(j)} ) \ ,
\end{equation*}
where $\mathbf{i} _{t}$ is the $t^{th}$ column unit vector of order $k$, and
\begin{align*}
    \mathbf{e}_{j} = \emph{vec}( \mathbf{C}^{T} _{j}\mathbf{A}^{T}_{2} )  & \Leftrightarrow  \mathbf{i} _{u(j)} \mathbf{i}^{T}_{t(j)} =  \mathbf{A}_{2}  \mathbf{C}_{j} \\
                                                                                                                & \Leftrightarrow  \mathbf{C}_{j} = \mathbf{A}^{-1}_{2} \mathbf{i} _{u(j)} \mathbf{i}^{T}_{t(j)} \ ,
\end{align*}
since $\mathbf{A}_{2}$ is nonsingular by hypothesis. Consequently, $ \forall j \in \{1,2, \cdots ,k.k\} $, we have $ \mathbf{M}(\mathbf{a})_{2} \mathbf{e}_{j} = \mathbf{M}(\mathbf{a})_{1} \mathbf{z}_{j}  $ with
\begin{equation*}
    \mathbf{z}_{j} =  - \emph{vec}(   \mathbf{C}^{T}_{j}   \mathbf{A}^{T}_{1} ) = - \emph{vec}( \mathbf{i} _{t(j)} \mathbf{i}^{T}_{u(j)}  \big( \mathbf{A}^{-1}_{2})^{T} \mathbf{A}^{T}_{1} \big) \ ,
\end{equation*}
which is the desired result.

We now demonstrate that, $ \forall j \in \{1,2, \cdots ,k.k\} $, we also have $ \mathbf{L}(\mathbf{a})_{2} \mathbf{e}_{j} = \mathbf{L}(\mathbf{a})_{1} \mathbf{z}_{j} $ with $ \mathbf{z}_{j} $ defined as above. To this end, we first observe that if $\mathbf{W}$ and $P_{\Omega}(\mathbf{X} -\mathbf{A}\mathbf{\widehat{B}})$ are partitioned as
\begin{equation*}
    \mathbf{W} = \begin{bmatrix}
            \mathbf{W}_{1}  \\
            \mathbf{W}_{2}
        \end{bmatrix} \text{ and }
    P_{\Omega}(\mathbf{X} -\mathbf{A}\mathbf{\widehat{B}})  = \begin{bmatrix}
            \mathbf{P}_{1}  \\
            \mathbf{P}_{2}
        \end{bmatrix} 
\end{equation*} \ ,
with
\begin{equation*}
    \mathbf{W}_{1}, \mathbf{P}_{1} \in \mathbb{R}^{(p-k) \times n}   \quad  \text{and}   \quad \mathbf{W}_{2}, \mathbf{P}_{2} \in \mathbb{R}^{k \times n} \ ,
\end{equation*} 
then, using equation~\eqref{eq:L_mat} in this subsection and equation~\eqref{eq:partmat_kronprod} in Subsection~\ref{multlin_alg:box},
\begin{align*}
    \mathbf{L}(\mathbf{a})  &   =  ( \mathbf{F}(\mathbf{a})^{+} )^{T} \big( ( \mathbf{W} \odot P_{\Omega}(\mathbf{X} -\mathbf{A}\mathbf{\widehat{B}}) )^{T} \otimes \mathbf{I}_{k} \big)  \\
                                         &   =  ( \mathbf{F}(\mathbf{a})^{+} )^{T} \big\lbrack ( \mathbf{W}_{1} \odot \mathbf{P}_{1} )^{T} \otimes \mathbf{I}_{k} , ( \mathbf{W}_{2} \odot \mathbf{P}_{2} )^{T} \otimes \mathbf{I}_{k} \big\rbrack \ ;
\end{align*} 
hence $\mathbf{L}(\mathbf{a})_{1}$ and $\mathbf{L}(\mathbf{a})_{2}$ are given by
\begin{equation*}
    \mathbf{L}(\mathbf{a})_{1} = ( \mathbf{F}(\mathbf{a})^{+} )^{T} \big( ( \mathbf{W}_{1} \odot \mathbf{P}_{1} )^{T} \otimes \mathbf{I}_{k} \big) \text{ and } \mathbf{L}(\mathbf{a})_{2} =  ( \mathbf{F}(\mathbf{a})^{+} )^{T} \big( (\mathbf{W}_{2} \odot \mathbf{P}_{2} )^{T} \otimes \mathbf{I}_{k} \big)  \ ,
\end{equation*}
 $\forall j \in \{1,2, \cdots ,k.k\} $, we then have the implication
\begin{equation*}
     \big( (\mathbf{W}_{2} \odot \mathbf{P}_{2} )^{T} \otimes \mathbf{I}_{k} \big) \mathbf{e}_{j} = \big( ( \mathbf{W}_{1} \odot \mathbf{P}_{1} )^{T} \otimes \mathbf{I}_{k} \big) \mathbf{z}_{j}  \Rightarrow   \mathbf{L}(\mathbf{a})_{2} \mathbf{e}_{j} = \mathbf{L}(\mathbf{a})_{1} \mathbf{z}_{j}  \ ,
\end{equation*}
and it suffices to show that
\begin{equation*}
\big( (\mathbf{W}_{2} \odot \mathbf{P}_{2} )^{T} \otimes \mathbf{I}_{k} \big) \mathbf{e}_{j} - \big( ( \mathbf{W}_{1} \odot \mathbf{P}_{1} )^{T} \otimes \mathbf{I}_{k} \big) \mathbf{z}_{j} = \mathbf{0}^{p.k} 
\end{equation*}
in order to obtain $\mathbf{L}(\mathbf{a})_{2} \mathbf{e}_{j} = \mathbf{L}(\mathbf{a})_{1} \mathbf{z}_{j}$. Defining, as above,
\begin{equation*}
    \mathbf{C}_{j} = \mathbf{A}^{-1}_{2} \mathbf{i} _{u(j)} \mathbf{i}^{T}_{t(j)}
\end{equation*}
and remembering that
\begin{equation*}
      \mathbf{z}_{j} =  - \emph{vec}(   \mathbf{C}^{T}_{j}   \mathbf{A}^{T}_{1} ) \text{ and } \mathbf{e}_{j} = \emph{vec}( \mathbf{C}^{T} _{j}\mathbf{A}^{T}_{2} )  \ ,
\end{equation*}
we have (using equation~\eqref{eq:vec_kronprod} in Subsection~\ref{multlin_alg:box})
\begin{align*}
    \big( ( \mathbf{W}_{1} \odot \mathbf{P}_{1} )^{T} \otimes \mathbf{I}_{k} \big) \mathbf{z}_{j}  &  =  - \big( ( \mathbf{W}_{1} \odot \mathbf{P}_{1} )^{T} \otimes \mathbf{I}_{k} \big)  \emph{vec}( \mathbf{C}^{T}_{j} \mathbf{A}^{T}_{1} )     \\
                                                                                                                                                  &  =  - \emph{vec}\big( \mathbf{C}^{T}_{j} \mathbf{A}^{T}_{1}  ( \mathbf{W}_{1} \odot \mathbf{P}_{1} ) \big)                                                     \\
                                                                                                                                                  &  =  - \Big( \big( ( \mathbf{W}_{1} \odot \mathbf{P}_{1} )^{T} \mathbf{A}_{1} \big) \otimes \mathbf{I}_{k} \Big) \emph{vec}( \mathbf{C}^{T}_{j} ) \ ,
\end{align*}
and also
\begin{align*}
    \big( ( \mathbf{W}_{2} \odot \mathbf{P}_{2} )^{T} \otimes \mathbf{I}_{k} \big) \mathbf{e}_{j}  &  =  \big( ( \mathbf{W}_{2} \odot \mathbf{P}_{2} )^{T} \otimes \mathbf{I}_{k} \big)  \emph{vec}( \mathbf{C}^{T}_{j} \mathbf{A}^{T}_{2} )     \\
                                                                                                                                                  &  = \emph{vec}\big( \mathbf{C}^{T}_{j} \mathbf{A}^{T}_{2}  ( \mathbf{W}_{2} \odot \mathbf{P}_{2} ) \big)                                                     \\
                                                                                                                                                  &  =  \Big( \big( ( \mathbf{W}_{2} \odot \mathbf{P}_{2} )^{T} \mathbf{A}_{2} \big) \otimes \mathbf{I}_{k} \Big) \emph{vec}( \mathbf{C}^{T}_{j} ) \ .
\end{align*}
From these equalities, using equations~\eqref{eq:partmat_kronprod},~\eqref{eq:mulmat_kronprod}  and~\eqref{eq:vec_kronprod} in Subsection~\ref{multlin_alg:box}, we deduce that
\begin{align*}
&    \big( (\mathbf{W}_{2} \odot \mathbf{P}_{2} )^{T} \otimes \mathbf{I}_{k} \big) \mathbf{e}_{j} - \big( ( \mathbf{W}_{1} \odot \mathbf{P}_{1} )^{T} \otimes \mathbf{I}_{k} \big) \mathbf{z}_{j}    \\
          &   \qquad =   \Big( \big( ( \mathbf{W}_{2} \odot \mathbf{P}_{2} )^{T} \mathbf{A}_{2}  \big) \otimes \mathbf{I}_{k} + \big( ( \mathbf{W}_{1} \odot \mathbf{P}_{1} )^{T} \mathbf{A}_{1}  \big) \otimes \mathbf{I}_{k} \Big) \emph{vec}( \mathbf{C}^{T}_{j} )   \\
          &   \qquad =   \Big( \big( ( \mathbf{W}_{2} \odot \mathbf{P}_{2} )^{T} \mathbf{A}_{2} +  ( \mathbf{W}_{1} \odot \mathbf{P}_{1} )^{T} \mathbf{A}_{1} \big) \otimes \mathbf{I}_{k}  \Big) \emph{vec}( \mathbf{C}^{T}_{j} )   \\
          &  \qquad =  \Big( \big( ( \mathbf{W} \odot P_{\Omega}(\mathbf{X} -\mathbf{A}\mathbf{\widehat{B}}) )^{T} \mathbf{A}  \big)  \otimes \mathbf{I}_{k}  \Big) \emph{vec}( \mathbf{C}^{T}_{j} )   \\
          &   \qquad =  \emph{vec}\big( \mathbf{C}^{T}_{j} \mathbf{A}^{T}  ( \mathbf{W} \odot P_{\Omega}(\mathbf{X} -\mathbf{A}\mathbf{\widehat{B}}) )  \big)   \\
          &   \qquad =  ( \mathbf{I}_{n}  \otimes \mathbf{C}^{T}_{j} )( \mathbf{I}_{n}  \otimes \mathbf{A}^{T} ) \emph{vec}(  \mathbf{W} \odot P_{\Omega}(\mathbf{X} -\mathbf{A}\mathbf{\widehat{B}})  )  \\
          &   \qquad =  ( \mathbf{I}_{n}  \otimes \mathbf{C}^{T}_{j} )( \mathbf{I}_{n}  \otimes \mathbf{A}^{T} ) \emph{diag}( \emph{vec}(  \sqrt{\mathbf{W}} ) ) \emph{vec}(  \sqrt{\mathbf{W}} \odot P_{\Omega}(\mathbf{X} -\mathbf{A}\mathbf{\widehat{B}})  )  \\
          &    \qquad =  ( \mathbf{I}_{n}  \otimes \mathbf{C}^{T}_{j} ) \mathbf{F}( \mathbf{a} )^{T} \mathbf{r}( \mathbf{a} )   \\
          &   \qquad =  \mathbf{0}^{p.k}  \  ,
\end{align*}
since $ \mathbf{F}( \mathbf{a} )^{T} \mathbf{r}( \mathbf{a} ) = \mathbf{0}^{n.k} $. This implies that, $ \forall j \in \{1,2, \cdots ,k.k\} $, we  also have $ \mathbf{L}(\mathbf{a})_{2} \mathbf{e}_{j} = \mathbf{L}(\mathbf{a})_{1} \mathbf{z}_{j} $ as claimed in the theorem.

Finally, since $  \mathit{J}( \mathbf{r}(\mathbf{a}) ) = - ( \mathbf{M}(\mathbf{a}) + \mathbf{L}(\mathbf{a}) ) $, the preceding results also imply that, $ \forall j \in \{1,2, \cdots ,k.k\} $, we  have
\begin{equation*}
\mathit{J}  \big( \mathbf{r}(\mathbf{a}) \big)_{2} \mathbf{e}_{j} = \mathit{J }\big( \mathbf{r}(\mathbf{a}) \big)_{1} \mathbf{z}_{j} \ .
\end{equation*}
\end{proof}

Theorems~\ref{theo5.4:box} shows that it is possible to compute $d\mathbf{a}_{gn}$ efficiently, and without using the linear constraint $\mathbf{N}^{T} d\mathbf{a} = \mathbf{0}^{k.k}$ or an orthonormal basis $\mathbf{O}^{\bot}$ of $\emph{ran}( \mathbf{A} )^{\bot}$, by using a simplified stable COD of $\mathit{J}( \mathbf{r}(\mathbf{a}) )$ (see equation~\eqref{eq:cod} in Subsection~\ref{lin_alg:box} for details). When $\emph{rank}( \mathit{J}( \mathbf{r}(\mathbf{a}) ) ) = k.( p - k)$, the first $k.(p-k)$ columns of $\mathit{J}( \mathbf{r}(\mathbf{a}) )$ are linearly independent and the last  $k.k$ columns of this matrix are linearly dependent upon these first $k.(p-k)$ columns as soon as $ \emph{rank}( \mathbf{A}_{2} ) = k$. This property is also verified for the matrix $-\mathbf{M}(\mathbf{a})$, which can be used as an approximation of $\mathit{J}( \mathbf{r}(\mathbf{a}) )$ in the optimization algorithms as we will discuss in Section~\ref{vpalg:box}. In these conditions, a simplified COD of $\mathit{J}( \mathbf{r}(\mathbf{a}) )$ (or alternatively of $-\mathbf{M}(\mathbf{a})$) is given by
\begin{equation*}
\mathit{J}( \mathbf{r}(\mathbf{a}) ) = \mathbf{Q}\mathbf{T}\mathbf{Z}^{T} =  \mathbf{Q}
\begin{bmatrix}
\mathbf{T}_{11}  & \mathbf{0}^{k(p-k) \times k.k} 
\end{bmatrix}
\mathbf{Z}^{T} \ ,
\end{equation*}
where $\mathbf{Q} \in \mathbb{O}^{n.p \times k.(p-k)}$, $\mathbf{Z} \in \mathbb{O}^{k.p \times k.p }$, $\mathbf{T} \in \mathbb{R}^{k.(p-k)  \times k.p}$ and $\mathbf{T}_{11} \in \mathbb{R}^{k.(p-k)  \times k.(p-k)}$ is a  full rank upper triangular matrix.
The advantage of this COD formulation for rank deficient matrices, such as $\mathit{J}( \mathbf{r}(\mathbf{a}) )$ or  $-\mathbf{M}(\mathbf{a})$, is the ability to compute directly the unique minimum 2-norm solution of the problem
\begin{equation*}
\min_{d\mathbf{a} \in \mathbb{R}^{p.k}}   \,   \frac{1}{2} \Vert \mathbf{r}(  \mathbf{a} ) + \mathit{J}( \mathbf{r}(\mathbf{a}) ) d\mathbf{a} \Vert^{2}_{2}
\end{equation*}
as follows
\begin{equation*}
d\mathbf{a}_{gn} = - \mathbf{Z}
\begin{bmatrix}
\mathbf{T}^{-1}_{11}  \mathbf{Q}^{T} \mathbf{r}(  \mathbf{a} ) \\ \mathbf{0}^{k.k  \times k(p-k) } 
\end{bmatrix} \ ,
\end{equation*}
without computing the SVD or the pseudo-inverse of $\mathit{J}( \mathbf{r}(\mathbf{a}) )$, but only this simplified COD of $\mathit{J}( \mathbf{r}(\mathbf{a}) )$, which is based on a simple QR factorization without column pivoting of the first $k.(p-k)$ columns of $\mathit{J}( \mathbf{r}(\mathbf{a}) )$ as a first step. Furthermore, as $\mathit{J}( \mathbf{r}(\mathbf{a}) )$ is a tall and skinny matrix, this QR factorization can be parallelized very efficiently to reduce the computational time and the memory requirements as we will illustrate in Section~\ref{vpalg:box}, see also~\cite{DGHL2012} for details. We also note that this COD can also be used to derive another variation of the Levenberg-Marquardt algorithm, which will also take care of the singularity of the Jacobian matrix or its approximation without using the linear constraint $\mathbf{N}^{T} d\mathbf{a} = \mathbf{0}^{k.k}$ or an orthonormal basis, $\mathbf{O}^{\bot}$, of $\emph{ran}( \mathbf{A} )^{\bot}$ as in the  formulations derived above; see again  Section~\ref{vpalg:box} for details.

Theorems~\ref{theo5.2:box} and~\ref{theo5.4:box} are valid for an arbitrary weight matrix $\mathbf{W}$. However, the ranks of $\mathbf{M}(\mathbf{a})$, $\mathbf{L}(\mathbf{a})$ and $\mathit{J}( \mathbf{r}(\mathbf{a}) )$ may also be altered by the choice of a particular weight matrix, as already illustrated by Theorem~\ref{theo5.3:box}. The two following theorems further illustrate the strong dependency of the ranks of $\mathbf{M}(\mathbf{a})$, $\mathbf{L}(\mathbf{a})$ and $\mathit{J}( \mathbf{r}(\mathbf{a}) )$ to the number and distribution of zero elements in $\mathbf{W}$, respectively.

Let us again give some definitions before stating the two theorems. First, we introduce the $p \times n$ incidence matrix, $\boldsymbol{\delta}$, associated with the matrix $\mathbf{W}$, defined by
\begin{equation} \label{eq:incidence_mat}
    \boldsymbol{\delta}_{ij} = 
              \begin{cases}
                    1 & \text{if }  \mathbf{W}_{ij} \ne 0   \\
                    0 & \text{if }  \mathbf{W}_{ij} = 0   \\
              \end{cases} \ .
\end{equation}
The number of "nonmissing" elements in $\mathbf{X}$ (or equivalently the number of non-zero elements in $\mathbf{W}$) is then equal to $ \displaystyle{ \sum_{ij} \boldsymbol{\delta}_{ij} } = nobs $. Since
\begin{equation*}
    \mathbf{F}(\mathbf{a})  =  \bigoplus_{j=1}^n  \left( \emph{diag}( \emph{vec}(  \sqrt{\mathbf{W}}_{.j} ) ) \mathbf{A} \right)  = \bigoplus_{j=1}^n  \mathbf{F}_{j}(\mathbf{a}) \ ,
\end{equation*}
we first observe that
\begin{equation*}
    \emph{rank}\big( \mathbf{F}(\mathbf{a}) \big) = \displaystyle{ \sum_{j=1}^n r_{j} }   \quad \text{where } r_{j} = \emph{rank}( \mathbf{F}_{j}(\mathbf{a}) )  \text{, for } j=1, \cdots, n \ .
\end{equation*}
Furthermore, using equation~\eqref{eq:rank2} in Subsection~\ref{lin_alg:box}, we have
\begin{equation*}
    r_{j} \leqslant \emph{min} \left( \displaystyle{ \sum_{i=1}^p \boldsymbol{\delta}_{ij} } , k \right) \text{, for } j=1, \cdots, n, \quad \text{and} \quad \emph{rank}\big( \mathbf{F}(\mathbf{a}) \big) \leqslant \emph{min} \left( nobs , k.n \right)  \ .
\end{equation*}
We now state the following simple bounds for the ranks of $\mathbf{M}(\mathbf{a})$, $\mathbf{L}(\mathbf{a})$ and $\mathit{J}( \mathbf{r}(\mathbf{a}) )$:
\\
 \begin{theo5.5} \label{theo5.5:box}
With these definitions and the same notations as in Theorem~\ref{theo5.2:box}, then the following inequalities hold:
 \begin{align*}
    \emph{rank}( \mathbf{M}(\mathbf{a}) )                       & \leqslant    \emph{min} \left( nobs - \emph{rank} \big( \mathbf{F}(\mathbf{a}) \big), k.( p - k_{ \mathbf{A} } ) \right)              \\
    \emph{rank}( \mathbf{L}(\mathbf{a}) )                        & \leqslant    \emph{min} \left( nobs , k.\emph{min}( n, p - k_{ \mathbf{A} } ) \right)   \\
    \emph{rank}( \mathit{J}( \mathbf{r}(\mathbf{a}) ) )     & \leqslant    \emph{min} \left( nobs , k.( p - k_{ \mathbf{A} } ) \right)   \ .
\end{align*}
\end{theo5.5}
\begin{proof}
Omitted.

\end{proof}
Let us further define for any integer $ i \in \{1,2, \cdots ,p\} $ and any $ p \times k $ matrix $ \mathbf{A} $, the finite subset of $\mathbb{N}$
 \begin{equation*}
    \xi ( i, \mathbf{A} ) = \Big\{ j \in  \{1,2, \cdots ,n\}  \text{ } / \text{ }  \boldsymbol{\delta}_{ij} \ne 0 \text{ and }   \displaystyle{ \sum_{l=1}^p  \boldsymbol{\delta}_{lj} } >  r_{j} \Big\} \  ,
\end{equation*}
where $ r_{j} = \emph{rank}( \mathbf{F}_{j}(\mathbf{a}) ) $ for $ j=1, \cdots, n $ as above. If now, we define the finite subset of $\mathbb{N}$
 \begin{equation*}
    \Omega( \mathbf{A} ) = \Big\{ i \in  \{1,2, \cdots ,p\}  \text{ } /   \text{ }  \xi ( i, \mathbf{A} ) =  \varnothing  \Big\} \ ,
\end{equation*}
the following theorem is valid.
 \begin{theo5.6} \label{theo5.6:box}
 With these definitions, let $ \emph{card}\big( \Omega( \mathbf{A} ) \big) $ be the number of elements of $  \Omega( \mathbf{A} ) $. If $ \Omega( \mathbf{A} ) $ is not empty then the matrices $ \mathbf{M}(\mathbf{a}) $, $ \mathbf{L}(\mathbf{a}) $ and $ \mathit{J}( \mathbf{r}(\mathbf{a}) ) $ have $ k.\emph{card}\big( \Omega( \mathbf{A} ) \big) $ columns equal to zero, moreover the indices of these columns in these matrices correspond.
\end{theo5.6}
\begin{proof}
We first consider the matrix $ \mathbf{M}(\mathbf{a}) $ and recall that this matrix has the following form
\begin{equation*}
\mathbf{M}(\mathbf{a}) = \mathbf{P}_{\mathbf{F}(\mathbf{a})}^{\bot} \emph{diag}(\emph{vec}(\sqrt{\mathbf{W}}))(\mathbf{\widehat{B}}^{T} \otimes \mathbf{I}_{p}) \mathbf{K}_{(k,p)} \  .
\end{equation*}
Then, we also recall that $ \mathbf{P}_{\mathbf{F}(\mathbf{a})}^{\bot} $ is a block-diagonal matrix of the form
\begin{equation*}
    \mathbf{P}_{\mathbf{F}(\mathbf{a})}^{\bot} = 
    \left\lbrack
\begin{array}{ccccc}
\mathbf{P}_{\mathbf{F}_1(\mathbf{a})}^{\bot} & 0          & \ldots             & 0          & 0          \\
0                      & \ddots &  0                    & \ldots   & 0          \\
\vdots              & 0         & \mathbf{P}_{\mathbf{F}_j(\mathbf{a})}^{\bot} & 0           & \vdots \\
0                      & \ldots  & 0                     & \ddots  & 0          \\
0                      & 0         & \ldots              & 0           & \mathbf{P}_{\mathbf{F}_n(\mathbf{a})}^{\bot}
\end{array} \right\rbrack = \bigoplus _{j=1}^n \mathbf{P}_{\mathbf{F}_j(\mathbf{a})}^{\bot} \ ,
\end{equation*}
with
\begin{equation*}
 \mathbf{P}_{\mathbf{F}_j(\mathbf{a})}^{\bot} =  \mathbf{I}_p -  \mathbf{F}_j(\mathbf{a}) \mathbf{F}_j(\mathbf{a})^{+} = \mathbf{I}_p - \big( \emph{diag}(\emph{vec}(\sqrt{\mathbf{W}}_{.j})) \mathbf{A} \big) \big( \emph{diag}(\emph{vec}(\sqrt{\mathbf{W}}_{.j}) ) \mathbf{A} \big)^{+} \ ,
\end{equation*}
for $ j \in  \{1,2, \cdots ,n\}$. Hence,
\begin{equation*}
    \mathbf{P}_{\mathbf{F}(\mathbf{a})}^{\bot} \emph{diag}(\emph{vec}(\sqrt{\mathbf{W}})) = \bigoplus _{j=1}^n \left( \mathbf{P}_{\mathbf{F}_j(\mathbf{a})}^{\bot} \emph{diag}(\emph{vec}( \sqrt{\mathbf{W}}_{.j}) ) \right) \ .
\end{equation*}
Consider now the following uniform blocking of the $ np \times pk $ matrix $ (\mathbf{\widehat{B}}^{T} \otimes \mathbf{I}_{p}) \mathbf{K}_{(k,p)} $ into $n.p$ submatrices $ \mathbf{O}_{ji} \in  \mathbb{R}^{p \times k} $, for $ j \in  \{1,2, \cdots ,n\} $ and $ i \in  \{1,2, \cdots ,p\} $,
\begin{equation*}
    (\mathbf{\widehat{B}}^{T} \otimes \mathbf{I}_{p}) \mathbf{K}_{(k,p)} = 
    \left\lbrack
\begin{array}{cccc}
\mathbf{O}_{11}  & \mathbf{O}_{12} & \ldots  & \mathbf{O}_{1p}   \\
\mathbf{O}_{21}  & \mathbf{O}_{22} & \ldots  & \mathbf{O}_{2p}   \\
\vdots                  & \vdots                 & \ddots & \vdots                   \\
\mathbf{O}_{n1}  & \mathbf{O}_{n2} & \ldots  & \mathbf{O}_{np}
\end{array} \right\rbrack , \mathbf{O}_{ji} \in \mathbb{R}^{ p \times k } \ .
\end{equation*}
Since $ (\mathbf{\widehat{B}}^{T} \otimes \mathbf{I}_{p}) \mathbf{K}_{(k,p)} = \mathbf{K}_{(n,p)} ( \mathbf{I}_{p} \otimes \mathbf{\widehat{B}}^{T}  ) $, $ (\mathbf{\widehat{B}}^{T} \otimes \mathbf{I}_{p}) \mathbf{K}_{(k,p)} $ is the matrix having as rows, every $n^{th}$ row of the matrix $ \mathbf{I}_{p} \otimes \mathbf{\widehat{B}}^{T} $ of order $ n.p \times p.k $, starting with the first, then every $n^{th}$ row starting with the second, and so on. Thus, keeping in mind that $ \mathbf{I}_{p} \otimes \mathbf{\widehat{B}}^{T} $ is a block-diagonal matrix, it is not hard to see that the $\mathbf{O}_{ji}$ submatrices are very sparse with the following simple structure
\begin{equation*}
    \mathbf{O}_{ji} = 
        \left\lbrack
        \begin{array}{c}
            \mathbf{0}^{(i-1) \times k}                             \\
            \lbrack \mathbf{\widehat{B}}^{T} \rbrack_{j.} \\
              \mathbf{0}^{(p-i) \times k} 
        \end{array}
        \right\rbrack
        \begin{array}{lc}
            \rbrace (i-1) \text{ rows } \\
            \rbrace 1     \text{ row }   &  , \text{ for } j=1,\ldots,n \text{ and } i=1,\ldots,p\\
            \rbrace (p-i) \text{ rows }  
        \end{array} \ .
\end{equation*}
Now, let $ i \in \Omega( \mathbf{A} ) $ and consider the submatrix $\mathbf{M}_i$ defined by the columns $(i-1).k + 1$ to $i.k$ of $ \mathbf{M}(\mathbf{a}) $. $\mathbf{M}_i$ is equal to
\begin{equation*}
    \mathbf{M}_i = \bigoplus_{j=1}^n \left( \mathbf{P}_{\mathbf{F}_j(\mathbf{a})}^{\bot}   \emph{diag}( \sqrt{\mathbf{W}}_{.j} ) \right)
        \left\lbrack
        \begin{array}{c}
            \mathbf{O}_{1i}  \\
            \vdots                 \\
            \mathbf{O}_{ni} 
        \end{array}
        \right\rbrack = 
        \left\lbrack
        \begin{array}{c}
            \mathbf{P}_{\mathbf{F}_1(\mathbf{a})}^{\bot}    \emph{diag}( \sqrt{\mathbf{W}}_{.1} ) \mathbf{O}_{1i}  \\
            \vdots                 \\
            \mathbf{P}_{\mathbf{F}_n(\mathbf{a})}^{\bot}  \emph{diag}( \sqrt{\mathbf{W}}_{.n} ) \mathbf{O}_{ni} 
        \end{array}
        \right\rbrack \ .
\end{equation*}
Now, for $ j=1, \ldots, n $, we have $ \mathbf{W}_{ij} = 0 $ or $ \mathbf{W}_{ij} \ne 0 $:

- If $ \mathbf{W}_{ij} = 0 $, then $ \emph{diag}( \sqrt{\mathbf{W}}_{.j} ) \mathbf{O}_{ji} = 0^{ p \times k} $ and the rows $(j-1).p + 1$ to $j.p$ of $\mathbf{M}_i$ are equal to zero.

- If $ \mathbf{W}_{ij} \ne 0 $, we have $ \boldsymbol{\delta}_{ij} = 1 $ and $  \displaystyle{ \sum_{l=1}^p  \boldsymbol{\delta}_{lj} } =  r_{j}$ since $ i \in \Omega( \mathbf{A} ) $, and it follows that 
\begin{equation*}
    \mathbf{P}_{\mathbf{F}_j(\mathbf{a})}^{\bot} = \mathbf{I}_p  - \emph{diag}( \boldsymbol{\delta}_{.j} )
\end{equation*}
and
\begin{align*}
    \mathbf{P}_{\mathbf{F}_j(\mathbf{a})}^{\bot}  \emph{diag}( \sqrt{\mathbf{W}}_{.j} ) & =  \big( \mathbf{I}_p  - \emph{diag}( \boldsymbol{\delta}_{.j} ) \big) \emph{diag}( \sqrt{\mathbf{W}}_{.j} )  \\
                                                                                                                                     & =   \emph{diag}( \sqrt{\mathbf{W}}_{.j} ) - \emph{diag}( \sqrt{\mathbf{W}}_{.j} )     \\
                                                                                                                                     & =  \mathbf{0}^{p},
\end{align*}
and the rows $ (j-1).p + 1 $ to $ j.p $ of $\mathbf{M}_i$ are also equal to zero if $ \mathbf{W}_{ij} \ne 0 $.

Hence, we finally obtain $ \mathbf{M}_i = 0^{ p.n \times k} $ if $ i \in \Omega( \mathbf{A} ) $ and we conclude that $ \mathbf{M}(\mathbf{a}) $ has $ k.\emph{card}\big( \Omega( \mathbf{A} ) \big) $ columns equal to zero, as claimed in the theorem.

Turning now our attention to $ \mathbf{L}(\mathbf{a}) $, which is equal to
\begin{equation*}
    \mathbf{L}(\mathbf{a}) = ( \mathbf{F}(\mathbf{a})^{+} )^{T} \big( ( \mathbf{W} \odot P_{\Omega}(\mathbf{X} -\mathbf{A}\mathbf{\widehat{B}}) )^{T} \otimes \mathbf{I}_{k} \big) \ ,
\end{equation*}
we observe that it is sufficient to show that the $i^{th}$ row of $ \mathbf{W} \odot P_{\Omega}(\mathbf{X} -\mathbf{A}\mathbf{\widehat{B}}) $ is equal to zero to establish that the columns $ (i-1).k + 1 $ to $ i.k $ of $ \mathbf{L}(\mathbf{a}) $ are equal to zero if $ i \in \Omega( \mathbf{A} ) $.

Now, for $ j=1, \ldots, n $, we have $ \mathbf{W}_{ij} = 0 $ or $ \mathbf{W}_{ij} \ne 0 $:

- If $ \mathbf{W}_{ij} = 0 $, then $ \left\lbrack \mathbf{W} \odot P_{\Omega}(\mathbf{X} -\mathbf{A}\mathbf{\widehat{B}}) \right\rbrack_{ij} = 0 $ for any $ p \times n $ matrix $ P_{\Omega}(\mathbf{X} -\mathbf{A}\mathbf{\widehat{B}}) $.

- If $ \mathbf{W}_{ij} \ne 0 $, as above, we have $ \boldsymbol{\delta}_{ij} = 1 $ and $  \displaystyle{ \sum_{l=1}^p  \boldsymbol{\delta}_{lj} } =  r_{j} $ since $ i \in \Omega( \mathbf{A} ) $, implying that
\begin{equation*}
    \mathbf{P}_{\mathbf{F}_j(\mathbf{a})}^{\bot}  = \mathbf{I}_p  - \emph{diag}( \boldsymbol{\delta}_{.j} )
\end{equation*}
and it follows that
\begin{align*}
    \left\lbrack \mathbf{W} \odot P_{\Omega}(\mathbf{X} -\mathbf{A}\mathbf{\widehat{B}}) \right\rbrack_{.j} & =   \sqrt{\mathbf{W}}_{.j}  \odot \left( \sqrt{\mathbf{W}}_{.j} \odot P_{\Omega}(\mathbf{X} -\mathbf{A}\mathbf{\widehat{B}})_{.j}  \right) \\
                                                            & =  \sqrt{\mathbf{W}}_{.j}  \odot \mathbf{P}_{\mathbf{F}_j(\mathbf{a})}^{\bot}  ( \sqrt{\mathbf{W}}_{.j} \odot \mathbf{X}_{.j} ) \\
                                                            & =  \sqrt{\mathbf{W}}_{.j}  \odot \big( \mathbf{I}_p  - \emph{diag}( \boldsymbol{\delta}_{.j} ) \big) ( \sqrt{\mathbf{W}}_{.j} \odot \mathbf{X}_{.j} ) \\
                                                            & =  \sqrt{\mathbf{W}}_{.j}  \odot \big( \sqrt{\mathbf{W}}_{.j} \odot \mathbf{X}_{.j} - \sqrt{\mathbf{W}}_{.j} \odot \mathbf{X}_{.j} \big) \\
                                                            & =  \mathbf{0}^{p} \ ,
\end{align*}
and $ \left\lbrack \mathbf{W} \odot P_{\Omega}(\mathbf{X} -\mathbf{A}\mathbf{\widehat{B}}) \right\rbrack_{ij} = 0 $ also in the case $ \mathbf{W}_{ij} \ne 0 $.

Hence, we conclude that the columns $ (i-1).k + 1 $ to $ i.k $ of $ \mathbf{L}(\mathbf{a}) $ are equal to zero if $ i \in \Omega( \mathbf{A} ) $ and $ \mathbf{L}(\mathbf{a}) $  has also $ k.\emph{card}\big( \Omega( \mathbf{A} ) \big) $ columns equal to zero, as claimed in the theorem. Finally, the same result holds for $ \mathit{J}( \mathbf{r}(\mathbf{a}) ) $ since $ \mathit{J}( \mathbf{r}(\mathbf{a}) ) = - ( \mathbf{M}(\mathbf{a}) + \mathbf{L}(\mathbf{a}) ) $.

\end{proof}
Clearly, any variable projection Gauss-Newton or Levenberg-Marquardt algorithm using the Jacobian matrix $ \mathit{J}( \mathbf{r}(\mathbf{a}) ) $, or approximating this Jacobian matrix with $ -\mathbf{M}(\mathbf{a}) $, will be incorrect if $ \Omega( \mathbf{A} ) $ is not empty as demonstrated in Theorem~\ref{theo5.6:box}. The fact that $ \Omega( \mathbf{A} ) $ is not empty, is symptomatic of the situation where we try to fit the matrix $ \mathbf{X} $ by a model with too many components with respect to the number of missing elements in this matrix. However, if we restrict the set of WLRA problems by imposing the condition $\displaystyle{ \sum_{l=1}^p  \boldsymbol{\delta}_{lj} } >  k $ for all $ j=1, \cdots, n$, we are sure that $ \Omega( \mathbf{A} ) $ will be empty and the events described in Theorem~\ref{theo5.6:box} will not occurred.

Finally, if $ \emph{card}\big( \Omega( \mathbf{A} ) \big) $ is not equal to zero, a much better alternative to find a solution of this particular WLRA problem, is to minimize the regularized cost function $g_\lambda(.)$ defined in equation~\eqref{eq:g_func} of Subsection~\ref{approx_wlra:box} with a variable projection algorithm, as this regularized minimization problem is always well-posed (see Subsection~\ref{approx_wlra:box} for details).

\subsection{Computations and properties of the gradient vector and Hessian matrix} \label{hess:box}

Using the properties of the Jacobian matrix $\mathit{J}( \mathbf{r}(\mathbf{a}) )$ demonstrated in Subsection~\ref{jacob:box}, we now derive a simple expression for the gradient of $\psi(.)$ and we give a detailed study of the Hessian matrix $\nabla^2 \psi ( \mathbf{a} ) $. These results may be used to formulate steepest descent or Newton algorithms for minimizing the variable projection functional $\psi(.)$~\cite{SJ2004}\cite{C2008b}\cite{SE2010}\cite{DMK2011}\cite{BA2011}\cite{BA2015}\cite{BL2020}. We also again illustrate the tight relationships between these Euclidean gradient and Hessian operators and the corresponding Riemannian gradient and Hessian operators when we consider the cost function $\psi(.)$  as operating on the Grassmann manifold $\text{Gr}(p,k)$~\cite{AMS2008}\cite{BA2015}\cite{B2023} extending the investigations of~\cite{HF2015b} on this topic.

Since $\mathit{J}( \mathbf{r}(\mathbf{a}) )$ is rank-deficient everywhere according to Theorem~\ref{theo5.2:box}, the smallest eigenvalue of $\mathit{J}( \mathbf{r}(\mathbf{a}) )^{T}\mathit{J}( \mathbf{r}(\mathbf{a}) )$ is always equal to zero and, in these conditions, we cannot expect that the Gauss-Newton term $\mathit{J}( \mathbf{r}(\mathbf{a}) )^{T}\mathit{J}( \mathbf{r}(\mathbf{a}) )$ will dominate the second term in the expression of $\nabla^2 \psi ( \mathbf{a} ) $ (see Subsection~\ref{opt:box}). This suggests that a full-Newton approach can perform much better than the Gauss-Newton or Levenberg-Marquardt methods for solving the~\eqref{eq:VP1} problem. However, the Gauss-Newton or Levenberg-Marquardt algorithms perform surprisingly better than different versions of the full Newton algorithm for solving the~\eqref{eq:VP1} problem in the comparative studies of Okatani et al.~\cite{OYD2011} or Hong et al.~\cite{HF2015}.These conclusions are thus very counterintuitive taking into account the excellent properties (e.g., fast convergence for NLLS problems with large residuals and quadratic local convergence in a neighborhood of a stationary point) of the full Newton approach compared to the Gauss-Newton or Levenberg-Marquardt methods for general or variable projection NLLS problems~\cite{DS1983}\cite{MN2010}\cite{B2009}\cite{HPS2012}.

The results of this Subsection will try to elucidate these contradictions and also will reveal the power of the variable projection framework in understanding the intrinsic difficulties associated with the WLRA problem. As an illustration, we will demonstrate below, that the Hessian matrix $\nabla^2 \psi (\mathbf{\widehat{a} })$ is  deficient at all first-order stationary points $\mathbf{\widehat{a} }$ of $\psi(.)$ (see Theorem~\ref{theo5.9:box} below) and, consequently, this Hessian matrix is expected to be nearly singular and ill-conditioned in a "small" neighborhood of a first-order stationary point $\mathbf{\widehat{a} }$ of  $\psi(.)$. This result is consistent with the fact that (local) minimizers of $\varphi^{*}(.)$ and $\psi(.)$ are never isolated, as already discussed in Subsection~\ref{noconv_wlra:box}, and that any neighborhood of a (local) minimizer of these cost functions contains also an infinite number of other minimizers, which attain the same minimum of $\psi(.)$, implying that the Hessian matrix $\nabla^2 \psi (\mathbf{\widehat{a} })$ at a local minimizer $\mathbf{\widehat{a}}$ is at best positive semi-definite, but never positive definite. When the weight matrix $\mathbf{W} \in \mathbb{R}^{p \times n}_{*+}$ and we consider the cost function $\psi(.)$ as defined on the (quotient) Grassmann manifold $\text{Gr}(n,k)$, because $\psi(.)$ is invariant on the equivalence classes of this quotient, the fact that the Hessian cannot  possibly be positive definite at first-order stationary points is already a known result, see Chapter 9 and Lemma 9.41 in~\cite{B2023}. Our Theorem~\ref{theo5.9:box} thus provides an extension of this result as it also applies to the case where  $\mathbf{W} \in \mathbb{R}^{p \times n}_{+}$. These results also imply that the Hessian matrices at points arbitrarily closed to minima have vanishingly small, possibly negative eigenvalues, leading to ill-conditioned and indefinite linear systems with a severe loss of accuracy in Newton and trust-regions methods when the iterates reach a neighborhood of a minimum~\cite{RB2024}.

In these conditions, many of the excellent properties of full-Newton and trust-regions approaches are lost, which may explain the degraded performance of these methods for minimizing $\psi(.)$ in the comparative studies of Okatani et al.~\cite{OYD2011} and Hong et al.~\cite{HF2015}. This poor performance concerns especially the second-order Riemannian trust-region method (RTRMC2)  operating on the Grassmann manifold developed initially in Boumal and Absil~\cite{BA2011} and already discussed in Subsection~\ref{approx_wlra:box}. In their Riemannian Newton approach, Boumal and Absil~\cite{BA2011}\cite{BA2015} have derived a compact directional derivative formulae for the Riemannian Hessian of the cost function of a regularized form of the WLRA problem (see the  regularized cost function $g_\lambda(.)$ defined in equation~\eqref{eq:g_func} of Subsection~\ref{approx_wlra:box}), which is then used in an inexact subproblem solver based on a truncated conjugate gradient method to compute an approximate (Riemannian) Newton step at each iteration of their RTRMC2 method. However, in order to ensure convergence to (local) minima, the truncated conjugate gradient theory assumes a positive definite system at these minima of the cost function, an hypothesis which is not verified here near first-order critical points, where the RTRMC2 method may be applied to ill-conditioned, indefinite systems, leading to an erratic behaviour for some WLRA problems as the truncated conjugate gradient subproblem solver may be highly sensitive to negative eigenvalues, even of small magnitude. This challenging question of the convergence of trust-regions methods for non-isolated minima of smooth functions defined on a (arbitrary) manifold has been revisited recently in~\cite{RB2024} in which the authors were still able to derive convergence results for  an inexact solver based on a truncated conjugate gradient method under some additional hypotheses and with a carefully designed inexact subproblem solver with had hoc stopping criteria for the truncated conjugate gradient iterations so that the subproblem solver is not affected by the small negative eigenvalues of the Hessian matrix.

These convergence problems of trust-region methods near (local) minima may explain why Boumal and Absil~\cite{BA2015} have subsequently introduced a pre-conditioner for the Hessian matrix in their RTRMC2 method originally proposed in~\cite{BA2011}. However, we are not aware of any new comparison studies, which evaluate the performance of the updated and preconditioned RTRMC2 algorithm proposed in~\cite{BA2015} with the Gauss-Newton or Levenberg-Marquardt approaches described in~\cite{OD2007}\cite{OYD2011}\cite{HF2015}\cite{HZF2017} to verify if this preconditioned version of  RTRMC2 performs now better than the variable projection Gauss-Newton or Levenberg-Marquardt methods for solving  WLRA problems. Furthermore, almost all these previous studies deal only of WLRA problems with binary weights (e.g., the missing value problem) excepted of the work of Boumal and Absil~\cite{BA2015}. This clearly shows the need of new extensive  comparison studies to clarify the respective performance of the most recent various first- or second-order methods proposed in the literature for solving the general WLRA problem, but such ambitious task is outside of the scope of this paper, which is mainly devoted to a better understanding of the theoretical properties of variable projections methods for solving the WLRA problem.

To start with, we give convenient different expressions for the Euclidean gradient of the variable projection functional $\psi(.)$. This proposition is mainly a reformulation of Theorem 2.1 in Golub and Pereyra~\cite{GP1973}; however, we give a direct proof using only linear algebra in order to be self-contained and because this formula is not well-known in the literature related to the WLRA problem.
\\
 \begin{theo5.7} \label{theo5.7:box}
 Let $\varphi^{*}( \mathbf{A},\mathbf{\widehat{B}} )$  and $\psi( \mathbf{a} )$ be defined, respectively, as in equations~\eqref{eq:P1} and~\eqref{eq:psi_func}  of Section~\ref{seppb:box} with $\mathbf{A}  \in \mathbb{R}^{p \times k}_{k}$, $\mathbf{a} = \emph{vec}( \mathbf{A}^{T} ) \in \mathbb{R}^{p.k}$, $\mathbf{\widehat{B}}  \in \mathbb{R}^{k \times n}$ is such that $\mathbf{\widehat{b}} = \emph{vec}( \mathbf{\widehat{B}} ) = \mathbf{F}(\mathbf{a})^{+} \mathbf{x}  \in \mathbb{R}^{k.n}$ and  $\mathbf{x} =  \emph{vec}(  \sqrt{ \mathbf{W} } \odot \mathbf{X} )  \in \mathbb{R}^{p.n}$, where $\mathbf{X}  \in \mathbb{R}^{ p \times n }$  and $\mathbf{W}  \in \mathbb{R}^{ p \times n }_{+}$ are, respectively, the data and weight matrices of the WLRA problems~\eqref{eq:P0} or~\eqref{eq:P1}.
 
 We further assume that $\mathbf{a}$ belongs to an open set $ \varOmega \subset  \mathbb{R}^{p.k} $  in which $ \mathbf{F}(.) $ has a constant rank, so that the Jacobian matrix $ \mathit{J}( \mathbf{r}(\mathbf{a}) ) $ derived in Subsection~\ref{jacob:box} and the Euclidean gradient $\nabla \psi( \mathbf{a} )$ are both well-defined. Then,
\begin{equation*}
    \nabla \psi( \mathbf{a} ) = - \mathbf{M}(\mathbf{a})^{T} \mathbf{r}(\mathbf{a}) = - \mathbf{M}(\mathbf{a})^{T} \big( \mathbf{x} -  \mathbf{F}( \mathbf{a} ) \mathbf{\widehat{b}}  \big) = \frac{ \partial \varphi^{*}( \mathbf{A},\mathbf{\widehat{B}} ) }{ \partial \mathbf{a} }
\end{equation*}
and
\begin{equation*}
\Vert \nabla \psi( \mathbf{a} ) \Vert_{2} = \big\Vert  \nabla  \varphi^{*}_{\mathbf{a}}( \mathbf{A},\mathbf{\widehat{B}} )  \big\Vert_{2} = \big\Vert \nabla  \varphi^{*}_{\mathbf{A}}( \mathbf{A},\mathbf{\widehat{B}} )  \big\Vert_{F} \ ,
\end{equation*}
where $\frac{ \partial \varphi^{*}( \mathbf{A},\mathbf{\widehat{B}} ) }{ \partial \mathbf{a} }$, $\nabla  \varphi^{*}_{\mathbf{a}}( \mathbf{A},\mathbf{\widehat{B}} )$ and $\nabla  \varphi^{*}_{\mathbf{A}}( \mathbf{A},\mathbf{\widehat{B}} )$ are defined, respectively, by equations~\eqref{eq:Da_varphi*} and~\eqref{eq:DAB_varphi*}  in Theorem~\ref{theo4.3:box} and equation~\eqref{eq:D_partial_grad_varphi*} of Subsection~\ref{landscape_wlra:box}.

In addition, the theorem remains valid if $\mathbf{\widehat{b}} = \mathbf{F}(\mathbf{a})^- \mathbf{x}$ where $\mathbf{F}(\mathbf{a})^-$ is a symmetric generalized inverse of $\mathbf{F}(\mathbf{a})$ as defined in equations~\eqref{eq:sginv} or~\eqref{eq:sginv_qrcp} of Subsection~\ref{lin_alg:box}.
\end{theo5.7}
\begin{proof}
Using the notations and results in Subsection~\ref{jacob:box}, we have
\begin{align*}
    \nabla \psi( \mathbf{a} )  & =  \mathit{J}\big( \mathbf{r}(\mathbf{a}) \big)^{T} \mathbf{r}(\mathbf{a})  \\
                                           & =  - \big( \mathbf{M}(\mathbf{a}) + \mathbf{L}(\mathbf{a}) \big)^{T} \mathbf{r}(\mathbf{a})             \\
                                           & =  - \mathbf{M}(\mathbf{a})^{T} \mathbf{r}(\mathbf{a}) - \mathbf{L}(\mathbf{a})^{T} \mathbf{r}(\mathbf{a})  \\
                                           & =  - \mathbf{M}(\mathbf{a})^{T} \mathbf{r}(\mathbf{a}) \ ,
\end{align*}
since  $\mathbf{r}(\mathbf{a}) \in \emph{ran}( \mathbf{F}(\mathbf{a}) )^{\bot}$ (see equation~\eqref{eq:r_vec}), $\emph{ran}( \mathbf{F}(\mathbf{a}) )^{\bot}  \subset \emph{ran}( \mathbf{L}(\mathbf{a}) )^{\bot}$ (see equation~\eqref{eq:L_mat2}) and $\emph{ran}( \mathbf{L}(\mathbf{a}) )^{\bot} = \emph{null}( \mathbf{L}(\mathbf{a})^{T} )$. The last equality resulting from equation~\eqref{eq:null_ran} in Subsection~\ref{lin_alg:box}.
 
 Hence, $\mathbf{L}(\mathbf{a})$, the second term of $\mathit{J}( \mathbf{r}(\mathbf{a}) )$, does not contribute to the gradient $\nabla \psi( \mathbf{a} )$. Now, using the second formulation of $\mathbf{M}(\mathbf{a})$  given in equation~\eqref{eq:M_mat2} of Subsection~\ref{jacob:box}, we deduce
\begin{align*}
    \nabla \psi( \mathbf{a} )
    & =  - \big( \mathbf{P}_{\mathbf{F}(\mathbf{a})}^{\bot} \mathbf{K}_{(n,p)} \mathbf{G}(\mathbf{\widehat{b}}) \big)^{T} \mathbf{r}(\mathbf{a})  \\
    & =  - \mathbf{G}(\mathbf{\widehat{b}})^{T} \mathbf{K}_{(p,n)} \mathbf{P}_{\mathbf{F}(\mathbf{a})}^{\bot} \mathbf{r}(\mathbf{a})  \\
    & =  - \mathbf{G}(\mathbf{\widehat{b}})^{T} \mathbf{K}_{(p,n)} \mathbf{r}(\mathbf{a}) \ ,
\end{align*}
since $ \mathbf{r}(\mathbf{a}) \in  \emph{ran}( \mathbf{F}(\mathbf{a}) )^{\bot} = \emph{ran}( \mathbf{P}_{\mathbf{F}(\mathbf{a})}^{\bot} ) $. Noting that
\begin{equation*}
    \mathbf{r}(\mathbf{a}) = \mathbf{x} - \mathbf{F}(\mathbf{a}) \mathbf{\widehat{b}} = \mathbf{K}_{(n,p)} \big( \mathbf{z} - \mathbf{G}(\mathbf{\widehat{b}}) \mathbf{a} \big) \ ,
\end{equation*}
where $ \mathbf{z} = \emph{vec}\big( ( \sqrt{ \mathbf{W} }  \odot \mathbf{X} )^{T} \big)$ (see equation~\eqref{eq:r_vec} in Section~\ref{seppb:box} for details), we finally obtain
\begin{align*}
        \nabla \psi( \mathbf{a} )
    & =  - \mathbf{G}(\mathbf{\widehat{b}})^{T} \mathbf{K}_{(p,n)} \mathbf{K}_{(n,p)} \big( \mathbf{z} - \mathbf{G}(\mathbf{\widehat{b}}) \mathbf{a} \big)  \\
    & =  - \mathbf{G}(\mathbf{\widehat{b}})^{T}  \big( \mathbf{z} - \mathbf{G}(\mathbf{\widehat{b}}) \mathbf{a} \big)  \\
    & =    \mathbf{G}(\mathbf{\widehat{b}})^{T} \mathbf{G}(\mathbf{\widehat{b}}) \mathbf{a} - \mathbf{G}(\mathbf{\widehat{b}})^{T} \mathbf{z} \\
    & =   \frac{ \partial \varphi^{*}( \mathbf{A},\mathbf{\widehat{B}} ) }{ \partial \mathbf{a} } \ ,
\end{align*}
where the last equality results from equation~\eqref{eq:Da_varphi*} in Theorem~\ref{theo4.3:box}.
Finally, using the fact that
\begin{equation*}
\mathbf{\widehat{B}} = \text{Arg} \min_{\mathbf{B}\in\mathbb{R}^{k \times n}} \,  \varphi^{*}(\mathbf{A},\mathbf{B}  ) = \frac{1}{2}   \Vert \mathbf{x} - \mathbf{F}(\mathbf{a})\mathbf{b} \Vert^{2}_{2} \ ,
\end{equation*}
we have
 \begin{equation*}
\frac{ \partial  \varphi^{*}( \mathbf{A},\mathbf{\widehat{B}} ) }{ \partial \mathbf{b} } = \mathbf{0}^{k.n} \ ,
\end{equation*}
and so the vectorized form of the Euclidean gradient of  $\varphi^{*}(.)$ at $(\mathbf{a},\mathbf{\widehat{b}} )$ is given by
\begin{equation*}
\nabla \varphi^{*}(\mathbf{A}, \mathbf{\widehat{B}}  )  = \begin{bmatrix}   \frac{ \partial \varphi^{*}(\mathbf{A},\mathbf{\widehat{B}} )   }{ \partial \mathbf{a} }  &  \mathbf{0}^{k.n} \end{bmatrix} \ ,
\end{equation*}
which implies that
\begin{equation*}
\Vert \nabla \psi( \mathbf{a} ) \Vert_{2} = \big\Vert \nabla  \varphi^{*}_{\mathbf{a}}( \mathbf{A},\mathbf{\widehat{B}} )  \big\Vert_{2} \ .
\end{equation*}
Next, the equality
\begin{equation*}
\Vert \nabla \psi( \mathbf{a} ) \Vert_{2} = \big\Vert \nabla  \varphi^{*}_{\mathbf{A}}( \mathbf{A},\mathbf{\widehat{B}} )  \big\Vert_{F} 
\end{equation*}
is a direct consequence of equation~\eqref{eq:DAB_varphi*}  in Theorem~\ref{theo4.3:box}.

Finally, the above demonstration remains valid if $\mathbf{\widehat{b}} = \mathbf{F}(\mathbf{a})^{-} \mathbf{x}$ because the differential formula~\eqref{eq:d_projmat2} for an orthogonal projector is unchanged if a symmetric generalized inverse is used in place of the pseudo-inverse and the residual vector $\mathbf{r}(\mathbf{a})$ is also identical since $\mathbf{r}(\mathbf{a}) = \mathbf{P}_{\mathbf{F}(\mathbf{a})}^{\bot} \mathbf{x}$. This concludes the proof of Theorem~\ref{theo5.7:box}.
\\
\end{proof}
Remembering that $\mathbf{G}(\mathbf{\widehat{b}})$ is a block-diagonal matrix, we see that the computation of $\nabla \psi (\mathbf{a})$ is easy, fast and may be efficiently parallelized using Theorem~\ref{theo5.7:box}, since $ \mathbf{G}(\mathbf{\widehat{b}})^{T} \mathbf{G}(\mathbf{\widehat{b}}) $ is also a block-diagonal matrix. We further highlight that this formulation of $\nabla \psi (\mathbf{a})$ is much more efficient than the formulae given in Chen~\cite{C2008b} (see its equations 24 and 27), which does not exploit the fact that the residual vector $\mathbf{r}(\mathbf{a}) $ is linear in both $\mathbf{a}$ and $\mathbf{b}$, and involved multiplication of matrices with many zeros.

With the help of Theorem~\ref{theo5.7:box}, it is also easy to demonstrate that the Gauss-Newton directions defined in Subsection~\ref{jacob:box} are in a descent direction for $\psi(.)$ despite the systematic rank deficiency of the Jacobian matrix $\mathit{J}( \mathbf{r}(\mathbf{a}) )$ or of its approximation $-\mathbf{M}(\mathbf{a})$ proved in Theorem~\ref{theo5.2:box}.
\begin{corol5.7} \label{corol5.7:box}
Let $\mathbf{A}  \in \mathbb{R}^{p \times k}_{k}$ and $\mathbf{a} = \emph{vec}( \mathbf{A}^{T} ) \in \mathbb{R}^{k.p}$.

If $\nabla \psi( \mathbf{a} ) \ne  \mathbf{0}^{p.k}$, e.g., if $\mathbf{a}$ is not a first-order stationary point of $\psi(.)$, the Gauss-Newton directions defined by
\begin{align*}
d\mathbf{a}_{gn} & = - \mathit{J} \big( \mathbf{r}(\mathbf{a}) \big)^{+} \mathbf{r}(  \mathbf{a} ) = \big( \mathbf{M}(\mathbf{a}) + \mathbf{L}(\mathbf{a}) \big)^{+} \mathbf{r}(  \mathbf{a} ) \ , \\
 d\mathbf{a}_{gn} & =  \mathbf{M}(\mathbf{a})^{+} \mathbf{r}(  \mathbf{a} ) 
\end{align*}
are in a descent direction for $\psi(.)$.
\end{corol5.7}
\begin{proof}
To prove the assertion if the Gauss-Newton direction $d\mathbf{a}_{gn}$ is defined from the pseudo-inverse of the full Jacobian matrix $\mathit{J}( \mathbf{r}(\mathbf{a}) )$, note that
\begin{align*}
d\mathbf{a}_{gn}^{T} \nabla \psi( \mathbf{a} ) & =  - \mathbf{r}(  \mathbf{a} )^{T} \mathit{J}\big( \mathbf{r}(\mathbf{a}) \big)^{+T} \mathit{J}\big( \mathbf{r}(\mathbf{a}) \big)^{T} \mathbf{r}(  \mathbf{a} ) \\
                                                                      & =  - \mathbf{r}(  \mathbf{a} )^{T} \mathit{J}\big( \mathbf{r}(\mathbf{a}) \big) \mathit{J}\big( \mathbf{r}(\mathbf{a}) \big)^{+} \mathbf{r}(  \mathbf{a} ) \\
                                                                      & =  - \mathbf{r}(  \mathbf{a} )^{T} \mathbf{P}_{\mathit{J}( \mathbf{r}(\mathbf{a}) )} \mathbf{r}(  \mathbf{a} ) \\
                                                                      & =  - \mathbf{r}(  \mathbf{a} )^{T} \mathbf{P}_{\mathit{J}( \mathbf{r}(\mathbf{a}) )}^{T} \mathbf{P}_{\mathit{J}( \mathbf{r}(\mathbf{a}) )} \mathbf{r}(  \mathbf{a} ) \\
                                                                      & =  - \Vert  \mathbf{P}_{\mathit{J}( \mathbf{r}(\mathbf{a}) )} \mathbf{r}(  \mathbf{a} ) \Vert_{2}^{2} \ ,
\end{align*}
since the projector $\mathbf{P}_{\mathit{J}( \mathbf{r}(\mathbf{a}) )}$ is a symmetric and idempotent matrix. Further, $\nabla \psi( \mathbf{a} ) \ne \mathbf{0}^{p.k}$ by hypothesis, then $\mathit{J}( \mathbf{r}(\mathbf{a}) )^{T} \mathbf{r}(  \mathbf{a} ) \ne \mathbf{0}^{p.k}$ and $\mathbf{r}(  \mathbf{a} )$ is not in the null space of $\mathit{J}( \mathbf{r}(\mathbf{a}) )^{T}$ and is not the zero-vector. Since
\begin{equation*}
\emph{null}\Big( \mathit{J}\big( \mathbf{r}(\mathbf{a}) \big)^{T} \Big) = \emph{ran}\Big( \mathit{J}\big( \mathbf{r}(\mathbf{a}) \big) \Big)^{\bot} = \emph{null}( \mathbf{P}_{\mathit{J}( \mathbf{r}(\mathbf{a}) )} ) ,
\end{equation*}
it follows that $\mathbf{P}_{\mathit{J}( \mathbf{r}(\mathbf{a}) )} \mathbf{r}(\mathbf{a}) \ne \mathbf{0}^{p.n}$ and $d\mathbf{a}_{gn}^{T} \nabla \psi( \mathbf{a} ) = - \Vert  \mathbf{P}_{\mathit{J}( \mathbf{r}(\mathbf{a}) )} \mathbf{r}(  \mathbf{a} ) \Vert_{2}^{2} < 0$.
This proves the first assertion. The second assertion if the  Gauss-Newton direction $d\mathbf{a}_{gn}$ is defined from the pseudo-inverse of the approximate Jacobian matrix $- \mathbf{M}(\mathbf{a})$ is verified by the same way since
\begin{equation*}
\nabla \psi( \mathbf{a} ) = - \mathbf{M}(\mathbf{a})^{T} \mathbf{r}(\mathbf{a}),
\end{equation*}
as proved in Theorem~\ref{theo5.7:box}.
\\
\end{proof}
These results show that it is appropriate to use a line search algorithm in the Gauss-Newton methods described in Subsection~\ref{jacob:box} in order to obtain global convergence even though the Jacobian matrix or its approximation are always singular. Moreover, similar results hold for the Levenberg-Marquardt directions defined in Subsection~\ref{jacob:box} as demonstrated in the following corollary.
\begin{corol5.8} \label{corol5.8:box}
Let $\mathbf{A}  \in \mathbb{R}^{p \times k}_{k}$, $\mathbf{a} = \emph{vec}( \mathbf{A}^{T} ) \in \mathbb{R}^{p.k}$, $\lambda \in \mathbb{R}_{+*}$ be the Marquardt damping parameter and $\mathbf{D}$ be a diagonal scaling matrix of order $p.k$ with diagonal elements $\mathbf{D}_{ii}>0$.

If $\nabla \psi( \mathbf{a} ) \ne \mathbf{0}^{k.p}$, e.g., if $\mathbf{a}$ is not a first-order stationary point of $\psi(.)$, the Levenberg-Marquardt directions
\begin{align*}
d\mathbf{a}_{lm} & = -  \left(  \mathit{J}\big( \mathbf{r}(\mathbf{a}) \big)^{T} \mathit{J}\big( \mathbf{r}(\mathbf{a}) \big) +  \lambda \mathbf{D}^{T} \mathbf{D}   \right)^{-1} \mathit{J}\big( \mathbf{r}(\mathbf{a}) \big)^{T} \mathbf{r}(  \mathbf{a} ) \ ,  \\
d\mathbf{a}_{lm} & = \left(  \mathbf{M}(\mathbf{a})^{T}  \mathbf{M}(\mathbf{a})  +  \lambda \mathbf{D}^{T} \mathbf{D}   \right)^{-1}  \mathbf{M}(\mathbf{a})^{T} \mathbf{r}(  \mathbf{a} )  \ , \\
d\mathbf{a}_{lm} & = \left(  \mathbf{M}(\mathbf{a})^{T}  \mathbf{M}(\mathbf{a})   +  \mathbf{N} \mathbf{N}^{T} +  \lambda \mathbf{D}^{T} \mathbf{D}  \right)^{-1}  \mathbf{M}(\mathbf{a})^{T} \mathbf{r}(  \mathbf{a} )  \ , \\
d\mathbf{a}_{lm} & = -  \left(  \mathit{J}\big( \mathbf{r}(\mathbf{a}) \big)^{T} \mathit{J}\big( \mathbf{r}(\mathbf{a}) \big)  +  \mathbf{N} \mathbf{N}^{T} +  \lambda \mathbf{D}^{T} \mathbf{D}  \right)^{-1} \mathit{J}\big( \mathbf{r}(\mathbf{a}) \big)^{T} \mathbf{r}(  \mathbf{a} ) \ ,
\end{align*}
where $\mathbf{N} = \mathbf{K}_{(p,k)} ( \mathbf{I}_{k} \otimes \mathbf{A} )$, are well defined and are also in a descent direction for $\psi(.)$.

In addition, if $\lambda = 0$ and $\emph{rank}( \mathit{J}( \mathbf{r}(\mathbf{a}) ) ) = \emph{rank}( \mathbf{M}(\mathbf{a}) ) = (p - k).k$, these two last Levenberg-Marquardt directions are also well-defined, again in a descent direction for $\psi(.)$ and equal to the corresponding Gauss-Newton directions defined in Corollary~\ref{corol5.7:box}.
\end{corol5.8}
\begin{proof}
To prove the first part of the Corollary, we note that  $\lambda \ne 0$ and all the elements of the diagonal matrix $\mathbf{D}$ are strictly positive by hypothesis. In these conditions, the matrices
\begin{align*}
& \left(  \mathit{J}\big( \mathbf{r}(\mathbf{a}) \big)^{T} \mathit{J}\big( \mathbf{r}(\mathbf{a}) \big) +  \lambda \mathbf{D}^{T} \mathbf{D}   \right) &\text{ , } & \left(  \mathbf{M}(\mathbf{a})^{T}  \mathbf{M}(\mathbf{a})  +  \lambda \mathbf{D}^{T} \mathbf{D}  \right) , \\
& \left(  \mathit{J}\big( \mathbf{r}(\mathbf{a}) \big)^{T} \mathit{J}\big( \mathbf{r}(\mathbf{a}) \big) +  \lambda \mathbf{D}^{T} \mathbf{D}  +  \mathbf{N} \mathbf{N}^{T} \right) &\text{ , }& \left(  \mathbf{M}(\mathbf{a})^{T}  \mathbf{M}(\mathbf{a})  +  \lambda \mathbf{D}^{T} \mathbf{D} +  \mathbf{N} \mathbf{N}^{T} \right) \ ,
\end{align*}
are all positive definite and the associated the Levenberg-Marquardt directions are thus well-defined. Furthermore, in these conditions, if $\mathbf{C}$ represents any of these  positive definite matrices, we have immediately
\begin{equation*}
d\mathbf{a}_{lm}^{T} \nabla \psi( \mathbf{a} )  =  d\mathbf{a}_{lm}^{T}  \mathit{J}\big( \mathbf{r}(\mathbf{a}) \big)^{T} \mathbf{r}(  \mathbf{a} )  = - d\mathbf{a}_{lm}^{T} \mathbf{C} d\mathbf{a}_{lm} < 0 \ ,
\end{equation*}
since the hypothesis $\nabla \psi( \mathbf{a} ) \ne \mathbf{0}^{p.k}$ implies that $d\mathbf{a}_{lm} \ne \mathbf{0}^{p.k}$ and $\mathbf{C}$ is positive definite. In other words, all the associated  Levenberg-Marquardt directions are in a descent direction for $\psi(.)$ if $\lambda \ne 0$, which proves the first part of the Corollary.

On the other hand, if $\lambda = 0$ and $\emph{rank}( \mathit{J}( \mathbf{r}(\mathbf{a}) ) ) = \emph{rank}( \mathbf{M}(\mathbf{a}) ) = (p - k).k$, we have
\begin{align*}
d\mathbf{a}_{lm} & = -  \left(  \mathit{J}\big( \mathbf{r}(\mathbf{a}) \big)^{T} \mathit{J}\big( \mathbf{r}(\mathbf{a}) \big) +   \mathbf{N} \mathbf{N}^{T} \right)^{-1} \mathit{J}\big( \mathbf{r}(\mathbf{a}) \big)^{T} \mathbf{r}(  \mathbf{a} ) \\
 &  \quad \text{  or  }  \\
d\mathbf{a}_{lm} & = \left(  \mathbf{M}(\mathbf{a})^{T}  \mathbf{M}(\mathbf{a})  +  \mathbf{N} \mathbf{N}^{T} \right)^{-1}  \mathbf{M}(\mathbf{a})^{T} \mathbf{r}(  \mathbf{a} ) \ .
\end{align*}
Further, noting that the matrices
\begin{equation*}
\begin{bmatrix} \mathit{J}\big( \mathbf{r}(\mathbf{a}) \big) \\   \mathbf{N}^{T} \end{bmatrix} \text{ and } \begin{bmatrix}  \mathbf{M}(\mathbf{a}) \\   \mathbf{N}^{T} \end{bmatrix}
\end{equation*}
are of full column rank, as demonstrated in equation~\eqref{eq:N_mat} of Subsection~\ref{jacob:box}, we deduce that the matrices
\begin{equation*}
\left(  \mathit{J}\big( \mathbf{r}(\mathbf{a}) \big)^{T} \mathit{J}\big( \mathbf{r}(\mathbf{a}) \big) +   \mathbf{N} \mathbf{N}^{T} \right) \text{ and } \left(  \mathbf{M}(\mathbf{a})^{T}  \mathbf{M}(\mathbf{a})  +  \mathbf{N} \mathbf{N}^{T} \right)
\end{equation*}
are both positive definite. This implies that the corresponding  Levenberg-Marquardt directions  are again well defined and still in a descent direction for $\psi(.)$ by similar arguments as used in the first part of the demonstration.

Finally, the linear systems
\begin{align*}
\left(  \mathit{J}\big( \mathbf{r}(\mathbf{a}) \big)^{T} \mathit{J}\big( \mathbf{r}(\mathbf{a}) \big) +   \mathbf{N} \mathbf{N}^{T} \right) d\mathbf{a}_{lm} & = -  \mathit{J}\big( \mathbf{r}(\mathbf{a}) \big)^{T} \mathbf{r}(  \mathbf{a} )  \\
 \text{  and  }  & \\
\left(  \mathbf{M}(\mathbf{a})^{T}  \mathbf{M}(\mathbf{a})  +  \mathbf{N} \mathbf{N}^{T} \right) d\mathbf{a}_{lm}  & =  \mathbf{M}(\mathbf{a})^{T} \mathbf{r}(  \mathbf{a} )
\end{align*}
are, respectively, the normal equations of the linear least-squares problems
\begin{equation*}
\min_{d\mathbf{a} \in \mathbb{R}^{p.k}}   \,   \frac{1}{2} \big\Vert \begin{bmatrix} \mathbf{r}(  \mathbf{a} )  \\  \mathbf{0}^{k_{\mathbf{A}}.k}   \end{bmatrix} + \begin{bmatrix} \mathit{J}\big( \mathbf{r}(\mathbf{a}) \big) \\   \mathbf{N}^{T} \end{bmatrix} d\mathbf{a} \big\Vert^{2}_{2}
\  \text{ and } \
\min_{d\mathbf{a} \in \mathbb{R}^{p.k}}   \,   \frac{1}{2} \big\Vert \begin{bmatrix} \mathbf{r}(  \mathbf{a} )  \\  \mathbf{0}^{k_{\mathbf{A}}.k}   \end{bmatrix} - \begin{bmatrix} \mathbf{M}(\mathbf{a})  \\   \mathbf{N}^{T} \end{bmatrix} d\mathbf{a} \big\Vert^{2}_{2} \ ,
\end{equation*}
which both have an unique solution, as the associated coefficient matrices are of full column rank, and these solutions are, respectively, the minimum 2-norm solutions of the rank deficient linear least-squares problems
\begin{equation*}
\min_{d\mathbf{a} \in \mathbb{R}^{p.k}}   \,   \frac{1}{2} \big\Vert \mathbf{r}(  \mathbf{a} )  + \mathit{J}\big( \mathbf{r}(\mathbf{a}) \big) d\mathbf{a} \big\Vert^{2}_{2}
\ \text{ and } \
\min_{d\mathbf{a} \in \mathbb{R}^{p.k}}   \,   \frac{1}{2} \big\Vert \mathbf{r}(  \mathbf{a} ) - \mathbf{M}(\mathbf{a})   d\mathbf{a} \big\Vert^{2}_{2} \ ,
\end{equation*}
as demonstrated in Subsection~\ref{jacob:box}. In other words, we have
\begin{equation*}
d\mathbf{a}_{lm} = -  \mathit{J}\big( \mathbf{r}(\mathbf{a}) \big)^{+} \mathbf{r}(  \mathbf{a} ) = d\mathbf{a}_{gn}
\ \text{ or } \
d\mathbf{a}_{lm} = \mathbf{M}(\mathbf{a})^{+} \mathbf{r}(  \mathbf{a} ) = d\mathbf{a}_{gn} \ ,
\end{equation*}
if an approximate Jacobian matrix $-\mathbf{M}(\mathbf{a})$ is used. This concludes the demonstration of the Corollary.
\\
\end{proof}

We now explore the relationships between the Euclidean gradient $\nabla \psi( \mathbf{a} )$ of $\psi(.)$  at $\mathbf{a}$, considered as a real function from  $\mathbb{R}^{p.k}$ into $\mathbb{R}$ (e.g., of the vectorized form of the $\mathbf{A}$ matrix, see equation~\eqref{eq:psi_func}), and the  Riemannian gradient of the unvectorized form of  $\psi(.)$ (e.g.,  $\psi  \circ  h^{-1}(.)$ where $h^{-1}(.)$ is defined in equation~\eqref{eq:h_func} of Subsection~\ref{varpro_wlra:box} with $h^{-1} ( \mathbf{A} ) = \emph{vec}(  \mathbf{A}^{T} ) = \mathbf{a}, \forall \mathbf{A} \in \mathbb{R}^{p \times k}_{k}$), when this cost function is considered as defined on the Grassmann manifold $\text{Gr}(p,k)$, as already discussed at the end of Subsection~\ref{jacob:box}.

To this end, we require that  $\mathbf{W} \in  \mathbb{R}^{ p \times n  }_{+*}$ (as otherwise $\psi  \circ  h^{-1}(.)$ is not smooth on its whole domain  $\mathbb{R}^{ p \times k  }_{k}$) and that each element of $\text{Gr}(p,k)$ is represented by an element of the compact Stiefel manifold $\text{St}(p,k) = \mathbb{O}^{p \times k }$ instead of $\mathbb{R}^{ p \times k  }_{k}$ as it is customary for simplicity and numerical reasons in previous works on Riemannian optimization on $\text{Gr}(p,k)$~\cite{DKM2012}\cite{BA2015}\cite{B2023}. With these requirements, we first observe that the restriction of the mapping $\psi  \circ  h^{-1}(.)$ to the domain $\text{St}(p,k) = \mathbb{O}^{p \times k }$ can be considered as a smooth map defined on  the Stiefel manifold and, by extension, also on the Grassmann manifold as the Grassmann manifold  $\text{Gr}(p,k)$ is a Riemannian quotient manifold of $\text{St}(p,k)$ by the action of the orthogonal group $\mathbb{O}^{k \times k}$, see Subsections~\ref{calculus:box} and~\ref{jacob:box},  and~\cite{AMS2008}\cite{B2023} for more details. In these conditions, the Riemannian gradient of this smooth map defined on the Grassmann manifold $\text{Gr}(p,k)$ at $\mathbf{O} \in \text{St}(p,k)$, considered as the matrix representation of $\mathring{\mathbf{O}} = \emph{ran}( \mathbf{O} )  \in \text{Gr}(p,k)$, is an element of the tangent space of $\text{Gr}(p,k)$ at  $\mathring{\mathbf{O}}$, $\mathcal{T}_{\mathring{\mathbf{O}}}  \text{Gr}(p,k)$, which is a linear subspace of dimension $\emph{dim} \big( \text{Gr}(p,k) \big) = k.( p - k )$. Moreover, any element of $\mathcal{T}_{\mathring{\mathbf{O}}}  \text{Gr}(p,k)$ can be represented uniquely by a matrix $\mathbf{D} \in \mathbb{R}^{p \times k }$ verifying $\mathbf{O}^{T}\mathbf{D} = \mathbf{0}^{ k \times k}$, as already noted in Subsection~\ref{jacob:box}. This subset of $\mathbb{R}^{ p \times k  }$ is noted $\mathcal{T}_{\mathbf{O}}  \text{Gr}(p,k)$ and is nothing else than the  horizontal space $\mathcal{H}_{\mathbf{O}} \mathbb{O}^{p \times k }$ of $\text{St}(p,k) = \mathbb{O}^{p \times k }$ at $\mathbf{O} \in \text{St}(p,k)$.

Thus, the Riemannian gradient of $\psi  \circ  h^{-1}(.)$ at $\mathring{\mathbf{O}}  \in \text{Gr}(p,k)$, noted $\nabla_{R} \psi  \circ  h^{-1} (  \mathring{\mathbf{O}} )$, can be represented uniquely by one element of $\mathcal{T}_{\mathbf{O}}  \text{Gr}(p,k)$ and, according to equation~\eqref{eq:D_rgrad_vec2} in Subsection~\ref{calculus:box},  this Riemannian gradient is given, with a slight abuse of notation, by
\begin{align} \label{eq:R_Grad}
\nabla_{R} \psi  \circ  h^{-1} (  \mathring{\mathbf{O}} ) & =     \mathbf{P}_{\mathcal{H}_{\mathbf{O}} \mathbb{O}^{p \times k }} \nabla_{F} \psi  \circ  h^{-1} ( \mathbf{O} )   \nonumber \\ 
& =   ( \mathbf{I}_{p} - \mathbf{O}\mathbf{O}^{T}) \nabla_{F} \psi  \circ  h^{-1} ( \mathbf{O} )  \ ,
\end{align}
where  $\mathbf{O} \in \text{St}(p,k) = \mathbb{O}^{p \times k }$, $\mathring{\mathbf{O}} \in \text{Gr}(p,k)$, $  \mathbf{P}_{\mathcal{H}_{\mathbf{O}} \mathbb{O}^{p \times k }} = \mathbf{I}_{p} - \mathbf{O}\mathbf{O}^{T}$ is the orthogonal projector onto the orthogonal of the range of $\mathbf{O}$ in $\mathbb{R}^{p}$ (e.g., $\emph{ran}( \mathbf{O} )^{\bot}$) and also the orthogonal projector on the horizontal space of $\text{St}(p,k)=\mathbb{O}^{p \times k }$ at  $\mathbf{O}$ ( in the linear spcace $\mathbb{R}^{ p \times k  }$) and, finally, $\nabla_{F} \psi  \circ  h^{-1} ( \mathbf{O} )$ is the usual Frobenius gradient of $\psi  \circ  h^{-1}(.)$ at  $\mathbf{O}$ when $\psi  \circ  h^{-1}(.)$ is considered as a real function defined on the linear space $\mathbb{R}^{p \times k}$. See~\cite{AMS2008}\cite{BA2015}\cite{B2023} for a derivation of this standard result on the geometry of the Stiefel and Grassmann manifolds.

Next, we first observe that
\begin{equation*}
\nabla \psi( \mathbf{o} ) =  \emph{vec}  \Big( \big( \nabla_{F} \psi  \circ  h^{-1} ( \mathbf{O} ) \big)^{T} \Big)  \in \mathbb{R}^{p.k} \ ,
\end{equation*}
since following the conventions used throughout the monograph, we have defined $\mathbf{a} = \emph{vec}(  \mathbf{A}^{T} )$, $\forall \mathbf{A} \in \mathbb{R}^{p  \times k}$, instead of $\mathbf{a} = \emph{vec}(  \mathbf{A} )$ as it is commonly used in past works on the WLRA problem. In words, $\nabla \psi( \mathbf{o} )$ is the vectorized form of the transpose of the Frobenius gradient of $\psi  \circ  h^{-1}(.)$ at $\mathbf{O}$. Using a similar and consistent convention, we now define
\begin{equation*}
\nabla_{R} \psi( \mathbf{o} ) =  \emph{vec} \Big(  \big ( \nabla_{R} \psi  \circ  h^{-1} (  \mathring{\mathbf{O}} ) \big) ^{T} \Big)  \in \mathbb{R}^{p.k} \ ,
\end{equation*}
e.g., $\nabla_{R} \psi( \mathbf{o} )$ is the vectorized form of the transpose of the Riemannian gradient of $\psi  \circ  h^{-1}(.)$ at $ \mathring{\mathbf{O}} \in \text{Gr}(p,k)$. In these conditions, the equality~\eqref{eq:R_Grad} defining the Riemannian gradient of $\psi  \circ  h^{-1}(.)$ at  $ \mathring{\mathbf{O}}$ is equivalent to the matrix equality
\begin{equation*}
 \big( \nabla_{R} \psi  \circ  h^{-1} (  \mathring{\mathbf{O}} ) \big)^{T} = \big( \nabla_{F} \psi  \circ  h^{-1} ( \mathbf{O} ) \big)^{T} \big( \mathbf{I}_{p} - \mathbf{O}\mathbf{O}^{T} \big) \ ,
\end{equation*}
as an orthogonal projector is a symmetric matrix (see equation~\eqref{eq:projector}), and also to the vector equality
\begin{align*}
\emph{vec} \Big(  \big( \nabla_{R} \psi  \circ  h^{-1} (  \mathring{\mathbf{O}} )  \big)^{T}  \Big)  &  = \emph{vec} \Big(  \big( \nabla_{F} \psi  \circ  h^{-1} ( \mathbf{O} ) \big)^{T}  ( \mathbf{I}_{p} - \mathbf{O}\mathbf{O}^{T} )  \Big)   \\
                                                                                                              &  =  \big(  ( \mathbf{I}_{p} - \mathbf{O}\mathbf{O}^{T} )     \otimes    \mathbf{I}_{k}   \big)    \emph{vec} \Big(   \big( \nabla_{F} \psi  \circ  h^{-1} ( \mathbf{O} )  \big)^{T}  \Big)   \ ,
\end{align*}
where we have used again the symmetry of the orthogonal projector and equation~\eqref{eq:vec_kronprod}. In other words, the vectorized form of the Riemannnian gradient of $\psi  \circ  h^{-1}(.)$ at $ \mathring{\mathbf{O}}$ is given by
\begin{equation*}
\nabla_{R} \psi ( \mathbf{o} ) =  \big(  ( \mathbf{I}_{p} - \mathbf{O}\mathbf{O}^{T}  )     \otimes    \mathbf{I}_{k}   \big) \nabla \psi ( \mathbf{o} ) \ .
\end{equation*}
Before proceeding further, we now precise  the nature of the  $p.k \times p.k$ symmetric matrix $( \mathbf{I}_{p} - \mathbf{O}\mathbf{O}^{T} )     \otimes    \mathbf{I}_{k}$. More precisely, we show that this symmetric matrix is the orthogonal projector onto $\emph{null}( \mathit{J}( \mathbf{r}(\mathbf{o}) ) )^{\bot} = \emph{null}( \mathbf{M}(\mathbf{o}) )^{\bot}$ when $\mathbf{W} \in  \mathbb{R}^{ p \times n  }_{+*}$ and $\mathbf{O} \in  \text{St}(p,k) = \mathbb{O}^{p \times k }$. To this end, we first recall that, in these conditions and according to Corollary~\ref{corol5.6:box}, the columns of the matrix $\mathbf{\bar{O}} =  \mathbf{K}_{(p,k)} ( \mathbf{I}_{k}  \otimes  \mathbf{O})$ form an orthonormal basis of $\emph{null}\big( \mathit{J}( \mathbf{r}(\mathbf{o}) ) \big) = \emph{null}( \mathbf{M}(\mathbf{o}) )$ and that the columns of $\mathbf{\bar{O}}^{\bot} =  \mathbf{K}_{(p,k)} ( \mathbf{I}_{k}  \otimes  \mathbf{O}^{\bot})$ are an orthonormal basis of $\emph{null}\big( \mathit{J}( \mathbf{r}(\mathbf{o}) ) \big)^{\bot} = \emph{null}( \mathbf{M}(\mathbf{o}) )^{\bot}$ as soon as the columns of $\mathbf{O}$ and $\mathbf{O}^{\bot}$ are orthonormal bases of $\emph{ran}( \mathbf{O} )$ and $\emph{ran}( \mathbf{O} )^{\bot}$, respectively. Next, using equation~\eqref{eq:bilin_kronprod} in Subsection~\ref{multlin_alg:box}, we have
\begin{equation*}
( \mathbf{I}_{p} - \mathbf{O}\mathbf{O}^{T}  )     \otimes    \mathbf{I}_{k} = \mathbf{I}_{p.k} -  ( \mathbf{O}\mathbf{O}^{T}   \otimes    \mathbf{I}_{k} ) \ .
\end{equation*}
Furthermore, using the identities~\eqref{eq:mulmat_kronprod} and~\eqref{eq:com_kron} again in Subsection~\ref{multlin_alg:box}, we deduce that
\begin{align*}
\mathbf{O}\mathbf{O}^{T}   \otimes    \mathbf{I}_{k} & = ( \mathbf{O}   \otimes    \mathbf{I}_{k} ) (  \mathbf{O}^{T}   \otimes    \mathbf{I}_{k} )   \\
                                                                                  & =  ( \mathbf{O}   \otimes    \mathbf{I}_{k} )  \mathbf{K}_{(k,k)} \mathbf{K}_{(k,k)}  (  \mathbf{O}^{T}   \otimes    \mathbf{I}_{k} )  \\
                                                                                  & =  \big( \mathbf{K}_{(p,k)} (   \mathbf{I}_{k} \otimes  \mathbf{O}  ) \big)  \big( \mathbf{K}_{(p,k)} (   \mathbf{I}_{k} \otimes  \mathbf{O}  ) \big)^{T}   \\
                                                                                  & =  \mathbf{\bar{O}} \mathbf{\bar{O}}^{T}  \ .
\end{align*}
Thus, $\mathbf{O}\mathbf{O}^{T} \otimes    \mathbf{I}_{k} = \mathbf{\bar{O}} \mathbf{\bar{O}}^{T}$ is the orthogonal projector onto $\emph{null}\big( \mathit{J}( \mathbf{r}(\mathbf{o}) ) \big) = \emph{null}( \mathbf{M}(\mathbf{o}) )$ and
\begin{equation*}
\mathbf{I}_{p.k} - ( \mathbf{O}\mathbf{O}^{T}   \otimes    \mathbf{I}_{k} ) = \mathbf{I}_{p.k} - \mathbf{\bar{O}} \mathbf{\bar{O}}^{T} = \mathbf{\bar{O}}^{\bot} ( \mathbf{\bar{O}}^{\bot} )^{T} 
\end{equation*}
is the orthogonal projector onto $\emph{null}\big( \mathit{J}( \mathbf{r}(\mathbf{o}) ) \big)^{\bot} = \emph{null}( \mathbf{M}(\mathbf{o}) )^{\bot}$ as stated above. Using these results, Theorem~\ref{theo5.7:box}, and remembering that $\mathbf{\bar{O}}$ is an orthogonal basis of $\emph{null}( \mathbf{M}(\mathbf{o}) )$ and, thus, that $\mathbf{M}(\mathbf{o}) \mathbf{\bar{O}} = \mathbf{0}^{p.n \times k.k}$ , we deduce that
\begin{align*}
\nabla_{R} \psi ( \mathbf{o} )  & =  ( \mathbf{I}_{p.k} - \mathbf{\bar{O}} \mathbf{\bar{O}}^{T}  ) \nabla \psi ( \mathbf{o} )  \\
                                               & =  \nabla \psi ( \mathbf{o} ) - \mathbf{\bar{O}} \mathbf{\bar{O}}^{T}  \nabla \psi ( \mathbf{o} )   \\
                                               & =    \nabla \psi ( \mathbf{o} ) +  \mathbf{\bar{O}} \mathbf{\bar{O}}^{T}  \mathbf{M}(\mathbf{o})^{T} \mathbf{r}(\mathbf{o})   \\ 
                                               & =    \nabla \psi ( \mathbf{o} ) +  \mathbf{\bar{O}} \big(  \mathbf{M}(\mathbf{o}) \mathbf{\bar{O}} \big)^{T} \mathbf{r}(\mathbf{o})   \\ 
                                               & =    \nabla \psi ( \mathbf{o} ) \ . 
\end{align*}
From the vectorized equality $\nabla_{R} \psi ( \mathbf{o} ) =  \nabla \psi ( \mathbf{o} )$, by identification and unicity of the Frobenius (or Riemaniann) gradient, we deduce finally that
\begin{equation*}
\nabla_{R} \psi  \circ  h^{-1} ( \mathring{\mathbf{O}} ) =  \nabla_{F} \psi  \circ  h^{-1} ( \mathbf{O} )  \ .
\end{equation*}
In other words, the Riemannian gradient of $\psi  \circ  h^{-1} (.)$ at $\mathring{\mathbf{O}}$ is simply equal to its Frobenius gradient at $\mathbf{O}$, which again illustrates the tight relationships between the variable projection method used here and the Riemannian optimization framework on the Grassmann manifold developed in Boumal and Absil~\cite{BA2011}\cite{BA2015} despite the derivations are completely different.

To conclude this paragraph on the comparison of the Euclidean and Riemannian gradients of the cost function $\psi(.)$, we now derive an explicit matrix form of $\nabla_{F}  \psi \circ  h^{-1} ( \mathbf{A} )$ for $\mathbf{A} \in  \mathbb{R}^{ p \times k }_{k}$ starting from the vectorized equality
\begin{equation*}
\nabla \psi ( \mathbf{a} ) =  \emph{vec} \Big(  \big( \nabla_{F} \psi  \circ  h^{-1} ( \mathbf{A} )  \big)^{T}  \Big)  \ .
\end{equation*}
Using Theorems~\ref{theo4.3:box} and~\ref{theo5.7:box}, we first recall that
\begin{equation*}
\nabla \psi ( \mathbf{a} ) =   \mathbf{G}(\widehat{\mathbf{b}})^{T} \big( \mathbf{G}(\widehat{\mathbf{b}})\mathbf{a} -  \mathbf{z}  \big) ,
\end{equation*}
where, as usual,
\begin{align*}
   \mathbf{a}                                    & = \emph{vec}(  \mathbf{A}^{T} ) \ ,  \\
   \widehat{\mathbf{b}}                    & = \mathbf{F}(\mathbf{a})^{+} \mathbf{x} = \Big( \emph{diag}\big( \emph{vec}( \sqrt{\mathbf{W}} ) \big)  \big(  \mathbf{I}_n  \otimes \mathbf{A}  \big) \Big)^{+} \mathbf{x} \ ,   \\
   \mathbf{z}                                    & = \emph{vec} \big( (\sqrt{\mathbf{W}}  \odot \mathbf{X})^{T} \big) =  \emph{diag} \big( \emph{vec}( \sqrt{ \mathbf{W} }^{T} ) \big)   \emph{vec} ( \mathbf{X}^{T} )  \  ,    \\
  \mathbf{G}(\widehat{\mathbf{b}}) & =  \emph{diag} \big( \emph{vec}( \sqrt{ \mathbf{W} }^{T} ) \big) (  \mathbf{I}_{p}  \otimes \widehat{\mathbf{B}}^{T} ) \ .
\end{align*}
Using these results and equation~\eqref{eq:vec_kronprod} in Subsection~\ref{multlin_alg:box}, we deduce that
\begin{align*}
\nabla \psi ( \mathbf{a} )  & =  (  \mathbf{I}_{p}  \otimes \widehat{\mathbf{B}} )  \emph{diag} \big( \emph{vec}( \mathbf{W}^{T} ) \big)  \big(  (  \mathbf{I}_{p}  \otimes \widehat{\mathbf{B}}^{T} ) \emph{vec}(  \mathbf{A}^{T} ) -   \emph{vec} ( \mathbf{X}^{T} )  \big)  \\
                                        & =  (  \mathbf{I}_{p}  \otimes \widehat{\mathbf{B}} )  \emph{diag} \big( \emph{vec}( \mathbf{W}^{T} ) \big)  \big(  \emph{vec}(  \widehat{\mathbf{B}}^{T} \mathbf{A}^{T} ) -   \emph{vec} ( \mathbf{X}^{T} )  \big) \\
                                        & =  (  \mathbf{I}_{p}  \otimes \widehat{\mathbf{B}} )  \emph{diag} \big( \emph{vec}( \mathbf{W}^{T} ) \big)  \emph{vec}(  \widehat{\mathbf{B}}^{T} \mathbf{A}^{T}  -  \mathbf{X}^{T} )   \\
                                        & =  (  \mathbf{I}_{p}  \otimes \widehat{\mathbf{B}} )  \emph{vec} \big(  \mathbf{W}^{T} \odot ( \mathbf{A}\widehat{\mathbf{B}}  -  \mathbf{X} )^{T}  \big)    \\
                                        & =  (  \mathbf{I}_{p}  \otimes \widehat{\mathbf{B}} )  \emph{vec} \Big(  \big( \mathbf{W} \odot ( \mathbf{A}\widehat{\mathbf{B}}  -  \mathbf{X} ) \big)^{T}  \Big)    \\
                                        & =    \emph{vec} \Big(  \widehat{\mathbf{B}} \big( \mathbf{W} \odot ( \mathbf{A}\widehat{\mathbf{B}}  -  \mathbf{X} ) \big)^{T}  \Big)  \ .
\end{align*}
By identification and the unicity of the (Frobenius) gradient and equation~\eqref{eq:D_partial_grad_varphi*} of Subsection~\ref{landscape_wlra:box}, we obtain, finally, the following unvectorized form of $\nabla \psi ( \mathbf{a} )$ as:
\begin{equation*}
\nabla_{F} \psi  \circ  h^{-1} ( \mathbf{A} ) =  \big( \mathbf{W} \odot ( \mathbf{A}\widehat{\mathbf{B}}  -  \mathbf{X} ) \big)   \widehat{\mathbf{B}}^{T} =  \nabla  \varphi^{*}_{\mathbf{A}}( \mathbf{A},\mathbf{\widehat{B}} )  \ ,
\end{equation*}
which is entirely consistent with the results in Theorem~\ref{theo5.7:box}.

Consistent with the fact that $\nabla_{R} \psi  \circ  h^{-1} ( \mathbf{O} ) = \nabla_{F} \psi  \circ  h^{-1} ( \mathbf{O} )$ derived above, when $\mathbf{O}  \in  \text{St}(p,k) = \mathbb{O}^{p \times k }$ and $\mathbf{W} \in  \mathbb{R}^{ p \times n  }_{+*}$,  it can be further verified that this formulae for $ \nabla_{F} \psi  \circ  h^{-1} ( \mathbf{A} )$ is equal to the Riemannian gradient $\text{grad} f( \mathbf{A} )$ (which corresponds to $\nabla_{R} \psi  \circ  h^{-1} ( \mathring{\mathbf{A}} )$ in our notations) derived in equations 23 and 25 of Boumal and Absil~\cite{BA2015}, if we assumed that $\mathbf{A}$ is an orthonormal matrix, $\mathbf{W} = \mathbf{W}_{\lambda} \in  \mathbb{R}^{ p \times n  }_{+*}$ and $\mathbf{X} = \mathbf{X}_{\Omega}$ as defined, respectively, in equations~\eqref{eq:weight_proj_op} and~\eqref{eq:proj_op}  of Subsection~\ref{approx_wlra:box}, e.g., when a regularization parameter $\lambda > 0$ is used in the Riemannian optimization method  on the Grassmann manifold developed by Boumal and Absil~\cite{BA2015} to minimize their cost function $f(.)$, which is nothing else than a variable projection form  of the cost function $g_{\lambda}(.)$ defined in equation~\eqref{eq:g_func} already discussed in Subsection~\ref{approx_wlra:box}.

We summarize these different results about the unvectorized forms of the Frobenius and Riemannian gradients of $\psi (.)$ in the following Theorem:
 \begin{theo5.8} \label{theo5.8:box}
Let  $\psi(.)$ and $h^{-1} (.)$ be defined, respectively, as in equations~\eqref{eq:psi_func} and~\eqref{eq:h_func}  of Subsection~\ref{varpro_wlra:box},  $\mathbf{A} \in \mathbb{R}^{ p \times k }_{k}$ and $\mathbf{a} = h^{-1} ( \mathbf{A} ) = \emph{vec}( \mathbf{A}^{T} ) \in \mathbb{R}^{ p.k }$, and further assume that $\mathbf{a}$ belongs to an open set $ \varOmega \subset  \mathbb{R}^{p.k} $  in which the matrix function $ \mathbf{F}( \mathbf{a} ) $, defined in equation~\eqref{eq:F_mat}, has a constant rank, so that the Jacobian matrix $ \mathit{J}( \mathbf{r}(\mathbf{a}) ) $ derived in Subsection~\ref{jacob:box}, the Euclidean gradient $\nabla \psi( \mathbf{a} )$ derived in Theorem~\ref{theo5.7:box} and the Frobenius gradient $\nabla_{F} \psi  \circ  h^{-1} ( \mathbf{A} ) $ are all well-defined. Then
\begin{equation*}
\nabla_{F} \psi  \circ  h^{-1} ( \mathbf{A} ) =  \big( \mathbf{W} \odot ( \mathbf{A}\widehat{\mathbf{B}}  -  \mathbf{X} ) \big)   \widehat{\mathbf{B}}^{T} =  \nabla  \varphi^{*}_{\mathbf{A}}( \mathbf{A},\mathbf{\widehat{B}} )  \ ,
\end{equation*}
where $\mathbf{X}  \in \mathbb{R}^{ p \times n }$  and $\mathbf{W}  \in \mathbb{R}^{ p \times n }_{+}$ are, respectively, the data and weight matrices of the WLRA problem~\eqref{eq:P1} and $\mathbf{\widehat{B}}  \in \mathbb{R}^{ k \times n }$ is such that $\emph{vec}( \mathbf{\widehat{B}} ) =  \mathbf{F}(\mathbf{a})^{+} \mathbf{x}$, with $\mathbf{x} = \emph{vec}(  \sqrt{ \mathbf{W} } \odot \mathbf{X} ) $.

Furthermore, if we assume that $\mathbf{W}  \in \mathbb{R}^{ p \times n }_{+*}$, $\mathbf{O}   \in  \text{St}(p,k) = \mathbb{O}^{p \times k }$ and $\mathring{\mathbf{O}} = \emph{ran}( \mathbf{O} )  \in \text{Gr}(p,k)$, we have
\begin{equation*}
\nabla_{R} \psi  \circ  h^{-1} ( \mathring{\mathbf{O}} ) =  \nabla_{F} \psi  \circ  h^{-1} ( \mathbf{O} ) \ ,
\end{equation*}
where $\nabla_{R} \psi  \circ  h^{-1} (  \mathring{\mathbf{O}} )$ is the Riemannian gradient of the smooth cost function $\psi \circ  h^{-1}(.)$ (defined now on the Grassmann manifold  $\text{Gr}(p,k)$) at $\mathring{\mathbf{O}} \in \text{Gr}(p,k)$.

$\Box$
\end{theo5.8}

To conclude, we highlight that the above results about the comparison of the gradients in the variable projection and Grassmann manifold optimization frameworks are a slight extension of results derived in a different way in~\cite{HF2015b}, who do not consider the case where a regularization parameter $\lambda$ is used in the Boumal and Absil algorithms~\cite{BA2011}\cite{BA2015}.

We now derive an explicit form for the Hessian matrix $\mathbf{H} = \nabla^2 \psi (\mathbf{a})$ where $\mathbf{a} = h^{-1} ( \mathbf{A} ) = \emph{vec}( \mathbf{A}^{T} ) \in \mathbb{R}^{ p.k }$ with $\mathbf{A} \in \mathbb{R}^{ p \times k }_{k}$. To this end, we follow the derivation of the Newton step presented in Borges~\cite{B2009} for solving a general variable projection NLLS problem. In the context of the WLRA problem, this Hessian matrix has been already explicitly derived in different forms by several authors~\cite{C2008b}\cite{SE2010}\cite{BA2015}, but these forms apply mainly to the case of binary weights or are not convenient to derive important properties of $\nabla^2 \psi (\mathbf{a})$. As an illustration, Chen~\cite{C2008b} has derived $\mathbf{H}$ (see its equation 25) for the case where $\mathbf{W}$ is a presence-absence matrix by expressing $\psi(\mathbf{a})$ as the sum of the $j^{th}$ atomic functions $\psi_j (.)$ defined in equation~\eqref{eq:psi_atomic_func} of Subsection~\ref{varpro_wlra:box}:
\begin{align*}
    \psi(\mathbf{a}) & =    \frac{ 1 }{ 2 } \sum_{j=1}^{n} { \psi_j (\mathbf{a}) }   \\
                              & =   \frac{ 1 }{ 2 } \sum_{j=1}^{n} \big\Vert   \big( \mathbf{I}_p - \mathbf{F}_{j} ( \mathbf{a} )\mathbf{F}_{j} ( \mathbf{a} )^{+}  \big) \mathbf{x}_{j}  \big\Vert^{2}_{2}  \\                                                                                       
                              & =   \frac{ 1 }{ 2 } \sum_{j=1}^{n} \big\Vert   \left( \mathbf{I}_p - \left(  \emph{diag}( \sqrt{ \mathbf{W} }_{.j}) \mathbf{A} \right) \left(  \emph{diag}(\sqrt{ \mathbf{W} }_{.j})  \mathbf{A} \right)^{+} \right) ( \sqrt{ \mathbf{W} }_{.j} \odot \mathbf{X}_{.j} ) \big\Vert^{2}_{2} \,
\end{align*}
and computing the Hessian matrices of $ \psi_j ( . ) $ at $\mathbf{a}$, for $j=1, \cdots, n $, and summing all these Hessian matrices. However, this formulation of $\nabla^2 \psi (\mathbf{a})$  does not give a compact form for  $\mathbf{H}$ and is not convenient to derive important properties of the Hessian matrix later in this subsection. The two other formulations of $\mathbf{H}$ derived in~\cite{SE2010}\cite{BA2015} apply, respectively, only to the binary weights case or only to the variable projection formulation of the regularized cost function $g_\lambda(.)$, defined in equation~\eqref{eq:g_func} of Subsection~\ref{approx_wlra:box}, and this one lacks also of generality as $\mathbf{H}$ is not obtained explicitly, but in the form of a directional derivative~\cite{BA2015}, which is used in an inexact subproblem solver based a truncated conjugate gradient method to compute an approximate (Riemannian) Newton step at each iteration~\cite{BA2015}.

To start with, we recall that the Hessian matrix is given formally  by
\begin{equation} \label{eq:hess_mat}
      \mathbf{H}  = \mathit{J} \big( \mathbf{r}(\mathbf{a}) \big)^{T} \mathit{J} \big( \mathbf{r}(\mathbf{a}) \big) + \sum_{l=1}^{n.p} \mathbf{r}_l (\mathbf{a}) \nabla^2 \mathbf{r}_l (\mathbf{a}) \ ,
\end{equation}
where $ \nabla^2 \mathbf{r}_l (\mathbf{a}) $ is the Hessian matrix of the functional $ \mathbf{r}_l (\mathbf{a})$ for $l=1, \cdots, n.p$ and is a $ p.k \times p.k $ matrix given by
\begin{equation*}
    \lbrack \nabla^2 \mathbf{r}_l (\mathbf{a}) \rbrack_{ij} = \frac{ \partial^2  \mathbf{r}_l (\mathbf{a}) }{ \partial \mathbf{a}_i  \partial \mathbf{a}_j } \text{ for } i=1, \cdots, p.k \text{ and } j=1, \cdots, p.k \ .
\end{equation*}
We first show that the calculation of the first term of $ \nabla^2 \psi (\mathbf{a}) $, involving only the Jacobian matrix $\mathit{J}( \mathbf{r}(\mathbf{a}) )$ and which corresponds exactly to the Gauss-Newton approximation of $ \nabla^2 \psi (\mathbf{a}) $, can be simplified as for the gradient of $\psi (.)$. Using the notations and equation~\eqref{eq:J_mat}  of Subsection~\ref{jacob:box}, we have
\begin{equation*}
    \mathit{J}( \mathbf{r}(\mathbf{a}) ) = - \big( \mathbf{M}(\mathbf{a}) + \mathbf{L}(\mathbf{a}) \big) \ ,
\end{equation*}
and it follows that
\begin{align} \label{eq:gn_hess_mat}
    \mathit{J} \big( \mathbf{r}(\mathbf{a}) \big)^{T} \mathit{J} \big( \mathbf{r}(\mathbf{a}) \big) \nonumber
               & =  \big( \mathbf{M}(\mathbf{a}) + \mathbf{L}(\mathbf{a}) \big)^{T} \big( \mathbf{M}(\mathbf{a}) + \mathbf{L}(\mathbf{a}) \big)  \nonumber  \\
               & =  \mathbf{M}(\mathbf{a})^{T} \mathbf{M}(\mathbf{a}) + \mathbf{M}(\mathbf{a})^{T} \mathbf{L}(\mathbf{a}) + \mathbf{L}(\mathbf{a})^{T}  \mathbf{M}(\mathbf{a}) + \mathbf{L}(\mathbf{a})^{T}\mathbf{L}(\mathbf{a}) \nonumber \\
               & =  \mathbf{M}(\mathbf{a})^{T} \mathbf{M}(\mathbf{a}) +  \mathbf{L}(\mathbf{a})^{T}\mathbf{L}(\mathbf{a}) \ ,
\end{align}
since $ \emph{ran}( \mathbf{M}(\mathbf{a}) ) \subset \emph{ran}( \mathbf{F}(\mathbf{a}) )^{\bot} $ and $ \emph{ran}( \mathbf{L}(\mathbf{a}) ) \subset \emph{ran}( \mathbf{F}(\mathbf{a}) )$, see Subsection~\ref{jacob:box} for more details.

We now derive an explicit expression for the second  term of the Hessian matrix,  given formally by
\begin{equation} \label{eq:S_hess_mat}
    \mathbf{S} = \sum_{l=1}^{n.p} \mathbf{r}_l (\mathbf{a}) \nabla^2 \mathbf{r}_l (\mathbf{a}) \ .
\end{equation}
To this end, we first recall that the Jacobian matrix may also be expressed as
\begin{equation*}
    \mathit{J}( \mathbf{r}(\mathbf{a}) ) = - \mathit{D}( \mathbf{P}_{\mathbf{F}(\mathbf{a})} ) \mathbf{x} \ ,
\end{equation*}
where $\mathbf{x} =  \emph{vec}(  \sqrt{\mathbf{W}}  \odot \mathbf{X})$  (see equation~\eqref{eq:r_jacob} in Subsection~\ref{jacob:box}) and this implies that the $i^{th}$ column of $ \mathit{J}( \mathbf{r}(\mathbf{a}) ) $ is given by
\begin{equation*}
   \lbrack \mathit{J}( \mathbf{r}(\mathbf{a}) ) \rbrack_{.j} = - \frac{ \partial \mathbf{P}_{\mathbf{F}(\mathbf{a})} }{ \partial \mathbf{a}_i }\mathbf{x} = \frac{ \partial  \mathbf{r} (\mathbf{a}) }{ \partial \mathbf{a}_i } \ ,
\end{equation*}
and it follows that
\begin{equation*}
    \frac{ \partial^2  \mathbf{r} (\mathbf{a}) }{ \partial \mathbf{a}_i \partial \mathbf{a}_j } = - \frac{ \partial^2 \mathbf{P}_{\mathbf{F}(\mathbf{a})} }{ \partial \mathbf{a}_i \partial \mathbf{a}_j }\mathbf{x}  \text{  for  } i=1, \cdots, k.p \text{ and }  j=1, \cdots, k.p \  .
\end{equation*}
Using this last equality and the definition of $\mathbf{S}$, it is not difficult to see that the $ij$ element of $ \mathbf{S} $ is given by
\begin{equation*}
    \mathbf{S}_{ij} = - \mathbf{r}(\mathbf{a})^{T} \frac{ \partial^2 \mathbf{P}_{\mathbf{F}(\mathbf{a})} }{ \partial \mathbf{a}_i \partial \mathbf{a}_j }\mathbf{x} \text{  for  } i=1, \cdots, k.p \text{ and }  j=1, \cdots, k.p \  .
\end{equation*}
Thus, in order to evaluate $ \mathbf{S} $, we need an explicit expression for all the second partial derivatives of $\mathbf{P}_{\mathbf{F}(\mathbf{a})}$
\begin{equation*}
\frac{ \partial^2 \mathbf{P}_{\mathbf{F}(\mathbf{a})} }{ \partial \mathbf{a}_i \partial \mathbf{a}_j } \in \pounds \big( \mathbb{R}, \pounds ( \mathbb{R}, \mathbb{R}^{n.p \times n.p} ) \big) \text{  for  } i=1, \cdots, k.p \text{ and }  j=1, \cdots, k.p \ .
\end{equation*}
From the expression of $\mathit{D}( \mathbf{P}_{\mathbf{F}(\mathbf{a})} )$ given in Corollary~\ref{corol5.1:box}, it follows that
\begin{align} \label{eq:part_dev_P_F(a)}
    \frac{ \partial \mathbf{P}_{\mathbf{F}(\mathbf{a}) } }{ \partial \mathbf{a}_i }
         & =  \mathbf{P}_{\mathbf{F}(\mathbf{a})}^{\bot}  \frac{  \partial  \mathbf{F}(\mathbf{a}) }{ \partial \mathbf{a}_i }  \mathbf{F}(\mathbf{a})^{+}   +  \Big( \mathbf{P}_{\mathbf{F}(\mathbf{a})}^{\bot}  \frac{  \partial  \mathbf{F}(\mathbf{a}) }{ \partial \mathbf{a}_i } \mathbf{F}(\mathbf{a})^{+} \Big)^{T}  \nonumber  \\
         & =  \mathbf{P}_{\mathbf{F}(\mathbf{a})}^{\bot}  \frac{  \partial  \mathbf{F}(\mathbf{a}) }{ \partial \mathbf{a}_i }  \mathbf{F}(\mathbf{a})^{+}   +  \big(\mathbf{F}(\mathbf{a})^{+} \big)^{T} \frac{  \partial  \mathbf{F}(\mathbf{a})^{T} }{ \partial \mathbf{a}_i } \mathbf{P}_{\mathbf{F}(\mathbf{a})}^{\bot} \ .
\end{align}
Using the product rule to take the partial derivative of $ \partial \mathbf{P}_{\mathbf{F}(\mathbf{a}) } / { \partial \mathbf{a}_i } $ with respect to $\mathbf{a}_j$~\cite{C2017}, we obtain
\begin{align*}
    \frac{ \partial^2 \mathbf{P}_{\mathbf{F}(\mathbf{a}) } }{ \partial \mathbf{a}_i \partial \mathbf{a}_j } & =  \frac{ \partial \mathbf{P}_{\mathbf{F}(\mathbf{a})}^{\bot} } { \partial \mathbf{a}_j } \frac{  \partial  \mathbf{F}(\mathbf{a}) }{ \partial \mathbf{a}_i }  \mathbf{F}(\mathbf{a})^{+}   & + \mathbf{P}_{\mathbf{F}(\mathbf{a})}^{\bot}  \frac{  \partial^2  \mathbf{F}(\mathbf{a}) }{ \partial \mathbf{a}_i \partial \mathbf{a}_j }  \mathbf{F}(\mathbf{a})^{+} & + \mathbf{P}_{\mathbf{F}(\mathbf{a})}^{\bot}  \frac{  \partial  \mathbf{F}(\mathbf{a}) }{ \partial \mathbf{a}_i }  \frac{ \partial \mathbf{F}(\mathbf{a})^{+} } { \partial \mathbf{a}_j }    \\
    & \: + \Big( \frac{ \partial \mathbf{P}_{\mathbf{F}(\mathbf{a})}^{\bot} } { \partial \mathbf{a}_j } \frac{  \partial  \mathbf{F}(\mathbf{a}) }{ \partial \mathbf{a}_i }  \mathbf{F}(\mathbf{a})^{+}   & + \mathbf{P}_{\mathbf{F}(\mathbf{a})}^{\bot}  \frac{  \partial^2  \mathbf{F}(\mathbf{a}) }{ \partial \mathbf{a}_i \partial \mathbf{a}_j }  \mathbf{F}(\mathbf{a})^{+}  & + \mathbf{P}_{\mathbf{F}(\mathbf{a})}^{\bot}  \frac{  \partial  \mathbf{F}(\mathbf{a}) }{ \partial \mathbf{a}_i }  \frac{ \partial \mathbf{F}(\mathbf{a})^{+} } { \partial \mathbf{a}_j } \Big)^{T} \ .
\end{align*}
For the sake of brevity, we define
\begin{equation*}
    \mathbf{D}_{(i,j)} = \frac{ \partial \mathbf{P}_{\mathbf{F}(\mathbf{a})}^{\bot} } { \partial \mathbf{a}_j } \frac{  \partial  \mathbf{F}(\mathbf{a}) }{ \partial \mathbf{a}_i }  \mathbf{F}(\mathbf{a})^{+}   + \mathbf{P}_{\mathbf{F}(\mathbf{a})}^{\bot}  \frac{  \partial^2  \mathbf{F}(\mathbf{a}) }{ \partial \mathbf{a}_i \partial \mathbf{a}_j }  \mathbf{F}(\mathbf{a})^{+} + \mathbf{P}_{\mathbf{F}(\mathbf{a})}^{\bot}  \frac{  \partial  \mathbf{F}(\mathbf{a}) }{ \partial \mathbf{a}_i }  \frac{ \partial \mathbf{F}(\mathbf{a})^{+} } { \partial \mathbf{a}_j }  \ ,
\end{equation*}
and, thus,
\begin{equation*}
    \frac{ \partial^2 \mathbf{P}_{\mathbf{F}(\mathbf{a}) } }{ \partial \mathbf{a}_i \partial \mathbf{a}_j } = \mathbf{D}_{(i,j)} + \mathbf{D}_{(i,j)}^{T}  \ .
\end{equation*}
We will now simplify the computation of $ \mathbf{D}_{(i,j)} $. First, using the fact that
\begin{equation*}
    \frac{  \partial \mathbf{P}_{\mathbf{F}(\mathbf{a}) }^{\bot} } { \partial \mathbf{a}_j } = - \frac{  \partial \mathbf{P}_{\mathbf{F}(\mathbf{a}) } } { \partial \mathbf{a}_j }  \ ,
\end{equation*}
we deduce that
\begin{equation*}
    \mathbf{D}_{(i,j)} = - \frac{ \partial \mathbf{P}_{\mathbf{F}(\mathbf{a})} } { \partial \mathbf{a}_j } \frac{  \partial  \mathbf{F}(\mathbf{a}) }{ \partial \mathbf{a}_i }  \mathbf{F}(\mathbf{a})^{+}   + \mathbf{P}_{\mathbf{F}(\mathbf{a})}^{\bot}  \frac{  \partial^2  \mathbf{F}(\mathbf{a}) }{ \partial \mathbf{a}_i \partial \mathbf{a}_j }  \mathbf{F}(\mathbf{a})^{+} + \mathbf{P}_{\mathbf{F}(\mathbf{a})}^{\bot}  \frac{  \partial  \mathbf{F}(\mathbf{a}) }{ \partial \mathbf{a}_i }  \frac{ \partial \mathbf{F}(\mathbf{a})^{+} } { \partial \mathbf{a}_j }  \ .
\end{equation*}
Next, we recall that $ \mathbf{F}(.) $ is a linear transformation from $ \mathbb{R}^{k.p} $ into $ \mathbb{R}^{n.p \times n.k} $ (see equation~\eqref{eq:F_mat} in Subsection~\ref{varpro_wlra:box}). This implies, in particular, that $ \mathit{D}\big( \mathbf{F}(\mathbf{a}) \big) = \mathbf{F}(.) $, for all $ \mathbf{a} \in \mathbb{R}^{k.p} $, meaning that $ \mathit{D}\big( \mathbf{F}(.) \big) $ is a constant function from $ \mathbb{R}^{k.p} $ into $ \pounds ( \mathbb{R}^{k.p}, \mathbb{R}^{n.p \times n.k} ) $. From this, we deduce that $ \mathit{D}^2 \big( \mathbf{F}(\mathbf{a}) \big) $ is the null application for every $ \mathbf{a} \in \mathbb{R}^{k.p} $, which implies that all the second partial derivatives $ \partial^2 \mathbf{F}(\mathbf{a}) / \partial \mathbf{a}_i \partial \mathbf{a}_j $ are identically zero.
\\
\\
This may be used to simplify $ \mathbf{D}_{(i,j)} $ again, giving
\begin{equation*}
    \mathbf{D}_{(i,j)} = - \frac{ \partial \mathbf{P}_{\mathbf{F}(\mathbf{a})} } { \partial \mathbf{a}_j } \frac{  \partial  \mathbf{F}(\mathbf{a}) }{ \partial \mathbf{a}_i }  \mathbf{F}(\mathbf{a})^{+}  + \mathbf{P}_{\mathbf{F}(\mathbf{a})}^{\bot}  \frac{  \partial  \mathbf{F}(\mathbf{a}) }{ \partial \mathbf{a}_i }  \frac{ \partial \mathbf{F}(\mathbf{a})^{+} } { \partial \mathbf{a}_j } \ .
\end{equation*}
Now, we need to compute the derivative of $\mathbf{F}(\mathbf{a})^{+}$ which is given in Theorem~\ref{theo5.1:box} under the assumption that $\mathbf{F}(\mathbf{a})$ is of local constant rank at any point $\mathbf{a}$ in which differentiation is to be performed (see Theorem~\ref{theo5.1:box} for details):
\begin{align*}
    \mathit{D}( \mathbf{F}(\mathbf{a})^{+} )  =  & - \mathbf{F}(\mathbf{a})^{+} \mathit{D}( \mathbf{F}(\mathbf{a}) ) \mathbf{F}(\mathbf{a})^{+} + \mathbf{F}(\mathbf{a})^{+} (\mathbf{F}(\mathbf{a})^{+})^{T} \mathit{D}( \mathbf{F}(\mathbf{a})^{T} ) \mathbf{P}_{\mathbf{F}(\mathbf{a})}^{\bot}  \\
    & +\: \big( \mathbf{I}_{n.k} - \mathbf{F}(\mathbf{a})^{+} \mathbf{F}(\mathbf{a}) \big) \mathit{D}( \mathbf{F}(\mathbf{a})^{T} ) (\mathbf{F}(\mathbf{a})^{+})^{T} \mathbf{F}(\mathbf{a})^{+} \ .
\end{align*}
In order to simplify the derivation, we will now further assume that $ \mathbf{F}(\mathbf{a}) $ has full column-rank, leaving the general case for future research. For practical applications, this hypothesis implies that the rank of $ \mathbf{F}(\mathbf{a}) $ is $k.n$ and that each column of $ \mathbf{X} $ must have at least $k$ "nonmissing" elements. In other words, using the $p \times n$ incidence matrix, $\boldsymbol{\delta}$, associated with the matrix $\mathbf{X}$ (see equation~\eqref{eq:incidence_mat} in Subsection~\ref{jacob:box}), we must have
\begin{equation*}
    \sum_{l=1}^p  \boldsymbol{\delta}_{lj} >  k \quad , \text{for all } j=1, \cdots, n \ .
\end{equation*}
Note also that this hypothesis is automatically verified if $\mathbf{W}\in\mathbb{R}^{p \times n}_{+*}$ and $\mathbf{A}$ is of full column rank since
\begin{equation*}
    \mathbf{F}(  \mathbf{a} ) = \bigoplus_{j=1}^n \mathbf{F}_{j}(  \mathbf{a} ) \text{  where  } \mathbf{F}_{j}(  \mathbf{a} ) = \emph{diag}(\sqrt{\mathbf{W}}_{.j})\mathbf{A} \ .
\end{equation*}
From this hypothesis, we deduce that
\begin{equation*}
    \mathbf{F}(\mathbf{a})^{+} \mathbf{F}(\mathbf{a}) = \mathbf{I}_{n.k} \ ,
\end{equation*}
yielding a simplified formula for the derivative of $ \mathbf{F}(\mathbf{a})^{+} $
\begin{equation*}
    \mathit{D}( \mathbf{F}(\mathbf{a})^{+} )  =  - \mathbf{F}(\mathbf{a})^{+} \mathit{D}( \mathbf{F}(\mathbf{a}) ) \mathbf{F}(\mathbf{a})^{+} + \mathbf{F}(\mathbf{a})^{+} (\mathbf{F}(\mathbf{a})^{+})^{T} \mathit{D}( \mathbf{F}(\mathbf{a})^{T} ) \mathbf{P}_{\mathbf{F}(\mathbf{a})}^{\bot} 
\end{equation*}
and its partial derivatives
\begin{align*}
    \frac{ \partial  \mathbf{F}(\mathbf{a})^{+} } { \partial \mathbf{a}_j } =  - \mathbf{F}(\mathbf{a})^{+}  \frac{ \partial  \mathbf{F}(\mathbf{a}) } { \partial \mathbf{a}_j } \mathbf{F}(\mathbf{a})^{+} + \mathbf{F}(\mathbf{a})^{+} (\mathbf{F}(\mathbf{a})^{+})^{T}  \frac{ \partial  \mathbf{F}(\mathbf{a})^{T} } { \partial \mathbf{a}_j } \mathbf{P}_{\mathbf{F}(\mathbf{a})}^{\bot}  \ .
\end{align*}
Substituting this definition into $ \mathbf{D}_{(i,j)} $ gives
\begin{align*}
    \mathbf{D}_{(i,j)}  = & - \frac{ \partial \mathbf{P}_{\mathbf{F}(\mathbf{a})} } { \partial \mathbf{a}_j } \frac{  \partial  \mathbf{F}(\mathbf{a}) }{ \partial \mathbf{a}_i }  \mathbf{F}(\mathbf{a})^{+}  \\
    & +\: \mathbf{P}_{\mathbf{F}(\mathbf{a})}^{\bot}  \frac{  \partial  \mathbf{F}(\mathbf{a}) }{ \partial \mathbf{a}_i } \Big( - \mathbf{F}(\mathbf{a})^{+}  \frac{ \partial  \mathbf{F}(\mathbf{a}) } { \partial \mathbf{a}_j } \mathbf{F}(\mathbf{a})^{+} + \mathbf{F}(\mathbf{a})^{+} (\mathbf{F}(\mathbf{a})^{+})^{T}  \frac{ \partial  \mathbf{F}(\mathbf{a})^{T} } { \partial \mathbf{a}_j } \mathbf{P}_{\mathbf{F}(\mathbf{a})}^{\bot} \Big) \\
     = & - \frac{ \partial \mathbf{P}_{\mathbf{F}(\mathbf{a})} } { \partial \mathbf{a}_j } \frac{  \partial  \mathbf{F}(\mathbf{a}) }{ \partial \mathbf{a}_i }  \mathbf{F}(\mathbf{a})^{+}  \\
    & +\: \mathbf{P}_{\mathbf{F}(\mathbf{a})}^{\bot}  \frac{  \partial  \mathbf{F}(\mathbf{a}) }{ \partial \mathbf{a}_i } \mathbf{F}(\mathbf{a})^{+} \Big( (\mathbf{F}(\mathbf{a})^{+})^{T}  \frac{ \partial  \mathbf{F}(\mathbf{a})^{T} } { \partial \mathbf{a}_j } \mathbf{P}_{\mathbf{F}(\mathbf{a})}^{\bot} -  \frac{ \partial  \mathbf{F}(\mathbf{a}) } { \partial \mathbf{a}_j } \mathbf{F}(\mathbf{a})^{+}  \Big)  \ .
\end{align*}
Using this result, we are now in the position to deduce an explicit formula for $\mathbf{S}_{ij}$
\begin{align*}
    \mathbf{S}_{ij} & =  - \mathbf{r}(\mathbf{a})^{T} \big( \mathbf{D}_{(i,j)} + \mathbf{D}_{(i,j)}^{T} \big)  \mathbf{x}  \\
                            & =  - \mathbf{r}(\mathbf{a})^{T} \mathbf{D}_{(i,j)} \mathbf{x} - \mathbf{r}(\mathbf{a})^{T} \mathbf{D}_{(i,j)}^{T} \mathbf{x}  \\
                            & =  - \mathbf{r}(\mathbf{a})^{T} \mathbf{D}_{(i,j)} \mathbf{x} -  \mathbf{x}^{T} \mathbf{D}_{(i,j)} \mathbf{r}(\mathbf{a}) \ .
\end{align*}
Using the facts that $ \mathbf{r}(\mathbf{a}) = \mathbf{P}_{\mathbf{F}(\mathbf{a})}^{\bot} \mathbf{x} $, $ \mathbf{P}_{\mathbf{F}(\mathbf{a})}^{\bot} \mathbf{r}(\mathbf{a}) = \mathbf{r}(\mathbf{a}) $ and $ \mathbf{b} = \mathbf{F}(\mathbf{a})^{+} \mathbf{x} $, the first term in the right hand side of this equation reduces to
\begin{align*}
    - \mathbf{r}(\mathbf{a})^{T} \mathbf{D}_{(i,j)} \mathbf{x}  =  \ &  \mathbf{r}(\mathbf{a})^{T} \frac{ \partial \mathbf{P}_{\mathbf{F}(\mathbf{a})} } { \partial \mathbf{a}_j } \frac{  \partial  \mathbf{F}(\mathbf{a}) }{ \partial \mathbf{a}_i } \mathbf{b}   \\
                 & -\: \mathbf{r}(\mathbf{a})^{T} \frac{  \partial  \mathbf{F}(\mathbf{a}) }{ \partial \mathbf{a}_i } \mathbf{F}(\mathbf{a})^{+} \Big( (\mathbf{F}(\mathbf{a})^{+})^{T}  \frac{ \partial  \mathbf{F}(\mathbf{a})^{T} } { \partial \mathbf{a}_j } \mathbf{r}(\mathbf{a}) -  \frac{ \partial  \mathbf{F}(\mathbf{a}) } { \partial \mathbf{a}_j } \mathbf{b}  \Big)  \\
                  =  \ & \mathbf{r}(\mathbf{a})^{T} \frac{ \partial \mathbf{P}_{\mathbf{F}(\mathbf{a})} } { \partial \mathbf{a}_j } \frac{  \partial  \mathbf{F}(\mathbf{a}) }{ \partial \mathbf{a}_i } \mathbf{b}  -  \mathbf{r}(\mathbf{a})^{T} \frac{  \partial  \mathbf{F}(\mathbf{a}) }{ \partial \mathbf{a}_i } \mathbf{F}(\mathbf{a})^{+} (\mathbf{F}(\mathbf{a})^{+})^{T}  \frac{ \partial  \mathbf{F}(\mathbf{a})^{T} } { \partial \mathbf{a}_j } \mathbf{r}(\mathbf{a}) \\
                 & +\:  \mathbf{r}(\mathbf{a})^{T} \frac{  \partial  \mathbf{F}(\mathbf{a}) }{ \partial \mathbf{a}_i } \mathbf{F}(\mathbf{a})^{+} \frac{ \partial  \mathbf{F}(\mathbf{a}) } { \partial \mathbf{a}_j } \mathbf{b} \ .
\end{align*}
Substituting $ \partial \mathbf{P}_{\mathbf{F}(\mathbf{a})} / \partial \mathbf{a}_j  $ (see equation~\eqref{eq:part_dev_P_F(a)} above) yields
\begin{align*}
    - \mathbf{r}(\mathbf{a})^{T} \mathbf{D}_{(i,j)} \mathbf{x} = \ &  \mathbf{r}(\mathbf{a})^{T} \Big(  \mathbf{P}_{\mathbf{F}(\mathbf{a})}^{\bot}  \frac{  \partial  \mathbf{F}(\mathbf{a}) }{ \partial \mathbf{a}_j }  \mathbf{F}(\mathbf{a})^{+}   +  \big(\mathbf{F}(\mathbf{a})^{+} \big)^{T} \frac{  \partial  \mathbf{F}(\mathbf{a})^{T} }{ \partial \mathbf{a}_j } \mathbf{P}_{\mathbf{F}(\mathbf{a})}^{\bot} \Big) \frac{  \partial  \mathbf{F}(\mathbf{a}) }{ \partial \mathbf{a}_i } \mathbf{b}   \\
                 & -\:  \mathbf{r}(\mathbf{a})^{T} \frac{  \partial  \mathbf{F}(\mathbf{a}) }{ \partial \mathbf{a}_i } \mathbf{F}(\mathbf{a})^{+} (\mathbf{F}(\mathbf{a})^{+})^{T}  \frac{ \partial  \mathbf{F}(\mathbf{a})^{T} } { \partial \mathbf{a}_j } \mathbf{r}(\mathbf{a}) \\
                 & +\:  \mathbf{r}(\mathbf{a})^{T} \frac{  \partial  \mathbf{F}(\mathbf{a}) }{ \partial \mathbf{a}_i } \mathbf{F}(\mathbf{a})^{+} \frac{ \partial  \mathbf{F}(\mathbf{a}) } { \partial \mathbf{a}_j } \mathbf{b} \ .
\end{align*}
Noting that $ \mathbf{P}_{\mathbf{F}(\mathbf{a})}^{\bot}  \mathbf{r}(\mathbf{a}) = \mathbf{r}(\mathbf{a})$ and $ \mathbf{F}(\mathbf{a})^{+} \mathbf{r}(\mathbf{a}) = \mathbf{0}^{k.n} $, this reduces to
\begin{align*}
    - \mathbf{r}(\mathbf{a})^{T} \mathbf{D}_{(i,j)} \mathbf{x}  = \ &  \mathbf{r}(\mathbf{a})^{T}   \frac{  \partial  \mathbf{F}(\mathbf{a}) }{ \partial \mathbf{a}_j }  \mathbf{F}(\mathbf{a})^{+}  \frac{  \partial  \mathbf{F}(\mathbf{a}) }{ \partial \mathbf{a}_i } \mathbf{b}   \\
                 & -\:  \mathbf{r}(\mathbf{a})^{T} \frac{  \partial  \mathbf{F}(\mathbf{a}) }{ \partial \mathbf{a}_i } \mathbf{F}(\mathbf{a})^{+} (\mathbf{F}(\mathbf{a})^{+})^{T}  \frac{ \partial  \mathbf{F}(\mathbf{a})^{T} } { \partial \mathbf{a}_j } \mathbf{r}(\mathbf{a}) \\
                 & +\:  \mathbf{r}(\mathbf{a})^{T} \frac{  \partial  \mathbf{F}(\mathbf{a}) }{ \partial \mathbf{a}_i } \mathbf{F}(\mathbf{a})^{+} \frac{ \partial  \mathbf{F}(\mathbf{a}) } { \partial \mathbf{a}_j } \mathbf{b} \ .
\end{align*}
Next, if we consider the second term in $ \mathbf{S}_{ij} $, it reduces to
\begin{align*}
     -  \mathbf{x}^{T} \mathbf{D}_{(i,j)} \mathbf{r}(\mathbf{a})  = \ &  \mathbf{x}^{T} \frac{ \partial \mathbf{P}_{\mathbf{F}(\mathbf{a})} } { \partial \mathbf{a}_j } \frac{  \partial  \mathbf{F}(\mathbf{a}) }{ \partial \mathbf{a}_i } \mathbf{F}(\mathbf{a})^{+} \mathbf{r}(\mathbf{a})   \\
                 & -\: \mathbf{x}^{T} \mathbf{P}_{\mathbf{F}(\mathbf{a})}^{\bot} \frac{  \partial  \mathbf{F}(\mathbf{a}) }{ \partial \mathbf{a}_i } \mathbf{F}(\mathbf{a})^{+} (\mathbf{F}(\mathbf{a})^{+})^{T}  \frac{ \partial  \mathbf{F}(\mathbf{a})^{T} } { \partial \mathbf{a}_j } \mathbf{P}_{\mathbf{F}(\mathbf{a})}^{\bot}  \mathbf{r}(\mathbf{a})   \\
                & +\: \mathbf{x}^{T} \mathbf{P}_{\mathbf{F}(\mathbf{a})}^{\bot} \frac{  \partial  \mathbf{F}(\mathbf{a}) }{ \partial \mathbf{a}_i } \mathbf{F}(\mathbf{a})^{+}  \frac{ \partial  \mathbf{F}(\mathbf{a}) } { \partial \mathbf{a}_j } \mathbf{F}(\mathbf{a})^{+} \mathbf{r}(\mathbf{a})   \\
                   = \ & - \mathbf{r}(\mathbf{a})^{T} \frac{  \partial  \mathbf{F}(\mathbf{a}) }{ \partial \mathbf{a}_i } \mathbf{F}(\mathbf{a})^{+}  (\mathbf{F}(\mathbf{a})^{+})^{T}  \frac{ \partial  \mathbf{F}(\mathbf{a})^{T} } { \partial \mathbf{a}_j }  \mathbf{r}(\mathbf{a}) \ ,
\end{align*}
using again the facts that  $ \mathbf{P}_{\mathbf{F}(\mathbf{a})}^{\bot}  \mathbf{r}(\mathbf{a}) = \mathbf{r}(\mathbf{a})$, $ \mathbf{F}(\mathbf{a})^{+} \mathbf{r}(\mathbf{a}) = \mathbf{0}^{k.n}  $ and $\mathbf{P}_{\mathbf{F}(\mathbf{a})}^{\bot} = (\mathbf{P}_{\mathbf{F}(\mathbf{a})}^{\bot})^{T} $.

Collecting all these results together, we deduce that
\begin{align*}
    \mathbf{S}_{ij}
           = \  & \mathbf{r}(\mathbf{a})^{T}   \frac{  \partial  \mathbf{F}(\mathbf{a}) }{ \partial \mathbf{a}_j }  \mathbf{F}(\mathbf{a})^{+}  \frac{  \partial  \mathbf{F}(\mathbf{a}) }{ \partial \mathbf{a}_i } \mathbf{b}   \\
                 & +\:  \mathbf{r}(\mathbf{a})^{T} \frac{  \partial  \mathbf{F}(\mathbf{a}) }{ \partial \mathbf{a}_i } \mathbf{F}(\mathbf{a})^{+} \frac{ \partial  \mathbf{F}(\mathbf{a}) } { \partial \mathbf{a}_j } \mathbf{b}  \\
                 & -\:  2 \mathbf{r}(\mathbf{a})^{T} \frac{  \partial  \mathbf{F}(\mathbf{a}) }{ \partial \mathbf{a}_i } \mathbf{F}(\mathbf{a})^{+} (\mathbf{F}(\mathbf{a})^{+})^{T}  \frac{ \partial  \mathbf{F}(\mathbf{a})^{T} } { \partial \mathbf{a}_j } \mathbf{r}(\mathbf{a}) \ .
\end{align*}
Now from equations~\eqref{eq: deriv_F(a)_delta_b} and~\eqref{eq: deriv_F(a)_delta_T_r(a)} in Subsection~\ref{jacob:box}, $\forall \thickspace \mathbf{a}, \triangle\mathbf{a}  \in \mathbb{R}^{p.k}$, we have
\begin{equation*}
\big(\mathit{D}( \mathbf{F}(\mathbf{a}) )(\triangle\mathbf{a})\big)\mathbf{b} = \mathbf{U} (\mathbf{a}) \triangle\mathbf{a} \quad \text{and} \quad \big( \mathit{D}( \mathbf{F}(\mathbf{a}) )( \triangle\mathbf{a} ) \big)^{T} \mathbf{r}(\mathbf{a}) =  \mathbf{V} (\mathbf{a}) \triangle\mathbf{a} \ ,
\end{equation*}
where $ \mathbf{U} (\mathbf{a}) $ and $ \mathbf{V} (\mathbf{a}) $ are defined in equations~\eqref{eq:U_mat} and~\eqref{eq:V_mat} of Subsection~\ref{jacob:box}. In these conditions, it is evident that
\begin{equation*}
    \frac{  \partial  \mathbf{F}(\mathbf{a}) }{ \partial \mathbf{a}_i } \mathbf{b} = \mathbf{U}(\mathbf{a})_{.i} \quad \text{and} \quad \frac{ \partial  \mathbf{F}(\mathbf{a})^{T} } { \partial \mathbf{a}_i } \mathbf{r}(\mathbf{a}) = \mathbf{V}(\mathbf{a})_{.i} \quad \text{for all } i=1, \cdots, p.k \ .
\end{equation*}
Thus,
\begin{align*}
    \mathbf{S}_{ij}
                & =  \big( \mathbf{V}(\mathbf{a})_{.j} \big)^{T} \mathbf{F}(\mathbf{a})^{+} \mathbf{U}(\mathbf{a})_{.i} + \big( \mathbf{V}(\mathbf{a})_{.i} \big)^{T} \mathbf{F}(\mathbf{a})^{+} \mathbf{U}(\mathbf{a})_{.j} - 2 . \big( \mathbf{V}(\mathbf{a})_{.i} \big)^{T} \mathbf{F}\big(\mathbf{a})^{+} \big(\mathbf{F}(\mathbf{a})^{+} \big)^{T} \mathbf{V}(\mathbf{a})_{.j}  \\
                 & =  \big( \mathbf{L}(\mathbf{a})_{.j} \big)^{T}  \mathbf{U}(\mathbf{a})_{.i} + \big( \mathbf{L}(\mathbf{a})_{.i} \big)^{T}  \mathbf{U}(\mathbf{a})_{.j} - 2 . \big( \mathbf{L}(\mathbf{a})_{.i} \big)^{T}  \mathbf{L}(\mathbf{a})_{.j}  \\
                 & =  \big(\mathbf{U}(\mathbf{a})_{.i}\big)^{T} \mathbf{L}(\mathbf{a})_{.j} + \big( \mathbf{L}(\mathbf{a})_{.i} \big)^{T}  \mathbf{U}(\mathbf{a})_{.j} - 2 . \big( \mathbf{L}(\mathbf{a})_{.i} \big)^{T}  \mathbf{L}(\mathbf{a})_{.j} \ ,
\end{align*}
where $\mathbf{L}(\mathbf{a})$ is  defined in equation~\eqref{eq:L_mat} of Subsection~\ref{jacob:box}. This implies, finally, that
\begin{equation}  \label{eq:S_hess_mat2}
    \mathbf{S} = \mathbf{U}(\mathbf{a})^{T} \mathbf{L}(\mathbf{a}) + \mathbf{L}(\mathbf{a})^{T} \mathbf{U}(\mathbf{a}) - 2 . \mathbf{L}(\mathbf{a})^{T} \mathbf{L}(\mathbf{a}) \ .
\end{equation}
\\
Using the fact that $ \mathit{J}( \mathbf{r}(\mathbf{a}) )^{T} \mathit{J}( \mathbf{r}(\mathbf{a}) ) = \mathbf{M}(\mathbf{a})^{T} \mathbf{M}(\mathbf{a}) + \mathbf{L}(\mathbf{a})^{T} \mathbf{L}(\mathbf{a}) $ established in equation~\eqref{eq:gn_hess_mat} above, we finally obtain the following explicit and compact expression for the Hessian matrix $\nabla^2 \psi (\mathbf{a})$ as the sum of three symmetric matrix terms
\begin{equation} \label{eq:H_hess_mat}
    \mathbf{H} = \mathbf{M}(\mathbf{a})^{T} \mathbf{M}(\mathbf{a}) - \mathbf{L}(\mathbf{a})^{T} \mathbf{L}(\mathbf{a}) + \big( \mathbf{U}(\mathbf{a})^{T} \mathbf{L} (\mathbf{a}) + \mathbf{L}(\mathbf{a})^{T} \mathbf{U}(\mathbf{a}) \big) \ ,
\end{equation}
under the hypothesis that $ \mathbf{F}(.) $ has full column-rank in a neighborhood of $ \mathbf{a} $.

We first note that $\mathbf{H}$ is relatively cheap to evaluate as all the matrix terms involved in the above expression of $\mathbf{H}$ are also needed to compute the Jacobian matrix $\mathit{J}( \mathbf{r}(\mathbf{a}) )$ and are, thus, available in most cases. This small difference between a full Newton and a Gauss-Newton approach in terms of computational load is due to the fact that all the second partial derivatives of $ \mathbf{F}(.) $ are identically zero everywhere and that, consequently, all the mixed partial derivatives matrix terms, which are normally present in the Hessian matrix, vanish here~\cite{B2009}. This small overhead of a full Newton method for solving the WLRA problem with a variable projection or  Grassmann manifold approaches has already been highlighted by Boumal and Absil~\cite{BA2015} when solving a regularized version of the WLRA problem when zero weights are present with a Riemannian method.

Next, we observe that the Gauss-Newton approximation of the Hessian matrix is
\begin{equation*}
 \nabla^2 \psi( \mathbf{a} ) \approx \mathbf{M}(\mathbf{a})^{T} \mathbf{M}(\mathbf{a}) + \mathbf{L}(\mathbf{a})^{T} \mathbf{L}(\mathbf{a}) \  ,
\end{equation*}
while the first two symmetric terms in the exact formulation of $\nabla^2  \psi( \mathbf{a} )$ derived above are 
\begin{equation*}
    \mathbf{M}(\mathbf{a})^{T} \mathbf{M}(\mathbf{a}) - \mathbf{L}(\mathbf{a})^{T} \mathbf{L}(\mathbf{a}) \  .
\end{equation*}
This suggests that  the term $\mathbf{M}(\mathbf{a})^{T} \mathbf{M}(\mathbf{a})$ can be eventually a much better estimate of the full Hessian matrix $\nabla^2 \psi (\mathbf{a})$ than its Gauss-Newton estimate surprisingly. This feature is not specific of the WLRA problem, but rather related to the variable projection framework, and remains true for any separable NLLS problem~\cite{B2009}. This is also verified in the comparison experiments of Hong et al.~\cite{HF2015} in which a Levenberg-Marquardt algorithm using $-\mathbf{M}(\mathbf{a})$ as a simplified Jacobian matrix performs equally or even better than the one using the full Jacobian matrix $-\left( \mathbf{M}(\mathbf{a}) + \mathbf{L}(\mathbf{a}) \right)$. As we will explain in Subsection~\ref{vp_gn_alg:box}, approximating the Jacobian matrix by the term $-\mathbf{M}(\mathbf{a})$ corresponds also to the simplification of the standard variable projection Gauss-Newton and Levenberg-Marquardt algorithms introduced by Kaufman~\cite{K1975} and Ruhe and Wedin~\cite{RW1980}. These results were also given given in Hong and Fitzgibbon~\cite{HF2015}, but their notations are different from those used here.

Based on similar arguments, we can also develop an efficient quasi-Newton method in which we drop the last symmetric term, e.g., $\mathbf{U}(\mathbf{a})^{T} \mathbf{L} (\mathbf{a})+ \mathbf{L}(\mathbf{a})^{T} \mathbf{U}(\mathbf{a})$ in the full Hessian, and use the symmetric matrix
\begin{equation} \label{eq:approx_hess_mat}
\bar{ \mathbf{H} } = \mathbf{M}(\mathbf{a})^{T} \mathbf{M}(\mathbf{a}) - \mathbf{L}(\mathbf{a})^{T} \mathbf{L}(\mathbf{a})  \ ,
\end{equation}
as an approximate Hessian matrix. This approach is new and has never been proposed  in the literature to the best of our knowledge. Obviously, we are not assure that this approximation is always positive semi-definite as in any Newton-like method, but this will be generally the rule as the term $\mathbf{M}(\mathbf{a})$ usually dominates the term $\mathbf{L}(\mathbf{a})$ in the Jacobian matrix as discussed above. However, an iterative algorithm based on this approximate Hessian $\bar{ \mathbf{H} }$ should perform similarly or  better than any Gauss-Newton- or Levenberg-Marquardt-like methods using the full Jacobian $-\left( \mathbf{M}(\mathbf{a}) + \mathbf{L}(\mathbf{a}) \right)$ or its approximation $-\mathbf{M}(\mathbf{a})$ as discussed above.

The approximate Hessian matrix $\bar{ \mathbf{H} }$ also inherits of the singularity of the Jacobian matrix and of its two terms (see Theorems~\ref{theo5.2:box} and~\ref{theo5.3:box}) since $\emph{null}( \mathit{J}( \mathbf{r}(\mathbf{a}) ) ) = \emph{null}( \mathbf{M}(\mathbf{a})  ) \cap \emph{null}( \mathbf{L}(\mathbf{a})  )$ is always included in the null space of $\bar{ \mathbf{H} }$ and this can be considered at first sight as a disadvantage. However, we will demonstrate below that a Newton method using the full Hessian matrix $\mathbf{H}$ has also to overcome similar singularity and ill-conditioning problems in a small neighborhood of a first-order stationary point of $\psi(.)$. Furthermore, we will finally illustrate that, in both cases, we can handle these difficulties for many practical cases with the same techniques developed for the Gauss-Newton or Levenberg-Marquardt algorithms in Subsection~\ref{jacob:box}), e.g., by restricting the  search directions for the Newton correction to $ \emph{ran}( \mathbf{A} )^{\bot}$ at each iteration.
 \\
 \\
The results in Subsection~\ref{jacob:box} show that the problem of minimizing $\psi(.)$ is an ill-posed NLLS problem in the sense that the Jacobian matrix, $ \mathit{J}( \mathbf{r}(\mathbf{a}) ) $, is exactly rank-deficient everywhere in the solution space $\mathbb{R}^{p \times k}_{k}$. We first derive some KKT theorems which give a characterization of a local minimum of $ \psi(.) $ for this very special class of NLLS problems.
 \\
\begin{theo5.9} \label{theo5.9:box}
Let  $ \mathbf{\widehat{a}} = \emph{vec}( \mathbf{\widehat{A}} )^{T}$, with $\mathbf{\widehat{A}} \in \mathbb{R}^{p \times k}_{k}$, be a first-order stationary point of  $ \psi(.) $ and $ \mathbf{F}(.) $ has full column-rank in a neighborhood of $ \mathbf{\widehat{a}} $. Then, the Hessian matrix of $ \psi(.) $ at $ \mathbf{\widehat{a}} $, $ \mathbf{\widehat{H}} = \nabla^2  \psi( \mathbf{\widehat{a}} ) $, is a matrix of rank less than or equal to $ ( p - k ).k $ as for the Jacobian matrix $ \mathit{J}( \mathbf{r}(\mathbf{\widehat{a}} ) )$ (see Corollary~\ref{corol5.3:box}).
\end{theo5.9}
\begin{proof}
First, the hypothesis that $ \mathbf{F}(\mathbf{a}) $ has full column-rank in a neighborhood of $ \mathbf{\widehat{a}} $ implies the existence of $ \mathbf{\widehat{H}} $ as demonstrated above. Furthermore, from equation~\eqref{eq:H_hess_mat}, $ \mathbf{\widehat{H}} $ is given by
\begin{equation*}
    \mathbf{\widehat{H}} = \mathbf{M}(\mathbf{\widehat{a}})^{T} \mathbf{M}(\mathbf{\widehat{a}}) - \mathbf{L}(\mathbf{\widehat{a}})^{T} \mathbf{L}(\mathbf{\widehat{a}} ) + \mathbf{U}(\mathbf{\widehat{a}} )^{T} \mathbf{L} (\mathbf{\widehat{a}}) + \mathbf{L}(\mathbf{\widehat{a}})^{T} \mathbf{U}(\mathbf{\widehat{a}} ) \  .
\end{equation*}
Now, as in Theorem~\ref{theo5.2:box}, we consider the matrix
\begin{equation*}
    \mathbf{\widehat{N}} = \mathbf{K}_{(p,k)} ( \mathbf{I}_{k} \otimes \mathbf{\widehat{A}} ) \  .
\end{equation*}
From Theorem~\ref{theo5.2:box}, we have $ \emph{rank} ( \mathbf{\widehat{N}} ) = k.k $  and the relation
\begin{equation*}
    \emph{ran}( \mathbf{\widehat{N}} ) \subset \emph{null}( \mathbf{M}(\mathbf{\widehat{a}}) ) \cap \emph{null}( \mathbf{L} (\mathbf{\widehat{a}}) ) \ .
\end{equation*}
We also observe that the theorem will be proved if the relation $ \emph{ran}( \mathbf{\widehat{N}} ) \subset \emph{null}( \mathbf{\widehat{H}} ) $ is verified.

To prove this inclusion, we first note that, if $ \mathbf{y} \in \emph{ran}( \mathbf{\widehat{N}} ) $, we have
\begin{equation*}
    \mathbf{M}(\mathbf{\widehat{a}}) \mathbf{y} = \mathbf{L} (\mathbf{\widehat{a}})  \mathbf{y} = \mathbf{0}^{n.p} \ .
\end{equation*}
From this, we deduce that
\begin{align*}
    \mathbf{\widehat{H}} \mathbf{y}
          & =  \mathbf{M}(\mathbf{\widehat{a}})^{T} \mathbf{M}(\mathbf{\widehat{a}}) \mathbf{y}  -  \mathbf{L} (\mathbf{\widehat{a}})  ^{T}  \mathbf{L} (\mathbf{\widehat{a}})  \mathbf{y} + \mathbf{U}(\mathbf{\widehat{a}} )^{T}  \mathbf{L} (\mathbf{\widehat{a}})  \mathbf{y} +  \mathbf{L} (\mathbf{\widehat{a}})^{T} \mathbf{U}(\mathbf{\widehat{a}} ) \mathbf{y} \\
          & =  \mathbf{L} (\mathbf{\widehat{a}})^{T}\mathbf{U}(\mathbf{\widehat{a}} )\mathbf{y} \ .
\end{align*}
Thus, to demonstrate that $ \mathbf{y} \in \emph{null}( \mathbf{\widehat{H}} ) $, it suffices to show that $ \mathbf{L} (\mathbf{\widehat{a}})^{T} \mathbf{U}(\mathbf{\widehat{a}} ) \mathbf{y} = \mathbf{0}^{k.p} $.

If $ \mathbf{y} \in \emph{ran}( \mathbf{\widehat{N}} ) $, then $  \exists \mathbf{Z} \in \mathbb{R}^{k \times k} $ such that
\begin{equation*}
    \mathbf{y} = \mathbf{\widehat{N}} \emph{vec} ( \mathbf{Z} ) = \mathbf{K}_{(p,k)} ( \mathbf{I}_{k} \otimes \mathbf{\widehat{A}} ) \emph{vec} ( \mathbf{Z} ) \ .
\end{equation*}
Consequently, using the definition of $\mathbf{U}( \mathbf{\widehat{a}} )$ (see Subsection~\ref{jacob:box}), we have
\begin{align*}
    \mathbf{U}( \mathbf{\widehat{a}} ) \mathbf{y}
           & =  \emph{diag} \big( \emph{vec}(\sqrt{\mathbf{W}})  \big) (\mathbf{\widehat{B}}^{T} \otimes \mathbf{I}_{p}) \mathbf{K}_{(k,p)} \mathbf{K}_{(p,k)} ( \mathbf{I}_{k} \otimes \mathbf{\widehat{A}} ) \emph{vec} ( \mathbf{Z} )   \\
           & =  \emph{diag} \big( \emph{vec}(\sqrt{\mathbf{W}})  \big) (\mathbf{\widehat{B}}^{T} \otimes \mathbf{I}_{p}) ( \mathbf{I}_{k} \otimes \mathbf{\widehat{A}} ) \emph{vec} ( \mathbf{Z} )   \\
           & =  \emph{diag} \big( \emph{vec}(\sqrt{\mathbf{W}})  \big) (\mathbf{\widehat{B}}^{T} \otimes \mathbf{\widehat{A}} )  \emph{vec} ( \mathbf{Z} )   \\
           & =  \emph{diag} \big( \emph{vec}(\sqrt{\mathbf{W}})  \big)  \emph{vec} ( \mathbf{\widehat{A}} \mathbf{Z} \mathbf{\widehat{B}})   \\
           & =  \emph{diag} \big( \emph{vec}(\sqrt{\mathbf{W}})  \big)  ( \mathbf{I}_{n} \otimes \mathbf{\widehat{A}} ) \emph{vec} ( \mathbf{Z} \mathbf{\widehat{B}} )   \\
           & =  \mathbf{F}( \mathbf{\widehat{a}} ) \emph{vec} ( \mathbf{Z} \mathbf{\widehat{B}} ) \ .
\end{align*}
Now, as in Subsection~\ref{jacob:box}, using the projection operator, $P_{\Omega}(.)$, associated with the $p \times n$ weight matrix $\mathbf{W}$, we have
\begin{equation*}
    \big \lbrack P_{\Omega}(\mathbf{X} - \mathbf{\widehat{A}} \mathbf{\widehat{B}} ) \big \rbrack_{ij} =
    \begin{cases}
        \displaystyle{ \mathbf{X}_{ij} - \sum_{l=1}^{k} \mathbf{\widehat{A}}_{il} \mathbf{\widehat{B}}_{lj} } & \text{if } \mathbf{W}_{ij} \ne 0,\\
         0                                                                                                                                                       & \text{if } \mathbf{W}_{ij}  =   0 .
    \end{cases} \ ,
\end{equation*}
where $\mathbf{\widehat{b}} = \mathbf{F}(\mathbf{\widehat{a}})^{+}\mathbf{x}$  and $\mathbf{\widehat{B}} = \emph{mat}( \mathbf{\widehat{b}} )$, and the variable projection residual vector of $\mathbf{X}$  at $\mathbf{\widehat{A}}$ can be written as
\begin{equation*}
    \mathbf{r}(\mathbf{\widehat{a}}) = \emph{vec} \big( \sqrt{\mathbf{W}} \odot P_{\Omega}(\mathbf{X} -\mathbf{\widehat{A}} \mathbf{\widehat{B}}) \big) \ .
\end{equation*}
Thus, using the definition of $\mathbf{L} (\mathbf{\widehat{a}})$ (see equation~\eqref{eq:L_mat} in Section~\ref{jacob:box}) and the fact that $ \mathbf{F}( \mathbf{\widehat{a}} )^{+}  \mathbf{F}( \mathbf{\widehat{a}}  ) = \mathbf{I}_{n.k} $, which results from the hypothesis that $ \mathbf{F}( \mathbf{\widehat{a}} ) $ has full column-rank, we have
\begin{align*}
   \mathbf{L} (\mathbf{\widehat{a}})^{T} \mathbf{U}(\mathbf{\widehat{a}} ) \mathbf{y}
          & =  \big( ( \mathbf{W} \odot P_{\Omega}(\mathbf{X} - \mathbf{\widehat{A}} \mathbf{\widehat{B}} )  ) \otimes \mathbf{I}_{k} \big)  \mathbf{F}( \mathbf{\widehat{a}} )^{+}  \mathbf{F}( \mathbf{\widehat{a}}  ) \emph{vec} ( \mathbf{Z} \mathbf{\widehat{B}} )   \\
          & =  \big( ( \mathbf{W} \odot P_{\Omega}(\mathbf{X} - \mathbf{\widehat{A}} \mathbf{\widehat{B}} )  ) \otimes \mathbf{I}_{k} \big)   \emph{vec} ( \mathbf{Z} \mathbf{\widehat{B}} ) \ ,
\end{align*}
Finally, this leads to the equalities
\begin{align*}
    \mathbf{\widehat{L}}^{T} \mathbf{U}(\mathbf{\widehat{a}} ) \mathbf{y}
          & =  \emph{vec} \big( \mathbf{Z} \mathbf{\widehat{B}} ( \mathbf{W} \odot P_{\Omega}(\mathbf{X} - \mathbf{\widehat{A}} \mathbf{\widehat{B}} ) )^{T} \big)   \\
          & =  ( \mathbf{I}_{p}  \otimes \mathbf{Z}  )    \emph{vec} \big(  \mathbf{\widehat{B}} ( \mathbf{W} \odot P_{\Omega}(\mathbf{X} - \mathbf{\widehat{A}} \mathbf{\widehat{B}} ) )^{T} \big)   \\          
          & =  ( \mathbf{I}_{p}  \otimes \mathbf{Z}  )  ( \mathbf{I}_{p}  \otimes \mathbf{\widehat{B}}  )  \emph{vec} \big(  ( \mathbf{W} \odot P_{\Omega}(\mathbf{X} - \mathbf{\widehat{A}} \mathbf{\widehat{B}} ) )^{T} \big)   \\
          & =  ( \mathbf{I}_{p}  \otimes \mathbf{Z}  )  ( \mathbf{I}_{p}  \otimes \mathbf{\widehat{B}}  )  \emph{diag} \big( \emph{vec}(\sqrt{\mathbf{W}}^{T}) \big) \emph{vec} \big(  ( \sqrt{\mathbf{W}} \odot P_{\Omega}(\mathbf{X} - \mathbf{\widehat{A}} \mathbf{\widehat{B}} ) )^{T} \big)   \\
          & =  ( \mathbf{I}_{p}  \otimes \mathbf{Z}  )  ( \mathbf{I}_{p}  \otimes \mathbf{\widehat{B}}  )  \emph{diag} \big( \emph{vec}(\sqrt{\mathbf{W}}^{T}) \big) \emph{vec} \big(  ( \sqrt{\mathbf{W}} \odot P_{\Omega}(\mathbf{X} - \mathbf{\widehat{A}} \mathbf{\widehat{B}} ) )^{T} \big)   \\         
          & =  ( \mathbf{I}_{p}  \otimes \mathbf{Z}  ) \Big(  \emph{diag} \big( \emph{vec}(\sqrt{\mathbf{W}}^{T}) \big) ( \mathbf{I}_{p}  \otimes \mathbf{\widehat{B}}  )^{T}  \Big)^{T} \emph{vec} \big(  ( \sqrt{\mathbf{W}} \odot P_{\Omega}(\mathbf{X} - \mathbf{\widehat{A}} \mathbf{\widehat{B}} ) )^{T} \big)   \\         
          & =  ( \mathbf{I}_{p}  \otimes \mathbf{Z}  ) \Big(  \emph{diag} \big( \emph{vec}(\sqrt{\mathbf{W}}^{T}) \big) ( \mathbf{I}_{p}  \otimes \mathbf{\widehat{B}}^{T}  )  \Big)^{T} \emph{vec} \big(  ( \sqrt{\mathbf{W}} \odot P_{\Omega}(\mathbf{X} - \mathbf{\widehat{A}} \mathbf{\widehat{B}} ) )^{T} \big)   \\         
          & =  ( \mathbf{I}_{p}  \otimes \mathbf{Z}  ) \mathbf{G}( \mathbf{\widehat{b}} )^{T}  \emph{vec} \big(  ( \sqrt{\mathbf{W}} \odot P_{\Omega}(\mathbf{X} - \mathbf{\widehat{A}} \mathbf{\widehat{B}} ) )^{T} \big)   \\         
          & =  ( \mathbf{I}_{p}  \otimes \mathbf{Z}  ) \mathbf{G}( \mathbf{\widehat{b}} )^{T}  \mathbf{K}_{(p,n)} \emph{vec} ( \sqrt{\mathbf{W}} \odot P_{\Omega}(\mathbf{X} - \mathbf{\widehat{A}} \mathbf{\widehat{B}} ) )   \\         
          & =  ( \mathbf{I}_{p}  \otimes \mathbf{Z}  ) \mathbf{G}( \mathbf{\widehat{b}} )^{T}  \mathbf{K}_{(p,n)}  \mathbf{r} ( \mathbf{\widehat{a}} ) \\         
          & =  - ( \mathbf{I}_{p}  \otimes \mathbf{Z}  )  \nabla \psi( \mathbf{\widehat{a}} )  \\
          & =  \mathbf{0}^{k.p} \ ,
\end{align*}
the two last equalities resulting from Theorem~\ref{theo5.7:box} and the hypothesis $ \nabla \psi( \mathbf{\widehat{a}} )  = \mathbf{0}^{k.p} $. Thus, $ \mathbf{\widehat{H}} \mathbf{y} = \mathbf{0}^{k.p} $ and $ \emph{ran}( \mathbf{\widehat{N}} ) \subset \emph{null}( \mathbf{\widehat{H}} ) $. This implies, finally, that $ k.k \leqslant  \emph{dim}( \emph{null}( \mathbf{\widehat{H}} ) ) $ and $ \emph{rank}( \mathbf{\widehat{H}} ) \leqslant  ( p - k ).k $ as claimed in the theorem.

\end{proof}
Furthermore, if $ \emph{rank}( \mathit{J}( \mathbf{r}( \mathbf{\widehat{a}} ) ) = ( p - k ).k $, for example if the hypotheses of Theorem~\ref{theo5.3:box} are satisfied, we have the following corollary:
\\
\begin{corol5.9} \label{corol5.9:box}
Let  $ \mathbf{\widehat{a}} = \emph{vec}( \mathbf{\widehat{A}} )^{T}$, with $\mathbf{\widehat{A}} \in \mathbb{R}^{p \times k}_{k}$,  be a first-order stationary point of  $ \psi (.)$ and suppose that $\mathbf{F}(\mathbf{a})$ has full column-rank in a neighborhood of $ \mathbf{\widehat{a}} $ and the rank of the Jacobian matrix $ \mathit{J}( \mathbf{r}(\mathbf{\widehat{a}}) ) $ is equal to $ r = ( p - k ).k $. Then, the Hessian matrix of $ \psi (.)$ at $ \mathbf{\widehat{a}} $, $ \mathbf{\widehat{H}} = \nabla^2 \psi( \mathbf{\widehat{a}} ) $, is a matrix of rank less than or equal to $ r $ with its null space containing the null space of $ \mathit{J}( \mathbf{r}(\widehat{a} ) )$.
\end{corol5.9}
\begin{proof}
Using the same notations as in Theorem~\ref{theo5.9:box}, the hypotheses imply the relation
\begin{equation*}
    \emph{dim} \Big( \emph{null}  \big( \mathit{J}( \mathbf{r}(\mathbf{\widehat{a}}) )  \big) \Big) =  k.k = \emph{rank}( \mathbf{\widehat{N}} ) \ .
\end{equation*}
From the relation $ \emph{ran}( \mathbf{\widehat{N}} ) \subset \emph{null}( \mathbf{M} (\mathbf{\widehat{a}}) ) \cap \emph{null}( \mathbf{L} (\mathbf{\widehat{a}}) ) = \emph{null}( \mathit{J} \big( \mathbf{r}(\mathbf{\widehat{a}}) ) \big) $, we then deduce that $ \emph{ran}( \mathbf{\widehat{N}} ) = \emph{null}( \mathit{J} \big( \mathbf{r}(\mathbf{\widehat{a}}) ) \big) $ and all the results are a direct consequence of Theorem~\ref{theo5.9:box}.

\end{proof}
\begin{remark5.3} \label{remark5.3:box}
Theorem~\ref{theo5.9:box} and Corollary~\ref{corol5.9:box} are similar to Theorem 2.2 of Eriksson et al.~\cite{EWGS2005} and its consequences, which deal with NLLS problems that are uniformly rank-deficient, i.e., with a Jacobian matrix that have the same deficient rank in the neighborhood of a solution. However, it is important to highlight that the strong assumption of uniform rank deficiency of the Jacobian matrix in the neighborhood of a solution has not been used here to derive Theorem~\ref{theo5.9:box}, while this hypothesis is central in the results of Eriksson et al.~\cite{EWGS2005} through the use of the Constant-Rank Theorem~\cite{C2017}. $\blacksquare$
\\
\end{remark5.3}
\begin{remark5.4} \label{remark5.4:box}
Theorem~\ref{theo5.9:box} and Corollary~\ref{corol5.9:box} also explain why full Newton variable projection algorithms (without any specific adaptation or regularization) perform poorly in the comparison experiments of Okatani et al.~\cite{OYD2011} and Hong et al.~\cite{HF2015}. As soon as we are in a small neighborhood of a first-order stationary point $\mathbf{\widehat{a} }$ of  $\psi(.)$, the Hessian matrix $ \nabla^2 \psi( \mathbf{a} ) $ may become nearly singular and ill-conditioned according to Theorem~\ref{theo5.9:box} and Corollary~\ref{corol5.9:box} leading to an erratic behaviour and a dramatic lost of accuracy in the final steps of the Newton method. $\blacksquare$

\end{remark5.4}

Importantly, Theorem~\ref{theo5.9:box} and Corollary~\ref{corol5.9:box} can also be used to illustrate that the sufficient condition for the existence of a strict local minima of $\psi(.)$ (see Subsection~\ref{calculus:box} for details) are never meet, but non strict local minima may still exist though. To see this, suppose that $\psi(.)$ is twice continuously differentiable at $\mathbf{\widehat{a}}$ and that $\mathbf{\widehat{a}}$ is a first-order stationary point of $\psi(.)$, then  $ \nabla \psi( \mathbf{\widehat{a}} ) = \mathbf{0}^{k.p}$ and  the second-order Taylor expansion of $\psi(.)$ at $\mathbf{\widehat{a}}$ becomes
\begin{equation*}
\psi( \mathbf{\widehat{a}} + d\mathbf{a} ) = \psi( \mathbf{\widehat{a}} ) + d\mathbf{a}^{T} \nabla^2 \psi( \mathbf{\widehat{a}} ) d\mathbf{a} +  \mathcal{O}( \Vert d\mathbf{a} \Vert^{3}_{2} ) \ .
\end{equation*}
However, for $d\mathbf{a} \in \emph{null}( \mathbf{\widehat{N}} )$, where $\mathbf{\widehat{N}}= \mathbf{K}_{(p,k)} ( \mathbf{I}_{k} \otimes \mathbf{\widehat{A}} )$, and $d\mathbf{a}$ sufficiently small, this Taylor expansion reduces to
\begin{equation*}
\psi( \mathbf{\widehat{a}} + d\mathbf{a} ) = \psi( \mathbf{\widehat{a}} ) +  \mathcal{O}( \Vert d\mathbf{a} \Vert^{3}_{2} ) \ ,
\end{equation*}
according to Theorem~\ref{theo5.9:box}. This last equation and the fact that the matrix $\nabla^2 \psi( \mathbf{\widehat{a}} )$ is never positive definite, are consistent with the property that, if the~\eqref{eq:VP1} problem admits a solution $\mathbf{\widehat{a}}$, then there exist infinitely many solutions near $\mathbf{\widehat{a}}$ (see Remark~\ref{remark3.5:box} for details). Thus, having $\mathit{J}( \mathbf{r}(\mathbf{a}) )$ rank deficient everywhere makes the minimization of $\psi(.)$ an ill-posed problem since if this problem has a solution $\mathbf{\widehat{a}}$ than this solution is never unique and, in addition, there are an infinite set of solutions in any neighborhood of  $\mathbf{\widehat{a}}$.

To conclude this section, we now give sufficient second-order KKT conditions for $\psi(.)$ to have a (non strict) minimum under the exact and constant rank-deficient condition of $\mathit{J}( \mathbf{r}(\mathbf{a}) )$ in a neighborhood of a first-order stationary point in the following theorem, which is  a reformulation and an extension in our context of results demonstrated in Eriksson et al.~\cite{EWGS2005}.
\\
\begin{theo5.10} \label{theo5.10:box}
Let  $\mathbf{\widehat{a}} \in \mathbb{R}^{k.p}$, with $ \widehat{\mathbf{A}} = ( \emph{mat}_{k \times p}( \widehat{\mathbf{a}} ) )^{T}  \in  \mathbb{R}^{p \times k}_{k}$, be a first-order stationary point of  $\psi (.)$, and suppose that $\mathbf{F}(\mathbf{a})$ has full column-rank and  the rank of the Jacobian matrix $\mathit{J}( \mathbf{r}(\mathbf{a}) )$ is equal to $r = (p - k).k$ in a neighborhood of $ \mathbf{\widehat{a}} $. Let further $\mathbf{Q}  \in \mathbb{R}^{k.p \times r}_{r}$ such that the columns of $\mathbf{Q}$ form a basis for $\emph{ran}( \mathit{J}( \mathbf{r}(\mathbf{\widehat{a}}) )^{T} ) = \emph{null}( \mathit{J}( \mathbf{r}(\mathbf{\widehat{a}}) ) )^{\bot}$.
Then, define the $r \times r$ symmetric matrix
\begin{align}  \label{eq:T_mat}
    \mathbf{\widehat{T}} & =  \mathbf{Q}^{T} \nabla^2 \psi( \mathbf{\widehat{a}} ) \mathbf{Q}  \nonumber \\
                                      & =  \mathbf{Q}^{T}  \mathbf{\widehat{H}} \mathbf{Q}  \nonumber \\                                      
                                      & =  \mathbf{Q}^{T} \left( \mathbf{M}(\mathbf{\widehat{a}})^{T} \mathbf{M}(\mathbf{\widehat{a}})  -  \mathbf{L} (\mathbf{\widehat{a}})  ^{T}  \mathbf{L} (\mathbf{\widehat{a}})  +  \mathbf{U}(\mathbf{\widehat{a}} )^{T}  \mathbf{L} (\mathbf{\widehat{a}}) +  \mathbf{L} (\mathbf{\widehat{a}})^{T}  \mathbf{U}(\mathbf{\widehat{a}} ) \right)  \mathbf{Q}  \ .
\end{align}
In these conditions, if $\psi (.)$ has a local minimum at $\mathbf{\widehat{a}}$ then  $\mathbf{\widehat{T}}$ is positive  semi-definite. Reciprocally, if $\mathbf{\widehat{T}}$ is positive definite then $\psi (.)$  has a (non strict) local minimizer at $\mathbf{\widehat{a}}$ and we have the equalities
\begin{equation*}
 \emph{null} \big(   \nabla^2 \psi( \mathbf{\widehat{a}} )  \big) = \emph{null} \big( \mathit{J}( \mathbf{r}(\mathbf{\widehat{a}}) )  \big) \quad \text{and} \quad  \emph{rank} \big(   \nabla^2 \psi( \mathbf{\widehat{a}} )  \big) = r  \ .
\end{equation*}

\end{theo5.10}
\begin{proof}
We first verify that the condition is necessary if $\mathbf{\widehat{a}} \in \mathbb{R}^{k.p}$ is a local minimizer of  $\psi (.)$. To this end, let $\mathbf{P}  \in \mathbb{R}^{k.p \times k.k}_{k.k}$ such that the columns of $\mathbf{P}$ form a basis for $\emph{null}( \mathit{J}( \mathbf{r}(\mathbf{\widehat{a}}) ) )$. We first note that orthonormal matrices  $\mathbf{P}$ and  $\mathbf{Q}$ can be easily constructed from the results of Corollary~\ref{corol5.6:box} and that, with this choice, the columns of the partitioned matrix
\begin{equation*}
\begin{bmatrix} \mathbf{P} &  \mathbf{Q}  \end{bmatrix} =  \begin{bmatrix} \mathbf{\bar{O}}  & \mathbf{\bar{O}}^{\bot} \end{bmatrix}
\end{equation*}
is an orthonormal basis of $\mathbb{R}^{k.p}$ since  $\emph{null}( \mathit{J}( \mathbf{r}(\mathbf{\widehat{a}}) ) )$ and $\emph{null}( \mathit{J}( \mathbf{r}(\mathbf{\widehat{a}}) ) )^{\bot}$ are orthogonal subspaces of $\mathbb{R}^{k.p}$ and
\begin{equation*}
   \mathbb{R}^{k.p}  = \emph{null}  \Big( \mathit{J} \big( \mathbf{r}(\mathbf{\widehat{a}}) \big) \Big)  \oplus \emph{null} \Big( \mathit{J} \big( \mathbf{r}(\mathbf{\widehat{a}}) \big) \Big)^{\bot}  \ ,
\end{equation*}
where $\oplus$ stands for the direct sum. Thus, $\forall \mathbf{y} \in \mathbb{R}^{k.p}$, $\mathbf{y}$ can be decomposed as
\begin{equation*}
 \mathbf{y} =  \begin{bmatrix} \mathbf{P} &  \mathbf{Q}  \end{bmatrix}  \begin{bmatrix} \mathbf{y}_{P} \\   \mathbf{y}_{Q}  \end{bmatrix} = \mathbf{P} \mathbf{y}_{P}  + \mathbf{Q} \mathbf{y}_{Q}  \text{ with } \mathbf{y}_{P}  \in \mathbb{R}^{k.k} \text{ and } \mathbf{y}_{Q}  \in \mathbb{R}^{(p-k).k}  \ .
\end{equation*}
Now, as $\mathbf{\widehat{a}}$ is a local minimizer of $\psi (.)$, $\nabla^2 \psi( \mathbf{\widehat{a}} )$ is a positive semi-definite matrix (see Subsection\ref{calculus:box}), which implies that, $\forall \mathbf{y} \in \mathbb{R}^{k.p}$, we have
\begin{equation*}
\mathbf{y}^{T} \nabla^2 \psi( \mathbf{\widehat{a}} ) \mathbf{y} \ge 0 \  .
\end{equation*}
On the other hand, since the null space of $ \mathit{J}( \mathbf{r}(\widehat{a} ) )$ is included in the null space of $\nabla^2 \psi( \mathbf{\widehat{a}} )$ according to Corollary~\ref{corol5.9:box} and taking into account the symmetry of $\nabla^2 \psi( \mathbf{\widehat{a}} )$, we have
\begin{align*}
\mathbf{y}^{T} \nabla^2 \psi( \mathbf{\widehat{a}} ) \mathbf{y} & =  \left(  \mathbf{P} \mathbf{y}_{P} + \mathbf{Q} \mathbf{y}_{Q} \right)^{T} \nabla^2 \psi( \mathbf{\widehat{a}} ) \left(  \mathbf{P} \mathbf{y}_{P} + \mathbf{Q} \mathbf{y}_{Q} \right) \\
& =   \left(  \mathbf{Q} \mathbf{y}_{Q}  \right)^{T} \nabla^2 \psi( \mathbf{\widehat{a}} ) \left(  \mathbf{Q} \mathbf{y}_{Q} \right) \\
& =  \mathbf{y}_{Q}^{T} \mathbf{\widehat{T}} \mathbf{y}_{Q} \ge 0 \  ,
\end{align*}
which demonstrates that $\mathbf{\widehat{T}}$ is a positive semi-definite matrix as claimed in the theorem.

We now give a sketch of the proof for the sufficient condition and we refer the interested reader to Eriksson et al.~\cite{E1996}\cite{EW1996}\cite{EWGS2005} for more details. First, the hypothesis that $\widehat{\mathbf{A}}$ is of full column rank implies that it exists an open neighborhood of $\mathbf{\widehat{a}}$, say $\Upsilon$, in which  all its elements have also full-column rank, as $\mathbb{R}^{p \times k}_{k}$ is open in  $\mathbb{R}^{p \times k}$ according to Theorem~\ref{theo2.3:box}. By eventually restricting this open neighborhood $\Upsilon$, the hypotheses that $\mathbf{F}(\mathbf{a})$ has full column-rank and that the rank of the Jacobian matrix $\mathit{J}( \mathbf{r}(\mathbf{a}) )$ is equal to $r = (p - k).k$ in a neighborhood of $ \mathbf{\widehat{a}} $ imply that the residual function $\mathbf{r}(.)$ from $\mathbb{R}^{k.p}$ to $\mathbb{R}^{p.n}$ is infinitely differentiable in an open neighborhood $\Upsilon$ of $ \mathbf{\widehat{a}} $. This allows us to apply the Constant-Rank Theorem (see Theorem 2.1 in~\cite{EWGS2005} or~\cite{C2017}) 
to demonstrate that there exist two functions, e.g.,
\begin{equation*}
\mathbf{z} : \mathbb{R}^{k.p} \longrightarrow \mathbb{R}^{r}  \quad \text{ and } \quad \mathbf{h} : \mathbb{R}^{r} \longrightarrow \mathbb{R}^{n.p}   \ ,
\end{equation*}
which are  twice continuously differentiable functions, respectively, over an open neighborhood $\Upsilon$ of $\mathbf{\widehat{a}}$ and over  $\mathbb{R}^{r}$, such that $\mathbf{r}(\mathbf{a}) = \mathbf{h}(\mathbf{z}(\mathbf{a}) )$.

Using the chain rule for computing derivatives~\cite{C2017}, we get, for all $\mathbf{a} \in  \Upsilon$, the relation
\begin{equation}  \label{eq:jacob_mat_hz}
\mathit{J} \big( \mathbf{r}(\mathbf{a}) \big) = \mathit{J} \big( \mathbf{h} ( \mathbf{z}(\mathbf{a}) ) \big) \mathit{J} \big( \mathbf{z}(\mathbf{a}) \big)  \ .
\end{equation}
In these conditions, using the hypothesis that the rank of the Jacobian matrix $\mathit{J}( \mathbf{r}(\mathbf{a}) )$ is equal to $r = (p - k).k, \forall \mathbf{a} \in \Upsilon$, and equation~\eqref{eq:rank2} in Subsection~\ref{lin_alg:box}, we deduce immediately that the Jacobian matrices, $\mathit{J}( \mathbf{z}(\mathbf{a}) )$ and $\mathit{J}( \mathbf{h}( \mathbf{z}(\mathbf{a}) ) )$,  have  also a constant (full) rank equals to $r = (p - k).k$  for all $\mathbf{a} \in \Upsilon$. Furthermore,  since both $\mathit{J}( \mathbf{h}( \mathbf{z}(\mathbf{a}) ) )$  and $\mathit{J}( \mathbf{z}(\mathbf{a}) )$ have full rank $r$, which is also equal to the rank of $\mathit{J}( \mathbf{r}(\mathbf{a}) )$, we have the equalities
\begin{equation*}
\emph{null} \big( \mathit{J} ( \mathbf{r}(\mathbf{a}) ) \big) = \emph{null} \big( \mathit{J}( \mathbf{z}(\mathbf{a}) ) \big) \quad \text{and} \quad\emph{ran} \big( \mathit{J} ( \mathbf{r}(\mathbf{a}) ) \big) = \emph{ran} \big( \mathit{J} (\mathbf{h}( \mathbf{z}(\mathbf{a}) )) \big)  \ ,
\end{equation*}
and also the "transposed" versions of these equalities
\begin{equation*}
\emph{null} \big( \mathit{J} ( \mathbf{r}(\mathbf{a}) )^{T} \big) = \emph{null} \big( \mathit{J}( \mathbf{h}( \mathbf{z}(\mathbf{a}) ) )^{T} \big) \quad \text{and} \quad \emph{ran} \big( \mathit{J} ( \mathbf{r}(\mathbf{a}) )^{T} \big) = \emph{ran} \big( \mathit{J}( \mathbf{z} (\mathbf{a}) )^{T} \big) \ ,
\end{equation*}
for all $\mathbf{a} \in  \Upsilon$.

Now, the equality $\emph{null} \left( \mathit{J} ( \mathbf{r}( \mathbf{\widehat{a}} ) )^{T} \right) = \emph{null} \left( \mathit{J}( \mathbf{h}( \mathbf{z}( \mathbf{\widehat{a}} ) ) )^{T} \right) $ and the hypothesis that $\mathbf{\widehat{a}}$ is a first-order stationary point of $\psi (.)$, i.e.,
\begin{equation*}
\nabla \psi( \mathbf{\widehat{a}} ) =  \mathit{J} \big( \mathbf{r}( \mathbf{\widehat{a}} ) \big)^{T} \mathbf{r}( \mathbf{\widehat{a}} ) = \mathbf{0}^{p.k} \ ,
\end{equation*}
imply that
\begin{equation*}
\mathit{J} \big( \mathbf{h}( \mathbf{z}( \mathbf{\widehat{a}} ) ) \big)^{T} \mathbf{h}(\mathbf{z}( \mathbf{\widehat{a}}) )  = \mathit{J} \big( \mathbf{h}( \mathbf{z}(  \mathbf{\widehat{a}} ) ) \big)^{T} \mathbf{r}( \mathbf{\widehat{a}} ) =  \mathbf{0}^{p.k} \ .
\end{equation*}
This demonstrates that  $\mathbf{z}(  \mathbf{\widehat{a}} )$ is a first-order stationary point of the  twice continuously differentiable real function $\phi(.)$ defined by
\begin{equation*}
\phi  : \mathbb{R}^{r} \longrightarrow \mathbb{R}  : \mathbf{o} \mapsto  \frac{1}{2} \Vert  \mathbf{h}( \mathbf{o}  )  \Vert^{2}_{2} =  \frac{1}{2}   \mathbf{h}( \mathbf{o}  )^{T}  \mathbf{h}( \mathbf{o}  ) \ .
\end{equation*}
In these conditions, a sufficient condition for the first-order stationary point $\mathbf{z}(  \mathbf{\widehat{a}} )$ to be a strict local minimizer of $\phi (.)$ is that the Hessian matrix $ \nabla^2 \phi( \mathbf{z}(  \mathbf{\widehat{a}} ) ) $ is positive definite (see Subsection~\ref{calculus:box}) and this will also imply that $\mathbf{\widehat{a}}$ is a non strict local minimizer of $\psi (.)$ since
\begin{equation*}
\psi ( \mathbf{a} ) = \frac{1}{2} \mathbf{r}( \mathbf{a}  )^{T}  \mathbf{r}( \mathbf{a}  ) = \frac{1}{2} \mathbf{h} \big(\mathbf{z}(\mathbf{a}) \big)^{T}  \mathbf{h} \big(\mathbf{z}(\mathbf{a}) \big) = \phi \big ( \mathbf{z}(\mathbf{a}) \big)  \ ,
\end{equation*}
for all $\mathbf{a} \in  \Upsilon$. The Hessian matrix $ \nabla^2 \phi( \mathbf{z}(  \mathbf{\widehat{a}} ) ) $ is positive definite if and only if
\begin{equation*}
 \mathbf{o}^{T} \left( \nabla^2 \phi( \mathbf{z}(  \mathbf{\widehat{a}} ) ) \right) \mathbf{o}  > 0  \text{ for all } \mathbf{o} \in \mathbb{R}^{r}  \text{ with } \mathbf{o} \ne  \mathbf{0}^{r}  \ .
\end{equation*}
Now, the equality $\emph{null} \left( \mathit{J} ( \mathbf{r}(  \mathbf{\widehat{a}} ) ) \right) = \emph{null} \left( \mathit{J}( \mathbf{z}(  \mathbf{\widehat{a}} ) ) \right)$ demonstrated above implies that
\begin{equation*}
\emph{null} \big( \mathit{J} ( \mathbf{r}(  \mathbf{\widehat{a}} ) ) \big)^{\bot} = \emph{null} \big( \mathit{J}( \mathbf{z}(  \mathbf{\widehat{a}} ) ) \big)^{\bot} \ ,
\end{equation*}
and this shows that the $r \times r$ matrix $\mathit{J}( \mathbf{z}(  \mathbf{\widehat{a}} ) )\mathbf{Q}$ is of full rank $r$ as the $r$ columns of $\mathbf{Q}$ form a basis of $\emph{null} \left( \mathit{J}( \mathbf{z}(  \mathbf{\widehat{a}}) ) \right)^{\bot}$. In these conditions, the columns of $\mathit{J}( \mathbf{z}( \mathbf{\widehat{a}} ) )\mathbf{Q}$ form a basis of $\mathbb{R}^{r}$ and any  $\mathbf{o} \in \mathbb{R}^{r}$ can be written as $\mathbf{o}= \mathit{J}( \mathbf{z}( \mathbf{\widehat{a}} ) )\mathbf{Q} \mathbf{w}$ for some $\mathbf{w}  \in \mathbb{R}^{r}$ and the proposition that $ \nabla^2 \phi( \mathbf{z}(  \mathbf{\widehat{a}} ) ) $ is positive definite, is equivalent to
\begin{equation*}
 \mathbf{w}^{T} \left(  \mathbf{Q}^{T} \mathit{J}\big( \mathbf{z}( \mathbf{\widehat{a}} ) \big)^{T} \Big ( \nabla^2 \phi \big( \mathbf{z}(  \mathbf{\widehat{a}} ) \big) \Big ) \mathit{J}\big( \mathbf{z}( \mathbf{\widehat{a}} ) \big) \mathbf{Q} \right) \mathbf{w}  > 0  \text{ for all } \mathbf{w} \in \mathbb{R}^{r}  \text{ with } \mathbf{w} \ne  \mathbf{0}^{r} \ ,
\end{equation*}
e.g., that the $r \times r$ symmetric matrix $\left(  \mathbf{Q}^{T} \mathit{J}\big( \mathbf{z}( \mathbf{\widehat{a}} ) \big)^{T} \Big ( \nabla^2 \phi \big( \mathbf{z}(  \mathbf{\widehat{a}} ) \big) \Big ) \mathit{J}\big( \mathbf{z}( \mathbf{\widehat{a}} ) \big) \mathbf{Q} \right)$ is definite positive. We will now show that this symmetric matrix is nothing else than the  $r \times r$ symmetric matrix $\mathbf{\widehat{T}}$, defined in equation~\eqref{eq:T_mat}, and which is assumed to be positive definite by hypothesis.

Using equations~\ref{eq:hessian_nlls}  and~\ref{eq:jacob_mat_hz} above, we have
\begin{align*}
&   \mathit{J} \big( \mathbf{z}( \mathbf{\widehat{a}} ) \big)^{T} \Big (  \nabla^2 \phi( \mathbf{z}(  \mathbf{\widehat{a}} ) ) \Big )  \mathit{J} \big( \mathbf{z}( \mathbf{\widehat{a}} ) \big)  \\
& =  \mathit{J} \big( \mathbf{z}( \mathbf{\widehat{a}} ) \big)^{T} \left(  \mathit{J}( \mathbf{h}( \mathbf{z}( \mathbf{\widehat{a}} ) ) )^{T} \mathit{J}( \mathbf{h}( \mathbf{z}( \mathbf{\widehat{a}} ) ) ) +  \sum_{l=1}^{n.p} \mathbf{h}_l ( \mathbf{z}( \mathbf{\widehat{a}} ) ) \nabla^2 \mathbf{h}_l ( \mathbf{z}( \mathbf{\widehat{a}} ) )  \right) \mathit{J} \big( \mathbf{z}( \mathbf{\widehat{a}} ) \big) \\
& =    \mathit{J} \big( \mathbf{r}( \mathbf{\widehat{a}} ) \big)^{T}  \mathit{J} \big( \mathbf{r}( \mathbf{\widehat{a}} ) \big) +   \mathit{J} \big( \mathbf{z}( \mathbf{\widehat{a}} ) \big)^{T} \left( \sum_{l=1}^{n.p} \mathbf{r}_l (  \mathbf{\widehat{a}} ) \nabla^2 \mathbf{h}_l ( \mathbf{z}( \mathbf{\widehat{a}} ) )  \right) \mathit{J} \big( \mathbf{z}( \mathbf{\widehat{a}} ) \big)  \ .
\end{align*}
Next, using again the chain rule for computing the second derivatives of $\mathbf{r}_l(.)$ and the fact that
\begin{equation*}
\mathit{J} \big( \mathbf{h}( \mathbf{z}(  \mathbf{\widehat{a}} ) ) \big)^{T} \mathbf{r}( \mathbf{\widehat{a}} ) =  \mathbf{0}^{p.k}  \ ,
\end{equation*}
demonstrated above, we obtain the following expression for the second  term of the Hessian matrix $ \mathbf{\widehat{H}} = \nabla^2 \psi( \mathbf{\widehat{a}} ) $, defined in equations~\ref{eq:S_hess_mat} and~\ref{eq:S_hess_mat2}, in terms of the derivatives of $\mathbf{z}(.)$ and $\mathbf{h}(.)$:
\begin{equation*}
    \mathbf{\widehat{S}} =  \sum_{l=1}^{n.p} \mathbf{r}_l (\mathbf{\widehat{a}} ) \nabla^2 \mathbf{r}_l ( \mathbf{\widehat{a}} )
                                      =  \mathit{J} \big( \mathbf{z}( \mathbf{\widehat{a}} ) \big)^{T} \left( \sum_{l=1}^{n.p} \mathbf{r}_l (  \mathbf{\widehat{a}} ) \nabla^2 \mathbf{h}_l ( \mathbf{z}( \mathbf{\widehat{a}} ) )  \right) \mathit{J} \big( \mathbf{z}( \mathbf{\widehat{a}} ) \big) \ .
\end{equation*}
This implies the equality
\begin{equation*}
\mathit{J}  \big( \mathbf{z}( \mathbf{\widehat{a}} )  \big)^{T} \Big( \nabla^2 \phi  \big( \mathbf{z}(  \mathbf{\widehat{a}} )  \big) \Big) \mathit{J}  \big( \mathbf{z}( \mathbf{\widehat{a}} )  \big) = \mathit{J}  \big( \mathbf{r}( \mathbf{\widehat{a}} )  \big)^{T}  \mathit{J}  \big( \mathbf{r}( \mathbf{\widehat{a}} )  \big) +  \mathbf{\widehat{S}} = \mathbf{\widehat{H}}  \ ,
\end{equation*}
from which we deduce that
\begin{equation*}
\mathbf{Q}^{T} \mathit{J} \big( \mathbf{z}( \mathbf{\widehat{a}} ) \big)^{T} \Big( \nabla^2 \phi  \big( \mathbf{z}(  \mathbf{\widehat{a}} )  \big) \Big) \mathit{J} \big( \mathbf{z}( \mathbf{\widehat{a}} ) \big) \mathbf{Q} = \mathbf{Q}^{T} \mathbf{\widehat{H}} \mathbf{Q} = \mathbf{\widehat{T}}  \ .
\end{equation*}
As the matrix $\mathbf{\widehat{T}}$ is positive definite by hypothesis, this also shows that $\nabla^2 \phi ( \mathbf{z}(  \mathbf{\widehat{a}} ) )$ is definite positive since the proposition that $\nabla^2 \phi ( \mathbf{z}(  \mathbf{\widehat{a}} ) )$  is definite positive is equivalent to the proposition that $\mathbf{\widehat{T}}$ is definite positive as noted above. Furthermore, this implies immediately that the first-order stationary point $ \mathbf{z}( \mathbf{\widehat{a}} )$ is a strict local minimizer of $\phi(.)$ and that the first-order stationary point $\mathbf{\widehat{a}}$ is a non strict local minimizer of $\psi(.)$ as claimed in the theorem.

Finally, since $\nabla^2 \phi ( \mathbf{z}(  \mathbf{\widehat{a}} ) )$ is definite positive and we have the equality
\begin{equation*}
\mathbf{\widehat{H}}  = \mathit{J}  \big( \mathbf{z}( \mathbf{\widehat{a}} )  \big)^{T} \Big( \nabla^2 \phi  \big( \mathbf{z}(  \mathbf{\widehat{a}} )  \big) \Big) \mathit{J}  \big( \mathbf{z}( \mathbf{\widehat{a}} )  \big) \ ,
\end{equation*}
it is easy to deduce that $\emph{null} ( \mathbf{\widehat{H}} ) = \emph{null} \big( \mathit{J} ( \mathbf{z}( \mathbf{\widehat{a}} ) )  \big) = \emph{null} \big( \mathit{J} ( \mathbf{r}( \mathbf{\widehat{a}} ) )  \big)$ and $\emph{dim}   \big( \emph{null} ( \mathbf{\widehat{H}} )  \big) = k.k$, using the hypothesis $\emph{rank} \big( \mathit{J}  \big( \mathbf{r}( \mathbf{\widehat{a}} ) )  \big) = r$ and the rank-nullity relationship~\ref{eq:rank}. This last equality implies, finally, that $\emph{rank} ( \mathbf{\widehat{H}} ) = r $ again by the rank-nullity theorem~\ref{eq:rank}, which concludes the demonstration of the theorem.

\end{proof}

Interestingly, under the same hypotheses as used in Theorem~\ref{theo5.10:box} and when, in addition, the $r \times r$ symmetric matrix $\mathbf{\widehat{T}}$ is definite positive and $\mathbf{\widehat{a}}$ is thus a (local) minimizer of $\psi(.)$, the continuum of (local) minimizers of $\psi(.)$ which achieve the same minimum value of $\psi(\mathbf{\widehat{a}})$, say $\mathcal{S} ( \mathbf{\widehat{a}} )$ defined formally as follows, 
\begin{equation} \label{eq:psi_minima}
\mathcal{S} ( \mathbf{\widehat{a}} ) = \left \lbrace  \mathbf{a} \in \mathbb{R}^{p.k} \, /  \,   \psi( \mathbf{a}) =  \psi( \mathbf{\widehat{a}}) \text{ and } \mathbf{A} = \big ( \emph{mat}_{k \times p}( \mathbf{a} ) \big)^{T}  \in  \mathbb{R}^{p \times k}_{k} \right \rbrace \ ,
\end{equation}
is locally well-behaved around $\mathbf{\widehat{a}}$ in the sense that this continuum forms a smooth (e.g., $C^{1}$) submanifold of $\mathbb{R}^{k.p}$ (locally) around $\mathbf{\widehat{a}}$ of dimension $k.k$ and the tangent linear space to $\mathcal{S} ( \mathbf{\widehat{a}} )$  at $\mathbf{\widehat{a}}$ is exactly the kernel of  $\nabla^2 \psi( \mathbf{\widehat{a}} )$, i.e.,
\begin{equation*}
 \emph{null} \big(   \nabla^2 \psi( \mathbf{\widehat{a}} )  \big) = \mathcal{T}_{\mathbf{\widehat{a}}} \mathcal{S} ( \mathbf{\widehat{a}} ) \ .
\end{equation*}
Note that, this is the best, we can hope for smooth functions like $\psi(.)$ with singular (local)  minima (because of over-parameterization) and this expresses the fact that the restriction of  $\nabla^2 \psi( \mathbf{\widehat{a}})$ to the normal linear space of  $\mathcal{S} ( \mathbf{\widehat{a}} )$  at $\mathbf{\widehat{a}}$, say  $\mathcal{N}_{\mathbf{\widehat{a}}} \mathcal{S} ( \mathbf{\widehat{a}} )$ (i.e., $\mathcal{N}_{\mathbf{\widehat{a}}} \mathcal{S} ( \mathbf{\widehat{a}} ) = \big( \mathcal{T}_{\mathbf{\widehat{a}}} \mathcal{S} ( \mathbf{\widehat{a}} ) \big)^{\bot}$), is definite positive (e.g., $\nabla^2 \psi( \mathbf{\widehat{a}})$ is definite positive along its normal directions) as demonstrated in Theorem~\ref{theo5.10:box}.

This nice geometrical property is called the Morse-Bott property in Rebjock and Boumal~\cite{RB2024b}, see their Definition 1.1. They further show that this local Morse-Bott property is essentially equivalent to various (local) structural assumptions like the Polyak–Lojasiewicz condition, Quadratic Growth and  Error Bound properties, which have been used in the past to explain the surprising local convergence with a fast local rate of (quasi-)Newton methods in some neighborhood of non-isolated optima; e.g., for minimizing  $C^{2}$ cost functions like $\psi(.)$ for which the Hessian at (local) minima is at best positive semi-definite, but never positive definite because local minima are never isolated due to over-parameterization.

In order to demonstrate that $\psi(.)$ satisfies effectively the Morse-Bott property at $\mathbf{\widehat{a}}$, let us first give a precise definition of a smooth (i.e., $C^{1}$) submanifold around one of its points in $\mathbb{R}^{n}$, taken from Example 6.8 in~\cite{RW1998}. 
\begin{def5.1} \label{def5.1:box}
Let $\mathcal{C}$ be a nonempty subset of $\mathbb{R}^{n}$.

We say that $\mathcal{C}$ is a $d$-dimensional smooth submanifold in $\mathbb{R}^{n}$ around the point $\mathbf{\widehat{a}} \in \mathcal{C}$, if $\mathcal{C}$ can be represented relative to an open neighborhood $U$ of $\mathbf{\widehat{a}}$ as the set of solutions to $g( \mathbf{a} ) = \mathbf{0}^{m}$ where $g(.)$ is a $C^{1}$ mapping from $U$ to $\mathbb{R}^{m}$ with its differential at $\mathbf{\widehat{a}}$, $g^{'} (\mathbf{\widehat{a}})$, being a surjective linear operator, which is equivalent to say that the Jacobian matrix  $\mathit{J}( g (\mathbf{\widehat{a}}) ) \in \mathbb{R}^{m \times n} $ is of full rank m, where $m = n - d $.

Thus, equivalently, $\mathcal{C}$ is a $d$-dimensional smooth submanifold in $\mathbb{R}^{n}$ around the point $\mathbf{\widehat{a}} \in \mathcal{C}$, if
\begin{equation*}
\mathcal{C} \cap U = \lbrace  \mathbf{a} \in U \, /  \, g( \mathbf{a} ) = \mathbf{0}^{m}  \rbrace \quad  \text{ and } \quad  \emph{rank} \Big( \mathit{J} \big( g (\mathbf{\widehat{a}}) \big) \Big ) = m = n - d \ .
\end{equation*}
Such function $g(.)$, if it exists, is called a local defining function for $\mathcal{C}$ at $\mathbf{\widehat{a}}$ as in the Definition~\ref{def2.4:box} of a classical $C^{p}$ embedded submanifold of $\mathbb{R}^{n}$ given in Subsection~\ref{calculus:box}.
\\
\end{def5.1}
Obviously, if $\mathcal{M}$ is a smooth embedded submanifold of $\mathbb{R}^{n}$ in the sense of Definition~\ref{def2.4:box}, $\mathcal{M}$ is also a $C^{1}$ submanifold around each of its points in the sense of Definition~\ref{def5.1:box}, but the converse is not true in general. Thus, a smooth submanifold around one of its points, is just the "local" version of a classical smooth submanifold embedded in $\mathbb{R}^{n}$.

Next, we will use the definition of tangent vectors to an arbitrary nonempty subset of a normed vector space, given in Definition~\ref{def2.5:box}, and a finite version of the Lyusternik Theorem~\cite{L1934} (given below), which gives a complete characterization of the set of tangent vectors to a subset of a normed vector space of the form
\begin{equation*}
\mathcal{S} = \lbrace  \mathbf{a} \in U \, /  \, g( \mathbf{a} ) = \mathbf{0}^{m}  \rbrace \ ,
\end{equation*}
where $U$ is an open set of $\mathbb{R}^{n}$ and $g(.)$ is a function from $U$ to  $\mathbb{R}^{m}$ of class $C^{1}$.
\begin{theo5.11} \label{theo5.11:box}
Let $U$ be an open subset of $\mathbb{R}^{n}$ and $g(.)$ a $C^{1}$ function from $U$ to  $\mathbb{R}^{m}$. Further, define the subset of $\mathbb{R}^{n}$
\begin{equation*}
\mathcal{S} = \lbrace  \mathbf{a} \in U \, /  \, g( \mathbf{a} ) = \mathbf{0}^{m}  \rbrace \ ,
\end{equation*}
and let $\mathbf{\widehat{a}} \in \mathcal{S}$.

If the differential of $g(.)$ at $\mathbf{\widehat{a}}$, $g^{'}(\mathbf{\widehat{a}})$, is a surjective linear operator from  $\mathbb{R}^{n}$ to $\mathbb{R}^{m}$, which is equivalent to say that $ \emph{rank} (  \mathit{J} ( g (\mathbf{\widehat{a}}) ) ) = m$, then 
\begin{equation*}
\emph{null} \big(  \mathit{J}( g(\mathbf{\widehat{a}}) ) \big) = \mathcal{T}_{\mathbf{\widehat{a}}} \mathcal{S} \ ,
\end{equation*}
where  $\mathit{J}( g(\mathbf{\widehat{a}}) )$ is the Jacobian matrix of  $g(.)$ at  $\mathbf{\widehat{a}}$ and $\mathcal{T}_{\mathbf{\widehat{a}}} \mathcal{S}$ is the set of tangent vectors to $\mathcal{S}$ at $\mathbf{\widehat{a}}$ in the sense of Definition~\ref{def2.5:box}, which is thus a linear subspace of $\mathbb{R}^{n}$ of dimension $d = n - m$.

\end{theo5.11}
\begin{proof}
Omitted. See~\cite{L1934} and~\cite{S1973} for proof.
\end{proof}
The Lyusternik Theorem~\ref{theo5.11:box} is the angular stone behind the Definition~\ref{def5.1:box} of an embedded $C^{1}$ submanifold $\mathcal{C}$ around one of its points in $\mathbb{R}^{n}$. Furthermore, it allows an easy identification of the tangent space to such "local" submanifolds as the kernel of the Jacobian of the local defining function for $\mathcal{C}$ at this particular point and also of the dimension of such smooth submanifold around one of its points. These "local" results are further consistent with those concerning classical smooth  submanifolds embedded in normed vector spaces discussed in Subsection~\ref{calculus:box}.

With these preliminaries, we are now in position to demonstrate that the continuum of (local) minima of  $\psi(.)$, $\mathcal{S} ( \mathbf{\widehat{a}} )$, defined in equation~\eqref{eq:psi_minima}, is effectively a smooth (e.g., $C^{1}$) submanifold of $\mathbb{R}^{k.p}$ (locally) around $\mathbf{\widehat{a}}$ of dimension $k.k$ in the sense of Definition~\ref{def5.1:box} and that
\begin{equation*}
 \emph{null} \big(   \nabla^2 \psi( \mathbf{\widehat{a}} )  \big) = \mathcal{T}_{\mathbf{\widehat{a}}} \mathcal{S} ( \mathbf{\widehat{a}} ) \ ,
\end{equation*}
if we assume, as in Theorem~\ref{theo5.10:box}, that:

- $\mathbf{\widehat{a}} \in \mathbb{R}^{k.p}$, with $ \widehat{\mathbf{A}} = ( \emph{mat}_{k \times p}( \widehat{\mathbf{a}} ) )^{T}  \in  \mathbb{R}^{p \times k}_{k}$, is a first-order stationary point of  $\psi (.)$;

- $\mathbf{F}(\mathbf{a})$ has full column-rank and  the rank of the Jacobian matrix $\mathit{J}( \mathbf{r}(\mathbf{a}) )$ is equal to $r = (p - k).k$ in a neighborhood of $ \mathbf{\widehat{a}} $.

- and, finally, the $r \times r$ symmetric matrix  $\mathbf{\widehat{T}}$ defined in equation~\eqref{eq:T_mat} is positive definite.

Under these hypotheses, we know from Theorem~\ref{theo5.10:box}, that  $\mathbf{\widehat{a}}$ is a (local) minimizer of $\psi(.)$ and we can assert that the set $\mathcal{S} ( \mathbf{\widehat{a}} )$ defined in equation~\eqref{eq:psi_minima} is nonempty and even forms a continuum of points. Furthermore, again as a direct consequence of Theorem~\ref{theo5.10:box}, we also have
\begin{equation*}
 \emph{null} \big(   \nabla^2 \psi( \mathbf{\widehat{a}} )  \big) = \emph{null} \big( \mathit{J}( \mathbf{r}(\mathbf{\widehat{a}}) )  \big) \quad \text{and} \quad  \emph{dim} \big(    \emph{null} \big(   \nabla^2 \psi( \mathbf{\widehat{a}} )  \big)  \big) = k.k  \ .
\end{equation*}
We now first demonstrate that $\mathcal{S} ( \mathbf{\widehat{a}} )$ is effectively an embedded  $k.k$-dimensional $C^{1}$ submanifold in $\mathbb{R}^{k.p}$ around the point $\mathbf{\widehat{a}} \in \mathcal{S} ( \mathbf{\widehat{a}} )$ in the sense of Definition~\ref{def5.1:box} and that the tangent space to $\mathcal{S} ( \mathbf{\widehat{a}} )$ around $\mathbf{\widehat{a}}$ is nothing else than the kernel of $\mathit{J}( \mathbf{r}(\mathbf{\widehat{a}}) )$ using the Lyusternik Theorem~\ref{theo5.11:box}.

To verify Definition~\ref{def5.1:box} for $\mathcal{S} ( \mathbf{\widehat{a}} )$, we need to find an open neighborhood $U$ of $\mathbf{\widehat{a}}$ in $\mathbb{R}^{k.p}$ and a  $C^{1}$ mapping $g(.)$ from $U$ to $\mathbb{R}^{r}$ with $r = k.(p - k)$ such that
\begin{equation*}
\mathcal{S} ( \mathbf{\widehat{a}} ) \cap U = \lbrace  \mathbf{a} \in U \, /  \, g( \mathbf{a} ) = \mathbf{0}^{r}  \rbrace \quad  \text{ and } \quad   \emph{rank} \big( \mathit{J}( \mathbf{g}(\mathbf{\widehat{a}}) )  \big) = r  \ .
\end{equation*}

However, as in the demonstration of Theorem~\ref{theo5.10:box}, the above hypotheses imply that there exist an open neighborhood $\Upsilon$ of $ \mathbf{\widehat{a}} $ in $\mathbb{R}^{k.p}$, a twice continuously differentiable function $\mathbf{z}(.)$ from $\Upsilon$ to $\mathbb{R}^{r}$ and a twice continuously differentiable function $\mathbf{h}(.)$ from $\mathbb{R}^{r}$ to $\mathbb{R}^{k.p}$ such that $\mathbf{r}(\mathbf{a}) = \mathbf{h}(\mathbf{z}(\mathbf{a}) )$ and the Jacobian matrices, $\mathit{J}( \mathbf{z}(\mathbf{a}) )$ and $\mathit{J}( \mathbf{h}( \mathbf{z}(\mathbf{a}) ) )$,  have  a constant rank equals to $r = (p - k).k$  for all $\mathbf{a} \in \Upsilon$.

Furthermore, we can also define a twice continuously differentiable real function $\phi(.)$ from $\mathbb{R}^{r}$ to $\mathbb{R}$ such that $\phi ( \mathbf{o} ) =  \frac{1}{2} \Vert  \mathbf{h}( \mathbf{o}  )  \Vert^{2}_{2} \, , \, \forall \mathbf{o} \in \mathbb{R}^{r}$, and, again according to the demonstration of Theorem~\ref{theo5.10:box}, the point $\mathbf{z}( \mathbf{\widehat{a}} )$ is a strict (local) minimizer of $\phi(.)$. This implies the existence of an open neighborhood $V$ of  $\mathbf{z}( \mathbf{\widehat{a}} )$ in $\mathbb{R}^{r}$ such that ,$\forall \mathbf{o} \in V$, we have $\phi( \mathbf{o} ) > \phi( \mathbf{z}( \mathbf{\widehat{a}} ) )$. As the function $\mathbf{z}(.)$ is continuous over $\Upsilon$, the set $U = \mathbf{z}^{-1} (V)$ is  of the form $U = \Upsilon  \cap W$, where $W$ is an open set  of  $\mathbb{R}^{k.p}$, and $U$ is thus an open  neighborhood of $ \mathbf{\widehat{a}} $ in $\mathbb{R}^{k.p}$. Clearly, the elements $\mathbf{a} \in \mathcal{S} ( \mathbf{\widehat{a}} ) \cap U$ verify $ \mathbf{z}( \mathbf{a} ) = \mathbf{z}( \mathbf{\widehat{a}} )$.

In these conditions, $U$ is an open neighborhood of $\mathbf{\widehat{a}}$  and we can define a  twice continuously differentiable function $g(.)$ from $U$ to  $\mathbb{R}^{r}$ by
\begin{equation*}
g :  U \longrightarrow \mathbb{R}^{r}  : \mathbf{a} \mapsto  g( \mathbf{a} ) =  \mathbf{z}( \mathbf{a} ) - \mathbf{z}( \mathbf{\widehat{a}} ) \ ,
\end{equation*}
and we have effectively
\begin{equation*}
\mathcal{S} ( \mathbf{\widehat{a}} )  \cap U = \left \lbrace  \mathbf{a} \in U \, /  \, g( \mathbf{a} ) = \mathbf{0}^{m}  \right \rbrace  \ .
\end{equation*}
Furthermore, the Jacobian matrix $ \mathit{J} \big( g (\mathbf{a} ) \big)  =  \mathit{J} \big(  \mathbf{z} (\mathbf{a} ) \big)  \in \mathbb{R}^{r \times k.p}$ verifies $\emph{rank} \Big(  \mathit{J} \big( g (\mathbf{a} ) \big) \Big ) =  r , \forall \mathbf{a}  \in  U$ (see the demonstration of Theorem~\ref{theo5.10:box}). This implies in particular that $ \emph{rank} \Big( \mathit{J} \big( g (\mathbf{\widehat{a}}) \big) \Big ) = r $ and $g^{'} ( \mathbf{\widehat{a}} )$ is a surjective linear mapping from $\mathbb{R}^{k.p}$ to $\mathbb{R}^{r}$ and we conclude that $\mathcal{S} ( \mathbf{\widehat{a}} )$ is effectively an embedded  $C^{1}$ submanifold  around the (local) minimizer $\mathbf{\widehat{a}}$ in $\mathbb{R}^{k.p}$ of dimension $k.p - r = k.k$.

Moreover, a closer look at the preceding demonstration further shows that the set $\mathcal{S} ( \mathbf{\widehat{a}} )  \cap U$ is in fact also  an embedded  $k.k$-dimensional $C^{2}$ submanifold in $\mathbb{R}^{k.p}$ in the sense of Definition~\ref{def2.4:box} because the mapping $g(.)$ from $U$ to $\mathbb{R}^{r}$, defined above, is of class $C^{2}$ and a valid local defining function for $\mathcal{S} ( \mathbf{\widehat{a}} )  \cap U$ at all $\mathbf{a} \in \mathcal{S} ( \mathbf{\widehat{a}} )  \cap U$.

Let us now determine the tangent space to $\mathcal{S} ( \mathbf{\widehat{a}} )  \cap U$ at $\mathbf{a} \in \mathcal{S} ( \mathbf{\widehat{a}} )  \cap U$. By applying the Lyusternik Theorem~\ref{theo5.11:box}  to the set $\mathcal{S} ( \mathbf{\widehat{a}} )  \cap U$, we have immediately
\begin{equation*}
\emph{null} \Big(  \mathit{J} \big( g (\mathbf{a} ) \big) \Big ) = \emph{null} \Big(  \mathit{J} \big(  \mathbf{z} (\mathbf{a} ) \big) \Big ) =  \mathcal{T}_{\mathbf{a}}  \big( \mathcal{S} ( \mathbf{\widehat{a}} )  \cap U   \big) \quad , \forall  \mathbf{a} \in \mathcal{S} ( \mathbf{\widehat{a}} )  \cap U \ .
\end{equation*}
Finally, using the equality
\begin{equation*}
\emph{null} \big( \mathit{J} ( \mathbf{r}(\mathbf{a}) ) \big) = \emph{null} \big( \mathit{J}( \mathbf{z}(\mathbf{a}) ) \big) \, , \forall  \mathbf{a} \in \Upsilon \ ,
\end{equation*}
established in the demonstration of Theorem~\ref{theo5.10:box}), we get the equality
\begin{equation*}
\emph{null} \big( \mathit{J} ( \mathbf{r}(\mathbf{a}) ) \big) = \mathcal{T}_{\mathbf{a}}  \big( \mathcal{S} ( \mathbf{\widehat{a}} )  \cap U   \big) \quad , \forall  \mathbf{a} \in \mathcal{S} ( \mathbf{\widehat{a}} )  \cap U \ .
\end{equation*}
This shows that the tangent space to $\mathcal{S} ( \mathbf{\widehat{a}} )  \cap U$ at $\mathbf{a} \in \mathcal{S} ( \mathbf{\widehat{a}} )  \cap U$ is nothing else than the kernel of the Jacobian matrix $\mathit{J} ( \mathbf{r}(\mathbf{a}) )$, $ \forall  \mathbf{a} \in \mathcal{S} ( \mathbf{\widehat{a}} )  \cap U$. This implies in particular that
\begin{equation*}
\emph{null} \big( \mathit{J} ( \mathbf{r}( \mathbf{\widehat{a}} ) ) \big) = \mathcal{T}_{\mathbf{\widehat{a}}}  \big( \mathcal{S} ( \mathbf{\widehat{a}} )  \cap U   \big)  \ .
\end{equation*}
Finally, using the hypothesis that $\mathbf{\widehat{T}}$ (defined in equation~\eqref{eq:T_mat}) is positive definite and Theorem~\ref{theo5.10:box}, we obtain the equality
\begin{equation*}
 \emph{null} \big(   \nabla^2 \psi( \mathbf{\widehat{a}} )  \big) = \emph{null} \big( \mathit{J}( \mathbf{r}(\mathbf{\widehat{a}}) )  \big) =  \mathcal{T}_{\mathbf{\widehat{a}}}  \big( \mathcal{S} ( \mathbf{\widehat{a}} )  \cap U   \big)   \ ,
\end{equation*}
which demonstrates that the cost function $\psi(.)$ effectively verifies the so-called Morse-Bott property at the (local) minimizer $\mathbf{\widehat{a}}$ under the hypotheses of Theorem~\ref{theo5.10:box} and when, in addition, the $r \times r$ symmetric matrix  $\mathbf{\widehat{T}}$ defined in equation~\eqref{eq:T_mat} is positive definite as claimed above. Note that if $\mathbf{\widehat{a}}$ is a (local) minimizer of $\psi(.)$, but we don't assume the hypothesis $\mathbf{\widehat{T}}$ is positive definite, we can still get the inclusion 
\begin{equation*}
\mathcal{T}_{\mathbf{\widehat{a}}}  \big( \mathcal{S} ( \mathbf{\widehat{a}} )  \cap U   \big) \subset  \emph{null} \big(   \nabla^2 \psi( \mathbf{\widehat{a}} )  \big)  \ ,
\end{equation*}
by using Theorem~\ref{theo5.9:box} and Corollary~\ref{corol5.9:box}, and the fact that $\mathbf{\widehat{a}}$ is a first-order stationary point of $\psi(.)$.

In order to show now that Theorem~\ref{theo5.10:box} covers a large body of real applications, we note that the matrix variable $\mathbf{A}$ is always requested to have full column rank to assure the differentiability of $\psi(.)$ and that $\mathbf{F}( \mathbf{a} )$ is of full column rank in most cases as soon as every line or column of the data matrix $\mathbf{X}$ has at least $k$ non-missing entries. Furthermore, the hypothesis that the rank of the Jacobian matrix $\mathit{J}( \mathbf{r}( \mathbf{a} ) )$ is constant and equal to $r = (p-k).k$ in a neighborhood of $\mathbf{\widehat{a}}$ is automatically verified if the hypotheses of Theorem~\ref{theo5.3:box} are fulfilled in a neighborhood of $\mathbf{\widehat{a}}$. Finally, we note that  orthonormal bases given, respectively, by the columns of $\mathbf{P}$ and $\mathbf{Q}$ for the null space of $\mathit{J}( \mathbf{r}( \mathbf{\widehat{a}} ) )$ and its orthogonal complement, used in Theorem~\ref{theo5.10:box}, can be easily computed from the results of Corollary~\ref{corol5.6:box}.

Obviously, Theorem~\ref{theo5.10:box} also justifies the extension of the regularization techniques, developed in Subsection~\ref{jacob:box} to overcome the systematic singularity of the Jacobian matrix (or its approximation) in the Gauss-Newton or Levenberg-Marquardt methods, to the Newton method using the full Hessian matrix $\mathbf{H}$ or its two-term approximation $\bar{\mathbf{H}}$ derived above.
In other words, to avoid the singularity or ill-conditioning of these two symmetric matrices near first-order stationary points of  $\psi(.)$, the Newton correction vector $d\mathbf{a}_{n} \in \mathbb{R}^{k.p}$ can be found in a two-step procedure at each iteration, as for the Gauss-Newton algorithm described in Subsection~\ref{jacob:box}. In a first step, we solve the $( p - k ).k \times ( p - k ).k$ symmetric linear system
\begin{equation}  \label{eq:newton_sys}
(\mathbf{\bar{O}}^{\bot} )^{T}\mathbf{H} \mathbf{\bar{O}}^{\bot} d\mathbf{\bar{a}}_{n} = -  ( \mathbf{\bar{O}}^{\bot}  )^{T}  \mathit{J} \big( \mathbf{r}(\mathbf{a})  \big)^{T}  \mathbf{r}(\mathbf{a}) =  \big( \mathbf{M}(\mathbf{a}) \mathbf{\bar{O}}^{\bot}   \big)^{T} \mathbf{r}(\mathbf{a})  \ ,
\end{equation}
or, if we use the two-term approximation of the Hessian, $\bar{\mathbf{H}}$ defined in equation~\eqref{eq:approx_hess_mat},
\begin{equation} \label{eq:approx_newton_sys}
(\mathbf{\bar{O}}^{\bot} )^{T} \bar{\mathbf{H}} \mathbf{\bar{O}}^{\bot} d\mathbf{\bar{a}}_{n} =  \big( \mathbf{M}(\mathbf{a}) \mathbf{\bar{O}}^{\bot}   \big)^{T} \mathbf{r}(\mathbf{a}) \ ,
\end{equation}
for  $d\mathbf{\bar{a}}_{n} \in \mathbb{R}^{(p-k).k}$ and where $\mathbf{\bar{O}}^{\bot} =  \mathbf{K}_{(p,k)} ( \mathbf{I}_{k}  \otimes  \mathbf{O}^{\bot}) \in \mathbb{O}^{k.p \times ( p - k ).k }$ and $\mathbf{O}^{\bot} \in \mathbb{O}^{p \times ( p - k ) }$  are orthonormal matrices whose columns form, respectively,  a basis of   $\emph{null} \big( \mathit{J}( \mathbf{r}(\mathbf{a} ) )  \big )^{\bot}$ and $ \emph{ran}( \mathbf{A} )^{\bot}$, see Corollary~\ref{corol5.6:box} and Theorem~\ref{theo5.7:box}  for details. In a second step,  to get the Newton correction vector we need to compute the following matrix-vector product
\begin{equation*}
d\mathbf{a}_{n} = \mathbf{\bar{O}}^{\bot} d\mathbf{\bar{a}}_{n} \in \mathbb{R}^{k.p} \ ,
\end{equation*}
or, equivalently, in matrix format,
\begin{equation*}
d\mathbf{A}_{n} = \mathbf{O}^{\bot} d\mathbf{\bar{A}}_{n} \in  \mathbb{R}^{p  \times k } \ .
\end{equation*}
This modification of the Newton algorithm in the context of the~\eqref{eq:VP1} problem, first suggested by Chen~\cite{C2008b}, has a strong theoretical justification as it can be interpreted as a Riemannian Newton operating on the (quotient) Grassmann manifold $\text{Gr}(p,k)$ as we will show below, but it is also computationally very expensive as $\mathbf{\bar{O}}^{\bot}$ is huge matrix in most cases.

Alternatively, we can use a cheaper alternative based on the orthogonality constraint
\begin{equation*}
 \mathbf{N}^{T} d\mathbf{a}_{n} =  \mathbf{0}^{k.k} \ ,
\end{equation*}
where $\mathbf{N} = \mathbf{K}_{(p,k)} ( \mathbf{I}_{k} \otimes \mathbf{A} )$ or $\mathbf{N} = \mathbf{K}_{(p,k)} ( \mathbf{I}_{k} \otimes \mathbf{O} )$, and the columns of the matrix $\mathbf{A}$ ($\mathbf{O}$) form a (orthonormal) basis of $\emph{ran}( \mathbf{A} )$ and the columns of $\mathbf{N}$ is a (orthonormal) basis of $\emph{null} \big( \mathit{J}( \mathbf{r}(\mathbf{a} ) )  \big) = \emph{null} \big( \mathbf{M}(\mathbf{a} ) \big)$, as demonstrated in Corollary~\ref{corol5.6:box} and discussed in Subsection~\ref{jacob:box}, and proceed in one step by solving the damped symmetric linear system
\begin{equation}  \label{eq:newton_sys2}
\big( \mathbf{H} +  \mathbf{N} \mathbf{N}^{T} \big)  d\mathbf{a}_{n}  =  - \nabla \psi( \mathbf{a} ) = \mathbf{M}(\mathbf{a})^{T} \mathbf{r}(\mathbf{a})
\end{equation}
or 
\begin{equation}  \label{eq:approx_newton_sys2}
\big(   \bar{\mathbf{H}} +  \mathbf{N} \mathbf{N}^{T} \big)  d\mathbf{a}_{n}  =  - \nabla \psi( \mathbf{a} ) = \mathbf{M}(\mathbf{a})^{T} \mathbf{r}(\mathbf{a}) \ ,
\end{equation}
if we use the two-term approximation of the Hessian matrix. Note that, in both cases, the term $\mathbf{N} \mathbf{N}^{T}$ can be efficiently evaluated as
\begin{equation}  \label{eq:NNt_mat}
\mathbf{N} \mathbf{N}^{T} = \mathbf{K}_{(p,k)} ( \mathbf{I}_{k} \otimes \mathbf{A} \mathbf{A}^{T}) \mathbf{K}_{(k,p)}  \quad \text{or}  \quad \mathbf{N} \mathbf{N}^{T} = \mathbf{K}_{(p,k)} ( \mathbf{I}_{k} \otimes \mathbf{O} \mathbf{O}^{T}) \mathbf{K}_{(k,p)} \ .
\end{equation}
This  last approach using a linear constraint is new in the context of Newton methods, but is simply an extension to (quasi-)Newton methods of the technique first proposed by Okatani et al.~\cite{OYD2011} for Gauss-Newton and Levenberg-Marquardt methods. The variant of the Newton method, defined by equation~\eqref{eq:approx_newton_sys2}, can also be interpreted as a Riemannian quasi-Newton operating on the (quotient) Grassmann manifold $\text{Gr}(p,k)$, but not the first one defined by  equation~\eqref{eq:newton_sys2} as we will  see below.

Assuming that the dimension of the null space of $\mathbf{H}$ or $\bar{\mathbf{H}}$ at  first-order stationary points of $\psi(.)$ is equal to the dimension of the null space of the Jacobian matrix and that both are equal to $k.k$, these different approaches will efficiently overcome the systematic singularity of $\bar{\mathbf{H}}$ or those of $\mathbf{H}$ at these first-order stationary points of $\psi(.)$. In these conditions, the above symmetric linear systems will have an unique solution and the (quasi-)Newton direction is thus well defined. Of course, as always in Newton methods, one common weakness of these two second-order approaches for solving the~\eqref{eq:VP1} problem is that $\mathbf{H}$ or its two-term approximation $\bar{\mathbf{H}}$ may not be positive definite at some points in a region of mixed curvature. In such condition, the Newton direction may not be in a descent direction and the above linear systems cannot be solved by a simple Cholesky factorization. However, many standard techniques  are available to deal with this classical problem~\cite{NW2006}\cite{MN2010} and can be applied here in our specific WLRA context without any modification as we will discussed in Subsection~\ref{vp_n_alg:box}.

Alternatively, we can again recast the WLRA problem in its variable projection formulation as an optimization problem on the (quotient) Grassmann manifold $\text{Gr}(p,k)$~\cite{AMS2008}\cite{B2023} and use a Riemannian Newton method to solve it as was done for example in~\cite{BA2011}\cite{BA2015}. To clarify the differences and similarities between the Riemannian Newton method operating on $\text{Gr}(p,k)$ and the above (quasi)-Newton algorithms operating on $\mathbb{R}^{p \times k}_{k}$ (or on $\text{St}(p,k) = \mathbb{O}^{p \times k}$), we now investigate the relationships between the Euclidean Hessian of $\psi(.)$  at $\mathbf{a}$, given in equation~\eqref{eq:H_hess_mat} (e.g., when $\psi(.)$ is considered as a real function from  $\mathbb{R}^{p.k}$ to $\mathbb{R}$, see equation~\eqref{eq:psi_func}), and the  Riemannian Hessian of the unvectorized form of  $\psi(.)$ at $\mathbf{A} \in \mathbb{R}^{p \times k}_{k}$ (e.g., the real function $\psi  \circ  h^{-1}(.)$ from  $\mathbb{R}^{p \times k}_{k}$ to $\mathbb{R}$, where $h^{-1}(.)$ is defined in equation~\eqref{eq:h_func} of Subsection~\ref{varpro_wlra:box} with $h^{-1} ( \mathbf{A} ) = \emph{vec}(  \mathbf{A}^{T} ) = \mathbf{a}, \forall \mathbf{A} \in \mathbb{R}^{p \times k}_{k}$), when this cost function is considered abusively as defined on the Grassmann manifold $\text{Gr}(p,k)$, as already discussed at the end of Subsection~\ref{jacob:box} and at the beginning of this subsection. The results we present now are a slight extension, with our notations, of those given in~\cite{HF2015b}.

As in our previous discussion on the connections between the Euclidean gradient of $\psi(.)$  and the Riemannian gradient of $\psi  \circ  h^{-1}(.)$, to simplify the presentation we require that  $\mathbf{W} \in  \mathbb{R}^{ p \times n  }_{+*}$ and that each element $\mathring{\mathbf{O}}  \in \text{Gr}(p,k)$ is represented by an element of the compact Stiefel manifold $\mathbf{O} \in \text{St}(p,k)$, in line with previous works on Riemannian optimization on $\text{Gr}(p,k)$~\cite{DKM2012}\cite{BA2015}\cite{B2023}. Recall also that  any element of the tangent space of $\text{Gr}(p,k)$ at $\mathring{\mathbf{O}}$,  $\mathcal{T}_{\mathring{\mathbf{O}}}  \text{Gr}(p,k)$, can be represented uniquely by an element of the following linear subspace of $\mathbb{R}^{p \times k }$ of dimension $(p-k).k$:
\begin{equation*}
\mathcal{T}_{\mathbf{O}} \text{Gr}(p,k) =  \big\{ \mathbf{D} \in \mathbb{R}^{ p \times k  }  \text{ } / \text{ } \mathbf{O}^{T} \mathbf{D} =  \mathbf{0}^{ k \times k  } \big\}  \ ,
\end{equation*}
as already noted in Subsection~\ref{jacob:box}. $\mathcal{T}_{\mathbf{O}} \text{Gr}(p,k)$ is nothing else than the horizontal space of $\text{St}(p,k)$ at $\mathbf{O} \in \text{St}(p,k)$, noted $\mathcal{H}_{\mathbf{O}} \text{St}(p,k)$ or $\mathcal{H}_{\mathbf{O}} \mathbb{O}^{p \times k}$ in Subsection~\ref{calculus:box}. In these conditions, using Theorem~\ref{theo5.8:box}, we have
\begin{equation*}
\nabla_{R} \psi  \circ  h^{-1} ( \mathring{\mathbf{O}} ) =  \left ( \mathbf{I}_{p} - \mathbf{O} \mathbf{O}^{T}  \right ) \nabla_{F} \psi  \circ  h^{-1} ( \mathbf{O} ) =  \nabla_{F} \psi  \circ  h^{-1} ( \mathbf{O} ) \ ,
\end{equation*}
where $\nabla_{R} \psi  \circ  h^{-1} (  \mathring{\mathbf{O}} )$ is the Riemannian gradient of the smooth cost function $\psi \circ  h^{-1}(.)$ (defined on the Grassmann manifold  $\text{Gr}(p,k)$) at $\mathring{\mathbf{O}} \in \text{Gr}(p,k)$ and $\nabla_{F} \psi  \circ  h^{-1} ( \mathbf{O} )$ is the standard Frobenius gradient of $\psi \circ  h^{-1}(.)$ at $\mathbf{O} \in \text{St}(p,k)$.

 In this framework, the Riemannian Hessian of the smooth map $\psi  \circ  h^{-1}(.)$ defined on the Grassmann manifold $\text{Gr}(p,k)$ at  $\mathring{\mathbf{O}} = \emph{ran}( \mathbf{O} )  \in \text{Gr}(p,k)$, is thus a linear transformation from $\mathcal{T}_{\mathbf{O}} \text{Gr}(p,k)$ to $\mathcal{T}_{\mathbf{O}} \text{Gr}(p,k)$. Furthermore, using the above equality between the Riemannian and Frobenius gradients and equation~\eqref{eq:D_rhess_mat2} in Subsection~\ref{calculus:box}, $\forall \mathbf{D} \in \mathcal{T}_{\mathbf{O}} \text{Gr}(p,k)$, the image of $\mathbf{D}$ by the Hessian $\nabla_{R}^{2} \psi  \circ  h^{-1} ( \mathring{\mathbf{O}} )$ considered as a linear operator from $\mathcal{T}_{\mathbf{O}} \text{Gr}(p,k)$ to $\mathcal{T}_{\mathbf{O}} \text{Gr}(p,k)$ is given by
\begin{equation} \label{eq:R_H_hess_mat}
\lbrack  \nabla_{R}^{2} \psi  \circ  h^{-1} ( \mathring{\mathbf{O}} ) \rbrack ( \mathbf{D} ) =   \left ( \mathbf{I}_{p} - \mathbf{O} \mathbf{O}^{T}  \right )  \lbrack  \nabla_{F}^{2} \psi  \circ  h^{-1} ( \mathbf{O} )  \rbrack  ( \mathbf{D} ) \ ,
\end{equation}
where $\nabla_{F}^{2} \psi  \circ  h^{-1} ( \mathbf{O} ) $ is the standard Frobenius Hessian of $\psi  \circ  h^{-1} (.)$ at $\mathbf{O} \in \text{St}(p,k)$ and  $\left ( \mathbf{I}_{p} - \mathbf{O} \mathbf{O}^{T}  \right )$ is the orthogonal projector onto $\emph{ran}(  \mathbf{O} )$ in $\mathbb{R}^{p}$, but also the orthogonal projector onto $\mathcal{T}_{\mathbf{O}} \text{Gr}(p,k)$ in $\mathbb{R}^{p \times k}$, when  $\mathcal{T}_{\mathbf{O}} \text{Gr}(p,k)$ is considered as a subspace of $\mathbb{R}^{p \times k}$ of dimension $( p - k).k$, which is also the horizontal space of $\text{St}(p,k)$ at $\mathbf{O} \in \text{St}(p,k)$ as noted above.

Thus, we have  $\mathbf{D} \in \mathcal{T}_{\mathbf{O}} \text{Gr}(p,k)$ and $\lbrack \nabla_{R}^{2} \psi  \circ  h^{-1} ( \mathring{\mathbf{O}} )  \rbrack ( \mathbf{D} )  \in \mathcal{T}_{\mathbf{O}} \text{Gr}(p,k)$, and both are considered as $p  \times k$ matrices in equation~\eqref{eq:R_H_hess_mat}. In these conditions, as noted in Subsection~\ref{jacob:box}, equivalently, we can vectorize both $\mathbf{D}$ and $\lbrack \nabla_{R}^{2} \psi  \circ  h^{-1} ( \mathring{\mathbf{O}} )\rbrack (  \mathbf{D} )$ as
\begin{equation*}
\mathbf{h} = \emph{vec} \big(  \mathbf{D}^{T} \big) \in \mathbb{R}^{p.k} \quad \text{and} \quad \nabla_{R}^{2} \psi ( \mathring{\mathbf{o}} )  \mathbf{d}  =  \emph{vec} \left(  \big( \lbrack \nabla_{R}^{2} \psi  \circ  h^{-1} ( \mathring{\mathbf{O}} ) \rbrack (  \mathbf{D} ) \big)^{T}  \right) \in \mathbb{R}^{p.k} \ ,
\end{equation*}
where now  $\nabla_{R}^{2} \psi  ( \mathring{\mathbf{o}} )$ is a $p.k \times p.k$ (asymmetric) matrix.

Furthermore, using equations~\eqref{eq:R_H_hess_mat},~\eqref{eq:mulmat_kronprod},~\eqref{eq:com_kron} and Corollary~\ref{corol5.6:box}, we have, $\forall \mathbf{d} \in \mathbb{R}^{p.k}$,
\begin{align*}
 \nabla_{R}^{2} \psi ( \mathring{\mathbf{o}} )  \mathbf{d}  & =  \emph{vec} \left(  \big(    \left ( \mathbf{I}_{p} - \mathbf{O} \mathbf{O}^{T}  \right )  \lbrack  \nabla_{F}^{2} \psi  \circ  h^{-1} ( \mathbf{O} )  \rbrack  ( \mathbf{D} )   \big)^{T}  \right) \\
  & =  \emph{vec} \left(  \big(  \lbrack  \nabla_{F}^{2} \psi  \circ  h^{-1} ( \mathbf{O} )  \rbrack  ( \mathbf{D} )   \big)^{T}  \left ( \mathbf{I}_{p} - \mathbf{O} \mathbf{O}^{T}  \right )  \right) \\
  & =  \left ( ( \mathbf{I}_{p} - \mathbf{O} \mathbf{O}^{T}  )   \otimes  \mathbf{I}_{k} \right ) \emph{vec} \left(  \big(  \lbrack  \nabla_{F}^{2} \psi  \circ  h^{-1} ( \mathbf{O} )  \rbrack  ( \mathbf{D} )   \big)^{T}  \right)  \\
  & =  \left ( \mathbf{I}_{p.k} - ( \mathbf{O} \mathbf{O}^{T}   \otimes  \mathbf{I}_{k} ) \right ) \emph{vec} \left(  \big(  \lbrack  \nabla_{F}^{2} \psi  \circ  h^{-1} ( \mathbf{O} )  \rbrack  ( \mathbf{D} )   \big)^{T}  \right)  \\
  & =  \left ( \mathbf{I}_{p.k} -  \bar{\mathbf{O}} \bar{\mathbf{O}}^{T}   \right ) \emph{vec} \left(  \big(  \lbrack  \nabla_{F}^{2} \psi  \circ  h^{-1} ( \mathbf{O} )  \rbrack  ( \mathbf{D} )   \big)^{T}  \right)  \\
  & =  \left ( \mathbf{I}_{p.k} -  \bar{\mathbf{O}} \bar{\mathbf{O}}^{T}   \right )  \nabla^{2} \psi ( \mathbf{o} ) \mathbf{d} \ ,
\end{align*}
where $\nabla^{2} \psi ( \mathbf{o} )$ is a $p.k \times p.k$ symmetric matrix and the orthogonal projector $\left ( \mathbf{I}_{p} - \mathbf{O} \mathbf{O}^{T}  \right )$ onto  $\mathcal{T}_{\mathbf{O}} \text{Gr}(p,k)$ is represented in vectorized form by the orthogonal projector onto  $\emph{null} \big( \mathit{J}( \mathbf{r}(\mathbf{o} ) )  \big )^{\bot} = \emph{null} \big( \mathbf{M}(\mathbf{o} )  \big )^{\bot}$ given by $ \left ( \mathbf{I}_{p.k} -  \bar{\mathbf{O}} \bar{\mathbf{O}}^{T}   \right )$, where the columns of $\mathbf{\bar{O}} =  \mathbf{K}_{(p,k)} ( \mathbf{I}_{k}  \otimes  \mathbf{O})$ are an orthonormal basis of $\emph{null} \big( \mathit{J}( \mathbf{r}(\mathbf{o} ) )  \big ) = \emph{null} \big( \mathbf{M}(\mathbf{o} )  \big )$,  as demonstrated in Corollary~\ref{corol5.6:box} of Subsection~\ref{jacob:box}.

This leads to the matrix equality
\begin{equation*}
 \nabla_{R}^{2} \psi ( \mathring{\mathbf{o}} ) =  \left ( \mathbf{I}_{p.k} -  \bar{\mathbf{O}} \bar{\mathbf{O}}^{T}   \right )  \nabla^{2} \psi ( \mathbf{o} ) \ ,
\end{equation*}
and, using the explicit form of $ \nabla^{2} \psi ( \mathbf{o} )$ given by equation~\eqref{eq:H_hess_mat}, we obtain:
\begin{align} \label{eq:R_H_hess_mat2}
 \nabla_{R}^{2} \psi ( \mathring{\mathbf{o}} ) & =  \left ( \mathbf{I}_{p.k} -  \bar{\mathbf{O}} \bar{\mathbf{O}}^{T}   \right )   \big (  \mathbf{M}(\mathbf{o})^{T} \mathbf{M}(\mathbf{o}) - \mathbf{L}(\mathbf{o})^{T} \mathbf{L}(\mathbf{o}) +  \mathbf{U}(\mathbf{o})^{T} \mathbf{L} (\mathbf{o}) + \mathbf{L}(\mathbf{o})^{T} \mathbf{U}(\mathbf{o}) \big)  \nonumber \\
  & =  \mathbf{M}(\mathbf{o})^{T} \mathbf{M}(\mathbf{o}) - \mathbf{L}(\mathbf{o})^{T} \mathbf{L}(\mathbf{o}) +    \left ( \mathbf{I}_{p.k} -  \bar{\mathbf{O}} \bar{\mathbf{O}}^{T}   \right )  \mathbf{U}(\mathbf{o})^{T} \mathbf{L} (\mathbf{o}) + \mathbf{L}(\mathbf{o})^{T} \mathbf{U}(\mathbf{o}) \ ,
\end{align}
since
\begin{align*}
 \left ( \mathbf{I}_{p.k} - \bar{\mathbf{O}} \bar{\mathbf{O}}^{T}   \right )  \mathbf{M}(\mathbf{o})^{T} & = \mathbf{M}(\mathbf{o})^{T} -  \bar{\mathbf{O}} ( \bar{\mathbf{O}}^{T}   \mathbf{M}(\mathbf{o})^{T} )  \\
  & = \mathbf{M}(\mathbf{o})^{T} -  \bar{\mathbf{O}} ( \mathbf{M}(\mathbf{o}) \bar{\mathbf{O}}  )^{T}  \\
  & = \mathbf{M}(\mathbf{o})^{T} \ ,
\end{align*}
as the columns of $\mathbf{\bar{O}}$ form an orthonormal basis of $\emph{null} \big( \mathbf{M}(\mathbf{o} )  \big )$ and, similarly,
\begin{align*}
 \left ( \mathbf{I}_{p.k} - \bar{\mathbf{O}} \bar{\mathbf{O}}^{T}   \right )  \mathbf{L}(\mathbf{o})^{T} & = \mathbf{L}(\mathbf{o})^{T} -  \bar{\mathbf{O}} ( \bar{\mathbf{O}}^{T}   \mathbf{L}(\mathbf{o})^{T} )  \\
  & = \mathbf{L}(\mathbf{o})^{T} -  \bar{\mathbf{O}} ( \mathbf{L}(\mathbf{o}) \bar{\mathbf{O}}  )^{T}  \\
  & = \mathbf{L}(\mathbf{o})^{T} \ ,
\end{align*}
as $\emph{null} \big( \mathit{J}( \mathbf{r}(\mathbf{o} ) )  \big ) = \emph{null} \big( \mathbf{M}(\mathbf{o} )  \big ) \cap  \emph{null} \big( \mathbf{L}(\mathbf{o} )  \big ) \subset  \emph{null} \big( \mathbf{L}(\mathbf{o} )  \big )$ and, thus, each column vector of  $\mathbf{\bar{O}}$ is also an element of $\emph{null} \big( \mathbf{L}(\mathbf{o} )  \big )$ when $\mathbf{W} \in  \mathbb{R}^{ p \times n  }_{+*}$ according to Corollary~\ref{corol5.5:box}.

Since $\nabla_{R}^{2} \psi  \circ  h^{-1} ( \mathring{\mathbf{O}} )$ is a linear operator from $\mathcal{T}_{\mathbf{O}} \text{Gr}(p,k)$ to $\mathcal{T}_{\mathbf{O}} \text{Gr}(p,k)$, which is of dimension $( p - k ).k$, the $p.k \times p.k$ matrix  $\nabla_{R}^{2} \psi ( \mathring{\mathbf{o}} )$ represents a linear mapping from $\emph{null} \big( \mathit{J}( \mathbf{r}(\mathbf{a} ) )  \big )^{\bot}$ to $\emph{null} \big( \mathit{J}( \mathbf{r}(\mathbf{a} ) )  \big )^{\bot}$ (or equivalently from $\emph{null} \big( \mathbf{M}(\mathbf{o} )  \big )^{\bot}$ to $\emph{null} \big( \mathbf{M}(\mathbf{o} )  \big )^{\bot}$) and we can represent this linear mapping in terms of the orthonormal basis of $\emph{null} \big( \mathit{J}( \mathbf{r}(\mathbf{a} ) )  \big )^{\bot}$ given by the columns of $\mathbf{\bar{O}}^{\bot} =  \mathbf{K}_{(p,k)} ( \mathbf{I}_{k}  \otimes  \mathbf{O}^{\bot}) \in \mathbb{O}^{k.p \times ( p - k ).k }$. When we do so the Riemannian Newton direction vector, $d\mathbf{\bar{o}}_{r-n}$, in $\emph{null} \big( \mathbf{M}(\mathbf{o} )  \big )^{\bot}$ can be computed as the solution of the following  $( p - k ).k \times ( p - k ).k$ symmetric linear system of equations
\begin{equation}  \label{eq:r_newton_sys}
(\mathbf{\bar{O}}^{\bot} )^{T}  \nabla_{R}^{2} \psi ( \mathring{\mathbf{o}} ) \mathbf{\bar{O}}^{\bot} d\mathbf{\bar{o}}_{r-n} = - ( \mathbf{\bar{O}}^{\bot} )^{T}   \nabla_{R} \psi ( \mathring{\mathbf{o}} ) \ ,
\end{equation}
which is exactly equivalent to the symmetric linear system given in equation~\ref{eq:newton_sys} since
\begin{align*}
(\mathbf{\bar{O}}^{\bot} )^{T}  \nabla_{R}^{2} \psi ( \mathring{\mathbf{o}} ) & = (\mathbf{\bar{O}}^{\bot} )^{T} \left ( \mathbf{I}_{p.k} -  \bar{\mathbf{O}} \bar{\mathbf{O}}^{T}   \right )  \nabla^{2} \psi ( \mathbf{o} )  \\
& = (\mathbf{\bar{O}}^{\bot} )^{T}  \mathbf{\bar{O}}^{\bot} (\mathbf{\bar{O}}^{\bot} )^{T} \nabla^{2} \psi ( \mathbf{o} ) \\
& = (\mathbf{\bar{O}}^{\bot} )^{T} \nabla^{2} \psi ( \mathbf{o} ) \ ,
\end{align*}
and
\begin{align*}
 - ( \mathbf{\bar{O}}^{\bot} )^{T}   \nabla_{R} \psi ( \mathring{\mathbf{o}} ) & = - ( \mathbf{\bar{O}}^{\bot} )^{T}   \nabla \psi ( \mathbf{o} ) \\
 & = -  ( \mathbf{\bar{O}}^{\bot}  )^{T}  \mathit{J}  \big( \mathbf{r}(\mathbf{o})  \big)^{T}  \mathbf{r}(\mathbf{o}) \\
 & =  \big( \mathbf{M}(\mathbf{o}) \mathbf{\bar{O}}^{\bot}  \big)^{T} \mathbf{r}(\mathbf{o}) \ .
\end{align*}
In a final step, we can also get  the Riemannian Newton direction vector, $d\mathbf{o}_{r-n}$, in $\mathbb{R}^{ k.p }$ as
\begin{equation*}
d\mathbf{o}_{r-n} = \mathbf{\bar{O}}^{\bot} d\mathbf{\bar{o}}_{r-n}  \ ,
\end{equation*}
and we have both $d\mathbf{o}_{r-n} = d\mathbf{o}_{n}$ and $d\mathbf{\bar{o}}_{r-n} = d\mathbf{\bar{o}}_{n}$. Thus, the Newton iteration defined by equation~\eqref{eq:newton_sys} can effectively be considered as a Riemannian Newton method if the next Newton iterate is computed with the help of a proper retraction onto the Stiefel manifold $\text{St}(p,k)$, as defined in equation~\eqref{eq:D_retraction} of Subsection~\ref{jacob:box}. Obviously, a similar argument shows that  the quasi-Newton iteration defined by equation~\eqref{eq:approx_newton_sys} can also be interpreted as a Riemannian quasi-Newton method.

We now show that the Newton iteration defined by equation~\eqref{eq:newton_sys2} does not share this nice property in general. Since $ \nabla_{R}^{2} \psi ( \mathring{\mathbf{o}} )$ is asymmetric and is considered to be a linear operator from $\mathcal{T}_{\mathbf{O}} \text{Gr}(p,k)$ to $\mathcal{T}_{\mathbf{O}} \text{Gr}(p,k)$ (more precisely from $\emph{null} \big( \mathit{J}( \mathbf{r}(\mathbf{a} ) )  \big )^{\bot}$ to $\emph{null} \big( \mathit{J}( \mathbf{r}(\mathbf{a} ) )  \big )^{\bot}$), it is first convenient to define the  $ p.k \times p.k$. symmetric projected Riemannian Hessian as 
\begin{align} \label{eq:sym_R_H_hess_mat}
& \nabla_{R}^{2} \psi ( \mathring{\mathbf{o}} )  \left ( \mathbf{I}_{p.k} -  \bar{\mathbf{O}} \bar{\mathbf{O}}^{T}   \right ) \nonumber \\
& =  \left ( \mathbf{I}_{p.k} -  \bar{\mathbf{O}} \bar{\mathbf{O}}^{T}   \right )   \big (  \mathbf{M}(\mathbf{o})^{T} \mathbf{M}(\mathbf{o}) - \mathbf{L}(\mathbf{o})^{T} \mathbf{L}(\mathbf{o}) +  \mathbf{U}(\mathbf{o})^{T} \mathbf{L} (\mathbf{o}) + \mathbf{L}(\mathbf{o})^{T} \mathbf{U}(\mathbf{o}) \big)   \left ( \mathbf{I}_{p.k} -  \bar{\mathbf{O}} \bar{\mathbf{O}}^{T}   \right ) \nonumber \\
  & =  \mathbf{M}(\mathbf{o})^{T} \mathbf{M}(\mathbf{o}) - \mathbf{L}(\mathbf{o})^{T} \mathbf{L}(\mathbf{o}) +    \left ( \mathbf{I}_{p.k} -  \bar{\mathbf{O}} \bar{\mathbf{O}}^{T}   \right )  \mathbf{U}(\mathbf{o})^{T} \mathbf{L} (\mathbf{o}) + \mathbf{L}(\mathbf{o})^{T} \mathbf{U}(\mathbf{o})  \left ( \mathbf{I}_{p.k} -  \bar{\mathbf{O}} \bar{\mathbf{O}}^{T}   \right ) \nonumber \\ 
   & =    \left ( \mathbf{I}_{p.k} -  \bar{\mathbf{O}} \bar{\mathbf{O}}^{T}   \right )   \nabla^{2} \psi ( \mathbf{o} )   \left ( \mathbf{I}_{p.k} -  \bar{\mathbf{O}} \bar{\mathbf{O}}^{T}   \right )  \ ,
\end{align}
which is now a symmetric linear mapping from $\mathbb{R}^{p.k}$ to $\mathbb{R}^{p.k}$, but this one has, however, the inconvenient to be always rank-deficient with its null space equals to $\emph{null}( \mathit{J}( \mathbf{r}(\mathbf{o}) ) ) = \emph{null}( \mathbf{M}(\mathbf{o}) )$ in regular cases, e.g., when the hypotheses of Corollary~\ref{corol5.5:box} are fulfilled and  the matrix function $\mathbf{F}(.)$ has full column rank in a neighborhood of $\mathbf{o}$ so that $\nabla^{2} \psi ( \mathbf{o} )$ exists.

Next, adding the term $\mathbf{N} \mathbf{N}^{T}$, defined in equation~\eqref{eq:NNt_mat}, to this  symmetric projected Riemannian Hessian matrix will remove its systematic rank degeneracy~\cite{HF2015b} and provides an alternative method to compute  the Riemannian Newton direction vector $d\mathbf{o}_{r-n} $, defined above, as the unique solution vector of the following  $ p.k \times p.k $ symmetric linear system in regular cases
\begin{equation*}
 \left ( \mathbf{I}_{p.k} -  \bar{\mathbf{O}} \bar{\mathbf{O}}^{T}   \right )   \nabla^{2} \psi ( \mathbf{o} )   \left ( \mathbf{I}_{p.k} -  \bar{\mathbf{O}} \bar{\mathbf{O}}^{T}   \right ) d\mathbf{o}_{r-n} = \mathbf{M}(\mathbf{o})^{T} \mathbf{r}(\mathbf{o})  \ .
\end{equation*}
However, we observe that Riemannian Newton iteration based on this damped version of the symmetric projected Riemannian Hessian will differ in general from the damped Newton iteration based on equation~\eqref{eq:newton_sys2} despite the damping term is the same. This is obvious from equation~\eqref{eq:sym_R_H_hess_mat} defining the symmetric projected Riemannian Hessian, which shows that the third symmetric term in the symmetric projected Riemannian Hessian differs from the third symmetric term in the formulation of the Euclidean Hessian given in equation~\eqref{eq:H_hess_mat}. On the other hand, as the first two-terms of the  symmetric projected Riemannian and Euclidean Hessians are identical, we can conclude that the damped quasi-Newton iteration based on a two-term approximation of the Euclidean Hessian $\bar{\mathbf{H}}$, defined in equation~\eqref{eq:approx_newton_sys2}, can still be interpreted as a Riemannian quasi-Newton method operating on the (quotient) Grassmann manifold $\text{Gr}(p,k)$ as claimed above.

To conclude that subsection, we now give an overview of the second-order trust-region method (RTRMC2) proposed by Boumal and Absil~\cite{BA2015} to minimize the cost function $\psi \circ  h^{-1}(.) = g_{\lambda}( \mathbf{.} ) $ on the Grassmann manifold  $\text{Gr}(p,k)$, where $g_{\lambda}( . )$ is defined in equation~\eqref{eq:g_func},  in order to contrast its advantages and drawbacks compared to the quasi-Newton methods just described above. To this end and as above, we assume that $\mathbf{W} \in  \mathbb{R}^{ p \times n  }_{+*}$, that each element $\mathring{\mathbf{O}}  \in \text{Gr}(p,k)$ is represented by an element of the compact Stiefel manifold $\mathbf{O} \in \text{St}(p,k)$ and that $\psi  \circ  h^{-1} (.)$ is a smooth function on $\text{St}(p,k)$ such that both the Riemannian gradient and Hessian of $\psi  \circ  h^{-1}(.)$ exist $\forall \mathbf{O} \in \text{St}(p,k)$. In this framework, the RTRMC2 method can be both interpreted as a variable projection method and a Riemannian second-order optimization method operating on the quotient Grassmann manifold $\text{Gr}(p,k)$.

First, Boumal and Absil have avoided the direct computation of the vectorized forms of the Euclidean or Riemannian Hessian of $\psi(.)$ as we did in equations~\eqref{eq:H_hess_mat} and~\eqref{eq:R_H_hess_mat2}, respectively. Instead, they have  derived only the directional derivative in the direction of  $\mathbf{D}  \in \mathcal{T}_{\mathbf{O}} \text{Gr}(p,k) = \mathcal{H}_{\mathbf{O}} \text{St}(p,k)$ of the Riemannian gradient of $\psi  \circ  h^{-1} (.)$ at $\mathring{\mathbf{O}} \in \text{Gr}(p,k)$ in a compact formulae,  see equation 27 in~\cite{BA2015}, which is relatively inexpensive and efficient compared to the evaluation of the full Hessians in equations~\eqref{eq:H_hess_mat} and~\eqref{eq:R_H_hess_mat2}. In other words and in our notations, this expression uses unvectorized variables and is essentially equivalent to equation~\eqref{eq:R_H_hess_mat} given above. Furthermore, this compact expression is sufficient for implementing efficiently an (inexact) iterative inner subsolver inside the RTRMC2 method, as we will describe now.

Besides, the RTRMC2 method generates a sequence of iterates $\mathbf{O}^{i} \in \text{St}(p,k)$ (more precisely of iterates $\mathring{\mathbf{O}}^{i} \in \text{Gr}(p,k)$) together with a sequence of trust-region radii $\delta^{i} > 0$. At each iteration $i$, a  subproblem solver computes a new step  $d\mathbf{O}^{i} \in   \mathcal{T}_{\mathbf{O}} \text{Gr}(p,k)  = \mathcal{H}_{\mathbf{O}} \text{St}(p,k)$ and the next iterate $\mathbf{O}^{i+1} \in \text{St}(p,k)$ is obtained by performing a retraction $\text{Retraction}_{\mathbf{O}^{i} } (  d\mathbf{O}^{i} )$, as defined in equation~\eqref {eq:D_retraction}, to go back to the correct manifold, e.g., $\text{St}(p,k)$ (or more precisely $\text{Gr}(p,k)$~\cite{BA2015}). The RTRMC2 method proceeds by computing $d\mathbf{O}^{i}$ via minimizing an approximate second-order Taylor expansion of the objective function  $\psi \circ  h^{-1}(.)$ within a "trust-region" $ \Vert d\mathbf{O}^{i}  \Vert_{F} \le \delta^{i}$ with adaptively chosen radii $\delta^{i}$~\cite{AMS2008}\cite{BA2015}\cite{B2023}. Depending on the performance of the inner subsolver, the outer RTRMC2  algorithm decides to accept or reject the proposed step $d\mathbf{O}^{i}$, and, possibly, decides to increase or reduce the trust-region radii. More precisely, in the inner subsolver, $d\mathbf{O}^{i}$ is computed to approximately minimize the following quadratic model function $m^{i} ( d\mathbf{O} ) :  \mathcal{T}_{\mathbf{O}} \text{Gr}(p,k)  \mapsto \mathbb{R}$ defined by
\begin{equation*}
m^{i} ( d\mathbf{O} ) = \psi  \circ  h^{-1} \left( \mathbf{O}^{i} \right) + \langle  d\mathbf{O}  , \nabla_{R} \psi  \circ  h^{-1} ( \mathring{\mathbf{O}} ) \rangle_{F} + \frac{1}{2}  \langle  d\mathbf{O} , \lbrack  \nabla_{R}^{2} \psi  \circ  h^{-1} ( \mathring{\mathbf{O}} ) \rbrack ( d\mathbf{O} ) \rangle_{F}  \ ,
\end{equation*} 
under the constraint $ \Vert d\mathbf{O} \Vert_{F} \le \delta^{i}$ and so that $m^{i} ( \mathbf{0}^{ p  \times k}  ) =  \psi  \circ  h^{-1} \left( \mathbf{O}^{i} \right)  $ and $m^{i} ( d\mathbf{O}^{i} ) \approx \psi  \circ  h^{-1} \left( \mathbf{O}^{i+1} \right)$. The selected inner subsolver is inexact and uses truncated conjugate gradient iterations to attempt to minimize $m^{i} ( . )$ over  $\mathcal{T}_{\mathbf{O}} \text{Gr}(p,k)$~\cite{BA2015}\cite{B2023} .

Despite the use of an iterative inner subsolver inside of an outer solver, the RTRMC2 method has a medium per-iteration complexity cost of the order of $\mathcal{O}( (p.k)^{2} )$, while the (quasi-)Newton methods given in equations~\eqref{eq:newton_sys}, ~\eqref{eq:approx_newton_sys},~\eqref{eq:newton_sys2}, ~\eqref{eq:approx_newton_sys2} have a much high per-iteration complexity of order at least $\mathcal{O}( (p.k)^{3} )$  for inverting the damped Hessian, giving a clear advantage to the RTRMC2 method in terms of speed and efficiency. On the other hand, the RTRMC2 method has more severe instability issues because its iterative inner solver based on (truncated) conjugate gradient iterations is not always robust when the Hessian is singular or ill-conditioned, which is always the case at first-order stationarity points and near the (local) non-isolated minima of $\psi  \circ  h^{-1} (.)$ where the Hessian has always vanishingly small, possibly negative eigenvalues~\cite{RB2024}. This is verified experimentally in Hong et al.~\cite{HF2015} and is explained by the theory developed in~\cite{RB2024}. Thus, with respect to accuracy and robustness, the Gauss-Newton, Levenberg-Marquardt and (quasi-)Newton methods developed in the previous and present subsections have a clear advantage as illustrated in the comparative studies of Okatani et al.~\cite{OYD2011} or Hong et al.~\cite{HF2015}.

This suggests finally that a fruitful area of future research to get the advantages of the two worlds, for solving efficiently and accurately difficult WLRA problems, may be to use a damped version of the Riemannian Hessian or of its two-term approximation (e.g., with a damping term equivalent to the one defined in equation~\eqref{eq:NNt_mat} for $\psi (.)$) in the quadratic model function $m^{i} (.)$ minimized by the inner solver of RTRMC2. This may eventually help to reduce its instability issues near the non-isolated minima of  $\psi \circ  h^{-1} (.)$, while keeping its lower per-iteration complexity.

\section{Implementation of variable projection NLLS methods for solving the WLRA problem} \label{vpalg:box}

This section is concerned with the formulation of practical and effective second-order algorithms for minimizing the cost function $\psi(.)$, i.e., the description of variable projection (pseudo-)second-order algorithms designed to solve the~\eqref{eq:VP1} formulation of the WLRA problem using the theoretical results established in the previous sections, especially the systematic rank-deficient nature of the Jacobian matrix $\mathit{J}( \mathbf{r}( \mathbf{a} ) )$ everywhere in the search space and the exactly rank-deficient nature of the Hessian matrix $\nabla^2 \psi( \mathbf{\widehat{a}} )$ if $\mathbf{\widehat{a}}$ is a stationary point of $\psi(.)$.

Standard results for the convergence of (quasi-)Newton methods suppose that the cost function is smooth, the target local minimum has a positive definite Hessian and that the algorithm is initialized in a neighborhood of this minimum~\cite{DS1983}\cite{NW2006}. However, this hypothesis is always violated here for the cost function $\psi(.)$ used in the~\eqref{eq:VP1} formulation of the WLRA problem as local minima of $\psi(.)$ are never isolated and can even form a continuum or a smooth manifold in some circumstances (e.g., when $\psi(.)$ locally verified the Morse-Bott property as illustrated in the previous section). In such deteriorated conditions, standard NLLS methods such as the Gauss-Newton or (quasi-)Newton methods would have difficulties~\cite{RB2024b}\cite{Z2024}. The Levenberg-Marquardt and trust-region Gauss-Newton methods can be used without modifications if the Jacobian $\mathit{J}( \mathbf{r}( \mathbf{a} ) )$ is not of full rank if the Marquardt damping parameter $\lambda$ or the radius $\Delta$ of the trust-region are controlled appropriately to not approach zero during the iterations and, especially, in the neighborhood of a critical point  $\mathbf{\widehat{a}}$. However, these methods may have very slow convergence if the Jacobian matrix $\mathit{J}( \mathbf{r}( \mathbf{a} ) )$ is singular everywhere in the search space and the damping factor is controlled so as not to tend to zero during the iterations, as in this scenario, we loose the quadratic or superlinear convergence speed of these methods near a critical point and the attainable accuracy is also limited~\cite{Z2024}. Moreover, most convergence results for all these methods essentially depend on the assumption that, at a solution point $\mathbf{\widehat{a}}$, the Jacobian $\mathit{J}( \mathbf{r}( \mathbf{\widehat{a}} ) )$ is nonsingular~\cite{DS1983}\cite{NW2006}.

However, some convergence results are also available for the exactly rank-deficient Jacobian case in the neighborhood of a solution point $\mathbf{\widehat{a}}$ for Gauss-Newton- or Levenberg-Marquardt-like methods~\cite{B1965}\cite{B1966}\cite{F1968}\cite{LH1974}\cite{B1976}\cite{H1986}\cite{Y1989}\cite{E1996}\cite{EW1996}\cite{YF2001}\cite{FY2005}\cite{EWGS2005}\cite{CKB2012}\cite{BM2015} or for a singular Hessian at  local minima $\mathbf{\widehat{a}}$ for  (quasi-)Newton methods~\cite{B1965}\cite{B1966}\cite{R1980}\cite{DK1980}\cite{GO1981}\cite{G1985}\cite{RR1986}\cite{NC1993}\cite{CNQ1997}\cite{ABHJ1999}\cite{DK2002}\cite{RB2024b}\cite{Z2024}\cite{DL2024}. Moreover, in practice, some of these (quasi-)second-order methods preserve their favourable faster convergence compared to first-order methods for non-isolated minima, a surprising result, which has been explained under mild assumptions like the Polyak–Lojasiewicz, Quadratic Growth and  Error Bound conditions~\cite{Z2024}\cite{RB2024}\cite{RB2024b}\cite{DL2024}. Recently, Rebjock and Boumal~\cite{RB2024b} unified these different results by demonstrated that  these three conditions  are essentially equivalent to the Morse-Bott property if the objective function is a $C^{2}$ function in a neighborhood of the target local minimum like the cost function $\psi(.)$ used in the~\eqref{eq:VP1} formulation of the WLRA problem under some circumstances (see the previous section for more details).

In the case of the Newton method, one of the most successful approaches in the case of a singular Hessian is the use of bordering techniques introduced by Griewank and co-workers~\cite{G1985}\cite{RR1986} and the adaptation of the Newton method described in the previous section to solve the~\eqref{eq:VP1} form of the WLRA problem falls in this category. On the other hand, if $\mathbf{a}_s$ is the current approximate solution, recall that the Gauss-Newton method for minimizing $\psi(.)$ is based on a linear approximation of the residual function $\mathbf{r}(.)$ from the Taylor's expansion around $\mathbf{a}_s$ 
\begin{equation*}
\mathbf{r}(  \mathbf{a} ) = \mathbf{r}(  \mathbf{a}_s ) + \mathit{J} \big ( \mathbf{r}(\mathbf{a}_s) \big)( \mathbf{a} - \mathbf{a}_s ) +  \mathcal{O}( \Vert  \mathbf{a} - \mathbf{a}_s  \Vert^{2}_{2} ) \ .
\end{equation*}
The Gauss-Newton step $d\mathbf{a}_{gn} = \mathbf{a} - \mathbf{a}_s$ is then the solution of the linear least-squares problem
\begin{equation*}
d\mathbf{a}_{gn} = \text{Arg}\min_{d\mathbf{a} \in \mathbb{R}^{p.k}}   \,   \frac{1}{2} \Vert \mathbf{r}(  \mathbf{a}_s ) + \mathit{J} \big( \mathbf{r}(\mathbf{a}_s) \big) d\mathbf{a} \Vert^{2}_{2} \ .
\end{equation*}
Since $\mathit{J}( \mathbf{r}( \mathbf{a}_s ) )$ is always rank deficient, the solution of this linear least-squares problem is not unique. However, a natural choice is to take $d\mathbf{a}_{gn}$ to be the unique minimum 2-norm solution, which is given by
\begin{equation*}
d\mathbf{a}_{gn} = - \mathit{J} \big( \mathbf{r}(\mathbf{a}_s) \big)^{+}  \mathbf{r}(\mathbf{a}_s) \ ,
\end{equation*}
since the linearization argument used to derive the Gauss-Newton iteration is only valid in a "small" neighborhood of $\mathbf{a}_s$. This leads to the generalized (Gauss-)Newton or Ben-Israel iterative method
\begin{equation*}
\mathbf{a}_{s+1} = \mathbf{a}_s +  d\mathbf{a}_{gn} =  \mathbf{a}_s - \mathit{J} \big( \mathbf{r}(\mathbf{a}_s) \big)^{+}  \mathbf{r}(\mathbf{a}_s) \ .
\end{equation*}
The local and global convergence properties of this algorithm when $\mathit{J}( \mathbf{r}( \mathbf{a} ) )$ does not have full rank in the neighborhood of a solution point $\mathbf{\widehat{a}}$ have been first investigated by Ben-Israel~\cite{B1965}\cite{B1966} and then by several authors afterward~\cite{B1976}\cite{H1986}\cite{NC1993}\cite{CNQ1997}\cite{DK2002}\cite{Z2024}. On assumptions like that $\mathit{J}( \mathbf{r}( \mathbf{a} ) )^{+}$ is Lipschitz-continuous, the Jacobian $\mathit{J}( \mathbf{r}( \mathbf{a} ) )$ is of constant rank in some neighborhood of $\mathbf{\widehat{a}}$ or the cost function $\psi(.)$ verifies a Morse-Bott-like property in a neighborhood of  $\mathbf{\widehat{a}}$ , they were able to show the convergence of this generalized (Gauss-)Newton method to a stationary point of $\psi(.)$. In practice, it is necessary to include some strategy to estimate the numerical rank of $\mathit{J}( \mathbf{r}( \mathbf{a}_s ) )$ and the assigned rank can have a decisive influence on the success of the method. Thus, a QR- or COD-factorization, or even a SVD-decomposition, of the matrix $\mathit{J}( \mathbf{r}( \mathbf{a}_s ) )$ are natural tools for the Ben-Israel iteration. This implies that the Ben-Israel method is relatively expensive. An alternative to this problem is to use a Tikhonov regularized version of the Ben-Israel iteration as in Levenberg-Marquardt or trust-region methods (see Subsection~\ref{opt:box}), but the choice of the regularization parameter is also a challenge in the rank-deficient case. Fortunately, the specific properties of the Jacobian $\mathit{J}( \mathbf{r}( \mathbf{a} ) )$ derived in Subsection~\ref{jacob:box} can be used for this purpose and also to reduce drastically the cost of the Ben-Israel iterative method or its regularized version as we will illustrate below.

Alternatively, Menzel~\cite{M1985} has proposed to reformulate the exactly rank deficient problem as an auxiliary least-squares problem of higher dimension which can be shown to be a well-posed one if the rank deficiency of $\mathit{J}( \mathbf{r}( \mathbf{a} ) )$ is small. Moreover, he was able to prove that for arbitrary rank deficiency in the consistent case $\psi( \mathbf{\widehat{ a}}  ) = 0$, the Gauss-Newton sequence for his auxiliary least-squares problem converges at least superlinearly to $\mathbf{\widehat{ a}}$. However, his technique, which expands considerably the size of the problem, is useful in practice only for small dimensions and small rank deficiency values and cannot be applied to our WLRA problem where both the dimensions and the rank deficiency values may be high.

More recently, Eriksson and Wedin~\cite{E1996}\cite{EW1996}\cite{EWGS2005} also considered the NLLS problem with an exactly rank-deficient Jacobian matrix for all points in a neighborhood of a local solution. As this situation corresponds exactly to the minimization of the cost function $\psi(.)$ in the~\eqref{eq:VP1} formulation of the WLRA problem, as demonstrated in the previous sections, we discuss now their method in some details. They suggested the following reformulation (in our notations) of the~\eqref{eq:VP1} variant of the WLRA problem in order to obtain a uniquely defined solution
 \begin{equation*}
  \mathbf{\widehat{a}} =
    \begin{cases}
        \text{Arg}\min_{\mathbf{a} \in \mathbb{R}^{p.k}} \,  \frac{1}{2}  \Vert \mathbf{a} \Vert^{2}_{2} \\
        \text{s.t. }  \text{Arg}\min_{\mathbf{a} \in \mathbb{R}^{p.k}} \,  \frac{1}{2} \psi( \mathbf{a} )  = \frac{1}{2}  \Vert \mathbf{r}( \mathbf{a} )  \Vert^{2}_{2}
    \end{cases} \ ,
\end{equation*}
\\
and they proposed two different iterative methods to solve this problem: a Gauss-Newton method and a Tikhonov regularization method. In words, what Eriksson and Wedin~\cite{E1996}\cite{EW1996}\cite{EWGS2005} suggested is to actually obtain the minimum Euclidean norm solution to the problem of minimizing $\psi(.)$. A similar approach was already mentioned by Boggs~\cite{B1976}. Using a Taylor series through two terms around $\mathbf{a}_s$, we get the linearized version of this problem and the following iterative method
 \begin{equation*}
  \mathbf{a}_{s+1} =
    \begin{cases}
        \text{Arg}\min_{\mathbf{a} \in \mathbb{R}^{p.k}} \,  \frac{1}{2}  \Vert \mathbf{a} \Vert^{2}_{2} \\
        \text{s.t. }  \text{Arg}\min_{\mathbf{a} \in \mathbb{R}^{p.k}} \,  \frac{1}{2}  \Vert \mathbf{r}( \mathbf{a}_{s} ) + \mathit{J} \big( \mathbf{r}( \mathbf{a}_s ) \big) (\mathbf{a} - \mathbf{a}_{s} ) \Vert^{2}_{2}
    \end{cases} \ .
\end{equation*}
\\
Following Pes and Rodriguez~\cite{PR2020}\cite{PR2022}, we will denote this as the Minimal-Norm Gauss-Newton (MNGN) method. It must be noted that the Ben-Israel iterative method has no predisposition toward the minimum 2-norm solution of minimizing $\psi(.)$ in the sense that any limit point generated by the Ben-Israel iteration is a least-squares solution, but not in general the minimum 2-norm solution and, consequently, the Ben-Israel and MNGN solutions will differ in general~\cite{B1976}\cite{EWGS2005}\cite{CKB2012}\cite{PR2020}. Since
\begin{equation*}
\mathbf{r}( \mathbf{a}_{s} ) + \mathit{J} \big( \mathbf{r}( \mathbf{a}_s ) \big) (\mathbf{a} - \mathbf{a}_{s} ) = \Big( \mathbf{r}( \mathbf{a}_{s} ) - \mathit{J} \big( \mathbf{r}( \mathbf{a}_s ) \big) \mathbf{a}_{s} \Big) + \mathit{J} \big( \mathbf{r}( \mathbf{a}_s ) \big) \mathbf{a} \ ,
\end{equation*}
the minimum 2-norm solution $\mathbf{a}_{s+1}$ of the problem
 \begin{equation*}
 \min_{\mathbf{a} \in \mathbb{R}^{p.k}} \,  \frac{1}{2}  \Vert \mathbf{r}( \mathbf{a}_{s} ) + \mathit{J} \big( \mathbf{r}( \mathbf{a}_s ) \big) (\mathbf{a} - \mathbf{a}_{s} ) \Vert^{2}_{2}
\end{equation*}
is
\begin{align*}
\mathbf{a}_{s+1} & =   - \mathit{J} \big( \mathbf{r}( \mathbf{a}_s ) \big)^{+} \Big( \mathbf{r}( \mathbf{a}_{s} ) - \mathit{J} \big( \mathbf{r}( \mathbf{a}_s ) \big) \mathbf{a}_{s} \Big)  \\
                            & =   - \mathit{J} \big( \mathbf{r}( \mathbf{a}_s ) \big)^{+} \mathbf{r}( \mathbf{a}_{s} ) +  \mathit{J} \big( \mathbf{r}( \mathbf{a}_s ) \big)^{+} \mathit{J} \big( \mathbf{r}( \mathbf{a}_s ) \big) \mathbf{a}_{s}  \\
                            & =   d\mathbf{a}_{gn} + \mathbf{P}_{\mathit{J} ( \mathbf{r}( \mathbf{a}_s ) )^{T}} \mathbf{a}_{s} \\
                            & =    \mathbf{P}_{\mathit{J} ( \mathbf{r}( \mathbf{a}_s ) )^{T}} \big(   \mathbf{a}_{s} + d\mathbf{a}_{gn}   \big) \ ,
\end{align*}
where $\mathbf{P}_{\mathit{J}( \mathbf{r}( \mathbf{a}_s ) )^{T}}$ is the orthogonal projector onto $\emph{ran}( \mathit{J}( \mathbf{r}( \mathbf{a}_s ) )^{T} ) = \emph{null}( \mathit{J}( \mathbf{r}( \mathbf{a}_s ) ) )^{\bot}$, e.g., onto the row space of $\mathit{J}( \mathbf{r}( \mathbf{a}_s ) )$ and the last equality results from equation~\eqref{eq:ginv}. Thus, to ensure computation of the minimal 2-norm solution, at the $s^{th}$ iteration of the MNGN method, the $standard$ Gauss-Newton approximation $\mathbf{a}_{s} + d\mathbf{a}_{gn}$, computed by the Ben-Israel method, is (orthogonally) projected onto the orthogonal of the null space of $\mathit{J}( \mathbf{r}( \mathbf{a}_s ) )$, e.g., $\emph{null}( \mathit{J}( \mathbf{r}( \mathbf{a}_s ) ) )^{\bot}$.
Alternatively, the MNGN iteration proposed by Eriksson and Wedin~\cite{EW1996}\cite{EWGS2005} can be written as
\begin{equation*}
\mathbf{a}_{s+1} = \mathbf{a}_{s} + d\mathbf{a}_{mngn}  \ ,
\end{equation*}
with the MNGN step computed as
\begin{align*}
d\mathbf{a}_{mngn}  & =  - \mathit{J} \big( \mathbf{r}( \mathbf{a}_s )  \big)^{+} \mathbf{r}( \mathbf{a}_{s} ) +  \mathbf{P}_{\mathit{J}( \mathbf{r}( \mathbf{a}_s ) )^{T}} \mathbf{a}_{s} - \mathbf{a}_s  \\
                                  & =  d\mathbf{a}_{gn} -  \mathbf{P}_{\mathit{J}( \mathbf{r}( \mathbf{a}_s ) )^{T}}^{\bot} \mathbf{a}_{s} \\
                                  & =  d\mathbf{a}_{gn} -  \mathbf{P}_{\emph{null}(\mathit{J}( \mathbf{r}( \mathbf{a}_s ) ) )} \mathbf{a}_{s}  \ ,
\end{align*}
and where $\mathbf{P}_{\emph{null}(\mathit{J}( \mathbf{r}( \mathbf{a}_s ) ) )}$ is the orthogonal projector onto the null space of $\mathit{J}( \mathbf{r}( \mathbf{a}_s ) )$. See Subsection~\ref{lin_alg:box} for details how the orthogonal projectors $\mathbf{P}_{\emph{null}(\mathit{J}( \mathbf{r}( \mathbf{a}_s ) ) )}$ and $\mathbf{P}_{\mathit{J}( \mathbf{r}( \mathbf{a}_s ) )^{T}}$ can be represented using the SVD of $\mathit{J}( \mathbf{r}( \mathbf{a}_s ) )$ or much more efficiently using a COD of this matrix. Furthermore, note that these orthogonal projectors  can be easily computed using the orthonormal bases of $\emph{null}(\mathit{J}( \mathbf{r}( \mathbf{a}_s ) ) )$ and its orthogonal complement identified in Corollary~\ref{corol5.6:box} (assuming that the rank of $\mathit{J}( \mathbf{r}( \mathbf{a}_s ) )$ is equal to $r = k.(p-k)$).

In order to ensure global convergence of the MNGN method, Campbell et al.~\cite{CKB2012} and Pes and Rodriguez~\cite{PR2022} have also considered the inclusion of relaxation (or damping) parameters  in the MNGN method, giving the iterative methods
\begin{equation} \label{eq:MNGN1} \tag{MNGN1}
\mathbf{a}_{s+1} = \mathbf{a}_{s} + \alpha_{s} d\mathbf{a}_{gn} - \mathbf{P}_{\emph{null}(\mathit{J}( \mathbf{r}( \mathbf{a}_s ) ) )} \mathbf{a}_{s}
\end{equation}
or
\begin{equation} \label{eq:MNGN2} \tag{MNGN2}
\mathbf{a}_{s+1} =  \mathbf{a}_{s} + d\mathbf{a}_{gn} - \beta_{s} \mathbf{P}_{\emph{null}(\mathit{J}( \mathbf{r}( \mathbf{a}_s ) ) )} \mathbf{a}_{s} 
\end{equation}
and also
\begin{equation}  \label{eq:MNGN3} \tag{MNGN3}
\mathbf{a}_{s+1} =  \mathbf{a}_{s} + \alpha_{s} d\mathbf{a}_{gn} - \beta_{s} \mathbf{P}_{\emph{null}(\mathit{J}( \mathbf{r}( \mathbf{a}_s ) ) )} \mathbf{a}_{s} \ ,
\end{equation}
where $\alpha_{s}$ and  $\beta_{s}$ are step length parameters, which can be determined by a line search and specific strategies described in~\cite{CKB2012}\cite{PR2022}.

As suggested by Eriksson and Wedin~\cite{E1996}\cite{EW1996}\cite{EWGS2005} and Pes and Rodriguez~\cite{PR2020}, a good approximate Minimal-Norm Levenberg-Marquardt step, $d\mathbf{a}_{mnlm}$, can also be found by solving the regularized linear least-squares problem
\begin{equation} \label{eq:MNLM_LLSQ}
d\mathbf{a}_{mnlm} = \text{Arg} \min_{d\mathbf{a} \in \mathbb{R}^{p.k}}   \,   \frac{1}{2} \Big\Vert \begin{bmatrix} \mathbf{r}(  \mathbf{a}_s )  \\  \mu \mathbf{a}_{s}   \end{bmatrix} + \begin{bmatrix} \mathit{J} \big( \mathbf{r}( \mathbf{a}_{s} )  \big)   \\   \mu \mathbf{I}_{k.p} \end{bmatrix} d\mathbf{a} \Big\Vert^{2}_{2} \ ,
\end{equation}
with a sufficiently small Tikhonov regularization parameter $\mu$ since this regularized problem has an unique solution
\begin{align*}
d\mathbf{a}_{mnlm} = & -  \left(    \mathit{J}  \big( \mathbf{r}( \mathbf{a}_{s} )  \big)^{T}  \mathit{J} \big( \mathbf{r}( \mathbf{a}_{s} )  \big) +  \mu^{2}  \mathbf{I}_{k.p}   \right)^{-1} \mathit{J}  \big( \mathbf{r}( \mathbf{a}_{s} )  \big)^{T} \mathbf{r}( \mathbf{a}_{s} )  \\
                                    & -  \left(    \mathit{J}  \big( \mathbf{r}( \mathbf{a}_{s} )  \big)^{T}  \mathit{J} \big( \mathbf{r}( \mathbf{a}_{s} )  \big) +  \mu^{2}  \mathbf{I}_{k.p}   \right)^{-1}  \mu^{2} \mathbf{a}_{s}
\end{align*}
and
\begin{equation*}
\begin{split}
        & \lim_{ \mu  \to 0} \, \left(    \mathit{J} \big( \mathbf{r}( \mathbf{a}_{s} ) \big)^{T}  \mathit{J} \big( \mathbf{r}( \mathbf{a}_{s} ) \big) +  \mu^{2}  \mathbf{I}_{k.p}   \right)^{-1} \mathit{J} \big( \mathbf{r}( \mathbf{a}_{s} ) \big)^{T} =  \mathit{J} \big( \mathbf{r}( \mathbf{a}_{s} ) \big)^{+} \ , \\
       & \lim_{ \mu  \to 0}  \,   \left(  \mathit{J} \big( \mathbf{r}( \mathbf{a}_{s} ) \big)^{T}  \mathit{J} \big( \mathbf{r}( \mathbf{a}_{s} ) \big) +  \mu^{2}  \mathbf{I}_{k.p}   \right)^{-1}  \mu^{2} =  \mathbf{P}_{\emph{null}(\mathit{J}( \mathbf{r}( \mathbf{a}_s ) ) )} \ ,
\end{split}
\end{equation*}
as demonstrated in~\cite{GW2000}. The Tikhonov method proposed by Eriksson and Wedin~\cite{EWGS2005}, and also considered by Pes and Rodriguez~\cite{PR2020}, then approximately solves a sequence of Tikhonov regularized NLLS problems for a sequence of decreasing regularization parameter $\mu_t$ where the index $t$ does not necessarily equal the iteration index $s$. The approximate solution of one regularized NLLS problem with Tikhonov parameter $\mu_t$ is taken as the starting point for the next regularized problem with Tikhonov parameter $\mu_{t+1}<\mu_t$. Eriksson and Wedin~\cite{EWGS2005} proved that both the MNGN and Tikhonov methods converge to a local minimum 2-norm solution if the iterations are started in a neighborhood of the solution.

These two methods are directly applicable to an exactly rank deficient problem such as the~\eqref{eq:VP1} formulation of the WLRA problem. However, we note that the different damped variants of the MNGN method may fail to converge in many cases because projecting the GN solution orthogonally to the null space of $\mathit{J}( \mathbf{r}( \mathbf{a}_s ) )$  may cause the residual to increase during the iterations for the MNGN1 variant or because the MNGN2 variant is equivalent to the application of the undamped Gauss-Newton method, whose convergence is not theoretically guaranteed~\cite{PR2022}\cite{CKB2012}. Note further that in the case of the MNGN3 variant, convergence can be generally obtained if $\alpha_{s}$ and $\beta_{s}$ are suitably chosen, but in that case the method does not converge to the minimum 2-norm solution unless $\beta_{s}=1$ for $s$ close to convergence~\cite{PR2022}. Moreover, it is not of interest, nor natural to find the minimum 2-norm least-squares solution of our WLRA problem due to its very special structure, separability properties and over-parameterization (e.g., because it is an optimization problem on manifolds or subspaces). Additionally, the MNGN and Tikhonov methods involve the extra-computations of the terms
\begin{equation*}
\mathbf{P}_{\emph{null}(\mathit{J}( \mathbf{r}( \mathbf{a}_s ) ) )} \mathbf{a}_{s} \quad \text{and} \quad \left(    \mathit{J} \big(( \mathbf{r}( \mathbf{a}_{s} ) \big()^{T}  \mathit{J} \big(( \mathbf{r}( \mathbf{a}_{s} ) \big() +  \mu^{2}  \mathbf{I}_{k.p}   \right)^{-1}  \mu^{2} \mathbf{a}_{s} \ ,
\end{equation*}
in each iteration, respectively. This suggests that simple Ben-Israel or regularized Tikhonov methods will be more appropriate to minimize $\psi(.)$.

With these considerations, we are now in position to give a full formal description of different variations of the Gauss-Newton, Levenberg-Marquardt and (quasi-)Newton algorithms which may be used to minimize $\psi(.)$ in practice.

\subsection{Variable projection Gauss-Newton algorithms} \label{vp_gn_alg:box}

Using similar notations and definitions as in previous sections, an outline of the variable projection Gauss-Newton algorithms is as follows:
\\
\begin{gn_alg} \label{gn_alg:box}
\end{gn_alg}
Choose starting matrix $\mathbf{A}_{1} \in \mathbb{R}^{p \times k}$, $\varepsilon_{1}, \varepsilon_{2}, \varepsilon_{3}  \in \mathbb{R}_{+*}$ and $i_{max}, j_{max} \in \mathbb{N}_{*}$,  appropriately

\textbf{For} $i=1, 2, \ldots$ \textbf{until convergence do}
\begin{enumerate}
\item[\textbf{(0)}]  Optionally, compute a QRCP of $\mathbf{A}_{i}$ (see equation~\eqref{eq:qrcp}) to determine $k_{i} = \emph{rank}( \mathbf{A}_{i} )$ and an orthonormal basis of $\emph{ran}( \mathbf{A}_{i} )$:
\begin{itemize}
\item[] $\mathbf{Q}_{i} \mathbf{A}_{i} \mathbf{P}_{i} =  \begin{bmatrix} \mathbf{R}_{i}    &  \mathbf{S}_{i}   \\  \mathbf{0}^{(p-k_{i}) \times k_{i} }  & \mathbf{0}^{(p-k_{i}) \times (k-k_{i}) }  \end{bmatrix}$ ,
\end{itemize}
where $\mathbf{Q}_{i}$ is an $p \times p$ orthogonal matrix, $\mathbf{P}_{i}$ is an $k \times k$ permutation matrix, $\mathbf{R}_{i}$ is an $k_{i} \times k_{i}$ nonsingular upper triangular matrix (with diagonal elements of decreasing absolute magnitude) and $\mathbf{S}_{i}$ an $k_{i} \times (k-k_{i})$ full matrix, which is vacuous if $k_{i} = k$. If $k_{i} < k$, complete the orthonormal basis of $\emph{ran}( \mathbf{A}_{i} )$ with $k - k_{i}$ orthonormal vectors by using the $p \times p$ orthogonal matrix $\mathbf{Q}_{i}$  computed implicitly during the QRCP of $\mathbf{A}_{i}$.

In other words, in all cases, compute a $p \times k$  matrix $\mathbf{O}_{i}$  with orthonormal columns as the first $k$ columns of $\mathbf{Q}_{i}$ (i.e., such that $\emph{ran}( \mathbf{A}_{i} ) \subset \emph{ran}( \mathbf{O}_{i} )$ if $k_{i} < k$ and  $\emph{ran}( \mathbf{A}_{i} ) = \emph{ran}( \mathbf{O}_{i} )$ if $k_{i} = k$) and set
\begin{itemize}
\item[] $\mathbf{A}_{i} = \mathbf{O}_{i}$ .
\end{itemize}
This optional orthogonalization step is a safe-guard as the condition $k_{i} = k$ is a necessary condition for the differentiability of $\psi(.)$ at a point $\mathbf{A}_{i}$ and is also useful to limit the occurrence of overflows and underflows in the next steps by enforcing that the matrix variable $\mathbf{A}_{i} \in \mathbb{O}^{p \times k}$.
\item[\textbf{(1)}]  Determine (implicitly) the block diagonal matrix
\begin{itemize}
\item[] $\mathbf{F}(\mathbf{a}_{i}) = \emph{diag}\big( \emph{vec}( \sqrt{\mathbf{W}} ) \big)  \big(  \mathbf{I}_n  \otimes \mathbf{A}_{i}  \big)$ ,
\end{itemize}
where $\mathbf{a}_{i} =  \emph{vec}( \mathbf{A}_{i}^{T} )$ .
\item[\textbf{(2)}]  Compute (implicitly) a QRCP of $\mathbf{F}(\mathbf{a}_{i})$ to determine $\mathbf{P}^{\bot}_{\mathbf{F}(  \mathbf{a}_{i} ) }$ and  $\mathbf{F}(\mathbf{a}_{i})^-$ (see equations~\eqref{eq:ginv_proj_ortho} and \eqref{eq:sginv_qrcp}) or, alternatively, a COD of  $\mathbf{F}(\mathbf{a}_{i})$ to determine $\mathbf{P}^{\bot}_{\mathbf{F}(  \mathbf{a}_{i} ) }$ and $\mathbf{F}(\mathbf{a}_{i})^{+}$ (see equations~\eqref{eq:ginv_proj_ortho} and \eqref{eq:ginv_cod}).

Note also that $\mathbf{F}(\mathbf{a}_{i})^- = \mathbf{F}(\mathbf{a}_{i})^{+}$ if $\mathbf{F}(\mathbf{a}_{i})$ is of full column rank and that $\mathbf{P}^{\bot}_{\mathbf{F}(  \mathbf{a}_{i} ) }$, $\mathbf{F}(\mathbf{a}_{i})^-$ and $\mathbf{F}(\mathbf{a}_{i})^{+}$ are also block diagonal matrices.
\item[\textbf{(3)}]  Solve the block diagonal linear least-squares problem
\begin{itemize}
\item[] $\mathbf{b}_{i} = \text{Arg}\min_{\mathbf{b}\in\mathbb{R}^{k.n}} \,   \Vert \mathbf{x} - \mathbf{F}(\mathbf{a}_{i})\mathbf{b} \Vert^{2}_{2}$ ,
\end{itemize}
e.g., compute
\begin{itemize}
\item[] $\mathbf{b}_{i} =
    \begin{cases}
         \mathbf{F}(\mathbf{a}_{i})^{-} \mathbf{x}   \quad\ \lbrace \text{if a QRCP of $\mathbf{F}(\mathbf{a}_{i})$ is used in step } \textbf{(2)}  \rbrace \\
         \mathbf{F}(\mathbf{a}_{i})^{+} \mathbf{x}  \quad\ \lbrace \text{if a COD of $\mathbf{F}(\mathbf{a}_{i})$ is used in step } \textbf{(2)} \rbrace
    \end{cases}$ .
\end{itemize}
\item[\textbf{(4)}] Determine:
\begin{itemize}
\item[] $\mathbf{r}(\mathbf{a}_{i}) = \mathbf{P}^{\bot}_{\mathbf{F}(  \mathbf{a}_{i} ) } \mathbf{x}$ $\lbrace \text{current residual vector} \rbrace$
\item[] $\psi(\mathbf{a}_{i} ) = \frac{1}{2} \Vert \mathbf{r}(\mathbf{a}_{i}) \Vert^{2}_{2}$ $\lbrace \text{current value of the cost function }  \rbrace$
\item[] $\nabla \psi( \mathbf{a}_{i} ) =  \mathbf{G}(\mathbf{b}_{i})^{T} \mathbf{G}(\mathbf{b}_{i})\mathbf{a}_{i} - \mathbf{G}(\mathbf{b}_{i})^{T} \mathbf{z}$ $\lbrace$see Theorems~\ref{theo4.3:box} and~\ref{theo5.7:box}$\rbrace$
\end{itemize}
Note that the steps \textbf{(1)} to \textbf{(4)} above can be very easily parallelized using the block diagonal structure of $\mathbf{F}(\mathbf{a}_{i})$.
\item[\textbf{(5)}]   Check for convergence. Relevant convergence criteria in the algorithms are of the form:
\begin{itemize}
\item $\Vert \nabla \psi(\mathbf{a}_{i} ) \Vert_{2} \le \varepsilon_{1}$
\item $\Vert \mathbf{a}_{i} - \mathbf{a}_{i-1} \Vert_{2} \le \varepsilon_{2} ( \varepsilon_{2} + \Vert \mathbf{a}_{i} \Vert_{2})$  $\lbrace$if $ i \ne 1\rbrace$

If step \textbf{(0)} is used, this last convergence condition can be simplified as:

$\Vert \mathbf{a}_{i} - \mathbf{a}_{i-1} \Vert_{2} \le \varepsilon_{2}  \Vert \mathbf{a}_{i} \Vert_{2} = \varepsilon_{2} \sqrt{k}$
\item $\vert \psi(\mathbf{a}_{i-1} ) - \psi(\mathbf{a}_{i} ) \vert \le \varepsilon_{3} ( \varepsilon_{3} +  \psi(\mathbf{a}_{i} ) )$  $\lbrace$if $ i \ne 1\rbrace$
\item $ i \ge i_{max}$ $\lbrace \text{e.g., give up if the number of iterations is too large} \rbrace$
\end{itemize}
where $\varepsilon_{1}, \varepsilon_{2}, \varepsilon_{3}$ and $i_{max}$ are constants chosen by the user.

\textbf{Exit if convergence}. \textbf{Otherwise, go to} step \textbf{(6)}
\item[\textbf{(6)}]  Compute the Gauss-Newton correction vector $d\mathbf{a}_{gn}$ as the minimum 2-norm solution of one of the following linear least-squares problems:
\begin{description}
\item[Golub-Pereyra step:] Golub and Pereyra~\cite{GP1973}, Ruhe and Wedin~\cite{RW1980}
\begin{align*}
d\mathbf{a}_{gp-gn} & =   \big( \mathbf{M}( \mathbf{a}_{i}  ) + \mathbf{L}( \mathbf{a}_{i}  ) \big)^{+} \mathbf{r}( \mathbf{a}_{i} )   \\
                                 & = 
    \begin{cases}
        \text{Arg}\min_{d\mathbf{a} \in \mathbb{R}^{p.k}} \,   \Vert d\mathbf{a} \Vert^{2}_{2} \\
        \text{s.t. }  \text{Arg}\min_{d\mathbf{a} \in \mathbb{R}^{p.k}} \,   \Vert \mathbf{r}( \mathbf{a}_{i} ) - \big( \mathbf{M}( \mathbf{a}_{i}  ) + \mathbf{L}( \mathbf{a}_{i}  )  \big) d\mathbf{a}   \Vert^{2}_{2}
    \end{cases}   \\
\end{align*}
\item[Kaufman step:] Kaufman~\cite{K1975}, Ruhe and Wedin~\cite{RW1980}
\begin{align*}
d\mathbf{a}_{k-gn} & =   \mathbf{M}( \mathbf{a}_{i}  )^{+} \mathbf{r}( \mathbf{a}_{i} )   \\
                               & = 
    \begin{cases}
        \text{Arg}\min_{d\mathbf{a} \in \mathbb{R}^{p.k}} \,   \Vert d\mathbf{a} \Vert^{2}_{2} \\
        \text{s.t. }  \text{Arg}\min_{d\mathbf{a} \in \mathbb{R}^{p.k}} \,   \Vert \mathbf{r}( \mathbf{a}_{i} ) - \mathbf{M}( \mathbf{a}_{i}  ) d\mathbf{a}   \Vert^{2}_{2}
    \end{cases}   \\
\end{align*}
\item[Gauss-Seidel step:] Ruhe and Wedin~\cite{RW1980}
\begin{align*}
d\mathbf{a}_{gs-gn} & =   \big( \mathbf{K}_{(n,p)}  \mathbf{G}( \mathbf{b}_{i}  ) \big)^{+} \mathbf{r}( \mathbf{a}_{i} )   \\
                                 & = 
    \begin{cases}
        \text{Arg}\min_{d\mathbf{a} \in \mathbb{R}^{p.k}} \,   \Vert d\mathbf{a} \Vert^{2}_{2} \\
        \text{s.t. }  \text{Arg}\min_{d\mathbf{a} \in \mathbb{R}^{p.k}} \,   \Vert \mathbf{r}( \mathbf{a}_{i} ) - \mathbf{K}_{(n,p)}  \mathbf{G}( \mathbf{b}_{i}  ) d\mathbf{a}   \Vert^{2}_{2}
    \end{cases}   \\
\end{align*}
\end{description}
\item[\textbf{(7)}]  Increment $\mathbf{a}_{i} = \emph{vec}( \mathbf{A}_{i}^{T} )$, e.g., compute $\mathbf{a}_{i+1} = \emph{vec}( \mathbf{A}_{i+1}^{T} )$ such that $\psi( \mathbf{a}_{i+1} ) < \psi( \mathbf{a}_{i} )$ in order to obtain global convergence.
\begin{enumerate}
\item[\textbf{(7.1)}]  To this end, first compute
\begin{itemize}
\item[]  $\mathbf{a}_{i+1} = \mathbf{a}_{i} +  d\mathbf{a}_{gn}$
\end{itemize}
and if a Golub-Pereyra or Kaufman step is used in step  \textbf{(6)}
\begin{itemize}
\item[] $\psi(\mathbf{a}_{i+1} ) = \frac{1}{2} \Vert \mathbf{r}(\mathbf{a}_{i+1}) \Vert^{2}_{2} = \frac{1}{2} \Vert \mathbf{P}^{\bot}_{\mathbf{F}(  \mathbf{a}_{i+1} ) } \mathbf{x} \Vert^{2}_{2}$ ,
\end{itemize}
using (implicitly) a QRCP of the block diagonal matrix $\mathbf{F}(  \mathbf{a}_{i+1} )$.
\item[\textbf{(7.2)}]  \textbf{If} a Golub-Pereyra or Kaufman step is used in step \textbf{(6)} and $\psi( \mathbf{a}_{i+1} ) > \psi( \mathbf{a}_{i} )$  \textbf{then} recompute $\mathbf{a}_{i+1}$ by one of the following methods:
\begin{description}
\item[Gauss-Seidel:]
$\mathbf{a}_{i+1} = \mathbf{a}_{i} +  d\mathbf{a}_{gs-gn}$ where $d\mathbf{a}_{gs-gn}$ is a Gauss-Seidel step~\cite{RW1980} defined as
\begin{align*}
d\mathbf{a}_{gs-gn} & =   \big( \mathbf{K}_{(n,p)}  \mathbf{G}( \mathbf{b}_{i}  ) \big)^{+} \mathbf{r}( \mathbf{a}_{i} )   \\
                                 & =
    \begin{cases}
        \text{Arg}\min_{d\mathbf{a} \in \mathbb{R}^{p.k}} \,   \Vert d\mathbf{a} \Vert^{2}_{2} \\
        \text{s.t. }  \text{Arg}\min_{d\mathbf{a} \in \mathbb{R}^{p.k}} \,   \Vert \mathbf{r}( \mathbf{a}_{i} ) - \mathbf{K}_{(n,p)}  \mathbf{G}( \mathbf{b}_{i}  ) d\mathbf{a}   \Vert^{2}_{2}
    \end{cases}   \\
\end{align*}
\item[Block alternating least-squares:]
\begin{align*}
  \mathbf{a}_{i+1} & =    \mathbf{G}( \mathbf{b}_{i}  )^{+} \mathbf{z}   \\
                             & = 
    \begin{cases}
        \text{Arg}\min_{\mathbf{a} \in \mathbb{R}^{p.k}} \,   \Vert \mathbf{a} \Vert^{2}_{2} \\
        \text{s.t. }  \text{Arg}\min_{\mathbf{a} \in \mathbb{R}^{p.k}} \,   \Vert \mathbf{z} -  \mathbf{G}( \mathbf{b}_{i}  ) \mathbf{a}   \Vert^{2}_{2}
    \end{cases}   \\
\end{align*}
\item[Line search:]
\begin{equation*}
\mathbf{a}_{i+1} = \mathbf{a}_{i} + \alpha_{i}  d\mathbf{a}_{gn} \ ,
\end{equation*}
where $\alpha_i < 1$ is determined by a line search to make the algorithm a descent method (i.e., such that $\psi( \mathbf{a}_{i+1} ) <  \psi( \mathbf{a}_{i} ) $). This is always possible as the correction vector $d\mathbf{a}_{gn}$ is in a descent direction for $\psi(.)$ if $\Vert \nabla \psi(\mathbf{a}_{i} ) \Vert_{2} \ne 0$, see Corollary~\ref{corol5.7:box}.

As an illustration, a simple, but still efficient, strategy is to first shorten the correction step to half the Gauss-Newton length, compute the new trial value for $\psi( \mathbf{a}_{i+1} )$ and, if it is still worse, continue to reduce the step until we get a step short enough such that $\psi( \mathbf{a}_{i+1} ) <  \psi( \mathbf{a}_{i} )$. The following loop incorporates this simple step-shortening algorithm:

\textbf{For} $j=1, 2, \ldots$ \textbf{while}$\big( \psi( \mathbf{a}_{i+1} ) > \psi( \mathbf{a}_{i} ) \big)$
\begin{enumerate}
\item[]  $d\mathbf{a}_{gn}  = \frac{1}{2}  d\mathbf{a}_{gn}$
\item[]  $\mathbf{a}_{i+1} = \mathbf{a}_{i} +  d\mathbf{a}_{gn}$
\item[]  $\psi(\mathbf{a}_{i+1} ) =  \frac{1}{2} \Vert \mathbf{P}^{\bot}_{\mathbf{F}(  \mathbf{a}_{i+1} ) } \mathbf{x} \Vert^{2}_{2}$ $\lbrace$using a QRCP of the matrix $\mathbf{F}(  \mathbf{a}_{i+1} ) \rbrace$
\item[]  \textbf{If} $j > j_{max}$ \textbf{exit} $\lbrace \text{e.g., give up if the number of iterations is too large} \rbrace$
\end{enumerate}
\textbf{End do}

\end{description}
\end{enumerate}
\end{enumerate}
\textbf{End do}

For the convenience of the reader, we first recall the shape and definition of the  vector and matrix variables used in these Gauss-Newton algorithms. We have:
$\mathbf{X}\in\mathbb{R}^{p \times n}$, $\mathbf{W}\in\mathbb{R}^{p \times n}_+$, $\mathbf{A}_{i}\in\mathbb{R}^{p \times k}$, $\mathbf{B}_{i}\in\mathbb{R}^{k \times n}$,  $\mathbf{O}_{i}\in\mathbb{O}^{p \times k}$ and
\begin{align*}
\mathbf{x} & =  \emph{vec}( \sqrt{\mathbf{W}} \odot \mathbf{X} ) \ , \\
\mathbf{z} & =  \emph{vec} \big( (\sqrt{\mathbf{W}}  \odot \mathbf{X})^{T} \big) \ , \\
\mathbf{a}_{i} & =   \emph{vec}( \mathbf{A}_{i}^{T} ) \ , \\
\mathbf{F}(\mathbf{a}_{i}) & =  \emph{diag}\big( \emph{vec}( \sqrt{\mathbf{W}} ) \big)  \big(  \mathbf{I}_n  \otimes \mathbf{A}_{i}  \big) \ , \\
\mathbf{b}_{i} & =   \emph{vec}( \mathbf{B}_{i} ) =
    \begin{cases}
         \mathbf{F}(\mathbf{a}_{i})^{-} \mathbf{x}   \quad\  \lbrace \text{if a QRCP is used in step } \textbf{(2)}  \rbrace  \\
         \mathbf{F}(\mathbf{a}_{i})^{+} \mathbf{x}  \quad\  \lbrace \text{if a COD is used in step } \textbf{(2)} \rbrace 
    \end{cases} \ ,
 \\
\mathbf{G}(\mathbf{b}_{i}) & =  \emph{diag}\big( \emph{vec}( \sqrt{\mathbf{W}}^{T} ) \big)  \big(  \mathbf{I}_p  \otimes \mathbf{B}_{i}^{T}  \big) \ , \\
\mathbf{U}(\mathbf{a}_{i})    & = \emph{diag}( \emph{vec}(\sqrt{\mathbf{W}} ) ) ( \mathbf{B}_{i}^{T} \otimes \mathbf{I}_{p}) \mathbf{K}_{(k,p)} \ ,  \\
\mathbf{V}(\mathbf{a}_{i})    & =  \big( ( \mathbf{W} \odot P_{\Omega}( \mathbf{X} -\mathbf{A}_{i}\mathbf{B}_{i} ) )^{T} \otimes \mathbf{I}_{k} \big)  \lbrace \text{where } P_{\Omega}(.) \text{ is defined in equation~\eqref{eq:proj_op}} \rbrace \ , \\
\mathbf{M}(\mathbf{a}_{i}) & =  \mathbf{P}_{\mathbf{F}(\mathbf{a}_{i})}^{\bot} \mathbf{U}(\mathbf{a}_{i}) = \mathbf{P}_{\mathbf{F}(\mathbf{a}_{i})}^{\bot} \mathbf{K}_{(n,p)} \mathbf{G}(\mathbf{b}_{i}) \ ,  \\
\mathbf{L}(\mathbf{a}_{i}) & = 
    \begin{cases}
        \mathbf{F}(\mathbf{a}_{i})^{-T} \mathbf{V}(\mathbf{a}_{i}) \quad\ \lbrace \text{if a QRCP is used in step } \textbf{(2)} \rbrace  \\
        \mathbf{F}(\mathbf{a}_{i})^{+T} \mathbf{V}(\mathbf{a}_{i})  \quad\ \lbrace \text{if a COD is used in step } \textbf{(2)} \rbrace
    \end{cases} \ .
\end{align*}
In these Gauss-Newton algorithms, the iterations are terminated either when one or several of the convergence criteria listed in step \textbf{(5)} are satisfied, or when the iteration count exceeds the predetermined number $i_{max}$.

The Golub-Pereyra step $d\mathbf{a}_{gp-gn}$ corresponds exactly to the standard Gauss-Newton step $d\mathbf{a}_{gn}$  applied to the minimization of the variable projection functional  $\psi(.)$, which is introduced in Subsection~\ref{opt:box}.

The following is a brief review of the basic ideas underlying the simplification introduced by the Kaufman step $d\mathbf{a}_{k-gn}$ in the Gauss-Newton algorithms~\eqref{gn_alg:box} and also in the Marquardt-Levenberg algorithms~\eqref{lm_alg1:box}  described in the next subsection. As stated in Subsection~\ref{opt:box}, the Gauss-Newton method may be interpreted as a variation of Newton's method to find a zero of the gradient of  $\psi(.)$. More precisely, dropping the iteration index of the algorithm in order to simplify the notation, the Gauss-Newton algorithm approximates the Hessian matrix $\mathbf{H} = \nabla^2 \psi( \mathbf{a} )$ with the cross-product matrix $\mathit{J}( \mathbf{r}(\mathbf{a}) )^{T} \mathit{J}( \mathbf{r}(\mathbf{a}) )$ at each iteration based on the assumption that the components of the residual vector $\vert \mathbf{r}_{l} (\mathbf{a}) \vert$ are small near a solution and using the fact that the second term of the Hessian matrix $\nabla^2 \psi( \mathbf{a} )$ given by
\begin{equation*}
\mathbf{S} =  \sum_{l=1}^{n.p} \mathbf{r}_l (\mathbf{a}) \nabla^2 \mathbf{r}_l (\mathbf{a})
\end{equation*}
is of order $\mathcal{O}( \Vert \mathbf{r}(\mathbf{a})  \Vert_{2} )$ (see equation~\eqref{eq:S_hess_mat2} for the explicit form of $\mathbf{S}$ in our WLRA context). Now, the Gauss-Newton approximation of the Hessian matrix is
\begin{equation*}
\nabla^2 \psi( \mathbf{a} ) \approx \mathit{J}( \mathbf{r}(\mathbf{a}) )^{T} \mathit{J}( \mathbf{r}(\mathbf{a})  = \mathbf{M}(\mathbf{a})^{T} \mathbf{M}(\mathbf{a}) + \mathbf{L}(\mathbf{a})^{T} \mathbf{L}(\mathbf{a}) \ ,
\end{equation*}
as demonstrated in Subsection~\ref{hess:box} (see equation~\eqref{eq:gn_hess_mat} for details). Hence, the term $-\mathbf{L}(\mathbf{a})$ does not contribute to the gradient $\nabla \psi( \mathbf{a} )$ (see Theorem~\ref{theo5.7:box}) and changes the Gauss-Newton approximation of the Hessian only by the term $\mathbf{L}(\mathbf{a})^{T} \mathbf{L}(\mathbf{a})$, which is of order $\mathcal{O}( \Vert \mathbf{r}(\mathbf{a})  \Vert^{2}_{2} )$, see equation~\eqref{eq:L_mat}. If $\Vert \mathbf{r}(\mathbf{a})  \Vert_{2}$ is small then this term is smaller than the second term of $\nabla^2 \psi( \mathbf{a} )$, $\mathbf{S}$, which is of order $\mathcal{O}( \Vert \mathbf{r}(\mathbf{a})  \Vert_{2} )$, and which is already dropped in the Gauss-Newton and Levenberg-Marquardt methods. Following the Gauss-Newton philosophy, it is then natural to drop this term here too. In addition, since the approximation of the Hessian matrix by its first two symmetric exact terms is given by
\begin{equation*}
   \nabla^2 \psi( \mathbf{a} ) \approx   \mathbf{M}(\mathbf{a})^{T} \mathbf{M}(\mathbf{a}) - \mathbf{L}(\mathbf{a})^{T} \mathbf{L}(\mathbf{a}) \  ,
\end{equation*}
see equation~\eqref{eq:approx_hess_mat}, we may expect that approximating the Hessian matrix by the cross-product matrix $\mathbf{M}(\mathbf{a})^{T} \mathbf{M}(\mathbf{a})$ (or equivalently the Jacobian matrix by $-\mathbf{M}(\mathbf{a})$) can perform even better than the Gauss-Newton approximation of the Hessian matrix, see Subsection~\ref{hess:box} and equation~\eqref{eq:H_hess_mat} for details. This leads to the following linear least-squares problems for computing the simplified Kaufman correction step in the Gauss-Newton and Levenberg-Marquardt methods, respectively,
\begin{align*}
d\mathbf{a}_{k-gn} & =  \text{Arg} \min_{d\mathbf{a} \in \mathbb{R}^{p.k}}   \,   \frac{1}{2} \Vert \mathbf{r}(  \mathbf{a} ) - \mathbf{M}(\mathbf{a}) d\mathbf{a} \Vert^{2}_{2}  \ ,  \\
d\mathbf{a}_{k-lm} & =  \text{Arg} \min_{d\mathbf{a} \in \mathbb{R}^{p.k}}   \,   \frac{1}{2} \Vert \mathbf{r}(  \mathbf{a} ) - \mathbf{M}(\mathbf{a}) d\mathbf{a} \Vert^{2}_{2} +  \frac{\lambda}{2}  \big\Vert \mathbf{D} d\mathbf{a} \big\Vert^{2}_{2}  \ , 
\end{align*}
where $\lambda$ is the damping Marquardt parameter and  $\mathbf{D}$ is a diagonal scaling matrix of dimension $k.p$, or their variants including a linear constraint to deal with the singularity of the $\mathbf{M}(\mathbf{a})$ matrix as described in Subsection~\ref{jacob:box}.

For more details in a general variable projection context, see Kaufman~\cite{K1975}, where this simplification has been derived for the first time using a QR factorization of the matrix $\mathbf{F}(\mathbf{a})$ and differentiation of orthogonal matrices, Ruhe and Wedin~\cite{RW1980} where a more direct derivation and generalizations are given, and also O'Leary and Rust~\cite{OR2013} for a recent discussion of the respective merits of this Jacobian approximation against the true Jacobian matrix, again for general variable projection NLLS algorithms. Note that, in the computer vision's community, this Kaufman variant of the Gauss-Newton or Levenberg-Marquardt algorithms has been already extensively used for solving Structure-From-Motion (SFM) tasks~\cite{OD2007}\cite{C2008a}\cite{C2008b}\cite{D2011}\cite{OYD2011}\cite{GM2011}\cite{HF2015}\cite{HZF2017} and is incorrectly called the Wiberg's algorithm in reference of the conference paper~\cite{W1976}. However, the first application of this algorithm to solve WLRA problems (with binary weights) is in fact due to Ruhe~\cite{R1974}. Furthermore, again in the computer vision community, this algorithm has frequently be assumed to be similar to the block ALS algorithm~\cite{SIR1995}\cite{BF2005} (described in Section~\ref{nipals:box}), which is also incorrect as shown above. The first correct derivation of this variant of the Gauss-Newton algorithm to solve WLRA problems with binary weights in the computer vision field is due to Okatani and Deguchi~\cite{OD2007}, see also the excellent Master thesis of Daskalov~\cite{D2011} in which this derivation is revisited.

Now, the Gauss-Seidel step, $d\mathbf{a}_{gs-gn}$, is closely related to the Kaufman step, $d\mathbf{a}_{k-gn}$, since it corresponds to applying the Gauss-Seidel iteration to the linear system appearing in the Kaufman-Gauss-Newton iteration, see~\cite{RW1980} for details. It is also closely related to the block ALS algorithm described in Section~\ref{nipals:box}, as we will demonstrate now.

The Gauss-Seidel step, $d\mathbf{a}_{gs-gn}$, is computed as the minimum 2-norm solution of the linear least-squares problem
\begin{equation*}
\min_{d\mathbf{a} \in \mathbb{R}^{p.k}} \,   \Vert \mathbf{r}( \mathbf{a}_{i} ) - \mathbf{K}_{(n,p)}  \mathbf{G}( \mathbf{b}_{i}  ) d\mathbf{a}   \Vert^{2}_{2}
\end{equation*}
and we have
\begin{equation*}
\mathbf{r}(\mathbf{a}_{i}) = \mathbf{P}^{\bot}_{\mathbf{F}(  \mathbf{a}_{i} ) } \mathbf{x}  =  \mathbf{x} - \mathbf{F}( \mathbf{a}_{i}  ) \mathbf{b}_{i} =  \mathbf{x} -  \mathbf{K}_{(n,p)}  \mathbf{G}( \mathbf{b}_{i}  ) \mathbf{a}_{i} \  ,
\end{equation*}
as demonstrated in Subsection~\ref{varpro_wlra:box}. Hence,
\begin{align*}
\mathbf{r}( \mathbf{a}_{i} ) - \mathbf{K}_{(n,p)}  \mathbf{G}( \mathbf{b}_{i}  ) d\mathbf{a} & =  \mathbf{x} -  \mathbf{K}_{(n,p)}  \mathbf{G}( \mathbf{b}_{i}  ) \mathbf{a}_{i} -  \mathbf{K}_{(n,p)}  \mathbf{G}( \mathbf{b}_{i}  ) d\mathbf{a} \\
                & =  \mathbf{x} -  \mathbf{K}_{(n,p)}  \mathbf{G}( \mathbf{b}_{i}  ) \big( \mathbf{a}_{i} + d\mathbf{a} \big)  \\
                & =  \mathbf{K}_{(n,p)}  \Big(   \mathbf{z}  - \mathbf{G}( \mathbf{b}_{i}  ) \big( \mathbf{a}_{i} + d\mathbf{a} \big) \Big)  \ .
\end{align*}
Furthermore,
\begin{equation*}
\Vert  \mathbf{r}( \mathbf{a}_{i} ) - \mathbf{K}_{(n,p)}  \mathbf{G}( \mathbf{b}_{i}  ) d\mathbf{a} \Vert^{2}_{2} = \Vert \mathbf{z}  - \mathbf{G}( \mathbf{b}_{i}  ) \big( \mathbf{a}_{i} + d\mathbf{a} \big)  \Vert^{2}_{2} \ ,
\end{equation*}
since $\mathbf{K}_{(n,p)}$ is an orthogonal matrix, and we see that the Gauss-Seidel step, $d\mathbf{a}_{gs-gn}$, solves the linear least-squares problem
\begin{equation*}
\min_{d\mathbf{a} \in \mathbb{R}^{p.k}} \,  \Vert \mathbf{z}  - \mathbf{G}( \mathbf{b}_{i}  ) \big( \mathbf{a}_{i} + d\mathbf{a} \big)  \Vert^{2}_{2} =   \Vert \big( \mathbf{z}  - \mathbf{G}( \mathbf{b}_{i}  ) \mathbf{a}_{i} \big) - \mathbf{G}( \mathbf{b}_{i}  ) d\mathbf{a}  \Vert^{2}_{2} \ ,
\end{equation*}
while the ALS iteration computes $\mathbf{a}_{i+1}$ as the minimum 2-norm solution of the linear least-squares problem
\begin{equation*}
\min_{\mathbf{a} \in \mathbb{R}^{p.k}} \,  \Vert \mathbf{z}  - \mathbf{G}( \mathbf{b}_{i}  ) \mathbf{a}   \Vert^{2}_{2} \  .
\end{equation*}
If $\mathbf{G}( \mathbf{b}_{i}  )$ is a full column-rank matrix, we then have
\begin{equation*}
\mathbf{a}_{i+1} =    \mathbf{G}( \mathbf{b}_{i}  )^{+}  \mathbf{z} = \big(  \mathbf{G}( \mathbf{b}_{i}  )^{T}  \mathbf{G}( \mathbf{b}_{i}  ) \big)^{-1} \mathbf{G}( \mathbf{b}_{i}  )^{T} \mathbf{z}
\end{equation*}
and
\begin{align*}
d\mathbf{a}_{gs-gn} & =   \mathbf{G}( \mathbf{b}_{i}  )^{+}   \big( \mathbf{z}  - \mathbf{G}( \mathbf{b}_{i}  ) \mathbf{a}_{i} \big)  \\
                                 & =  \big(  \mathbf{G}( \mathbf{b}_{i}  )^{T}  \mathbf{G}( \mathbf{b}_{i}  ) \big)^{-1} \mathbf{G}( \mathbf{b}_{i}  )^{T}  \big( \mathbf{z}  - \mathbf{G}( \mathbf{b}_{i}  ) \mathbf{a}_{i} \big) \\
                                 & =   \big(  \mathbf{G}( \mathbf{b}_{i}  )^{T}  \mathbf{G}( \mathbf{b}_{i}  ) \big)^{-1} \mathbf{G}( \mathbf{b}_{i}  )^{T}  \mathbf{z} - \mathbf{a}_{i} \\
                                 & =  \mathbf{a}_{i+1} - \mathbf{a}_{i} \  ,
\end{align*}
and, in these conditions, the ALS and Gauss-Seidel-Gauss-Newton algorithms generate exactly the same iterates. On the other hand, if $\mathbf{G}( \mathbf{b}_{i}  )$ is not a full column-rank matrix, we still have the equality
\begin{equation*}
 \Vert \mathbf{z}  - \mathbf{G}( \mathbf{b}_{i}  ) \mathbf{a}_{i+1}   \Vert^{2}_{2} =  \Vert \big( \mathbf{z}  - \mathbf{G}( \mathbf{b}_{i}  ) \mathbf{a}_{i} \big) - \mathbf{G}( \mathbf{b}_{i}  ) d\mathbf{a}_{gs-gn}   \Vert^{2}_{2} \ .
\end{equation*}
However, in general $\mathbf{a}_{i+1}  \ne \mathbf{a}_{i} + d\mathbf{a}_{gs-gn}$ since the Gauss-Seidel iteration produces the minimum 2-norm correction vector $d\mathbf{a}_{gs-gn}$ to the above linear least-squares problem  while the ALS algorithm obtains the minimum 2-norm solution $\mathbf{a}_{i+1}$ of this linear least-squares problem. Thus, unless $\mathbf{a}_{i+1}$  belongs to the correct manifold, the Gauss-Seidel and ALS steps do not produce the same iterate when the matrix $\mathbf{G}( \mathbf{b}_{i}  )$ is not of  full column-rank.

We finally observe that a line search in step \textbf{(7.2)} of the Gauss-Newton algorithms~\eqref{gn_alg:box}  is not required for the Gauss-Seidel correction $d\mathbf{a}_{gs-gn}$ in order to obtain the inequality $\psi( \mathbf{a}_{i+1} )  \le \psi( \mathbf{a}_{i} )$ and global convergence of the iterations, see Section~\ref{nipals:box} for details. On the other hand, for both the Golub-Pereyra correction $d\mathbf{a}_{gp-gn}$ and the Kaufman correction $d\mathbf{a}_{k-gn}$, it may happen that $\psi( \mathbf{a}_{i+1} )  > \psi( \mathbf{a}_{i} )$ meaning that the $\psi(.)$ surface is not reliably approximated by a quadratic function. In other words, the quadratic approximation is only good in the local neighborhood of $\mathbf{a}_{i}$, not at the bottom of the quadratic valley that the Gauss-Newton approach uses. In such cases, a line search algorithm to determine $\alpha_i$ at each iteration such that $\psi( \mathbf{a}_{i+1} )  < \psi( \mathbf{a}_{i} )$ must be incorporated in step \textbf{(7.2)} of the Gauss-Newton algorithms~\eqref{gn_alg:box} in order to obtain global convergence. Note also that this is always possible despite the singularity of the Jacobian matrix $- \big( \mathbf{M}( \mathbf{a}_{i}  ) + \mathbf{L}( \mathbf{a}_{i}  )  \big)$ or its Kaufman approximation $- \mathbf{M}( \mathbf{a}_{i}  )$ (see Theorem~\ref{theo5.2:box}) as the correction vectors $d\mathbf{a}_{gp-gn}$ and $d\mathbf{a}_{k-gn}$ are in a descent direction for $\psi(.)$ if $\nabla \psi( \mathbf{a}_{i} ) \ne  \mathbf{0}^{k.p}$ (see Corollary~\ref{corol5.7:box}). However, to develop damped versions of these Gauss-Newton algorithms by implementing a line search algorithm, we have  to perform the second part of step \textbf{(4)}  of the Gauss-Newton algorithms~\eqref{gn_alg:box}, every time we want to get  $\psi( \mathbf{a}_{i+1} )$ for a new trial value of $\mathbf{a}_{i+1}$, since
\begin{equation*}
\psi(\mathbf{a}_{i+1} ) = \frac{1}{2} \Vert \mathbf{P}^{\bot}_{\mathbf{F}(  \mathbf{a}_{i+1} ) }\mathbf{x} \Vert^{2}_{2}  \ ,
\end{equation*}
and a line search can involve many extra evaluations of $\psi(.)$, which do not get us closer to the solution. Furthermore, if a line search is incorporated, the Gauss-Newton algorithms~\eqref{gn_alg:box} must be slightly reorganized to avoid duplicate computations in steps \textbf{(4)} and \textbf{(7)}, but we omit these details here.

 In these conditions, to obtain global convergence, it is tempting to perform one iteration with a Gauss-Seidel step $d\mathbf{a}_{gs-gn}$ or even several iterations with the fast block ALS method described in Section~\ref{nipals:box} to compute $\mathbf{a}_{i+1}$  in step \textbf{(7.2)} (e.g., if $\psi( \mathbf{a}_{i} + d\mathbf{a}_{gn} ) > \psi( \mathbf{a}_{i})$) instead of using a more costly line search. In other words, if a full Gauss-Newton step gives a sufficient decrease of $\psi(.)$, we accept this point as the new iterate. Otherwise we switch to the fast Gauss-Seidel or block ALS methods.

We now explain how the matrices $\mathbf{M}( \mathbf{a}_{i} )$ and $-\mathit{J}( \mathbf{r}(\mathbf{a}_{i}) ) = \mathbf{M}( \mathbf{a}_{i} ) + \mathbf{L}( \mathbf{a}_{i} )$ and their QR factorizations can be computed efficiently and with reduced storage in order to obtain the correction vectors $d\mathbf{a}_{k-gn}$  or $d\mathbf{a}_{gp-gn}$ at each iteration of the Gauss-Newton algorithms~\eqref{gn_alg:box}. To this end, we first recall from the results of Subsection~\ref{jacob:box} that we have the following explicit expressions for these matrices:
\begin{align*}
    \mathbf{M}(\mathbf{a})  & =  \mathbf{P}_{\mathbf{F}(\mathbf{a})}^{\bot} \mathbf{U(\mathbf{a})}  \  , \\
    \mathbf{L}(\mathbf{a})  & =   ( \mathbf{F}(\mathbf{a})^{+}  )^{T}  \mathbf{V(\mathbf{a})}  \  , \\
    - \mathit{J}( \mathbf{r}(\mathbf{a}) ) & =  \mathbf{M}(\mathbf{a}) + \mathbf{L}(\mathbf{a}) = \mathbf{P}_{\mathbf{F}(\mathbf{a})}^{\bot} \mathbf{U}(\mathbf{a}) + ( \mathbf{F}(\mathbf{a})^{+}  )^{T}  \mathbf{V} (\mathbf{a}) \  ,
\end{align*}
with
\begin{align*}
    \mathbf{U}(\mathbf{a}) & =  \emph{diag}( \emph{vec}(\sqrt{\mathbf{W}} ) ) (\mathbf{\widehat{B}}^{T} \otimes \mathbf{I}_{p}) \mathbf{K}_{(k,p)} = \mathbf{K}_{(n,p)} \mathbf{G}(\mathbf{\widehat{b}})  \  , \\
    \mathbf{V}(\mathbf{a}) & =  ( \mathbf{W} \odot P_{\Omega}(\mathbf{X} -\mathbf{A}\mathbf{\widehat{B}}) )^{T} \otimes \mathbf{I}_{k}  \  ,
\end{align*}
as stated in equations~\eqref{eq:U_mat},~\eqref{eq:V_mat} and~\eqref{eq:J_mat}. The notations here are exactly the same as in Subsection~\ref{jacob:box} and we have drop again the iteration index of the Gauss-Newton iterations in order to simplify the notations in the rest of this section. Note  that, in the definition of $\mathbf{L}(\mathbf{a})$, $\mathbf{F}(\mathbf{a})^{+}$ can be replaced by a symmetric generalized inverse $\mathbf{F}(\mathbf{a})^{-}$ at our convenience. We also recall that the matrices $-\mathit{J}( \mathbf{r}(\mathbf{a}) )$, $\mathbf{M}( \mathbf{a} )$ and $\mathbf{L}( \mathbf{a} )$ have $n.p$ rows if $\mathbf{W} \in \mathbb{R}^{p \times n}_{+*}$, but only $nobs$ rows if $\mathbf{W} \in \mathbb{R}^{p \times n}_{+}$, where $nobs$  is the number of non-zero rows of $\mathbf{F}( \mathbf{a} )$. As explained in Subsection~\ref{jacob:box}, $nobs$ is simply the number of "non-missing" elements in the data matrix $\mathbf{X}$ or, equivalently, the number of non-zero weights in the weight matrix $\mathbf{W}$, namely,
\begin{equation*}
nobs = \displaystyle{ \sum_{ij} \boldsymbol{\delta}_{ij} }  \  ,
\end{equation*}
where $\boldsymbol{\delta}$ is the incidence matrix associated the weight matrix (also defined in Subsection~\ref{jacob:box}). From the above equations, we see that the matrix $-\mathit{J}( \mathbf{r}(\mathbf{a}) )$ and its two matrix components are tall and skinny in most cases, even if the number of missing values is high, as $k$ is expected to be much smaller than min$(p,n)$ and that their evaluations require the computations of the orthogonal projector $\mathbf{P}_{\mathbf{F}(\mathbf{a})}^{\bot}$, and,  of $\mathbf{F}(\mathbf{a})^{+}$ and  $\mathbf{\widehat{b}} = \mathbf{F}( \mathbf{a} )^{+}  \mathbf{x}$ if a COD is used in step \textbf{(2)} of the Gauss-Newton algorithms~\eqref{gn_alg:box} or, alternatively,  of $\mathbf{F}(\mathbf{a})^{-}$ and $\mathbf{\widehat{b}} = \mathbf{F}( \mathbf{a} )^{-}  \mathbf{x}$ if a QRCP is used in this step \textbf{(2)}.

The key-observation to compute efficiently these different matrices is to remember that $\mathbf{F}(\mathbf{a})$ is a block-diagonal matrix, namely,
\begin{equation*}
\mathbf{F}(\mathbf{a}) =  \bigoplus_{j=1}^n \mathbf{F}_{j}(  \mathbf{a} )  \text{ , where }  \mathbf{F}_{j}(  \mathbf{a} ) =  \emph{diag}(\sqrt{\mathbf{W}}_{.j})\mathbf{A}  \  .
\end{equation*}
In these conditions, as already discussed in Subsection~\ref{jacob:box}, we have
\begin{equation*}
\emph{rank} ( \mathbf{F}( \mathbf{a} ) ) = \sum^{n}_{j=1} \emph{rank} ( \mathbf{F}_{j}( \mathbf{a} ) ) =  \sum^{n}_{j=1} r_{j} = {r}_{\mathbf{F}( \mathbf{a} )}
\end{equation*}
and
\begin{equation*}
\mathbf{P}_{\mathbf{F}(\mathbf{a})}^{\bot} =   \bigoplus^{n}_{j=1} \mathbf{P}_{\mathbf{F}_{j}(\mathbf{a})}^{\bot} \text{ , } \mathbf{F}( \mathbf{a} )^{+} =  \bigoplus^{n}_{j=1}  \mathbf{F}_{j}( \mathbf{a} )^{+} \text{ and } \mathbf{F}( \mathbf{a} )^{-} =  \bigoplus^{n}_{j=1}  \mathbf{F}_{j}( \mathbf{a} )^{-}.
\end{equation*}
Taking advantage of these block-structures of $\mathbf{P}_{\mathbf{F}(\mathbf{a})}^{\bot}$ and  of the generalized inverses of $\mathbf{F}(\mathbf{a})$ can reduce drastically the required storage and allows us to use efficient parallelization techniques for reducing the computing time needed to solve large WLRA problems using variable projection second-order algorithms as we will illustrate below. Interestingly, the techniques used for this purpose are very similar to those developed for solving large and dense structured linear least-squares problems arising in the context of separable NLLS problems with multiple right hand sides (e.g., NLLS problems in which a linear combination of nonlinear functions is fit linearly to data in many datasets); see Kaufman and Silvester~\cite{KS1992},  Kaufman et al.~\cite{KSW1994} and Kaufman~\cite{K2010} for more details.

Taking into account the block-structures of $\mathbf{P}_{\mathbf{F}(\mathbf{a})}^{\bot}$, $\mathbf{F}( \mathbf{a} )^{+}$ and $\mathbf{F}( \mathbf{a} )^{-}$, we first observe that the $n.p  \times k.p$ matrices $-\mathit{J}( \mathbf{r}(\mathbf{a}) )$, $\mathbf{M}( \mathbf{a} )$ and $\mathbf{L}( \mathbf{a} )$ can be divided into $n$ blocks, each of shape $p  \times k.p$ if $\mathbf{W} \in \mathbb{R}^{p \times n}_{+*}$:
\begin{equation} \label{eq:JML_mat_blocking}
-\mathit{J}( \mathbf{r}(\mathbf{a}) ) = \begin{bmatrix} \mathbf{J}_{1}  \\ \vdots \\  \mathbf{J}_{j}  \\ \vdots   \\ \mathbf{J}_{n} \end{bmatrix} \text{, }  \mathbf{M}( \mathbf{a} ) = \begin{bmatrix} \mathbf{M}_{1}  \\ \vdots \\  \mathbf{M}_{j}  \\ \vdots   \\ \mathbf{M}_{n} \end{bmatrix}   \text{ and }  \mathbf{L}( \mathbf{a} ) = \begin{bmatrix} \mathbf{L}_{1}  \\ \vdots \\  \mathbf{L}_{j}  \\ \vdots   \\ \mathbf{L}_{n} \end{bmatrix}  \  .
\end{equation}
If $\mathbf{W} \in \mathbb{R}^{p \times n}_{+}$, these matrices can also be divided  into $n$ blocks, but the number of rows in each block  $\mathbf{J}_{j}$, $\mathbf{M}_{j}$ and  $\mathbf{L}_{j}$ will differ and will be equal to the number of non-missing elements in the corresponding column of the data matrix $\mathbf{X}$. In order to simplify the exposition, but without loss of generality, we will assume in the rest of this section that $\mathbf{W} \in \mathbb{R}^{p \times n}_{+*}$. Obviously, if $\mathbf{W} \in \mathbb{R}^{p \times n}_{+}$, the zero-rows of these submatrices should be eliminated in real computations.

For the same reasons, the matrices $\mathbf{U}(\mathbf{a})$ and $\mathbf{V}(\mathbf{a})$ involved, respectively, in the definitions of $\mathbf{M}( \mathbf{a} )$ and $\mathbf{L}( \mathbf{a} )$ can also be considered as stacks of $n$ blocks.
We also observe that these two matrices are very sparse as $\mathbf{U}(\mathbf{a})$ is a row-permuted block diagonal matrix and $\mathbf{V}(\mathbf{a})$ is the Kronecker product of a matrix with $\mathbf{I}_{k}$, the identity matrix of order $k$. However, they have also a well-defined regular structure for the positions of their non-zero elements, which can be exploited in practical computations. As an illustration, it is easily checked, using the equality, $\mathbf{U}(\mathbf{a}) = \mathbf{K}_{(n,p)} \mathbf{G}(\mathbf{\widehat{b}})$, that the matrix $\mathbf{U}(\mathbf{a})$ as the following block structure with at most $k$ non-zero elements in each of its rows and at most $n$ non-zero elements in each of its columns:
\begin{equation} \label{eq:blk_U_mat}
\mathbf{U}(\mathbf{a}) = \begin{bmatrix} \mathbf{U}_{1}  \\ \vdots \\  \mathbf{U}_{j}  \\ \vdots   \\ \mathbf{U}_{n} \end{bmatrix} \text{ with } \mathbf{U}_{j} \in \mathbb{R}^{p \times k.p} \text{ for }  j = 1, \cdots, n \text{ and }
\end{equation}
\begin{equation*}
\mathbf{U}_{j} =
\left\lbrack  \begin{array}{ccccc}
   \sqrt{\mathbf{W}}_{1j} (\mathbf{\widehat{B}}_{.j})^{T}  & 0    & \ldots     & 0          & 0          \\
    0   &  \sqrt{\mathbf{W}}_{2j} (\mathbf{\widehat{B}}_{.j})^{T}  &  0          & \ldots   & 0          \\
    \vdots   & \ddots                                                                 & \ddots    & \ddots   & \vdots     \\
    0   & \ldots  & 0       & \sqrt{\mathbf{W}}_{(p-1)j} (\mathbf{\widehat{B}}_{.j})^{T}   & 0          \\
    0          & 0          &\ldots    & 0   & \sqrt{\mathbf{W}}_{pj} (\mathbf{\widehat{B}}_{.j})^{T}
\end{array} \right\rbrack \ .
\end{equation*}
Similarly, it is easily verified that the matrix $\mathbf{V}(\mathbf{a})$ has the following block structure with at most $n$ non-zero elements in each of its columns and at most  $p$ non-zero elements in each of its rows:
\begin{equation} \label{eq:blk_V_mat}
\mathbf{V}(\mathbf{a}) = \begin{bmatrix} \mathbf{V}_{1}  \\ \vdots \\  \mathbf{V}_{j}  \\ \vdots   \\ \mathbf{V}_{n} \end{bmatrix} \text{ with } \mathbf{V}_{j} \in \mathbb{R}^{k \times k.p} \text{ for }  j = 1, \cdots, n  \text{ and } \mathbf{V}_{j} = \begin{bmatrix} \mathbf{Z}_{1}  & \cdots &  \mathbf{Z}_{i}  & \cdots  & \mathbf{Z}_{p} \end{bmatrix} \ ,
\end{equation}
with $\mathbf{Z}_{i} \in \mathbb{R}^{k \times k}$ and  $\mathbf{Z}_{i} = \beta_{i} \mathbf{I}_{k}$ for $ i = 1, \cdots, p$, where $\beta_{i} = \mathbf{W}_{ij} ( \mathbf{\bar{X}}_{ij} -  \sum_{l=1}^k { \mathbf{A}_{il} \mathbf{\widehat{B}}_{lj} } )$,  $\mathbf{I}_{k}$ is the identity matrix of order $k$ and 
\begin{equation*}
     \mathbf{\bar{X}}_{ij} =
    \begin{cases}
        \mathbf{X}_{ij}  & \text{if } \mathbf{W}_{ij} \ne 0 \\
         0                      & \text{if } \mathbf{W}_{ij}  =   0
    \end{cases} \ .
\end{equation*}

Finally, if $\mathbf{W} \in \mathbb{R}^{p \times n}_{+*}$, the residual vector $\mathbf{r}(\mathbf{a}) = \mathbf{P}^{\bot}_{\mathbf{F}(  \mathbf{a} ) } \mathbf{x}$ can also be considered as a stack of $n$ p-vectors with
\begin{equation}  \label{eq:res_vec_blocking}
\mathbf{r}(\mathbf{a}) = \begin{bmatrix} \mathbf{r}_{1}(\mathbf{a}) \\ \vdots \\  \mathbf{r}_{j}(\mathbf{a})  \\ \vdots   \\ \mathbf{r}_{n}(\mathbf{a})  \end{bmatrix}  \text{ with } \mathbf{r}_{j}(\mathbf{a}) =   \mathbf{P}^{\bot}_{\mathbf{F}_{j}(  \mathbf{a} ) }  \mathbf{x}_{j} \text{ and } \mathbf{x}_{j} = \sqrt{\mathbf{W}}_{.j} \odot \mathbf{X}_{.j} \ .
\end{equation}
The subvector $\mathbf{r}_{j}(\mathbf{a})$ is the $j^{th}$ residual vector associated with the $j^{th}$ atomic function $\psi_{j}(.)$ defined in equation~\eqref{eq:psi_atomic_func} of Subsection~\ref{varpro_wlra:box}. Using these different block-structures, we can use the following strategy for computing $\mathbf{M}( \mathbf{a} )$ and  $-\mathit{J}( \mathbf{r}(\mathbf{a}) )$ in $n$ independent steps.

At the $j^{th}$ step, we compute the blocks $\mathbf{M}_j$, $\mathbf{L}_j$ and $\mathbf{J}_j$ defined above, namely,
\begin{align*}
    \mathbf{J}_{j}  & =  \mathbf{M}_{j} + \mathbf{L}_{j}  \ , \\
    \mathbf{M}_{j} & =  \mathbf{P}_{\mathbf{F}_{j} (\mathbf{a})}^{\bot} \mathbf{U}_{j}  \ , \\
    \mathbf{L}_{j} & =  \mathbf{F}_{j}(\mathbf{a})^{+T}  \mathbf{V}_{j}   \text{ or }  \mathbf{L}_{j} =  \mathbf{F}_{j}(\mathbf{a})^{-T}  \mathbf{V}_{j}  \ .
\end{align*}
To this end, we need to process the $j^{th}$ columns of $\mathbf{X}$  and $\mathbf{W}$, and we first compute the matrix $\mathbf{F}_{j}(\mathbf{a})$ as
\begin{equation*}
 \mathbf{F}_{j}(  \mathbf{a} ) =  \emph{diag}(\sqrt{\mathbf{W}}_{.j})\mathbf{A}  \ ,
\end{equation*}
where $\mathbf{A}$ is the current estimate for this matrix variable of the factor model and we eliminate eventually the zero rows in $\mathbf{F}_{j}(  \mathbf{a} )$ if some elements of the weight column-vector $\mathbf{W}_{.j}$ are equal to zero. Next, from the above equations, we see that the computation of $\mathbf{J}_{j}$ and its two matrix components requires the computations of $\mathbf{P}_{\mathbf{F}_{j} (\mathbf{a})}^{\bot}$, $\mathbf{F}_{j}(\mathbf{a})^{+}$ and $\mathbf{\widehat{B}}_{.j}=\mathbf{F}_{j}(\mathbf{a})^{+} \mathbf{x}_{j}$, or, $\mathbf{F}_{j}(\mathbf{a})^{-}$ and $\mathbf{\widehat{B}}_{.j}=\mathbf{F}_{j}(\mathbf{a})^{-} \mathbf{x}_{j}$, where $\mathbf{x}_{j} = \sqrt{\mathbf{W}}_{.j} \odot \mathbf{X}_{.j}$. This implies that the two matrix components of $\mathbf{J}_{j}$ can be obtained from a QRCP or a COD (see equations~\ref{eq:qrcp} and~\ref{eq:cod} in Subsection~\ref{lin_alg:box}) of the $p \times k$ matrix $\mathbf{F}_{j}(\mathbf{a})$. See also Golub and Pereyra~\cite{GP1973}~\cite{GP1976}, Krogh~\cite{K1974}, Kaufman~\cite{K1975} and  Gay and Kaufman~\cite{GK1991} for more details in a more general framework dealing with general separable NLLS problems. Thus, we first compute the QRCP of $\mathbf{F}_{j}(\mathbf{a})$ as
\begin{equation*}
\mathbf{F}_{j}(\mathbf{a})  =  \mathbf{Q}^{T}_{j} \begin{bmatrix} \mathbf{R}_{j}    &  \mathbf{S}_{j}   \\  \mathbf{0}^{(p-r_{j}) \times r_{j}}  & \mathbf{0}^{ (p-r_{j}) \times (k-r_{j}) }  \end{bmatrix} \mathbf{P}^{T}_{j} \ ,
\end{equation*}
where $\mathbf{Q}_{j}$ is an $p \times p$ orthogonal matrix, $\mathbf{R}_{j}$ is an $r_{j} \times r_{j}$ nonsingular upper triangular matrix with $r_{j} =  \emph{rank}( \mathbf{F}_{j}(\mathbf{a}) )$, which can be estimated during the QRCP, $\mathbf{S}_{j}$ is vacuous unless $\mathbf{F}_{j}(\mathbf{a})$ is rank deficient and $\mathbf{P}_{j}$ is a $k \times k$ permutation matrix. Note that $\mathbf{F}_{j}(\mathbf{a})$ is a dense matrix, so that its QRCP can be computed efficiently and cheaply with the help of standard dense methods, see Subsection~\ref{lin_alg:box} and~\cite{GVL1996}\cite{B2015} for details.
From the above QRCP of $\mathbf{F}_{j}(\mathbf{a})$, we can compute $\mathbf{P}_{\mathbf{F}_{j} (\mathbf{a})}^{\bot}$ as
\begin{equation*}
\mathbf{P}_{\mathbf{F}_{j} (\mathbf{a})}^{\bot} = \mathbf{Q}^{T}_{j} \begin{bmatrix} \mathbf{0}^{ r_{j} \times r_{j} }     &  \mathbf{0}^{ r_{j} \times (p - r_{j} ) }   \\  \mathbf{0}^{(p-r_{j}) \times r_{j}}  & \mathbf{I}^{ p-r_{j} }  \end{bmatrix} \mathbf{Q}_{j}
\end{equation*}
and also a symmetric generalized inverse of $\mathbf{F}_{j} (\mathbf{a})$ as
\begin{equation*}
\mathbf{F}_{j}(\mathbf{a})^{-}  = \mathbf{P}_{j}   \begin{bmatrix} \mathbf{R}^{-1}_{j}    &  \mathbf{0}_{  r_{j} \times  (p-r_{j}) }   \\  \mathbf{0}^{(k-r_{j}) \times r_{j}}  & \mathbf{0}^{ (k-r_{j}) \times (p-r_{j}) }  \end{bmatrix}  \mathbf{Q}_{j} \ .
\end{equation*}
See again Subsection~\ref{lin_alg:box} for more information. Next, from this formulation of $\mathbf{F}_{j}(\mathbf{a})^{-}$, the $j^{th}$ column of $\mathbf{\widehat{B}}$ can be computed as
\begin{equation*}
\mathbf{\widehat{B}}_{.j} = \mathbf{F}_{j}(\mathbf{a})^{-} \mathbf{x}_{j} =  \mathbf{P}^{1}_{j}  \mathbf{R}^{-1}_{j}  \mathbf{Q}^{1}_{j}  \mathbf{x}_{j} \ .
\end{equation*}
In this last equation, the orthogonal matrices $\mathbf{Q}_{j}$ and $\mathbf{P}_{j}$ have been partitioned as
\begin{equation*}
\mathbf{Q}_{j} = \begin{bmatrix}    \mathbf{Q}^{1}_{j} \\   \mathbf{Q}^{2}_{j}  \end{bmatrix}    \text{ and } \begin{bmatrix}    \mathbf{P}^{1}_{j} \ &  \mathbf{P}^{2}_{j}  \end{bmatrix} \ ,
\end{equation*}
where
\begin{itemize}
\item $\mathbf{Q}^{1}_{j}$ and $\mathbf{Q}^{2}_{j}$ have, respectively, $r_{j}$ and $p - r_{j}$ rows \ ,
\item $\mathbf{P}^{1}_{j}$ and $\mathbf{P}^{2}_{j}$ have, respectively, $r_{j}$ and $k - r_{j}$ columns \ .
\end{itemize}

 If $r_{j} = k$ then $\mathbf{F}_{j}(\mathbf{a})^{+} = \mathbf{F}_{j}(\mathbf{a})^{-}$. On the other hand,  if $\mathbf{F}_{j}(\mathbf{a})$ is singular, we can optionally compute its COD from its QRCP in order to obtain $\mathbf{F}_{j}(\mathbf{a})^{+}$ and compute $\mathbf{\widehat{B}}_{.j}$ as $\mathbf{F}_{j}(\mathbf{a})^{+} \mathbf{x}_{j}$. However, taking into account the special structure and indeterminacy associated with the solutions $(\mathbf{\widehat{A}},\mathbf{\widehat{B}})$ of the WLRA problem (see Section~\ref{seppb:box}), we don't really need to compute $\mathbf{F}_{j}(\mathbf{a})^{+}$ even if $\mathbf{F}_{j}(\mathbf{a})$ is rank-deficient, so we omit here the details of the optional computation of the COD of $\mathbf{F}_{j}(\mathbf{a})$.
 
 Once $\mathbf{\widehat{B}}_{.j}$ has been computed, the submatrices $\mathbf{U}_{j}$ and $\mathbf{V}_{j}$ defined above can then be evaluated, or more precisely are available, to compute $\mathbf{M}_{j}$ and $\mathbf{L}_{j}$. Finally, $\mathbf{r}_{j}(\mathbf{a})$ can be computed as follows
\begin{equation*}
\mathbf{r}_{j}(\mathbf{a}) = \mathbf{P}^{\bot}_{\mathbf{F}_{j}(  \mathbf{a}_{i} ) } \mathbf{x}_{j} =  \mathbf{Q}_{j}^{T} \begin{bmatrix}  \mathbf{0}^{ r_{j} } \\  \mathbf{Q}^{2}_{j}  \mathbf{x}_{j}   \end{bmatrix} = ( \mathbf{Q}^{2}_{j} )^{T} ( \mathbf{Q}^{2}_{j}  \mathbf{x}_{j} ) \ ,
\end{equation*}
using the above results, and similarly if a QRCP or COD of $\mathbf{F}_{j}(\mathbf{a})$ is available. Inserting now the above expressions for $\mathbf{P}_{\mathbf{F}_{j} (\mathbf{a})}^{\bot}$ and $\mathbf{F}_{j}(\mathbf{a})^{-}$ in the definition of $\mathbf{J}_{j}$, we obtain
\begin{align*}
    \mathbf{J}_{j} & =   \mathbf{P}_{\mathbf{F}_{j}(\mathbf{a})}^{\bot} \mathbf{U}_{j} + ( \mathbf{F}_{j}(\mathbf{a})^{-}  )^{T}  \mathbf{V}_{j}  \\
                          & =   \mathbf{Q}^{T}_{j} \begin{bmatrix} \mathbf{0}^{ r_{j} \times r_{j} }     &  \mathbf{0}^{ r_{j} \times (p - r_{j} ) }   \\  \mathbf{0}^{(p-r_{j}) \times r_{j}}  & \mathbf{I}_{ p-r_{j} }  \end{bmatrix}  \mathbf{Q}_{j} \mathbf{U}_{j}  +    \mathbf{Q}^{T}_{j}   \begin{bmatrix} \mathbf{R}^{-T}_{j}    &  \mathbf{0}^{  r_{j} \times  (k-r_{j}) }   \\  \mathbf{0}^{(p-r_{j}) \times r_{j}}  & \mathbf{0}^{ (p-r_{j}) \times (k-r_{j}) }  \end{bmatrix}  \mathbf{P}^{T}_{j}  \mathbf{V}_{j} \\
                          & =   \mathbf{Q}^{T}_{j} \Big( \begin{bmatrix} \mathbf{0}^{ r_{j} \times r_{j} }     &  \mathbf{0}^{ r_{j} \times (p - r_{j} ) }   \\  \mathbf{0}^{(p-r_{j}) \times r_{j}}  & \mathbf{I}_{ p-r_{j} }  \end{bmatrix}  \mathbf{Q}_{j} \mathbf{U}_{j}  +  \begin{bmatrix} \mathbf{R}^{-T}_{j}    &  \mathbf{0}^{  r_{j} \times  (k-r_{j}) }   \\  \mathbf{0}^{(p-r_{j}) \times r_{j}}  & \mathbf{0}^{ (p-r_{j}) \times (k-r_{j}) }  \end{bmatrix}  \mathbf{P}^{T}_{j}  \mathbf{V}_{j} \Big) \\
                          & =  \mathbf{Q}^{T}_{j} \Big( \begin{bmatrix}  \mathbf{0}^{ r_{j} \times k.p }   \\   \mathbf{Q}^{2}_{j} \mathbf{U}_{j}    \end{bmatrix}  +  \begin{bmatrix} \mathbf{R}^{-T}_{j}  (\mathbf{P}^{1}_{j})^{T}  \mathbf{V}_{j}     \\  \mathbf{0}^{(p-r_{j}) \times k.p  }   \end{bmatrix}  \Big) \\
                          & = \mathbf{Q}^{T}_{j}  \begin{bmatrix} \mathbf{R}^{-T}_{j}  (\mathbf{P}^{1}_{j})^{T}  \mathbf{V}_{j}  \\   \mathbf{Q}^{2}_{j} \mathbf{U}_{j}    \end{bmatrix} .
\end{align*}
In these conditions, for $d\mathbf{a} \in \mathbb{R}^{k.p}$, we have
\begin{equation*}
\mathbf{r}_{j}(\mathbf{a}) - \mathbf{J}_{j} d\mathbf{a} =  \mathbf{Q}^{T}_{j}  \Big(   \begin{bmatrix}  \mathbf{0}^{ r_{j} } \\  \mathbf{Q}^{2}_{j}  \mathbf{x}_{j}   \end{bmatrix}  -   \begin{bmatrix} \mathbf{R}^{-T}_{j}  (\mathbf{P}^{1}_{j})^{T}  \mathbf{V}_{j}  \\   \mathbf{Q}^{2}_{j} \mathbf{U}_{j}    \end{bmatrix} d\mathbf{a}  \Big) \ .
\end{equation*}
At this point, we introduce several new block matrix and vector definitions again to simplify the notation going forward:
\begin{equation} \label{eq:tj_mat}
\widetilde{\mathit{J}}( \mathbf{r}(\mathbf{a}) ) = \begin{bmatrix} \widetilde{\mathbf{J}}_{1}  \\ \vdots \\  \widetilde{\mathbf{J}}_{j}  \\ \vdots   \\ \widetilde{\mathbf{J}}_{n} \end{bmatrix} \text{ with }   \widetilde{\mathbf{J}}_{j} =  \begin{bmatrix} \mathbf{R}^{-T}_{j}  (\mathbf{P}^{1}_{j})^{T}  \mathbf{V}_{j}  \\   \mathbf{Q}^{2}_{j} \mathbf{U}_{j}    \end{bmatrix} 
\end{equation}
and
\begin{equation} \label{eq:tr_vec}
\widetilde{\mathbf{r}}(\mathbf{a}) = \begin{bmatrix} \widetilde{\mathbf{r}}_{1}  \\ \vdots \\  \widetilde{\mathbf{r}}_{j}  \\ \vdots   \\ \widetilde{\mathbf{r}}_{n}  \end{bmatrix}  \text{ with } \widetilde{\mathbf{r}}_{j} = \begin{bmatrix}  \mathbf{0}^{ r_{j} } \\  \mathbf{Q}^{2}_{j}  \mathbf{x}_{j}   \end{bmatrix}.
\end{equation}
Thus, $\widetilde{\mathbf{r}}(\mathbf{a})$ is an $n.p$-vector, which is a stack of $n$ $(p-r_{j})$-subvectors, separated by $r_{j}$ zero elements in sequential order. Finally, we conceptually define the orthogonal block diagonal matrix
\begin{equation} \label{eq:Q_mat_blocking}
\mathbf{Q}_{\mathbf{F}} = \bigoplus_{j=1}^n \mathbf{Q}_{j} \ .
\end{equation}
With these new notations and the preceding results, we have
\begin{equation} \label{eq:Q_J_mat}
- \mathit{J}( \mathbf{r}(\mathbf{a}) ) = \mathbf{Q}_{\mathbf{F}} \widetilde{\mathit{J}}( \mathbf{r}(\mathbf{a}) )  \text{  and  } \mathbf{r}(\mathbf{a}) = \mathbf{Q}_{\mathbf{F}} \widetilde{\mathbf{r}}(\mathbf{a}) \ .
\end{equation}
Now, as the 2-norm is unitarily invariant, $\forall d\mathbf{a} \in \mathbb{R}^{k.p}$, we have
\begin{equation*}
\Vert \mathbf{r}( \mathbf{a} ) +  \mathit{J} \big ( \mathbf{r}(\mathbf{a}) \big ) d\mathbf{a}   \Vert_{2} = \Vert \mathbf{r}( \mathbf{a} ) - \big( \mathbf{M}( \mathbf{a}  ) + \mathbf{L}( \mathbf{a}  )  \big) d\mathbf{a}   \Vert_{2} = \Vert \widetilde{\mathbf{r}}(\mathbf{a}) -  \widetilde{\mathit{J}} \big ( \mathbf{r}(\mathbf{a}) \big ) d\mathbf{a}   \Vert_{2} \ ,
\end{equation*}
as $\mathbf{Q}_{\mathbf{F}}$ is an $p.n \times p.n$ orthogonal matrix. In other words, the linear least-squares problem
\begin{equation}  \label{eq:llsq_j_mat}
\min_{d\mathbf{a} \in \mathbb{R}^{p.k}} \,   \Vert \mathbf{r}( \mathbf{a} ) - \big( \mathbf{M}( \mathbf{a}  ) + \mathbf{L}( \mathbf{a}  )  \big) d\mathbf{a}   \Vert^{2}_{2} \ ,
\end{equation}
which must be solved at each iteration of the Gauss-Newton algorithm if a Golub-Pereyra step $d\mathbf{a}_{gp-gn}$  is used, is equivalent to the linear least-squares problem
\begin{equation} \label{eq:llsq_tj_mat}
\min_{d\mathbf{a} \in \mathbb{R}^{p.k}} \,   \Vert \widetilde{\mathbf{r}}(\mathbf{a}) -  \widetilde{\mathit{J}}( \mathbf{r}(\mathbf{a}) ) d\mathbf{a}   \Vert^{2}_{2} \ .
\end{equation}
Similarly, if we use the Kaufman variant at each iteration of the Gauss-Newton algorithm~\eqref{gn_alg:box}, it is not difficult to verify using similar arguments that computing a QRCP of $\mathbf{F}_{j} (\mathbf{a})$ at  step \textbf{(2)} of this algorithm is again sufficient and that the associated linear least-squares problem to solve for computing the correction vector $d\mathbf{a}_{k-gn}$, namely,
\begin{equation} \label{eq:llsq_m_mat}
\min_{d\mathbf{a} \in \mathbb{R}^{p.k}} \,   \Vert \mathbf{r}( \mathbf{a} ) - \mathbf{M}( \mathbf{a}  ) d\mathbf{a}   \Vert^{2}_{2} \ ,
\end{equation}
is equivalent to the linear least-squares problem
\begin{equation*}
\min_{d\mathbf{a} \in \mathbb{R}^{p.k}} \,   \Vert \widetilde{\mathbf{r}}(\mathbf{a}) - \widetilde{\mathbf{M}}( \mathbf{a}  ) d\mathbf{a}   \Vert^{2}_{2} \ ,
\end{equation*}
where $\widetilde{\mathbf{M}}( \mathbf{a}  )$ has the following block structure
\begin{equation}\label{eq:tm_mat}
\widetilde{\mathbf{M}}( \mathbf{a} ) = \begin{bmatrix} \widetilde{\mathbf{M}}_{1}  \\ \vdots \\  \widetilde{\mathbf{M}}_{j}  \\ \vdots   \\ \widetilde{\mathbf{M}}_{n} \end{bmatrix} \text{ with }   \widetilde{\mathbf{M}}_{j} = \begin{bmatrix} \mathbf{0}^{ r_{j} \times k.p} \\   \mathbf{Q}^{2}_{j} \mathbf{U}_{j}    \end{bmatrix}  \ ,
\end{equation}
and we also have the matrix equality
\begin{equation} \label{eq:Q_M_mat}
\mathbf{M}( \mathbf{a} ) = \mathbf{Q}_{\mathbf{F}} \widetilde{\mathbf{M}}( \mathbf{a} ) \ .
\end{equation}
Furthermore, as zero rows appearing in the coefficient matrix of a linear least-squares problem do not affect the solution of this linear least-squares problem, these zero rows can be deleted and the linear least-squares problem to be solved at each iteration of the Kaufman variant of the Gauss-Newton algorithm~\eqref{gn_alg:box} reduces, finally, to
\begin{equation} \label{eq:llsq_tm_mat}
\min_{d\mathbf{a} \in \mathbb{R}^{p.k}} \,   \Vert \bar{\mathbf{r}}(\mathbf{a}) - \bar{\mathbf{M}}( \mathbf{a}  ) d\mathbf{a}   \Vert^{2}_{2} \ ,
\end{equation}
where
\begin{equation}\label{eq:ttm_mat}
\bar{\mathbf{M}}( \mathbf{a}  ) = \begin{bmatrix} \bar{\mathbf{M}}_{1}  \\ \vdots \\  \bar{\mathbf{M}}_{j}  \\ \vdots   \\ \bar{\mathbf{M}}_{n} \end{bmatrix} \text{ with }   \bar{\mathbf{M}}_{j} =  \mathbf{Q}^{2}_{j} \mathbf{U}_{j}   
\end{equation}
and
\begin{equation}\label{eq:ttr_vec}
\bar{\mathbf{r}}(\mathbf{a}) = \begin{bmatrix} \bar{\mathbf{r}}_{1}  \\ \vdots \\  \bar{\mathbf{r}}_{j}  \\ \vdots   \\ \bar{\mathbf{r}}_{n}  \end{bmatrix}  \text{ with } \bar{\mathbf{r}}_{j} =  \mathbf{Q}^{2}_{j}  \mathbf{x}_{j} \ .
\end{equation}
$\bar{\mathbf{M}}( \mathbf{a}  )$ and $\bar{\mathbf{r}}(\mathbf{a})$ have, respectively, only $n.p - {r}_{\mathbf{F}( \mathbf{a} )}$ rows and elements if $\mathbf{W} \in \mathbb{R}^{p \times n}_{+*}$, and, $nobs - {r}_{\mathbf{F}( \mathbf{a} )}$ rows and elements if $\mathbf{W} \in \mathbb{R}^{p \times n}_{+}$. Thus, the work involved in solving this linear least-squares problem is further reduced in addition to the simplifications introduced by the use of the approximated Jacobian matrix $-\mathbf{M}( \mathbf{a})$ as the coefficient matrix of the linear least-squares problem.

In all the above alternative formulations of the linear least-squares problems involving the Jacobian matrix or its approximation, which must be solved at each iteration of the Gauss-Newton algorithms~\eqref{gn_alg:box}, we observe that both the coefficient matrix and the right hand-side vector of the associated linear least-squares problems can be computed independently in $n$ steps, which may offer some important speed-up in a parallel environment.

The next critical step is to solve the linear least-squares problems involving the huge, but tall and skinny, matrices $\widetilde{\mathit{J}}( \mathbf{r}(\mathbf{a}) )$ or $\bar{\mathbf{M}}( \mathbf{a}  )$ in a computationally responsible manner. If $\mathbf{W} \in \mathbb{R}^{p \times n}_{+*}$,  $\widetilde{\mathit{J}}( \mathbf{r}(\mathbf{a}) )$ and $\bar{\mathbf{M}}( \mathbf{a}  )$ will have, respectively, $n.p$ and $n.p - {r}_{\mathbf{F}( \mathbf{a} )}$ rows and $k.p$ columns (with $k \ll \text{min}(n,n)$) and it is not conceivable to store such huge matrices in main memory to compute their SVD, QRCP or COD by standard methods as soon as both $p$ and $n$ are relatively large numbers.
We suggest two different strategies to alleviate this problem. The first one uses a QR decomposition of the transformed matrices $\widetilde{\mathit{J}}( \mathbf{r}(\mathbf{a}) )$ or $\bar{\mathbf{M}}( \mathbf{a}  )$ as a preliminary step to reduce the size of the problems and the second one consists in solving the normal equations associated with the linear least-squares problems~\eqref{eq:llsq_j_mat}  and~\eqref{eq:llsq_m_mat} or their transformed versions~\eqref{eq:llsq_tj_mat}  and~\eqref{eq:llsq_tm_mat}.

For most NLLS problems, separable or not, using a QR decomposition of the Jacobian matrix takes at least twice as long as using the normal equations, but gives improved accuracy when the Jacobian matrix is ill-conditioned. As we already know that  $-\mathit{J}( \mathbf{r}(\mathbf{a}) )$ and $\mathbf{M}( \mathbf{a}  )$, or their transformed versions, are always rank-deficient matrices (see Theorem~\ref{theo5.2:box}), this may favor a QR approach for more reliable results. However, computing a Cholesky factor of the Gauss-Newton approximations of the Hessian matrix is substantially faster than computing a QR factorization of the (approximated) Jacobian matrix, especially when this matrix is tall and skinny. This explains why the normal equations approach has been favored in past studies~\cite{OYD2011}\cite{HF2015}\cite{HZF2017} despite the inherent difficulties in this approach to deal efficiently and accurately with the singularity of the Gauss-Newton approximations of the Hessian matrix.

We first describe the iterative methods, which aim at computing the thin QR decomposition of the $n.p \times k.p$ transformed Jacobian matrix
\begin{equation} \label{eq:QR_tjmat}
\widetilde{\mathit{J}}( \mathbf{r}(\mathbf{a}) ) = \widetilde{\mathbf{Q}}_{\mathit{J}} \mathbf{R}_{\mathit{J}}  \  ,
\end{equation}
where $\widetilde{\mathbf{Q}}_{\mathit{J}}$ is an $n.p \times k.p$ matrix (if  $\mathbf{W} \in \mathbb{R}^{p \times n}_{+*}$) with orthonormal columns and $\mathbf{R}_{\mathit{J}}$ is an  $k.p \times k.p$  (singular) upper-triangular matrix, or of its Kaufman approximation
\begin{equation} \label{eq:QR_tmmat}
 \bar{\mathbf{M}}( \mathbf{a}  ) =  \bar{\mathbf{Q}}_{\mathbf{M}}  \mathbf{R}_{\mathbf{M}}  \  ,
\end{equation}
where $\bar{\mathbf{Q}}_{\mathbf{M}} $ and  $\mathbf{R}_{\mathbf{M}}$ have similar shapes than $\bar{\mathbf{M}}( \mathbf{a}  )$ and $\mathbf{R}_{\mathit{J}}$, respectively. In most practical applications $k$ is much smaller than min$(p,n)$ and the matrices $\widetilde{\mathit{J}}( \mathbf{r}(\mathbf{a}) )$ and $\bar{\mathbf{M}}( \mathbf{a}  )$ are tall and skinny, as discussed above, for which highly efficient parallel QR algorithms have been proposed in the literature~\cite{DGHL2012}. These dedicated Tall and Skinny QR (TSQR) algorithms are often referred as communication avoiding algorithms and outperform significantly the conventional Householder  QR algorithm for an $m \times n$ matrix with $m \gg n$. Furthermore, these TSQR algorithms can be combined or more precisely "fused" with the $n$ independent steps described above for the computation of $\widetilde{\mathit{J}}( \mathbf{r}(\mathbf{a}) )$ or $\bar{\mathbf{M}}( \mathbf{a}  )$. This allows us to get the triangular factors $\mathbf{R}_{\mathit{J}}$ or $\mathbf{R}_{\mathbf{M}}$ in the standard QR factorizations of $\widetilde{\mathit{J}}( \mathbf{r}(\mathbf{a}) )$ or $\bar{\mathbf{M}}( \mathbf{a}  )$ without storing in main memory these huge matrices or explicitly computing the orthonormal matrices $\widetilde{\mathbf{Q}}_{\mathit{J}}$ or $\bar{\mathbf{Q}}_{\mathbf{M}}$.

We now focus on two TSQR algorithms for computing these triangular factors and also the matrix-vector products $\widetilde{\mathbf{Q}}_{\mathit{J}}^{T} \widetilde{\mathbf{r}}(\mathbf{a})$ or $\bar{\mathbf{Q}}_{\mathbf{M}}^{T} \bar{\mathbf{r}}(\mathbf{a})$, which are needed to solve the associated linear least-squares problems involving $\widetilde{\mathit{J}}( \mathbf{r}(\mathbf{a}) )$ and $\bar{\mathbf{M}}( \mathbf{a}  )$ in a final step. The first TSQR method is a serial algorithm and the second one is a parallel algorithm. Both of them processes the rows of $\widetilde{\mathit{J}}( \mathbf{r}(\mathbf{a}) )$ or $\bar{\mathbf{M}}( \mathbf{a}  )$ in $n$ steps as described above, but in different order and with different computational kernels as we will see now.

For the sake of convenience, we first describe the serial TSQR, which proceeds the $n$ steps in sequential order from $j=1$ to $n$. Again, to simplify the presentation, but without loss of generality, we will also assume that $\mathbf{W} \in \mathbb{R}^{p \times n}_{+*}$, so that $\widetilde{\mathit{J}}( \mathbf{r}(\mathbf{a}) )$ have $n.p$ rows and each of the $n$ steps processes exactly $p$ rows of  $\widetilde{\mathit{J}}( \mathbf{r}(\mathbf{a}) )$, namely, the $j^{th}$ step processes the submatrix $\widetilde{\mathit{J}}_{j}$ defined above.

In the first step, the submatrix $\widetilde{\mathit{J}}_{1}$ and the subvector $\widetilde{\mathbf{r}}_{1}$ are evaluated exactly as described above. Then, a standard Householder QR algorithm is applied to $\widetilde{\mathit{J}}_{1}$ to transform this target matrix into an upper triangular (if $k=1$) or trapezoidal (if $k>1$) matrix $\widetilde{\mathbf{R}}_{1}$:
\begin{equation*}
\widetilde{\mathit{J}}_{1} = \widetilde{\mathbf{Q}}_{1} \widetilde{\mathbf{R}}_{1}  \ ,
\end{equation*}
where $\widetilde{\mathbf{Q}}_{1}$ is an $p \times p$ orthogonal matrix and $\widetilde{\mathbf{R}}_{1}$ is an $p \times k.p$ upper triangular or trapezoidal matrix. This is performed by the application of a sequence of $p$ Householder transformations, whose product implicitly represents the orthogonal matrix $\widetilde{\mathbf{Q}}_{1}$, see Subsection~\ref{lin_alg:box} for more details. The algorithm consists of the iterations of two steps: generation of the Householder transformation from the target column vector of $\widetilde{\mathit{J}}_{1}$ and application of this Householder transformation to the trailing part of $\widetilde{\mathit{J}}_{1}$. Next, the right hand-side subvector $\widetilde{\mathbf{r}}_{1}$ is pre-multiplied by the transpose of $\widetilde{\mathbf{Q}}_{1}$. Note that, in a practical implementation of this TSQR algorithm, it is convenient to concatenate and store $\widetilde{\mathit{J}}_{1}$ and $\widetilde{\mathbf{r}}_{1}$ in the same matrix array (with $p$ rows and $k.p+1$ columns if $\mathbf{W} \in \mathbb{R}^{p \times n}_{+*}$) so that the matrix-vector product $\widetilde{\mathbf{Q}}_{1}^{T} \widetilde{\mathbf{r}}_{1}$ is directly computed when the Householder transformations are applied to the trailing part of $\widetilde{\mathit{J}}_{1}$ during its QR factorization.

In the second step of the TSQR algorithm, the submatrix $\widetilde{\mathit{J}}_{2}$ and the subvector $\widetilde{\mathbf{r}}_{2}$ are first evaluated and then combined with the results of the first step as follows:
\begin{equation*}
\widetilde{\widetilde{\mathit{J}}}_{2} =  \begin{bmatrix}  \widetilde{\mathbf{R}}_{1}  \\  \widetilde{\mathit{J}}_{2} \end{bmatrix} \text{ and }  \widetilde{\widetilde{\mathbf{r}}}_{2} = \begin{bmatrix}  \widetilde{\mathbf{Q}}_{1}^{T} \widetilde{\mathbf{r}}_{1} \\  \widetilde{\mathbf{r}}_{2} \end{bmatrix}  \ .
\end{equation*}
Then, a structured and thin QR factorization of the matrix $\widetilde{\widetilde{\mathit{J}}}_{2}$ is performed
\begin{equation*}
\widetilde{\widetilde{\mathit{J}}}_{2} = \widetilde{\mathbf{Q}}_{2} \widetilde{\mathbf{R}}_{2}  \ ,
\end{equation*}
where $\widetilde{\mathbf{Q}}_{2}$ is a matrix with min$(2.p,k.p)$ orthonormal columns and $\widetilde{\mathbf{R}}_{2}$ is an upper triangular or trapezoidal matrix. This reduction to triangular or trapezoidal form of $\widetilde{\widetilde{\mathit{J}}}_{2}$ can be accomplished by a special sequence of min$(2.p,k.p)$ Householder transformations in which the $i^{th}$ transformation is designed to annihilate the nonzero subdiagonal elements in the $i^{th}$ column of $\widetilde{\widetilde{\mathit{J}}}_{2}$. Note that no fill-in occurs during this process because the columns of $\widetilde{\widetilde{\mathit{J}}}_{2}$ are reduced from left to right. Finally, the vector $\widetilde{\widetilde{\mathbf{r}}}_{2}$ is pre-multiplied by the transpose of $\widetilde{\mathbf{Q}}_{2}$ and this ends the second step.

At the $j^{th}$ step, the submatrix $\widetilde{\mathit{J}}_{j}$ and subvector $\widetilde{\mathbf{r}}_{j}$ are evaluated and concatenated with the outputs of the $j-1$ step as
\begin{equation*}
\widetilde{\widetilde{\mathit{J}}}_{j} =  \begin{bmatrix}  \widetilde{\mathbf{R}}_{j-1}  \\  \widetilde{\mathit{J}}_{j} \end{bmatrix} \text{ and }  \widetilde{\widetilde{\mathbf{r}}}_{j} = \begin{bmatrix}  \widetilde{\mathbf{Q}}_{j-1}^{T} \widetilde{\mathbf{r}}_{j-1} \\  \widetilde{\mathbf{r}}_{j} \end{bmatrix}  \ ,
\end{equation*}
and a structured and thin QR factorization of $\widetilde{\widetilde{\mathit{J}}}_{j}$ is performed as
\begin{equation*}
\widetilde{\widetilde{\mathit{J}}}_{j} = \widetilde{\mathbf{Q}}_{j} \widetilde{\mathbf{R}}_{j}  \ ,
\end{equation*}
where $\widetilde{\mathbf{Q}}_{j}$ is a matrix with min$(j.p,k.p)$ orthonormal columns and $\widetilde{\mathbf{R}}_{j}$ is an upper triangular or trapezoidal matrix. The vector $\widetilde{\widetilde{\mathbf{r}}}_{j}$ is also pre-multiplied by the transpose of $\widetilde{\mathbf{Q}}_{j}$ after this new structured QR factorization.

Then, the following steps are exactly similar to the $j^{th}$ step and this process continues in blocks of $p$ rows of $\widetilde{\mathit{J}}( \mathbf{r}(\mathbf{a}) )$ until there are no more rows of  $\widetilde{\mathit{J}}( \mathbf{r}(\mathbf{a}) )$ left. In exact arithmetic without roundoff errors, it can be shown that the final triangular factor $\widetilde{\mathbf{R}}_{n}$ obtained by this recursive algorithm is the upper triangular factor $\mathbf{R}_{\mathit{J}}$ of the standard thin QR factorization of $\widetilde{\mathit{J}}( \mathbf{r}(\mathbf{a}) )$ defined in equation~\eqref{eq:QR_tjmat}. Furthermore, the associated right hand-side vector $\widetilde{\mathbf{Q}}_{n}^{T} \widetilde{\widetilde{\mathbf{r}}}_{n}$ in output of the recursive process is also equal the matrix-vector product $\widetilde{\mathbf{Q}}_{\mathit{J}}^{T} \widetilde{\mathbf{r}}(\mathbf{a})$, where $\widetilde{\mathbf{Q}}_{\mathit{J}}$ is also defined in equation~\eqref{eq:QR_tjmat} and $\widetilde{\mathbf{r}}(\mathbf{a})$ in equation~\eqref{eq:tr_vec}.

For the description of the parallel TSQR algorithm, we assume that $t$ processors are available with their own memory. Then, the $n$ columns of the matrices $\mathbf{X}$ and $\mathbf{W}$ and the $n$ steps are distributed equally among these $t$ processors (or eventually such that the rows of $\widetilde{\mathit{J}}( \mathbf{r}(\mathbf{a}) )$ are partitioned equally among the $t$ processors if $\mathbf{W} \in \mathbb{R}^{p \times n}_{+}$).

Next, each processor $i$ processes its own $n_i$ steps independently without any communication between the processors as in the serial TSQR algorithm. The obtained triangular factors and transformed right hand-side vectors, say,
\begin{equation*}
\widetilde{\mathbf{R}(i)} = \widetilde{\mathbf{R}(i)}_{n_i}  \text{ and }   \widetilde{\mathbf{r}(i)}(\mathbf{a}) =  \widetilde{\mathbf{Q}(i)}_{n_i}^{T} \widetilde{\widetilde{\mathbf{r}(i)}}_{n_i}  \text{ for } i=1 \text{ to } t \ ,
\end{equation*}
are then reduced into an unique triangular factor $\mathbf{R}_{\mathit{J}}$ and an unique transformed right hand-side vector $\widetilde{\mathbf{Q}}_{\mathit{J}}^{T} \widetilde{\mathbf{r}}(\mathbf{a})$ by the QR factorizations (in parallel) of a sequence of matrices built by coupling two upper-triangular factors $\widetilde{\mathbf{R}(i)}$ and $\widetilde{\mathbf{R}(j)}$ on top of each other. We also call such special QR factorization, a structured QR factorization, and this QR factorization can also be performed by a special sequence of Householder transformations~\cite{LH1974}\cite{DGHL2012}. In this second part of the parallel TSQR algorithm, we note that several reduction trees are available to obtain the final triangular factor $\mathbf{R}_{\mathit{J}}$ and transformed right hand-side vector $\widetilde{\mathbf{Q}}_{\mathit{J}}^{T} \widetilde{\mathbf{r}}(\mathbf{a})$~\cite{DGHL2012}. But, for simplicity, we further assume that $t$ is a power of two and that a binary reduction tree is used. In other words, the pairs of processors used in the second part of the parallel TSQR algorithm are given by simply grouping together, initially, a processor and its neighbour, e.g., $(0,1), (2,3), \cdots, (t-2,t-1)$. Then, in the next recursion level, we group the left most processor of a pair, e.g., $(0,2), (4,6), \cdots$ and proceed in this fashion until the last pair is $(0,t/2)$ with the result that the triangular factor and right hand-side vector stored in processor $0$ contains $\mathbf{R}_{\mathit{J}}$ and the matrix-vector product $\widetilde{\mathbf{Q}}_{\mathit{J}}^{T} \widetilde{\mathbf{r}}(\mathbf{a})$, e.g., the actual triangular factor in the QR factorization of $\widetilde{\mathit{J}}( \mathbf{r}(\mathbf{a}) )$ and the associated transformed right hand-side vector.

Thus, in the parallel version of the TSQR algorithm, after each processor $i$ has computed its own triangular factor $\widetilde{\mathbf{R}(i)}$, one repeats sending, receiving one's $\widetilde{\mathbf{R}(j)}$ to/from one's processor neighbour and calculating structured QR factorizations in parallel. Then, another round of a reduced number of structured QR factorizations is performed in parallel and this process repeats until the final  $\mathbf{R}_{\mathit{J}}$ factor is obtained. Obviously, other groupings of processors can be used to take advantage of a given topology of a network of processors. As noted in Demmel et al.~\cite{DGHL2012}, any sequence of tree of (structured) QR factorizations between the serial and parallel versions of the TSQR algorithm will work. Our binary tree version shows how to compute the TSQR factorization with maximum parallelism while the serial version does not exhibit any parallelism, but requires less memory. Note, finally, that the two structured QR factorizations used in the serial and parallel TSQR algorithms, namely,
\begin{equation*}
\begin{bmatrix}   \mathbf{R}_1 \\  \mathbf{R}_2 \end{bmatrix} =  \mathbf{Q} \mathbf{R} \quad \text{and}  \quad \begin{bmatrix}  \mathbf{R}_1  \\   \mathbf{X} \end{bmatrix} = \mathbf{Q} \mathbf{R} \ ,
\end{equation*}
where $\mathbf{R}_1, \mathbf{R}_2$ and $\mathbf{R}$ are upper triangular or trapezoidal matrices, $\mathbf{Q}$ is an orthogonal matrix  and $\mathbf{X}$ is a full matrix, can be computed very efficiently using BLAS3 algorithms exploiting the triangular structure of the $\mathbf{R}_1$ and $\mathbf{R}_2$ matrices appearing in these structured QR factorizations~\cite{DGHL2012}. As an illustration, such computational kernels are already available in the recent versions of the LAPACK library.

Furthermore, the orthonormal matrix $\widetilde{\mathbf{Q}}_{\mathit{J}}$ in the thin QR factorization of $\widetilde{\mathit{J}}( \mathbf{r}(\mathbf{a}) )$ (see equation~\eqref{eq:QR_tjmat}) is never explicitly formed in both the serial and parallel TSQR algorithms, but pre-multiplying the vector $\widetilde{\mathbf{r}}(\mathbf{a})$ by $\widetilde{\mathbf{Q}}_{\mathit{J}}^{T}$ can be done also recursively and efficiently during the TSQR algorithms (as illustrated above) and this is sufficient to solve the linear least-squares problems involving $\widetilde{\mathit{J}}( \mathbf{r}(\mathbf{a}) )$ as a coefficient matrix, or as a block of the coefficient matrix to deal with its singularity as we will illustrate below. Obviously, the same recursive TSQR methods can be used to compute the thin QR factorization of $\bar{\mathbf{M}}( \mathbf{a}  )$ if a Kaufman step is used in the Gauss-Newton algorithms~\eqref{gn_alg:box}, but we omit the details here as the steps are essentially the same.

In the previous paragraphs, we show how variable projection WLRA solvers, which request the Jacobian matrix (or its approximation in the case of the Kaufman variant) and perform a QR decomposition of it for solving the linear least-squares problems~\eqref{eq:llsq_j_mat} or~\eqref{eq:llsq_m_mat} at each iteration, can be implemented efficiently and with reduced memory requirements by exploiting the block diagonal structure of $\mathbf{F}(\mathbf{a})$ and using a parallel two-step TSQR algorithm.
We now consider variable projection WLRA solvers, which, at each iteration, solve the linear least-squares problems~\eqref{eq:llsq_j_mat} or~\eqref{eq:llsq_m_mat} by computing the $k.p \times k.p$ cross-product positive semi-definite matrices
\begin{equation*}
\Delta =  \mathit{J}( \mathbf{r}(\mathbf{a}) )^{T} \mathit{J}( \mathbf{r}(\mathbf{a}) )    \text{  or  }    \Lambda =  \mathbf{M}(\mathbf{a} )^{T} \mathbf{M}(\mathbf{a} )
\end{equation*}
and solving the associated normal equations, e.g.,
\begin{equation*}
\Delta d\mathbf{a} =  - \nabla \psi( \mathbf{a} )  \text{  or  }    \Lambda d\mathbf{a} =  - \nabla \psi( \mathbf{a} )  \ ,
\end{equation*}
where $\nabla \psi( \mathbf{a} ) = \mathit{J}( \mathbf{r}(\mathbf{a}) )^{T} \mathbf{r}(\mathbf{a})$.

In order to present in more details this Cholesky approach, we first recall from the results in Section~\ref{jacob:box} that
\begin{equation*}
\mathit{J}( \mathbf{r}(\mathbf{a}) ) = - \big( \mathbf{M}(\mathbf{a}) + \mathbf{L}(\mathbf{a}) \big)
\end{equation*}
and that the columns  of $\mathbf{L}(\mathbf{a})$ lie in $\emph{ran}(  \mathbf{F}(\mathbf{a}) )$, the range of  $\mathbf{F}(\mathbf{a})$, and those of $\mathbf{M}(\mathbf{a})$ lie in $\emph{ran}(  \mathbf{F}(\mathbf{a}) )^{\bot}$. This implies the equalities
\begin{equation*}
\mathbf{M}(\mathbf{a})^{T} \mathbf{L}(\mathbf{a}) = \mathbf{0}^{k.p \times k.p}  \quad  \text{and} \quad  \mathbf{L}(\mathbf{a})^{T} \mathbf{r}(\mathbf{a}) = \mathbf{0}^{k.p}  \ ,
\end{equation*}
from which we deduce that
\begin{equation*}
\Delta =   \mathbf{M}(\mathbf{a})^{T} \mathbf{M}(\mathbf{a}) + \mathbf{L}(\mathbf{a})^{T} \mathbf{L}(\mathbf{a}) \quad  \text{and} \quad \nabla \psi( \mathbf{a} ) = \mathbf{M}( \mathbf{a})^{T} \mathbf{r}(\mathbf{a}) \  .
\end{equation*} 
Thus, $\Delta$ is the sum of two positive semi-definite matrices and $\mathbf{L}(\mathbf{a})$ does not contribute in $\nabla \psi( \mathbf{a} )$, see Section~\ref{hess:box} for details. As in the QR approach, to speed up and parallelize the computations, it is convenient to consider the matrices $-\mathit{J}( \mathbf{r}(\mathbf{a}) )$, $\mathbf{M}(\mathbf{a})$ and  $\mathbf{L}(\mathbf{a})$ as stacks of $n$ submatrices (each of $p$ rows if $\mathbf{W} \in \mathbb{R}^{p \times n}_{+*}$) as in equation~\eqref{eq:JML_mat_blocking}. As before, the $j^{th}$ blocks  $\mathbf{J}_{j}$, $\mathbf{M}_{j}$ and $\mathbf{L}_{j}$ are defined by
\begin{align*}
    \mathbf{J}_{j}  & =  \mathbf{M}_{j} + \mathbf{L}_{j}  \ , \\
    \mathbf{M}_{j} & =  \mathbf{P}_{\mathbf{F}_{j} (\mathbf{a})}^{\bot} \mathbf{U}_{j}  \  , \\
    \mathbf{L}_{j} & =   \mathbf{L}_{j} =  \mathbf{F}_{j}(\mathbf{a})^{-T}  \mathbf{V}_{j}  \ .
\end{align*}
and can be computed from the $j^{th}$ columns of the matrices $\mathbf{X}$ and $\mathbf{W}$. Similarly, it is useful to consider the $n.p$-vector $\mathbf{r}(\mathbf{a})$ as a stack of $n$ subvectors of dimension $p$ (if $\mathbf{W} \in \mathbb{R}^{p \times n}_{+*}$) as in equation~\eqref{eq:res_vec_blocking}. With these block-structures of $-\mathit{J}( \mathbf{r}(\mathbf{a}) )$, $\mathbf{M}(\mathbf{a})$, $\mathbf{L}(\mathbf{a})$ and $\mathbf{r}(\mathbf{a})$, we have the equalities
\begin{equation*}
\Lambda = \mathbf{M}(\mathbf{a})^{T} \mathbf{M}(\mathbf{a}) = \sum_{j=1}^n  \mathbf{M}_{j}^{T} \mathbf{M}_{j} \quad  \text{and} \quad   \mathbf{L}(\mathbf{a})^{T} \mathbf{L}(\mathbf{a}) = \sum_{j=1}^n  \mathbf{L}_{j}^{T} \mathbf{L}_{j}
\end{equation*}
and also
\begin{equation*}
\Delta = \big(  \sum_{j=1}^n  \mathbf{M}_{j}^{T} \mathbf{M}_{j} +  \sum_{j=1}^n  \mathbf{L}_{j}^{T} \mathbf{L}_{j}  \big) \quad  \text{and} \quad   \nabla \psi( \mathbf{a} ) =  \sum_{j=1}^n  \mathbf{M}_{j}^{T}  \mathbf{r}_{j}(\mathbf{a}) \  ,
\end{equation*}
which show that the computations of $\Delta$, $\Lambda$ and $\nabla \psi( \mathbf{a} )$ can be easily parallelized if several processors are available.

As in the QR approach, we can also use the transformed versions of $-\mathit{J}( \mathbf{r}(\mathbf{a}) )$ and $\mathbf{M}(\mathbf{a})$ (see equations~\eqref{eq:Q_mat_blocking},~\eqref{eq:Q_J_mat} and~\eqref{eq:Q_M_mat}) defined as
\begin{equation*}
\widetilde{\mathit{J}} \big ( \mathbf{r}(\mathbf{a}) \big ) =  \mathbf{Q}_{\mathbf{F}}^{T} \mathit{J} \big ( \mathbf{r}(\mathbf{a}) \big ) \quad  \text{and} \quad  \widetilde{\mathbf{M}}( \mathbf{a} ) = \mathbf{Q}_{\mathbf{F}}^{T} \mathbf{M}( \mathbf{a} )
\end{equation*}
since 
\begin{equation*}
\Delta = \widetilde{\mathit{J}} \big( \mathbf{r}(\mathbf{a}) \big)^{T}  \widetilde{\mathit{J}} \big( \mathbf{r}(\mathbf{a}) \big) \quad  \text{and} \quad  \Lambda = \widetilde{\mathbf{M}}( \mathbf{a}  )^{T} \widetilde{\mathbf{M}}( \mathbf{a}  ) = \bar{\mathbf{M}}( \mathbf{a}  )^{T} \bar{\mathbf{M}}( \mathbf{a}  ) \ .
\end{equation*}
Furthermore, taking again into account the respective block-structures of $\widetilde{\mathit{J}}( \mathbf{r}(\mathbf{a}) )$, $\widetilde{\mathbf{M}}( \mathbf{a}  )$  and $\bar{\mathbf{M}}$, defined, respectively, in equations~\eqref{eq:tj_mat},~\eqref{eq:tm_mat} and~\eqref{eq:ttm_mat}, we have the equalities
\begin{equation*}
\Delta =   \sum_{j=1}^n  \widetilde{\mathbf{J}}_{j} ^{T} \widetilde{\mathbf{J}}_{j}   \quad  \text{and} \quad \Lambda =  \sum_{j=1}^n  \widetilde{\mathbf{M}}_{j}^{T}    \widetilde{\mathbf{M}}_{j}  =  \sum_{j=1}^n  \bar{\mathbf{M}}_{j}^{T}  \bar{\mathbf{M}}_{j} \  ,
\end{equation*}
which demonstrate that the evaluation of these transformed forms of $\Delta$  and $\Lambda$ can also be easily parallelized. Similarly, using the block-structures of $\widetilde{\mathbf{r}}(\mathbf{a})$ and $\bar{\mathbf{r}}(\mathbf{a})$ (defined in equations~\eqref{eq:tr_vec} and~\eqref{eq:ttr_vec})  and equations~\eqref{eq:tm_mat},~\eqref{eq:Q_M_mat} and~\eqref{eq:ttm_mat}, we can also express $\nabla \psi( \mathbf{a} )$ as
\begin{equation*}
\nabla \psi( \mathbf{a} ) =  \sum_{j=1}^n  \widetilde{\mathbf{M}}_{j}^{T} \widetilde{\mathbf{r}}_{j}   =  \sum_{j=1}^n  \bar{\mathbf{M}}_{j}^{T} \bar{\mathbf{r}}_{j} \  .
\end{equation*}
This  demonstrates that the evaluation of $\nabla \psi( \mathbf{a} )$ can also be  easily evaluated and parallelized in the normal-equation framework. An alternative approach to evaluate $\nabla \psi( \mathbf{a} )$ in parallel is to use Theorems~\ref{theo4.3:box} and~\ref{theo5.7:box} and the equality
\begin{equation*}
\nabla \psi( \mathbf{a} ) = \frac{ \partial \varphi^{*}( \mathbf{A},\mathbf{\widehat{B}} ) }{ \partial \mathbf{a} } = \mathbf{G}(\mathbf{\widehat{b}})^{T} \mathbf{G}(\mathbf{\widehat{b}}) \mathbf{a} - \mathbf{G}(\mathbf{\widehat{b}})^{T} \mathbf{z} \  ,
\end{equation*}
where $ \mathbf{z} = \emph{vec}\big( ( \sqrt{ \mathbf{W} }  \odot \mathbf{X} )^{T} \big)$, as already indicated in the formal description of the variable projection Gauss-Newton algorithms~\eqref{gn_alg:box} at the beginning of this section. Taking into account the diagonal block-structure of $\mathbf{G}(\mathbf{\widehat{b}})$ (see equation~\eqref{eq:G_mat}), the evaluation of $\nabla \psi( \mathbf{a} )$ using this last formulation can also be parallelized very efficiently.

We now explain how to compute the correction vectors $d\mathbf{a}_{gp-gn}$ and $d\mathbf{a}_{k-gn}$ in the Gauss-Newton algorithms~\eqref{gn_alg:box} using both the TSQR and normal-equation approaches described above.

In the TSQR approach, we compute implicitly a thin QR decomposition of $-\mathit{J}( \mathbf{r}(\mathbf{a}) )$ or $\mathbf{M}(\mathbf{a})$ in two stages. For the Jacobian matrix $\mathit{J}( \mathbf{r}(\mathbf{a}) )$, these two steps are as follows
\begin{itemize}
\item[$(a)$] $- \mathit{J} \big( \mathbf{r}(\mathbf{a}) \big) = \mathbf{Q}_{\mathbf{F}} \widetilde{\mathit{J}} \big( \mathbf{r}(\mathbf{a}) \big)$ ,
\item[$(b)$] $\widetilde{\mathit{J}} \big( \mathbf{r}(\mathbf{a}) \big) = \widetilde{\mathbf{Q}}_{\mathit{J}} \mathbf{R}_{\mathit{J}}$ ,
\end{itemize}
giving the thin QR factorization of $-\mathit{J}( \mathbf{r}(\mathbf{a}) ) = \mathbf{M}(\mathbf{a}) + \mathbf{L}(\mathbf{a})$ as
\begin{equation} \label{eq:QR_jmat}
\mathbf{M}(\mathbf{a}) + \mathbf{L}(\mathbf{a}) =  (\mathbf{Q}_{\mathbf{F}} \widetilde{\mathbf{Q}}_{\mathit{J}}) \mathbf{R}_{\mathit{J}} = \mathbf{Q}_{\mathit{J}} \mathbf{R}_{\mathit{J}} \  ,
\end{equation}
where $ \mathbf{Q}_{\mathit{J}}$ is an $p.n \times k.p$ matrix with orthonormal columns and $\mathbf{R}_{\mathit{J}}$ is an $k.p \times k.p$ upper triangular matrix. Similarly, if we use the approximate Jacobian matrix $-\mathbf{M}(\mathbf{a})$, we have the two steps:
\begin{itemize}
\item[$(a)$] $\mathbf{M}(\mathbf{a}) = \mathbf{Q}_{\mathbf{F}} \widetilde{\mathbf{M}}(\mathbf{a})$ ,
\item[$(b)$] $\widetilde{\mathbf{M}}(\mathbf{a}) = \widetilde{\mathbf{Q}}_{\mathbf{M}} \mathbf{R}_{\mathbf{M}}$ ,
\end{itemize}
giving the thin QR factorization of $\mathbf{M}(\mathbf{a})$ as
\begin{equation}  \label{eq:QR_mmat}
\mathbf{M}(\mathbf{a}) =  (\mathbf{Q}_{\mathbf{F}}  \widetilde{\mathbf{Q}}_{\mathbf{M}} ) \mathbf{R}_{\mathbf{M}} = \mathbf{Q}_{\mathbf{M}} \mathbf{R}_{\mathbf{M}} \  ,
\end{equation}
where, again, $\mathbf{Q}_{\mathbf{M}}$ is an $p.n \times k.p$ matrix with orthonormal columns and $\mathbf{R}_{\mathbf{M}}$ is an $k.p \times k.p$ upper triangular matrix.

As noted above, $\widetilde{\mathbf{M}}(\mathbf{a})$ has ${r}_{\mathbf{F}( \mathbf{a} )}$ zero rows and these zero rows can be eliminated in the second step, giving the following simplified step
\begin{itemize}
\item[$(b')$] $\bar{\mathbf{M}}(\mathbf{a}) = \bar{\mathbf{Q}}_{\mathbf{M}} \mathbf{R}_{\mathbf{M}}$ ,
\end{itemize}
for the computation of the upper triangular factor $\mathbf{R}_{\mathbf{M}}$ in the Kaufman variant of the Gauss-Newton algorithm. Note that the matrices $\mathbf{Q}_{\mathit{J}}$ and $\mathbf{Q}_{\mathbf{M}}$ are never explicitly computed and all that is needed to solve the associated linear least-squares problems~\eqref{eq:llsq_j_mat} and~\eqref{eq:llsq_m_mat}, and compute the correction vectors $d\mathbf{a}_{gp-gn}$ and $d\mathbf{a}_{k-gn}$ are the triangular factor $\mathbf{R}_{\mathit{J}}$ and the matrix-vector product $\mathbf{Q}_{\mathit{J}}^{T} \mathbf{r}(\mathbf{a})$ in the case of the Golub-Pereyra variant or the triangular factor $\mathbf{R}_{\mathbf{M}}$ and the matrix-vector product $\mathbf{Q}_{\mathbf{M}}^{T} \mathbf{r}(\mathbf{a})$ in the case of the Kaufman variant. These two matrix-vector products are given, respectively, by
\begin{equation}  \label{eq:jr_product}
\mathbf{Q}_{\mathit{J}}^{T} \mathbf{r}(\mathbf{a}) = \widetilde{\mathbf{Q}}_{\mathit{J}}^{T}   \mathbf{Q}_{\mathbf{F}}^{T} \mathbf{r}(\mathbf{a}) = \widetilde{\mathbf{Q}}_{\mathit{J}}^{T} \widetilde{\mathbf{r}}(\mathbf{a})
\end{equation}
and
\begin{equation} \label{eq:mr_product}
\mathbf{Q}_{\mathbf{M}}^{T} \mathbf{r}(\mathbf{a}) = \widetilde{\mathbf{Q}}_{\mathbf{M}}^{T}   \mathbf{Q}_{\mathbf{F}}^{T} \mathbf{r}(\mathbf{a}) = \widetilde{\mathbf{Q}}_{\mathbf{M}}^{T} \widetilde{\mathbf{r}}(\mathbf{a}) = \bar{\mathbf{Q}}_{\mathbf{M}} ^{T} \bar{\mathbf{r}}(\mathbf{a}) \ .
\end{equation}
As explained in the previous paragraphs, both $\mathbf{Q}_{\mathit{J}}^{T} \mathbf{r}(\mathbf{a})$ and $\mathbf{Q}_{\mathbf{M}}^{T} \mathbf{r}(\mathbf{a})$ can be computed recursively, as the triangular factors $\mathbf{R}_{\mathit{J}}$ and and $\mathbf{R}_{\mathbf{M}}$, and without explicitly computing the matrices $\mathbf{Q}_{\mathit{J}}$ and $\mathbf{Q}_{\mathbf{M}}$ with the help of the (parallel) TSQR algorithm.

Now, if we assume that the triangular matrices $\mathbf{R}_{\mathit{J}}$ and $\mathbf{R}_{\mathbf{M}}$ are of full rank (which is equivalent to assume that $- \mathit{J}( \mathbf{r}(\mathbf{a}) )$ and $\mathbf{M}(\mathbf{a} )$ are of full column-rank), the unique solutions of the associated linear least-squares problems~\eqref{eq:llsq_j_mat} and~\eqref{eq:llsq_m_mat} are simply found by solving the following upper triangular systems, for $d\mathbf{a} \in  \mathbb{R}^{k.p}$,
\begin{equation*}
\mathbf{R}_{\mathit{J}} d\mathbf{a} = \mathbf{Q}_{\mathit{J}}^{T}\mathbf{r}(\mathbf{a}) \quad \text{and} \quad \mathbf{R}_{\mathbf{M}} d\mathbf{a} = \mathbf{Q}_{\mathbf{M}}^{T}\mathbf{r}(\mathbf{a}) \ ,
\end{equation*}
since the 2-norm is unitarily invariant. Thus, in this case, we obtain
\begin{equation*}
d\mathbf{a}_{gp-gn} = \mathbf{R}_{\mathit{J}}^{-1}  \mathbf{Q}_{\mathit{J}}^{T}\mathbf{r}(\mathbf{a}) \quad \text{and}  \quad d\mathbf{a}_{k-gn} = \mathbf{R}_{\mathbf{M}}^{-1} \mathbf{Q}_{\mathbf{M}}^{T}\mathbf{r}(\mathbf{a}) \ .
\end{equation*}
Of course, in our case,  this simple approach cannot be used as we already know that the triangular factors $\mathbf{R}_{\mathit{J}}$ and $\mathbf{R}_{\mathbf{M}}$ are rank deficient since $- \mathit{J}( \mathbf{r}(\mathbf{a}) )$ and $\mathbf{M}(\mathbf{a})$ are always column rank deficient (see Theorem~\ref{theo5.2:box}). Thus, we will discuss how to proceed in order to compute the minimum  2-norm solutions of these triangular systems (and, thus, of the associated linear least-squares problems~\eqref{eq:llsq_j_mat} and~\eqref{eq:llsq_m_mat}) using the theoretical results of Subsection~\ref{jacob:box} once we have presented how the correction vectors $d\mathbf{a}_{gp-gn}$ and $d\mathbf{a}_{k-gn}$ can be computed in the normal-equation approach.

In the normal-equation approach, the correction vectors  $d\mathbf{a}_{gp-gn}$ and $d\mathbf{a}_{k-gn}$ can be found by computing, respectively, the matrices $\Delta$ and $\Lambda$, as described above, and solving the associated normal equations, namely,
\begin{equation*}
\Delta d\mathbf{a} =  - \nabla \psi( \mathbf{a} )  \quad \text{or}  \quad  \Lambda d\mathbf{a} =  - \nabla \psi( \mathbf{a} ) \ ,
\end{equation*}
where $\nabla \psi( \mathbf{a} ) = \mathit{J}( \mathbf{r}(\mathbf{a}) )^{T} \mathbf{r}(\mathbf{a})$ is also evaluated by one of the methods described in the preceding paragraphs. Since $\Delta$  and $\Lambda$ are symmetric and positive semi-definite matrices, these normal equations are usually solved by computing the Cholesky factorization of $\Delta$  and $\Lambda$, namely,
\begin{equation*}
\Delta =   \mathbf{R}^{T}_{\Delta}   \mathbf{R}_{\Delta}    \quad  \text{or}  \quad  \Lambda =  \mathbf{R}^{T}_{\Lambda}   \mathbf{R}_{\Lambda} \ ,
\end{equation*}
where  $\mathbf{R}_{\Delta}$ and $\mathbf{R}_{\Lambda}$ are $k.p \times k.p$ upper-triangular matrices. Then, if we assume that $\Delta$ and $\Lambda$ are of full rank and, thus, positive definite (which is again equivalent to assume that $\mathit{J}( \mathbf{r}(\mathbf{a}) )$ and $\mathbf{M}(\mathbf{a} )$ are of full column-rank), the normal equations can be solved by backward and forward substitutions using these Cholesky triangular factors. For example, assuming that $\Delta$ is nonsingular, we first solve, for $d\mathbf{a}_{\Delta} \in \mathbb{R}^{k.p}$, the triangular system
\begin{equation*}
\mathbf{R}_{\Delta}^{T} d\mathbf{a}_{\Delta} = - \mathit{J} \big ( \mathbf{r}(\mathbf{a}) \big )^{T} \mathbf{r}(\mathbf{a})  \Longrightarrow   d\mathbf{a}_{\Delta} = - \mathbf{R}_{\Delta}^{-T} \mathit{J} \big( \mathbf{r}(\mathbf{a}) \big)^{T} \mathbf{r}(\mathbf{a}) \ ,
\end{equation*}
and, then, solve for $d\mathbf{a}_{gp-gn} \in \mathbb{R}^{k.p}$,
\begin{equation*}
\mathbf{R}_{\Delta} d\mathbf{a}_{gp-gn} = d\mathbf{a}_{\Delta}   \Longrightarrow   d\mathbf{a}_{gp-gn} = \mathbf{R}_{\Delta}^{-1} d\mathbf{a}_{\Delta}  \ .
\end{equation*}
Note that, up to the sign of the rows of the matrices, we have the equality
\begin{equation*}
\mathbf{R}_{\Delta} = \mathbf{R}_{\mathit{J}}  \ ,
\end{equation*}
and, also up to the sign of the elements of the vectors, the equality
\begin{equation*}
d\mathbf{a}_{\Delta} =  \mathbf{Q}_{\mathit{J}}^{T} \mathbf{r}(\mathbf{a})  \ ,
\end{equation*}
where $\mathbf{Q}_{\mathit{J}}$ and $\mathbf{R}_{\mathit{J}}$ are, respectively, a $n.p \times k.p$ matrix with orthonormal columns and n $k.p \times k.p$ upper triangular matrix, which define the two matrix factors of the QR decomposition of $-\mathit{J}( \mathbf{r}(\mathbf{a}) )$ in equation~\eqref{eq:QR_jmat}. Obviously, similar results are valid for $\Lambda$.
In such conditions, we have, thus, the equivalences
\begin{equation*}
\mathbf{R}_{\Delta} d\mathbf{a}_{gp-gn}  = d\mathbf{a}_{\Delta}  \Longleftrightarrow  \mathbf{R}_{\mathit{J}} d\mathbf{a}_{gp-gn}  =  \mathbf{Q}_{\mathit{J}}^{T} \mathbf{r}(\mathbf{a})
\end{equation*}
and
\begin{equation*}
\mathbf{R}_{\Lambda} d\mathbf{a}_{k-gn}  = d\mathbf{a}_{\Lambda}  \Longleftrightarrow  \mathbf{R}_{\mathbf{M}} d\mathbf{a}_{k-gn}  =  \mathbf{Q}_{\mathbf{M}}^{T} \mathbf{r}(\mathbf{a})  \ ,
\end{equation*}
which establish the equivalence between the QR and Cholesky approaches when we assume that $-\mathit{J}( \mathbf{r}(\mathbf{a}) )$ and $\mathbf{M}(\mathbf{a})$ are of full column rank and that the computations are performed with exact arithmetic without any roundoff errors. However, again, $d\mathbf{a}_{gp-gn}$ and $d\mathbf{a}_{k-gn}$ cannot be found in this simple way with the normal-equation approach as we already know that the symmetric matrices $\Delta$ and $\Lambda$ are only positive semi-definite since $-\mathit{J}( \mathbf{r}(\mathbf{a}) )$ and $\mathbf{M}(\mathbf{a})$ are always rank-deficient and we will also come back to this problem after discussing the respective merits of the QR and normal-equation approaches in more details.

Because $n.p$ is in most cases much larger than $k.p$ (assuming that $p<n$ and $k \ll p$) and the computation of the QR decomposition of the (approximated) Jacobian matrix is the major portion of the time for the variable projection WLRA solvers, which use this QR approach (despite of the use of a parallel TSQR algorithm and BLAS3 kernels for the structured QR decompositions in it), a standard normal-equation approach can be much faster than the QR approach. However, while the normal-equation approach is extremely efficient in terms of work and speed, it may be also numerically unreliable for ill-conditioned linear least-squares problems as it is well know~\cite{LH1974}\cite{GVL1996}\cite{HPS2012}\cite{B2015}. Furthermore, as $-\mathit{J}( \mathbf{r}(\mathbf{a}) )$ and $\mathbf{M}(\mathbf{a})$ are always rank-deficient, but their precise rank is not known if $\mathbf{W} \in \mathbb{R}^{p \times n}_{+}$ and the number of zero weights in $\mathbf{W}$ is large (see Section~\ref{varpro:box}), it is important to have reliable information about the rank of these matrices, which can be obtained theoretically from the triangular factors $\mathbf{R}_{\mathit{J}}$ and $\mathbf{R}_{\mathbf{M}}$  in the QR approach or $\mathbf{R}_{\Delta}$ and $\mathbf{R}_{\Lambda}$ in the normal equations approach. However, the errors in the computation of $\mathbf{R}_{\mathit{J}}$ and $\mathbf{R}_{\mathbf{M}}$ depend on the condition number of $-\mathit{J}( \mathbf{r}(\mathbf{a}) )$ and $\mathbf{M}(\mathbf{a})$, while that of  $\mathbf{R}_{\Delta}$ and $\mathbf{R}_{\Lambda}$ depend on the square of the condition numbers of $-\mathit{J}( \mathbf{r}(\mathbf{a}) )$ and $\mathbf{M}(\mathbf{a})$~\cite{LH1974}\cite{GVL1996}\cite{HPS2012}\cite{B2015}. Thus, the QR approach can still be preferred for stability and accuracy reasons, especially when the number of zero weights in $\mathbf{W}$ is large.

More precisely, from the results of Subsection~\ref{jacob:box}, we know that the matrices $-\mathit{J}( \mathbf{r}(\mathbf{a}) )$ and $\mathbf{M}(\mathbf{a})$ have a rank $r$ at most equal to $(p-k).k$ if $\emph{rank}(\mathbf{A}) = k$ (see Theorem~\ref{theo5.2:box}  and Corollary~\ref{corol5.3:box}; recall also that the condition $\emph{rank}(\mathbf{A}) = k$ is required for the continuity and differentiability of $\mathbf{P}^{\bot}_{\mathbf{F}(.) }$, $\mathbf{F}(.)^{+}$ and $\mathbf{F}(.)^-$ at $\mathbf{a}$), but that $r$ can be smaller than $(p-k).k$ depending on the number of zero weights or missing values in $\mathbf{X}$ (see Theorems~\ref{theo5.5:box} and~\ref{theo5.6:box} for details). Standard tools for solving linear least-squares problems with such deficient matrices are the SVD or the COD, as outlined in Subsection~\ref{lin_alg:box}, which both allow to estimate the solution vectors $d\mathbf{a}_{gp-gn}$ and $d\mathbf{a}_{k-gn}$ of minimum 2-norm of problems~\eqref{eq:llsq_tj_mat} and~\eqref{eq:llsq_tm_mat} with the help of the generalized inverse of $\mathbf{R}_{\mathit{J}}$ or $\mathbf{R}_{\mathbf{M}}$ in the QR approach.

As an illustration, if the SVD of the upper triangular matrix $\mathbf{R}_{\mathit{J}}$ is given by
\begin{equation*}
\mathbf{R}_{\mathit{J}} = \boldsymbol{\mathit{U}_{\mathit{J}}} \boldsymbol{\mathit{D}_{\mathit{J}}} \boldsymbol{\mathit{V}_{\mathit{J}}}^{T} \ ,
\end{equation*}
where $\boldsymbol{\mathit{U}_{\mathit{J}}}$ and $\boldsymbol{\mathit{V}_{\mathit{J}}}$ are $k.p \times k.p$ orthogonal matrices and $\boldsymbol{\mathit{D}_{\mathit{J}}}$ is a $k.p \times k.p$ diagonal matrix with its diagonal elements equal to the singular values of $\mathbf{R}_{\mathit{J}}$ in decreasing order with at least its $k.k$ last diagonal elements  equal to zero according to Theorem~\ref{theo5.2:box}. As stated in equation~\eqref{eq:ginv_svd}, we have
\begin{equation*}
\mathbf{R}_{\mathit{J}}^{+} = \boldsymbol{\mathit{V}_{\mathit{J}}} \boldsymbol{\mathit{D}_{\mathit{J}}}^{+} \boldsymbol{\mathit{U}_{\mathit{J}}}^{T} \ ,
\end{equation*}
where $\lbrack  \boldsymbol{\mathit{D}_{\mathit{J}}}^{+} \rbrack_{ii} = \lbrack  \boldsymbol{\mathit{D}_{\mathit{J}}} \rbrack_{ii}^{-1}$ if $\lbrack  \boldsymbol{\mathit{D}_{\mathit{J}}} \rbrack_{ii}  \ne 0$ and  $\lbrack  \boldsymbol{\mathit{D}_{\mathit{J}}}^{+} \rbrack_{ii} = 0$ if  $\lbrack  \boldsymbol{\mathit{D}_{\mathit{J}}} \rbrack_{ii} = 0$ and the correction vector $d\mathbf{a}_{gp-gn}$ of minimum 2-norm can be computed as
\begin{equation*}
d\mathbf{a}_{gp-gn} = \mathbf{R}_{\mathit{J}}^{+}  \mathbf{Q}_{\mathit{J}}^{T} \mathbf{r}(\mathbf{a})
\end{equation*}
and, similarly, if we use the Kaufman variant of the Gauss-Newton algorithm, the correction vector  $d\mathbf{a}_{k-gn}$ of minimum 2-norm can be estimated as
\begin{equation*}
d\mathbf{a}_{k-gn} = \mathbf{R}_{\mathbf{M}}^{+}  \mathbf{Q}_{\mathbf{M}}^{T} \mathbf{r}(\mathbf{a}) =  \mathbf{R}_{\mathbf{M}}^{+}  \bar{\mathbf{Q}}_{\mathbf{M}}^{T} \bar{\mathbf{r}}(\mathbf{a}) \ .
\end{equation*}
Note that if any  $\lbrack \boldsymbol{\mathit{D}_{\mathit{J}}} \rbrack_{ii}$ or  $\lbrack \boldsymbol{\mathit{D}_{\mathbf{M}}} \rbrack_{ii}$ is small, but non-zero, these computations can be numerically unstable, which makes important to consider approximate methods, which can provide control over the size of $d\mathbf{a}_{gp-gn}$ or $d\mathbf{a}_{k-gn}$. This leads to consider  low rank estimates of $\mathbf{R}_{\mathit{J}}$ and $\mathbf{R}_{\mathbf{M}}$ by considering only their singular values, which are above a suitable threshold $ \nu \in \mathbb{R}_{+*}$ in the SVDs of $\mathbf{R}_{\mathit{J}}$ and $\mathbf{R}_{\mathbf{M}}$, and, finally, in the computations of $\mathbf{R}_{\mathit{J}}^{+}$ and $d\mathbf{a}_{gp-gn}$ or, alternatively, of  $\mathbf{R}_{\mathit{M}}^{+}$ and $d\mathbf{a}_{k-gn}$.

As an illustration, for the Golub-Pereyra variant of the Gauss-Newton algorithm, without a threshold, we will have
\begin{equation*}
d\mathbf{a}_{gp-gn} = \sum_{ \lbrack \boldsymbol{\mathit{D}_{\mathit{J}}} \rbrack_{ii}  > 0} \frac{ \big( \mathbf{Q}_{\mathit{J}}^{T} \mathbf{r}(\mathbf{a})  \big)^{T}  \lbrack \boldsymbol{\mathit{U}_{\mathit{J}}} \rbrack_{.i} }{ \lbrack \boldsymbol{\mathit{D}_{\mathit{J}}} \rbrack_{ii}}   \lbrack \boldsymbol{\mathit{V}_{\mathit{J}}} \rbrack_{.i}  \  ,
\end{equation*}
but, using the threshold $\nu$, we will get
\begin{equation*}
d\mathbf{a}_{gp-gn} = \sum_{ \lbrack \boldsymbol{\mathit{D}_{\mathit{J}}} \rbrack_{ii}  > \nu} \frac{ \big( \mathbf{Q}_{\mathit{J}}^{T} \mathbf{r}(\mathbf{a})  \big)^{T}  \lbrack \boldsymbol{\mathit{U}_{\mathit{J}}} \rbrack_{.i} }{ \lbrack \boldsymbol{\mathit{D}_{\mathit{J}}} \rbrack_{ii}}   \lbrack \boldsymbol{\mathit{V}_{\mathit{J}}} \rbrack_{.i}  \ ,
\end{equation*}
which obviously limits the potential occurrence of large elements in $d\mathbf{a}_{gp-gn}$. Alternatively, we can estimate $\mathbf{R}_{\mathit{J}}^{+}$ and $\mathbf{R}_{\mathbf{M}}^{+}$ by a COD of $\mathbf{R}_{\mathit{J}}$ and $\mathbf{R}_{\mathbf{M}}$ as described in Subsection~\ref{lin_alg:box}. This will be less time consuming, but also less reliable in estimating precisely the ranks of  $-\mathit{J}( \mathbf{r}(\mathbf{a}) )$ and $\mathbf{M}(\mathbf{a})$. However, recent investigations in the context of general NLLS problems suggest that using a COD can even be more reliable than the truncated SVD approach described above, see~\cite{IKP2011} for details.

Similarly, in the normal-equation approach, we can compute the Eigenvalue-Vector Decomposition (EVD) of the positive semi-definite matrices $\Delta$ (or $\Lambda$) and use a truncated EVD to estimate its pseudo-inverses $\Delta^{+}$ (or $\Lambda^{+}$) and, finally, $d\mathbf{a}_{gp-gn}$ (or $d\mathbf{a}_{k-gn}$) as
\begin{equation*}
d\mathbf{a}_{gp-gn} = - \sum_{ \lbrack \boldsymbol{\mathit{D}_{\mathit{J}}} \rbrack_{ii}^{2}  > \nu^{2} } \frac{ \big( \mathit{J}( \mathbf{r}(\mathbf{a}) )^{T} \mathbf{r}(\mathbf{a})  \big)^{T}  \lbrack \boldsymbol{\mathit{V}_{\mathit{J}}} \rbrack_{.i} }{ \lbrack \boldsymbol{\mathit{D}_{\mathit{J}}} \rbrack_{ii}^{2}}   \lbrack \boldsymbol{\mathit{V}_{\mathit{J}}} \rbrack_{.i} \ ,
\end{equation*}
where the EVD of $\Delta$ is given by
\begin{equation*}
\Delta = \boldsymbol{\mathit{V}_{\mathit{J}}} \boldsymbol{\mathit{D}_{\mathit{J}}}^{2} \boldsymbol{\mathit{V}_{\mathit{J}}}^{T} \ ,
\end{equation*}
where the matrices $\boldsymbol{\mathit{V}_{\mathit{J}}}$ and $\boldsymbol{\mathit{D}_{\mathit{J}}}$ have the same meaning as in the SVD of $\mathbf{R}_{\mathit{J}}$.

However, if $p$ is large and the rank $k$ of the WLRA matrix approximation we are seeking is also not small, the dimensions of $\mathbf{R}_{\mathit{J}}$ and $\mathbf{R}_{\mathbf{M}}$ in the TSQR approach, or, $\Delta$ and $\Lambda$ in the normal-equation approach, can be very large. In these conditions, computing a SVD, or even a COD, of  $\mathbf{R}_{\mathit{J}}$ (or $\mathbf{R}_{\mathbf{M}}$) in the TSQR approach or an EVD of $\Delta$ (or $\Lambda$) in the normal-equation approach, at each iteration of the variable projection Gauss-Newton algorithms~\eqref{gn_alg:box}, can be very costly and, consequently, must be avoided as much as possible.

In many practical applications, for example if $\mathbf{W} \in \mathbb{R}_{+*}^{p \times n}$ or if the number of observed values in each column and row of the data matrix $\mathbf{X}$ is larger than the rank $k$ of the matrix approximation we are seeking, we also know that $r = \emph{rank}( \mathit{J}( \mathbf{r}(\mathbf{a}) ) )$  will be exactly equal to $(p-k).k$ with high probability (see Theorem~\ref{theo5.3:box} and Corollaries~\ref{corol5.4:box} and~\ref{corol5.5:box}), the last $k.k$ columns of $\mathit{J}( \mathbf{r}(\mathbf{a}) )$ and $\mathbf{M}(\mathbf{a})$ are linearly dependent upon the first $(p-k).k$ columns of these matrices (see Theorem~\ref{theo5.4:box}) and, finally, that it is easy to compute an orthonormal basis of the null-space of $\mathit{J}( \mathbf{r}(\mathbf{a}) )$ and $\mathbf{M}(\mathbf{a})$ from an orthonormal basis of $\mathbf{A}$ if  $\emph{rank}( \mathbf{A} ) = k$ (see Corollary~\ref{corol5.6:box}). Collectively, these different results suggest that the minimum 2-norm solutions of the linear least-squares problems~\eqref{eq:llsq_j_mat} and~\eqref{eq:llsq_m_mat} involving the rank deficient matrices $-\mathit{J}( \mathbf{r}(\mathbf{a}) )$ and $\mathbf{M}(\mathbf{a})$ can still be found accurately and much more efficiently without resorting to costly techniques like the SVD or a full COD in the TSQR approach or the EVD in the normal-equation framework as already illustrated in Subsection~\ref{jacob:box}.

We first reconsider the TSQR approach in which we want to find the minimum 2-norm solutions $d\mathbf{a}_{gp-gn}$ or $d\mathbf{a}_{k-gn}$ of the rank-deficient, but consistent upper triangular systems,
\begin{equation*}
\mathbf{R}_{\mathit{J}} d\mathbf{a}_{gp-gn} = \mathbf{Q}_{\mathit{J}}^{T}\mathbf{r}(\mathbf{a}) \text{ or } \mathbf{R}_{\mathbf{M}}  d\mathbf{a}_{k-gn} = \mathbf{Q}_{\mathbf{M}}^{T}\mathbf{r}(\mathbf{a})  \ ,
\end{equation*}
assuming that
\begin{equation*}
\emph{rank}( \mathbf{R}_{\mathit{J}} ) = \emph{rank}( \mathit{J}( \mathbf{r}(\mathbf{a}) ) ) = (p-k).k
\end{equation*}
and, similarly, that
\begin{equation*}
\emph{rank}( \mathbf{R}_{\mathbf{M}} ) = \emph{rank}( \mathbf{M}(\mathbf{a}) ) = (p-k).k \  .
\end{equation*}
Under the hypotheses of Theorem~\ref{theo5.4:box}, we know that the first $(p-k).k$ columns of $-\mathit{J}( \mathbf{r}(\mathbf{a}) )$ and $\mathbf{M}(\mathbf{a})$ are linearly independent and that the last $k.k$ columns of these matrices are linearly dependent onto the first $(p-k).k$ columns of these matrices. Obviously, the same relationships hold for $\mathbf{R}_{\mathit{J}}$ and $\mathbf{R}_{\mathbf{M}}$ are these matrices are, respectively, the triangular factors in the QR factorizations of $-\mathit{J}( \mathbf{r}(\mathbf{a}) )$ and $\mathbf{M}(\mathbf{a})$. Then, it is possible to compute the vectors $d\mathbf{a}_{gp-gn}$ or $d\mathbf{a}_{k-gn}$ efficiently as follows.

First, we define the following partitions of the two matrix factors in the thin QR decompositions of $-\mathit{J}( \mathbf{r}(\mathbf{a}) )$ and $\mathbf{M}(\mathbf{a})$, defined in equations~\eqref{eq:QR_jmat}  and~\eqref{eq:QR_mmat}, which correspond to the partitions of $\mathit{J}( \mathbf{r}(\mathbf{a}) )$ and $\mathbf{M}(\mathbf{a})$ used in Theorem~\ref{theo5.4:box}:
\begin{equation*}
\mathbf{Q}_{\mathit{J}} = \begin{bmatrix} \mathbf{Q}_{\mathit{J}}^{1}  & \mathbf{Q}_{\mathit{J}}^{2}  \end{bmatrix} \text{ , }  \mathbf{R}_{\mathit{J}} = \begin{bmatrix} \mathbf{R}_{\mathit{J}}^{11}  & \mathbf{R}_{\mathit{J}}^{12}  \\   \mathbf{0}^{k.k  \times (p-k).k}  & \mathbf{R}_{\mathit{J}}^{22} \end{bmatrix} \ ,
\end{equation*}
and
\begin{equation*}
\mathbf{Q}_{\mathbf{M}} = \begin{bmatrix} \mathbf{Q}_{\mathbf{M}}^{1}  & \mathbf{Q}_{\mathbf{M}}^{2}  \end{bmatrix} \text{ , }  \mathbf{R}_{\mathbf{M}} = \begin{bmatrix} \mathbf{R}_{\mathbf{M}}^{11}  & \mathbf{R}_{\mathbf{M}}^{12}  \\   \mathbf{0}^{k.k  \times (p-k).k}  & \mathbf{R}_{\mathbf{M}}^{22} \end{bmatrix} \ ,
\end{equation*}
where
\begin{itemize}
\item $\mathbf{Q}^{1}_{\mathit{J}}$ and $\mathbf{Q}^{1}_{\mathbf{M}} \in \mathbb{O}^{p.n \times (p-k).k}$ ,
\item $\mathbf{Q}^{2}_{\mathit{J}}$ and $\mathbf{Q}^{2}_{\mathbf{M}}  \in \mathbb{O}^{p.n \times k.k}$ ,
\item $\mathbf{R}_{\mathit{J}}^{11}$ and $\mathbf{R}_{\mathbf{M}}^{11}$ are $(p-k).k \times (p-k).k$ upper triangular matrices ,
\item $\mathbf{R}_{\mathit{J}}^{22}$ and $\mathbf{R}_{\mathbf{M}}^{22}$ are $k.k \times k.k$ upper triangular matrices ,
\item $\mathbf{R}_{\mathit{J}}^{12}$ and $\mathbf{R}_{\mathbf{M}}^{12}$ are $(p-k).k \times k.k$ full matrices .
\end{itemize}
From Theorem~\ref{theo5.4:box}, we deduce immediately that
\begin{equation*}
\mathbf{R}_{\mathit{J}}^{22} = \mathbf{R}_{\mathbf{M}}^{22} = \mathbf{0}^{k.k  \times k.k} \ ,
\end{equation*}
as $\mathbf{R}_{\mathit{J}}$ and $\mathbf{R}_{\mathbf{M}}$ exhibit the same linear dependencies as $\mathit{J}( \mathbf{r}(\mathbf{a}) )$ and $\mathbf{M}(\mathbf{a})$.

Next, applying $(p-k).k$ Householder transformations to the right of $\begin{bmatrix} \mathbf{R}_{\mathit{J}}^{11}  & \mathbf{R}_{\mathit{J}}^{12}  \end{bmatrix}$ and $\begin{bmatrix} \mathbf{R}_{\mathbf{M}}^{11}  & \mathbf{R}_{\mathbf{M}}^{12}  \end{bmatrix}$ to annihilate $\mathbf{R}_{\mathit{J}}^{12}$ and $\mathbf{R}_{\mathbf{M}}^{12}$, we obtain the following simplified CODs of $\mathbf{R}_{\mathit{J}}$ and $\mathbf{R}_{\mathbf{M}}$:
\begin{equation*}
\mathbf{R}_{\mathit{J}} = \begin{bmatrix} \mathbf{R}_{\mathit{J}}^{11}  & \mathbf{R}_{\mathit{J}}^{12}  \\   \mathbf{0}^{k.k  \times (p-k).k}  &\mathbf{0}^{k.k  \times k.k} \end{bmatrix} = \begin{bmatrix} \mathbf{T}_{\mathit{J}}  & \mathbf{0}^{(p-k).k  \times k.k}  \\   \mathbf{0}^{k.k  \times (p-k).k}  & \mathbf{0}^{k.k  \times k.k} \end{bmatrix} \mathbf{Z}_{\mathit{J}}^{T} 
\end{equation*}
and
\begin{equation*}
\mathbf{R}_{\mathbf{M}} = \begin{bmatrix} \mathbf{R}_{\mathbf{M}}^{11}  & \mathbf{R}_{\mathbf{M}}^{12}  \\   \mathbf{0}^{k.k  \times (p-k).k}  &\mathbf{0}^{k.k  \times k.k} \end{bmatrix} = \begin{bmatrix} \mathbf{T}_{\mathbf{M}}  & \mathbf{0}^{(p-k).k  \times k.k}  \\   \mathbf{0}^{k.k  \times (p-k).k}  & \mathbf{0}^{k.k  \times k.k} \end{bmatrix} \mathbf{Z}_{\mathbf{M}}^{T} \ ,
\end{equation*}
where $\mathbf{T}_{\mathit{J}}$ and $\mathbf{T}_{\mathbf{M}}$ are $(p-k).k  \times (p-k).k$ nonsingular upper triangular matrices, and, $\mathbf{Z}_{\mathit{J}}$ and  $\mathbf{Z}_{\mathbf{M}}$ are $k.p  \times k.p$ orthogonal matrices, which are the product of $(p-k).k$ Householder transformations designed to annihilate $\mathbf{R}_{\mathit{J}}^{12}$ and $\mathbf{R}_{\mathbf{M}}^{12}$, respectively.
In doing so, we implicitly obtain the following CODs of $-\mathit{J}( \mathbf{r}(\mathbf{a}) )$ or $\mathbf{M}(\mathbf{a})$ from their thin QR decompositions (defined in equations~\eqref{eq:QR_jmat}  and~\eqref{eq:QR_mmat}) and computed by the TSQR algorithm:
\begin{equation*}
- \mathit{J}( \mathbf{r}(\mathbf{a}) ) =  \mathbf{Q}_{\mathit{J}}  \mathbf{R}_{\mathit{J}} =  \mathbf{Q}_{\mathit{J}} \begin{bmatrix} \mathbf{T}_{\mathit{J}}  & \mathbf{0}^{(p-k).k  \times k.k}  \\   \mathbf{0}^{k.k  \times (p-k).k}  & \mathbf{0}^{k.k  \times k.k} \end{bmatrix} \mathbf{Z}_{\mathit{J}}^{T}
\end{equation*}
and
\begin{equation*}
\mathbf{M}(\mathbf{a}) = \mathbf{Q}_{\mathbf{M}}  \mathbf{R}_{\mathbf{M}} =  \mathbf{Q}_{\mathbf{M}}  \begin{bmatrix} \mathbf{T}_{\mathbf{M}}  & \mathbf{0}^{(p-k).k  \times k.k}  \\   \mathbf{0}^{k.k  \times (p-k).k}  & \mathbf{0}^{k.k  \times k.k} \end{bmatrix} \mathbf{Z}_{\mathbf{M}}^{T} \ ,
\end{equation*}
and we can express $-\mathit{J}( \mathbf{r}(\mathbf{a}) )^{+}$ or $\mathbf{M}(\mathbf{a})^{+}$ as
\begin{equation*}
- \mathit{J}( \mathbf{r}(\mathbf{a}) )^{+} =  \mathbf{Z}_{\mathit{J}}   \begin{bmatrix} \mathbf{T}_{\mathit{J}}^{-1}  & \mathbf{0}^{(p-k).k  \times k.k}  \\   \mathbf{0}^{k.k  \times (p-k).k}  & \mathbf{0}^{k.k  \times k.k} \end{bmatrix}   \mathbf{Q}_{\mathit{J}}^{T} = \mathbf{Z}_{\mathit{J}}   \begin{bmatrix} \mathbf{T}_{\mathit{J}}^{-1} (\mathbf{Q}_{\mathit{J}}^{1})^{T} \\   \mathbf{0}^{k.k  \times (p-k).k} \end{bmatrix} 
\end{equation*}
and
\begin{equation*}
\mathbf{M}(\mathbf{a})^{+} =  \mathbf{Z}_{\mathbf{M}}   \begin{bmatrix} \mathbf{T}_{\mathbf{M}}^{-1}  & \mathbf{0}^{(p-k).k  \times k.k}  \\   \mathbf{0}^{k.k  \times (p-k).k}  & \mathbf{0}^{k.k  \times k.k} \end{bmatrix}   \mathbf{Q}_{\mathbf{M}}^{T} = \mathbf{Z}_{\mathbf{M}}   \begin{bmatrix} \mathbf{T}_{\mathbf{M}}^{-1} (\mathbf{Q}_{\mathbf{M}}^{1})^{T} \\   \mathbf{0}^{k.k  \times (p-k).k} \end{bmatrix} \ .
\end{equation*}
Finally, $d\mathbf{a}_{gp-gn}$ and $d\mathbf{a}_{k-gn}$ can be computed as
\begin{equation*}
d\mathbf{a}_{gp-gn} = - \mathit{J}( \mathbf{r}(\mathbf{a}) )^{+} \mathbf{r}(\mathbf{a}) =  \mathbf{Z}_{\mathit{J}}   \begin{bmatrix} \mathbf{T}_{\mathit{J}}^{-1} (\mathbf{Q}_{\mathit{J}}^{1})^{T} \mathbf{r}(\mathbf{a}) \\   \mathbf{0}^{k.k  \times (p-k).k} \end{bmatrix} = \mathbf{Z}_{\mathit{J}}   \begin{bmatrix} \mathbf{T}_{\mathit{J}}^{-1} (\widetilde{\mathbf{Q}}_{\mathit{J}}^{1})^{T} \widetilde{\mathbf{r}}(\mathbf{a}) \\   \mathbf{0}^{k.k  \times (p-k).k} \end{bmatrix}
\end{equation*}
and
\begin{equation*}
d\mathbf{a}_{k-gn} = \mathbf{M}(\mathbf{a})^{+} \mathbf{r}(\mathbf{a}) = \mathbf{Z}_{\mathbf{M}}   \begin{bmatrix} \mathbf{T}_{\mathbf{M}}^{-1} (\mathbf{Q}_{\mathbf{M}}^{1})^{T}  \mathbf{r}(\mathbf{a}) \\   \mathbf{0}^{k.k  \times (p-k).k} \end{bmatrix} = \mathbf{Z}_{\mathbf{M}}   \begin{bmatrix} \mathbf{T}_{\mathbf{M}}^{-1} (\bar{\mathbf{Q}}_{\mathbf{M}}^{1})^{T}  \bar{\mathbf{r}}(\mathbf{a}) \\   \mathbf{0}^{k.k  \times (p-k).k} \end{bmatrix},
\end{equation*}
where we have used equations~\eqref{eq:jr_product}  and~\eqref{eq:mr_product}, and, in both cases, the matrix expressions on the right hand-side of these equalities are available on output of the TSQR algorithm.

Alternatively, the correction vectors $d\mathbf{a}_{gp-gn}$ and $d\mathbf{a}_{k-gn}$ can be computed with the help of Theorem~\ref{theo5.3:box} and Corollary~\ref{corol5.6:box}, which state that the matrix
\begin{equation*}
\mathbf{N} = \mathbf{K}_{(p,k)} ( \mathbf{I}_{k} \otimes \mathbf{A} )
\end{equation*}
is a matrix of full column rank and that the columns of $\mathbf{N}$ form a basis of  $\emph{null}\big( \mathit{J}( \mathbf{r}(\mathbf{a}) ) \big) = \emph{null}( \mathbf{M}(\mathbf{a}) )$, if $\emph{rank}( \mathbf{A} ) = k$ and the hypotheses of Theorem~\ref{theo5.3:box} are verified. Note that the condition $\emph{rank}( \mathbf{A} ) = k$  is always verified if step \textbf{(0)} of the Gauss-Newton algorithms~\eqref{gn_alg:box} is performed at each iteration. On the other hand, the hypotheses of Theorem~\ref{theo5.3:box} can be violated or, more generally, the condition $\emph{rank}( \mathit{J}( \mathbf{r}(\mathbf{a}) ) ) = \emph{rank}( \mathbf{M}(\mathbf{a}) ) = (p-k).k$ can be not verified, especially, in the case of a very large number of missing values in $\mathbf{X}$ or zero weights in $\mathbf{W}$ as demonstrated in Theorems~\ref{theo5.5:box} and~\ref{theo5.6:box}. However, as noted at the end of Subsection~\ref{jacob:box} , if we restrict the set of WLRA problems by imposing the condition $\displaystyle{ \sum_{l=1}^p  \boldsymbol{\delta}_{lj} } >  k $ for all $ j=1, \cdots, n$, where $\boldsymbol{\delta}$ is the incidence matrix associated with the matrix $\mathbf{X}$, the  condition $\emph{rank}( \mathit{J}( \mathbf{r}(\mathbf{a}) ) ) = \emph{rank}( \mathbf{M}(\mathbf{a}) ) = (p-k).k$ will be also verified in the majority of practical applications. In this scenario, the columns of $\mathbf{N}$ form also an orthonormal basis of $\emph{null}\big( \mathit{J}( \mathbf{r}(\mathbf{a}) ) \big) = \emph{null}( \mathbf{M}(\mathbf{a}) )$ if  step \textbf{(0)} of the Gauss-Newton algorithms~\eqref{gn_alg:box} is used at each iteration according to Corollary~\ref{corol5.6:box}, since $\mathbf{A}$ has orthonormal columns after step \textbf{(0)} is performed. Then, using the results of Section~\ref{jacob:box}, it is not difficult to see that the vector $d\mathbf{a}_{gp-gn}$ is the unique solution of the structured linear system
\begin{equation*}
\begin{bmatrix}   \mathbf{R}_{\mathit{J}} \\  \mathbf{N}^{T}  \end{bmatrix} d\mathbf{a}_{gp-gn} = \begin{bmatrix}  \widetilde{\mathbf{Q}}_{\mathit{J}}^{T} \widetilde{\mathbf{r}}(\mathbf{a}) \\ \mathbf{0}^{k.k  \times k.p} \end{bmatrix} \ ,
\end{equation*}
which can be solved  by computing the structured thin QR decomposition of $\begin{bmatrix}   \mathbf{R}_{\mathit{J}} \\  \mathbf{N}^{T}  \end{bmatrix}$ in a first step as
\begin{equation*}
\begin{bmatrix}   \mathbf{R}_{\mathit{J}} \\  \mathbf{N}^{T}  \end{bmatrix} = \mathbf{Q}_{\mathit{J}}(\mathbf{N}) \mathbf{R}_{\mathit{J}} (\mathbf{N}) \ ,
\end{equation*}
where $\mathbf{Q}_{\mathit{J}}(\mathbf{N})$ is an $k.(p+k) \times k.p$ matrix with orthonormal columns and $\mathbf{R}_{\mathit{J}}(\mathbf{N})$ is an $k.p \times k.p$ nonsingular upper triangular matrix. In these conditions, $d\mathbf{a}_{gp-gn}$ is the unique solution of the upper triangular system
\begin{equation*}
\mathbf{R}_{\mathit{J}}(\mathbf{N})   d\mathbf{a}_{gp-gn} =  \mathbf{Q}_{\mathit{J}}(\mathbf{N})^{T}  \begin{bmatrix}  \widetilde{\mathbf{Q}}_{\mathit{J}}^{T} \widetilde{\mathbf{r}}(\mathbf{a}) \\ \mathbf{0}^{k.k  \times k.p} \end{bmatrix},
\end{equation*}
which can be easily solved by backward substitution. In a similar fashion and using the same notations, $d\mathbf{a}_{k-gn}$ can be evaluated by computing the structured thin QR decomposition of $\begin{bmatrix}   \mathbf{R}_{\mathbf{M}} \\  \mathbf{N}^{T}  \end{bmatrix}$ in a first step
\begin{equation*}
\begin{bmatrix}   \mathbf{R}_{\mathbf{M}} \\  \mathbf{N}^{T}  \end{bmatrix} = \mathbf{Q}_{\mathbf{M}}(\mathbf{N})  \mathbf{R}_{\mathbf{M}} (\mathbf{N}) \ ,
\end{equation*}
and by solving the following nonsingular upper triangular system in the second step
\begin{equation*}
\mathbf{R}_{\mathbf{M}}(\mathbf{N})  d\mathbf{a}_{k-gn} =  \mathbf{Q}_{\mathbf{M}}(\mathbf{N})^{T}  \begin{bmatrix}  \bar{\mathbf{Q}}_{\mathbf{M}}^{T} \bar{\mathbf{r}}(\mathbf{a}) \\ \mathbf{0}^{k.k  \times k.p} \end{bmatrix}.
\end{equation*}
Recall, finally, that solving these upper triangular systems for $d\mathbf{a}_{gp-gn}$ and $d\mathbf{a}_{k-gn}$ in output of the TSQR algorithm is equivalent to find, respectively, the unique solutions of the following "constrained" linear least-squares problems
\begin{equation} \label{eq:j_llsq}
d\mathbf{a}_{gp-gn} = \text{Arg} \min_{d\mathbf{a} \in \mathbb{R}^{p.k}}   \,   \frac{1}{2} \big\Vert \begin{bmatrix} \mathbf{r}(  \mathbf{a} )  \\  \mathbf{0}^{k.k}   \end{bmatrix} - \begin{bmatrix} \mathbf{M}(\mathbf{a}) + \mathbf{L}(\mathbf{a}) \\   \mathbf{N}^{T} \end{bmatrix} d\mathbf{a} \big\Vert^{2}_{2}
\end{equation}
and
\begin{equation} \label{eq:m_llsq}
d\mathbf{a}_{k-gn} = \text{Arg} \min_{d\mathbf{a} \in \mathbb{R}^{p.k}}   \,   \frac{1}{2} \big\Vert \begin{bmatrix} \mathbf{r}(  \mathbf{a} )  \\  \mathbf{0}^{k.k}   \end{bmatrix} - \begin{bmatrix} \mathbf{M}(\mathbf{a})  \\   \mathbf{N}^{T} \end{bmatrix} d\mathbf{a} \big\Vert^{2}_{2},
\end{equation}
in three steps at each iteration of the Golub-Pereyra or Kaufman variants of the Gauss-Newton algorithms~\eqref{gn_alg:box} (as discussed in Section~\ref{jacob:box}).
\\
\begin{remark6.1} \label{remark6.1:box}
A third solution for computing the correction vectors $d\mathbf{a}_{gp-gn}$ and $d\mathbf{a}_{k-gn}$ in the Gauss-Newton algorithms~\eqref{gn_alg:box}, if we assumed again that $\emph{rank}( \mathbf{A} ) = k$ and $\emph{rank}( \mathit{J}( \mathbf{r}(\mathbf{a}) ) ) = \emph{rank}( \mathbf{M}(\mathbf{a}) ) = (p-k).k$, is to apply the TSQR algorithm to the matrices $- \mathit{J}( \mathbf{r}(\mathbf{a}) )\mathbf{\bar{O}}^{\bot}$ and $\mathbf{M}(\mathbf{a})\mathbf{\bar{O}}^{\bot}$ defined in Corollary~\ref{corol5.6:box} instead to $- \mathit{J}( \mathbf{r}(\mathbf{a})$ and $\mathbf{M}(\mathbf{a})$ as described above. In these conditions, the upper triangular matrices obtained in the output of the TSQR algorithm are also nonsingular and can, thus, be directly solved by backward substitution. Finally, the correction vectors $d\mathbf{a}_{gp-gn}$ and $d\mathbf{a}_{k-gn}$ can be computed by a simple matrix-vector product as described in Subsection~\ref{jacob:box}. $\blacksquare$
\end{remark6.1}

Similarly, in the normal-equation approach, several faster alternative methods  can be used to find the correction vectors $d\mathbf{a}_{gp-gn}$ and $d\mathbf{a}_{k-gn}$ from the symmetric positive semi-definite matrices 
\begin{equation*}
\Delta = \mathit{J} \big ( \mathbf{r}(\mathbf{a})  \big )^{T} \mathit{J}  \big ( \mathbf{r}(\mathbf{a})  \big ) = \widetilde{\mathit{J}}  \big ( \mathbf{r}(\mathbf{a})  \big )^{T} \widetilde{\mathit{J}}  \big ( \mathbf{r}(\mathbf{a})  \big )
\end{equation*}
and
\begin{equation*}
\Lambda = \mathbf{M}(\mathbf{a})^{T} \mathbf{M}(\mathbf{a}) = \bar{\mathbf{M}}(\mathbf{a})^{T}  \bar{\mathbf{M}}(\mathbf{a})  \ ,
\end{equation*}
if we assume that $\emph{rank}( \mathbf{A} ) = k$ and $\emph{rank}( \mathit{J}( \mathbf{r}(\mathbf{a}) ) ) = \emph{rank}( \mathbf{M}(\mathbf{a})  ) = (p-k).k$. Remember again that the first condition is always true if step \textbf{(0)} of the Gauss-Newton algorithms~\eqref{gn_alg:box} is performed at each iteration.
Under these hypotheses, we deduce immediately that  $\emph{rank}( \Delta ) = \emph{rank}( \Lambda  ) = (p-k).k$ and to cope with this uniform rank deficiency of $\Delta$ and $\Lambda$, we can again use the results of Corollary~\ref{corol5.6:box}, which state that the matrices $\mathbf{N} = \mathbf{K}_{(p,k)} ( \mathbf{I}_{k} \otimes \mathbf{A} )$ and $\mathbf{\bar{O}} =  \mathbf{K}_{(p,k)} ( \mathbf{I}_{k}  \otimes  \mathbf{O})$ (where $\mathbf{O}$ is an orthonormal basis of $ \emph{ran}( \mathbf{A} ) $) are, respectively, a basis and an orthonormal basis of $\emph{null}\big( \mathit{J}( \mathbf{r}(\mathbf{a}) ) \big) = \emph{null}( \mathbf{M}(\mathbf{a}) )$ and, thus, also of the null spaces of the matrices $\Delta$ and $\Lambda$.

In these conditions, as first suggested by Okatani et al.~\cite{OYD2011}, we can compute
\begin{equation*}
\mathbf{N}\mathbf{N}^{T} =  \mathbf{K}_{(p,k)} ( \mathbf{I}_{k} \otimes \mathbf{A} \mathbf{A}^{T} ) \mathbf{K}_{(k,p)}
\end{equation*}
or
\begin{equation*}
\mathbf{\bar{O}}\mathbf{\bar{O}}^{T} =  \mathbf{K}_{(p,k)} ( \mathbf{I}_{k} \otimes \mathbf{O} \mathbf{O}^{T} ) \mathbf{K}_{(k,p)}  \ ,
\end{equation*}
and add these positive semi-definite matrices of rank $k.k$ to $\Delta$ and $\Lambda$. Again, note that, if step \textbf{(0)} of the Gauss-Newton algorithms~\eqref{gn_alg:box} is performed, we have $\mathbf{N} = \mathbf{\bar{O}}$ and, thus, $\mathbf{N}\mathbf{N}^{T} = \mathbf{\bar{O}}\mathbf{\bar{O}}^{T}$. Next, we can compute
\begin{equation*}
\Delta (\mathbf{N}) = \Delta + \mathbf{N}\mathbf{N}^{T} =  \begin{bmatrix}   \mathit{J}( \mathbf{r}(\mathbf{a}) )   \\  \mathbf{N}^{T}  \end{bmatrix}^{T}  \begin{bmatrix}   \mathit{J}( \mathbf{r}(\mathbf{a}) )   \\  \mathbf{N}^{T}  \end{bmatrix} \text{ or } \Lambda(\mathbf{N}) = \Lambda + \mathbf{N}\mathbf{N}^{T} =  \begin{bmatrix}    \mathbf{M}(\mathbf{a})   \\  \mathbf{N}^{T}  \end{bmatrix}^{T}  \begin{bmatrix}   \mathbf{M}(\mathbf{a})   \\  \mathbf{N}^{T}  \end{bmatrix}.
\end{equation*}
Since $\emph{ran}( \mathbf{N}\mathbf{N}^{T} ) = \emph{ran}(  \Delta )^{\bot} = \emph{ran}(  \Lambda )^{\bot} $ if  $\emph{rank}( \mathit{J}( \mathbf{r}(\mathbf{a}) ) ) = \emph{rank}( \mathbf{M}(\mathbf{a})  ) = (p-k).k$, we have the relationships
\begin{equation*}
\emph{dim}( \Delta (\mathbf{N}) ) = \emph{dim}( \Delta ) + \emph{dim}( \mathbf{N}\mathbf{N}^{T} ) = (p-k).k + k.k = k.p   \  ,
\end{equation*}
and, similarly, $\emph{dim}( \Lambda (\mathbf{N}) ) = k.p$. In other words, the matrices $\Delta (\mathbf{N})$ and $\Lambda(\mathbf{N})$ are positive definite and, thus, of full rank. In these conditions, the normal equations
\begin{equation*}
\Delta (\mathbf{N}) d\mathbf{a}_{gp-gn} = - \nabla \psi( \mathbf{a} )  \text{  or  }\Lambda(\mathbf{N}) d\mathbf{a}_{k-gn} = - \nabla \psi( \mathbf{a} ) 
\end{equation*}
have an unique solution, which are also the solutions of the associated linear least-square problems~\eqref{eq:j_llsq} and~\eqref{eq:m_llsq} solved  in the TSQR approach in order to find $d\mathbf{a}_{gp-gn}$ and $d\mathbf{a}_{k-gn}$. Thus, when using the normal-equation approach, we finally need to compute the Cholesky factorizations of $\Delta (\mathbf{N})$ or $\Lambda(\mathbf{N})$ and solve the above positive definite systems by forward and backward substitutions using these triangular Cholesky factors as described in Okatani et al.~\cite{OYD2011}.
\\
\begin{remark6.2} \label{remark6.2:box}
As for the TSQR method, an alternative solution for computing the correction vectors $d\mathbf{a}_{gp-gn}$ and $d\mathbf{a}_{k-gn}$ in the Gauss-Newton algorithms~\eqref{gn_alg:box}, if we assumed that $\emph{rank}( \mathit{J}( \mathbf{r}(\mathbf{a}) ) ) = \emph{rank}( \mathbf{M}(\mathbf{a}) ) = (p-k).k$, is to apply the normal-equation algorithm to the matrices $-\mathit{J}( \mathbf{r}(\mathbf{a}) )\mathbf{\bar{O}}^{\bot}$ and $\mathbf{M}(\mathbf{a})\mathbf{\bar{O}}^{\bot}$ defined in Corollary~\ref{corol5.6:box} instead to $-\mathit{J}( \mathbf{r}(\mathbf{a})$ and $\mathbf{M}(\mathbf{a})$ as described above. In these conditions, the upper triangular matrices obtained in the output of the Cholesky factorization are nonsingular and can, thus, be directly solved by forward and backward substitutions. Finally, the correction vectors $d\mathbf{a}_{gp-gn}$ and $d\mathbf{a}_{k-gn}$ can also be computed by a simple matrix-vector product as described in Subsection~\ref{jacob:box}. $\blacksquare$
\\
\end{remark6.2}

\subsection{Variable projection Levenberg-Marquardt algorithms} \label{vp_lm_alg:box}

This subsection describes and investigates variable projection Levenberg-Marquardt methods for the solution of the WLRA problem. Using similar notations as in the Gauss-Newton algorithms~\eqref{gn_alg:box} described in the previous subsection, an outline of a first version of the variable projection Levenberg-Marquardt algorithms is as follows:
\\
\begin{lm_alg1} \label{lm_alg1:box}
\end{lm_alg1}
Choose starting matrix $\mathbf{A}_{1} \in \mathbb{R}^{p \times k}$ , $\varepsilon_{1}, \varepsilon_{2}, \varepsilon_{3} , \beta, \Vert \nabla\psi \Vert_{min} \in \mathbb{R}_{+*}$ and $i_{max}, j_{max} \in \mathbb{N}_{*}$, appropriately

\textbf{For} $i=1, 2, \ldots$ \textbf{until convergence do}
\begin{enumerate}
\item[\textbf{(0)}] Optionally, compute a QRCP of $\mathbf{A}_{i}$ (see equation~\eqref{eq:qrcp}) to determine $k_{i} = \emph{rank}( \mathbf{A}_{i} )$ and an orthonormal basis of $\emph{ran}( \mathbf{A}_{i} )$:
\begin{itemize}
\item[] $\mathbf{Q}_{i} \mathbf{A}_{i} \mathbf{P}_{i} =  \begin{bmatrix} \mathbf{R}_{i}    &  \mathbf{S}_{i}   \\  \mathbf{0}^{(p-k_{i}) \times k_{i} }  & \mathbf{0}^{(p-k_{i}) \times (k-k_{i}) }  \end{bmatrix}$ ,
\end{itemize}
where $\mathbf{Q}_{i}$ is an $p \times p$ orthogonal matrix, $\mathbf{P}_{i}$ is an $k \times k$ permutation matrix, $\mathbf{R}_{i}$ is an $k_{i} \times k_{i}$ nonsingular upper triangular matrix (with diagonal elements of decreasing absolute magnitude) and $\mathbf{S}_{i}$ an $k_{i} \times (k-k_{i})$ full matrix, which is vacuous if $k_{i} = k$.

In all cases, compute an $p \times k$  matrix $\mathbf{O}_{i}$  with orthonormal columns as the first $k$ columns of $\mathbf{Q}_{i}$ (i.e., such that $\emph{ran}( \mathbf{A}_{i} ) \subset \emph{ran}( \mathbf{O}_{i} )$ if $k_{i} < k$ and  $\emph{ran}( \mathbf{A}_{i} ) = \emph{ran}( \mathbf{O}_{i} )$ if $k_{i} = k$) and set
\begin{itemize}
\item[] $\mathbf{A}_{i} = \mathbf{O}_{i}$ .
\end{itemize}
This optional orthogonalization step is a safe-guard as the condition $k_{i} = k$ is a necessary condition for the differentiability of $\psi(.)$ at a point $\mathbf{A}_{i}$ and is also to limit the occurrence of overflows and underflows in the next steps by enforcing that the matrix variable $\mathbf{A}_{i} \in \mathbb{O}^{p \times k}$.
\item[\textbf{(1)}] Determine (implicitly) the block diagonal matrix
\begin{itemize}
\item[] $\mathbf{F}(\mathbf{a}_{i}) = \emph{diag}\big( \emph{vec}( \sqrt{\mathbf{W}} ) \big)  \big(  \mathbf{I}_n  \otimes \mathbf{A}_{i}  \big)$ ,
\end{itemize}
where $\mathbf{a}_{i} =  \emph{vec}( \mathbf{A}_{i}^{T} )$.
\item[\textbf{(2)}] Compute (implicitly) a QRCP of $\mathbf{F}(\mathbf{a}_{i})$ to determine $\mathbf{P}^{\bot}_{\mathbf{F}(  \mathbf{a}_{i} ) }$ and  $\mathbf{F}(\mathbf{a}_{i})^-$ (see equations~\eqref{eq:ginv_proj_ortho} and \eqref{eq:sginv_qrcp}) or, alternatively, a COD of  $\mathbf{F}(\mathbf{a}_{i})$ to determine $\mathbf{P}^{\bot}_{\mathbf{F}(  \mathbf{a}_{i} ) }$ and $\mathbf{F}(\mathbf{a}_{i})^{+}$ (see equations~\eqref{eq:ginv_proj_ortho} and \eqref{eq:ginv_cod}).

Note also that $\mathbf{F}(\mathbf{a}_{i})^- = \mathbf{F}(\mathbf{a}_{i})^{+}$ if $\mathbf{F}(\mathbf{a}_{i})$ is of full column rank and that $\mathbf{P}^{\bot}_{\mathbf{F}(  \mathbf{a}_{i} ) }$, $\mathbf{F}(\mathbf{a}_{i})^-$ and $\mathbf{F}(\mathbf{a}_{i})^{+}$ are also block diagonal matrices.
\item[\textbf{(3)}] Solve the block diagonal linear least-squares problem
\begin{itemize}
\item[] $\mathbf{b}_{i} = \text{Arg}\min_{\mathbf{b}\in\mathbb{R}^{k.n}} \,   \Vert \mathbf{x} - \mathbf{F}(\mathbf{a}_{i})\mathbf{b} \Vert^{2}_{2}$ ,
\end{itemize}
e.g., compute
\begin{itemize}
\item[] $\mathbf{b}_{i} =
    \begin{cases}
         \mathbf{F}(\mathbf{a}_{i})^{-} \mathbf{x}   \quad\ \lbrace \text{if a QRCP of $\mathbf{F}(\mathbf{a}_{i})$ is used in step } \textbf{(2)}  \rbrace \\
         \mathbf{F}(\mathbf{a}_{i})^{+} \mathbf{x}  \quad\ \lbrace \text{if a COD of $\mathbf{F}(\mathbf{a}_{i})$ is used in step } \textbf{(2)} \rbrace
    \end{cases}$ .
\end{itemize}
\item[\textbf{(4)}] Determine:
\begin{itemize}
\item[] $\mathbf{r}(\mathbf{a}_{i}) = \mathbf{P}^{\bot}_{\mathbf{F}(  \mathbf{a}_{i} ) } \mathbf{x}$ $\lbrace \text{current residual vector} \rbrace$
\item[] $\psi(\mathbf{a}_{i} ) = \frac{1}{2} \Vert \mathbf{r}(\mathbf{a}_{i}) \Vert^{2}_{2}$ $\lbrace \text{current value of the cost function }  \rbrace$
\item[] $\nabla \psi( \mathbf{a}_{i} ) =  \mathbf{G}(\mathbf{b}_{i})^{T} \mathbf{G}(\mathbf{b}_{i})\mathbf{a}_{i} - \mathbf{G}(\mathbf{b}_{i})^{T} \mathbf{z}$ $\lbrace$see Theorems~\ref{theo4.3:box} and~\ref{theo5.7:box}$\rbrace$
\item[] $\lambda_{i} = \beta \Vert  \nabla \psi( \mathbf{a}_{i} ) \Vert^{2}_{2}$ $\lbrace \text{set ridge parameter proportional to the squared 2-norm of the gradient} \rbrace$
\end{itemize}
Note that the steps \textbf{(1)} to \textbf{(4)} above can be very easily parallelized using the block diagonal structure of $\mathbf{F}(\mathbf{a}_{i})$.
\item[\textbf{(5)}]  Check for convergence. Relevant convergence criteria in the algorithms are of the form:
\begin{itemize}
\item $\Vert \nabla \psi(\mathbf{a}_{i} ) \Vert_{2} \le \varepsilon_{1}$
\item $\Vert \mathbf{a}_{i} - \mathbf{a}_{i-1} \Vert_{2} \le \varepsilon_{2} ( \varepsilon_{2} + \Vert \mathbf{a}_{i} \Vert_{2})$  $\lbrace$if $ i \ne 1\rbrace$

If step \textbf{(0)} is used, this last convergence condition can be simplified as:

$\Vert \mathbf{a}_{i} - \mathbf{a}_{i-1} \Vert_{2} \le \varepsilon_{2}  \Vert \mathbf{a}_{i} \Vert_{2} = \varepsilon_{2} \sqrt{k}$
\item $\vert \psi(\mathbf{a}_{i-1} ) - \psi(\mathbf{a}_{i} ) \vert \le \varepsilon_{3} ( \varepsilon_{3} +  \psi(\mathbf{a}_{i} ) )$   $\lbrace$if $ i \ne 1\rbrace$
\item $ i \ge i_{max}$ $\lbrace \text{e.g., give up if the number of iterations is too large} \rbrace$
\end{itemize}
where $\varepsilon_{1}, \varepsilon_{2}, \varepsilon_{3}$ and $i_{max}$ are constants chosen by the user.

\textbf{Exit if convergence}. \textbf{Otherwise, go to} step \textbf{(6)}
\item[\textbf{(6)}] Compute the Levenberg-Marquardt correction vector $d\mathbf{a}_{lm}$ as the (minimum 2-norm) solution of one of the following (regularized) linear least-squares problems:
\begin{enumerate}
\item[\textbf{(6.1)}] \textbf{If} $\Vert \nabla \psi(\mathbf{a}_{i} ) \Vert_{2} \ge \Vert \nabla\psi \Vert_{min}$ \textbf{then}
\begin{description}
\item[Golub-Pereyra Levenberg-Marquardt step:] Golub and Pereyra~\cite{GP1973}
\begin{align*}
d\mathbf{a}_{gp-lm} & =   \begin{bmatrix}  \mathbf{M}( \mathbf{a}_{i} ) + \mathbf{L}( \mathbf{a}_{i}  ) \\   \sqrt{\lambda_{i}} \mathbf{D}_{i}  \end{bmatrix}^{+}  \begin{bmatrix}  \mathbf{r}( \mathbf{a}_{i} )   \\   \mathbf{0}^{k.p} \end{bmatrix}  \\
                                 & =  \text{Arg}\min_{d\mathbf{a} \in \mathbb{R}^{p.k}} \,   \Vert \mathbf{r}( \mathbf{a}_{i} ) - \big( \mathbf{M}( \mathbf{a}_{i}  ) + \mathbf{L}( \mathbf{a}_{i}  )  \big) d\mathbf{a}   \Vert^{2}_{2} + \lambda_{i} \Vert \mathbf{D}_{i}  d\mathbf{a} \Vert^{2}_{2} 
\end{align*}
\item[Kaufman Levenberg-Marquardt  step:] Kaufman~\cite{K1975}
\begin{align*}
d\mathbf{a}_{k-lm} & =  \begin{bmatrix}  \mathbf{M}( \mathbf{a}_{i} ) \\   \sqrt{\lambda_{i}} \mathbf{D}_{i}  \end{bmatrix}^{+}  \begin{bmatrix}  \mathbf{r}( \mathbf{a}_{i} )   \\   \mathbf{0}^{k.p} \end{bmatrix}  \\
                               & =  \text{Arg}\min_{d\mathbf{a} \in \mathbb{R}^{p.k}} \,  \Vert \mathbf{r}( \mathbf{a}_{i} ) - \mathbf{M}( \mathbf{a}_{i}  ) d\mathbf{a}   \Vert^{2}_{2} + \lambda_{i} \Vert \mathbf{D}_{i}  d\mathbf{a} \Vert^{2}_{2}
\end{align*}
\end{description}
\item[\textbf{(6.2)}] \textbf{Else}
\begin{description}
\item[Golub-Pereyra Gauss-Newton step:] Golub and Pereyra~\cite{GP1973}, Ruhe and Wedin~\cite{RW1980}
\begin{align*}
d\mathbf{a}_{gp-gn} & =   \big( \mathbf{M}( \mathbf{a}_{i}  ) + \mathbf{L}( \mathbf{a}_{i}  ) \big)^{+} \mathbf{r}( \mathbf{a}_{i} )   \\
                                 & = 
    \begin{cases}
        \text{Arg}\min_{d\mathbf{a} \in \mathbb{R}^{p.k}} \,   \Vert d\mathbf{a} \Vert^{2}_{2} \\
        \text{s.t. }  \text{Arg}\min_{d\mathbf{a} \in \mathbb{R}^{p.k}} \,   \Vert \mathbf{r}( \mathbf{a}_{i} ) - \big( \mathbf{M}( \mathbf{a}_{i}  ) + \mathbf{L}( \mathbf{a}_{i}  )  \big) d\mathbf{a}   \Vert^{2}_{2}
    \end{cases}   \\
\end{align*}
\item[Kaufman Gauss-Newton step:] Kaufman~\cite{K1975}, Ruhe and Wedin~\cite{RW1980}
\begin{align*}
d\mathbf{a}_{k-gn} & =  \mathbf{M}( \mathbf{a}_{i}  )^{+} \mathbf{r}( \mathbf{a}_{i} )   \\
                               & =
    \begin{cases}
        \text{Arg}\min_{d\mathbf{a} \in \mathbb{R}^{p.k}} \,   \Vert d\mathbf{a} \Vert^{2}_{2} \\
        \text{s.t. }  \text{Arg}\min_{d\mathbf{a} \in \mathbb{R}^{p.k}} \,   \Vert \mathbf{r}( \mathbf{a}_{i} ) - \mathbf{M}( \mathbf{a}_{i}  ) d\mathbf{a}   \Vert^{2}_{2}
    \end{cases}   \\
\end{align*}
\end{description}
\end{enumerate}
\item[\textbf{(7)}] Increment $\mathbf{a}_{i} = \emph{vec}( \mathbf{A}_{i}^{T} )$, e.g., compute $\mathbf{a}_{i+1} = \emph{vec}( \mathbf{A}_{i+1}^{T} )$ such that $\psi( \mathbf{a}_{i+1} ) < \psi( \mathbf{a}_{i} )$ in order to obtain global convergence.
\begin{enumerate}
\item[\textbf{(7.1)}] To this end, first compute
\begin{itemize}
\item[] $\mathbf{a}_{i+1} = \mathbf{a}_{i} +  d\mathbf{a}_{lm}$
\item[] $\psi(\mathbf{a}_{i+1} ) = \frac{1}{2} \Vert \mathbf{r}(\mathbf{a}_{i+1}) \Vert^{2}_{2} = \frac{1}{2} \Vert \mathbf{P}^{\bot}_{\mathbf{F}(  \mathbf{a}_{i+1} ) } \mathbf{x} \Vert^{2}_{2}$ ,
\end{itemize}
using (implicitly) a QRCP of the block diagonal matrix $\mathbf{F}(  \mathbf{a}_{i+1} )$.
\item[\textbf{(7.2)}] \textbf{If}  $\psi( \mathbf{a}_{i+1} ) > \psi( \mathbf{a}_{i} )$ \textbf{then} recompute $\mathbf{a}_{i+1}$ by one of the following methods:
\begin{description}
\item[Gauss-Seidel:]
$\mathbf{a}_{i+1} = \mathbf{a}_{i} +  d\mathbf{a}_{gs-gn}$ where $d\mathbf{a}_{gs-gn}$ is a Gauss-Seidel step~\cite{RW1980}
\begin{align*}
d\mathbf{a}_{gs-gn} & =  \big( \mathbf{K}_{(n,p)}  \mathbf{G}( \mathbf{b}_{i}  ) \big)^{+} \mathbf{r}( \mathbf{a}_{i} )   \\
                                 & =
    \begin{cases}
        \text{Arg}\min_{d\mathbf{a} \in \mathbb{R}^{p.k}} \,   \Vert d\mathbf{a} \Vert^{2}_{2} \\
        \text{s.t. }  \text{Arg}\min_{d\mathbf{a} \in \mathbb{R}^{p.k}} \,   \Vert \mathbf{r}( \mathbf{a}_{i} ) - \mathbf{K}_{(n,p)}  \mathbf{G}( \mathbf{b}_{i}  ) d\mathbf{a}   \Vert^{2}_{2}
    \end{cases}   \\
\end{align*}
\item[Block alternating least-squares:]
\begin{align*}
  \mathbf{a}_{i+1} & =    \mathbf{G}( \mathbf{b}_{i}  )^{+} \mathbf{z}   \\
                             & = 
    \begin{cases}
        \text{Arg}\min_{\mathbf{a} \in \mathbb{R}^{p.k}} \,   \Vert \mathbf{a} \Vert^{2}_{2} \\
        \text{s.t. }  \text{Arg}\min_{\mathbf{a} \in \mathbb{R}^{p.k}} \,   \Vert \mathbf{z} -  \mathbf{G}( \mathbf{b}_{i}  ) \mathbf{a}   \Vert^{2}_{2}
    \end{cases}   \\
\end{align*}
\item[Line search:]
\begin{equation*}
\mathbf{a}_{i+1} = \mathbf{a}_{i} + \alpha_{i}  d\mathbf{a}_{lm}
\end{equation*}
where $\alpha_i < 1$ is determined by a line search to make the algorithm a descent method (i.e, such that $\psi( \mathbf{a}_{i+1} ) <  \psi( \mathbf{a}_{i} ) $). This is always possible as the correction vector $d\mathbf{a}_{lm}$ is in a descent direction for $\psi(.)$ if $\Vert \nabla \psi(\mathbf{a}_{i} ) \Vert_{2} \ne 0$, see Corollaries~\ref{corol5.7:box} and~\ref{corol5.8:box}.

As an illustration, a simple, but still efficient, strategy is to first shorten the correction step to half the Levenberg-Marquardt length (or Gauss-Newton length if $\Vert \nabla \psi(\mathbf{a}_{i} ) \Vert_{2} < \Vert \nabla\psi \Vert_{min}$), compute the new trial value for $\psi( \mathbf{a}_{i+1} )$ and, if it is still worse, continue to reduce the step until we get a step short enough such that $\psi( \mathbf{a}_{i+1} ) <  \psi( \mathbf{a}_{i} )$. The following loop incorporates this simple step-shortening algorithm:

\textbf{For} $j=1, 2, \ldots$ \textbf{while}$\big( \psi( \mathbf{a}_{i+1} ) > \psi( \mathbf{a}_{i} ) \big)$
\begin{enumerate}
\item[]  $d\mathbf{a}_{lm}  = \frac{1}{2}  d\mathbf{a}_{lm}$
\item[]  $\mathbf{a}_{i+1} = \mathbf{a}_{i} +  d\mathbf{a}_{lm}$
\item[]  $\psi(\mathbf{a}_{i+1} ) =  \frac{1}{2} \Vert \mathbf{P}^{\bot}_{\mathbf{F}(  \mathbf{a}_{i+1} ) } \mathbf{x} \Vert^{2}_{2}$ $\lbrace$using a QRCP of the matrix $\mathbf{F}(  \mathbf{a}_{i+1} ) \rbrace$
\item[]  \textbf{If} $j > j_{max}$ \textbf{exit} $\lbrace \text{give up if the number of iterations is too large} \rbrace$
\end{enumerate}
\textbf{End do}

\end{description}
\end{enumerate}
\end{enumerate}
\textbf{End do}

In this version of the Levenberg-Marquardt algorithms~\eqref{lm_alg1:box}, the shape and  definition of the different vector and matrix variables are exactly the same as in the Gauss-Newton algorithms~\eqref{gn_alg:box} described in the previous subsection. Furthermore, if, during the iterations, $\Vert \nabla \psi(\mathbf{a}_{i} ) \Vert_{2} <  \Vert \nabla\psi \Vert_{min}$ where $\Vert \nabla\psi \Vert_{min}$ is a positive real constant greater than $\varepsilon_{1}$ chosen by the user, we use a Gauss-Newton correction step as also described in the previous subsection. In addition, as in the Gauss-Newton algorithm,  the computations in the above Levenberg-Marquardt algorithms~\eqref{lm_alg1:box} are terminated either when one or several of the convergence criteria listed in step \textbf{(5)} are satisfied, or when the iteration count exceeds the predetermined number $i_{max}$. Obviously, this version~\eqref{lm_alg1:box} of the  Levenberg-Marquardt algorithms is, thus, similar to the Gauss-Newton algorithms~\eqref{gn_alg:box}, except in step \textbf{(6.1)}, when $\Vert \nabla \psi(\mathbf{a}_{i} ) \Vert_{2}  \ge \Vert \nabla\psi \Vert_{min}$ .

This approach is first justified by the fact that, when $\varepsilon_{1} < \Vert \nabla \psi(\mathbf{a}_{i} ) \Vert_{2} < \Vert \nabla\psi \Vert_{min}$, we are near a stationary or local solution point of our minimization problem, in which case, we want to benefit from the faster convergence of the Gauss-Newton method, see Subsection~\ref{opt:box} for more details. On the other hand, if $\Vert \nabla \psi(\mathbf{a}_{i} ) \Vert_{2}  \ge \Vert \nabla\psi \Vert_{min}$, we consider that we are far away from a stationary or solution point  in which case the Gauss-Newton method may be much less satisfactorily and we prefer to use a more robust correction step, which will be more in the steepest descent direction, in order to widen the basin of convergence of the method. With these considerations in mind, when $\Vert \nabla \psi(\mathbf{a}_{i} ) \Vert_{2}  \ge \Vert \nabla\psi \Vert_{min}$, we introduce both a strictly positive damping parameter $\lambda_{i}$ (e.g. the Marquardt parameter), which takes into account how far we are from a solution and, optionally, a strictly positive scaling diagonal matrix  $\mathbf{D}_{i} \in \mathbb{R}^{k.p \times k.p}_{+}$, which may be useful to render the algorithm invariant under diagonal scaling of the solution vector $\hat{\mathbf{a}}$ and even more robust when $\lambda_{i}$ becomes very large as discussed also in Subsection~\ref{opt:box}.

The choice of $\lambda_{i}$ influences both the direction and the size of the correction vector $d\mathbf{a}_{gp-lm}$ and $d\mathbf{a}_{k-lm}$. If $\lambda_{i}$ tends to zero, $d\mathbf{a}_{gp-lm}$ and $d\mathbf{a}_{k-lm}$ will tend, respectively, to the corresponding Gauss-Newton steps $d\mathbf{a}_{gp-gn}$ and $d\mathbf{a}_{k-gn}$. On the other hand, if $\lambda_{i}$ tends to infinity, then $d\mathbf{a}_{gp-lm}$ and $d\mathbf{a}_{k-lm}$ will tend to a short step in the steepest descent direction, e.g., $- \frac{1}{\lambda_{i}} \nabla \psi( \mathbf{a}_{i} )$, see Subsection~\ref{opt:box} for more information. Thus, the choice of the Marquardt parameter $\lambda_{i}$ is based on the following considerations: if we are close to a local solution then we want the faster convergence of the Gauss-Newton method while it is safe to choose the steepest descent method when we are far from the solution. In other words, the selection procedure
\begin{equation*}
\lambda_{i} = \beta \Vert  \nabla \psi( \mathbf{a}_{i} ) \Vert^{2}_{2} \ ,
\end{equation*}
used in step \textbf{(6.1)} of our Levenberg-Marquardt algorithms~\eqref{lm_alg1:box}, is first motivated by the fact that the method of steepest descent has global convergence not held by the Gauss-Newton method. When one is far away from the solution (i.e. $\Vert  \nabla \psi( \mathbf{a}_{i} ) \Vert_{2}$ is large),  $\lambda_{i}$ is chosen to be large in order to weight the descent part of the correction. As the iterates proceed toward the solution (i.e. $\Vert  \nabla \psi( \mathbf{a}_{i} )  \Vert_{2}$ is small), $\lambda_{i}$ is decreased to weight the Gauss-Newton part of the correction. When we are far from the solution we are interested in the stability of the steepest descent method; when we are close, we strive for the rapidity of convergence of the Gauss-Newton method.

However, taking into account the systematic rank deficiency of the Jacobian matrix  $\mathit{J}( \mathbf{r}(\mathbf{a}_{i}) )$ or its Kaufman approximation  $-\mathbf{M}( \mathbf{a}_{i}  )$ demonstrated in the previous sections, we cannot let $\lambda_{i}$  tends to zero freely and we need to control it appropriately in order to avoid numerical instability when computing the correction steps $d\mathbf{a}_{gp-lm}$  or $d\mathbf{a}_{k-lm}$ if $\lambda_{i}$ approaches zero. Thus, in an actual computer implementation, the condition test $\lambda_{i} = 0$ must be replaced by the condition $\lambda_{i} \le \lambda_{min}$ (with $\lambda_{min} \in  \mathbb{R}_{+*}$) to switch to the Gauss-Newton method, where $\lambda_{min}$ is a suitably chosen real constant such that the matrices
\begin{equation*}
\begin{bmatrix}  \mathbf{M}( \mathbf{a}_{i} ) + \mathbf{L}( \mathbf{a}_{i}  ) \\   \sqrt{\lambda_{i}} \mathbf{D}_{i}  \end{bmatrix} \text{ and } \begin{bmatrix}  \mathbf{M}( \mathbf{a}_{i} )  \\   \sqrt{\lambda_{i}} \mathbf{D}_{i}  \end{bmatrix}
\end{equation*}
do not become nearly singular or ill-conditioned when $\lambda_{i}$ approaches zero as it is expected after some iterations of the Levenberg-Marquardt algorithms~\eqref{lm_alg1:box}. Equivalently, in our version~\eqref{lm_alg1:box} of the Levenberg-Marquardt algorithms, such numerical test is also performed in step \textbf{(6.1)}, on $\Vert \nabla \psi(\mathbf{a} ) \Vert_{2}$ using the user defined threshold $\Vert \nabla\psi \Vert_{min}  \in  \mathbb{R}_{+*}$ rather on $\lambda_{i}$. This is justified by the fact that $\lambda_{i} = \beta \Vert \nabla \psi(\mathbf{a}_{i} ) \Vert_{2}$, where $\beta$ is a strictly positive real constant also chosen by the user. Obviously, the choice of $\Vert \nabla\psi \Vert_{min}$ (or alternatively $\lambda_{min}$) can be tricky as it depends on the scaling of the problem. Moreover, it must be done with care to avoid numerical instabilities when computing $d\mathbf{a}_{gp-lm}$ or  $d\mathbf{a}_{k-lm}$, and, at the same time, maintain the global convergence properties of the Levenberg-Marquardt algorithms~\eqref{lm_alg1:box} .

In addition, we have introduced  a strictly positive diagonal matrix  $\mathbf{D}_{i} \in \mathbb{R}^{k.p \times k.p}_{+}$ in step \textbf{(6.1)} of the Levenberg-Marquardt algorithms. A common simple choice for this scaling diagonal matrix  is to set $\mathbf{D}_{i} = \mathbf{I}_{k.p}$, the identity matrix of order $k.p$. This choice together with a suitable strategy to update $\lambda_{i}$ across the iterations gives the original Levenberg algorithm~\cite{L1944}. Note that $\mathbf{D}_{i}$ can also vary during the iterations and permits for example to introduce some scaling in order to take into account the relative sizes of the columns of the Jacobian matrix  or its Kaufman approximation. Thus, as first suggested by Marquardt~\cite{M1963}, we can also set
\begin{equation*}
\big \lbrack \mathbf{D}_{i} \big \rbrack_{jj} =
    \begin{cases}
 \Vert  \big \lbrack  \mathbf{M}( \mathbf{a}_{i}  ) +  \mathbf{L}( \mathbf{a}_{i}  ) \rbrack_{.j}  \Vert_{2}    &\lbrace \text{if  $d\mathbf{a}_{lm} = d\mathbf{a}_{gp-lm}$ in step \textbf{(6.1)}} \rbrace \\
 \Vert  \big \lbrack  \mathbf{M}( \mathbf{a}_{i}  ) \rbrack_{.j}  \Vert_{2}    &\lbrace \text{if  $d\mathbf{a}_{lm} = d\mathbf{a}_{k-lm}$ in step \textbf{(6.1)}} \rbrace
    \end{cases} \ ,
\end{equation*}
for $j = 1, \cdots, k.p$, see Subsection~\ref{opt:box} for further details. Note, however, that this last choice for the scaling matrix $\mathbf{D}_{i}$ implies that the conditions stated in Theorem~\ref{theo5.6:box} are not verified as otherwise some of the elements of the diagonal of $\mathbf{D}_{i}$ will be equal to zero during the whole iterative process.
\\
\\
The Golub-Pereyra step $d\mathbf{a}_{gp-lm}$ corresponds exactly to the standard  Levenberg-Marquardt  step $d\mathbf{a}_{lm}$  applied to the minimization of the variable projection functional  $\psi(.)$, which is introduced in Subsection~\ref{opt:box}. The philosophy behind the Kaufman step  $d\mathbf{a}_{k-lm}$ is exactly similar to the one detailed for the Gauss-Newton algorithms~\eqref{lm_alg1:box}: in most cases, approximating the Jacobian matrix by $-\mathbf{M}( \mathbf{a} )$ can perform even better than to use the exact Jacobian matrix, taking into account the particular form of the Hessian matrix $\nabla^2 \psi( \mathbf{a} )$ derived in Subsection~\ref{hess:box}.

The Golub-Pereyra and Kaufman variants in the Gauss-Newton algorithms~\eqref{gn_alg:box} generate a sequence $\lbrace \mathbf{a}_{i} \rbrace$ by setting $\mathbf{a}_{i+1} = \mathbf{a}_{i} + \alpha_{i} d\mathbf{a}_{i}$, where $d\mathbf{a}_{i}$ is the minimum 2-norm solution of one of the linearized subproblems
\begin{equation*}
d\mathbf{a}_{i} =  \text{Arg}\min_{d\mathbf{a} \in \mathbb{R}^{p.k}} \,  \frac{1}{2} \Vert l( d\mathbf{a} )  \Vert^{2}_{2} \ ,
\end{equation*}
with
\begin{equation} \label{eq:approx_rfunc}
 l( d\mathbf{a} ) = 
     \begin{cases}
  \mathbf{r}( \mathbf{a}_{i} ) - \big( \mathbf{M}( \mathbf{a}_{i}  ) + \mathbf{L}( \mathbf{a}_{i}  )  \big) d\mathbf{a}    &\lbrace  \text{if } d\mathbf{a}_{gp-gn} \text{ is used in step \textbf{(6)}}  \rbrace \\
  \mathbf{r}( \mathbf{a}_{i} ) - \mathbf{M}( \mathbf{a}_{i}  ) d\mathbf{a}  &\lbrace  \text{if } d\mathbf{a}_{k-gn} \text{ is used in step \textbf{(6)}}  \rbrace 
    \end{cases} \ ,
\end{equation}
as explained in the previous subsection. However, we know from the results of the previous sections that a solution of the WLRA problem, if it exists, is never unique and isolated. Furthermore, since we also know that the above linear least-squares subproblems are always rank-deficient, Levenberg-Marquardt methods, which replace them by regularized linearized subproblems of the form
\begin{equation*}
d\mathbf{a}_{i} =  \text{Arg}\min_{d\mathbf{a} \in \mathbb{R}^{p.k}} \,  \frac{1}{2} \Vert l( d\mathbf{a} )  \Vert^{2}_{2}  + \frac{\lambda_{i}}{2}  \Vert \mathbf{D}_{i}  d\mathbf{a} \Vert^{2}_{2} \ ,
\end{equation*}
where $\lambda_{i}$ is a strictly positive parameter and $\mathbf{D}_{i}$ is a (positive) diagonal matrix, are an interesting alternative. Equivalently, this means that the correction step $d\mathbf{a}_{i}$ in the Levenberg-Marquardt approach minimizes one of the following quadratic models
\begin{equation*}
L_{\lambda_{i}}( d\mathbf{a} ) = 
    \begin{cases}
       \psi( \mathbf{a}_{i} ) + d\mathbf{a}^{T}  \nabla \psi( \mathbf{a}_{i} ) +  \frac{1}{2} d\mathbf{a}^{T}  \big( \mathbf{M}( \mathbf{a}_{i}  )^{T} \mathbf{M}( \mathbf{a}_{i}  ) + \mathbf{L}( \mathbf{a}_{i}  )^{T} \mathbf{L}( \mathbf{a}_{i}  )  +  \lambda_{i} \mathbf{D}_{i}^{2}  \big) d\mathbf{a}     \\
       \psi( \mathbf{a}_{i} ) + d\mathbf{a}^{T} \nabla \psi( \mathbf{a}_{i} )  +  \frac{1}{2} d\mathbf{a}^{T}  \big( \mathbf{M}( \mathbf{a}_{i}  )^{T} \mathbf{M}( \mathbf{a}_{i}  )  +  \lambda_{i} \mathbf{D}_{i}^{2} \big) d\mathbf{a}
    \end{cases} \ ,
\end{equation*}
since $\nabla \psi( \mathbf{a}_{i} )  =  - \big( \mathbf{M}( \mathbf{a}_{i}  ) + \mathbf{L}( \mathbf{a}_{i}  ) \big)^{T}  \mathbf{r}( \mathbf{a}_{i} ) =  - \mathbf{M}( \mathbf{a}_{i}  )^{T} \mathbf{r}( \mathbf{a}_{i} )$.
On the other hand, the correction step $d\mathbf{a}_{i}$ in the Gauss-Newton methods is based on the simpler quadratic models
\begin{equation*}
G( d\mathbf{a} ) = L_{0}( d\mathbf{a} ) = 
    \begin{cases}
       \psi( \mathbf{a}_{i} ) + d\mathbf{a}^{T}  \nabla \psi( \mathbf{a}_{i} ) +  \frac{1}{2} d\mathbf{a}^{T}  \big( \mathbf{M}( \mathbf{a}_{i}  )^{T} \mathbf{M}( \mathbf{a}_{i}  ) + \mathbf{L}( \mathbf{a}_{i}  )^{T} \mathbf{L}( \mathbf{a}_{i}  )   \big) d\mathbf{a}     \\
       \psi( \mathbf{a}_{i} ) + d\mathbf{a}^{T} \nabla \psi( \mathbf{a}_{i} )  +  \frac{1}{2} d\mathbf{a}^{T} \mathbf{M}( \mathbf{a}_{i}  )^{T} \mathbf{M}( \mathbf{a}_{i}  )  d\mathbf{a}
    \end{cases} \ .
\end{equation*}
The gradients of the quadratic functions $L_{\lambda_{i}}(.)$ are, respectively,
\begin{equation*}
\nabla L_{\lambda_{i}}( d\mathbf{a} ) =
    \begin{cases}
       \nabla \psi( \mathbf{a}_{i} )   +  \big( \mathbf{M}( \mathbf{a}_{i}  )^{T} \mathbf{M}( \mathbf{a}_{i}  ) + \mathbf{L}( \mathbf{a}_{i}  )^{T} \mathbf{L}( \mathbf{a}_{i}  )  +  \lambda_{i} \mathbf{D}_{i}^{2}  \big) d\mathbf{a}     \\
       \nabla \psi( \mathbf{a}_{i} )     +  \big( \mathbf{M}( \mathbf{a}_{i}  )^{T} \mathbf{M}( \mathbf{a}_{i}  )  +  \lambda_{i} \mathbf{D}_{i}^{2}  \big) d\mathbf{a}
    \end{cases} \ ,
\end{equation*}
and, by setting these gradients equal to zero, we get $d\mathbf{a}_{gp-lm}$ and $d\mathbf{a}_{k-lm}$ as the solutions to the linear systems
\begin{equation*}
     \big( \mathbf{M}( \mathbf{a}_{i}  )^{T} \mathbf{M}( \mathbf{a}_{i}  ) + \mathbf{L}( \mathbf{a}_{i}  )^{T} \mathbf{L}( \mathbf{a}_{i}  )  +  \lambda_{i} \mathbf{D}_{i}^{2}  \big) d\mathbf{a}_{gp-lm}  = - \nabla \psi( \mathbf{a}_{i} ) = \mathbf{M}( \mathbf{a}_{i}  )^{T} \mathbf{r}( \mathbf{a}_{i} )
\end{equation*}
and
\begin{equation*}
     \big( \mathbf{M}( \mathbf{a}_{i}  )^{T} \mathbf{M}( \mathbf{a}_{i}  ) +  \lambda_{i} \mathbf{D}_{i}^{2}  \big) d\mathbf{a}_{k-lm}   = - \nabla \psi( \mathbf{a}_{i} ) = \mathbf{M}( \mathbf{a}_{i}  )^{T} \mathbf{r}( \mathbf{a}_{i} ) \ ,
\end{equation*}
which are, respectively, the normal equations for the damped linear least-squares problems
\begin{equation*}
\min_{d\mathbf{a} \in \mathbb{R}^{p.k}}   \,    \big\Vert  \begin{bmatrix} \mathbf{r}(  \mathbf{a}_{i} ) \\ \mathbf{0}^{p.k} \end{bmatrix}  -  \begin{bmatrix}   \mathbf{M}( \mathbf{a}_{i}  ) + \mathbf{L}( \mathbf{a}_{i}  )    \\  \sqrt{\lambda_{i}} \mathbf{D}_{i}  \end{bmatrix}  d\mathbf{a}  \big\Vert^{2}_{2}
\end{equation*}
and
\begin{equation*}
\min_{d\mathbf{a} \in \mathbb{R}^{p.k}}   \,    \big\Vert  \begin{bmatrix} \mathbf{r}(  \mathbf{a}_{i} ) \\ \mathbf{0}^{p.k} \end{bmatrix}  -  \begin{bmatrix}   \mathbf{M}( \mathbf{a}_{i}  )   \\  \sqrt{\lambda_{i}} \mathbf{D}_{i}  \end{bmatrix}  d\mathbf{a}  \big\Vert^{2}_{2} \ .
\end{equation*}
Furthermore, as the coefficient matrices of the above normal equations are always positive definite if $\lambda_{i}>0$, these linear systems have always an unique solution, which are the global minimizers  of the associated linear least-squares problems or quadratic model functions, and these quadratic functions are also strictly convex. These nice properties are important numerically and are also an another advantage of the Levenberg-Marquardt methods over a simple Gauss-Newton approach. These results also show that the correction vectors $d\mathbf{a}_{gp-lm}$ and $d\mathbf{a}_{k-lm}$ can be computed, alternatively, by a normal-equation or a more stable QR method as for the Gauss-Newton correction vectors $d\mathbf{a}_{gp-gn}$ and $d\mathbf{a}_{k-gn}$ and we will discuss this matter in more details below after we derive a second and third versions of the Levenberg-Marquardt algorithms.

One disadvantage with the simple strategy used in the Levenberg-Marquardt algorithms~\eqref{lm_alg1:box} for updating the Marquardt parameter $\lambda_{i}$ is, however, that strict descent (i.e., $\psi( \mathbf{a}_{i+1}) < \psi( \mathbf{a}_{i})$) is not guaranteed if a line search or alternative strategies for computing $\mathbf{a}_{i+1}$ are not incorporated in step \textbf{(7.2)} of the algorithms. However, as for the Gauss-Newton algorithms~\eqref{gn_alg:box}, in order to implement a line search algorithm, we have to perform the second part of step \textbf{(4)} of the algorithm, every time we want to get $\psi( \mathbf{a}_{i+1} )$ for a new trial value of $\alpha_{i}$ since
\begin{equation*}
\psi(\mathbf{a}_{i+1} ) = \frac{1}{2} \Vert \mathbf{P}^{\bot}_{\mathbf{F}(  \mathbf{a}_{i+1} ) }\mathbf{x} \Vert^{2}_{2}  \ ,
\end{equation*}
and this line search can involve many extra evaluations of $\psi(.)$, which do not get us closer to an acceptable solution. In these conditions, it is again tempting to perform one or several iterations with the fast Gauss-Seidel or  block ALS methods to compute $\mathbf{a}_{i+1}$  in step \textbf{(7.2)} in case we have $\psi( \mathbf{a}_{i} + d\mathbf{a}_{lm} ) > \psi( \mathbf{a}_{i})$ instead of using a more costly line search algorithm. In other words, if a full Levenberg-Marquardt step gives a sufficient decrease of $\psi(.)$, we accept the point $\mathbf{a}_{i} + d\mathbf{a}_{lm}$ as the new iterate. Otherwise we switch to the fast Gauss-Seidel or block ALS methods.

Alternatively, it is well-known that a line search can be completely avoided in Levenberg-Marquardt methods by using a more sophisticated strategy for updating $\lambda$  during the iterations since the choice of the Marquardt parameter influences both the direction and the size of the correction vector $d\mathbf{a}_{lm}$~\cite{NW2006}\cite{MN2010}.
Furthermore, it is always possible to find a $\lambda$ such that $\psi( \mathbf{a}_{i} + d\mathbf{a}_{lm} ) < \psi( \mathbf{a}_{i})$ if $\Vert  \nabla \psi( \mathbf{a}_{i} ) \Vert_{2} \ne 0$~\cite{NW2006}\cite{MN2010}. Thus, by a proper adjustment of the damping parameter $\lambda$ we have also a direct method for ensuring the descent condition $\psi( \mathbf{a}_{i} + d\mathbf{a}_{lm} ) < \psi( \mathbf{a}_{i})$.

As first suggested by Marquardt~\cite{M1963}, one such strategy is to start with $\lambda$ sets at a small value, $10^{-8}$ for example. Whenever a step is unsuccessful, $\lambda$ gets multiplied by $10$ to force smaller steps until $\psi( \mathbf{a}_{i} + d\mathbf{a}_{lm} )< \psi( \mathbf{a}_{i})$. On the other hand, when the steps become successful $\lambda$ is divided by $10$. This simple strategy results in a fully adaptive technique that behaves just like Gauss-Newton when Gauss-Newton is successful, but shifts in the steepest descent direction and shortens steps when the steps are not successful.

More sophisticated strategies for updating the Marquardt parameter  $\lambda$ during the iterations are based on the so-called gain factor
\begin{equation*}
\rho = \frac{ \psi( \mathbf{a}_{i}) - \psi( \mathbf{a}_{i} + d\mathbf{a} )}{ G(\mathbf{0}^{k.p}) -  G( d\mathbf{a} ) } = \frac{ \psi( \mathbf{a}_{i}) - \psi( \mathbf{a}_{i} + d\mathbf{a} )}{ L_{0}(\mathbf{0}^{k.p}) -  L_{0}( d\mathbf{a} ) } \ ,
\end{equation*}
where
\begin{equation*}
G( d\mathbf{a} ) = L_{0}( d\mathbf{a} ) = \frac{1}{2} l( d\mathbf{a} )^{T} l( d\mathbf{a} ) \ ,
\end{equation*}
and $l(.)$ is defined in equation~\eqref{eq:approx_rfunc}, see~\cite{O1976}\cite{M1978}\cite{NW2006}\cite{MN2010} for a discussion of this strategy in a general NLLS context. $G( d\mathbf{a} )$ is assumed to be a good approximation to $\psi( \mathbf{a}_{i} + d\mathbf{a} )$ when $d\mathbf{a}$ is sufficiently small since $l( d\mathbf{a} )$ is based on first order Taylor's expansions for the residual function $\mathbf{r}( \mathbf{a} )$ around the current iterate $\mathbf{a}_{i}$.
Note further that $G( d\mathbf{a} ) = L_{0}( d\mathbf{a} )$ is the quadratic model, which is used to approximate $\psi( \mathbf{a})$ in the neighborhood of the current iterate $\mathbf{a}_{i}$ and is minimized at each iteration of the variable projection Gauss-Newton algorithms~\eqref{gn_alg:box} as discussed above. Now, we have for the Kaufman variant of the variable projection Levenberg-Marquardt  algorithm the equalities
\begin{align*}
G(\mathbf{0}^{k.p}) -  G( d\mathbf{a}_{k-lm} ) & =   \frac{1}{2} l( \mathbf{0}^{k.p} )^{T} l( \mathbf{0}^{k.p} )  - \frac{1}{2} l( d\mathbf{a}_{k-lm} )^{T} l( d\mathbf{a}_{k-lm} ) \\
                                                                        & =    \psi( \mathbf{a}_{i} ) - \frac{1}{2} \big(  \mathbf{r}( \mathbf{a}_{i} ) - \mathbf{M}( \mathbf{a}_{i}  ) d\mathbf{a}_{k-lm} \big)^{T} \big(  \mathbf{r}( \mathbf{a}_{i} ) - \mathbf{M}( \mathbf{a}_{i}  ) d\mathbf{a}_{k-lm} \big) \\
                                                                    & =   - d\mathbf{a}_{k-lm}^{T}  \nabla \psi( \mathbf{a}_{i} ) -  \frac{1}{2} d\mathbf{a}_{k-lm}^{T} \mathbf{M}( \mathbf{a}_{i}  )^{T} \mathbf{M}( \mathbf{a}_{i}  ) d\mathbf{a}_{k-lm} \\
                                                                    & =    - \frac{1}{2} d\mathbf{a}_{k-lm}^{T}  \big( 2 \nabla \psi( \mathbf{a}_{i} )  +   \mathbf{M}( \mathbf{a}_{i}  )^{T} \mathbf{M}( \mathbf{a}_{i}  ) d\mathbf{a}_{k-lm}  \big)  \\
                                                                    & =    - \frac{1}{2} d\mathbf{a}_{k-lm}^{T}  \Big( 2 \nabla \psi( \mathbf{a}_{i} )  +  \big(  \mathbf{M}( \mathbf{a}_{i}  )^{T} \mathbf{M}( \mathbf{a}_{i}  ) +  \lambda_{i} \mathbf{D}_{i}^{2}  -  \lambda_{i}  \mathbf{D}_{i}^{2}  \big) d\mathbf{a}_{k-lm}  \Big)  \\
                                                                    & =  - \frac{1}{2} d\mathbf{a}_{k-lm}^{T}  \Big( 2 \nabla \psi( \mathbf{a}_{i} )  - \nabla \psi( \mathbf{a}_{i} ) - \lambda_{i} \mathbf{D}_{i}^{2} d\mathbf{a}_{k-lm}  \Big)  \\
                                                                    & =  \frac{1}{2} d\mathbf{a}_{k-lm}^{T}  \Big( \lambda_{i} \mathbf{D}_{i}^{2} d\mathbf{a}_{k-lm} - \nabla \psi( \mathbf{a}_{i} ) \Big)  \\
                                                                    & =   \frac{1}{2} \Big(  \lambda_{i} \Vert \mathbf{D}_{i}  d\mathbf{a}_{k-lm} \Vert^{2}_{2} - d\mathbf{a}_{k-lm}^{T} \nabla \psi( \mathbf{a}_{i} )   \Big) \ ,
\end{align*}
where we have used the equalities
\begin{equation*}
\nabla \psi( \mathbf{a}_{i} ) = - \mathbf{M}( \mathbf{a}_{i}  )^{T} \mathbf{r}( \mathbf{a}_{i} )  \text{ and } \big(  \mathbf{M}( \mathbf{a}_{i}  )^{T} \mathbf{M}( \mathbf{a}_{i}  ) +  \lambda_{i} \mathbf{D}_{i}^{2}   \big) d\mathbf{a}_{k-lm} = - \nabla \psi( \mathbf{a}_{i} ) \ ,
\end{equation*}
derived in Theorem~\ref{theo5.7:box} and Corollary~\ref{corol5.8:box}.
Moreover, the same equality holds for the Golub-Pereyra variant since
\begin{equation*}
\emph{ran}( \mathbf{M}(\mathbf{a}_{i}) ) \subset \emph{ran}( \mathbf{F}( \mathbf{a}_{i} ) )^{\bot} \text{ and } \emph{ran}( \mathbf{L}(\mathbf{a}_{i}) )  \subset \emph{ran}( \mathbf{F}( \mathbf{a}_{i} ) )  \ .
\end{equation*}
Thus, in both cases, $G(\mathbf{0}^{k.p}) -  G( d\mathbf{a}_{lm} )$ and the gain factor $\rho$ can be easily computed at each iteration of the Levenberg-Marquardt  algorithms as the terms  $\lambda_{i} \Vert \mathbf{D}_{i}  d\mathbf{a}_{lm} \Vert^{2}_{2}$ and $\nabla \psi( \mathbf{a}_{i} )$ are already available before the gain factor must be evaluated.

Furthermore, $G(\mathbf{0}^{k.p}) - G( d\mathbf{a}_{lm} )$ is always guaranteed to be positive as  both $\Vert \mathbf{D}_{i}  d\mathbf{a}_{lm} \Vert^{2}_{2}$ and $- d\mathbf{a}_{lm}^{T} \nabla \psi( \mathbf{a}_{i} )$ are positive if $\Vert \nabla \psi( \mathbf{a}_{i} ) \Vert_{2} \ne 0$ since $d\mathbf{a}_{lm}$ is in a descent direction for $\psi(.)$, as demonstrated in Corollary~\ref{corol5.7:box}. In these conditions, it follows that the condition $\rho > 0$ is equivalent to the descending condition $\psi( \mathbf{a}_{i} + d\mathbf{a}_{lm} ) < \psi( \mathbf{a}_{i})$ if $\Vert \nabla \psi( \mathbf{a}_{i} ) \Vert_{2}  \ne 0$ for both the Golub-Pereyra and Kaufman variants of the Levenberg-Marquardt algorithm. Using these results, the following clever strategy proposed by Madsen and Nielsen~\cite{MN2010} may be used to update $\lambda$ at each iteration of Levenberg-Marquardt methods:

\textbf{For} $i=1, 2, \ldots$ \textbf{until convergence do}
\begin{enumerate}
\item[\textbf{(0)}]   $\cdots$
\item[]  $\vdots$
\item[\textbf{(6)}]  Compute the Levenberg-Marquardt correction vector $d\mathbf{a}_{lm}$ as the solution of one of the following constrained and damped linear least-squares problems:
\begin{description}
\item[Golub-Pereyra Levenberg-Marquardt step:]
\begin{align*}
d\mathbf{a}_{gp-lm} & =  \text{Arg}\min_{d\mathbf{a} \in \mathbb{R}^{p.k}} \,   \Vert \mathbf{r}( \mathbf{a}_{i} ) - \big( \mathbf{M}( \mathbf{a}_{i}  ) + \mathbf{L}( \mathbf{a}_{i}  )  \big) d\mathbf{a}   \Vert^{2}_{2} +   \lambda \Vert \mathbf{D}_{i}  d\mathbf{a} \Vert^{2}_{2}  \ ,
\end{align*}
\item[Kaufman Levenberg-Marquardt  step:]
\begin{align*}
d\mathbf{a}_{k-lm}  & =  \text{Arg}\min_{d\mathbf{a} \in \mathbb{R}^{p.k}} \,  \Vert \mathbf{r}( \mathbf{a}_{i} ) - \mathbf{M}( \mathbf{a}_{i}  ) d\mathbf{a}  \Vert^{2}_{2}  + \lambda \Vert \mathbf{D}_{i}  d\mathbf{a} \Vert^{2}_{2}  \ ,
\end{align*}
\end{description}
\item[\textbf{(7)}]  Increment $\mathbf{a}_{i} = \emph{vec}( \mathbf{A}_{i}^{T} )$, e.g., compute $\mathbf{a}_{i+1} = \emph{vec}( \mathbf{A}_{i+1}^{T} )$ such that $\psi( \mathbf{a}_{i+1} ) < \psi( \mathbf{a}_{i} )$ in order to obtain global convergence. To this end, first compute the gain factor
\begin{itemize}
\item[] $\rho = \frac{ \psi( \mathbf{a}_{i}) - \psi( \mathbf{a}_{i} + d\mathbf{a}_{lm} )}{ G(\mathbf{0}^{k.p}) -  G( d\mathbf{a}_{lm} ) } =  \frac{ \psi( \mathbf{a}_{i}) - \psi( \mathbf{a}_{i} + d\mathbf{a}_{lm} )}{  \frac{1}{2}  \big(  \lambda \Vert \mathbf{D}_{i}  d\mathbf{a}_{lm} \Vert^{2}_{2} - d\mathbf{a}_{lm}^{T} \nabla \psi( \mathbf{a}_{i} )  \big) } $
\end{itemize}
\item[] \textbf{If}  $\rho > 0$ \textbf{then}
\begin{itemize}
\item[] $\mathbf{a}_{i+1} = \mathbf{a}_{i} + d\mathbf{a}_{lm}$
\item[] $\lambda = \lambda .\emph{max}\big( \frac{1}{3}, 1-(2.\rho-1)^3 \big)$
\item[] $\nu = 2$
\end{itemize}
\item[] \textbf{Else}
\begin{itemize}
\item[] $\lambda = \nu.\lambda$
\item[] $\nu = 2.\nu$
\item[] \textbf{Go to} step \textbf{(6)} 
\end{itemize}
\end{enumerate}
\textbf{End do}

In this algorithm, the factor $\nu$ is initialized to $2$ and the Marquardt parameter $\lambda$ is again initialized to a small value like $10^{-8}$ (see below for more details).
With this updating strategy, if $\rho \le 0$ then $\mathbf{A}_{i}$ is kept fixed, but we increase $\lambda$ quickly with the twofold purpose of getting closer to the steepest descent direction and reduce the step length. On the other hand, if $\rho > 0$ we accept the new point $\mathbf{a}_{i+1} = \mathbf{a}_{i} + d\mathbf{a}_{lm}$. However, if $\rho$ is small, the quadratic model $G(d\mathbf{a})$ is considered not to be a good approximation to the function $\psi(\mathbf{a}_{i} + d\mathbf{a})$ in the neighborhood of $\mathbf{a}_{i}$ and $\lambda$ is increased. On the other hand, if $\rho$ is large the quadratic model $G(d\mathbf{a})$ is considered to be a good approximation to the function $\psi(\mathbf{a}_{i} + d\mathbf{a})$ and $\lambda$ is decreased in order to get closer to the Gauss-Newton direction in the next iteration step. Thus, this more sophisticated strategy results also in a fully adaptive technique that behaves just like Gauss-Newton when it is successful, but shifts smoothly in the steepest descent direction and shortens steps when the steps are not successful~\cite{MN2010}.

However, in our specific WLRA context, both updating schemes of the Marquardt  parameter must be adapted to take care of the systematic rank deficiency of the coefficient matrices $\mathbf{M}( \mathbf{a}_{i} ) + \mathbf{L}( \mathbf{a}_{i}  )$ or $\mathbf{M}( \mathbf{a}_{i} )$ across the iterations and this without using a specific threshold like $\lambda_{min}$ or $\Vert \nabla\psi \Vert_{min}$ to shift to a Gauss-Newton step when $\lambda$ approaches zero as this will break the incremental nature of the updating schemes if we set $\lambda = 0$ when a full Gauss-Newton step is used at one particular iteration. Furthermore, as before, we must also avoid the near singularity and ill-conditioning of the regularized coefficient matrices
\begin{equation*}
\begin{bmatrix}   \mathbf{M}( \mathbf{a}_{i}  ) + \mathbf{L}( \mathbf{a}_{i}  )    \\  \sqrt{\lambda} \mathbf{D}_{i}  \end{bmatrix}  \quad \text{and}  \quad
 \begin{bmatrix}   \mathbf{M}( \mathbf{a}_{i}  )   \\  \sqrt{\lambda} \mathbf{D}_{i}  \end{bmatrix}  \ ,
\end{equation*}
when $\lambda$ approaches zero.

A clever way to avoid  the use of a threshold and to deal efficiently with the uniform rank deficiency of the Jacobian matrix $\mathit{J}( \mathbf{r}(\mathbf{a}_{i}) )$ or its Kaufman approximation $-\mathbf{M}(\mathbf{a}_{i})$, if we want to incorporate these updating schemes in the variable projection Levenberg-Marquardt algorithms, is to solve at step \textbf{(6)} the constrained and regularized problems
\begin{align*}
d\mathbf{a}_{gp-lm} & =   \text{Arg}\min_{d\mathbf{a} \in \mathbb{R}^{p.k}} \,   \Vert \mathbf{r}( \mathbf{a}_{i} ) - \big( \mathbf{M}( \mathbf{a}_{i}  ) + \mathbf{L}( \mathbf{a}_{i}  )  \big) d\mathbf{a}   \Vert^{2}_{2} +   \Vert \mathbf{N}_{i}^{T}  d\mathbf{a} \Vert^{2}_{2} + \lambda_{i} \Vert \mathbf{D}_{i}  d\mathbf{a} \Vert^{2}_{2}
\end{align*}
or
\begin{align*}
d\mathbf{a}_{k-lm}  & =  \text{Arg}\min_{d\mathbf{a} \in \mathbb{R}^{p.k}} \,  \Vert \mathbf{r}( \mathbf{a}_{i} ) - \mathbf{M}( \mathbf{a}_{i}  ) d\mathbf{a}  \Vert^{2}_{2}  +   \Vert \mathbf{N}_{i}^{T}  d\mathbf{a} \Vert^{2}_{2} + \lambda_{i} \Vert \mathbf{D}_{i}  d\mathbf{a} \Vert^{2}_{2}  \ ,
\end{align*}
where the columns of $\mathbf{N}_{i} = \mathbf{K}_{(p,k)} ( \mathbf{I}_{k} \otimes \mathbf{A}_{i} )$ are a (orthonormal) basis of  $\emph{null}\big( \mathbf{M}(\mathbf{a}_{i}) + \mathbf{L}(\mathbf{a}_{i}) \big) = \emph{null}( \mathbf{M}(\mathbf{a}_{i}) )$ if  $\emph{rank}(\mathbf{A}_{i}) = k$ and $\emph{rank}( \mathit{J}( \mathbf{r}(\mathbf{a}_{i}) ) ) = \emph{rank}( \mathbf{M}(\mathbf{a}_{i}) ) = k.(p-k)$ (see Corollary~\ref{corol5.6:box} for details). The associated quadratic models are
\begin{equation*}
L_{\lambda_{i}}^{\mathbf{N}_{i}} ( d\mathbf{a} ) = \psi( \mathbf{a}_{i} ) + d\mathbf{a}^{T}  \nabla \psi( \mathbf{a}_{i} ) +  \frac{1}{2} d\mathbf{a}^{T}  \big( \mathbf{M}( \mathbf{a}_{i}  )^{T} \mathbf{M}( \mathbf{a}_{i}  ) + \mathbf{L}( \mathbf{a}_{i}  )^{T} \mathbf{L}( \mathbf{a}_{i}  )  +  \mathbf{N}_{i} \mathbf{N}_{i}^{T} +  \lambda_{i} \mathbf{D}_{i}^{2}  \big) d\mathbf{a}
\end{equation*}
or
\begin{equation*}
L_{\lambda_{i}}^{\mathbf{N}_{i}}( d\mathbf{a} ) = \psi( \mathbf{a}_{i} ) + d\mathbf{a}^{T} \nabla \psi( \mathbf{a}_{i} )  +  \frac{1}{2} d\mathbf{a}^{T}  \big( \mathbf{M}( \mathbf{a}_{i}  )^{T} \mathbf{M}( \mathbf{a}_{i}  )  +  \mathbf{N}_{i} \mathbf{N}_{i}^{T} +  \lambda_{i} \mathbf{D}_{i}^{2} \big) d\mathbf{a}  \  ,
\end{equation*}
which can be minimized by solving the symmetric linear systems
\begin{equation*}
     \big( \mathbf{M}( \mathbf{a}_{i}  )^{T} \mathbf{M}( \mathbf{a}_{i}  ) + \mathbf{L}( \mathbf{a}_{i}  )^{T} \mathbf{L}( \mathbf{a}_{i}  ) + \mathbf{N}_{i} \mathbf{N}_{i}^{T} +  \lambda_{i} \mathbf{D}_{i}^{2}  \big) d\mathbf{a}_{gp-lm}  = \mathbf{M}( \mathbf{a}_{i}  )^{T} \mathbf{r}( \mathbf{a}_{i} )
\end{equation*}
or
\begin{equation*}
     \big( \mathbf{M}( \mathbf{a}_{i}  )^{T} \mathbf{M}( \mathbf{a}_{i}  ) + \mathbf{N}_{i} \mathbf{N}_{i}^{T} +  \lambda_{i} \mathbf{D}_{i}^{2}  \big) d\mathbf{a}_{k-lm}  = \mathbf{M}( \mathbf{a}_{i}  )^{T} \mathbf{r}( \mathbf{a}_{i} ) \ ,
\end{equation*}
which, in turn, are the normal equations for the constrained and damped linear systems
\begin{equation*}
\min_{d\mathbf{a} \in \mathbb{R}^{p.k}}   \,    \big\Vert  \begin{bmatrix} \mathbf{r}(  \mathbf{a}_{i} ) \\ \mathbf{0}^{k.k} \\ \mathbf{0}^{p.k} \end{bmatrix}  -  \begin{bmatrix}   \mathbf{M}( \mathbf{a}_{i}  ) + \mathbf{L}( \mathbf{a}_{i}  )  \\  \mathbf{N}_{i}^{T}  \\  \sqrt{\lambda_{i}} \mathbf{D}_{i}  \end{bmatrix}  d\mathbf{a}  \big\Vert^{2}_{2}
\end{equation*}
and
\begin{equation*}
\min_{d\mathbf{a} \in \mathbb{R}^{p.k}}   \,    \big\Vert  \begin{bmatrix} \mathbf{r}(  \mathbf{a}_{i} ) \\ \mathbf{0}^{k.k} \\ \mathbf{0}^{p.k} \end{bmatrix}  -  \begin{bmatrix}   \mathbf{M}( \mathbf{a}_{i}  )  \\  \mathbf{N}_{i}^{T}  \\  \sqrt{\lambda_{i}} \mathbf{D}_{i}  \end{bmatrix}  d\mathbf{a}  \big\Vert^{2}_{2} \ .
\end{equation*}
This approach was first suggested by Okatani  et al.~\cite{OYD2011}. These two linear least-squares problems  always have an unique solution independently of the value of $\lambda$ if we assume that
\begin{equation*}
\emph{rank}(\mathbf{A}_{i}) = k   \quad   \text{and}  \quad \emph{rank}( \mathit{J}( \mathbf{r}(\mathbf{a}_{i}) ) ) = \emph{rank}( \mathbf{M}(\mathbf{a}_{i}) ) = k.(p-k) \ ,
\end{equation*}
during the iterations or, at least, that these conditions are verified as soon as $\lambda = 0$. If $\lambda = 0$, we again obtained $d\mathbf{a}_{gp-lm} = d\mathbf{a}_{gp-gn}$ and $d\mathbf{a}_{k-lm} = d\mathbf{a}_{k-gn}$ as in version~\eqref{lm_alg1:box} of the Levenberg-Marquardt algorithms if the above assumptions are verified. However, the key-difference with this version~\eqref{lm_alg1:box} is that these two linear least-squares problems remain nonsingular and well-conditioned when $\lambda$ approaches zero if the above assumptions are verified thanks to the inclusion of the block $\mathbf{N}_{i}^{T}$ in the coefficient matrix of these linear least-squares problems or to the addition of the term  $\mathbf{N}_{i} \mathbf{N}_{i}^{T}$ in the associated normal equations if we use a normal-equation approach to solve them.

Finally, when $\lambda  \gg 0$, and we want to shift $d\mathbf{a}_{gp-lm}$ or $d\mathbf{a}_{k-lm}$ in the steepest descent direction in order to benefit of the good global convergence ability of the gradient descent method, the inclusion of the block $\mathbf{N}_{i}^{T}$ or the term $\mathbf{N}_{i} \mathbf{N}_{i}^{T}$ is of secondary importance and will not impair the performance of the algorithm as they just try to constrain the columns of the perturbation matrices $d\mathbf{A}_{gp-lm}$ or $d\mathbf{A}_{k-lm}$ to belong to $\emph{ran}( \mathbf{A}_{i}  )^{\bot}$, see the discussion after Corollary~\ref{corol5.6:box} for details.

Finally, note that the computations of the gain factor $\rho$ at step \textbf{(7)} must also be slightly modified as follows
\begin{equation*}
\rho = \frac{ \psi( \mathbf{a}_{i}) - \psi( \mathbf{a}_{i} + d\mathbf{a}_{lm} )}{ G(\mathbf{0}^{k.p}) -  G( d\mathbf{a}_{lm} ) } =  \frac{ \psi( \mathbf{a}_{i}) - \psi( \mathbf{a}_{i} + d\mathbf{a}_{lm} )}{  \frac{1}{2}   \big(  \Vert \mathbf{N}_{i}^{T}  d\mathbf{a}_{lm} \Vert^{2}_{2} + \lambda \Vert \mathbf{D}_{i}  d\mathbf{a}_{lm} \Vert^{2}_{2} - d\mathbf{a}_{lm}^{T} \nabla \psi( \mathbf{a}_{i} )  \big) }  \  ,
\end{equation*}
when the above constrained and regularized linear least-squares problems are solved in step \textbf{(6)}. This is justified by the fact that we cannot assume that $\Vert \mathbf{N}_{i}^{T}  d\mathbf{a}_{lm} \Vert_{2} = 0$ if a damping term $\lambda \Vert \mathbf{D}_{i}  d\mathbf{a}_{lm} \Vert^{2}_{2}$ with $\lambda > 0$ is also used in the algorithms.

These different considerations lead to our second and third versions of the Levenberg-Marquardt algorithms, which use, respectively, the simple and more sophisticated updating strategies of $\lambda$ discussed above:
\\
\begin{lm_alg2} \label{lm_alg2:box}
\end{lm_alg2}
Choose starting matrix $\mathbf{A}_{1} \in \mathbb{R}^{p \times k}$ , $\varepsilon_{1}, \varepsilon_{2}, \varepsilon_{3} , \lambda \in \mathbb{R}_{+*}$ and $i_{max}, j_{max} \in \mathbb{N}_{*}$, appropriately

\textbf{For} $i=1, 2, \ldots$ \textbf{until convergence do}
\begin{enumerate}
\item[\textbf{(0)}] Optionally, compute a QRCP of $\mathbf{A}_{i}$ (see equation~\eqref{eq:qrcp}) to determine $k_{i} = \emph{rank}( \mathbf{A}_{i} )$ and an orthonormal basis of $\emph{ran}( \mathbf{A}_{i} )$:
\begin{itemize}
\item[] $\mathbf{Q}_{i} \mathbf{A}_{i} \mathbf{P}_{i} =  \begin{bmatrix} \mathbf{R}_{i}    &  \mathbf{S}_{i}   \\  \mathbf{0}^{(p-k_{i}) \times k_{i} }  & \mathbf{0}^{(p-k_{i}) \times (k-k_{i}) }  \end{bmatrix}$ ,
\end{itemize}
where $\mathbf{Q}_{i}$ is an $p \times p$ orthogonal matrix, $\mathbf{P}_{i}$ is an $k \times k$ permutation matrix, $\mathbf{R}_{i}$ is an $k_{i} \times k_{i}$ nonsingular upper triangular matrix (with diagonal elements of decreasing absolute magnitude) and $\mathbf{S}_{i}$ an $k_{i} \times (k-k_{i})$ full matrix, which is vacuous if $k_{i} = k$.

In all cases, compute an $p \times k$  matrix $\mathbf{O}_{i}$  with orthonormal columns as the first $k$ columns of $\mathbf{Q}_{i}$ (i.e., such that $\emph{ran}( \mathbf{A}_{i} ) \subset \emph{ran}( \mathbf{O}_{i} )$ if $k_{i} < k$ and  $\emph{ran}( \mathbf{A}_{i} ) = \emph{ran}( \mathbf{O}_{i} )$ if $k_{i} = k$) and set
\begin{itemize}
\item[] $\mathbf{A}_{i} = \mathbf{O}_{i}$ .
\end{itemize}
This optional orthogonalization step is a safe-guard as the condition $k_{i} = k$ is a necessary condition for the differentiability of $\psi(.)$ at a point $\mathbf{A}_{i}$ and also to limit the occurence of overflows and underflows in the next steps by enforcing that the matrix variable $\mathbf{A}_{i} \in \mathbb{O}^{p \times k}$.
\item[\textbf{(1)}] Determine (implicitly) the block diagonal matrix
\begin{itemize}
\item[] $\mathbf{F}(\mathbf{a}_{i}) = \emph{diag}\big( \emph{vec}( \sqrt{\mathbf{W}} ) \big)  \big(  \mathbf{I}_n  \otimes \mathbf{A}_{i}  \big)$ ,
\end{itemize}
where $\mathbf{a}_{i} =  \emph{vec}( \mathbf{A}_{i}^{T} )$.
\item[\textbf{(2)}] Compute (implicitly) a QRCP of $\mathbf{F}(\mathbf{a}_{i})$ to determine $\mathbf{P}^{\bot}_{\mathbf{F}(  \mathbf{a}_{i} ) }$ and  $\mathbf{F}(\mathbf{a}_{i})^-$ (see equations~\eqref{eq:ginv_proj_ortho} and \eqref{eq:sginv_qrcp}) or, alternatively, a COD of  $\mathbf{F}(\mathbf{a}_{i})$ to determine $\mathbf{P}^{\bot}_{\mathbf{F}(  \mathbf{a}_{i} ) }$ and $\mathbf{F}(\mathbf{a}_{i})^{+}$ (see equations~\eqref{eq:ginv_proj_ortho} and \eqref{eq:ginv_cod}).

Note also that $\mathbf{F}(\mathbf{a}_{i})^- = \mathbf{F}(\mathbf{a}_{i})^{+}$ if $\mathbf{F}(\mathbf{a}_{i})$ is of full column rank and that $\mathbf{P}^{\bot}_{\mathbf{F}(  \mathbf{a}_{i} ) }$, $\mathbf{F}(\mathbf{a}_{i})^-$ and $\mathbf{F}(\mathbf{a}_{i})^{+}$ are also block diagonal matrices.
\item[\textbf{(3)}] Solve the block diagonal linear least-squares problem
\begin{itemize}
\item[] $\mathbf{b}_{i} = \text{Arg}\min_{\mathbf{b}\in\mathbb{R}^{k.n}} \,   \Vert \mathbf{x} - \mathbf{F}(\mathbf{a}_{i})\mathbf{b} \Vert^{2}_{2}$ ,
\end{itemize}
e.g., compute
\begin{itemize}
\item[] $\mathbf{b}_{i} =
    \begin{cases}
         \mathbf{F}(\mathbf{a}_{i})^{-} \mathbf{x}   \quad\ \lbrace \text{if a QRCP of $\mathbf{F}(\mathbf{a}_{i})$ is used in step } \textbf{(2)}  \rbrace \\
         \mathbf{F}(\mathbf{a}_{i})^{+} \mathbf{x}  \quad\ \lbrace \text{if a COD of $\mathbf{F}(\mathbf{a}_{i})$ is used in step } \textbf{(2)} \rbrace
    \end{cases}$  \ .
\end{itemize}
\item[\textbf{(4)}] Determine and set:
\begin{itemize}
\item[] $\mathbf{r}(\mathbf{a}_{i}) = \mathbf{P}^{\bot}_{\mathbf{F}(  \mathbf{a}_{i} ) } \mathbf{x}$  $\lbrace \text{current residual vector} \rbrace$
\item[] $\psi(\mathbf{a}_{i} ) = \frac{1}{2} \Vert \mathbf{r}(\mathbf{a}_{i}) \Vert^{2}_{2}$ $\lbrace \text{current value of the cost function }  \rbrace$
\item[] $\nabla \psi( \mathbf{a}_{i} ) =  \mathbf{G}(\mathbf{b}_{i})^{T} \mathbf{G}(\mathbf{b}_{i})\mathbf{a}_{i} - \mathbf{G}(\mathbf{b}_{i})^{T} \mathbf{z}$ $\lbrace$see Theorems~\ref{theo4.3:box} and~\ref{theo5.7:box}$\rbrace$
\item[] $j = 0$ $\lbrace \text{initialize counter for the ridge scaling subiterations} \rbrace$
\end{itemize}
Note that the steps  \textbf{(1)} to \textbf{(4)}  above can be very easily parallelized using the block diagonal structure of $\mathbf{F}(\mathbf{a}_{i})$.
\item[\textbf{(5)}]  Check for convergence. Relevant convergence criteria in the algorithms are of the form:
\begin{itemize}
\item $\Vert \nabla \psi(\mathbf{a}_{i} ) \Vert_{2} \le \varepsilon_{1}$
\item $\Vert \mathbf{a}_{i} - \mathbf{a}_{i-1} \Vert_{2} \le \varepsilon_{2} ( \varepsilon_{2} + \Vert \mathbf{a}_{i} \Vert_{2})$  $\lbrace$if $ i \ne 1\rbrace$

If step \textbf{(0)} is used, this convergence condition can be simplified as:

$\Vert \mathbf{a}_{i} - \mathbf{a}_{i-1} \Vert_{2} \le \varepsilon_{2}  \Vert \mathbf{a}_{i} \Vert_{2} = \varepsilon_{2} \sqrt{k}$
\item $\vert \psi(\mathbf{a}_{i-1} ) - \psi(\mathbf{a}_{i} ) \vert \le \varepsilon_{3} ( \varepsilon_{3} +  \psi(\mathbf{a}_{i} ) )$  $\lbrace$if $ i \ne 1\rbrace$
\item $ i \ge i_{max}$ $\lbrace \text{e.g., give up if the number of iterations is too large} \rbrace$
\end{itemize}
where $\varepsilon_{1}, \varepsilon_{2}, \varepsilon_{3}$ and $i_{max}$ are constants chosen by the user.

\textbf{Exit if convergence}. \textbf{Otherwise, go to} step \textbf{(6)}
\item[\textbf{(6)}]Compute the Levenberg-Marquardt correction vector $d\mathbf{a}_{lm}$ as the solution of one of the following constrained and damped linear least-squares problems:
\begin{description}
\item[Golub-Pereyra Levenberg-Marquardt step:]
\begin{align*}
d\mathbf{a}_{gp-lm} & =    \begin{bmatrix}  \mathbf{M}( \mathbf{a}_{i} ) + \mathbf{L}( \mathbf{a}_{i}  ) \\  \mathbf{N}_{i}^{T}  \\ \sqrt{\lambda} \mathbf{D}_{i}  \end{bmatrix}^{+}  \begin{bmatrix}  \mathbf{r}( \mathbf{a}_{i} )  \\   \mathbf{0}^{k.k}   \\   \mathbf{0}^{k.p} \end{bmatrix}  \\
                                 & =   \text{Arg}\min_{d\mathbf{a} \in \mathbb{R}^{p.k}} \,   \Vert \mathbf{r}( \mathbf{a}_{i} ) - \big( \mathbf{M}( \mathbf{a}_{i}  ) + \mathbf{L}( \mathbf{a}_{i}  )  \big) d\mathbf{a}   \Vert^{2}_{2} +  \Vert \mathbf{N}_{i}^{T}  d\mathbf{a} \Vert^{2}_{2} + \lambda \Vert \mathbf{D}_{i}  d\mathbf{a} \Vert^{2}_{2} 
\end{align*}
\item[Kaufman Levenberg-Marquardt  step:] Okatani  et al.~\cite{OYD2011}
\begin{align*}
d\mathbf{a}_{k-lm} & =  \begin{bmatrix}  \mathbf{M}( \mathbf{a}_{i} ) \\  \mathbf{N}_{i}^{T}  \\   \sqrt{\lambda} \mathbf{D}_{i}  \end{bmatrix}^{+}  \begin{bmatrix}  \mathbf{r}( \mathbf{a}_{i} )  \\   \mathbf{0}^{k.k}  \\   \mathbf{0}^{k.p} \end{bmatrix}  \\
                               & =  \text{Arg}\min_{d\mathbf{a} \in \mathbb{R}^{p.k}} \,  \Vert \mathbf{r}( \mathbf{a}_{i} ) - \mathbf{M}( \mathbf{a}_{i}  ) d\mathbf{a}  \Vert^{2}_{2} +  \Vert \mathbf{N}_{i}^{T}  d\mathbf{a} \Vert^{2}_{2} + \lambda \Vert \mathbf{D}_{i}  d\mathbf{a} \Vert^{2}_{2} 
\end{align*}
\end{description}
where the columns of $\mathbf{N}_{i} = \mathbf{K}_{(p,k)} ( \mathbf{I}_{k} \otimes \mathbf{A}_{i} )$ are a (orthonormal) basis of  $\emph{null}\big( \mathbf{M}(\mathbf{a}_{i}) + \mathbf{L}(\mathbf{a}_{i}) \big) = \emph{null}( \mathbf{M}(\mathbf{a}_{i}) )$, see Corollary~\ref{corol5.6:box}.
\item[\textbf{(7)}] Compute $\mathbf{a}_{i+1} = \emph{vec}( \mathbf{A}_{i+1}^{T} )$ such that $\psi( \mathbf{a}_{i+1} ) < \psi( \mathbf{a}_{i} )$ in order to obtain global convergence.
\begin{enumerate}
\item[\textbf{(7.1)}]  To this end, first compute
\begin{itemize}
\item[] $\psi( \mathbf{a}_{i} +  d\mathbf{a}_{lm} ) = \frac{1}{2} \Vert \mathbf{r}(\mathbf{a}_{i} +  d\mathbf{a}_{lm}) \Vert^{2}_{2} = \frac{1}{2} \Vert \mathbf{P}^{\bot}_{\mathbf{F}(  \mathbf{a}_{i} +  d\mathbf{a}_{lm} ) } \mathbf{x} \Vert^{2}_{2}$ ,
\end{itemize}
using (implicitly) a QRCP of the block diagonal matrix $\mathbf{F}(  \mathbf{a}_{i} +  d\mathbf{a}_{lm} )$.
\item[\textbf{(7.2)}] \textbf{If}  $\psi( \mathbf{a}_{i} +  d\mathbf{a}_{lm} ) > \psi( \mathbf{a}_{i} )$ \textbf{then}  $\lbrace \text{step rejected} \rbrace$
\begin{itemize}
\item[] $j = j + 1$
\item[] $\lambda = 10.\lambda$ $\lbrace \text{scale up the ridge parameter} \rbrace$
\item[] \textbf{If} $j  \le j_{max}$ \textbf{go to} step  \textbf{(6)} $\lbrace \text{recompute } d\mathbf{a}_{lm} \text{ with inflated diagonal}  \rbrace$
\end{itemize}
\item[\textbf{(7.3)}] \textbf{Else} $\lbrace \text{step acceptable} \rbrace$
\begin{itemize}
\item[] \textbf{If} $j  = 0$   \textbf{then} $\lambda = \lambda /10$  $\lbrace \text{scale down the ridge parameter if step is successful} \rbrace$
\end{itemize}
\item[\textbf{(7.4)}]  Increment $\mathbf{a}_{i}$:

$\mathbf{a}_{i+1} = \mathbf{a}_{i} +  d\mathbf{a}_{lm} \lbrace \text{compute new iterate} \rbrace $
\end{enumerate}
\end{enumerate}
\textbf{End do}
\\
\\
\begin{lm_alg3} \label{lm_alg3:box}
\end{lm_alg3}
Choose starting matrix $\mathbf{A}_{1} \in \mathbb{R}^{p \times k}$ , $\varepsilon_{1}, \varepsilon_{2}, \varepsilon_{3} , \lambda  \in \mathbb{R}_{+*}$ and $i_{max}, j_{max} \in \mathbb{N}_{*}$, appropriately, and initialize $\nu = 2$

\textbf{For} $i=1, 2, \ldots$ \textbf{until convergence do}
\begin{enumerate}
\item[\textbf{(0)}]  Optionally, compute a QRCP of $\mathbf{A}_{i}$ (see equation~\eqref{eq:qrcp}) to determine $k_{i} = \emph{rank}( \mathbf{A}_{i} )$ and an orthonormal basis of $\emph{ran}( \mathbf{A}_{i} )$:
\begin{itemize}
\item[] $\mathbf{Q}_{i} \mathbf{A}_{i} \mathbf{P}_{i} =  \begin{bmatrix} \mathbf{R}_{i}    &  \mathbf{S}_{i}   \\  \mathbf{0}^{(p-k_{i}) \times k_{i} }  & \mathbf{0}^{(p-k_{i}) \times (k-k_{i}) }  \end{bmatrix}$ ,
\end{itemize}
where $\mathbf{Q}_{i}$ is an $p \times p$ orthogonal matrix, $\mathbf{P}_{i}$ is an $k \times k$ permutation matrix, $\mathbf{R}_{i}$ is an $k_{i} \times k_{i}$ nonsingular upper triangular matrix (with diagonal elements of decreasing absolute magnitude) and $\mathbf{S}_{i}$ an $k_{i} \times (k-k_{i})$ full matrix, which is vacuous if $k_{i} = k$.

In all cases, compute an $p \times k$  matrix $\mathbf{O}_{i}$  with orthonormal columns as the first $k$ columns of $\mathbf{Q}_{i}$ (i.e., such that $\emph{ran}( \mathbf{A}_{i} ) \subset \emph{ran}( \mathbf{O}_{i} )$ if $k_{i} < k$ and  $\emph{ran}( \mathbf{A}_{i} ) = \emph{ran}( \mathbf{O}_{i} )$ if $k_{i} = k$) and set
\begin{itemize}
\item[] $\mathbf{A}_{i} = \mathbf{O}_{i}$ .
\end{itemize}
This optional orthogonalization step is a safe-guard as the condition $k_{i} = k$ is a necessary condition for the differentiability of $\psi(.)$ at a point $\mathbf{A}_{i}$ and also to limit the occurence of overflows and underflows in the next steps by enforcing that the matrix variable $\mathbf{A}_{i} \in \mathbb{O}^{p \times k}$.
\item[\textbf{(1)}]  Determine (implicitly) the block diagonal matrix
\begin{itemize}
\item[] $\mathbf{F}(\mathbf{a}_{i}) = \emph{diag}\big( \emph{vec}( \sqrt{\mathbf{W}} ) \big)  \big(  \mathbf{I}_n  \otimes \mathbf{A}_{i}  \big)$ ,
\end{itemize}
where $\mathbf{a}_{i} =  \emph{vec}( \mathbf{A}_{i}^{T} )$.
\item[\textbf{(2)}] Compute (implicitly) a QRCP of $\mathbf{F}(\mathbf{a}_{i})$ to determine $\mathbf{P}^{\bot}_{\mathbf{F}(  \mathbf{a}_{i} ) }$ and  $\mathbf{F}(\mathbf{a}_{i})^-$ (see equations~\eqref{eq:ginv_proj_ortho} and \eqref{eq:sginv_qrcp}) or, alternatively, a COD of  $\mathbf{F}(\mathbf{a}_{i})$ to determine $\mathbf{P}^{\bot}_{\mathbf{F}(  \mathbf{a}_{i} ) }$ and $\mathbf{F}(\mathbf{a}_{i})^{+}$ (see equations~\eqref{eq:ginv_proj_ortho} and \eqref{eq:ginv_cod}).

Note also that $\mathbf{F}(\mathbf{a}_{i})^- = \mathbf{F}(\mathbf{a}_{i})^{+}$ if $\mathbf{F}(\mathbf{a}_{i})$ is of full column rank and that $\mathbf{P}^{\bot}_{\mathbf{F}(  \mathbf{a}_{i} ) }$, $\mathbf{F}(\mathbf{a}_{i})^-$ and $\mathbf{F}(\mathbf{a}_{i})^{+}$ are also block diagonal matrices.
\item[\textbf{(3)}]  Solve the block diagonal linear least-squares problem
\begin{itemize}
\item[] $\mathbf{b}_{i} = \text{Arg}\min_{\mathbf{b}\in\mathbb{R}^{k.n}} \,   \Vert \mathbf{x} - \mathbf{F}(\mathbf{a}_{i})\mathbf{b} \Vert^{2}_{2}$,
\end{itemize}
e.g., compute
\begin{itemize}
\item[] $\mathbf{b}_{i} =
    \begin{cases}
         \mathbf{F}(\mathbf{a}_{i})^{-} \mathbf{x}   \quad\ \lbrace \text{if a QRCP of $\mathbf{F}(\mathbf{a}_{i})$ is used in step } \textbf{(2)}  \rbrace \\
         \mathbf{F}(\mathbf{a}_{i})^{+} \mathbf{x}  \quad\ \lbrace \text{if a COD of $\mathbf{F}(\mathbf{a}_{i})$ is used in step } \textbf{(2)} \rbrace
    \end{cases}$ .
\end{itemize}
\item[\textbf{(4)}]  Determine and set:
\begin{itemize}
\item[] $\mathbf{r}(\mathbf{a}_{i}) = \mathbf{P}^{\bot}_{\mathbf{F}(  \mathbf{a}_{i} ) } \mathbf{x}$ $\lbrace \text{current residual vector} \rbrace$
\item[] $\psi(\mathbf{a}_{i} ) = \frac{1}{2} \Vert \mathbf{r}(\mathbf{a}_{i}) \Vert^{2}_{2}$ $\lbrace \text{current value of the cost function }  \rbrace$
\item[] $\nabla \psi( \mathbf{a}_{i} ) =  \mathbf{G}(\mathbf{b}_{i})^{T} \mathbf{G}(\mathbf{b}_{i})\mathbf{a}_{i} - \mathbf{G}(\mathbf{b}_{i})^{T} \mathbf{z}$ $\lbrace$see Theorems~\ref{theo4.3:box} and~\ref{theo5.7:box}$\rbrace$
\item[] $j = 0$ $\lbrace \text{initialize counter for the ridge scaling subiterations} \rbrace$
\end{itemize}
Note that the steps \textbf{(1)} to \textbf{(4)} above can be very easily parallelized using the block diagonal structure of $\mathbf{F}(\mathbf{a}_{i})$.
\item[\textbf{(5)}]   Check for convergence. Relevant convergence criteria in the algorithms are of the form:
\begin{itemize}
\item $\Vert \nabla \psi(\mathbf{a}_{i} ) \Vert_{2} \le \varepsilon_{1}$
\item $\Vert \mathbf{a}_{i} - \mathbf{a}_{i-1} \Vert_{2} \le \varepsilon_{2} ( \varepsilon_{2} + \Vert \mathbf{a}_{i} \Vert_{2})$  $\lbrace$if $ i \ne 1\rbrace$

If step \textbf{(0)} is used, this convergence condition can be simplified as:

$\Vert \mathbf{a}_{i} - \mathbf{a}_{i-1} \Vert_{2} \le \varepsilon_{2}  \Vert \mathbf{a}_{i} \Vert_{2} = \varepsilon_{2} \sqrt{k}$
\item $\vert \psi(\mathbf{a}_{i-1} ) - \psi(\mathbf{a}_{i} ) \vert \le \varepsilon_{3} ( \varepsilon_{3} +  \psi(\mathbf{a}_{i} ) )$  $\lbrace$if $ i \ne 1\rbrace$
\item $ i \ge i_{max}$ $\lbrace \text{e.g., give up if the number of iterations is too large} \rbrace$
\end{itemize}
where $\varepsilon_{1}, \varepsilon_{2}, \varepsilon_{3}$ and $i_{max}$ are constants chosen by the user.

\textbf{Exit if convergence}. \textbf{Otherwise, go to} step \textbf{(6)}
\item[\textbf{(6)}]  Compute the Levenberg-Marquardt correction vector $d\mathbf{a}_{lm}$ as the solution of one of the following constrained and damped linear least-squares problems:
\begin{description}
\item[Golub-Pereyra Levenberg-Marquardt step:]
\begin{align*}
d\mathbf{a}_{gp-lm} & =    \begin{bmatrix}  \mathbf{M}( \mathbf{a}_{i} ) + \mathbf{L}( \mathbf{a}_{i}  ) \\  \mathbf{N}_{i}^{T}  \\ \sqrt{\lambda} \mathbf{D}_{i}  \end{bmatrix}^{+}  \begin{bmatrix}  \mathbf{r}( \mathbf{a}_{i} )  \\   \mathbf{0}^{k.k}   \\   \mathbf{0}^{k.p} \end{bmatrix}  \\
                                 & =   \text{Arg}\min_{d\mathbf{a} \in \mathbb{R}^{p.k}} \,   \Vert \mathbf{r}( \mathbf{a}_{i} ) - \big( \mathbf{M}( \mathbf{a}_{i}  ) + \mathbf{L}( \mathbf{a}_{i}  )  \big) d\mathbf{a}   \Vert^{2}_{2} +  \Vert \mathbf{N}_{i}^{T}  d\mathbf{a} \Vert^{2}_{2} + \lambda \Vert \mathbf{D}_{i}  d\mathbf{a} \Vert^{2}_{2}
\end{align*}
\item[Kaufman Levenberg-Marquardt  step:] Okatani  et al.~\cite{OYD2011}
\begin{align*}
d\mathbf{a}_{k-lm} & =  \begin{bmatrix}  \mathbf{M}( \mathbf{a}_{i} ) \\  \mathbf{N}_{i}^{T}  \\   \sqrt{\lambda} \mathbf{D}_{i}  \end{bmatrix}^{+}  \begin{bmatrix}  \mathbf{r}( \mathbf{a}_{i} )  \\   \mathbf{0}^{k.k}  \\   \mathbf{0}^{k.p} \end{bmatrix}  \\
                               & =  \text{Arg}\min_{d\mathbf{a} \in \mathbb{R}^{p.k}} \,  \Vert \mathbf{r}( \mathbf{a}_{i} ) - \mathbf{M}( \mathbf{a}_{i}  ) d\mathbf{a}  \Vert^{2}_{2} +  \Vert \mathbf{N}_{i}^{T}  d\mathbf{a} \Vert^{2}_{2} + \lambda \Vert \mathbf{D}_{i}  d\mathbf{a} \Vert^{2}_{2}
\end{align*}
\end{description}
where the columns of $\mathbf{N}_{i} = \mathbf{K}_{(p,k)} ( \mathbf{I}_{k} \otimes \mathbf{A}_{i} )$ are a (orthonormal) basis of  $\emph{null}\big( \mathbf{M}(\mathbf{a}_{i}) + \mathbf{L}(\mathbf{a}_{i}) \big) = \emph{null}( \mathbf{M}(\mathbf{a}_{i}) )$, see Corollary~\ref{corol5.6:box}.
\item[\textbf{(7)}]  Compute $\mathbf{a}_{i+1} = \emph{vec}( \mathbf{A}_{i+1}^{T} )$ such that $\psi( \mathbf{a}_{i+1} ) < \psi( \mathbf{a}_{i} )$ in order to obtain global convergence.
\begin{enumerate}
\item[\textbf{(7.1)}]  To this end, first compute
\begin{itemize}
\item[] $\psi( \mathbf{a}_{i} +  d\mathbf{a}_{lm} ) = \frac{1}{2} \Vert \mathbf{P}^{\bot}_{\mathbf{F}(  \mathbf{a}_{i} +  d\mathbf{a}_{lm} ) } \mathbf{x} \Vert^{2}_{2}$ ,
\end{itemize}
using (implicitly) a QRCP of the block diagonal matrix $\mathbf{F}(  \mathbf{a}_{i} +  d\mathbf{a}_{lm} )$, and the gain factor
\begin{itemize}
\item[] $\rho =  \frac{ \psi( \mathbf{a}_{i}) - \psi( \mathbf{a}_{i} + d\mathbf{a}_{lm} )}{ G(\mathbf{0}^{k.p}) -  G( d\mathbf{a}_{lm} ) } = \frac{ \psi( \mathbf{a}_{i}) - \psi( \mathbf{a}_{i} + d\mathbf{a}_{lm} )}{  \frac{1}{2}   \big(  \Vert \mathbf{N}_{i}^{T}  d\mathbf{a}_{lm} \Vert^{2}_{2} + \lambda \Vert \mathbf{D}_{i}  d\mathbf{a}_{lm} \Vert^{2}_{2} - d\mathbf{a}_{lm}^{T} \nabla \psi( \mathbf{a}_{i} )  \big) }$
\end{itemize}
\item[\textbf{(7.2)}]  \textbf{If}  $\rho > 0$ \textbf{then} $\lbrace \text{step acceptable} \rbrace$
\begin{itemize}
\item[] $\lambda = \lambda .\emph{max}\big( \frac{1}{3}, 1-(2.\rho-1)^3 \big)$ $\lbrace \text{scale down the ridge parameter} \rbrace$
\item[] $\nu = 2$ $\lbrace \text{reinitialize the growth factor of the ridge parameter} \rbrace$
\end{itemize}
\item[\textbf{(7.3)}]  \textbf{Else} $\lbrace \text{step rejected} \rbrace$
\begin{itemize}
\item[] $j = j + 1$
\item[] $\lambda = \nu.\lambda$ $\lbrace \text{scale up the ridge parameter} \rbrace$
\item[] $\nu = 2.\nu$ $\lbrace \text{increase the growth factor of the ridge parameter} \rbrace$
\item[] \textbf{If} $j  \le j_{max}$ \textbf{go to} step \textbf{(6)} $\lbrace \text{recompute } d\mathbf{a}_{lm} \text{ with inflated diagonal} \rbrace$
\end{itemize}
\item[\textbf{(7.4)}]  Increment $\mathbf{a}_{i}$:

$\mathbf{a}_{i+1} = \mathbf{a}_{i} +  d\mathbf{a}_{lm}  \lbrace \text{compute new iterate} \rbrace $
\end{enumerate}
\end{enumerate}
\textbf{End do}

In the Levenberg-Marquardt algorithms~\eqref{lm_alg2:box} and~\eqref{lm_alg3:box}, the Marquardt parameter $\lambda$ is initialized to $\lambda = \tau$ if
\begin{equation*}
\big \lbrack \mathbf{D}_{i} \big \rbrack_{jj} = 
    \begin{cases}
 \Vert  \big \lbrack  \mathbf{M}( \mathbf{a}_{i}  ) +  \mathbf{L}( \mathbf{a}_{i}  ) \rbrack_{.j}  \Vert_{2}    &\lbrace \text{when } d\mathbf{a}_{lm} = d\mathbf{a}_{gp-lm}  \text{ in step \textbf{(6)}} \rbrace \\
 \Vert  \big \lbrack  \mathbf{M}( \mathbf{a}_{i}  ) \rbrack_{.j}  \Vert_{2}    &\lbrace \text{when }   d\mathbf{a}_{lm} = d\mathbf{a}_{k-lm}  \text{ in step \textbf{(6)}} \rbrace,
    \end{cases}
\end{equation*}
for $j = 1, \cdots, k.p$ during the iterations, or to
\begin{equation*}
\lambda = \tau.
    \begin{cases}
  \emph{max}_{ j = 1, \cdots, k.p} \Vert  \big \lbrack ( \mathbf{M}( \mathbf{a}_{1}  ) +  \mathbf{L}( \mathbf{a}_{1}  ) \rbrack_{.j}  \Vert_{2}^{2} )  &\lbrace \text{when } d\mathbf{a}_{lm} = d\mathbf{a}_{gp-lm}  \text{ in step \textbf{(6)}} \rbrace \\\
  \emph{max}_{ j = 1, \cdots, k.p} \Vert  \big \lbrack ( \mathbf{M}( \mathbf{a}_{1}  ) \rbrack_{.j}  \Vert_{2}^{2})     &\lbrace \text{when }   d\mathbf{a}_{lm} = d\mathbf{a}_{k-lm}  \text{ in step \textbf{(6)}} \rbrace,
    \end{cases}
\end{equation*}
if $\mathbf{D}_{i}$ is set to the identity matrix during the iterations. In both cases $\tau$ is taken in the interval $\lbrack  10^{-8}  \,  1 \rbrack$ and a small value of $\tau$ is selected if we believe that $\mathbf{A}_{1}$ is close to a solution (say $\tau = 10^{-6}$).  Otherwise, we can use $\tau=10^{-3}$ or even $1$. The algorithms are not very sensitive to this initial choice of  $\tau$ as $\lambda$ is quickly updated during the iterations in both Levenberg-Marquardt algorithms~\eqref{lm_alg2:box} and~\eqref{lm_alg3:box}. Version~\eqref{lm_alg3:box} of the Levenberg-Marquardt algorithms also uses a growth factor $\nu$ for the ridge parameter, which is  initialized to $2$ at the start of the algorithm and reinitialized to this initial value in step \textbf{(7.2)} when a step is successful.
\\
\begin{remark6.3} \label{remark6.3:box}
An alternative for computing the correction vectors $d\mathbf{a}_{gp-lm}$ and $d\mathbf{a}_{k-lm}$ in step \textbf{(6)} of the Levenberg-Marquardt algorithms~\eqref{lm_alg2:box} and~\eqref{lm_alg3:box}, if we assume again that $\emph{rank}( \mathbf{A}_{i} ) = k$ and $\emph{rank}( \mathit{J}( \mathbf{r}(\mathbf{a}_{i}) ) ) = \emph{rank}( \mathbf{M}(\mathbf{a}_{i}) ) = (p-k).k$, is to first find the unique solutions of the following "reduced" and damped linear least-squares problems
\begin{equation*}
 d\bar{\mathbf{a}}_{gp-lm} = \text{Arg}\min_{d\bar{\mathbf{a}} \in \mathbb{R}^{(p-k).k}} \,   \Vert \mathbf{r}( \mathbf{a}_{i} ) - \big( \mathbf{M}( \mathbf{a}_{i}  ) + \mathbf{L}( \mathbf{a}_{i}  )  \big)\mathbf{\bar{O}}_{i}^{\bot} d\bar{\mathbf{a}}   \Vert^{2}_{2} + \lambda \Vert \bar{\mathbf{D}}_{i}  d\bar{\mathbf{a}} \Vert^{2}_{2} \ ,
\end{equation*}
or
\begin{equation*}
d\bar{\mathbf{a}}_{k-lm} = \text{Arg}\min_{d\bar{\mathbf{a}} \in \mathbb{R}^{(p-k).k}} \,   \Vert \mathbf{r}( \mathbf{a}_{i} ) - \mathbf{M}( \mathbf{a}_{i}  ) \mathbf{\bar{O}}_{i}^{\bot} d\bar{\mathbf{a}}   \Vert^{2}_{2} + \lambda \Vert \bar{\mathbf{D}}_{i}  d\bar{\mathbf{a}} \Vert^{2}_{2} \ ,
\end{equation*}
where $\mathbf{\bar{O}}_{i}^{\bot}$ is an orthonormal basis of  $\emph{null}\big( \mathit{J}( \mathbf{r}(\mathbf{a}_{i}) ) \big)^{\bot} = \emph{null}( \mathbf{M}(\mathbf{a}_{i}) )^{\bot}$ and $\bar{\mathbf{D}}_{i}$ is now a positive diagonal matrix of order $(p-k).k$, see Corollary~\ref{corol5.6:box} for more information. Finally, in an additional step just before incrementing $\mathbf{a}_{i}$ in step \textbf{(7.4)}, the correction vectors $d\mathbf{a}_{gp-lm}$ and $d\mathbf{a}_{k-lm}$ can be computed by the matrix-vector products
\begin{equation*}
d\mathbf{a}_{gp-lm} =  \mathbf{\bar{O}}_{i}^{\bot} d\bar{\mathbf{a}}_{gp-lm}   \quad \text{and}  \quad d\mathbf{a}_{k-lm} = \mathbf{\bar{O}}_{i}^{\bot} d\bar{\mathbf{a}}_{k-lm} \ ,
\end{equation*}
or, equivalently, the matrix-matrix products
\begin{equation*}
d\mathbf{A}_{gp-lm} =  \mathbf{O}_{i}^{\bot} d\bar{\mathbf{A}}_{gp-lm}    \quad \text{and}  \quad d\mathbf{A}_{k-lm} = \mathbf{O}_{i}^{\bot} d\bar{\mathbf{A}}_{k-lm} \ ,
\end{equation*}
where $\mathbf{O}_{i}^{\bot}$ is an orthonormal basis of $\emph{ran}( \mathbf{A}_{i} )^{\bot}$, as also described in Subsection~\ref{jacob:box}. $\blacksquare$
\\
\end{remark6.3}

We now consider in more details how to compute the correction vectors $d\mathbf{a}_{gp-lm}$ and $d\mathbf{a}_{k-lm}$ in the variable projection Levenberg-Marquardt algorithms~\eqref{lm_alg1:box},~\eqref{lm_alg2:box} and~\eqref{lm_alg3:box} using the normal-equation or TSQR approaches.

Using the normal-equation framework, the first step to obtain the correction vectors  $d\mathbf{a}_{gp-lm}$ or $d\mathbf{a}_{k-lm}$ in all the Levenberg-Marquardt algorithms is to form the cross-product positive semi-definite matrices  (e.g., the Gauss-Newton approximations of the Hessian matrix)
\begin{equation*}
\Delta =  \mathit{J} \big ( \mathbf{r}(\mathbf{a}) \big)^{T} \mathit{J} \big( \mathbf{r}(\mathbf{a}) \big) = \mathbf{M}(\mathbf{a} )^{T} \mathbf{M}(\mathbf{a} ) + \mathbf{L}(\mathbf{a} )^{T} \mathbf{L}(\mathbf{a} )  \quad  \text{or}  \quad  \Lambda =  \mathbf{M}(\mathbf{a} )^{T} \mathbf{M}(\mathbf{a} )  \ ,
\end{equation*}
exactly as in the Gauss-Newton algorithms~\eqref{gn_alg:box} described in the previous subsection. Note that, in these last equations and the rest of this subsection, we have drop  the iteration index $i$ of the Levenberg-Marquardt algorithms for notational convenience and the notations are exactly similar as in the Gauss-Newton methods described in the previous subsection. Furthermore, this first and costly step can be easily parallelized as for the Gauss-Newton algorithms~\eqref{gn_alg:box}.

For the Levenberg-Marquardt algorithms~\eqref{lm_alg1:box}, in a second stage, we just need to regularize (or damp) these positive semi-definite matrices by adding the diagonal matrix $\lambda \mathbf{D}^{2}$ to these Gauss-Newton approximations of the Hessian matrix:
\begin{equation*}
\Delta( \lambda ) = \Delta + \lambda \mathbf{D}^{2} \quad  \text{or}  \quad   \Lambda( \lambda ) = \Lambda + \lambda \mathbf{D}^{2} \ ,
\end{equation*}
where $\lambda > 0$ is a ridge parameter, which will control both the magnitude and the direction of the correction vectors $d\mathbf{a}_{gp-lm}$ and $d\mathbf{a}_{k-lm}$, and $\mathbf{D}$ is a strictly positive diagonal matrix. Finally, in the third and last step of the Levenberg-Marquardt algorithms~\eqref{lm_alg1:box}, we have to solve the consistent linear systems
\begin{equation} \label{eq:neq_gp_lm_alg1}
\Delta( \lambda ) d\mathbf{a}_{gp-lm} = - \mathit{J} \big( \mathbf{r}(\mathbf{a}) \big)^{T} \mathbf{r}(\mathbf{a}) = \mathbf{M}(\mathbf{a} )^{T} \mathbf{r}(\mathbf{a})
\end{equation}
or
\begin{equation}  \label{eq:neq_k_lm_alg1}
\Lambda( \lambda ) d\mathbf{a}_{k-lm} = - \mathit{J} \big( \mathbf{r}(\mathbf{a}) \big)^{T} \mathbf{r}(\mathbf{a}) = \mathbf{M}(\mathbf{a} )^{T} \mathbf{r}(\mathbf{a}) \ ,
\end{equation}
in order to obtain the correction vectors $d\mathbf{a}_{gp-lm}$ and $d\mathbf{a}_{k-lm}$. 
As $\lambda > 0$, $\Delta( \lambda )$ and $\Lambda( \lambda )$ are positive definite matrices and we can simply compute their Cholesky factorizations as
\begin{equation*}
\Delta( \lambda ) = \mathbf{R}_{\Delta}( \lambda )^{T} \mathbf{R}_{\Delta}( \lambda )  \quad  \text{or}  \quad    \Lambda( \lambda ) = \mathbf{R}_{\Lambda}( \lambda )^{T} \mathbf{R}_{\Lambda}( \lambda )  \ ,
\end{equation*}
where $\mathbf{R}_{\Delta}( \lambda )$ and $\mathbf{R}_{\Lambda}( \lambda )$ are $k.p \times k.p$ nonsingular upper triangular matrices. Then, in a final step, $d\mathbf{a}_{gp-lm}$ and $d\mathbf{a}_{k-lm}$ can be obtained by forward and backward substitutions in the usual manner, using these Cholesky factors, as
\begin{equation*}
d\mathbf{a}_{gp-lm} = \mathbf{R}_{\Delta}( \lambda )^{-1}  \mathbf{R}_{\Delta}( \lambda )^{-T} \mathbf{M}(\mathbf{a} )^{T} \mathbf{r}(\mathbf{a})  \quad  \text{and}  \quad d\mathbf{a}_{k-lm} =  \mathbf{R}_{\Lambda}( \lambda )^{-1}  \mathbf{R}_{\Lambda}( \lambda )^{-T}  \mathbf{M}(\mathbf{a} )^{T} \mathbf{r}(\mathbf{a}) \ .
\end{equation*}

On the other hand, if we use a normal-equation approach in both the Levenberg-Marquardt algorithms~\eqref{lm_alg2:box} and~\eqref{lm_alg3:box}, we have first to compute the constrained and damped cross-product (approximated) Jacobian matrices
\begin{equation*}
\Delta( \mathbf{N}, \lambda ) = \Delta + \mathbf{N} \mathbf{N}^{T} + \lambda \mathbf{D}^{2}
 \quad \text{or}  \quad
 \Lambda( \mathbf{N}, \lambda ) = \Lambda + \mathbf{N} \mathbf{N}^{T} + \lambda \mathbf{D}^{2} \  ,
\end{equation*}
perform their Cholesky decompositions as
\begin{equation*}
\Delta( \mathbf{N}, \lambda ) = \mathbf{R}_{\Delta}( \mathbf{N}, \lambda )^{T} \mathbf{R}_{\Delta}( \mathbf{N}, \lambda )   \quad \text{or}   \quad  \Lambda( \mathbf{N}, \lambda ) = \mathbf{R}_{\Lambda}( \mathbf{N}, \lambda )^{T} \mathbf{R}_{\Lambda}( \mathbf{N}, \lambda ) \  ,
\end{equation*}
where $\mathbf{R}_{\Delta}( \mathbf{N}, \lambda )$ and $\mathbf{R}_{\Lambda}( \mathbf{N}, \lambda )$ are $k.p \times k.p$ upper triangular matrices. and, finally, solve the consistent linear systems
\begin{equation} \label{eq:neq_gp_lm_alg23}
\Delta( \mathbf{N}, \lambda ) d\mathbf{a}_{gp-lm} =  \mathbf{M}( \mathbf{a}  )^{T} \mathbf{r}(\mathbf{a})
 \end{equation}
and
\begin{equation} \label{eq:neq_k_lm_alg23}
\Lambda( \mathbf{N},  \lambda ) d\mathbf{a}_{k-lm} =  \mathbf{M}( \mathbf{a}  )^{T} \mathbf{r}(\mathbf{a}) \  ,
\end{equation}
using these Cholesky factorizations. These normal equations  always have an unique solution independently of the value of $\lambda$ if we assume that $\emph{rank}(\mathbf{A}) = k$ and $\emph{rank}( \mathit{J}( \mathbf{r}(\mathbf{a}) ) ) = \emph{rank}( \mathbf{M}(\mathbf{a}) ) = k.(p-k)$ during the iterations (or at least that these conditions are verified as soon as $\lambda = 0$). If $\lambda = 0$, we simply obtained $d\mathbf{a}_{gp-lm} = d\mathbf{a}_{gp-gn}$ and $d\mathbf{a}_{k-lm} = d\mathbf{a}_{k-gn}$ as in version~\eqref{lm_alg1:box} of the Levenberg-Marquardt algorithm when $\Vert \nabla \psi(\mathbf{a} ) \Vert_{2} < \Vert \nabla\psi \Vert_{min}$. However, the key-difference with version~\eqref{lm_alg1:box} of the Levenberg-Marquardt algorithms is that these two linear systems remain nonsingular and well-conditioned when $\lambda$ approaches zero thanks to the addition of the term  $\mathbf{N} \mathbf{N}^{T}$ in the coefficient matrix of these normal equations. In other words, in these conditions,  $\mathbf{R}_{\Delta}( \mathbf{N}, \lambda )$ and $\mathbf{R}_{\Lambda}( \mathbf{N}, \lambda )$ are nonsingular upper triangular matrices and the normal equations can be solved by forward and backward substitutions to get $d\mathbf{a}_{gp-lm}$ and $d\mathbf{a}_{k-lm}$ in step \textbf{(6)} of the Levenberg-Marquardt algorithms~\eqref{lm_alg2:box} and~\eqref{lm_alg3:box}, respectively,:
\begin{equation*}
d\mathbf{a}_{gp-lm} = \mathbf{R}_{\Delta}(\mathbf{N},  \lambda )^{-1}  \mathbf{R}_{\Delta}( \mathbf{N}, \lambda )^{-T} \mathbf{M}(\mathbf{a} )^{T} \mathbf{r}(\mathbf{a})
\end{equation*}
or
\begin{equation*}
d\mathbf{a}_{k-lm} =  \mathbf{R}_{\Lambda}( \mathbf{N},  \lambda )^{-1}  \mathbf{R}_{\Lambda}( \mathbf{N},  \lambda )^{-T}  \mathbf{M}(\mathbf{a} )^{T} \mathbf{r}(\mathbf{a})  \ .
\end{equation*}
Note that, as this step \textbf{(6)} of the Levenberg-Marquardt algorithms~\eqref{lm_alg2:box} and~\eqref{lm_alg3:box} has to be performed several times with different values of $\lambda$, but the same matrix $\mathbf{N}$, for some particular iterations $i$ of these algorithms, it is convenient to first add the cross-product matrix $\mathbf{N} \mathbf{N}^{T}$ to $\Delta$ or $\Lambda$ at each iteration and, then update only the diagonal of these intermediate matrices with $\lambda \mathbf{D}^{2}$ before computing the Cholesky decomposition for a new $\lambda$ value at each subiteration $j$ of the algorithms.

If we want to use a more accurate QR approach in the Levenberg-Marquardt algorithms~\eqref{lm_alg1:box},~\eqref{lm_alg2:box} and~\eqref{lm_alg3:box}, the key-observation is to recognize that the linear systems~\eqref{eq:neq_gp_lm_alg1},~\eqref{eq:neq_k_lm_alg1},~\eqref{eq:neq_gp_lm_alg23} and~\eqref{eq:neq_k_lm_alg23} solved in the above Cholesky approach are, respectively, the normal equations of the damped linear least-squares problems
\begin{equation} \label{eq:llsq_D_J}
d\mathbf{a}_{gp-lm} = \text{Arg} \min_{d\mathbf{a} \in \mathbb{R}^{p.k}}   \,    \big\Vert  \begin{bmatrix} \mathbf{r}(  \mathbf{a} ) \\ \mathbf{0}^{p.k} \end{bmatrix}  -  \begin{bmatrix}   \mathbf{M}( \mathbf{a} ) + \mathbf{L}( \mathbf{a}  )    \\  \sqrt{\lambda} \mathbf{D}  \end{bmatrix}  d\mathbf{a}  \big\Vert^{2}_{2}
\end{equation}
\begin{equation}  \label{eq:llsq_D_M}
d\mathbf{a}_{k-lm} = \text{Arg} \min_{d\mathbf{a} \in \mathbb{R}^{p.k}}   \,    \big\Vert  \begin{bmatrix} \mathbf{r}(  \mathbf{a} ) \\ \mathbf{0}^{p.k} \end{bmatrix}  -  \begin{bmatrix}   \mathbf{M}( \mathbf{a}  )   \\  \sqrt{\lambda} \mathbf{D} \end{bmatrix}  d\mathbf{a}  \big\Vert^{2}_{2}  \ , 
\end{equation}
for the Levenberg-Marquardt algorithms~\eqref{lm_alg1:box}, and of the constrained and damped linear least-squares problems
\begin{equation} \label{eq:llsq_N_D_J}
d\mathbf{a}_{gp-lm} = \text{Arg} \min_{d\mathbf{a} \in \mathbb{R}^{p.k}}   \,    \big\Vert  \begin{bmatrix} \mathbf{r}(  \mathbf{a} ) \\ \mathbf{0}^{k.k} \\ \mathbf{0}^{p.k} \end{bmatrix}  -  \begin{bmatrix}   \mathbf{M}( \mathbf{a} ) + \mathbf{L}( \mathbf{a}  )   \\  \mathbf{N}^{T}  \\  \sqrt{\lambda} \mathbf{D}  \end{bmatrix}  d\mathbf{a}  \big\Vert^{2}_{2}
\end{equation}
\begin{equation}  \label{eq:llsq_N_D_M}
d\mathbf{a}_{k-lm} = \text{Arg} \min_{d\mathbf{a} \in \mathbb{R}^{p.k}}   \,    \big\Vert  \begin{bmatrix} \mathbf{r}(  \mathbf{a} ) \\ \mathbf{0}^{k.k} \\ \mathbf{0}^{p.k} \end{bmatrix}  -  \begin{bmatrix}   \mathbf{M}( \mathbf{a}  )  \\  \mathbf{N}^{T}  \\  \sqrt{\lambda} \mathbf{D} \end{bmatrix}  d\mathbf{a}  \big\Vert^{2}_{2}  \ ,
\end{equation}
for the Levenberg-Marquardt algorithms~\eqref{lm_alg2:box} and~\eqref{lm_alg3:box}.

Thus, the linear least-squares problems~\eqref{eq:llsq_D_J},~\eqref{eq:llsq_D_M},~\eqref{eq:llsq_N_D_J} and~\eqref{eq:llsq_N_D_M} can be solved by computing, respectively, a thin QR decomposition of the "damped" Jacobian matrices
\begin{equation} \label{eq:mat_D_J}
\mathit{J}( \mathbf{r}(\mathbf{a}) ) ( \lambda ) = \begin{bmatrix}   \mathbf{M}( \mathbf{a}  ) + \mathbf{L}( \mathbf{a}  )    \\  \sqrt{\lambda} \mathbf{D}  \end{bmatrix}  \quad \text{or}  \quad \mathbf{M}( \mathbf{a}  ) ( \lambda ) = \begin{bmatrix}   \mathbf{M}( \mathbf{a}  )   \\  \sqrt{\lambda} \mathbf{D}  \end{bmatrix} \ ,
\end{equation}
for the Levenberg-Marquardt algorithms~\eqref{lm_alg1:box} and a thin QR decomposition of the "constrained" and damped" Jacobian matrices
\begin{equation} \label{eq:mat_N_D_J}
\mathit{J}( \mathbf{r}(\mathbf{a}) ) (\mathbf{N}, \lambda ) = \begin{bmatrix}   \mathbf{M}( \mathbf{a}  ) + \mathbf{L}( \mathbf{a}  )   \\  \mathbf{N}^{T}  \\  \sqrt{\lambda} \mathbf{D}  \end{bmatrix}  \quad \text{or}  \quad \mathbf{M}( \mathbf{a}  ) ( \mathbf{N}, \lambda ) = \begin{bmatrix}   \mathbf{M}( \mathbf{a}  )   \\  \mathbf{N}^{T} \\  \sqrt{\lambda} \mathbf{D}  \end{bmatrix} \ ,
\end{equation}
for the Levenberg-Marquardt algorithms~\eqref{lm_alg2:box} and~\eqref{lm_alg3:box}. Furthermore, these different thin QR decompositions can again be done in several steps if we take into account the block-column structure of these constrained and damped Jacobian matrices in order to reduce the memory footprint of the algorithms.

In the first stage, for all the algorithms, we compute the QR decomposition of $\mathbf{M}( \mathbf{a}  ) + \mathbf{L}( \mathbf{a}  )$ or $\mathbf{M}( \mathbf{a}  )$ (without column pivoting) with the same TSQR algorithms as used in the Gauss-Newton methods. This produces implicitly the thin QR factorizations
\begin{equation*}
  \mathbf{M}( \mathbf{a}  ) + \mathbf{L}( \mathbf{a}  ) =  \mathbf{Q}_{\mathit{J}} \mathbf{R}_{\mathit{J}} \text{ or } \mathbf{M}( \mathbf{a}  ) =  \mathbf{Q}_{\mathbf{M}} \mathbf{R}_{\mathbf{M}}  \ ,
\end{equation*}
where $\mathbf{Q}_{\mathit{J}}$ and $\mathbf{Q}_{\mathbf{M}}$ are $n.p \times k.p$ matrices with orthonormal columns, and, $\mathbf{R}_{\mathit{J}}$ and $\mathbf{R}_{\mathbf{M}}$ are $k.p \times k.p$ singular upper triangular matrices as discussed in the previous subsection.

Next, in step \textbf{(6.1)} of the Levenberg-Marquardt algorithms~\eqref{lm_alg1:box} in which $\lambda > 0$, we first note, using the above thin QR decomposition of $\mathbf{M}( \mathbf{a}  ) + \mathbf{L}( \mathbf{a}  )$ or $\mathbf{M}( \mathbf{a}  )$, that we have
\begin{equation*}
\mathit{J} \big ( \mathbf{r}(\mathbf{a}) \big ) (\lambda)  =  \begin{bmatrix}    \mathbf{Q}_{\mathit{J}} \mathbf{R}_{\mathit{J}}  \\    \sqrt{\lambda} \mathbf{D} \end{bmatrix}  = \begin{bmatrix}    \mathbf{Q}_{\mathit{J}} &   \mathbf{0}^{ n.p \times k.p }  \\   \mathbf{0}^{k.p \times k.p }  &   \mathbf{I}_{k.p} \end{bmatrix}  \begin{bmatrix}  \mathbf{R}_{\mathit{J}} \\    \sqrt{\lambda} \mathbf{D}  \end{bmatrix}
\end{equation*}
and 
\begin{equation*}
\mathbf{M}( \mathbf{a}  ) (\lambda)  =  \begin{bmatrix}     \mathbf{Q}_{\mathbf{M}}  \mathbf{R}_{\mathbf{M}}  \\    \sqrt{\lambda} \mathbf{D} \end{bmatrix}  = \begin{bmatrix}    \mathbf{Q}_{\mathbf{M}} &   \mathbf{0}^{ n.p \times k.p }  \\   \mathbf{0}^{k.p \times k.p }  &   \mathbf{I}_{k.p} \end{bmatrix}  \begin{bmatrix}  \mathbf{R}_{\mathbf{M}} \\    \sqrt{\lambda} \mathbf{D}  \end{bmatrix}  \ .
\end{equation*}
Thus, in a second stage, we can factorize the block-column matrices $\begin{bmatrix}  \mathbf{R}_{\mathit{J}} \\    \sqrt{\lambda} \mathbf{D}  \end{bmatrix}$ or $\begin{bmatrix} \mathbf{R}_{\mathbf{M}} \\    \sqrt{\lambda} \mathbf{D}  \end{bmatrix}$ on the right hand side of these equations into the product of an orthogonal matrix times a rectangular matrix in upper triangular form. This can be done efficiently with a sequence of $k.p.(k.p+1)$ Givens rotations applied to the left of these block-column matrices to annihilate their bottom diagonal block, $\sqrt{\lambda} \mathbf{D}$. These Givens rotations use the diagonal elements of $\mathbf{R}_{\mathit{J}}$ and $\mathbf{R}_{\mathbf{M}}$ to eliminate the diagonal elements of $\sqrt{\lambda} \mathbf{D}$ and reduce the fill-in in this process. Note further that, in this recursive process, the bands of zeros introduced in the previous stages are unaffected by the subsequent stages thanks to the use of Givens rotations in the calculations; see Section 10.3 of Nocedal and Wright~\cite{NW2006} for a good account of this computing scheme. At the end, we get  the matrix equations
\begin{equation*}
\mathbf{W}_{\mathit{J}} (\lambda) \begin{bmatrix}  \mathbf{R}_{\mathit{J}} \\    \sqrt{\lambda} \mathbf{D}  \end{bmatrix} = \begin{bmatrix}  \mathbf{R}_{\mathit{J}} (\lambda) \\   \mathbf{0}^{k.p \times k.p }   \end{bmatrix}
 \text{ or }
\mathbf{W}_{\mathbf{M}} (\lambda) \begin{bmatrix}  \mathbf{R}_{\mathbf{M}} \\    \sqrt{\lambda} \mathbf{D}  \end{bmatrix} = \begin{bmatrix}  \mathbf{R}_{\mathbf{M}} (\lambda) \\   \mathbf{0}^{k.p \times k.p }   \end{bmatrix} \ .
\end{equation*}
Here, $\mathbf{W}_{\mathit{J}} (\lambda)$  and $\mathbf{W}_{\mathbf{M}} (\lambda)$ are $2.k.p  \times 2.k.p$ orthogonal matrices, which are the products of $k.p.(k.p+1)$ Givens rotations, and, $\mathbf{R}_{\mathit{J}} (\lambda)$ and $\mathbf{R}_{\mathbf{M}} (\lambda)$ are nonsingular $k.p  \times k.p$ upper triangular matrices as $\lambda > 0$ and $\mathbf{D}$ has no zero elements on its diagonal.

Proceeding in this way, we implicitly build up thin QR decompositions of $\mathit{J}( \mathbf{r}(\mathbf{a}) ) (\lambda)$ or $\mathbf{M}( \mathbf{a}  ) (\lambda)$ in several stages since
\begin{align*}
\mathit{J} \big ( \mathbf{r}(\mathbf{a}) \big ) (\lambda)  & =  \begin{bmatrix}    \mathbf{Q}_{\mathit{J}} &   \mathbf{0}^{ n.p \times k.p }  \\   \mathbf{0}^{k.p \times k.p }  &   \mathbf{I}_{k.p} \end{bmatrix}  \mathbf{W}_{\mathit{J}} (\lambda)^{T} \begin{bmatrix}  \mathbf{R}_{\mathit{J}} (\lambda) \\   \mathbf{0}^{k.p \times k.p }   \end{bmatrix}  \\
                                      & =  \begin{bmatrix}    \mathbf{Q}_{\mathit{J}} \mathbf{W}_{\mathit{J}}^{11} (\lambda)^{T}  \\    \mathbf{W}_{\mathit{J}}^{12} (\lambda)^{T}   \end{bmatrix}   \mathbf{R}_{\mathit{J}} (\lambda) \\
                                      & =  \mathbf{Q}_{\mathit{J}} (\lambda)  \mathbf{R}_{\mathit{J}} (\lambda)
\end{align*}
and 
\begin{align*}
\mathbf{M}( \mathbf{a}  ) (\lambda)  & =   \begin{bmatrix}    \mathbf{Q}_{\mathbf{M}} &   \mathbf{0}^{ n.p \times k.p }  \\   \mathbf{0}^{k.p \times k.p }  &   \mathbf{I}_{k.p} \end{bmatrix} \mathbf{W}_{\mathbf{M}} (\lambda)^{T} \begin{bmatrix}  \mathbf{R}_{\mathbf{M}} (\lambda) \\   \mathbf{0}^{k.p \times k.p }   \end{bmatrix}   \\
                 & =   \begin{bmatrix}    \mathbf{Q}_{\mathbf{M}} \mathbf{W}_{\mathbf{M}}^{11} (\lambda)^{T}  \\   \mathbf{W}_{\mathbf{M}}^{12} (\lambda)^{T} \end{bmatrix} \mathbf{R}_{\mathbf{M}} (\lambda)  \\
                 & =  \mathbf{Q}_{\mathbf{M}} (\lambda) \mathbf{R}_{\mathbf{M}} (\lambda) \ ,
\end{align*}
where the $2.k.p \times 2.k.p$ orthogonal matrices  $\mathbf{W}_{\mathit{J}} (\lambda)$ and $\mathbf{W}_{\mathbf{M}} (\lambda)$ have been partitioned in four blocks of $k.p$ rows and columns each:
\begin{equation*}
\mathbf{W}_{\mathit{J}} (\lambda) = \begin{bmatrix}   \mathbf{W}_{\mathit{J}}^{11} (\lambda)  &   \mathbf{W}_{\mathit{J}}^{12} (\lambda)   \\    \mathbf{W}_{\mathit{J}}^{21} (\lambda) & \mathbf{W}_{\mathit{J}}^{22} (\lambda) \end{bmatrix}
 \text{ , }
\mathbf{W}_{\mathbf{M}} (\lambda) = \begin{bmatrix}  \mathbf{W}_{\mathbf{M}}^{11} (\lambda)  & \mathbf{W}_{\mathbf{M}}^{12} (\lambda)  \\    \mathbf{W}_{\mathbf{M}}^{21} (\lambda)  & \mathbf{W}_{\mathbf{M}}^{22} (\lambda) \end{bmatrix},
\end{equation*}
and the $(p.n + k.p) \times k.p$  matrices $\mathbf{Q}_{\mathit{J}} (\lambda)$ and $\mathbf{Q}_{\mathbf{M}} (\lambda)$ have orthonormal columns since $\mathbf{W}_{\mathit{J}} (\lambda)$ and $\mathbf{W}_{\mathbf{M}} (\lambda)$ are orthogonal matrices, and, the matrices $ \mathbf{Q}_{\mathit{J}}$ and $\mathbf{Q}_{\mathbf{M}}$ have orthonormal columns.

Finally, using these thin QR decompositions of $\mathit{J}( \mathbf{r}(\mathbf{a}) ) (\lambda)$ and $\mathbf{M}( \mathbf{a}  ) (\lambda)$, the solutions $d\mathbf{a}_{gp-lm}$ and $d\mathbf{a}_{k-lm}$ of the damped linear least-square problems~\eqref{eq:llsq_D_J}  and~\eqref{eq:llsq_D_M}, which must be solved in step \textbf{(6.1)} of the Levenberg-Marquardt algorithms~\eqref{lm_alg1:box}, can be easily computed in a last step as $\mathbf{R}_{\mathit{J}} (\lambda)$ and $\mathbf{R}_{\mathbf{M}} (\lambda)$ are nonsingular upper triangular matrices, e.g.,
\begin{equation*}
d\mathbf{a}_{gp-lm} =   \mathbf{R}_{\mathit{J}} (\lambda)^{-1} \mathbf{Q}_{\mathit{J}} (\lambda)^{T} \begin{bmatrix} \mathbf{r}(\mathbf{a})   \\  \mathbf{0}^{ k.p }  \end{bmatrix}  \text{ or }  d\mathbf{a}_{k-lm} = \mathbf{R}_{\mathbf{M}} (\lambda)^{-1} \mathbf{Q}_{\mathbf{M}} (\lambda) ^{T}\begin{bmatrix} \mathbf{r}(\mathbf{a})   \\  \mathbf{0}^{ k.p }  \end{bmatrix} \ .
\end{equation*}

Next, if we want to use a QR approach in the Levenberg-Marquardt algorithms~\eqref{lm_alg2:box} and~\eqref{lm_alg3:box}, we have to solve the constrained and damped linear least-squares problems~\eqref{eq:llsq_N_D_J} or~\eqref{eq:llsq_N_D_M} at each iteration. This can also be done by computing the thin QR decomposition of the block-column matrices $\mathit{J}( \mathbf{r}(\mathbf{a}) ) ( \mathbf{N},  \lambda )$ and $\mathbf{M}( \mathbf{a}  ) ( \mathbf{N}, \lambda )$ defined in equation~\eqref{eq:mat_N_D_J}  in several steps to reduce the memory footprint of the algorithms. More precisely with one more steps compared to the structured QR algorithm used in step \textbf{(6.1)} of the Levenberg-Marquardt algorithms~\eqref{lm_alg1:box} to reduce the matrices $\mathit{J}( \mathbf{r}(\mathbf{a}) ) ( \lambda )$ and $\mathbf{M}( \mathbf{a}  ) (\lambda )$ (defined in equation~\eqref{eq:mat_D_J}) to triangular form.

Thus, after using the same first stage as before to get the QR decomposition of $\mathbf{M}( \mathbf{a}  ) + \mathbf{L}( \mathbf{a}  )$ or $\mathbf{M}( \mathbf{a}  )$ with the TSQR algorithms, we next compute implicitly the thin QR factorizations of
\begin{equation*}
\mathit{J} \big ( \mathbf{r}(\mathbf{a}) \big ) ( \mathbf{N} ) = \begin{bmatrix}   \mathbf{M}( \mathbf{a}  ) + \mathbf{L}( \mathbf{a}  )    \\  \mathbf{N}^{T}  \end{bmatrix} = \mathbf{Q}_{\mathit{J}} (\mathbf{N})  \mathbf{R}_{\mathit{J}} (\mathbf{N})
\end{equation*}
or
\begin{equation*}
\mathbf{M}( \mathbf{a}  ) ( \mathbf{N} ) = \begin{bmatrix}   \mathbf{M}( \mathbf{a}  )    \\  \mathbf{N}^{T}   \end{bmatrix} = \mathbf{Q}_{\mathbf{M}} (\mathbf{N}) \mathbf{R}_{\mathbf{M}} (\mathbf{N}) \ ,
\end{equation*}
where $\mathbf{Q}_{\mathit{J}} (\mathbf{N})$ and  $\mathbf{Q}_{\mathbf{M}} (\mathbf{N})$ are $(p.n + k.k) \times k.p$ matrices with orthonormal columns and $\mathbf{R}_{\mathit{J}} (\mathbf{N})$ and $\mathbf{R}_{\mathbf{M}} (\mathbf{N})$ are $k.p \times k.p$ upper triangular matrices. Since
\begin{equation*}
\mathit{J} \big ( \mathbf{r}(\mathbf{a}) \big ) ( \mathbf{N} ) =  \begin{bmatrix}   \mathbf{Q}_{\mathit{J}}  \mathbf{R}_{\mathit{J}}   \\  \mathbf{N}^{T}  \end{bmatrix}  = \begin{bmatrix}  \mathbf{Q}_{\mathit{J}} &  \mathbf{0}^{ n.p \times k.k}   \\   \mathbf{0}^{k.k \times k.p}   & \mathbf{I}_{k.k}  \end{bmatrix} \begin{bmatrix}  \mathbf{R}_{\mathit{J}} \\  \mathbf{N}^{T} \end{bmatrix}
\end{equation*}
and 
\begin{equation*}
\mathbf{M}( \mathbf{a}  ) ( \mathbf{N} ) = \begin{bmatrix} \mathbf{Q}_{\mathbf{M}} \mathbf{R}_{\mathbf{M}}   \\  \mathbf{N}^{T}  \end{bmatrix} = \begin{bmatrix}  \mathbf{Q}_{\mathbf{M}}  &  \mathbf{0}^{ n.p \times k.k}  \\   \mathbf{0}^{k.k \times k.p}  & \mathbf{I}_{k.k}   \end{bmatrix} \begin{bmatrix}  \mathbf{R}_{\mathbf{M}}  \\  \mathbf{N}^{T} \end{bmatrix} \ ,
\end{equation*}
this can be done by performing structured and thin QR factorizations of the matrices $\begin{bmatrix}  \mathbf{R}_{\mathit{J}} \\  \mathbf{N}^{T} \end{bmatrix}$ and $\begin{bmatrix}  \mathbf{R}_{\mathbf{M}} \\  \mathbf{N}^{T} \end{bmatrix}$; more precisely, by applying a sequence of $k.p$ dedicated Householder transformations on the left of the matrices. These Householder transformations are designed to annihilate the lower block $\mathbf{N}^{T}$ of these matrices, giving the matrix equalities,
\begin{equation*}
\mathbf{W}_{\mathit{J}} (  \mathbf{N} ) \begin{bmatrix}  \mathbf{R}_{\mathit{J}}  \\   \mathbf{N}^{T}  \end{bmatrix} = \begin{bmatrix}  \mathbf{R}_{\mathit{J}} ( \mathbf{N}) \\   \mathbf{0}^{k.k \times k.p }   \end{bmatrix}
\quad  \text{or}  \quad
\mathbf{W}_{\mathbf{M}} (  \mathbf{N} ) \begin{bmatrix}  \mathbf{R}_{\mathbf{M}}  \\   \mathbf{N}^{T}  \end{bmatrix} = \begin{bmatrix}  \mathbf{R}_{\mathbf{M}} ( \mathbf{N} ) \\   \mathbf{0}^{k.k \times k.p }   \end{bmatrix} \ ,
\end{equation*}
where $\mathbf{W}_{\mathit{J}} (  \mathbf{N} )$ and $\mathbf{W}_{\mathbf{M}} (  \mathbf{N} )$ are $k.(p + k ) \times k.(p + k)$ orthogonal matrices composed of the product of these $k.p$ elementary Householder transformations, and $\mathbf{R}_{\mathit{J}} ( \mathbf{N})$ and $\mathbf{R}_{\mathbf{M}} ( \mathbf{N} )$ are nonsingular upper triangular matrices if $\emph{rank}( \mathbf{A}) = k$ and  $\emph{rank}( \mathbf{M}( \mathbf{a}  ) + \mathbf{L}( \mathbf{a}  ) ) = \emph{rank}( \mathbf{M}( \mathbf{a}  ) ) = k.(p-k)$. This reduction to triangular from is exactly similar to one of the steps of the serial TSQR algorithm described in Subsection~\ref{vp_gn_alg:box}, after the first one. Note also that, in this recursive process, the band of zeros introduced in the previous stages or preceding Householder transformations are unaffected by the subsequent stages if dedicated Householder transformations are used in the calculations.

Using this computational sequence, we finally obtain a thin QR factorization of $\mathit{J}( \mathbf{r}(\mathbf{a}) ) ( \mathbf{N} )$ or $\mathbf{M}( \mathbf{a}  ) ( \mathbf{N} ) $ since
\begin{align*}
\mathit{J} \big( \mathbf{r}(\mathbf{a}) \big) ( \mathbf{N} ) & =   \begin{bmatrix}  \mathbf{Q}_{\mathit{J}}  &  \mathbf{0}^{ n.p  \times k.k}   \\   \mathbf{0}^{k.k \times k.p}   & \mathbf{I}_{k.k}  \end{bmatrix} \mathbf{W}_{\mathit{J}} (  \mathbf{N} )^{T}  \begin{bmatrix}  \mathbf{R}_{\mathit{J}} ( \mathbf{N}) \\   \mathbf{0}^{k.k \times k.p }   \end{bmatrix}   \\
             & =    \begin{bmatrix}    \mathbf{Q}_{\mathit{J}} \mathbf{W}_{\mathit{J}}^{11} ( \mathbf{N} )^{T}  \\    \mathbf{W}_{\mathit{J}}^{12} ( \mathbf{N} )^{T}   \end{bmatrix}    \mathbf{R}_{\mathit{J}} ( \mathbf{N})  \\
            & =  \mathbf{Q}_{\mathit{J}} ( \mathbf{N} ) \mathbf{R}_{\mathit{J}} ( \mathbf{N})
\end{align*}
and
\begin{align*}
\mathbf{M}( \mathbf{a}  ) ( \mathbf{N} ) & =   \begin{bmatrix}  \mathbf{Q}_{\mathbf{M}} &  \mathbf{0}^{ n.p \times k.k}  \\   \mathbf{0}^{k.k \times k.p}  & \mathbf{I}_{k.k}   \end{bmatrix} \mathbf{W}_{\mathbf{M}} ( \mathbf{N} )^{T} \begin{bmatrix}  \mathbf{R}_{\mathbf{M}} ( \mathbf{N} ) \\   \mathbf{0}^{k.k \times k.p }   \end{bmatrix}  \\
        & =    \begin{bmatrix}    \mathbf{Q}_{\mathbf{M}} \mathbf{W}_{\mathbf{M}}^{11} ( \mathbf{N} )^{T}  \\   \mathbf{W}_{\mathbf{M}}^{12} ( \mathbf{N} )^{T} \end{bmatrix}  \mathbf{R}_{\mathbf{M}} ( \mathbf{N} )   \\
        & =  \mathbf{Q}_{\mathbf{M}} (  \mathbf{N} ) \mathbf{R}_{\mathbf{M}} ( \mathbf{N} ) \ ,
\end{align*}
where the $k.(p + k) \times k.(p + k)$ orthogonal matrices $\mathbf{W}_{\mathit{J}} (  \mathbf{N} )$ and $\mathbf{W}_{\mathbf{M}} (  \mathbf{N} )$ have been partitioned in four blocks as
\begin{equation*}
\mathbf{W}_{\mathit{J}} (\mathbf{N}) = \begin{bmatrix}   \mathbf{W}_{\mathit{J}}^{11} (\mathbf{N})  &   \mathbf{W}_{\mathit{J}}^{12} (\mathbf{N})   \\    \mathbf{W}_{\mathit{J}}^{21} (\mathbf{N}) & \mathbf{W}_{\mathit{J}}^{22} (\mathbf{N}) \end{bmatrix}
 \text{ and }
\mathbf{W}_{\mathbf{M}} (\mathbf{N}) = \begin{bmatrix}  \mathbf{W}_{\mathbf{M}}^{11} (\mathbf{N})  & \mathbf{W}_{\mathbf{M}}^{12} (\mathbf{N})  \\    \mathbf{W}_{\mathbf{M}}^{21} (\mathbf{N})  & \mathbf{W}_{\mathbf{M}}^{22} (\mathbf{N}) \end{bmatrix}  \ ,
\end{equation*}
where
\begin{itemize}
\item[]  $\mathbf{W}_{\mathit{J}}^{11} (\mathbf{N}), \mathbf{W}_{\mathbf{M}}^{11} (\mathbf{N}) \in  \mathbb{R}^{k.p \times k.p}$ ,
\item[]  $\mathbf{W}_{\mathit{J}}^{12} (\mathbf{N}), \mathbf{W}_{\mathbf{M}}^{12} (\mathbf{N}) \in  \mathbb{R}^{k.p \times k.k}$ ,
\item[]  $\mathbf{W}_{\mathit{J}}^{21} (\mathbf{N}), \mathbf{W}_{\mathbf{M}}^{21} (\mathbf{N}) \in  \mathbb{R}^{k.k \times k.p}$ ,
\item[]  $\mathbf{W}_{\mathit{J}}^{22} (\mathbf{N}), \mathbf{W}_{\mathbf{M}}^{22} (\mathbf{N}) \in  \mathbb{R}^{k.k \times k.k}$ ,
\end{itemize}
and the $(n.p + k.k) \times k.p$ matrices $\mathbf{Q}_{\mathit{J}} ( \mathbf{N} )$ and $\mathbf{Q}_{\mathbf{M}} ( \mathbf{N})$ have orthonormal columns since $\mathbf{W}_{\mathit{J}} ( \mathbf{N})$ and $\mathbf{W}_{\mathbf{M}} ( \mathbf{N} )$ are  orthogonal matrices and the $n.p \times k.p$ matrices $\mathbf{Q}_{\mathit{J}}$ and $\mathbf{Q}_{\mathbf{M}}$ have orthonormal columns.

Using these results, we can write, finally,
\begin{equation*}
\mathit{J} \big ( \mathbf{r}(\mathbf{a}) \big) (\mathbf{N}, \lambda ) =  \begin{bmatrix}  \mathit{J}( \mathbf{r}(\mathbf{a}) ) ( \mathbf{N} )  \\  \sqrt{\lambda} \mathbf{D}  \end{bmatrix} = \begin{bmatrix}  \mathbf{Q}_{\mathit{J}} ( \mathbf{N} ) \mathbf{R}_{\mathit{J}} ( \mathbf{N}) \\  \sqrt{\lambda} \mathbf{D}  \end{bmatrix} = \begin{bmatrix}    \mathbf{Q}_{\mathit{J}}  ( \mathbf{N} )&   \mathbf{0}^{ ( n.p + k.k ) \times k.p }  \\   \mathbf{0}^{k.p \times k.p }  &   \mathbf{I}_{k.p} \end{bmatrix}  \begin{bmatrix}  \mathbf{R}_{\mathit{J}}  ( \mathbf{N} ) \\    \sqrt{\lambda} \mathbf{D}  \end{bmatrix}
\end{equation*}
and
\begin{equation*}
\mathbf{M}( \mathbf{a}  ) ( \mathbf{N}, \lambda ) = \begin{bmatrix}  \mathbf{M}( \mathbf{a}  ) ( \mathbf{N} ) \\  \sqrt{\lambda} \mathbf{D}  \end{bmatrix} = \begin{bmatrix}  \mathbf{Q}_{\mathbf{M}} (  \mathbf{N} ) \mathbf{R}_{\mathbf{M}} ( \mathbf{N} ) \\  \sqrt{\lambda} \mathbf{D}  \end{bmatrix}  = \begin{bmatrix}    \mathbf{Q}_{\mathbf{M}} (  \mathbf{N} ) &   \mathbf{0}^{ ( n.p + k.k ) \times k.p }  \\   \mathbf{0}^{k.p \times k.p }  &   \mathbf{I}_{k.p} \end{bmatrix}  \begin{bmatrix}  \mathbf{R}_{\mathbf{M}} (  \mathbf{N} ) \\    \sqrt{\lambda} \mathbf{D}  \end{bmatrix} \ ,
\end{equation*}
and the thin QR decompositions of $\mathit{J}( \mathbf{r}(\mathbf{a}) ) (\mathbf{N}, \lambda )$ and $\mathbf{M}( \mathbf{a}  ) ( \mathbf{N}, \lambda )$ can be obtained by computing the structured and thin QR decompositions of the column-block matrices $\begin{bmatrix}  \mathbf{R}_{\mathit{J}}  ( \mathbf{N} ) \\    \sqrt{\lambda} \mathbf{D}  \end{bmatrix}$ and $\begin{bmatrix}  \mathbf{R}_{\mathbf{M}} (  \mathbf{N} ) \\    \sqrt{\lambda} \mathbf{D}  \end{bmatrix}$ using the same computing sequence as performed in the last stage of step \textbf{(6.1)} of the Levenberg-Marquardt algorithm~\eqref{lm_alg1:box} described above.

In other words, the elements of the diagonal matrix  $\sqrt{\lambda} \mathbf{D}$ can be eliminated by a sequence of $k.p.(k.p+1)$ Givens rotations and, at the end of this process, we get
\begin{equation*}
\mathbf{W}_{\mathit{J}} (\mathbf{N}, \lambda) \begin{bmatrix}  \mathbf{R}_{\mathit{J}} (\mathbf{N}) \\    \sqrt{\lambda} \mathbf{D}  \end{bmatrix} = \begin{bmatrix}  \mathbf{R}_{\mathit{J}} (\mathbf{N}, \lambda) \\   \mathbf{0}^{k.p \times k.p }   \end{bmatrix}
 \text{ or } 
\mathbf{W}_{\mathbf{M}} (\mathbf{N}, \lambda) \begin{bmatrix}  \mathbf{R}_{\mathbf{M}} (\mathbf{N}) \\    \sqrt{\lambda} \mathbf{D}  \end{bmatrix} = \begin{bmatrix}  \mathbf{R}_{\mathbf{M}} (\mathbf{N}, \lambda) \\   \mathbf{0}^{k.p \times k.p }   \end{bmatrix} \ .
\end{equation*}
where $\mathbf{W}_{\mathit{J}} (\mathbf{N}, \lambda)$  and $\mathbf{W}_{\mathbf{M}} (\mathbf{N}, \lambda)$ are $2.k.p  \times 2.k.p$ orthogonal matrices, which are the products of these $k.p.(k.p+1)$ Givens rotations, and, $\mathbf{R}_{\mathit{J}} (\mathbf{N}, \lambda)$ and $\mathbf{R}_{\mathbf{M}} (\mathbf{N}, \lambda)$ are nonsingular $k.p  \times k.p$ upper triangular matrices. Finally, using these matrices equalities, we obtain the thin QR factorizations of  $\mathit{J}( \mathbf{r}(\mathbf{a}) ) ( \mathbf{N}, \lambda )$ and $\mathbf{M}( \mathbf{a}  ) ( \mathbf{N}, \lambda )$ since
\begin{align*}
\mathit{J}  \big ( \mathbf{r}(\mathbf{a})  \big ) (\mathbf{N}, \lambda)  & =  \begin{bmatrix}    \mathbf{Q}_{\mathit{J}}(\mathbf{N}) &   \mathbf{0}^{ (n.p + k.k) \times k.p }  \\   \mathbf{0}^{k.p \times k.p }  &   \mathbf{I}_{k.p} \end{bmatrix}  \mathbf{W}_{\mathit{J}} (\mathbf{N}, \lambda)^{T} \begin{bmatrix}  \mathbf{R}_{\mathit{J}} (\mathbf{N}, \lambda) \\   \mathbf{0}^{k.p \times k.p }   \end{bmatrix}  \\
                                      & =  \begin{bmatrix}    \mathbf{Q}_{\mathit{J}} (\mathbf{N}) \mathbf{W}_{\mathit{J}}^{11} (\mathbf{N}, \lambda)^{T}  \\    \mathbf{W}_{\mathit{J}}^{12} (\mathbf{N}, \lambda)^{T}   \end{bmatrix}   \mathbf{R}_{\mathit{J}} (\mathbf{N}, \lambda) \\
                                      & =  \mathbf{Q}_{\mathit{J}} (\mathbf{N}, \lambda)  \mathbf{R}_{\mathit{J}} (\mathbf{N}, \lambda)
\end{align*}
and 
\begin{align*}
\mathbf{M}( \mathbf{a}  ) (\mathbf{N}, \lambda)  & =   \begin{bmatrix}    \mathbf{Q}_{\mathbf{M}} (\mathbf{N}) &   \mathbf{0}^{ (n.p + k.k) \times k.p }  \\   \mathbf{0}^{k.p \times k.p }  &   \mathbf{I}_{k.p} \end{bmatrix} \mathbf{W}_{\mathbf{M}} (\mathbf{N}, \lambda)^{T} \begin{bmatrix}  \mathbf{R}_{\mathbf{M}} (\mathbf{N}, \lambda) \\   \mathbf{0}^{k.p \times k.p }   \end{bmatrix}   \\
                 & =   \begin{bmatrix}    \mathbf{Q}_{\mathbf{M}} (\mathbf{N}) \mathbf{W}_{\mathbf{M}}^{11} (\mathbf{N}, \lambda)^{T}  \\   \mathbf{W}_{\mathbf{M}}^{12} (\mathbf{N}, \lambda)^{T} \end{bmatrix} \mathbf{R}_{\mathbf{M}} (\mathbf{N}, \lambda)  \\
                 & =  \mathbf{Q}_{\mathbf{M}} (\mathbf{N}, \lambda) \mathbf{R}_{\mathbf{M}} (\mathbf{N}, \lambda) \ ,
\end{align*}
where the $2.k.p \times 2.k.p$ orthogonal matrices  $\mathbf{W}_{\mathit{J}} (\mathbf{N}, \lambda)$ and $\mathbf{W}_{\mathbf{M}} (\mathbf{N}, \lambda)$ have been partitioned in four blocks of $k.p$ rows and columns each, and the $(p.n + k.k + k.p) \times k.p$  matrices $\mathbf{Q}_{\mathit{J}} (\mathbf{N}, \lambda)$ and $\mathbf{Q}_{\mathbf{M}} (\mathbf{N}, \lambda)$ have orthonormal columns since $\mathbf{W}_{\mathit{J}} (\mathbf{N}, \lambda)$ and $\mathbf{W}_{\mathbf{M}} (\mathbf{N}, \lambda)$ are orthogonal matrices, and, the matrices $ \mathbf{Q}_{\mathit{J}} (\mathbf{N})$ and $\mathbf{Q}_{\mathbf{M}} (\mathbf{N})$ have orthonormal columns.

Using these thin QR factorizations of $\mathit{J}( \mathbf{r}(\mathbf{a}) ) ( \mathbf{N}, \lambda )$ and $\mathbf{M}( \mathbf{a}  ) ( \mathbf{N}, \lambda )$, the correction vectors $d\mathbf{a}_{gp-lm}$ and $d\mathbf{a}_{k-lm}$ in step \textbf{(6)} of the Levenberg-Marquardt algorithms~\eqref{lm_alg2:box} and~\eqref{lm_alg3:box} can then be computed as
\begin{equation*}
d\mathbf{a}_{gp-lm} =   \mathbf{R}_{\mathit{J}} (\mathbf{N}, \lambda)^{-1} \mathbf{Q}_{\mathit{J}} (\mathbf{N}, \lambda)^{T}  \begin{bmatrix} \mathbf{r}(\mathbf{a})   \\  \mathbf{0}^{ k.k + k.p }  \end{bmatrix}
\end{equation*}
and
\begin{equation*}
d\mathbf{a}_{k-lm} = \mathbf{R}_{\mathbf{M}} (\mathbf{N}, \lambda)^{-1} \mathbf{Q}_{\mathbf{M}} (\mathbf{N}, \lambda)^{T} \begin{bmatrix} \mathbf{r}(\mathbf{a})   \\  \mathbf{0}^{ k.k + k.p }  \end{bmatrix} \ ,
\end{equation*}
if we use a QR approach in these algorithms.

Note, finally, that if, at a particular iteration $i$ of the Levenberg-Marquardt algorithms~\eqref{lm_alg2:box} and~\eqref{lm_alg3:box}, we have to solve several times the constrained and damped linear least squares problems~\eqref{eq:llsq_N_D_J} and~\eqref{eq:llsq_N_D_M} for the same value of $\mathbf{a}_{i}$, but different values of $\lambda$ in order to ensure the descending condition $\psi(\mathbf{a}_{i+1} ) < \psi(\mathbf{a}_{i} )$ and the convergence of the algorithms, only the last stage of the structured and thin QR factorizations of $\mathit{J}( \mathbf{r}(\mathbf{a}_{i}) ) ( \mathbf{N}_{i}, \lambda )$ and $\mathbf{M}( \mathbf{a}_{i}  ) ( \mathbf{N}_{i}, \lambda )$ involving the damping parameter $\lambda$ has to be performed again as the first two stages remain identical if $\mathbf{a}_{i}$ is not changed. This obviously can save a lot of computing time as $p.n \gg k.p$ for most WLRA problems encountered in practice. This is an interesting feature of the QR approach, see Section 10.3 of Nocedal and Wright~\cite{NW2006} for further discussion on these algorithmic, but important, details in a more general NLLS context.

\subsection{Variable projection Newton and quasi-Newton algorithms} \label{vp_n_alg:box}

For large residuals WLRA problems, the variable projection Gauss-Newton and Levenberg-Marquardt algorithms described in Subsections~\ref{vp_gn_alg:box} and~\ref{vp_lm_alg:box} can be much less efficient as they do not include second-order derivative information from the Hessian matrix $\nabla^2 \psi( \mathbf{a} )$ in their associated quadratic models for the variations of $\psi(.)$ in a neighborhood of $\mathbf{a}$. Thus, their asymptotic convergence rates are expected to be only linear in these conditions~\cite{DS1983}\cite{NW2006}. Fortunately, from the results of Subsection~\ref{hess:box}, we have a compact expression for the Hessian matrix $\nabla^2 \psi( \mathbf{a} )$ (see equation~\eqref{eq:H_hess_mat})
\begin{equation*}
    \mathbf{H} = \mathbf{M}(\mathbf{a})^{T} \mathbf{M}(\mathbf{a}) - \mathbf{L}(\mathbf{a})^{T} \mathbf{L}(\mathbf{a}) + \big( \mathbf{U}(\mathbf{a})^{T} \mathbf{L} (\mathbf{a}) + \mathbf{L}(\mathbf{a})^{T} \mathbf{U}(\mathbf{a}) \big) \ ,
\end{equation*}
under the hypothesis that $ \mathbf{F}(.) $ has full column-rank in a neighborhood of $ \mathbf{a} $, and, also, all the machinery to implement full Newton algorithms based on this exact three-term expression of $\nabla^2 \psi( \mathbf{a} )$ or quasi-Newton methods based on its two-term approximation (see equation~\eqref{eq:approx_hess_mat})
\begin{equation*}
    \bar{\mathbf{H}} = \mathbf{M}(\mathbf{a})^{T} \mathbf{M}(\mathbf{a}) - \mathbf{L}(\mathbf{a})^{T} \mathbf{L}(\mathbf{a})  \ .
\end{equation*}
Here $\mathbf{M}(\mathbf{a}), \mathbf{L}(\mathbf{a})$ and  $\mathbf{U}(\mathbf{a})$ are defined, respectively, in equations~\eqref{eq:M_mat},~\eqref{eq:L_mat} and~\eqref{eq:U_mat}. Note further that, as we assume that $\mathbf{F}(\mathbf{a})$ is of full column rank, e.g.,  ${r}_{\mathbf{F}( \mathbf{a} )} = \emph{rank}( \mathbf{F}(  \mathbf{a} ) ) =k.p$, we have $\mathbf{F}(\mathbf{a})^{+} = \mathbf{F}(\mathbf{a})^-$ and the computation of  the three matrices $\mathbf{M}(\mathbf{a}), \mathbf{L}(\mathbf{a})$ and $\mathbf{U}(\mathbf{a})$ can be simplified accordingly. Finally, the cost of the above quasi-Newton methods is similar to those of the Golub-Pereyra Gauss-Newton algorithms using a Cholesky approach as both methods compute exactly the same symmetric matrix terms, e.g.,  $\mathbf{M}(\mathbf{a})^{T} \mathbf{M}(\mathbf{a})$ and $\mathbf{L}(\mathbf{a})^{T} \mathbf{L}(\mathbf{a})$.

Furthermore, as demonstrated in the Cholesky approach of the Gauss-Newton algorithms~\eqref{gn_alg:box}, the computations of these two symmetric terms can be easily parallelized since
\begin{equation*}
\mathbf{M}(\mathbf{a})^{T} \mathbf{M}(\mathbf{a}) = \sum_{j=1}^n  \mathbf{M}_{j}^{T} \mathbf{M}_{j} \text{ and }  \mathbf{L}(\mathbf{a})^{T} \mathbf{L}(\mathbf{a}) = \sum_{j=1}^n  \mathbf{L}_{j}^{T} \mathbf{L}_{j}  \ ,
\end{equation*}
where $\mathbf{M}(\mathbf{a})$ and  $\mathbf{L}(\mathbf{a})$ have been divided into $n$ blocks $\mathbf{M}_{j}$ and $\mathbf{L}_{j}$ as defined in equation~\eqref{eq:JML_mat_blocking}. The last symmetric term
 $\mathbf{U}(\mathbf{a}) ^{T} \mathbf{L} (\mathbf{a}) + \mathbf{L}(\mathbf{a})^{T} \mathbf{U}(\mathbf{a}) $ in the above expression of $\nabla^2 \psi( \mathbf{a} )$ can also be efficiently computed in two stages. First, by computing in $n$ parallel steps, the term
\begin{equation*}
\mathbf{L}(\mathbf{a})^{T} \mathbf{U}(\mathbf{a})  =   \mathbf{V}(\mathbf{a}) ^{T} \mathbf{F}(\mathbf{a})^{-} \mathbf{U}(\mathbf{a})  = \sum_{j=1}^n \mathbf{V}_{j}^{T} \mathbf{F}_{j}(\mathbf{a})^{-} \mathbf{U}_{j}  \ ,
\end{equation*}
where $\mathbf{V}(\mathbf{a})$ is defined  in equation~\eqref{eq:V_mat} and the matrices $\mathbf{U}(\mathbf{a})$ and  $\mathbf{V}(\mathbf{a})$  have also been partitioned into $n$ blocks $\mathbf{V}_{j}$ and $\mathbf{U}_{j}$ as defined in equations~\eqref{eq:blk_U_mat} and~\eqref{eq:blk_V_mat}. Note that, in this last equation, we have used again a symmetric generalized inverse of $\mathbf{F}(\mathbf{a})$ instead of the Moore-Penrose inverse $\mathbf{F}(\mathbf{a})^{+}$ as we assume here that $\mathbf{F}(\mathbf{a})$ is of full column rank, e.g., that ${r}_{\mathbf{F}( \mathbf{a} )} = \emph{rank}( \mathbf{F}(  \mathbf{a} ) ) =k.p$. In a second stage, we just transpose this squared matrix and sum this squared matrix and its transpose to get the third symmetric term in the above expression of $\nabla^2 \psi( \mathbf{a} )$.

Using the matrix $\mathbf{H} $ or its two-term approximation $\bar{\mathbf{H}}$, a basic variable projection Newton or quasi-Newton algorithm for the WLRA problem has to solve at each iteration, respectively, the constrained linear systems (see equations~\eqref{eq:newton_sys2} and~\eqref{eq:approx_newton_sys2})
\begin{equation*}
\big( \mathbf{H} +  \mathbf{N} \mathbf{N}^{T} \big)  d\mathbf{a}_{n}  = - \nabla \psi( \mathbf{a} ) = \mathbf{M}(\mathbf{a})^{T} \mathbf{r}(\mathbf{a})
\end{equation*}
and
\begin{equation*}
\big(   \bar{\mathbf{H}} +  \mathbf{N} \mathbf{N}^{T} \big)  d\mathbf{a}_{n}  = - \nabla \psi( \mathbf{a} ) = \mathbf{M}(\mathbf{a})^{T} \mathbf{r}(\mathbf{a})  \ ,
\end{equation*}
where the columns of the matrix $\mathbf{N} = \mathbf{K}_{(p,k)} ( \mathbf{I}_{k} \otimes \mathbf{A} )$ form a (orthonormal) basis of the null space of $\mathit{J}( \mathbf{r}(\mathbf{a}) )$. As already discussed in the previous subsections, such a basis can be easily computed with the help of Corollary~\ref{corol5.6:box} under the hypothesis that  $\emph{rank}(\mathbf{A}) = k$ and $\emph{rank}( \mathit{J}( \mathbf{r}(\mathbf{a}) ) ) = \emph{rank}( \mathbf{M}(\mathbf{a}) ) = k.(p-k)$. Furthermore, the cross-product matrix $\mathbf{N} \mathbf{N}^{T}$ can be evaluated efficiently with equation~\eqref{eq:NNt_mat} and the matrix-product $\mathbf{M}(\mathbf{a})^{T} \mathbf{r}(\mathbf{a}) = - \nabla \psi( \mathbf{a} )$ can also be computed in $n$ parallel steps as the matrices $\mathbf{H}$ and $\bar{\mathbf{H}}$ (see Subsection~\ref{vp_gn_alg:box} for details).

As discussed at the end of Subsection~\ref{hess:box}, adding the symmetric matrix $\mathbf{N} \mathbf{N}^{T}$ to $\mathbf{H}$ and $\bar{\mathbf{H}}$ will guarantee that the minimum 2-norm solutions of the consistent linear systems
\begin{equation*}
\mathbf{H}  d\mathbf{a}_{n}  =  \mathbf{M}(\mathbf{a})^{T} \mathbf{r}(\mathbf{a})
\quad \text{and} \quad
 \bar{\mathbf{H}}   d\mathbf{a}_{n}  = \mathbf{M}(\mathbf{a})^{T} \mathbf{r}(\mathbf{a})
\end{equation*}
are computed if $\emph{rank}(\mathbf{A}) = k$ and $\emph{rank}( \mathit{J}( \mathbf{r}(\mathbf{a}) ) ) = \emph{rank}( \mathbf{M}(\mathbf{a}) ) = k.(p-k)$ and, thus, will overcome the systematic singularity and ill-conditioning of  $\bar{\mathbf{H}}$ or those of $\mathbf{H}$ at the stationary points of $\psi(.)$ in most cases without using any pre-conditioner for the Hessian matrix or its approximation as suggested in~\cite{BA2015}. Note that these two constrained linear systems are based, respectively, on the quadratic approximation models
\begin{equation*}
 \psi( \mathbf{a} + d\mathbf{a} )   \approx N^{\mathbf{N}} ( d\mathbf{a} ) = \psi( \mathbf{a} ) + d\mathbf{a}^{T}  \nabla \psi( \mathbf{a} ) +  \frac{1}{2} d\mathbf{a}^{T}  \big(  \mathbf{H}  +  \mathbf{N} \mathbf{N}^{T}  \big) d\mathbf{a}
\end{equation*}
and
\begin{equation*}
 \psi( \mathbf{a} + d\mathbf{a} )   \approx N^{\mathbf{N}}( d\mathbf{a} ) = \psi( \mathbf{a} ) + d\mathbf{a}^{T} \nabla \psi( \mathbf{a} )  +  \frac{1}{2} d\mathbf{a}^{T}  \big(  \bar{\mathbf{H}} +  \mathbf{N} \mathbf{N}^{T}  \big) d\mathbf{a}   \ .
\end{equation*}
These quadratic functions are more accurate then  the corresponding Golub-Pereyra and Kaufman quadratic models
\begin{equation*}
G^{\mathbf{N}} ( d\mathbf{a} ) = \psi( \mathbf{a} ) + d\mathbf{a}^{T}  \nabla \psi( \mathbf{a} ) +  \frac{1}{2} d\mathbf{a}^{T}  \big(  \mathbf{M}( \mathbf{a}  )^{T} \mathbf{M}( \mathbf{a} ) + \mathbf{L}( \mathbf{a} )^{T} \mathbf{L}( \mathbf{a} )   +  \mathbf{N} \mathbf{N}^{T}  \big) d\mathbf{a}
\end{equation*}
and
\begin{equation*}
G^{\mathbf{N}}( d\mathbf{a} ) = \psi( \mathbf{a} ) + d\mathbf{a}^{T} \nabla \psi( \mathbf{a} )  +  \frac{1}{2} d\mathbf{a}^{T}  \big(  \mathbf{M}( \mathbf{a}  )^{T} \mathbf{M}( \mathbf{a} ) +  \mathbf{N} \mathbf{N}^{T}  \big) d\mathbf{a}   \ ,
\end{equation*}
 which are used in the variable projection Gauss-Newton methods, as information on second-order derivatives of $\psi(.)$ is included in $N^{\mathbf{N}} (.)$, but not in $G^{\mathbf{N}} (.)$.

However, if the addition of the term $\mathbf{N} \mathbf{N}^{T}$ to the Hessian matrix $\mathbf{H}$, or to its two-term approximation $\bar{\mathbf{H}}$, is useful to ensure that these matrices stay nonsingular and well-conditioned at each iteration, it is not sufficient to guarantee that these two matrices stay positive definite in all the steps and , thus, that $d\mathbf{a}_{n}$ is always in a descent direction for $\psi(.)$. As explained in Subsection~\ref{opt:box}, this property can be achieved by adding an (another) damping term $\lambda \mathbf{I}_{k.p}$ (where  $\lambda>0$) in the above Newton quadratic models
\begin{equation*}
N_{\lambda}^{\mathbf{N}} ( d\mathbf{a} ) = \psi( \mathbf{a} ) + d\mathbf{a}^{T}  \nabla \psi( \mathbf{a} ) +  \frac{1}{2} d\mathbf{a}^{T}  \big(  \mathbf{H} +  \mathbf{N} \mathbf{N}^{T} +  \lambda \mathbf{I}_{k.p}  \big) d\mathbf{a}
\end{equation*}
or
\begin{equation*}
N_{\lambda}^{\mathbf{N}}( d\mathbf{a} ) = \psi( \mathbf{a} ) + d\mathbf{a}^{T} \nabla \psi( \mathbf{a} )  +  \frac{1}{2} d\mathbf{a}^{T}  \big(   \bar{\mathbf{H}} +  \mathbf{N} \mathbf{N}^{T} +  \lambda \mathbf{I}_{k.p} \big) d\mathbf{a}   \ ,
\end{equation*}
which lead, respectively, to the damped and constrained linear systems
\begin{equation*}
\big( \mathbf{H} +  \mathbf{N} \mathbf{N}^{T}  +  \lambda \mathbf{I}_{k.p} \big)  d\mathbf{a}_{n}  = \mathbf{M}(\mathbf{a})^{T} \mathbf{r}(\mathbf{a})
\end{equation*}
and
\begin{equation*}
\big(   \bar{\mathbf{H}} +  \mathbf{N} \mathbf{N}^{T}  +  \lambda \mathbf{I}_{k.p} \big)  d\mathbf{a}_{n}  = \mathbf{M}(\mathbf{a})^{T} \mathbf{r}(\mathbf{a})  \ ,
\end{equation*}
for the computation of the Newton or quasi-Newton correction step at each iteration.
\\
\begin{remark6.4} \label{remark6.4:box}
Alternatively, as already discussed at the end of Subsection~\ref{hess:box}, the correction vectors in the variable projection Newton or quasi-Newton algorithms can be computed in a two-step procedure as first suggested by Chen~\cite{C2008b}. First, by solving the reduced  damped linear system
\begin{equation*}
(\mathbf{\bar{O}}^{\bot} )^{T} \big( \mathbf{H} +  \lambda \mathbf{I}_{k.p} \big) \mathbf{\bar{O}}^{\bot} d\mathbf{\bar{a}}_{n} =   \big( \mathbf{M}(\mathbf{a}) \mathbf{\bar{O}}^{\bot}   \big)^{T} \mathbf{r}(\mathbf{a})
\end{equation*}
or
\begin{equation*}
(\mathbf{\bar{O}}^{\bot} )^{T} \big( \bar{\mathbf{H}} +  \lambda \mathbf{I}_{k.p} \big) \mathbf{\bar{O}}^{\bot} d\mathbf{\bar{a}}_{n} =   \big( \mathbf{M}(\mathbf{a}) \mathbf{\bar{O}}^{\bot}   \big)^{T} \mathbf{r}(\mathbf{a})   \ ,
\end{equation*}
for $d\mathbf{\bar{a}}_{n}  \in \mathbb{R}^{(p-k).k}$ and where $\mathbf{\bar{O}}^{\bot} =  \mathbf{K}_{(p,k)} ( \mathbf{I}_{k}  \otimes  \mathbf{O}^{\bot}) \in \mathbb{O}^{k.p \times ( p - k ).k }$ and $\mathbf{O}^{\bot} \in \mathbb{O}^{p \times ( p - k ) }$ is an orthonormal matrix whose columns form a basis of $ \emph{ran}( \mathbf{A} )^{\bot}$. Next, $d\mathbf{a}_{n}$ is computed by
\begin{equation*}
d\mathbf{A}_{n} = \mathbf{O}^{\bot} d\mathbf{\bar{A}}_{n}  \ ,
\end{equation*}
in a second step. $\blacksquare$
\\
\end{remark6.4}
In other words, we can develop variable projection $Levenberg$-$Marquardt$-type Newton and quasi-Newton algorithms with a wider basin of convergence by using the same strategies as used in the Levenberg-Marquardt algorithms described in Subsection~\ref{vp_lm_alg:box}.

Note that we can also define and use a gain factor $\rho$ in the context of these  variable projection $Levenberg$-$Marquardt$-type Newton and quasi-Newton methods, e.g.,
\begin{equation*}
\rho = \frac{ \psi( \mathbf{a}_{i}) - \psi( \mathbf{a}_{i} + d\mathbf{a} )}{ N(\mathbf{0}^{k.p}) -  N( d\mathbf{a} ) }  \ ,
\end{equation*}
where $N( . )$ is the more accurate quadratic model used by the (quasi-)Newton iterations and defined by
\begin{equation*}
N( d\mathbf{a} ) = \psi( \mathbf{a} ) + d\mathbf{a}^{T} \nabla \psi( \mathbf{a} )  +  \frac{1}{2} d\mathbf{a}^{T}  \mathbf{H} d\mathbf{a}  
\end{equation*}
or
\begin{equation*}
N( d\mathbf{a} ) = \psi( \mathbf{a} ) + d\mathbf{a}^{T} \nabla \psi( \mathbf{a} )  +  \frac{1}{2} d\mathbf{a}^{T}  \bar{\mathbf{H}}  d\mathbf{a}   \ ,
\end{equation*}
if we use a quasi-Newton algorithm. Furthermore, in both cases, a direct computation shows that we can compute cheaply the difference $N(\mathbf{0}^{k.p}) -  N( d\mathbf{a}_{n} )$  at each iteration of the (quasi-)Newton algorithms as
\begin{equation*}
N(\mathbf{0}^{k.p}) -  N( d\mathbf{a}_{n} ) = \frac{1}{2}   \big(  \Vert \mathbf{N}_{i}^{T}  d\mathbf{a}_{n} \Vert^{2}_{2} + \lambda \Vert d\mathbf{a}_{n} \Vert^{2}_{2} - d\mathbf{a}_{n}^{T} \nabla \psi( \mathbf{a}_{i} )  \big) \ ,
\end{equation*}
similarly as for the Levenberg-Marquardt algorithms~\eqref{lm_alg3:box} developed in the last subsection. Furthermore, $N(\mathbf{0}^{k.p}) -  N( d\mathbf{a}_{n} )$ will be guaranteed to be positive, if $d\mathbf{a}_{n}$ is in a descent direction for $\psi( . )$ and $\Vert \nabla \psi( \mathbf{a}_{i} ) \Vert_{2} \ne 0$.

Using these considerations and inspired by the Levenberg-Marquardt algorithms~\eqref{lm_alg1:box},~\eqref{lm_alg2:box} and~\eqref{lm_alg3:box}, an outline of three different versions of the variable projection (quasi-)Newton algorithms is detailed below. The definitions and variables used in these Newton and quasi-Newton algorithms have the same meaning as in the previous Gauss-Newton and Levenberg-Marquardt algorithms.

Note that, in all these algorithms, we compute the Newton correction step $d\mathbf{a}_{n}$ with the exact Hessian matrix $\mathbf{H}$. However, at the user convenience, for example to reduce the computing time per iteration and the memory footprint in all the algorithms, we can eliminate the computation of the third symmetric term of $\nabla^2 \psi( \mathbf{a} )$ and use the approximate Hessian matrix $\bar{\mathbf{H}}$ instead to compute the Newton correction step $d\mathbf{a}_{n}$ at each iteration. Thus, this simple modification defines the quasi-Newton variant for all the algorithms.
\\
\begin{n_alg1} \label{n_alg1:box}
\end{n_alg1}
Choose starting matrix $\mathbf{A}_{1} \in \mathbb{R}^{p \times k}$ , $\varepsilon_{1}, \varepsilon_{2}, \varepsilon_{3} , \beta \in \mathbb{R}_{+*}$ and $i_{max}, j_{max} \in \mathbb{N}_{*}$, appropriately

\textbf{For} $i=1, 2, \ldots$ \textbf{until convergence do}
\begin{enumerate}
\item[\textbf{(0)}] Optionally, compute a QRCP of $\mathbf{A}_{i}$ (see equation~\eqref{eq:qrcp}) to determine $k_{i} = \emph{rank}( \mathbf{A}_{i} )$ and an orthonormal basis of $\emph{ran}( \mathbf{A}_{i} )$:
\begin{itemize}
\item[] $\mathbf{Q}_{i} \mathbf{A}_{i} \mathbf{P}_{i} =  \begin{bmatrix} \mathbf{R}_{i}    &  \mathbf{S}_{i}   \\  \mathbf{0}^{(p-k_{i}) \times k_{i} }  & \mathbf{0}^{(p-k_{i}) \times (k-k_{i}) }  \end{bmatrix}$ ,
\end{itemize}
where $\mathbf{Q}_{i}$ is an $p \times p$ orthogonal matrix, $\mathbf{P}_{i}$ is an $k \times k$ permutation matrix, $\mathbf{R}_{i}$ is an $k_{i} \times k_{i}$ nonsingular upper triangular matrix (with diagonal elements of decreasing absolute magnitude) and $\mathbf{S}_{i}$ an $k_{i} \times (k-k_{i})$ full matrix, which is vacuous if $k_{i} = k$.

In all cases, compute an $p \times k$  matrix $\mathbf{O}_{i}$  with orthonormal columns as the first $k$ columns of $\mathbf{Q}_{i}$ (i.e., such that $\emph{ran}( \mathbf{A}_{i} ) \subset \emph{ran}( \mathbf{O}_{i} )$ if $k_{i} < k$ and  $\emph{ran}( \mathbf{A}_{i} ) = \emph{ran}( \mathbf{O}_{i} )$ if $k_{i} = k$) and set
\begin{itemize}
\item[] $\mathbf{A}_{i} = \mathbf{O}_{i}$ .
\end{itemize}
This optional orthogonalization step is a safe-guard as the condition $k_{i} = k$ is a necessary condition for the differentiability of $\psi(.)$ at a point $\mathbf{A}_{i}$ and also to limit the occurrence of overflows and underflows in the next steps by enforcing that the matrix variable $\mathbf{A}_{i} \in \mathbb{O}^{p \times k}$.
\item[\textbf{(1)}] Determine (implicitly) the block diagonal matrix
\begin{itemize}
\item[] $\mathbf{F}(\mathbf{a}_{i}) = \emph{diag}\big( \emph{vec}( \sqrt{\mathbf{W}} ) \big)  \big(  \mathbf{I}_n  \otimes \mathbf{A}_{i}  \big)$ ,
\end{itemize}
where $\mathbf{a}_{i} =  \emph{vec}( \mathbf{A}_{i}^{T} )$.
\item[\textbf{(2)}] Compute (implicitly) a QRCP of $\mathbf{F}(\mathbf{a}_{i})$ to determine  ${r}_{\mathbf{F}( \mathbf{a}_{i} )}  = \emph{rank}( \mathbf{F}(  \mathbf{a}_{i} ) )$, $\mathbf{P}^{\bot}_{\mathbf{F}(  \mathbf{a}_{i} ) }$ and  $\mathbf{F}(\mathbf{a}_{i})^-$ (see equations~\eqref{eq:ginv_proj_ortho} and \eqref{eq:sginv_qrcp}). Note that  $\mathbf{P}^{\bot}_{\mathbf{F}(  \mathbf{a}_{i} ) }$ and $\mathbf{F}(\mathbf{a}_{i})^-$ are also block diagonal matrices and that the Newton and quasi-Newton algorithms assume that ${r}_{\mathbf{F}( \mathbf{a}_{i} )} = k.p$; see the derivation of the Hessian matrix  $\nabla^2 \psi( \mathbf{a}_{i} )$ in Subsection~\ref{hess:box} for more details. Note that the validity of this hypothesis can be checked here in output of the QRCP of $\mathbf{F}(\mathbf{a}_{i})$.

\item[\textbf{(3)}] Solve the block diagonal linear least-squares problem
\begin{itemize}
\item[] $\mathbf{b}_{i} = \text{Arg}\min_{\mathbf{b}\in\mathbb{R}^{k.n}} \,   \Vert \mathbf{x} - \mathbf{F}(\mathbf{a}_{i})\mathbf{b} \Vert^{2}_{2}$ ,
\end{itemize}
e.g., compute $\mathbf{b}_{i} = \mathbf{F}(\mathbf{a}_{i})^{-} \mathbf{x}$.
\item[\textbf{(4)}] Determine:
\begin{itemize}
\item[] $\mathbf{r}(\mathbf{a}_{i}) = \mathbf{P}^{\bot}_{\mathbf{F}(  \mathbf{a}_{i} ) } \mathbf{x}$ $\lbrace \text{current residual vector} \rbrace$
\item[] $\psi(\mathbf{a}_{i} ) = \frac{1}{2} \Vert \mathbf{r}(\mathbf{a}_{i}) \Vert^{2}_{2}$ $\lbrace \text{current value of the cost function }  \rbrace$
\item[] $\nabla \psi( \mathbf{a}_{i} ) =  \mathbf{G}(\mathbf{b}_{i})^{T} \mathbf{G}(\mathbf{b}_{i})\mathbf{a}_{i} - \mathbf{G}(\mathbf{b}_{i})^{T} \mathbf{z}$ $\lbrace \text{see Theorems~\ref{theo4.3:box} and~\ref{theo5.7:box}} \rbrace$
\item[] $\lambda_{i} = \beta \Vert  \nabla \psi( \mathbf{a}_{i} ) \Vert^{2}_{2}$ $\lbrace \text{set ridge parameter proportional to the squared 2-norm of the gradient} \rbrace$
\end{itemize}
Note that the steps \textbf{(1)} to \textbf{(4)} above can be very easily parallelized using the block diagonal structure of $\mathbf{F}(\mathbf{a}_{i})$.
\item[\textbf{(5)}]  Check for convergence. Relevant convergence criteria in the algorithms are of the form:
\begin{itemize}
\item $\Vert \nabla \psi(\mathbf{a}_{i} ) \Vert_{2} \le \varepsilon_{1}$
\item $\Vert \mathbf{a}_{i} - \mathbf{a}_{i-1} \Vert_{2} \le \varepsilon_{2} ( \varepsilon_{2} + \Vert \mathbf{a}_{i} \Vert_{2})$  $\lbrace$if $ i \ne 1\rbrace$.

If step \textbf{(0)} is used, this convergence condition can be simplified as:

$\Vert \mathbf{a}_{i} - \mathbf{a}_{i-1} \Vert_{2} \le \varepsilon_{2}  \Vert \mathbf{a}_{i} \Vert_{2} = \varepsilon_{2} \sqrt{k}$
\item $\vert \psi(\mathbf{a}_{i-1} ) - \psi(\mathbf{a}_{i} ) \vert \le \varepsilon_{3} ( \varepsilon_{3} +  \psi(\mathbf{a}_{i} ) )$  $\lbrace$if $ i \ne 1\rbrace$
\item $ i \ge i_{max}$ $\lbrace \text{e.g., give up if the number of iterations is too large} \rbrace$
\end{itemize}
where $\varepsilon_{1}, \varepsilon_{2}, \varepsilon_{3}$ and $i_{max}$ are constants chosen by the user.

\textbf{Exit if convergence}. \textbf{Otherwise, go to} step \textbf{(6)}
\item[\textbf{(6)}] Compute the Newton correction vector $d\mathbf{a}_{n}$.
\begin{enumerate}
\item[\textbf{(6.1)}] To this end, first compute Hessian matrix or its two-term approximation:
\begin{itemize}
\item[] $\mathbf{H}_{i}^{1} = \mathbf{M}(\mathbf{a}_{i})^{T} \mathbf{M}(\mathbf{a}_{i}) + \mathbf{N}_{i} \mathbf{N}_{i}^{T}  + \lambda_{i} \mathbf{I}_{k.p} $
\item[] $\mathbf{H}_{i}^{2} = \mathbf{L}(\mathbf{a}_{i})^{T} \mathbf{L}(\mathbf{a}_{i})$
\item[] $\mathbf{H}_{i}^{3} = \mathbf{U}(\mathbf{a}_{i})^{T} \mathbf{L} (\mathbf{a}_{i}) + \mathbf{L}(\mathbf{a}_{i})^{T} \mathbf{U}(\mathbf{a}_{i})$ $\lbrace$only if a full Newton step is wanted$\rbrace$
\item[] $\mathbf{H}_{i} = 
    \begin{cases}
        \mathbf{H}_{i}^{1} - \mathbf{H}_{i}^{2} + \mathbf{H}_{i}^{3}    & \lbrace \text{for a full Newton step}\rbrace \\
        \mathbf{H}_{i}^{1} - \mathbf{H}_{i}^{2}    & \lbrace \text{for a quasi-Newton step} \rbrace
     \end{cases} 
$
\end{itemize}
where the columns of $\mathbf{N}_{i} = \mathbf{K}_{(p,k)} ( \mathbf{I}_{k} \otimes \mathbf{A}_{i} )$ are a (orthonormal) basis of  $\emph{null}\big( \mathbf{M}(\mathbf{a}_{i}) + \mathbf{L}(\mathbf{a}_{i}) \big) = \emph{null}( \mathbf{M}(\mathbf{a}_{i}) )$, see Corollary~\ref{corol5.6:box}.
\item[\textbf{(6.2)}] \textbf{If} $\mathbf{H}_{i}$ is positive definite \textbf{then} $\lbrace \text{use Cholesky factorization} \rbrace$
\begin{description}
\item[Newton step:] get Newton or quasi-Newton step as the solution of the linear system
\begin{equation*}
\mathbf{H}_{i}  d\mathbf{a}_{n} =  - \nabla \psi( \mathbf{a}_{i} ) = \mathbf{M}(\mathbf{a}_{i})^{T} \mathbf{r}(\mathbf{a}_{i})
\end{equation*}
\end{description}
\item[\textbf{(6.3)}] \textbf{Else}
\begin{description}
\item[Gauss-Newton step:] get (Gauss-)Newton step as the solution of the linear system
\begin{equation*}
\mathbf{H}_{i}^{1}  d\mathbf{a}_{n} =  - \nabla \psi( \mathbf{a}_{i} ) =  \mathbf{M}(\mathbf{a}_{i})^{T} \mathbf{r}(\mathbf{a}_{i})
\end{equation*}
\end{description}
\end{enumerate}
\item[\textbf{(7)}] Increment $\mathbf{a}_{i} = \emph{vec}( \mathbf{A}_{i}^{T} )$, e.g., compute $\mathbf{a}_{i+1} = \emph{vec}( \mathbf{A}_{i+1}^{T} )$ such that $\psi( \mathbf{a}_{i+1} ) < \psi( \mathbf{a}_{i} )$ in order to obtain global convergence.
\begin{enumerate}
\item[\textbf{(7.1)}]  To this end, first compute:
\begin{itemize}
\item[] $\mathbf{a}_{i+1} = \mathbf{a}_{i} +  d\mathbf{a}_{n}$
\item[] $\psi(\mathbf{a}_{i+1} ) = \frac{1}{2} \Vert \mathbf{r}(\mathbf{a}_{i+1}) \Vert^{2}_{2} = \frac{1}{2} \Vert \mathbf{P}^{\bot}_{\mathbf{F}(  \mathbf{a}_{i+1} ) } \mathbf{x} \Vert^{2}_{2}$ ,
\end{itemize}
using (implicitly) a QRCP of the block diagonal matrix $\mathbf{F}(  \mathbf{a}_{i+1} )$.
\item[\textbf{(7.2)}] \textbf{If}  $\psi( \mathbf{a}_{i+1} ) > \psi( \mathbf{a}_{i} )$ \textbf{then} recompute $\mathbf{a}_{i+1}$ by one of the following methods:
\begin{description}
\item[Gauss-Seidel:]
$\mathbf{a}_{i+1} = \mathbf{a}_{i} +  d\mathbf{a}_{gs-gn}$ where $d\mathbf{a}_{gs-gn}$ is a Gauss-Seidel step~\cite{RW1980}
\begin{align*}
d\mathbf{a}_{gs-gn} & =   \big( \mathbf{K}_{(n,p)}  \mathbf{G}( \mathbf{b}_{i}  ) \big)^{+} \mathbf{r}( \mathbf{a}_{i} )   \\
                                 & = 
    \begin{cases}
        \text{Arg}\min_{d\mathbf{a} \in \mathbb{R}^{p.k}} \,   \Vert d\mathbf{a} \Vert^{2}_{2} \\
        \text{s.t. }  \text{Arg}\min_{d\mathbf{a} \in \mathbb{R}^{p.k}} \,   \Vert \mathbf{r}( \mathbf{a}_{i} ) - \mathbf{K}_{(n,p)}  \mathbf{G}( \mathbf{b}_{i}  ) d\mathbf{a}   \Vert^{2}_{2}
    \end{cases}   \\
\end{align*}
\item[Block alternating least-squares:]
\begin{align*}
  \mathbf{a}_{i+1} & =    \mathbf{G}( \mathbf{b}_{i}  )^{+} \mathbf{z}   \\
                             & = 
    \begin{cases}
        \text{Arg}\min_{\mathbf{a} \in \mathbb{R}^{p.k}} \,   \Vert \mathbf{a} \Vert^{2}_{2} \\
        \text{s.t. }  \text{Arg}\min_{\mathbf{a} \in \mathbb{R}^{p.k}} \,   \Vert \mathbf{z} -  \mathbf{G}( \mathbf{b}_{i}  ) \mathbf{a}   \Vert^{2}_{2}
    \end{cases}   \\
\end{align*}
\item[Line search:]
\begin{equation*}
\mathbf{a}_{i+1} = \mathbf{a}_{i} + \alpha_{i}  d\mathbf{a}_{n}
\end{equation*}
where $\alpha_i < 1$ is determined by a line search to make the algorithm a descent method (i.e., such that $\psi( \mathbf{a}_{i+1} ) <  \psi( \mathbf{a}_{i} ) $). This is always possible as the correction vector $d\mathbf{a}_{n}$ is in a descent direction for $\psi(.)$ if $\Vert \nabla \psi(\mathbf{a}_{i} ) \Vert_{2} \ne 0$.

A simple strategy is to first shorten the correction step to half the Newton length (or Gauss-Newton length if $\mathbf{H}_{i}$ is not positive definite), compute the new trial value for $\psi( \mathbf{a}_{i+1} )$ and, if it is still worse, continue to reduce the step until we get a step short enough such that $\psi( \mathbf{a}_{i+1} ) <  \psi( \mathbf{a}_{i} )$. The following loop incorporates this simple step-shortening algorithm:

\textbf{For} $j=1, 2, \ldots$ \textbf{while}$\big( \psi( \mathbf{a}_{i+1} ) > \psi( \mathbf{a}_{i} ) \big)$
\begin{enumerate}
\item[]  $d\mathbf{a}_{n}  = \frac{1}{2}  d\mathbf{a}_{n}$
\item[]  $\mathbf{a}_{i+1} = \mathbf{a}_{i} +  d\mathbf{a}_{n}$
\item[]  $\psi(\mathbf{a}_{i+1} ) =  \frac{1}{2} \Vert \mathbf{P}^{\bot}_{\mathbf{F}(  \mathbf{a}_{i+1} ) } \mathbf{x} \Vert^{2}_{2}$ $\lbrace$using a QRCP of the matrix $\mathbf{F}(  \mathbf{a}_{i+1} )\rbrace$
\item[]  \textbf{If} $j > j_{max}$ \textbf{exit} $\lbrace \text{give up if the number of iterations is too large} \rbrace$
\end{enumerate}
\textbf{End do}

\end{description}
\end{enumerate}
\end{enumerate}
\textbf{End do}
\\
\begin{n_alg2} \label{n_alg2:box}
\end{n_alg2}
Choose starting matrix $\mathbf{A}_{1} \in \mathbb{R}^{p \times k}$ , $\varepsilon_{1}, \varepsilon_{2}, \varepsilon_{3} , \lambda \in \mathbb{R}_{+*}$ and $i_{max}, j_{max} \in \mathbb{N}_{*}$, appropriately

\textbf{For} $i=1, 2, \ldots$ \textbf{until convergence do}
\begin{enumerate}
\item[\textbf{(0)}] Optionally, compute a QRCP of $\mathbf{A}_{i}$ (see equation~\eqref{eq:qrcp}) to determine $k_{i} = \emph{rank}( \mathbf{A}_{i} )$ and an orthonormal basis of $\emph{ran}( \mathbf{A}_{i} )$:
\begin{itemize}
\item[] $\mathbf{Q}_{i} \mathbf{A}_{i} \mathbf{P}_{i} =  \begin{bmatrix} \mathbf{R}_{i}    &  \mathbf{S}_{i}   \\  \mathbf{0}^{(p-k_{i}) \times k_{i} }  & \mathbf{0}^{(p-k_{i}) \times (k-k_{i}) }  \end{bmatrix}$ ,
\end{itemize}
where $\mathbf{Q}_{i}$ is an $p \times p$ orthogonal matrix, $\mathbf{P}_{i}$ is an $k \times k$ permutation matrix, $\mathbf{R}_{i}$ is an $k_{i} \times k_{i}$ nonsingular upper triangular matrix (with diagonal elements of decreasing absolute magnitude) and $\mathbf{S}_{i}$ an $k_{i} \times (k-k_{i})$ full matrix, which is vacuous if $k_{i} = k$.

In all cases, compute an $p \times k$  matrix $\mathbf{O}_{i}$  with orthonormal columns as the first $k$ columns of $\mathbf{Q}_{i}$ (i.e., such that $\emph{ran}( \mathbf{A}_{i} ) \subset \emph{ran}( \mathbf{O}_{i} )$ if $k_{i} < k$ and  $\emph{ran}( \mathbf{A}_{i} ) = \emph{ran}( \mathbf{O}_{i} )$ if $k_{i} = k$) and set
\begin{itemize}
\item[] $\mathbf{A}_{i} = \mathbf{O}_{i}$ .
\end{itemize}
This optional orthogonalization step is a safe-guard as the condition $k_{i} = k$ is a necessary condition for the differentiability of $\psi(.)$ at a point $\mathbf{A}_{i}$ and also to limit the occurrence of overflows and underflows in the next steps by enforcing that the matrix variable $\mathbf{A}_{i} \in \mathbb{O}^{p \times k}$.
\item[\textbf{(1)}] Determine (implicitly) the block diagonal matrix
\begin{itemize}
\item[] $\mathbf{F}(\mathbf{a}_{i}) = \emph{diag}\big( \emph{vec}( \sqrt{\mathbf{W}} ) \big)  \big(  \mathbf{I}_n  \otimes \mathbf{A}_{i}  \big)$ ,
\end{itemize}
where $\mathbf{a}_{i} =  \emph{vec}( \mathbf{A}_{i}^{T} )$
\item[\textbf{(2)}] Compute (implicitly) a QRCP of $\mathbf{F}(\mathbf{a}_{i})$ to determine  ${r}_{\mathbf{F}( \mathbf{a}_{i} )}  = \emph{rank}( \mathbf{F}(  \mathbf{a}_{i} ) )$, $\mathbf{P}^{\bot}_{\mathbf{F}(  \mathbf{a}_{i} ) }$ and  $\mathbf{F}(\mathbf{a}_{i})^-$ (see equations~\eqref{eq:ginv_proj_ortho} and \eqref{eq:sginv_qrcp}). Note that  $\mathbf{P}^{\bot}_{\mathbf{F}(  \mathbf{a}_{i} ) }$ and $\mathbf{F}(\mathbf{a}_{i})^-$ are also block diagonal matrices and that the Newton and quasi-Newton algorithms assume that ${r}_{\mathbf{F}( \mathbf{a}_{i} )}  =  k.p$; see the derivation of the Hessian matrix  $\nabla^2 \psi( \mathbf{a}_{i} )$ in Subsection~\ref{hess:box} for more details. Note that, optionally,  the validity of this hypothesis can be checked here in output of the QRCP of $\mathbf{F}(\mathbf{a}_{i})$.

\item[\textbf{(3)}] Solve the block diagonal linear least-squares problem
\begin{itemize}
\item[] $\mathbf{b}_{i} = \text{Arg}\min_{\mathbf{b}\in\mathbb{R}^{k.n}} \,   \Vert \mathbf{x} - \mathbf{F}(\mathbf{a}_{i})\mathbf{b} \Vert^{2}_{2}$ ,
\end{itemize}
e.g., compute $\mathbf{b}_{i} = \mathbf{F}(\mathbf{a}_{i})^{-} \mathbf{x}$.
\item[\textbf{(4)}] Determine and set:
\begin{itemize}
\item[] $\mathbf{r}(\mathbf{a}_{i}) = \mathbf{P}^{\bot}_{\mathbf{F}(  \mathbf{a}_{i} ) } \mathbf{x}$ $\lbrace \text{current residual vector} \rbrace$
\item[] $\psi(\mathbf{a}_{i} ) = \frac{1}{2} \Vert \mathbf{r}(\mathbf{a}_{i}) \Vert^{2}_{2}$ $\lbrace \text{current value of the cost function }  \rbrace$
\item[] $\nabla \psi( \mathbf{a}_{i} ) =  \mathbf{G}(\mathbf{b}_{i})^{T} \mathbf{G}(\mathbf{b}_{i})\mathbf{a}_{i} - \mathbf{G}(\mathbf{b}_{i})^{T} \mathbf{z}$ $\lbrace \text{see Theorems~\ref{theo4.3:box} and~\ref{theo5.7:box}} \rbrace$
\item[] $j = 0$ $\lbrace \text{initialize counter for the ridge scaling subiterations} \rbrace$
\end{itemize}
Note that the steps  \textbf{(1)} to \textbf{(4)}  above can be very easily parallelized using the block diagonal structure of $\mathbf{F}(\mathbf{a}_{i})$.
\item[\textbf{(5)}]  Check for convergence. Relevant convergence criteria in the algorithms are of the form:
\begin{itemize}
\item $\Vert \nabla \psi(\mathbf{a}_{i} ) \Vert_{2} \le \varepsilon_{1}$
\item $\Vert \mathbf{a}_{i} - \mathbf{a}_{i-1} \Vert_{2} \le \varepsilon_{2} ( \varepsilon_{2} + \Vert \mathbf{a}_{i} \Vert_{2})$  $\lbrace$if $ i \ne 1\rbrace$

If step \textbf{(0)} is used, this convergence condition can be simplified as:

$\Vert \mathbf{a}_{i} - \mathbf{a}_{i-1} \Vert_{2} \le \varepsilon_{2}  \Vert \mathbf{a}_{i} \Vert_{2} = \varepsilon_{2} \sqrt{k}$
\item $\vert \psi(\mathbf{a}_{i-1} ) - \psi(\mathbf{a}_{i} ) \vert \le \varepsilon_{3} ( \varepsilon_{3} +  \psi(\mathbf{a}_{i} ) )$  $\lbrace$if $ i \ne 1\rbrace$
\item $ i \ge i_{max}$ $\lbrace \text{e.g., give up if the number of iterations is too large} \rbrace$
\end{itemize}
where $\varepsilon_{1}, \varepsilon_{2}, \varepsilon_{3}$ and $i_{max}$ are constants chosen by the user.

\textbf{Exit if convergence}. \textbf{Otherwise, go to} step \textbf{(6)}
\item[\textbf{(6)}] Compute the Newton correction vector $d\mathbf{a}_{n}$.
\begin{enumerate}
\item[\textbf{(6.1)}] To this end, first compute Hessian matrix or its two-term approximation:
\begin{itemize}
\item[] $\mathbf{H}_{i}^{1} = \mathbf{M}(\mathbf{a}_{i})^{T} \mathbf{M}(\mathbf{a}_{i})$
\item[] $\mathbf{H}_{i}^{2} = \mathbf{L}(\mathbf{a}_{i})^{T} \mathbf{L}(\mathbf{a}_{i})$
\item[] $\mathbf{H}_{i}^{3} = \mathbf{U}(\mathbf{a}_{i})^{T} \mathbf{L} (\mathbf{a}_{i}) + \mathbf{L}(\mathbf{a}_{i})^{T} \mathbf{U}(\mathbf{a}_{i})$ $\lbrace$only if a full Newton step is wanted$\rbrace$
\item[] $\mathbf{H}_{i} =
    \begin{cases}
        \mathbf{H}_{i}^{1} - \mathbf{H}_{i}^{2} + \mathbf{H}_{i}^{3}   & \lbrace \text{for a full Newton step}\rbrace \\
        \mathbf{H}_{i}^{1} - \mathbf{H}_{i}^{2}                                   & \lbrace \text{for a quasi-Newton step} \rbrace
     \end{cases} 
$
\item[] $\mathbf{H}_{i} = \mathbf{H}_{i} + \mathbf{N}_{i} \mathbf{N}_{i}^{T}$
\end{itemize}
where the columns of $\mathbf{N}_{i} = \mathbf{K}_{(p,k)} ( \mathbf{I}_{k} \otimes \mathbf{A}_{i} )$ are a (orthonormal) basis of  $\emph{null}\big( \mathbf{M}(\mathbf{a}_{i}) + \mathbf{L}(\mathbf{a}_{i}) \big) = \emph{null}( \mathbf{M}(\mathbf{a}_{i}) )$, see Corollary~\ref{corol5.6:box}.
\item[\textbf{(6.2)}] Check first diagonal elements of $\mathbf{H}_{i}$:
\begin{itemize}
\item[] $h_{min} = \emph{min}_{ l = 1, \cdots, k.p} \lbrack \mathbf{H}_{i} \rbrack_{ll}$
\item[] \textbf{If} $h_{min}  < 0$ \textbf{then} $\lambda = \lambda - h_{min}$ $\lbrace \text{scale up the ridge parameter} \rbrace$
\end{itemize}
\item[\textbf{(6.3)}] \textbf{Do while} $\mathbf{H}_{i} +  \lambda \mathbf{I}_{k.p}$ is not positive definite $\lbrace \text{use Cholesky factorization} \rbrace$
\begin{itemize}
\item[] $j = j + 1$
\item[] $\lambda = 10.\lambda$ $\lbrace \text{scale up the ridge parameter} \rbrace$
\end{itemize}
\item[\textbf{(6.4)}] 
\begin{description}
\item[Newton step:] get Newton step as the solution of the positive definite linear system
\begin{equation*}
\big(  \mathbf{H}_{i}  +  \lambda \mathbf{I}_{k.p} \big) d\mathbf{a}_{n} =  - \nabla \psi( \mathbf{a}_{i} ) = \mathbf{M}(\mathbf{a}_{i})^{T} \mathbf{r}(\mathbf{a}_{i})
\end{equation*}
\end{description}
\end{enumerate}
\item[\textbf{(7)}] Compute next iterate $\mathbf{a}_{i+1} = \emph{vec}( \mathbf{A}_{i+1}^{T} )$ such that $\psi( \mathbf{a}_{i+1} ) < \psi( \mathbf{a}_{i} )$ in order to obtain global convergence.
\begin{enumerate}
\item[\textbf{(7.1)}]  To this end, first compute
\begin{itemize}
\item[] $\psi( \mathbf{a}_{i} +  d\mathbf{a}_{n} ) = \frac{1}{2} \Vert \mathbf{r}(\mathbf{a}_{i} +  d\mathbf{a}_{n}) \Vert^{2}_{2} = \frac{1}{2} \Vert \mathbf{P}^{\bot}_{\mathbf{F}(  \mathbf{a}_{i} +  d\mathbf{a}_{n} ) } \mathbf{x} \Vert^{2}_{2}$
\end{itemize}
using (implicitly) a QRCP of the block diagonal matrix $\mathbf{F}(  \mathbf{a}_{i} +  d\mathbf{a}_{n} )$.
\item[\textbf{(7.2)}] \textbf{If}  $\psi( \mathbf{a}_{i} +  d\mathbf{a}_{n} ) > \psi( \mathbf{a}_{i} )$ \textbf{then} $\lbrace \text{step rejected} \rbrace$
\begin{itemize}
\item[] $j = j + 1$
\item[] $\lambda = 10.\lambda$ $\lbrace \text{scale up the ridge parameter} \rbrace$
\item[] \textbf{If} $j  \le j_{max}$ \textbf{go to} step  \textbf{(6.3)} $\lbrace \text{recompute } d\mathbf{a}_{n} \text{ with inflated diagonal of  } \mathbf{H}_{i} \rbrace$
\end{itemize}
\item[\textbf{(7.3)}] \textbf{Else} $\lbrace \text{step acceptable} \rbrace$
\begin{itemize}
\item[] \textbf{If} $j  = 0$   \textbf{then} $\lambda = \lambda /10$ $\lbrace \text{scale down the ridge parameter if step is successful} \rbrace$
\end{itemize}
\item[\textbf{(7.4)}]  Increment $\mathbf{a}_{i}$:

$\mathbf{a}_{i+1} = \mathbf{a}_{i} +  d\mathbf{a}_{n} \lbrace \text{compute new iterate} \rbrace $
\end{enumerate}
\end{enumerate}
\textbf{End do}
\\
\\
\begin{n_alg3} \label{n_alg3:box}
\end{n_alg3}
Choose starting matrix $\mathbf{A}_{1} \in \mathbb{R}^{p \times k}$ , $\varepsilon_{1}, \varepsilon_{2}, \varepsilon_{3} , \lambda  \in \mathbb{R}_{+*}$ and $i_{max}, j_{max} \in \mathbb{N}_{*}$, appropriately, and initialize $\nu = 2$

\textbf{For} $i=1, 2, \ldots$ \textbf{until convergence do}
\begin{enumerate}
\item[\textbf{(0)}]  Optionally, compute a QRCP of $\mathbf{A}_{i}$ (see equation~\eqref{eq:qrcp}) to determine $k_{i} = \emph{rank}( \mathbf{A}_{i} )$ and an orthonormal basis of $\emph{ran}( \mathbf{A}_{i} )$:
\begin{itemize}
\item[] $\mathbf{Q}_{i} \mathbf{A}_{i} \mathbf{P}_{i} =  \begin{bmatrix} \mathbf{R}_{i}    &  \mathbf{S}_{i}   \\  \mathbf{0}^{(p-k_{i}) \times k_{i} }  & \mathbf{0}^{(p-k_{i}) \times (k-k_{i}) }  \end{bmatrix}$ ,
\end{itemize}
where $\mathbf{Q}_{i}$ is an $p \times p$ orthogonal matrix, $\mathbf{P}_{i}$ is an $k \times k$ permutation matrix, $\mathbf{R}_{i}$ is an $k_{i} \times k_{i}$ nonsingular upper triangular matrix (with diagonal elements of decreasing absolute magnitude) and $\mathbf{S}_{i}$ an $k_{i} \times (k-k_{i})$ full matrix, which is vacuous if $k_{i} = k$.

In all cases, compute an $p \times k$  matrix $\mathbf{O}_{i}$  with orthonormal columns as the first $k$ columns of $\mathbf{Q}_{i}$ (i.e., such that $\emph{ran}( \mathbf{A}_{i} ) \subset \emph{ran}( \mathbf{O}_{i} )$ if $k_{i} < k$ and  $\emph{ran}( \mathbf{A}_{i} ) = \emph{ran}( \mathbf{O}_{i} )$ if $k_{i} = k$) and set
\begin{itemize}
\item[] $\mathbf{A}_{i} = \mathbf{O}_{i}$ .
\end{itemize}
This optional orthogonalization step is a safe-guard as the condition $k_{i} = k$ is a necessary condition for the differentiability of $\psi(.)$ at a point $\mathbf{A}_{i}$ and also to limit the occurrence of overflows and underflows in the next steps by enforcing that the matrix variable $\mathbf{A}_{i} \in \mathbb{O}^{p \times k}$.
\item[\textbf{(1)}]  Determine (implicitly) the block diagonal matrix
\begin{itemize}
\item[] $\mathbf{F}(\mathbf{a}_{i}) = \emph{diag}\big( \emph{vec}( \sqrt{\mathbf{W}} ) \big)  \big(  \mathbf{I}_n  \otimes \mathbf{A}_{i}  \big)$ ,
\end{itemize}
where $\mathbf{a}_{i} =  \emph{vec}( \mathbf{A}_{i}^{T} )$.
\item[\textbf{(2)}] Compute (implicitly) a QRCP of $\mathbf{F}(\mathbf{a}_{i})$ to determine  ${r}_{\mathbf{F}( \mathbf{a}_{i} )}  = \emph{rank}( \mathbf{F}(  \mathbf{a}_{i} ) )$, $\mathbf{P}^{\bot}_{\mathbf{F}(  \mathbf{a}_{i} ) }$ and  $\mathbf{F}(\mathbf{a}_{i})^-$ (see equations~\eqref{eq:ginv_proj_ortho} and \eqref{eq:sginv_qrcp}). Note that  $\mathbf{P}^{\bot}_{\mathbf{F}(  \mathbf{a}_{i} ) }$ and $\mathbf{F}(\mathbf{a}_{i})^-$ are also block diagonal matrices and that the Newton and quasi-Newton algorithms assume that ${r}_{\mathbf{F}( \mathbf{a}_{i} )} = k.p$; see the derivation of the Hessian matrix  $\nabla^2 \psi( \mathbf{a}_{i} )$ in Subsection~\ref{hess:box} for more details. Note that, optionally, the validity of this hypothesis can be checked here in output of the QRCP of $\mathbf{F}(\mathbf{a}_{i})$.

\item[\textbf{(3)}] Solve the block diagonal linear least-squares problem
\begin{itemize}
\item[] $\mathbf{b}_{i} = \text{Arg}\min_{\mathbf{b}\in\mathbb{R}^{k.n}} \,   \Vert \mathbf{x} - \mathbf{F}(\mathbf{a}_{i})\mathbf{b} \Vert^{2}_{2}$ ,
\end{itemize}
e.g., compute $\mathbf{b}_{i} = \mathbf{F}(\mathbf{a}_{i})^{-} \mathbf{x}$.

\item[\textbf{(4)}]  Determine and set:
\begin{itemize}
\item[] $\mathbf{r}(\mathbf{a}_{i}) = \mathbf{P}^{\bot}_{\mathbf{F}(  \mathbf{a}_{i} ) } \mathbf{x}$ $\lbrace \text{current residual vector} \rbrace$
\item[] $\psi(\mathbf{a}_{i} ) = \frac{1}{2} \Vert \mathbf{r}(\mathbf{a}_{i}) \Vert^{2}_{2}$ $\lbrace \text{current value of the cost function }  \rbrace$
\item[] $\nabla \psi( \mathbf{a}_{i} ) =  \mathbf{G}(\mathbf{b}_{i})^{T} \mathbf{G}(\mathbf{b}_{i})\mathbf{a}_{i} - \mathbf{G}(\mathbf{b}_{i})^{T} \mathbf{z}$ $\lbrace \text{see Theorems~\ref{theo4.3:box} and~\ref{theo5.7:box}} \rbrace$
\item[] $j = 0$ $\lbrace \text{initialize counter for the ridge scaling subiterations} \rbrace$
\end{itemize}
Note that the steps \textbf{(1)} to \textbf{(4)} above can be very easily parallelized using the block diagonal structure of $\mathbf{F}(\mathbf{a}_{i})$.
\item[\textbf{(5)}]   Check for convergence. Relevant convergence criteria in the algorithms are of the form:
\begin{itemize}
\item $\Vert \nabla \psi(\mathbf{a}_{i} ) \Vert_{2} \le \varepsilon_{1}$
\item $\Vert \mathbf{a}_{i} - \mathbf{a}_{i-1} \Vert_{2} \le \varepsilon_{2} ( \varepsilon_{2} + \Vert \mathbf{a}_{i} \Vert_{2})$  $\lbrace$if $ i \ne 1\rbrace$

If step \textbf{(0)} is used, this convergence condition can be simplified as:

$\Vert \mathbf{a}_{i} - \mathbf{a}_{i-1} \Vert_{2} \le \varepsilon_{2}  \Vert \mathbf{a}_{i} \Vert_{2} = \varepsilon_{2} \sqrt{k}$
\item $\vert \psi(\mathbf{a}_{i-1} ) - \psi(\mathbf{a}_{i} ) \vert \le \varepsilon_{3} ( \varepsilon_{3} +  \psi(\mathbf{a}_{i} ) )$  $\lbrace$if $ i \ne 1\rbrace$
\item $ i \ge i_{max}$ $\lbrace \text{e.g., give up if the number of iterations is too large} \rbrace$
\end{itemize}
where $\varepsilon_{1}, \varepsilon_{2}, \varepsilon_{3}$ and $i_{max}$ are constants chosen by the user.

\textbf{Exit if convergence}. \textbf{Otherwise, go to} step \textbf{(6)}
\item[\textbf{(6)}] Compute the Newton correction vector $d\mathbf{a}_{n}$.
\begin{enumerate}
\item[\textbf{(6.1)}] To this end, first compute Hessian matrix or its two-term approximation:
\begin{itemize}
\item[] $\mathbf{H}_{i}^{1} = \mathbf{M}(\mathbf{a}_{i})^{T} \mathbf{M}(\mathbf{a}_{i})$
\item[] $\mathbf{H}_{i}^{2} = \mathbf{L}(\mathbf{a}_{i})^{T} \mathbf{L}(\mathbf{a}_{i})$
\item[] $\mathbf{H}_{i}^{3} = \mathbf{U}(\mathbf{a}_{i})^{T} \mathbf{L} (\mathbf{a}_{i}) + \mathbf{L}(\mathbf{a}_{i})^{T} \mathbf{U}(\mathbf{a}_{i})$ $\lbrace$only if a full Newton step is wanted$\rbrace$
\item[] $\mathbf{H}_{i} =
    \begin{cases}
        \mathbf{H}_{i}^{1} - \mathbf{H}_{i}^{2} + \mathbf{H}_{i}^{3}   & \lbrace \text{for a full Newton step}\rbrace \\
        \mathbf{H}_{i}^{1} - \mathbf{H}_{i}^{2}                                   & \lbrace \text{for a quasi-Newton step} \rbrace
     \end{cases} 
$
\item[] $\mathbf{H}_{i} = \mathbf{H}_{i} + \mathbf{N}_{i} \mathbf{N}_{i}^{T}$
\end{itemize}
where the columns of $\mathbf{N}_{i} = \mathbf{K}_{(p,k)} ( \mathbf{I}_{k} \otimes \mathbf{A}_{i} )$ are a (orthonormal) basis of  $\emph{null}\big( \mathbf{M}(\mathbf{a}_{i}) + \mathbf{L}(\mathbf{a}_{i}) \big) = \emph{null}( \mathbf{M}(\mathbf{a}_{i}) )$, see Corollary~\ref{corol5.6:box}.
\item[\textbf{(6.2)}] Check first diagonal elements of $\mathbf{H}_{i}$:
\begin{itemize}
\item[] $h_{min} = \emph{min}_{ l = 1, \cdots, k.p} \lbrack \mathbf{H}_{i} \rbrack_{ll}$
\item[] \textbf{If} $h_{min}  < 0$ \textbf{then} $\lambda = \lambda - h_{min}$ $\lbrace \text{scale up the ridge parameter} \rbrace$
\end{itemize}
\item[\textbf{(6.3)}] \textbf{Do while} $\mathbf{H}_{i} +  \lambda \mathbf{I}_{k.p}$ is not positive definite $\lbrace \text{use Cholesky factorization} \rbrace$
\begin{itemize}
\item[] $j = j + 1$
\item[] $\lambda = \nu.\lambda$ $\lbrace \text{scale up the ridge parameter} \rbrace$
\item[] $\nu = 2.\nu$ $\lbrace \text{increase the growth factor of the ridge parameter} \rbrace$
\end{itemize}
\item[\textbf{(6.4)}] 
\begin{description}
\item[Newton step:] get Newton step as the solution of the positive definite linear system
\begin{equation*}
\big(  \mathbf{H}_{i}  +  \lambda \mathbf{I}_{k.p} \big) d\mathbf{a}_{n} =  - \nabla \psi( \mathbf{a}_{i} ) = \mathbf{M}(\mathbf{a}_{i})^{T} \mathbf{r}(\mathbf{a}_{i})
\end{equation*}
\end{description}
\end{enumerate}
\item[\textbf{(7)}]  Compute next iterate $\mathbf{a}_{i+1} = \emph{vec}( \mathbf{A}_{i+1}^{T} )$ such that $\psi( \mathbf{a}_{i+1} ) < \psi( \mathbf{a}_{i} )$ in order to obtain global convergence.
\begin{enumerate}
\item[\textbf{(7.1)}]  To this end, first compute
\begin{itemize}
\item[] $\psi( \mathbf{a}_{i} +  d\mathbf{a}_{n} ) = \frac{1}{2} \Vert \mathbf{P}^{\bot}_{\mathbf{F}(  \mathbf{a}_{i} +  d\mathbf{a}_{n} ) } \mathbf{x} \Vert^{2}_{2}$,
\end{itemize}
using (implicitly) a QRCP of the block diagonal matrix $\mathbf{F}(  \mathbf{a}_{i} +  d\mathbf{a}_{n} )$, and the gain factor
\begin{itemize}
\item[] $\rho =  \frac{ \psi( \mathbf{a}_{i}) - \psi( \mathbf{a}_{i} + d\mathbf{a}_{n} )}{ N(\mathbf{0}^{k.p}) -  N( d\mathbf{a}_{n} ) } = \frac{ \psi( \mathbf{a}_{i}) - \psi( \mathbf{a}_{i} + d\mathbf{a}_{n} )}{  \frac{1}{2}   \big(  \Vert \mathbf{N}_{i}^{T}  d\mathbf{a}_{n} \Vert^{2}_{2} + \lambda \Vert d\mathbf{a}_{n} \Vert^{2}_{2} - d\mathbf{a}_{n}^{T} \nabla \psi( \mathbf{a}_{i} )  \big) }$
\end{itemize}
\item[\textbf{(7.2)}]  \textbf{If}  $\rho > 0$ \textbf{then} $\lbrace \text{step acceptable} \rbrace$
\begin{itemize}
\item[] $\lambda = \lambda .\emph{max}\big( \frac{1}{3}, 1-(2.\rho-1)^3 \big)$ $\lbrace \text{scale down the ridge parameter} \rbrace$
\item[] $\nu = 2$ $\lbrace \text{reinitialize the growth factor of the ridge parameter} \rbrace$
\end{itemize}
\item[\textbf{(7.3)}]  \textbf{Else} $\lbrace \text{step rejected} \rbrace$
\begin{itemize}
\item[] $j = j + 1$
\item[] $\lambda = \nu.\lambda$  $\lbrace \text{scale up the ridge parameter} \rbrace$
\item[] $\nu = 2.\nu$ $\lbrace \text{increase the growth factor of the ridge parameter} \rbrace$
\item[] \textbf{If} $j  \le j_{max}$ \textbf{go to} step  \textbf{(6.3)} $\lbrace \text{recompute } d\mathbf{a}_{n} \text{ with inflated diagonal of  } \mathbf{H}_{i} \rbrace$
\end{itemize}
\item[\textbf{(7.4)}]  Increment $\mathbf{a}_{i}$:

$\mathbf{a}_{i+1} = \mathbf{a}_{i} +  d\mathbf{a}_{n} \lbrace \text{compute new iterate} \rbrace $
\end{enumerate}
\end{enumerate}
\textbf{End do}

As in the Gauss-Newton or Levenberg-Marquardt algorithms, the computations in the above Newton algorithms are terminated either when one or several of the convergence criteria listed in step \textbf{(5)} are satisfied, or when the iteration count exceeds the predetermined number $i_{max}$.

In both the Newton algorithms~\eqref{n_alg2:box} and~\eqref{n_alg3:box}, the Marquardt parameter $\lambda$ is taken in the interval $\lbrack  10^{-8}  \,  1 \rbrack$ and a small value of $\lambda$ is selected if we believe that $\mathbf{A}_{1}$ is close to a solution (say $\lambda = 10^{-6}$). Otherwise, we can use $\lambda=10^{-3}$ or $10^{-4}$, or even $1$. The algorithms are not very sensitive to this initial choice of  $\lambda$ as this parameter is quickly updated during the iterations in both Newton algorithms~\eqref{n_alg2:box} and~\eqref{n_alg3:box}. Version~\eqref{n_alg3:box} of the Newton algorithms also uses a growth factor $\nu$ for the ridge parameter, which is  initialized to $2$ at the start of the algorithm and reinitialized to this initial value in step \textbf{(7.2)} when a Newton step is successful.

Note that it may happens during some iterations that the exact Hessian matrix $\mathbf{H}_i$ or its two-term approximation $\bar{\mathbf{H}}_i$ may become not positive definite in which case the computed Newton correction vector $d\mathbf{a}_{n}$ is not in a descent direction for $\psi(.)$. The Newton algorithms~\eqref{n_alg1:box} overcome this difficulty by using immediately a Gauss-Newton step, more precisely a Kaufman Gauss-Newton step, which is always in a descent direction if $\mathbf{a}_i$ is not a stationary point as demonstrated in Corollary~\ref{corol5.7:box}. On the other hand, the Newton algorithms~\eqref{n_alg2:box} and~\eqref{n_alg3:box} use the simple strategy of adding a multiple of the identity to the Hessian matrix until this modified Hessian matrix becomes positive definite; see step \textbf{(6.3)} in these algorithms. The drawback in this simple approach is that each time we add a multiple of the identity to the Hessian matrix, a new Cholesky factorization of $\mathbf{H}_{i}  +  \lambda \mathbf{I}_{k.p}$ must be attempted. This can become very expensive if the process is repeated many times across the iterations. This explains why $\lambda$ is multiplied by a value as large as $10$ in step \textbf{(6.3)} of the Newton algorithm~\eqref{n_alg2:box} and that the updating strategy of $\lambda$ in step \textbf{(6.3)} of  the Newton algorithm~\eqref{n_alg3:box} is modified so that consecutive failures of getting a positive definite modified Hessian matrix give a very fast growth of $\lambda$. The same modified strategy with a faster growth of $\lambda$ is also used in step  \textbf{(7.3)} of the Newton algorithms~\eqref{n_alg3:box} when the gain factor is negative and a trial Newton step is rejected. However, in both Newton algorithms~\eqref{n_alg2:box} and~\eqref{n_alg3:box}, if the number of ridge scaling subiterations in steps \textbf{(6.3)} and \textbf{(7.2)} (or \textbf{(7.3)} for Newton algorithms~\eqref{n_alg3:box}) become too important, we stop this ridge scaling process at the end of step \textbf{(7.2)} (or \textbf{(7.3)} for Newton algorithms~\eqref{n_alg3:box})  and we try to decrease the cost function $\psi(.)$ through the main loop of the algorithms instead.

\subsection{Variable projection hybrid algorithms} \label{vp_h_alg:box}

In the previous subsections, we have presented a large variety of variable projection WLRA solvers based on Gauss-Newton, Levenberg-Marquardt and Newton algorithms, which can all be considered as second-order or pseudo second-order methods. The complexity of these algorithms is relatively high, especially when $p \approx n$, even if $\emph{min}(p,n) \gg k$, as these algorithms have to solve a large linear system or a tall and skinny linear least-squares problem at each iteration as discussed in the previous subsections. Furthermore, the pre-processing of the Hessian matrix and its different approximations in a normal-equation approach or of the coefficient matrix of the linear least-squares problem in a QR approach for the Gauss-Newton and Levenberg-Marquardt methods is also expensive and has a high memory footprint~\cite{C2008b}\cite{OYD2011}\cite{HF2015}\cite{HZF2017}. In the previous section, we have provided different parallel techniques to reduce both the computing time and the memory storage requirements for this expensive pre-processing step included in all these  variable projection (pseudo) second-order algorithms.

However, all these (pseudo) second-order algorithms have a much higher complexity than the block ALS method (e.g., NIPALS, presented in Section~\ref{nipals:box}) or  other first-order algorithms based on majorization or  Expectation-Maximization (EM) methods, steepest, conjugate or stochastic gradient descent (or combinations of some of them), which have been proposed in the literature  to solve the WLRA problem, since  most of them are based on relatively inexpensive iterative algorithms that monotonically improve the function value by sequential repetition of local optimizations as the ALS method~\cite{K1997}\cite{K2002}\cite{JHJ2009}\cite{IR2010}\cite{HMLZ2015}\cite{MHT2010}\cite{RSW2016}\cite{BWZ2019}\cite{BL2020}\cite{TH2021}\cite{OUV2023}. Moreover, the block ALS method greatly reduces the memory footprint requirements and can also be parallelized easily and very efficiently taking into account the block diagonal structures of the matrices $\mathbf{F}(\mathbf{a})$ and $\mathbf{G}(\mathbf{b})$ as explained in Section~\ref{nipals:box}. On the other hand, it is known that the block ALS method and its variants fail frequently to converge to an acceptable optimal solution for difficult WLRA problems without the use of a proper regularization, and are prone to flattening, especially when the percentage of missing data or/and the level of noise are high~\cite{BF2005}\cite{C2008b}\cite{D2011}\cite{OYD2011}\cite{HF2015}\cite{HZF2017}.

These different features suggest that hybrid methods combining the fast block ALS algorithm with any of the (pseudo) second-order variable projection algorithms detailed in the last subsections can performed much better than any of the previous individual algorithms as first suggested by~\cite{BF2005} and confirmed by Chen~\cite{C2008b} in the specific case of variable projection Newton algorithms. In this way, we can benefit of the remarkable efficiency of the block ALS method in terms of speed and, at the same time, of the good convergence ability of the (pseudo) second-order methods to reduce drastically the computing time without sacrificing the overall performance of the (pseudo) second-order WLRA solvers.

As an illustration, one simple, but still very efficient, choice of such hybrid methods could be to add an additional step consisting of a fixed number of iterations (say between $2$ and $20$) of the block ALS algorithm before step \textbf{(0)} in all our (pseudo) second-order variable projection algorithms. Preliminary tests (not shown) suggest that this simple modification is always very beneficial in terms of speed without impairing the global convergence performance of the (pseudo) second-order variable projection WLRA solvers discussed here. Of course, many variations other than this simple hybrid scheme are possible~\cite{BF2005}\cite{C2008b}.

\section{Conclusions and discussion} \label{conclu:box}

In this monograph, we consider the difficult WLRA problem, an extension of the well-known matrix completion problem~\cite{DR2016}. The WLRA problem is NP-hard and has no closed solution in general~\cite{GG2011}, but has an increasing number of important applications in practice.

We survey many different approaches which have been used to solve it, with a particular focus on variable projection second-order methods~\cite{R1974}\cite{GP1973}\cite{RW1980}\cite{B2009}\cite{OYD2011}\cite{OR2013}\cite{HF2015}. A large variety of low-complexity first-order methods have been already proposed in the past to solve the WLRA problem~\cite{WWY2012}\cite{V2013}\cite{BA2015}\cite{HMLZ2015}\cite{BL2020} \cite{TH2021}\cite{OUV2023}, but only a few second-order methods, as these second-order methods have a very high per-iteration complexity, which preclude their use for very large datasets commonly found in recent applications. However, variable projection second-order methods perform better than first-order methods for badly conditioned WLRA problems, which are very common, and these more costly methods are thus still of interest in this context~\cite{OYD2011}\cite{HF2015}.

First, we review in detail the connections between manifold, variety, factorization and variable projection formulations of the WLRA problem, which are most often treated as disconnected approaches in the literature and establish relationships between (local) minima, first- and second-order critical points of the objective functions used in the these different formulations of the WLRA problem. These results are an illustration and slight extension in the context of the WLRA problem of recent results presented in~\cite{SU2015}\cite{HLU2019}\cite{HLB2020}\cite{LLZ2024}\cite{LKB2025} about the near equivalences of first- and second-order critical points of the objective functions in nonconvex factorization, manifold and variety formulations in more general low-rank matrix optimization.

Second, we provide an extended and original overview of the variable projection formulation of the WLRA problem both from theoretical and algorithmic perspectives. In particular, we study in detail the non-smoothness of the variable projection cost function of the WLRA problem when some weights are equal to zero (e.g., in presence of missing values) and when this cost function is not regularized. We characterize precisely its discontinuities generalizing the preliminary investigations of Dai et al.~\cite{DMK2011}\cite{DKM2012}  on this topic. These points of discontinuity form barrier sets in the feasible space of solutions, which prevent low-complexity algorithms like gradient descent or alternating least squares to converge to first-order critical points for some WLRA problems when missing values are present. Up to now, most of variable projection algorithms proposed in the literature for solving the WLRA problem simply "ignore" these discontinuities when missing values are present or use a regularized objective function to eliminate them. However, such regularization may degrade significantly the compression quality of the computed low-rank matrix approximation. It is thus interesting to explore whether  it is possible to detect if we move in front of these discontinuity points in the unregularized variable projection methods when missing values are present and to find ways of escaping from them during the course of the computations as was discussed in Dai et al.~\cite{DMK2011} for a gradient descent algorithm. Finally, it is also worth to explore in more depth the landscape connections between these discontinuity points in the variable projection formulation and the corresponding points in the manifold and factorization formulations of the WLRA problem. As an illustration, the matrix incoherency hypotheses, which are often used to prove the convergence of  gradient descent or alternating least-squares for matrix completion problems~\cite{VMS2016} mainly lead to the exclusion of the above barrier sets.

Next, we derive new formulae for the variable projection gradient vector, and Jacobian and Hessian matrices, which are pivotal in all the second-order variable projection methods. These new formulae also allow a better understanding of the geometric landscape of the objective function associated with the variable projection formulation. In particular, they allow to characterize precisely the systematic rank degeneracy of the variable projection Jacobian matrix and also the singularities of the variable projection Hessian at first-order stationarity points of the variable projection objective function. Furthermore, these new formulae allow us to demonstrate that most of the first- and second-order variable projection methods, which can be used to solve the WLRA problem, can be viewed as as Riemannian optimization methods operating on the Grassmann manifold when the variable projection objective function is smooth over all the points of the Grassmannian.

Our new formulae for the variable projection gradient vector, Jacobian and Hessian matrices also allow us to formulate more accurate and robust variable projection second-order algorithms based on stable orthogonal kernels to tackle the systematic rank degeneracy of the Jacobian, the singularity of the Hessian at first-order stationarity points of the variable projection functional and, finally, the ill-conditioning and indefinite nature of the Hessian at points arbitrarily closed to local minima of this functional. These singularities and instabilities are inescapable here as the (local) minima of the variable projection objective function are always non-isolated  and form a differentiable submanifold of the ambient linear space around around each of the minima under some regularity hypotheses of the variable projection functional as demonstrated in Subsection~\ref{hess:box}.

From an algorithmic point of view, these formulae also allow us to formulate more efficient variable projection second-order algorithms to tackle WLRA problems of larger dimensions, which are now the rule in many applications, despite these variable projection second-order algorithms have still a high per-iteration complexity. In particular, we improve significantly the scalability and efficiency of the proposed variable projection second-order algorithms compared to previous studies by taking better into account the sparseness of the  different matrix variables involved in each iteration of the algorithms and by using large-scale parallelization techniques and highly optimized BLAS3 kernels in the sensible parts of the algorithms. However, some parts of the variable projection algorithms still do not explicitly take into account the sparsity of the matrix variables. This concerns especially the recursive and parallel implementations of the QR decomposition of the the tall and skinny matrices generated by blocks and used in the iterations of both the Gauss-Newton and Levenberg-Marquardt algorithms. Devising more efficient recursive and parallel QR decompositions of these sparse, tall and skinny matrices is thus an important topic for future research, including possibly the use of very fast randomized QR methods~\cite{MDME2023}.

Finally, here, we have assumed that the rank of the low-rank matrix approximation we are seeking is given or bounded before hand by the user. Of course, in practice, this rank is often unknown. An interesting continuation of this work would thus be to devise efficient variable projection algorithms specifically designed to adaptively select or change the rank of the low-rank matrix solution of the WLRA problem during the computations as was done for Riemannian descent methods on low-rank matrix varieties in~\cite{UV2014}\cite{UV2015}. Such rank-adaptive optimization strategies in which local minima of smaller rank are used as starting points for improved approximation with a larger rank will be very useful to select the best rank in practical applications and can also lead to improve efficiency and accuracy in the computations of larger, but fixed low-rank matrix approximation, especially for ill-conditioned WLRA problems. 
More generally, selecting accurate, but cheap, initial low-rank matrix approximations as a first guess of the costly variable projection second-order methods described here is another useful and important continuation of this work.

Obviously, the variable projection second-order algorithms described here must be compared to the many other state-of-the-art first-order methods already proposed in the literature for solving the WLRA problem~\cite{WWY2012}\cite{V2013}\cite{BA2015}\cite{HMLZ2015}\cite{BL2020} \cite{TH2021}\cite{OUV2023} in comprehensive benchmark experiments. This benchmark must be based on both synthetic and real datasets and not only restricted to the matrix completion problem as in many past experiments. In particular, we expect that the variable projection second-order algorithms grow more efficient, robust and accurate than gradient-based or alternating least-squares methods as the amount of missing values (e.g., zero weights) increase in the WLRA problem, but this must be objectively validated in comprehensive experiments. Of course, if speed is the priority, first-order (e.g. gradient descent or alternating least squares) or hybrid algorithms as described in Subsection~\ref{vp_h_alg:box} will be the methods of choice.

Fortran90 codes with OpenMP and BLAS supports, for the variable projection second-order algorithms described here, will be later available in the open source STATPACK library available at:

 \url{https://pagesperso.locean-ipsl.upmc.fr/terray/statpack2.3/index.html}.

\end{document}